\documentclass[12pt]{amsbook}
\usepackage{amssymb}
\usepackage[all]{xy}

\usepackage[colorlinks=true,linkcolor=magenta,citecolor=magenta]{hyperref}
\makeindex

\newtheorem{theorem}{Theorem}[chapter]
\newtheorem{answer}[theorem]{Answer}
\newtheorem{claim}[theorem]{Claim}
\newtheorem{conclusion}[theorem]{Conclusion}
\newtheorem{conjecture}[theorem]{Conjecture}
\newtheorem{definition}[theorem]{Definition}
\newtheorem{exercise}[theorem]{Exercise}
\newtheorem{problem}[theorem]{Problem}
\newtheorem{proposition}[theorem]{Proposition}
\newtheorem{question}[theorem]{Question}
\newtheorem{solution}[theorem]{Solution}
\newtheorem{speculation}[theorem]{Speculation}

\begin{document}

\title{Quantum permutation groups}

\author{Teo Banica}
\address{Department of Mathematics, University of Cergy-Pontoise, F-95000 Cergy-Pontoise, France. {\tt teo.banica@gmail.com}}


\begin{abstract}
The permutation group $S_N$ has a quantum analogue $S_N^+$, which is infinite at $N\geq4$. We review the known facts regarding $S_N^+$, and notably its easiness property, Weingarten calculus, and the isomorphism $S_4^+=SO_3^{-1}$ and its consequences. We discuss then the structure of the closed subgroups $G\subset S_N^+$, and notably of the quantum symmetry groups of finite graphs $G^+(X)\subset S_N^+$, with particular attention to the quantum reflection groups $H_N^{s+}$. We also discuss, more generally, the quantum symmetry groups $S_Z^+$ of the finite quantum spaces $Z$, and their closed subgroups $G\subset S_Z^+$, with particular attention to the quantum graph case, and to quantum reflection groups.
\end{abstract}

\maketitle

\chapter*{Preface}

One of the most puzzling discoveries in quantum algebra, going back to work of Wang from the late 90s, in answer to a question of Connes, is that the set $\{1,\ldots,N\}$ with $N\geq4$ has an infinity of quantum permutations. At the first glance, this looks as one of these physicists' crazy things, which might be worth attention or not. But please don't go away, and listen to what I have so say. Yes, all this is related to physics. But the mathematics behind it is extremely simple, and worth some attention.

\bigskip

Let us first look at the symmetric group $S_N$. When regarding it as group of permutations of the $N$ coordinate axes of $\mathbb R^N$, so as subgroup $S_N\subset O_N$, the standard coordinates $u_{ij}\in C(S_N)$ are given by a very simple formula, $u_{ij}(\sigma)=\delta_{\sigma(j)i}$. It follows that these coordinates $u_{ij}\in C(S_N)$ form a matrix $u=(u_{ij})$ which is ``magic'', in the sense that its entries are projections, $p^2=p^*=p$, which sum up to 1 on each row and each column. Moreover, by Stone-Weierstrass we have $C(S_N)=<u_{ij}>$, and with a bit more work, by using the Gelfand theorem, we can see that $C(S_N)$ is isomorphic to the universal commutative $C^*$-algebra generated by the entries of a $N\times N$ magic matrix.

\bigskip

This suggests looking at the universal $C^*$-algebra $C(S_N^+)$ generated by the entries of a $N\times N$ magic matrix. In analogy with what happens for $C(S_N)$, this algebra has a comultiplication $\Delta$, a counit $\varepsilon$ and an antipode $S$, so according to general results of Woronowicz, its spectrum $S_N^+$ is a compact quantum group, called quantum permutation group. And the point is that the inclusion $S_N\subset S_N^+$ is not an isomorphism at $N\geq4$, because diagonally joining magic matrices of size $\geq2$ shows that $S_N^+$ is a non-classical, infinite compact quantum group, substantially bigger than $S_N$.

\bigskip

Summarizing, some interesting mathematics going on here, and with this digested, the first thought goes to physics. Can such beasts be of help in connection with statistical mechanics, along the lines suggested by Jones? What about quarks and the Standard Model, along the lines suggested by Connes? What about statistical mechanics and quantum physics alike, via random matrices, and freeness in the sense of Voiculescu? And also, in tune with our times, what about quantum information?

\bigskip

These questions are all old and difficult, going back to Wang's discovery of $S_N^+$ in the late 90s, and the few years that followed. So, perhaps more modestly, we should start with some pure mathematics. The point indeed is that $S_N$ and its subgroups $G\subset S_N$ have a lot of interesting mathematics, so before anything, we should understand what the analogous theory of $S_N^+$ and of its closed quantum subgroups $G\subset S_N^+$ is. And with a bit of luck, we will get in this way precisely into the mathematics of Connes, Jones, Voiculescu, putting us on the right track for doing some physics afterwards.

\bigskip

The present book is precisely about this, the mathematics of $S_N^+$ and of its closed subgroups $G\subset S_N^+$, with physics motivations in mind. We will be interested in mathematics in a large sense, mixing algebra, geometry, analysis and probability. At the level of subgroups, we will mostly insist on the quantum reflection groups, which are quite fundamental objects, mathematically speaking, and which are expected as well to be useful in physics. Finally, we will discuss also the more general case of the quantum permutation groups $S_Z^+$ of the arbitrary ``finite quantum spaces'' $Z$, and their closed subgroups $G\subset S_Z^+$, with special attention to the quantum reflection subgroups.

\bigskip

All in all, many things to be discussed. We will assume some familiarity with basic graduate level mathematics, such as operator algebras and quantum groups. In case you are not familiar with this, have a copy of my quantum group book \cite{ba4} handy. However, the present book is for the most self-contained, up to undergraduate mathematics, and if you're not afraid of theorems coming with very short proofs, it is for you as such. 

\bigskip

Finally, regarding physics, there will be not much of it in this book, which is meant to be a purely mathematical text. However, the relation with the work of Connes, Jones, Voiculescu, potentially leading to physics, will be carefully explained. In short, you will learn in this way everything that is needed for doing good physics afterwards. With this applying to myself too, my lattice model book \cite{ba5} being long overdue. 

\bigskip

This book is heavily based on a number of research papers on quantum permutations and reflections, and I am particularly grateful to Julien Bichon, Beno\^\i t Collins and Steve Curran, for substantial joint work on the subject. Many thanks go as well to my cats, for precious support and advice, during the preparation of the present book.

\bigskip

\

{\em Cergy, August 2024}

\smallskip

{\em Teo Banica}

\baselineskip=15.95pt
\tableofcontents
\baselineskip=14pt

\part{Quantum permutations}

\ \vskip50mm

\begin{center}
{\em Carry me, caravan, take me away

Take me to Portugal, take me to Spain

Andalucia with fields full of grain

I have to see you again and again}
\end{center}

\chapter{Quantum permutations}

\section*{1a. Quantum spaces}

Welcome to quantum permutations. In this first part of the present book we discuss the construction and basic properties of the quantum permutation groups $S_N^+$, and of their generalizations $S_Z^+$, with $Z$ being a finite quantum space. The story here involves the foundational papers of Woronowicz \cite{wo1}, \cite{wo2}, from the end of the 80s, then the key papers of Wang \cite{wa1}, \cite{wa2}, from the mid 90s, then my own papers \cite{ba1}, \cite{ba2}, \cite{ba3}, from the late 90s and early 00s, and finally some more specialized papers from the mid and late 00s, including \cite{bb2}, \cite{bbs}, \cite{bc1}, \cite{bc2}, \cite{bcu}, \cite{bsp}, containing a few fundamentals too.

\bigskip

In short, cavern man mathematics from about 20 years ago, but lots of things to be learned. We will provide here 100 pages on the subject, with a decent presentation of what is known about $S_N^+$ and $S_X^+$, of fundamental type, coming in the form of theorems accompanied by short proofs. For further details on all this, you have my graduate textbook on quantum groups \cite{ba4}, along with the original papers cited above.

\bigskip

Getting started now, at the beginning of everything, we have:

\begin{question}[Connes]
What is a quantum permutation group?
\end{question}

This question is more tricky than it might seem. For solving it you need a good formalism of quantum groups, and there is a bewildering number of choices here, with most of these formalisms leading nowhere, in connection with the above question. So, we are into philosophy, and for truly getting started, we have to go back in time, with:

\begin{question}[Heisenberg]
What is a quantum space?
\end{question}

Regarding this latter question, there are as many answers as quantum physicists, starting with Heisenberg himself in the early 1920s, then Schr\"odinger and Dirac short after, with each coming with his own answer to the question. Not to forget Einstein, who labeled all these solutions as ``nice, but probably fundamentally wrong''. 

\bigskip

In short, we are now into controversy, and a look at more modern physics does not help much, with the controversy basically growing instead of diminishing, over the time. So, in the lack of a good answer, let us take as starting point something nice and mathematical, rather agreed upon in the 1930s, coming from Dirac's work, namely:

\begin{answer}[von Neumann]
A quantum space is the dual of an operator algebra.
\end{answer}

Fast forward now to the 90s and to Connes' question, this remains something non-trivial, even when knowing what a quantum space is, and this for a myriad technical reasons. You have to work a bit on that question, try all sorts of things which do not work, until you hit the good answer. With this good answer being as follows:

\begin{answer}[Wang]
The quantum permutation group $S_N^+$ is the biggest compact quantum group acting on $\{1,\ldots,N\}$, by leaving the counting measure invariant.
\end{answer}

To be more precise, the idea is that $\{1,\ldots,N\}$ has all sorts of quantum permutations, and even when restricting the attention to the ``correct'' ones, namely those leaving invariant the counting measure, there is still an infinity of such quantum permutations, and the quantum group formed by this infinity of quantum permutations is compact.

\bigskip

This was for the story of the subject, very simplified, and as a final ingredient, two answers to two natural questions that you might have:

\bigskip

(1) Isn't the conclusion $|S_N^+|=\infty$ a bit too speculatory, not to say crazy? Certainly not, I would say, because in quantum mechanics particles do not have clear positions and speeds, and once you're deep into this viewpoint, ``think quantum'', a bit fuzzy about everything, why the set $\{1,\ldots,N\}$ not being allowed to have an infinity of quantum permutations, after all. So, no contradiction, philosophically speaking.

\bigskip

(2) Why was the theory of $S_N^+$ developed so late? Good question, and in answer, looking retrospectively, quantum groups and permutations should have been developed by von Neumann and Weyl, sometimes in the 1940s, perhaps with some help from Gelfand. But that never happened. As for the story after WW2, with mathematics, physics, and mankind in general: that was sex, drugs and rock and roll, forget about it.

\bigskip

Getting started now for good, we have the whole remainder of this chapter for understanding what Question 1.1 is about, and what its Answer 1.4 says. But before that, Question 1.2 and Answer 1.3 coming first. Leaving aside physics, we must first talk about operator algebras, and the starting definition here is as follows:

\index{Banach algebra}
\index{operator algebra}
\index{involutive algebra}

\begin{definition}
A $C^*$-algebra is a complex algebra $A$, having a norm $||.||$ making it a Banach algebra, and an involution $*$, related to the norm by the formula 
$$||aa^*||=||a||^2$$
which must hold for any $a\in A$.
\end{definition}

As a basic example, the algebra $M_N(\mathbb C)$ of the complex $N\times N$ matrices is a $C^*$-algebra, with the usual matrix norm and involution of matrices, namely:
$$||M||=\sup_{||x||=1}||Mx||\quad,\quad 
(M^*)_{ij}=\bar{M}_{ji}$$

More generally, any $*$-subalgebra $A\subset M_N(\mathbb C)$ is automatically closed, and so is a $C^*$-algebra. In fact, in finite dimensions, the situation is as follows:

\index{finite dimensional algebra}
\index{multimatrix algebra}

\begin{proposition}
The finite dimensional $C^*$-algebras are exactly the algebras
$$A=M_{n_1}(\mathbb C)\oplus\ldots\oplus M_{n_k}(\mathbb C)$$
with norm $||(a_1,\ldots,a_k)||=\sup_i||a_i||$, and involution $(a_1,\ldots,a_k)^*=(a_1^*,\ldots,a_k^*)$.
\end{proposition}

\begin{proof}
In one sense this is clear. In the other sense, this comes by splitting the unit of our algebra $A$ as a sum of central minimal projections, $1=p_1+\ldots+p_k$. Indeed, when doing so, each of the $*$-algebras $A_i=p_iAp_i$ follows to be a matrix algebra, $A_i\simeq M_{n_i}(\mathbb C)$, and this gives the direct sum decomposition in the statement.
\end{proof}

In general now, a main theoretical result about $C^*$-algebras, due to Gelfand, Naimark and Segal, and called GNS representation theorem, is as follows:

\index{GNS theorem}
\index{bounded operator}
\index{Hilbert space}
\index{operator algebra}
\index{linear operator}

\begin{theorem}
Given a complex Hilbert space $H$, finite dimensional or not, the algebra $B(H)$ of linear operators $T:H\to H$ which are bounded, in the sense that
$$||T||=\sup_{||x||=1}||Tx||$$
is finite, is a $C^*$-algebra, with the above norm, and with involution given by:
$$<Tx,y>=<x,T^*y>$$
More generally, and norm closed $*$-subalgebra of this full operator algebra
$$A\subset B(H)$$
is a $C^*$-algebra. Any $C^*$-algebra appears in this way, for a certain Hilbert space $H$.
\end{theorem}

\begin{proof}
There are several statements here, with the first ones being standard operator theory, and with the last one being the GNS theorem, the idea being as follows:

\medskip

(1) First of all, the full operator algebra $B(H)$ is a Banach algebra. Indeed, given a Cauchy sequence $\{T_n\}$ inside $B(H)$, we can set $Tx=\lim_{n\to\infty}T_nx$, for any $x\in H$. It is then routine to check that we have $T\in B(H)$, and that $T_n\to T$ in norm.

\medskip 

(2) Regarding the involution, the point is that we must have $<Tx,y>=<x,T^*y>$, for a certain vector $T^*y\in H$. But this can serve as a definition for $T^*$, and the fact that $T^*$ is indeed linear, and bounded, with the bound $||T^*||=||T||$, is routine. As for the formula $||TT^*||=||T||^2$, this is elementary as well, coming by double inequality.

\medskip

(3) Finally, the fact that any $C^*$-algebra appears as $A\subset B(H)$, for a certain Hilbert space $H$, is advanced. The idea is that each $a\in A$ acts on $A$ by multiplication, $T_a(b)=ab$. Thus, we are more or less led to the result, provided that we are able to convert our algebra $A$, regarded as a complex vector space, into a Hilbert space $H=L^2(A)$. But this latter conversion can be done, by taking some inspiration from abstract measure theory.
\end{proof}

As a third and last basic result about $C^*$-algebras, which will be of particular interest for us, we have the following well-known theorem of Gelfand:

\index{Gelfand theorem}
\index{commutative algebra}

\begin{theorem}
Given a compact space $X$, the algebra $C(X)$ of continuous functions $f:X\to\mathbb C$ is a $C^*$-algebra, with norm and involution as follows:
$$||f||=\sup_{x\in X}|f(x)|\quad,\quad 
f^*(x)=\overline{f(x)}$$
This algebra is commutative, and any commutative $C^*$-algebra $A$ is of this form, with $X=Spec(A)$ appearing as the space of Banach algebra characters $\chi:A\to\mathbb C$.
\end{theorem}

\begin{proof}
Once again, there are several statements here, some of them being trivial, and some of them being advanced, the idea being as follows:

\medskip

(1) First of all, the fact that $C(X)$ is indeed a Banach algebra is clear, because a uniform limit of continuous functions is continuous.

\medskip

(2) Regarding now for the formula $||ff^*||=||f||^2$, this is something trivial for functions, because on both sides we obtain $\sup_{x\in X}|f(x)|^2$.

\medskip

(3) Given a commutative $C^*$-algebra $A$, the character space $X=\{\chi:A\to\mathbb C\}$ is compact, and we have an evaluation morphism $ev:A\to C(X)$. 

\medskip

(4) The tricky point, which follows from basic spectral theory in Banach algebras, is to prove that $ev$ is indeed isometric. This gives the last assertion.
\end{proof}

In what follows, we will be mainly using Definition 1.5 and Theorem 1.8, as general framework. To be more precise, in view of Theorem 1.8, let us formulate:

\index{quantum space}
\index{compact quantum space}
\index{noncommutative space}

\begin{definition}
Given an arbitrary $C^*$-algebra $A$, we agree to write 
$$A=C(X)$$
and call the abstract space $X$ a compact quantum space.
\end{definition}

In other words, we can define the category of compact quantum spaces $X$ as being the category of the $C^*$-algebras $A$, with the arrows reversed. A morphism $f:X\to Y$ corresponds by definition to a morphism $\Phi:C(Y)\to C(X)$, a product of spaces $X\times Y$ corresponds by definition to a product of algebras $C(X)\otimes C(Y)$, and so on.

\bigskip

All this is of course quite speculative, and as a first result regarding these compact quantum spaces, coming from Proposition 1.6, we have:

\index{finite quantum space}

\begin{proposition}
The finite quantum spaces are exactly the disjoint unions of type
$$X=M_{n_1}\sqcup\ldots\sqcup M_{n_k}$$
where $M_n$ is the finite quantum space given by $C(M_n)=M_n(\mathbb C)$.
\end{proposition}

\begin{proof}
This is a reformulation of Proposition 1.6, by using the above philosophy. Indeed, for a compact quantum space $X$, coming from a $C^*$-algebra $A$ via the formula $A=C(X)$, being finite can only mean that the following number is finite:
$$|X|=\dim_\mathbb CA<\infty$$

Thus, by using Proposition 1.6, we are led to the conclusion that we must have:
$$C(X)=M_{n_1}(\mathbb C)\oplus\ldots\oplus M_{n_k}(\mathbb C)$$

But since direct sums of algebras $A$ correspond to disjoint unions of quantum spaces $X$, via the correspondence $A=C(X)$, this leads to the conclusion in the statement.
\end{proof}

This was for the basic theory of $C^*$-algebras, the idea being that we have some basic operator theory results, that can be further learned from any standard book, such as Blackadar \cite{bla}, and then we can talk about reformulations of these results in quantum space terms, by using Definition 1.9 and some basic common sense.

\bigskip

Finally, no discussion would be complete without a word about the von Neumann algebras. These are operator algebras of more advanced type, as follows:

\index{bicommutant}
\index{von Neumann algebra}
\index{weak topology}
\index{weak operator topology}
\index{bicommutant}
\index{commutative algebra}

\begin{theorem}
For a $*$-algebra $A\subset B(H)$ the following conditions are equivalent, and if they are satisfied, we say that $A$ is a von Neumann algebra:
\begin{enumerate}
\item $A$ is closed with respect to the weak topology, making each $T\to Tx$ continuous.

\item $A$ is equal to its algebraic bicommutant, $A=A''$, computed inside $B(H)$.
\end{enumerate}
As basic examples, we have the algebras $A=L^\infty(X)$, acting on $H=L^2(X)$. Such algebras are commutative, any any commutative von Neumann algebra is of this form.
\end{theorem}

\begin{proof}
There are several assertions here, the idea being as follows:

\medskip

(1) The equivalence $(1)\iff(2)$ is the well-known bicommutant theorem of von Neumann, which can be proved by using an amplification trick, $H\to\mathbb C^N\otimes H$.

\medskip

(2) Given a measured space $X$, we have indeed an emdedding $L^\infty(X)\subset B(L^2(X))$, with weakly closed image, given by $T_f:g\to fg$, as in the proof of the GNS theorem.

\medskip

(3) Given a commutative von Neumann algebra $A\subset B(H)$ we can write  $A=<T>$ with $T$ being a normal operator, and the Spectral Theorem gives $A\simeq L^\infty(X)$.
\end{proof}

In the context of a $C^*$-algebra representation $A\subset B(H)$ we can consider the weak closure, or bicommutant $A''\subset B(H)$, which is a von Neumann algebra. In the commutative case, $C(X)\subset B(L^2(X))$, the weak closure is $L^\infty(X)$. In general, we agree to write:
$$A''=L^\infty(X)$$

For more on all this, basic theory of the $C^*$-algebras and von Neumann algebras, we refer to any standard operator algebra book, such as Blackadar \cite{bla}.

\section*{1b. Quantum groups}

We are ready now to introduce the compact quantum groups. The axioms here, due to Woronowicz \cite{wo1}, and slightly modified for our present purposes, are as follows:

\index{Woronowicz algebra}
\index{comultiplication}
\index{counit}
\index{antipode}

\begin{definition}
A Woronowicz algebra is a $C^*$-algebra $A$, given with a unitary matrix $u\in M_N(A)$ whose coefficients generate $A$, such that the formulae
$$\Delta(u_{ij})=\sum_ku_{ik}\otimes u_{kj}\quad,\quad
\varepsilon(u_{ij})=\delta_{ij}\quad,\quad 
S(u_{ij})=u_{ji}^*$$
define morphisms of $C^*$-algebras $\Delta:A\to A\otimes A$, $\varepsilon:A\to\mathbb C$ and $S:A\to A^{opp}$, called comultiplication, counit and antipode. 
\end{definition}

In this definition the tensor product needed for $\Delta$ can be any $C^*$-algebra tensor product. In order to get rid of redundancies, coming from this and from amenability issues, we will divide everything by an equivalence relation, as follows:

\index{equivalence relation}
\index{amenability}

\begin{definition}
We agree to identify two Woronowicz algebras, $(A,u)=(B,v)$, when we have an isomorphism of $*$-algebras
$$<u_{ij}>\simeq<v_{ij}>$$
mapping standard coordinates to standard coordinates, $u_{ij}\to v_{ij}$.
\end{definition}

We say that $A$ is cocommutative when $\Sigma\Delta=\Delta$, where $\Sigma(a\otimes b)=b\otimes a$ is the flip. We have then the following key result, from \cite{wo1}, providing us with examples:

\index{cocommutative algebra}
\index{commutative algebra}
\index{compact Lie group}
\index{finitely generated group}

\begin{proposition}
The following are Woronowicz algebras, which are commutative, respectively cocommutative:
\begin{enumerate}
\item $C(G)$, with $G\subset U_N$ compact Lie group. Here the structural maps are:
$$\Delta(\varphi)=\big[(g,h)\to \varphi(gh)\big]\quad,\quad 
\varepsilon(\varphi)=\varphi(1)\quad,\quad
S(\varphi)=\big[g\to\varphi(g^{-1})\big]$$

\item $C^*(\Gamma)$, with $F_N\to\Gamma$ finitely generated group. Here the structural maps are:
$$\Delta(g)=g\otimes g\quad,\quad 
\varepsilon(g)=1\quad,\quad
S(g)=g^{-1}$$
\end{enumerate}
Moreover, we obtain in this way all the commutative/cocommutative algebras.
\end{proposition}

\begin{proof}
In both cases, we first have to exhibit a certain matrix $u$, and then prove that we have indeed a Woronowicz algebra. The constructions are as follows:

\medskip

(1) For the first assertion, we can use the matrix $u=(u_{ij})$ formed by the standard matrix coordinates of $G$, which is by definition given by:
$$g=\begin{pmatrix}
u_{11}(g)&\ldots&u_{1N}(g)\\
\vdots&&\vdots\\
u_{N1}(g)&\ldots&u_{NN}(g)
\end{pmatrix}$$

(2) For the second assertion, we can use the diagonal matrix formed by generators:
$$u=\begin{pmatrix}
g_1&&0\\
&\ddots&\\
0&&g_N
\end{pmatrix}$$

Finally, regarding the last assertion, in the commutative case this follows from the Gelfand theorem, and in the cocommutative case, we will be back to this.
\end{proof}

In order to get now to quantum groups, we will need as well:

\index{Pontrjagin duality}
\index{dual of group}
\index{Fourier transform}

\begin{proposition}
Assuming that $G\subset U_N$ is abelian, we have an identification of Woronowicz algebras $C(G)=C^*(\Gamma)$, with $\Gamma$ being the Pontrjagin dual of $G$:
$$\Gamma=\big\{\chi:G\to\mathbb T\big\}$$
Conversely, assuming that $F_N\to\Gamma$ is abelian, we have an identification of Woronowicz algebras $C^*(\Gamma)=C(G)$, with $G$ being the Pontrjagin dual of $\Gamma$:
$$G=\big\{\chi:\Gamma\to\mathbb T\big\}$$
Thus, the Woronowicz algebras which are both commutative and cocommutative are exactly those of type $A=C(G)=C^*(\Gamma)$, with $G,\Gamma$ being abelian, in Pontrjagin duality.
\end{proposition}

\begin{proof}
This follows from the Gelfand theorem applied to $C^*(\Gamma)$, and from the fact that the characters of a group algebra come from the characters of the group.
\end{proof}

In view of this result, and of the findings from Proposition 1.14 too, we have the following definition, complementing Definition 1.12 and Definition 1.13:

\index{compact quantum group}
\index{discrete quantum group}

\begin{definition}
Given a Woronowicz algebra, we write it as follows, and call $G$ a compact quantum Lie group, and $\Gamma$ a finitely generated discrete quantum group:
$$A=C(G)=C^*(\Gamma)$$
Also, we say that $G,\Gamma$ are dual to each other, and write $G=\widehat{\Gamma},\Gamma=\widehat{G}$.
\end{definition}

Let us discuss now some tools for studying the Woronowicz algebras, and the underlying quantum groups. First, we have the following result:

\index{square of antipode}
\index{Hopf algebra axioms}

\begin{proposition}
Let $(A,u)$ be a Woronowicz algebra.
\begin{enumerate} 
\item $\Delta,\varepsilon$ satisfy the usual axioms for a comultiplication and a counit, namely:
$$(\Delta\otimes id)\Delta=(id\otimes \Delta)\Delta$$
$$(\varepsilon\otimes id)\Delta=(id\otimes\varepsilon)\Delta=id$$

\item $S$ satisfies the antipode axiom, on the $*$-algebra generated by entries of $u$: 
$$m(S\otimes id)\Delta=m(id\otimes S)\Delta=\varepsilon(.)1$$

\item In addition, the square of the antipode is the identity, $S^2=id$.
\end{enumerate}
\end{proposition}

\begin{proof}
As a first observation, the result holds in the commutative case, $A=C(G)$ with $G\subset U_N$. Indeed, here we know from Proposition 1.14 that $\Delta,\varepsilon,S$ appear as functional analytic transposes of the multiplication, unit and inverse maps $m,u,i$:
$$\Delta=m^t\quad,\quad 
\varepsilon=u^t\quad,\quad 
S=i^t$$

Thus, in this case, the various conditions in the statement on $\Delta,\varepsilon,S$ simply come by transposition from the group axioms satisfied by $m,u,i$, namely:
$$m(m\times id)=m(id\times m)$$
$$m(u\times id)=m(id\times u)=id$$
$$m(i\times id)\delta=m(id\times i)\delta=1$$

Here $\delta(g)=(g,g)$. Observe also that the result holds as well in the cocommutative case, $A=C^*(\Gamma)$ with $F_N\to\Gamma$, trivially. In general now, the first axiom follows from:
$$(\Delta\otimes id)\Delta(u_{ij})=(id\otimes \Delta)\Delta(u_{ij})=\sum_{kl}u_{ik}\otimes u_{kl}\otimes u_{lj}$$

As for the other axioms, the verifications here are similar.
\end{proof}

In order to reach to more advanced results, the idea will be that of doing representation theory. Following Woronowicz \cite{wo1}, let us start with the following definition:

\index{representation}
\index{corepresentation}

\begin{definition}
Given $(A,u)$, we call corepresentation of it any unitary matrix $v\in M_n(\mathcal A)$, with $\mathcal A=<u_{ij}>$, satisfying the same conditions as $u$, namely:
$$\Delta(v_{ij})=\sum_kv_{ik}\otimes v_{kj}\quad,\quad
\varepsilon(v_{ij})=\delta_{ij}\quad,\quad
S(v_{ij})=v_{ji}^*$$
We also say that $v$ is a representation of the underlying compact quantum group $G$.
\end{definition}

In the commutative case, $A=C(G)$ with $G\subset U_N$, we obtain in this way the finite dimensional unitary smooth representations $v:G\to U_n$, via the following formula:
$$v(g)=\begin{pmatrix}
v_{11}(g)&\ldots&v_{1n}(g)\\
\vdots&&\vdots\\
v_{n1}(g)&\ldots&v_{nn}(g)
\end{pmatrix}$$

In the cocommutative case, $A=C^*(\Gamma)$ with $F_N\to\Gamma$, we will see in a moment that we obtain in this way the formal sums of elements of $\Gamma$, possibly rotated by a unitary. As a first result now regarding the corepresentations, we have:

\index{sum of representations}
\index{tensor product of representations}
\index{conjugate representation}
\index{spinned representation}

\begin{proposition}
The corepresentations are subject to the following operations:
\begin{enumerate}
\item Making sums, $v+w=diag(v,w)$.

\item Making tensor products, $(v\otimes w)_{ia,jb}=v_{ij}w_{ab}$.

\item Taking conjugates, $(\bar{v})_{ij}=v_{ij}^*$.

\item Rotating by a unitary, $v\to UvU^*$.
\end{enumerate}
\end{proposition}

\begin{proof}
We first check the fact that the matrices in the statement are unitaries:

\medskip

(1) The fact that $v+w$ is unitary is clear. 

\medskip

(2) Regarding now $v\otimes w$, this can be written in standard leg-numbering notation as $v\otimes w=v_{13}w_{23}$, and with this interpretation in mind, the unitarity is clear.

\medskip

(3) In order to check that $\bar{v}$ is unitary, we can use the antipode. Indeed, by regarding the antipode as an antimultiplicative map $S:A\to A$, we have:
$$(\bar{v}v^t)_{ij}
=\sum_kv_{ik}^*v_{jk}
=\sum_kS(v_{kj}^*v_{ki})
=S((v^*v)_{ji})
=\delta_{ij}$$
$$(v^t\bar{v})_{ij}
=\sum_kv_{ki}v_{kj}^*
=\sum_kS(v_{jk}v_{ik}^*)
=S((vv^*)_{ji})
=\delta_{ij}$$

(4) Finally, the fact that $UvU^*$ is unitary is clear. As for the verification of the comultiplicativity axioms, involving $\Delta,\varepsilon,S$, this is routine, in all cases.
\end{proof}

As a consequence of the above result, we can formulate:

\index{Peter-Weyl representation}
\index{colored integer}

\begin{definition}
We denote by $u^{\otimes k}$, with $k=\circ\bullet\bullet\circ\ldots$ being a colored integer, the various tensor products between $u,\bar{u}$, indexed according to the rules 
$$u^{\otimes\emptyset}=1\quad,\quad
u^{\otimes\circ}=u\quad,\quad 
u^{\otimes\bullet}=\bar{u}$$
and multiplicativity, $u^{\otimes kl}=u^{\otimes k}\otimes u^{\otimes l}$, and call them Peter-Weyl corepresentations. 
\end{definition}

Here are a few examples of such corepresentations, namely those coming from the colored integers of length 2, to be often used in what follows:
$$u^{\otimes\circ\circ}=u\otimes u\quad,\quad 
u^{\otimes\circ\bullet}=u\otimes\bar{u}$$
$$u^{\otimes\bullet\circ}=\bar{u}\otimes u\quad,\quad
u^{\otimes\bullet\bullet}=\bar{u}\otimes\bar{u}$$

In order to do representation theory, we first need to know how to integrate over $G$. And we have here the following key result, due to Woronowicz \cite{wo1}:

\index{Haar integration}
\index{Ces\`aro limit}
\index{fixed points}
\index{convolution}

\begin{theorem}
Any Woronowicz algebra $A=C(G)$ has a unique Haar integration, 
$$\left(\int_G\otimes\, id\right)\Delta=\left(id\otimes\int_G\right)\Delta=\int_G(.)1$$
which can be constructed by starting with any faithful positive form $\varphi\in A^*$, and setting
$$\int_G=\lim_{n\to\infty}\frac{1}{n}\sum_{k=1}^n\varphi^{*k}$$
where $\phi*\psi=(\phi\otimes\psi)\Delta$. Moreover, for any corepresentation $v\in M_n(\mathbb C)\otimes A$ we have
$$\left(id\otimes\int_G\right)v=P$$
where $P$ is the orthogonal projection onto $Fix(v)=\{\xi\in\mathbb C^n|v\xi=\xi\}$.
\end{theorem}

\begin{proof}
Following \cite{wo1}, this can be done in 3 steps, as follows:

\medskip

(1) Given $\varphi\in A^*$, our claim is that the following limit converges, for any $a\in A$:
$$\int_\varphi a=\lim_{n\to\infty}\frac{1}{n}\sum_{k=1}^n\varphi^{*k}(a)$$

Indeed, by linearity we can assume that $a$ is the coefficient of certain corepresentation, $a=(\tau\otimes id)v$. But in this case, an elementary computation gives the following formula, with $P_\varphi$ being the orthogonal projection onto the $1$-eigenspace of $(id\otimes\varphi)v$:
$$\left(id\otimes\int_\varphi\right)v=P_\varphi$$

(2) Since $v\xi=\xi$ implies $[(id\otimes\varphi)v]\xi=\xi$, we have $P_\varphi\geq P$, where $P$ is the orthogonal projection onto the following fixed point space:
$$Fix(v)=\left\{\xi\in\mathbb C^n\Big|v\xi=\xi\right\}$$

The point now is that when $\varphi\in A^*$ is faithful, by using a standard positivity trick, we can prove that we have $P_\varphi=P$. Assume indeed $P_\varphi\xi=\xi$, and let us set:
$$a=\sum_i\left(\sum_jv_{ij}\xi_j-\xi_i\right)\left(\sum_kv_{ik}\xi_k-\xi_i\right)^*$$

A straightforward computation shows then that $\varphi(a)=0$, and so $a=0$, as desired.

\medskip

(3) With this in hand, the left and right invariance of $\int_G=\int_\varphi$ is clear on coefficients, and so in general, and this gives all the assertions. See \cite{wo1}. 
\end{proof}

We can now develop a Peter-Weyl type theory for the corepresentations, in analogy with the theory from the classical case. We will need:

\index{Hom space}
\index{End space}
\index{Fix space}
\index{irreducible representation}
\index{equivalent representations}
\index{intertwining spaces}

\begin{definition}
Given two corepresentations $v\in M_n(A),w\in M_m(A)$, we set 
$$Hom(v,w)=\left\{T\in M_{m\times n}(\mathbb C)\Big|Tv=wT\right\}$$
and we use the following conventions:
\begin{enumerate}
\item We use the notations $Fix(v)=Hom(1,v)$, and $End(v)=Hom(v,v)$.

\item We write $v\sim w$ when $Hom(v,w)$ contains an invertible element.

\item We say that $v$ is irreducible, and write $v\in Irr(G)$, when $End(v)=\mathbb C1$.
\end{enumerate}
\end{definition}

In the classical case, where $A=C(G)$ with $G\subset U_N$ being a closed subgroup, we obtain in this way the usual notions regarding the representation intertwiners. Observe also that in the group dual case we have $g\sim h$ when $g=h$. Finally, observe that $v\sim w$ means that $v,w$ are conjugated by an invertible matrix. 

\bigskip

Here are now a few basic results, regarding the above linear spaces:

\index{tensor category}

\begin{proposition}
We have the following results:
\begin{enumerate}
\item $T\in Hom(u,v),S\in Hom(v,w)\implies ST\in Hom(u,w)$.

\item $S\in Hom(u,v),T\in Hom(w,z)\implies S\otimes T\in Hom(u\otimes w,v\otimes z)$.

\item $T\in Hom(v,w)\implies T^*\in Hom(w,v)$.
\end{enumerate}
In other words, the Hom spaces form a tensor $*$-category.
\end{proposition}

\begin{proof}
The proofs are all elementary, as follows:

\medskip

(1) Assume indeed that we have $Tu=vT$, $Sv=Ws$. We obtain, as desired:
$$STu
=SvT
=wST$$

(2) Assuming that we have $Su=vS$, $Tw=zT$, we obtain, as desired:
$$(S\otimes T)(u\otimes w)
=(Su)_{13}(Tw)_{23}
=(vS)_{13}(zT)_{23}
=(v\otimes z)(S\otimes T)$$

(3) By conjugating, and then using the unitarity of $v,w$, we obtain:
\begin{eqnarray*}
Tv=wT
&\implies&v^*T^*=T^*w^*\\
&\implies&vv^*T^*w=vT^*w^*w\\
&\implies&T^*w=vT^*
\end{eqnarray*}

Finally, the last assertion follows from definitions, and from the obvious fact that, in addition to (1,2,3), the Hom spaces are linear spaces, and contain the units.
\end{proof}

Finally, in order to formulate the Peter-Weyl results, we will need as well:

\index{character of representation}

\begin{proposition}
The characters of the corepresentations, given by 
$$\chi_v=\sum_iv_{ii}$$ 
behave as follows, in respect to the various operations:
$$\chi_{v+w}=\chi_v+\chi_w\quad,\quad
\chi_{v\otimes w}=\chi_v\chi_w\quad,\quad 
\chi_{\bar{v}}=\chi_v^*$$
In addition, given two equivalent corepresentations, $v\sim w$, we have $\chi_v=\chi_w$.
\end{proposition}

\begin{proof}
The three formulae in the statement are all clear from definitions. Regarding now the last assertion, assuming that we have $v=T^{-1}wT$, we obtain:
$$\chi_v
=Tr(v)
=Tr(T^{-1}wT)
=Tr(w)
=\chi_w$$

We conclude that $v\sim w$ implies $\chi_v=\chi_w$, as claimed.
\end{proof}

Consider the dense $*$-subalgebra $\mathcal A\subset A$ generated by the coefficients of the fundamental corepresentation $u$, and endow it with the following scalar product: 
$$<a,b>=\int_Gab^*$$

With this convention, we have the following fundamental result, from \cite{wo1}:

\index{Peter-Weyl theorem}
\index{Frobenius isomorphism}
\index{finite dimensional algebra}

\begin{theorem}
We have the following Peter-Weyl type results:
\begin{enumerate}
\item Any corepresentation decomposes as a sum of irreducible corepresentations.

\item Each irreducible corepresentation appears inside a certain $u^{\otimes k}$.

\item $\mathcal A=\bigoplus_{v\in Irr(A)}M_{\dim(v)}(\mathbb C)$, the summands being pairwise orthogonal.

\item The characters of irreducible corepresentations form an orthonormal system.
\end{enumerate}
\end{theorem}

\begin{proof}
All these results are from Woronowicz \cite{wo1}, the idea being as follows:

\medskip

(1) Given a corepresentation $v\in M_n(A)$, we know from Proposition 1.23 that $End(v)$ is a finite dimensional $C^*$-algebra, and by using Proposition 1.6, we obtain:
$$End(v)=M_{n_1}(\mathbb C)\oplus\ldots\oplus M_{n_k}(\mathbb C)$$

But this decomposition allows us to define subcorepresentations $v_i\subset v$, which are irreducible, so we obtain, as desired, a decomposition $v=v_1+\ldots+v_k$.

\medskip

(2) To any corepresentation $v\in M_n(A)$ we associate its space of coefficients, given by $C(v)=span(v_{ij})$. The construction $v\to C(v)$ is then functorial, in the sense that it maps subcorepresentations into subspaces. Observe also that we have:
$$\mathcal A=\sum_{k\in\mathbb N*\mathbb N}C(u^{\otimes k})$$

Now given an arbitrary corepresentation $v\in M_n(A)$, the corresponding coefficient space is a finite dimensional subspace $C(v)\subset\mathcal A$, and so we must have, for certain positive integers $k_1,\ldots,k_p$, an inclusion of vector spaces, as follows:
$$C(v)\subset C(u^{\otimes k_1}\oplus\ldots\oplus u^{\otimes k_p})$$

Thus we have $v\subset u^{\otimes k_1}\oplus\ldots\oplus u^{\otimes k_p}$, and by (1) we obtain the result.

\medskip

(3) As a first observation, which follows from an elementary computation, for any two corepresentations $v,w$ we have a Frobenius type isomorphism, as follows:
$$Hom(v,w)\simeq Fix(\bar{v}\otimes w)$$

Now assume $v\not\sim w$, and let us set $P_{ia,jb}=\int_Gv_{ij}w_{ab}^*$. According to Theorem 1.21, the matrix $P$ is the orthogonal projection onto the following vector space:
$$Fix(v\otimes\bar{w})
\simeq Hom(\bar{v},\bar{w})
=\{0\}$$

Thus we have $P=0$, and so $C(v)\perp C(w)$, which gives the result.

\medskip

(4) The fact that the characters form indeed an orthogonal system follows from (3). Regarding now the norm 1 assertion, consider the following integrals:
$$P_{ik,jl}=\int_Gv_{ij}v_{kl}^*$$

We know from Theorem 1.21 that these integrals form the orthogonal projection onto $Fix(v\otimes\bar{v})\simeq End(\bar{v})=\mathbb C1$. By using this fact, we obtain the following formula:
$$\int_G\chi_v\chi_v^*
=\sum_{ij}\int_Gv_{ii}v_{jj}^*
=\sum_i\frac{1}{N}
=1$$

Thus the characters have indeed norm 1, and we are done.
\end{proof}

Observe that in the cocommutative case, we obtain from (4) that our algebra must be of the form $A=C^*(\Gamma)$, for some discrete group $\Gamma$, as mentioned in Proposition 1.14. As another consequence of the above results, following Woronowicz \cite{wo1}, we have:

\index{full algebra}
\index{reduced algebra}
\index{amenable quantum group}
\index{amenable Woronowicz algebra}
\index{coamenable quantum group}
\index{Kesten amenability}

\begin{theorem}
Let $A_{full}$ be the enveloping $C^*$-algebra of $\mathcal A$, and $A_{red}$ be the quotient of $A$ by the null ideal of the Haar integration. The following are then equivalent:
\begin{enumerate}
\item The Haar functional of $A_{full}$ is faithful.

\item The projection map $A_{full}\to A_{red}$ is an isomorphism.

\item The counit map $\varepsilon:A_{full}\to\mathbb C$ factorizes through $A_{red}$.

\item We have $N\in\sigma(Re(\chi_u))$, the spectrum being taken inside $A_{red}$.
\end{enumerate}
If this is the case, we say that the underlying discrete quantum group $\Gamma$ is amenable.
\end{theorem}

\begin{proof}
This is well-known in the group dual case, $A=C^*(\Gamma)$, with $\Gamma$ being a usual discrete group. In general, the result follows by adapting the group dual case proof:

\medskip

$(1)\iff(2)$ This simply follows from the fact that the GNS construction for the algebra $A_{full}$ with respect to the Haar functional produces the algebra $A_{red}$.

\medskip

$(2)\iff(3)$ Here $\implies$ is trivial, and conversely, a counit $\varepsilon:A_{red}\to\mathbb C$ produces an isomorphism $\Phi:A_{red}\to A_{full}$, by slicing the map $\widetilde{\Delta}:A_{red}\to A_{red}\otimes A_{full}$.

\medskip

$(3)\iff(4)$ Here $\implies$ is clear, coming from $\varepsilon(N-Re(\chi (u)))=0$, and the converse can be proved by doing some functional analysis. See \cite{wo1}.
\end{proof}

With these results in hand, we can formulate, as a refinement of Definition 1.16:

\index{GNS construction}
\index{compact quantum group}
\index{discrete quantum group}

\begin{definition}
Given a Woronowicz algebra $A$, we formally write as before
$$A=C(G)=C^*(\Gamma)$$
and by GNS construction with respect to the Haar functional, we write as well
$$A''=L^\infty(G)=L(\Gamma)$$
with $G$ being a compact quantum group, and $\Gamma$ being a discrete quantum group.
\end{definition}

Now back to Theorem 1.26, as in the discrete group case, the most interesting criterion for amenability, leading to some interesting mathematics and physics, is the Kesten one, (4) there. This leads us into computing character laws:

\index{main character}
\index{amenability}
\index{coamenability}
\index{Kesten amenability}
\index{law of character}

\begin{theorem}
Given a Woronowicz algebra $(A,u)$, consider its main character:
$$\chi=\sum_iu_{ii}$$
\begin{enumerate}
\item The moments of $\chi$ are the numbers $M_k=\dim(Fix(u^{\otimes k}))$.

\item When $u\sim\bar{u}$ the law of $\chi$ is a real measure, supported by $\sigma(\chi)$.

\item The notion of coamenability of $A$ depends only on $law(\chi)$. 
\end{enumerate}
\end{theorem}

\begin{proof}
All this follows from the above results, the idea being as follows:

\medskip

(1) This follows indeed from Peter-Weyl theory.

\medskip

(2) When $u\sim\bar{u}$ we have $\chi=\chi^*$, which gives the result.

\medskip

(3) This follows from Theorem 1.26 (4), and from (2) applied to $u+\bar{u}$.
\end{proof}

This was for the basic theory of compact and discrete quantum groups. For more on all this, we refer to Woronowicz \cite{wo1} and related papers, or to the book \cite{ba4}.

\section*{1c. Quantum rotations}

We know so far that the compact quantum groups include the usual compact Lie groups, $G\subset U_N$, and the abstract duals $G=\widehat{\Gamma}$ of the finitely generated groups $F_N\to\Gamma$. We can combine these examples by performing basic operations, as follows:

\index{tensor product}
\index{free product}
\index{dual free product}

\begin{proposition}
The class of Woronowicz algebras is stable under taking:
\begin{enumerate}
\item Tensor products, $A=A'\otimes A''$, with $u=u'+u''$. At the quantum group level we obtain usual products, $G=G'\times G''$ and $\Gamma=\Gamma'\times\Gamma''$.

\item Free products, $A=A'*A''$, with $u=u'+u''$. At the quantum group level we obtain dual free products $G=G'\,\hat{*}\,G''$ and free products $\Gamma=\Gamma'*\Gamma''$.
\end{enumerate}
\end{proposition}

\begin{proof}
Everything here is clear from definitions. In addition to this, let us mention as well that we have $\int_{A'\otimes A''}=\int_{A'}\otimes\int_{A''}$ and $\int_{A'*A''}=\int_{A'}*\int_{A''}$. Also, the corepresentations of the products can be explicitely computed. See Wang \cite{wa1}.
\end{proof}

Here are some further basic operations, once again from Wang \cite{wa1}:

\index{Woronowicz subalgebra}
\index{quotient algebra}
\index{quantum subgroup}
\index{quotient quantum group}

\begin{proposition}
The class of Woronowicz algebras is stable under taking:
\begin{enumerate}
\item Subalgebras $A'=<u'_{ij}>\subset A$, with $u'$ being a corepresentation of $A$. At the quantum group level we obtain quotients $G\to G'$ and subgroups $\Gamma'\subset\Gamma$.

\item Quotients $A\to A'=A/I$, with $I$ being a Hopf ideal, $\Delta(I)\subset A\otimes I+I\otimes A$. At the quantum group level we obtain subgroups $G'\subset G$ and quotients $\Gamma\to\Gamma'$.
\end{enumerate}
\end{proposition}

\begin{proof}
Once again, everything is clear, and we have as well some straightforward supplementary results, regarding integration and corepresentations. See \cite{wa1}.
\end{proof}

Finally, here are two more operations, which are of key importance:

\index{projective version}
\index{free complexification}

\begin{proposition}
The class of Woronowicz algebras is stable under taking:
\begin{enumerate}
\item Projective versions, $PA=<w_{ia,jb}>\subset A$, where $w=u\otimes\bar{u}$. At the quantum group level we obtain projective versions, $G\to PG$ and $P\Gamma\subset\Gamma$.

\item Free complexifications, $\widetilde{A}=<zu_{ij}>\subset C(\mathbb T)*A$. At the quantum group level we obtain free complexifications, denoted $\widetilde{G}$ and $\widetilde{\Gamma}$.
\end{enumerate}
\end{proposition}

\begin{proof}
This is clear from the previous results. For details here, we refer to \cite{wa1}.
\end{proof}

Once again following Wang \cite{wa1} and related papers, let us discuss now a number of truly ``new'' quantum groups, obtained by liberating. We first have:

\index{free orthogonal group}
\index{free unitary group}

\begin{theorem}
The following universal algebras are Woronowicz algebras,
$$C(O_N^+)=C^*\left((u_{ij})_{i,j=1,\ldots,N}\Big|u=\bar{u},u^t=u^{-1}\right)$$
$$C(U_N^+)=C^*\left((u_{ij})_{i,j=1,\ldots,N}\Big|u^*=u^{-1},u^t=\bar{u}^{-1}\right)$$
so the underlying quantum spaces $O_N^+,U_N^+$ are compact quantum groups.
\end{theorem}

\begin{proof}
This comes from the elementary fact that if a matrix $u=(u_{ij})$ is orthogonal or biunitary, then so must be the following matrices:
$$(u^\Delta)_{ij}=\sum_ku_{ik}\otimes u_{kj}\quad,\quad 
(u^\varepsilon)_{ij}=\delta_{ij}\quad,\quad 
(u^S)_{ij}=u_{ji}^*$$

Thus we can define $\Delta,\varepsilon,S$ by using the universal property of $C(O_N^+)$, $C(U_N^+)$.
\end{proof}

Now with this done, we can look for various intermediate subgroups $O_N\subset O_N^\times\subset O_N^+$ and $U_N\subset U_N^\times\subset U_N^+$. Following \cite{bsp}, a basic construction here is as follows:

\index{half-classical orthogonal group}
\index{half-classical unitary group}
\index{half-commutation}

\begin{theorem}
The following quotient algebras are Woronowicz algebras,
$$C(O_N^*)=C(O_N^+)\Big/\left<abc=cba\Big|\forall a,b,c\in\{u_{ij}\}\right>$$
$$C(U_N^*)=C(U_N^+)\Big/\left<abc=cba\Big|\forall a,b,c\in\{u_{ij},u_{ij}^*\}\right>$$
so the underlying quantum spaces $O_N^*,U_N^*$ are compact quantum groups.
\end{theorem}

\begin{proof}
This follows as in the proof of Theorem 1.32, because if the entries of $u$ satisfy the half-commutation relations $abc=cba$, then so do the entries of $u^\Delta,u^\varepsilon,u^S$.
\end{proof}

Obviously, there are many more things that can be done here, with the above constructions being just the tip of the iceberg. But instead of discussing this, let us first verify that Theorem 1.32 and Theorem 1.33 provide us indeed with new quantum groups. For this purpose, we can use the notion of diagonal torus, which is as follows:

\index{diagonal torus}
\index{group dual}

\begin{proposition}
Given a closed subgroup $G\subset U_N^+$, consider its diagonal torus, which is the closed subgroup $T\subset G$ constructed as follows:
$$C(T)=C(G)\Big/\left<u_{ij}=0\Big|\forall i\neq j\right>$$
This torus is then a group dual, $T=\widehat{\Lambda}$, where $\Lambda=<g_1,\ldots,g_N>$ is the discrete group generated by the elements $g_i=u_{ii}$, which are unitaries inside $C(T)$.
\end{proposition}

\begin{proof}
Since $u$ is unitary, its diagonal entries $g_i=u_{ii}$ are unitaries inside $C(T)$. Moreover, from $\Delta(u_{ij})=\sum_ku_{ik}\otimes u_{kj}$ we obtain, when passing inside the quotient:
$$\Delta(g_i)=g_i\otimes g_i$$

It follows that we have $C(T)=C^*(\Lambda)$, modulo identifying as usual the $C^*$-completions of the various group algebras, and so that we have $T=\widehat{\Lambda}$, as claimed.
\end{proof}

We can now distinguish between our various quantum groups, as follows:

\begin{theorem}
The diagonal tori of the basic unitary quantum groups, namely
$$\xymatrix@R=15mm@C=15mm{
U_N\ar[r]&U_N^*\ar[r]&U_N^+\\
O_N\ar[r]\ar[u]&O_N^*\ar[r]\ar[u]&O_N^+\ar[u]}$$
are the following discrete group duals,
$$\xymatrix@R=14.5mm@C=15mm{
\widehat{\mathbb Z^N}\ar[r]&\widehat{\mathbb Z^{\circ N}}\ar[r]&\widehat{F_N}\\
\widehat{\mathbb Z_2^N}\ar[r]\ar[u]&\widehat{\mathbb Z_2^{\circ N}}\ar[r]\ar[u]&\widehat{\mathbb Z_2^{*N}}\ar[u]}$$
with $\circ$ standing for the half-classical product operation for groups.
\end{theorem}

\begin{proof}
This is clear for $U_N^+$, where on the diagonal we obtain the biggest possible group dual, namely $\widehat{F_N}$. For the other quantum groups this follows by taking quotients, which correspond to taking quotients as well, at the level of the groups $\Lambda=\widehat{T}$.
\end{proof}

As a consequence of the above result, the quantum groups that we have are indeed distinct. There are many more things that can be said about these quantum groups, and about further versions of these quantum groups that can be constructed. More later.

\section*{1d. Quantum permutations}

Eventually. Following Wang \cite{wa2}, let us discuss now the construction and basic properties of the quantum permutation group $S_N^+$. Let us first look at $S_N$. We have:

\index{magic matrix}
\index{magic unitary}
\index{symmetric group}

\begin{proposition}
Consider the symmetric group $S_N$, viewed as permutation group of the $N$ coordinate axes of $\mathbb R^N$. The coordinate functions on $S_N\subset O_N$ are given by
$$u_{ij}=\chi\left(\sigma\in G\Big|\sigma(j)=i\right)$$
and the matrix $u=(u_{ij})$ that these functions form is magic, in the sense that its entries are projections $(p^2=p^*=p)$, summing up to $1$ on each row and each column.
\end{proposition}

\begin{proof}
The action of $S_N$ on the standard basis $e_1,\ldots,e_N\in\mathbb R^N$ being given by $\sigma:e_j\to e_{\sigma(j)}$, this gives the formula of $u_{ij}$ in the statement. As for the fact that the matrix $u=(u_{ij})$ that these functions form is magic, this is clear.
\end{proof}

With a bit more effort, we obtain the following nice characterization of $S_N$:

\index{Gelfand theorem}
\index{commutative algebra}

\begin{theorem}
The algebra of functions on $S_N$ has the following presentation,
$$C(S_N)=C^*_{comm}\left((u_{ij})_{i,j=1,\ldots,N}\Big|u={\rm magic}\right)$$
and the multiplication, unit and inversion map of $S_N$ appear from the maps
$$\Delta(u_{ij})=\sum_ku_{ik}\otimes u_{kj}\quad,\quad 
\varepsilon(u_{ij})=\delta_{ij}\quad,\quad 
S(u_{ij})=u_{ji}$$
defined at the algebraic level, of functions on $S_N$, by transposing.
\end{theorem}

\begin{proof}
The universal algebra $A$ in the statement being commutative, by the Gelfand theorem it must be of the form $A=C(X)$, with $X$ being a certain compact space. Now since we have coordinates $u_{ij}:X\to\mathbb R$, we have an embedding $X\subset M_N(\mathbb R)$. Also, since we know that these coordinates form a magic matrix, the elements $g\in X$ must be 0-1 matrices, having exactly one 1 entry on each row and each column, and so $X=S_N$. Thus we have proved the first assertion, and the second assertion is clear as well.
\end{proof}

Following now Wang \cite{wa2}, we can liberate $S_N$, as follows:

\index{quantum permutation group}
\index{magic unitary}
\index{free symmetric group}

\begin{theorem}
The following universal $C^*$-algebra, with magic meaning as usual formed by projections $(p^2=p^*=p)$, summing up to $1$ on each row and each column,
$$C(S_N^+)=C^*\left((u_{ij})_{i,j=1,\ldots,N}\Big|u={\rm magic}\right)$$
is a Woronowicz algebra, with comultiplication, counit and antipode given by:
$$\Delta(u_{ij})=\sum_ku_{ik}\otimes u_{kj}\quad,\quad 
\varepsilon(u_{ij})=\delta_{ij}\quad,\quad 
S(u_{ij})=u_{ji}$$
Thus the space $S_N^+$ is a compact quantum group, called quantum permutation group.
\end{theorem}

\begin{proof}
As a first observation, the universal $C^*$-algebra in the statement is indeed well-defined, because the conditions $p^2=p^*=p$ satisfied by the coordinates give:
$$||u_{ij}||\leq1$$

In order to prove now that we have a Woronowicz algebra, we must construct maps $\Delta,\varepsilon,S$ given by the formulae in the statement. Consider the following matrices:
$$u^\Delta_{ij}=\sum_ku_{ik}\otimes u_{kj}\quad,\quad 
u^\varepsilon_{ij}=\delta_{ij}\quad,\quad 
u^S_{ij}=u_{ji}$$

Our claim is that, since $u$ is magic, so are these three matrices. Indeed, regarding $u^\Delta$, its entries are idempotents, as shown by the following computation:
$$(u_{ij}^\Delta)^2
=\sum_{kl}u_{ik}u_{il}\otimes u_{kj}u_{lj}
=\sum_{kl}\delta_{kl}u_{ik}\otimes\delta_{kl}u_{kj}
=u_{ij}^\Delta$$

These elements are self-adjoint as well, as shown by the following computation:
$$(u_{ij}^\Delta)^*
=\sum_ku_{ik}^*\otimes u_{kj}^*
=\sum_ku_{ik}\otimes u_{kj}
=u_{ij}^\Delta$$

The row and column sums for the matrix $u^\Delta$ can be computed as follows:
$$\sum_ju_{ij}^\Delta
=\sum_{jk}u_{ik}\otimes u_{kj}
=\sum_ku_{ik}\otimes 1
=1$$
$$\sum_iu_{ij}^\Delta
=\sum_{ik}u_{ik}\otimes u_{kj}
=\sum_k1\otimes u_{kj}
=1$$

Thus, $u^\Delta$ is magic. Regarding now $u^\varepsilon,u^S$, these matrices are magic too, and this for obvious reasons. Thus, all our three matrices $u^\Delta,u^\varepsilon,u^S$ are magic, so we can define $\Delta,\varepsilon,S$ by the formulae in the statement, by using the universality property of $C(S_N^+)$.
\end{proof}

Our first task now is to make sure that Theorem 1.38 produces indeed new quantum groups, which do not collapse to $S_N$. Following Wang \cite{wa2}, we have:

\index{liberation}

\begin{theorem}
We have an embedding $S_N\subset S_N^+$, given at the algebra level by: 
$$u_{ij}\to\chi\left(\sigma\in S_N\Big|\sigma(j)=i\right)$$
This is an isomorphism at $N\leq3$, but not at $N\geq4$, where $S_N^+$ is not classical, nor finite.
\end{theorem} 

\begin{proof}
The fact that we have indeed an embedding as above follows from Theorem 1.37. Observe that in fact more is true, because Theorems 1.37 and 1.38 give:
$$C(S_N)=C(S_N^+)\Big/\Big<ab=ba\Big>$$

Thus, the inclusion $S_N\subset S_N^+$ is a ``liberation'', in the sense that $S_N$ is the classical version of $S_N^+$. We will often use this basic fact, in what follows. Regarding now the second assertion, we can prove this in four steps, as follows:

\medskip

\underline{Case $N=2$}. The fact that $S_2^+$ is indeed classical, and hence collapses to $S_2$, is trivial, because the $2\times2$ magic matrices are as follows, with $p$ being a projection:
$$U=\begin{pmatrix}p&1-p\\1-p&p\end{pmatrix}$$

Indeed, this shows that the entries of $U$ commute. Thus $C(S_2^+)$ is commutative, and so equals its biggest commutative quotient, which is $C(S_2)$. Thus, $S_2^+=S_2$.

\medskip

\underline{Case $N=3$}. By using the same argument as in the $N=2$ case, and the symmetries of the problem, it is enough to check that $u_{11},u_{22}$ commute. But this follows from:
\begin{eqnarray*}
u_{11}u_{22}
&=&u_{11}u_{22}(u_{11}+u_{12}+u_{13})\\
&=&u_{11}u_{22}u_{11}+u_{11}u_{22}u_{13}\\
&=&u_{11}u_{22}u_{11}+u_{11}(1-u_{21}-u_{23})u_{13}\\
&=&u_{11}u_{22}u_{11}
\end{eqnarray*}

Indeed, by applying the involution to this formula, we obtain that we have as well $u_{22}u_{11}=u_{11}u_{22}u_{11}$. Thus, we obtain $u_{11}u_{22}=u_{22}u_{11}$, as desired.

\medskip

\underline{Case $N=4$}. Consider the following matrix, with $p,q$ being projections:
$$U=\begin{pmatrix}
p&1-p&0&0\\
1-p&p&0&0\\
0&0&q&1-q\\
0&0&1-q&q
\end{pmatrix}$$ 

This matrix is magic, and we can choose $p,q\in B(H)$ as for the algebra $<p,q>$ to be noncommutative and infinite dimensional. We conclude that $C(S_4^+)$ is noncommutative and infinite dimensional as well, and so $S_4^+$ is non-classical and infinite, as claimed.

\medskip

\underline{Case $N\geq5$}. Here we can use the standard embedding $S_4^+\subset S_N^+$, obtained at the level of the corresponding magic matrices in the following way:
$$u\to\begin{pmatrix}u&0\\ 0&1_{N-4}\end{pmatrix}$$

Indeed, with this in hand, the fact that $S_4^+$ is a non-classical, infinite compact quantum group implies that $S_N^+$ with $N\geq5$ has these two properties as well.
\end{proof}

The above result is quite surprising. How on Earth can the set $\{1,2,3,4\}$ have an infinity of quantum permutations, and will us be able to fully understand this, one day. But do not worry, the remainder of the present book will be here for that.

\bigskip

As a first observation, as a matter of doublechecking our findings, we are not wrong with our formalism, because as explained once again in \cite{wa2}, we have as well:

\index{quantum permutation}
\index{counting measure}
\index{coaction}
\index{standard coaction}

\begin{theorem}
The quantum permutation group $S_N^+$ acts on the set $X=\{1,\ldots,N\}$, the corresponding coaction map $\Phi:C(X)\to C(X)\otimes C(S_N^+)$ being given by:
$$\Phi(e_i)=\sum_je_j\otimes u_{ji}$$
In fact, $S_N^+$ is the biggest compact quantum group acting on $X$, by leaving the counting measure invariant, in the sense that $(tr\otimes id)\Phi=tr(.)1$, where $tr(e_i)=\frac{1}{N},\forall i$.
\end{theorem}

\begin{proof}
Our claim is that given a compact matrix quantum group $G$, the following formula defines a morphism of algebras, which is a coaction map, leaving the trace invariant, precisely when the matrix $u=(u_{ij})$ is a magic corepresentation of $C(G)$: 
$$\Phi(e_i)=\sum_je_j\otimes u_{ji}$$

Indeed, let us first determine when $\Phi$ is multiplicative. We have:
$$\Phi(e_i)\Phi(e_k)
=\sum_{jl}e_je_l\otimes u_{ji}u_{lk}
=\sum_je_j\otimes u_{ji}u_{jk}$$
$$\Phi(e_ie_k)
=\delta_{ik}\Phi(e_i)
=\delta_{ik}\sum_je_j\otimes u_{ji}$$

We conclude that the multiplicativity of $\Phi$ is equivalent to the following conditions:
$$u_{ji}u_{jk}=\delta_{ik}u_{ji}\quad,\quad\forall i,j,k$$

Similarly, $\Phi$ is unital when $\sum_iu_{ji}=1$, $\forall j$. Finally, the fact that $\Phi$ is a $*$-morphism translates into $u_{ij}=u_{ij}^*$, $\forall i,j$. Summing up, in order for $\Phi(e_i)=\sum_je_j\otimes u_{ji}$ to be a morphism of $C^*$-algebras, the elements $u_{ij}$ must be projections, summing up to 1 on each row of $u$. Regarding now the preservation of the trace, observe that we have:
$$(tr\otimes id)\Phi(e_i)=\frac{1}{N}\sum_ju_{ji}$$

Thus the trace is preserved precisely when the elements $u_{ij}$ sum up to 1 on each of the columns of $u$. We conclude from this that $\Phi(e_i)=\sum_je_j\otimes u_{ji}$ is a morphism of $C^*$-algebras preserving the trace precisely when $u$ is magic, and since the coaction conditions on $\Phi$ are equivalent to the fact that $u$ must be a corepresentation, this finishes the proof of our claim. But this claim proves all the assertions in the statement.
\end{proof}

As a technical comment here, the invariance of the counting measure is a key assumption in Theorem 1.40, in order to have an universal object $S_N^+$. That is, this condition is automatic for classical group actions, but not for quantum group actions, and when dropping it, there is no universal object of type $S_N^+$. This explains the main difficulty behind Question 1.1, and the credit for this discovery goes to Wang \cite{wa2}.

\bigskip

In order to study now $S_N^+$, we can use the technology that we have, which gives:

\index{dual free product}
\index{infinite dihedral group}
\index{half-classical symmetric group}
\index{coamenability}

\begin{theorem}
The quantum groups $S_N^+$ have the following properties:
\begin{enumerate}
\item We have $S_N^+\,\hat{*}\,S_M^+\subset S_{N+M}^+$, for any $N,M$.

\item In particular, we have an embedding $\widehat{D_\infty}\subset S_4^+$. 

\item $S_4\subset S_4^+$ are distinguished by their spinned diagonal tori.

\item If $\mathbb Z_{N_1}*\ldots*\mathbb Z_{N_k}\to\Gamma$, with $N=\sum N_i$, then $\widehat{\Gamma}\subset S_N^+$.

\item The quantum groups $S_N^+$ with $N\geq5$ are not coamenable.

\item The half-classical version $S_N^*=S_N^+\cap O_N^*$ collapses to $S_N$.
\end{enumerate}
\end{theorem}

\begin{proof}
These results follow from what we have, the proofs being as follows:

\medskip

(1) If we denote by $u,v$ the fundamental corepresentations of $C(S_N^+),C(S_M^+)$, the fundamental corepresentation of $C(S_N^+\,\hat{*}\,S_M^+)$ is by definition:
$$w=\begin{pmatrix}u&0\\0&v\end{pmatrix}$$

But this matrix is magic, because both $u,v$ are magic, and this gives the result.

\medskip

(2) This result, which refines our $N=4$ trick from the proof of Theorem 1.39, follows from (1) with $N=M=2$. Indeed, we have the following computation:
\begin{eqnarray*}
S_2^+\,\hat{*}\,S_2^+
&=&S_2\,\hat{*}\, S_2
=\mathbb Z_2\,\hat{*}\, \mathbb Z_2\\
&\simeq&\widehat{\mathbb Z_2}\,\hat{*}\, \widehat{\mathbb Z_2}
=\widehat{\mathbb Z_2*\mathbb Z_2}\\
&=&\widehat{D_\infty}
\end{eqnarray*}

(3) Observe first that $S_4\subset S_4^+$ are not distinguished by their diagonal torus, which is $\{1\}$ for both of them. However, according to the Peter-Weyl theory applied to the group duals, the group dual $\widehat{D_\infty}\subset S_4^+$ from (2) must be a subgroup of the diagonal torus of $(S_4^+,FuF^*)$, for a certain unitary $F\in U_4$, and this gives the result.

\medskip

(4) This result, which generalizes (2), can be deduced as follows:
\begin{eqnarray*}
\widehat{\Gamma}
&\subset&\widehat{\mathbb Z_{N_1}*\ldots*\mathbb Z_{N_k}}
=\widehat{\mathbb Z_{N_1}}\,\hat{*}\,\ldots\,\hat{*}\,\widehat{\mathbb Z_{N_k}}\\
&\simeq&\mathbb Z_{N_1}\,\hat{*}\,\ldots\,\hat{*}\,\mathbb Z_{N_k}
\subset S_{N_1}\,\hat{*}\,\ldots\,\hat{*}\,S_{N_k}\\
&\subset&S_{N_1}^+\,\hat{*}\,\ldots\,\hat{*}\,S_{N_k}^+
\subset S_N^+
\end{eqnarray*}

(5) This follows from (4), because at $N=5$ the dual of the group $\Gamma=\mathbb Z_2*\mathbb Z_3$, which is well-known not to be amenable, embeds into $S_5^+$. As for the general case, that of $S_N^+$ with $N\geq5$, here the result follows by using the embedding $S_5^+\subset S_N^+$.

\medskip

(6) We must prove that $S_N^*=S_N^+\cap O_N^*$ is classical. But here, we can use the fact that for a magic matrix, the entries on each row sum up to 1. Indeed, by making $c$ vary over a full row of $u$, we obtain $abc=cba\implies ab=ba$, as desired.
\end{proof}

The above results are all quite interesting, notably with (2) providing us with a better understanding of why $S_4^+$ is infinite, and with (4) telling us that $S_5^+$ is not only infinite, but just huge. We have as well (6), suggesting that $S_N^+$ might be the only liberation of $S_N$. We will be back to these observations, with further results, in due time.

\section*{1e. Exercises}

There has been a lot of mathematics in this chapter, and as a best exercise on all this, read more about operator algebras and quantum groups, in general. Regarding now the quantum rotations and permutations, as a first exercise about them, we have:

\begin{exercise}
Prove that the quantum group inclusion
$$\widetilde{O_N^+}\subset U_N^+$$
is an isomorphism at the level of the corresponding diagonal tori.
\end{exercise}

To be more precise, the fact that we have an inclusion $\widetilde{O_N^+}\subset U_N^+$ is clear, and it can be actually proved that this inclusion is an isomorphism, but this is non-trivial. So, more on this later, and in the meantime, you can have some fun with the above exercise.

\begin{exercise}
Find the shortest proof ever for the equality
$$S_3^+=S_3$$
by doing some original manipulations with the $3\times3$ magic matrices.
\end{exercise}

To be more precise, we have already seen such a proof in the above, and the problem is that of finding a new proof, a bit in the same spirit, based on new computations.

\begin{exercise}
Prove that we have $S_3^+=S_3$ by looking at the coaction
$$\Phi:\mathbb C^3\to\mathbb C^3\otimes C(S_3^+)$$
written in terms of the Fourier basis of $\mathbb C^3$.
\end{exercise}

To be more precise, the question here is that of changing the basis of $\mathbb C^3$, by using the Fourier transform over $\mathbb Z_3$, and then deducing that the coefficients must commute.

\begin{exercise}
Prove that the following quantum groups are not coamenable:
$$O_N^+\ (N\geq3)\quad,\quad U_N^+\ (N\geq2)$$
For a bonus point, prove as well that $O_2^+$ and $S_4^+$ are coamenable.
\end{exercise}

Here the first part can only be quite standard, by using the same type of ideas as for $S_N^+$ with $N\geq5$. As for the bonus point question, this is something quite difficult, and I don't know myself a simple proof for that. Which of course does not mean anything, old man and I might just miss something, and this is a good problem for you, reader.

\chapter{Diagrams, easiness}

\section*{2a. Some philosophy}

We have seen the definition and basic properties of $S_N^+$, and a number of more advanced results as well, such as the non-isomorphism of $S_N\subset S_N^+$ at $N\geq4$, obtained by using suitable group duals $\widehat{\Gamma}\subset S_N^+$. It is possible to further build along these lines, but all this remains quite amateurish. For strong results, we must do representation theory. 

\bigskip

So, let us first go back to the general closed subgroups $G\subset U_N^+$. We have seen in chapter 1 that such quantum groups have a Haar measure, and that by using this, a Peter-Weyl theory can be developed for them. However, all this is just a beginning, and many more things can be said, at the general level, which are all useful. We will present now this material, and go back afterwards to our problems regarding $S_N^+$.

\bigskip

Let us start with a claim, which is quite precise, and advanced, and which will stand as a guiding principle for this chapter, and in fact for the remainder of this book:

\index{main character}
\index{moments}

\begin{claim}
Given a closed subgroup $G\subset_uU_N^+$, no matter what you want to do with it, of algebraic or analytic type, you must compute the following spaces:
$$F_k=Fix(u^{\otimes k})$$
Moreover, for most questions, the computation of the dimensions $M_k=\dim F_k$, which are the moments of the main character $\chi=\sum_iu_{ii}$, will do.
\end{claim}

This might look like a quite bold claim, so let us explain this. Assuming first that you are interested in doing representation theory for $G$, you will certainly run into the spaces $F_k$, via Peter-Weyl theory. In fact, Peter-Weyl tells you that the irreducible representations appear as $r\subset u^{\otimes k}$, so for finding them, you must compute the algebras $C_k=End(u^{\otimes k})$. But the knowledge of these algebras $C_k$ is more or less the same thing as the knowledge of the spaces $F_k$, due to Frobenius duality, as follows:

\begin{proposition}
Given a closed subgroup $G\subset_uU_N^+$, consider the following spaces:
$$F_k=Fix(u^{\otimes k})\quad,\quad
C_k=End(u^{\otimes k})\quad,\quad
C_{kl}=Hom(u^{\otimes k},u^{\otimes l})$$
Then knowing the sequence $\{F_k\}$ is the same as knowing the double sequence $\{C_{kl}\}$, and in the case $1\in u$, this is the same as knowing the sequence $\{C_k\}$.
\end{proposition} 

\begin{proof}
In the particular case of the Peter-Weyl corepresentations, the Frobenius isomorphism $Hom(v,w)\simeq Fix(\bar{v}\otimes w)$, that we know from chapter 1, reads:
$$C_{kl}=Hom(u^{\otimes k},u^{\otimes l})=Fix(u^{\otimes \bar{k}l})=F_{\bar{k}l}$$

But this gives the equivalence in the statement. Regarding now the last assertion, assuming $1\in u$ we have $1\in u^{\otimes k}$ for any colored integer $k$, and so:
$$F_k=Hom(1,u^{\otimes k})\subset Hom(u^{\otimes k},u^{\otimes k})=C_k$$

Thus the spaces $F_k$ can be identified inside the algebras $C_k$, and we are done.
\end{proof}

Summarizing, we have now good algebraic motivations for Claim 2.1. Before going further, however, let us point out that looking at Proposition 2.2 leads us a bit into a dillema, on which spaces are the best to use. And the traditional answer here is that the spaces $C_{kl}$ are the best, due to Tannakian duality, which is as follows:

\index{tensor category}
\index{Tannakian category}
\index{Tannakian duality}
\index{bicommutant}

\begin{theorem}
The following operations are inverse to each other:
\begin{enumerate}
\item The construction $G\to C$, which associates to a closed subgroup $G\subset_uU_N^+$ the tensor category formed by the intertwiner spaces $C_{kl}=Hom(u^{\otimes k},u^{\otimes l})$.

\item The construction $C\to G$, associating to a tensor category $C$ the closed subgroup $G\subset_uU_N^+$ coming from the relations $T\in Hom(u^{\otimes k},u^{\otimes l})$, with $T\in C_{kl}$.
\end{enumerate}
\end{theorem}

\begin{proof}
This is something quite deep, going back to Woronowicz \cite{wo2} in a slightly different form, and to Malacarne \cite{mal} in the simplified form above. The idea is that we have indeed a construction $G\to C_G$, whose output is a tensor $C^*$-subcategory with duals of the tensor $C^*$-category of finite dimensional Hilbert spaces, as follows:
$$(C_G)_{kl}=Hom(u^{\otimes k},u^{\otimes l})$$

We have as well a construction $C\to G_C$, obtained by setting:
$$C(G_C)=C(U_N^+)\big/\left<T\in Hom(u^{\otimes k},u^{\otimes l})\Big|\forall k,l,\forall T\in C_{kl}\right>$$

Regarding now the bijection claim, some elementary algebra shows that $C=C_{G_C}$ implies $G=G_{C_G}$, and that $C\subset C_{G_C}$ is automatic. Thus we are left with proving:
$$C_{G_C}\subset C$$

But this latter inclusion can be proved indeed, by doing some algebra, and using von Neumann's bicommutant theorem, in finite dimensions. See Malacarne \cite{mal}. 
\end{proof}

The above result is something quite abstract, yet powerful. We will see applications of it in a moment, in the form of Brauer theorems for $U_N,O_N,S_N$ and $U_N^+,O_N^+,S_N^+$.

\bigskip

All this is very good, providing us with strong motivations for Claim 2.1. However, algebra is of course not everything, and we must comment now on analysis as well. As an analyst you would like to know how to integrate over $G$, and here, we have:

\index{Weingarten formula}
\index{Tannakian category}
\index{Gram matrix}

\begin{theorem}
The integration over $G\subset_uU_N^+$ is given by the Weingarten formula
$$\int_Gu_{i_1j_1}^{e_1}\ldots u_{i_kj_k}^{e_k}=\sum_{\pi,\sigma\in D_k}\delta_\pi(i)\delta_\sigma(j)W_k(\pi,\sigma)$$
for any colored integer $k=e_1\ldots e_k$ and indices $i,j$, where $D_k$ is a linear basis of $Fix(u^{\otimes k})$, 
$$\delta_\pi(i)=<\pi,e_{i_1}\otimes\ldots\otimes e_{i_k}>$$
and $W_k=G_k^{-1}$, with $G_k(\pi,\sigma)=<\pi,\sigma>$.
\end{theorem}

\begin{proof}
We know from chapter 1 that the integrals in the statement form altogether the orthogonal projection $P^k$ onto the following space:
$$Fix(u^{\otimes k})=span(D_k)$$

Consider now the following linear map, with $D_k=\{\xi_k\}$ being as in the statement:
$$E(x)=\sum_{\pi\in D_k}<x,\xi_\pi>\xi_\pi$$

By a standard linear algebra computation, it follows that we have $P=WE$, where $W$ is the inverse on $span(T_\pi|\pi\in D_k)$ of the restriction of $E$. But this restriction is the linear map given by $G_k$, and so $W$ is the linear map given by $W_k$, and this gives the result.
\end{proof}

As a conclusion, regardless on whether you're an algebraist or an analyst, if you want to study $G\subset_uU_N$ you are led into the computation of the spaces $F_k=Fix(u^{\otimes k})$. However, the story is not over here, because you might say that you are a functional analyst, interested in the fine analytic properties of the dual $\Gamma=\widehat{G}$. But here, I would strike back with the following statement, based on the Kesten amenability criterion:

\index{main character}
\index{amenability}
\index{coamenability}
\index{Kesten amenability}

\begin{proposition}
Given a closed subgroup $G\subset_uU_N^+$, consider its main character:
$$\chi=\sum_iu_{ii}$$
\begin{enumerate}
\item The moments of $\chi$ are the numbers $M_k=\dim(Fix(u^{\otimes k}))$.

\item When $u\sim\bar{u}$ the law of $\chi$ is a real measure, supported by $\sigma(\chi)$.

\item The notion of amenability of $\Gamma=\widehat{G}$ depends only on $law(\chi)$. 
\end{enumerate}
\end{proposition}

\begin{proof}
This is something that we know from chapter 1, the idea being that (1) comes from Peter-Weyl theory, that (2) comes from $u\sim\bar{u}\implies\chi=\chi^*$, and that (3) comes from the Kesten amenability criterion, and from (2) applied to $u+\bar{u}$.
\end{proof}

Finally, you might argue that you are in fact a pure mathematician interested in the combinatorial beauty of the dual $\Gamma=\widehat{G}$. But I have an answer to this too, as follows, again urging you to look at the spaces $F_k=Fix(u^{\otimes k})$, before getting into $\Gamma$:

\index{Cayley graph}
\index{geodesic distance}
\index{main character}

\begin{proposition}
Consider a closed subgroup $G\subset_uU_N^+$, and assume, by enlarging if necessary $u$, that we have $1\in u=\bar{u}$. The formula
$$d(v,w)=\min\left\{k\in\mathbb N\Big|1\subset\bar{v}\otimes w\otimes u^{\otimes k}\right\}$$
defines then a distance on $Irr(G)$, which coincides with the geodesic distance on the associated Cayley graph. Moreover, the moments of the main character, 
$$\int_G\chi^k=\dim\left(Fix(u^{\otimes k})\right)$$
count the loops based at $1$, having lenght $k$, on the corresponding Cayley graph.
\end{proposition}

\begin{proof}
Observe first the result holds indeed in the group dual case, where the Woronowicz algebra is $A=C^*(\Gamma)$, with $\Gamma=<S>$ being a finitely generated discrete group. In general, the fact that the lengths are finite follows from Peter-Weyl theory. The symmetry axiom is clear as well, and the triangle inequality is elementary to establish too. Finally, the last assertion, regarding the moments, is elementary too.
\end{proof}

As a conclusion, looks like I won the debate, with Claim 2.1 reigning over both the compact and discrete quantum group worlds, without opposition. Before getting further, let us record a result in relation with the second part of that claim, as follows:

\index{main character}
\index{Kesten measure}
\index{amenability}

\begin{theorem}
Given a closed subgroup $G\subset_uU_N^+$, the law of its main character
$$\chi=\sum_iu_{ii}$$
with respect to the Haar integration has the following properties:
\begin{enumerate}
\item The moments of $\chi$ are the numbers $M_k=\dim(Fix(u^{\otimes k}))$.

\item $M_k$ counts the lenght $k$ loops at $1$, on the Cayley graph of $\Gamma=\widehat{G}$.

\item $law(\chi)$ is the Kesten measure of the discrete quantum group $\Gamma=\widehat{G}$.

\item When $u\sim\bar{u}$ the law of $\chi$ is a usual measure, supported on $[-N,N]$.

\item $\Gamma=\widehat{G}$ is amenable precisely when $N\in supp(law(Re(\chi)))$.

\item Any inclusion $G\subset_u H\subset_vU_N^+$ must decrease the numbers $M_k$.

\item Such an inclusion is an isomorphism when $law(\chi_u)=law(\chi_v)$.
\end{enumerate}
\end{theorem}

\begin{proof}
All this is very standard, coming from the Peter-Weyl theory developed by Woronowicz in \cite{wo1}, and explained in chapter 1, the idea being as follows:

\medskip

(1) This comes from the Peter-Weyl type theory, which tells us the number of fixed points of $v=u^{\otimes k}$ can be recovered by integrating the character $\chi_v=\chi_u^k$.

\medskip

(2) This is something true, and well-known, for $G=\widehat{\Gamma}$ with $\Gamma=<g_1,\ldots,g_N>$ being a discrete group. In general, the proof is quite similar.

\medskip

(3) This is actually the definition of the Kesten measure, in the case $G=\widehat{\Gamma}$, with $\Gamma=<g_1,\ldots,g_N>$ being a discrete group. In general, this follows from (2).

\medskip

(4) The equivalence $u\sim\bar{u}$ translates into $\chi_u=\chi_u^*$, and this gives the first assertion. As for the support claim, this follows from $uu^*=1\implies||u_{ii}||\leq1$, for any $i$.

\medskip

(5) This is the Kesten amenability criterion, which can be established as in the group dual case, $G=\widehat{\Gamma}$, with $\Gamma=<g_1,\ldots,g_N>$ being a discrete group.

\medskip

(6) This is something elementary, which follows from (1), and from the fact that the inclusions of closed subgroups of $U_N^+$ decrease the spaces of fixed points.

\medskip

(7) This follows by using (6), and the Peter-Weyl type theory, the idea being that if $G\subset H$ is not injective, then it must strictly decrease one of the spaces $Fix(u^{\otimes k})$.
\end{proof}

As a conclusion to all this, somewhat improving Claim 2.1, given a closed subgroup $G\subset_uU_N^+$, regardless of our precise motivations, be that algebra, analysis or other, computing the law of $\chi=\sum_iu_{ii}$ is the ``main problem'' to be solved. Good to know.

\section*{2b. Diagrams, easiness}

Let us discuss now the representation theory of $S_N^+$, and the computation of the law of the main character. Our main result here, which will be something quite conceptual, will be the fact that $S_N\subset S_N^+$ is a liberation of ``easy quantum groups''. 

\bigskip

Looking at what has been said above, as a main tool, at the general level, we only have Tannakian duality. So, inspired by that, and following \cite{bsp}, let us formulate:

\index{category of partitions}
\index{concatenation of partitions}
\index{semicircle partition}
\index{matching pairing}

\begin{definition}
Let $P(k,l)$ be the set of partitions between an upper row of $k$ points, and a lower row of $l$ points. A collection of sets
$$D=\bigsqcup_{k,l}D(k,l)$$
with $D(k,l)\subset P(k,l)$ is called a category of partitions when it has the following properties:
\begin{enumerate}
\item Stability under the horizontal concatenation, $(\pi,\sigma)\to[\pi\sigma]$.

\item Stability under the vertical concatenation, $(\pi,\sigma)\to[^\sigma_\pi]$.

\item Stability under the upside-down turning, $\pi\to\pi^*$.

\item Each set $P(k,k)$ contains the identity partition $||\ldots||$.

\item The sets $P(\emptyset,\circ\bullet)$ and $P(\emptyset,\bullet\circ)$ both contain the semicircle $\cap$.
\end{enumerate}
\end{definition} 

As a basic example, we have the category of all partitions $P$ itself. Other basic examples are the category of pairings $P_2$, and the categories $NC,NC_2$ of noncrossing partitions, and pairings. We have as well the category $\mathcal P_2$ of pairings which are ``matching'', in the sense that they connect $\circ-\circ$, $\bullet-\bullet$ on the vertical, and $\circ-\bullet$ on the horizontal, and its subcategory $\mathcal{NC}_2\subset\mathcal P_2$ consisting of the noncrossing matching pairings.

\bigskip

There are many other examples, and we will be back to this. Following \cite{bsp}, the relation with the Tannakian categories and duality comes from:

\index{Kronecker symbol}
\index{maps associated to partitions}

\begin{proposition}
Each partition $\pi\in P(k,l)$ produces a linear map
$$T_\pi:(\mathbb C^N)^{\otimes k}\to(\mathbb C^N)^{\otimes l}$$
given by the following formula, with $e_1,\ldots,e_N$ being the standard basis of $\mathbb C^N$,
$$T_\pi(e_{i_1}\otimes\ldots\otimes e_{i_k})=\sum_{j_1\ldots j_l}\delta_\pi\begin{pmatrix}i_1&\ldots&i_k\\ j_1&\ldots&j_l\end{pmatrix}e_{j_1}\otimes\ldots\otimes e_{j_l}$$
and with the Kronecker type symbols $\delta_\pi\in\{0,1\}$ depending on whether the indices fit or not. The assignement $\pi\to T_\pi$ is categorical, in the sense that we have
$$T_\pi\otimes T_\sigma=T_{[\pi\sigma]}\quad,\quad 
T_\pi T_\sigma=N^{c(\pi,\sigma)}T_{[^\sigma_\pi]}\quad,\quad 
T_\pi^*=T_{\pi^*}$$
where $c(\pi,\sigma)$ are certain integers, coming from the erased components in the middle.
\end{proposition}

\begin{proof}
The concatenation axiom follows from the following computation:
\begin{eqnarray*}
&&(T_\pi\otimes T_\sigma)(e_{i_1}\otimes\ldots\otimes e_{i_p}\otimes e_{k_1}\otimes\ldots\otimes e_{k_r})\\
&=&\sum_{j_1\ldots j_q}\sum_{l_1\ldots l_s}\delta_\pi\begin{pmatrix}i_1&\ldots&i_p\\j_1&\ldots&j_q\end{pmatrix}\delta_\sigma\begin{pmatrix}k_1&\ldots&k_r\\l_1&\ldots&l_s\end{pmatrix}e_{j_1}\otimes\ldots\otimes e_{j_q}\otimes e_{l_1}\otimes\ldots\otimes e_{l_s}\\
&=&\sum_{j_1\ldots j_q}\sum_{l_1\ldots l_s}\delta_{[\pi\sigma]}\begin{pmatrix}i_1&\ldots&i_p&k_1&\ldots&k_r\\j_1&\ldots&j_q&l_1&\ldots&l_s\end{pmatrix}e_{j_1}\otimes\ldots\otimes e_{j_q}\otimes e_{l_1}\otimes\ldots\otimes e_{l_s}\\
&=&T_{[\pi\sigma]}(e_{i_1}\otimes\ldots\otimes e_{i_p}\otimes e_{k_1}\otimes\ldots\otimes e_{k_r})
\end{eqnarray*}

As for the composition and involution axioms, their proof is similar.
\end{proof}

In relation with quantum groups, we have the following result, from \cite{bsp}:

\index{Tannakian duality}

\begin{theorem}
Each category of partitions $D=(D(k,l))$ produces a family of compact quantum groups $G=(G_N)$, one for each $N\in\mathbb N$, via the formula
$$Hom(u^{\otimes k},u^{\otimes l})=span\left(T_\pi\Big|\pi\in D(k,l)\right)$$
which produces a Tannakian category, and so a closed subgroup $G_N\subset_uU_N^+$.
\end{theorem}

\begin{proof}
Let call $C_{kl}$ the spaces on the right. By using the axioms in Definition 2.8, and the categorical properties of the operation $\pi\to T_\pi$, from Proposition 2.9, we see that $C=(C_{kl})$ is a Tannakian category. Thus Theorem 2.3 applies, and gives the result.
\end{proof}

We can now formulate a key definition, as follows:

\index{easy quantum group}

\begin{definition}
A compact quantum group $G_N$ is called easy when we have
$$Hom(u^{\otimes k},u^{\otimes l})=span\left(T_\pi\Big|\pi\in D(k,l)\right)$$
for any colored integers $k,l$, for a certain category of partitions $D\subset P$.
\end{definition}

In other words, a compact quantum group is called easy when its Tannakian category appears in the simplest possible way: from a category of partitions. The terminology is quite natural, because Tannakian duality is basically our only serious tool. In relation now with quantum permutation groups, and with the orthogonal and unitary quantum groups too, here is our main result, coming from \cite{ba1}, \cite{bc1}, \cite{bsp}:

\index{noncrossing partition}
\index{noncrossing pairing}
\index{Brauer theorem}
\index{easy quantum group}
\index{easiness}
\index{Schur-Weyl duality}

\begin{theorem}
The basic quantum permutation and rotation groups,
$$\xymatrix@R=15mm@C=15mm{
S_N^+\ar[r]&O_N^+\ar[r]&U_N^+\\
S_N\ar[r]\ar[u]&O_N\ar[r]\ar[u]&U_N\ar[u]}$$
are all easy, the corresponding categories of partitions being as follows,
$$\xymatrix@R=16mm@C=15mm{
NC\ar[d]&NC_2\ar[l]\ar[d]&\mathcal{NC}_2\ar[l]\ar[d]\\
P&P_2\ar[l]&\mathcal P_2\ar[l]}$$
with $2$ standing for pairings, NC for noncrossing, and calligraphic for matching. 
\end{theorem}

\begin{proof}
This is something quite fundamental, the proof being as follows:

\medskip

(1) The quantum group $U_N^+$ is defined via the following relations:
$$u^*=u^{-1}\quad,\quad 
u^t=\bar{u}^{-1}$$ 

But, by doing some elementary computations, these relations tell us precisely that the following two operators must be in the associated Tannakian category $C$:
$$T_\pi\quad:\quad \pi={\ }^{\,\cap}_{\circ\bullet}\ ,\ 
{\ }^{\,\cap}_{\bullet\circ}$$

Thus, the associated Tannakian category is $C=span(T_\pi|\pi\in D)$, with:
$$D
=<{\ }^{\,\cap}_{\circ\bullet}\,\,,{\ }^{\,\cap}_{\bullet\circ}>
={\mathcal NC}_2$$

(2) The subgroup $O_N^+\subset U_N^+$ is defined by imposing the following relations:
$$u_{ij}=\bar{u}_{ij}$$

Thus, the following operators must be in the associated Tannakian category $C$:
$$T_\pi\quad:\quad\pi=|^{\hskip-1.32mm\circ}_{\hskip-1.32mm\bullet}\ ,\ 
|_{\hskip-1.32mm\circ}^{\hskip-1.32mm\bullet}$$

We conclude that the Tannakian category is $C=span(T_\pi|\pi\in D)$, with:
$$D
=<\mathcal{NC}_2,|^{\hskip-1.32mm\circ}_{\hskip-1.32mm\bullet},|_{\hskip-1.32mm\circ}^{\hskip-1.32mm\bullet}>
=NC_2$$

(3) The subgroup $U_N\subset U_N^+$ is defined via the following relations:
$$[u_{ij},u_{kl}]=0\quad,\quad 
[u_{ij},\bar{u}_{kl}]=0$$

Thus, the following operators must be in the associated Tannakian category $C$:
$$T_\pi\quad:\quad \pi={\slash\hskip-2.1mm\backslash}^{\hskip-2.5mm\circ\circ}_{\hskip-2.5mm\circ\circ}\ ,\ 
{\slash\hskip-2.1mm\backslash}^{\hskip-2.5mm\circ\bullet}_{\hskip-2.5mm\bullet\circ}$$

Thus the associated Tannakian category is $C=span(T_\pi|\pi\in D)$, with:
$$D
=<\mathcal{NC}_2,{\slash\hskip-2.1mm\backslash}^{\hskip-2.5mm\circ\circ}_{\hskip-2.5mm\circ\circ},{\slash\hskip-2.1mm\backslash}^{\hskip-2.5mm\circ\bullet}_{\hskip-2.5mm\bullet\circ}>
=\mathcal P_2$$

(4) In order to deal now with $O_N$, we can simply use the following formula: 
$$O_N=O_N^+\cap U_N$$

At the categorical level, this tells us that $O_N$ is indeed easy, coming from:
$$D
=<NC_2,\mathcal P_2>
=P_2$$

(5) We know that the subgroup $S_N^+\subset O_N^+$ appears as follows:
$$C(S_N^+)=C(O_N^+)\Big\slash\Big<u={\rm magic}\Big>$$

In order to interpret the magic condition, consider the fork partition:
$$Y\in P(2,1)$$

Given a corepresentation $u$, we have the following formulae:
$$(T_Yu^{\otimes 2})_{i,jk}
=\sum_{lm}(T_Y)_{i,lm}(u^{\otimes 2})_{lm,jk}
=u_{ij}u_{ik}$$
$$(uT_Y)_{i,jk}
=\sum_lu_{il}(T_Y)_{l,jk}
=\delta_{jk}u_{ij}$$

We conclude that we have the following equivalence:
$$T_Y\in Hom(u^{\otimes 2},u)\iff u_{ij}u_{ik}=\delta_{jk}u_{ij},\forall i,j,k$$

The condition on the right being equivalent to the magic condition, we obtain:
$$C(S_N^+)=C(O_N^+)\Big\slash\Big<T_Y\in Hom(u^{\otimes 2},u)\Big>$$

Thus $S_N^+$ is indeed easy, the corresponding category of partitions being:
$$D=<Y>=NC$$

(6) Finally, in order to deal with $S_N$, we can use the following formula: 
$$S_N=S_N^+\cap O_N$$

At the categorical level, this tells us that $S_N$ is indeed easy, coming from:
$$D
=<NC,P_2>
=P$$

Thus, we are led to the conclusions in the statement.
\end{proof}

The above result is something quite deep, and we will see in what follows countless applications of it. As a first such application, which is rather philosophical, we have:

\index{liberation}
\index{easy liberation}

\begin{theorem}
The constructions $G_N\to G_N^+$ with $G=U,O,S$ are easy quantum group liberations, in the sense that they come from the construction
$$D\to D\cap NC$$
at the level of the associated categories of partitions.
\end{theorem}

\begin{proof}
This is clear indeed from Theorem 2.12, and from the following trivial equalities, connecting the categories found there:
$$\mathcal{NC}_2=\mathcal{P}_2\cap NC\quad,\quad NC_2=P_2\cap NC\quad,\quad NC=P\cap NC$$

Thus, we are led to the conclusion in the statement.
\end{proof}

The above result is quite nice, because the various constructions $G_N\to G_N^+$ that we saw in chapter 1, although natural, were something quite ad-hoc. Now all this is no longer ad-hoc, and the next time that we will have to liberate a subgroup $G_N\subset U_N$, we know what the recipe is, namely check if $G_N$ is easy, and if so, simply define $G_N^+\subset U_N^+$ as being the easy quantum group coming from the category $D=D_G\cap NC$.

\section*{2c. Laws of characters}

Let us discuss now some more advanced applications of Theorem 2.12, this time to the computation of the law of the main character, in the spirit of Claim 2.1. First, we have the following result, valid in the general easy quantum group setting:

\index{main character}
\index{moments}

\begin{proposition}
For an easy quantum group $G=(G_N)$, coming from a category of partitions $D=(D(k,l))$, the moments of the main character are given by
$$\int_{G_N}\chi^k=\dim\left( span\left(\xi_\pi\Big|\pi\in D(k)\right)\right)$$
where $D(k)=D(\emptyset,k)$, and with the notation $\xi_\pi=T_\pi$, for partitions $\pi\in D(k)$.
\end{proposition}

\begin{proof}
According to the Peter-Weyl theory, and to the definition of easiness, the moments of the main character are given by the following formula:
\begin{eqnarray*}
\int_{G_N}\chi^k
&=&\int_{G_N}\chi_{u^{\otimes k}}\\
&=&\dim\left(Fix(u^{\otimes k})\right)\\
&=&\dim\left( span\left(\xi_\pi\Big|\pi\in D(k)\right)\right)
\end{eqnarray*}

Thus, we obtain the formula in the statement.
\end{proof}

With the above result in hand, you would probably say very nice, so in practice, this is just a matter of counting the partitions appearing in Theorem 2.12, and then recovering the measures having these numbers as moments. However, this is wrong, because such a computation would lead to a law of $\chi$ which is independent on $N\in\mathbb N$, and for the classical groups at least, $S_N,O_N,U_N$, we obviously cannot have such a result.

\bigskip

The mistake comes from the fact that the vectors $\xi_\pi$ are not necessarily linearly independent. Let us record this finding, which will be of key importance for us:

\begin{conclusion}
The vectors associated to the partitions $\pi\in P(k)$, namely
$$\xi_\pi=\sum_{i_1\ldots i_k}\delta_\pi(i_1,\ldots,i_k)\,e_{i_1}\otimes\ldots\otimes e_{i_k}$$
are not linearly independent, with this making the main character moments for $S_N$,
$$\int_{S_N}\chi^k=\dim\left( span\left(\xi_\pi\Big|\pi\in P(k)\right)\right)$$
depend on $N\in\mathbb N$. Moreover, the same phenomenon happens for $O_N,U_N$.
\end{conclusion}

All this suggests by doing some linear algebra for the vectors $\xi_\pi$, but this looks rather complicated, and let's keep that for later. What we can do right away, instead, is that of studying $S_N$ with alternative, direct techniques. And here we have:

\index{Poisson law}
\index{inclusion-exclusion}
\index{derangement}

\begin{theorem}
Consider the symmetric group $S_N$, regarded as a compact group of matrices, $S_N\subset O_N$, via the standard permutation matrices.
\begin{enumerate}
\item The main character $\chi\in C(S_N)$, defined as usual as $\chi=\sum_iu_{ii}$, counts the number of fixed points, $\chi(\sigma)=\#\{i|\sigma(i)=i\}$.

\item The probability for a permutation $\sigma\in S_N$ to be a derangement, meaning to have no fixed points at all, becomes, with $N\to\infty$, equal to $1/e$.

\item The law of the main character $\chi\in C(S_N)$ becomes with $N\to\infty$ the Poisson law $p_1=\frac{1}{e}\sum_k\delta_k/k!$, with respect to the counting measure.
\end{enumerate}
\end{theorem}

\begin{proof}
This is something very classical, the proof being as follows:

\medskip

(1) We have indeed the following computation, which gives the result:
$$\chi(\sigma)
=\sum_iu_{ii}(\sigma)
=\sum_i\delta_{\sigma(i)i}
=\#\left\{i\Big|\sigma(i)=i\right\}$$

(2) We use the inclusion-exclusion principle. Consider the following sets:
$$S_N^i=\left\{\sigma\in S_N\Big|\sigma(i)=i\right\}$$

The probability that we are interested in is then given by:
\begin{eqnarray*}
P(\chi=0)
&=&\frac{1}{N!}\left(|S_N|-\sum_i|S_N^i|+\sum_{i<j}|S_N^i\cap S_N^j|-\sum_{i<j<k}|S_N^i\cap S_N^j\cap S_N^k|+\ldots\right)\\
&=&\frac{1}{N!}\sum_{r=0}^N(-1)^r\sum_{i_1<\ldots<i_r}(N-r)!\\
&=&\frac{1}{N!}\sum_{r=0}^N(-1)^r\binom{N}{r}(N-r)!\\
&=&\sum_{r=0}^N\frac{(-1)^r}{r!}
\end{eqnarray*}

Since we have here the expansion of $1/e$, this gives the result.

\medskip

(3) This follows by generalizing the computation in (2). To be more precise, a similar application of the inclusion-exclusion principle gives the following formula:
$$\lim_{N\to\infty}P(\chi=k)=\frac{1}{k!e}$$

Thus, we obtain in the limit a Poisson law of parameter 1, as stated.
\end{proof}

The above result is quite interesting, and tells us what to do next. As a first goal, we can try to recover (3) there by using Proposition 2.14, and easiness. Then, once this understood, we can try to look at $S_N^+$, and then at $O_N,U_N$ and $O_N^+,U_N^+$ too, with the same objective, namely finding $N\to\infty$ results for the law of $\chi$, using easiness.

\bigskip

So, back to Proposition 2.14 and Conclusion 2.15, and we have now to courageously attack the main problem, namely the linear independence question for the vectors $\xi_\pi$. This will be quite technical. Let us begin with some standard combinatorics:

\index{lattice of partitions}
\index{order of partitions}
\index{supremum of partitions}

\begin{definition}
Let $P(k)$ be the set of partitions of $\{1,\ldots,k\}$, and $\pi,\sigma\in P(k)$.
\begin{enumerate}
\item We write $\pi\leq\sigma$ if each block of $\pi$ is contained in a block of $\sigma$.

\item We let $\pi\vee\sigma\in P(k)$ be the partition obtained by superposing $\pi,\sigma$.
\end{enumerate}
Also, we denote by $|.|$ the number of blocks of the partitions $\pi\in P(k)$.
\end{definition}

As an illustration here, at $k=2$ we have $P(2)=\{||,\sqcap\}$, and we have:
$$||\leq\sqcap$$

Also, at $k=3$ we have $P(3)=\{|||,\sqcap|,\sqcap\hskip-3.2mm{\ }_|\,,|\sqcap,\sqcap\hskip-0.7mm\sqcap\}$, and the order relation is as follows:
$$|||\ \leq\ \sqcap|\ ,\ \sqcap\hskip-3.2mm{\ }_|\ ,\ |\sqcap\ \leq\ \sqcap\hskip-0.7mm\sqcap$$

In relation with our linear independence questions, the idea will be that of using:

\begin{proposition}
The Gram matrix of the vectors $\xi_\pi$ is given by the formula
$$<\xi_\pi,\xi_\sigma>=N^{|\pi\vee\sigma|}$$
where $\vee$ is the superposition operation, and $|.|$ is the number of blocks.
\end{proposition}

\begin{proof}
According to the formula of the vectors $\xi_\pi$, we have:
\begin{eqnarray*}
<\xi_\pi,\xi_\sigma>
&=&\sum_{i_1\ldots i_k}\delta_\pi(i_1,\ldots,i_k)\delta_\sigma(i_1,\ldots,i_k)\\
&=&\sum_{i_1\ldots i_k}\delta_{\pi\vee\sigma}(i_1,\ldots,i_k)\\
&=&N^{|\pi\vee\sigma|}
\end{eqnarray*}

Thus, we have obtained the formula in the statement.
\end{proof}

In order to study the Gram matrix $G_k(\pi,\sigma)=N^{|\pi\vee\sigma|}$, and more specifically to compute its determinant, we will use several standard facts about the partitions. We have:

\index{Gram matrix}
\index{M\"obius function}

\begin{definition}
The M\"obius function of any lattice, and so of $P$, is given by
$$\mu(\pi,\sigma)=\begin{cases}
1&{\rm if}\ \pi=\sigma\\
-\sum_{\pi\leq\tau<\sigma}\mu(\pi,\tau)&{\rm if}\ \pi<\sigma\\
0&{\rm if}\ \pi\not\leq\sigma
\end{cases}$$
with the construction being performed by recurrence.
\end{definition}

As an illustration here, for $P(2)=\{||,\sqcap\}$, we have by definition:
$$\mu(||,||)=\mu(\sqcap,\sqcap)=1$$

Also, $||<\sqcap$, with no intermediate partition in between, so we obtain:
$$\mu(||,\sqcap)=-\mu(||,||)=-1$$

Finally, we have $\sqcap\not\leq||$, and so we have as well the following formula:
$$\mu(\sqcap,||)=0$$

Thus, as a conclusion, we have computed the M\"obius matrix $M_2(\pi,\sigma)=\mu(\pi,\sigma)$ of the lattice $P(2)=\{||,\sqcap\}$, the formula being as follows:
$$M_2=\begin{pmatrix}1&-1\\ 0&1\end{pmatrix}$$

Back to the general case now, the main interest in the M\"obius function comes from the M\"obius inversion formula, which states that the following happens:
$$f(\sigma)=\sum_{\pi\leq\sigma}g(\pi)\quad
\implies\quad g(\sigma)=\sum_{\pi\leq\sigma}\mu(\pi,\sigma)f(\pi)$$

In linear algebra terms, the statement and proof of this formula are as follows:

\index{M\"obius inversion}

\begin{theorem}
The inverse of the adjacency matrix of $P(k)$, given by
$$A_k(\pi,\sigma)=\begin{cases}
1&{\rm if}\ \pi\leq\sigma\\
0&{\rm if}\ \pi\not\leq\sigma
\end{cases}$$
is the M\"obius matrix of $P$, given by $M_k(\pi,\sigma)=\mu(\pi,\sigma)$.
\end{theorem}

\begin{proof}
This is well-known, coming for instance from the fact that $A_k$ is upper triangular. Indeed, when inverting, we are led into the recurrence from Definition 2.19.
\end{proof}

As an illustration, for $P(2)$ the formula $M_2=A_2^{-1}$ appears as follows:
$$\begin{pmatrix}1&-1\\ 0&1\end{pmatrix}=
\begin{pmatrix}1&1\\ 0&1\end{pmatrix}^{-1}$$

Now back to our Gram matrix considerations, we have the following key result:

\begin{proposition}
The Gram matrix of the vectors $\xi_\pi$ with $\pi\in P(k)$,
$$G_{\pi\sigma}=N^{|\pi\vee\sigma|}$$
decomposes as a product of upper/lower triangular matrices, $G_k=A_kL_k$, where
$$L_k(\pi,\sigma)=
\begin{cases}
N(N-1)\ldots(N-|\pi|+1)&{\rm if}\ \sigma\leq\pi\\
0&{\rm otherwise}
\end{cases}$$
and where $A_k$ is the adjacency matrix of $P(k)$.
\end{proposition}

\begin{proof}
We have the following computation, based on Proposition 2.18:
\begin{eqnarray*}
G_k(\pi,\sigma)
&=&N^{|\pi\vee\sigma|}\\
&=&\#\left\{i_1,\ldots,i_k\in\{1,\ldots,N\}\Big|\ker i\geq\pi\vee\sigma\right\}\\
&=&\sum_{\tau\geq\pi\vee\sigma}\#\left\{i_1,\ldots,i_k\in\{1,\ldots,N\}\Big|\ker i=\tau\right\}\\
&=&\sum_{\tau\geq\pi\vee\sigma}N(N-1)\ldots(N-|\tau|+1)
\end{eqnarray*}

According now to the definition of $A_k,L_k$, this formula reads:
\begin{eqnarray*}
G_k(\pi,\sigma)
&=&\sum_{\tau\geq\pi}L_k(\tau,\sigma)\\
&=&\sum_\tau A_k(\pi,\tau)L_k(\tau,\sigma)\\
&=&(A_kL_k)(\pi,\sigma)
\end{eqnarray*}

Thus, we are led to the formula in the statement.
\end{proof}

As an illustration for the above result, at $k=2$ we have $P(2)=\{||,\sqcap\}$, and the above decomposition $G_2=A_2L_2$ appears as follows:
$$\begin{pmatrix}N^2&N\\ N&N\end{pmatrix}
=\begin{pmatrix}1&1\\ 0&1\end{pmatrix}
\begin{pmatrix}N^2-N&0\\N&N\end{pmatrix}$$

We are led in this way to the following formula, due to Lindst\"om \cite{lin}:

\index{Gram matrix}
\index{Gram determinant}
\index{Lindst\"om formula}
\index{linear independence}

\begin{theorem}
The determinant of the Gram matrix $G_k$ is given by
$$\det(G_k)=\prod_{\pi\in P(k)}\frac{N!}{(N-|\pi|)!}$$
with the convention that in the case $N<k$ we obtain $0$.
\end{theorem}

\begin{proof}
If we order $P(k)$ as usual, with respect to the number of blocks, and then lexicographically, $A_k$ is upper triangular, and $L_k$ is lower triangular. Thus, we have:
\begin{eqnarray*}
\det(G_k)
&=&\det(A_k)\det(L_k)\\
&=&\det(L_k)\\
&=&\prod_\pi L_k(\pi,\pi)\\
&=&\prod_\pi N(N-1)\ldots(N-|\pi|+1)
\end{eqnarray*}

Thus, we are led to the formula in the statement.
\end{proof}

Now back to easiness and laws of characters, we can formulate:

\index{main character}
\index{stationary convergence}
\index{asymptotic moments}

\begin{theorem}
For an easy quantum group $G=(G_N)$, coming from a category of partitions $D=(D(k,l))$, the asymptotic moments of the main character are given by
$$\lim_{N\to\infty}\int_{G_N}\chi^k=|D(k)|$$
where $D(k)=D(\emptyset,k)$, with the limiting sequence on the left consisting of certain integers, and being stationary at least starting from the $k$-th term.
\end{theorem}

\begin{proof}
We know from Proposition 2.14 that we have the following formula:
$$\int_{G_N}\chi^k
=\dim\left( span\left(\xi_\pi\Big|\pi\in D(k)\right)\right)$$

Now since by Theorem 2.22 the vectors $\xi_\pi$ are linearly independent with $N\geq k$, and in particular with $N\to\infty$, we obtain the formula in the statement.
\end{proof}

This is very nice, and as a first application, we can recover as promised the Poisson law result from Theorem 2.16, this time by using easiness, as follows:

\index{Bell numbers}
\index{Poisson law}

\begin{theorem}
For the symmetric group $S_N$, the main character becomes Poisson
$$\chi\sim p_1$$
in the $N\to\infty$ limit.
\end{theorem}

\begin{proof}
As already mentioned, this is something that we already know, from Theorem 2.16. Alternatively, according to Theorem 2.23, we have the following formula:
$$\lim_{N\to\infty}\int_{S_N}\chi^k=|P(k)|$$

Now since a partition of $\{1,\ldots,k+1\}$ appears by choosing $s$ neighbors for $1$, among the $k$ numbers available, and then partitioning the $k-s$ elements left, the numbers on the right $B_k=|P(k)|$, called Bell numbers, satisfy the following recurrence:
$$B_{k+1}=\sum_s\binom{k}{s}B_{k-s}$$

On the other hand, the moments $M_k$ of the Poisson law $p_1=\frac{1}{e}\sum_r\delta_r/r!$ are subject to the same recurrence formula, as shown by the following computation:
\begin{eqnarray*}
M_{k+1}
&=&\frac{1}{e}\sum_r\frac{(r+1)^k}{r!}\\
&=&\frac{1}{e}\sum_r\frac{r^k}{r!}\left(1+\frac{1}{r}\right)^k\\
&=&\frac{1}{e}\sum_r\frac{r^k}{r!}\sum_s\binom{k}{s}r^{-s}\\
&=&\sum_s\binom{k}{s}\cdot\frac{1}{e}\sum_r\frac{r^{k-s}}{r!}\\
&=&\sum_s\binom{k}{s}M_{k-s}
\end{eqnarray*}

As for the initial values, at $k=1,2$, these are $1,2$, for both the Bell numbers $B_k$, and the Poisson moments $M_k$. Thus we have $B_k=M_k$, which gives the result.
\end{proof}

\section*{2d. Free probability} 

Moving ahead, we have now to work out free analogues of Theorem 2.24 for the other easy quantum groups that we know. A bit of thinking at traces of unitary matrices suggests that for the groups $O_N,U_N$ we should get the real and complex normal laws. As for $O_N^+,U_N^+,S_N^+$, we are a bit in the dark here, and we can only say that we can expect to have ``free versions'' of the real and complex normal laws, and of the Poisson law.

\bigskip

Long story short, the combinatorics ahead looks quite complicated, and we are in need of a crash course on probability. So, let us start with that, classical and free probability, and we will come back later to combinatorics and quantum groups. We first have:

\index{random variable}
\index{laws of variables}

\begin{definition}
Let $A$ be a $C^*$-algebra, given with a trace $tr:A\to\mathbb C$.
\begin{enumerate}
\item The elements $a\in A$ are called random variables.

\item The moments of such a variable are the numbers $M_k(a)=tr(a^k)$.

\item The law of such a variable is the functional $\mu:P\to tr(P(a))$.
\end{enumerate}
\end{definition}

Here $k=\circ\bullet\bullet\circ\ldots$ is by definition a colored integer, and the corresponding powers $a^k$ are defined by the following formulae, and multiplicativity: 
$$a^\emptyset=1\quad,\quad
a^\circ=a\quad,\quad
a^\bullet=a^*$$

As for the polynomial $P$, this is a noncommuting $*$-polynomial in one variable:
$$P\in\mathbb C<X,X^*>$$

Observe that the law is uniquely determined by the moments, because we have:
$$P(X)=\sum_k\lambda_kX^k\implies\mu(P)=\sum_k\lambda_kM_k(a)$$

Generally speaking, the above definition is something quite abstract, but there is no other way of doing things, at least at this level of generality. However, in certain special cases, the formalism simplifies, and we recover more familiar objects, as follows:

\index{normal element}
\index{spectral measure}

\begin{proposition}
Assuming that $a\in A$ is normal, $aa^*=a^*a$, its law corresponds to a probability measure on its spectrum $\sigma(a)\subset\mathbb C$, according to the following formula:
$$tr(P(a))=\int_{\sigma(a)}P(x)d\mu(x)$$
When the trace is faithful we have $supp(\mu)=\sigma(a)$. Also, in the particular case where the variable is self-adjoint, $a=a^*$, this law is a real probability measure.
\end{proposition}

\begin{proof}
This is something very standard, coming from the Gelfand theorem, applied to the algebra $<a>$, which is commutative, and then the Riesz theorem.
\end{proof}

Following Voiculescu \cite{vo1}, we have the following two notions of independence:

\index{independence}
\index{freeness}

\begin{definition}
Two subalgebras $A,B\subset C$ are called independent when
$$tr(a)=tr(b)=0\implies tr(ab)=0$$
holds for any $a\in A$ and $b\in B$, and free when
$$tr(a_i)=tr(b_i)=0\implies tr(a_1b_1a_2b_2\ldots)=0$$
holds for any $a_i\in A$ and $b_i\in B$.
\end{definition}
 
In short, we have here a straightforward extension of the usual notion of independence, in the framework of Definition 2.25, along with a quite natural free analogue of it. In order to understand what is going on, let us first discuss some basic models for independence and freeness. We have the following result, from \cite{vo1}, which clarifies things:

\index{tensor product}
\index{free product}

\begin{proposition}
Given two algebras $(A,tr)$ and $(B,tr)$, the following hold:
\begin{enumerate}
\item $A,B$ are independent inside their tensor product $A\otimes B$.

\item $A,B$ are free inside their free product $A*B$.
\end{enumerate}
\end{proposition}

\begin{proof}
Both the assertions are clear from definitions, after some standard discussion regarding the tensor product and free product trace. See Voiculescu \cite{vo1}.
\end{proof}

In relation with groups, we have the following result:

\index{independence}
\index{freeness}
\index{group algebra}
\index{tensor product}
\index{free product}

\begin{proposition}
We have the following results, valid for group algebras:
\begin{enumerate}
\item $C^*(\Gamma),C^*(\Lambda)$ are independent inside $C^*(\Gamma\times\Lambda)$.

\item $C^*(\Gamma),C^*(\Lambda)$ are free inside $C^*(\Gamma*\Lambda)$.
\end{enumerate}
\end{proposition}

\begin{proof}
This follows from the general results in Proposition 2.28, along with the following two isomorphisms, which are both standard:
$$C^*(\Gamma\times\Lambda)=C^*(\Lambda)\otimes C^*(\Gamma)\quad,\quad 
C^*(\Gamma*\Lambda)=C^*(\Lambda)*C^*(\Gamma)$$

Alternatively, we can prove this directly, by using the fact that each algebra is spanned by the corresponding group elements, and checking the result on group elements.
\end{proof}

In order to study independence and freeness, our main tool will be: 

\index{Fourier transform}
\index{R transform}

\begin{theorem}
The convolution is linearized by the log of the Fourier transform,
$$F_f(x)=\mathbb E(e^{ixf})$$
and the free convolution is linearized by the $R$-transform, given by:
$$G_\mu(\xi)=\int_\mathbb R\frac{d\mu(t)}{\xi-t}\implies G_\mu\left(R_\mu(\xi)+\frac{1}{\xi}\right)=\xi$$
\end{theorem}

\begin{proof}
In what regards the first assertion, if $f,g$ are independent, we have indeed:
\begin{eqnarray*}
F_{f+g}(x)
&=&\int_\mathbb Re^{ixz}d(\mu_f*\mu_g)(z)\\
&=&\int_{\mathbb R\times\mathbb R}e^{ix(z+t)}d\mu_f(z)d\mu_g(t)\\
&=&\int_\mathbb Re^{ixz}d\mu_f(z)\int_\mathbb Re^{ixt}d\mu_g(t)\\
&=&F_f(x)F_g(x)
\end{eqnarray*}

As for the second assertion, here we need a good model for free convolution, and the best is to use the semigroup algebra of the free semigroup on two generators:
$$A=C^*(\mathbb N*\mathbb N)$$

Indeed, we have some freeness in the semigroup setting, a bit in the same way as for the group algebras $C^*(\Gamma*\Lambda)$, from Proposition 2.29, and in addition to this fact, and to what happens in the group algebra case, the following two key things happen:

\medskip

(1) The variables of type $S^*+f(S)$, with $S\in C^*(\mathbb N)$ being the shift, and with $f\in\mathbb C[X]$ being a polynomial, model in moments all the distributions $\mu:\mathbb C[X]\to\mathbb C$. This is indeed something elementary, which can be checked via a direct algebraic computation.

\medskip

(2) Given $f,g\in\mathbb C[X]$, the variables $S^*+f(S)$ and $T^*+g(T)$, where $S,T\in C^*(\mathbb N*\mathbb N)$ are the shifts corresponding to the generators of $\mathbb N*\mathbb N$, are free, and their sum has the same law as $S^*+(f+g)(S)$. This follows indeed by using a $45^\circ$ argument.

\medskip

With this in hand, we can see that the operation $\mu\to f$ linearizes the free convolution. We are therefore left with a computation inside $C^*(\mathbb N)$, whose conclusion is that $R_\mu=f$ can be recaptured from $\mu$ via the Cauchy transform $G_\mu$, as stated. See \cite{vo1}.
\end{proof}

As a first result now, which is central and classical and free probability, we have:

\index{CLT}
\index{Wigner law}

\begin{theorem}[CLT]
Given self-adjoint variables $x_1,x_2,x_3,\ldots$ which are i.i.d./f.i.d., centered, with variance $t>0$, we have, with $n\to\infty$, in moments,
$$\frac{1}{\sqrt{n}}\sum_{i=1}^nx_i\sim g_t/\gamma_t$$
where $g_t/\gamma_t$ are the normal and Wigner semicircle law of parameter $t$, given by:
$$g_t=\frac{1}{\sqrt{2\pi t}}e^{-x^2/2t}dx\quad,\quad 
\gamma_t=\frac{1}{2\pi t}\sqrt{4t^2-x^2}dx$$
\end{theorem}

\begin{proof}
This is routine, by using the Fourier transform and the $R$-transform.
\end{proof}

Next, we have the following complex version of the CLT:

\index{CCLT}
\index{Voiculescu law}

\begin{theorem}[CCLT]
Given variables $x_1,x_2,x_3,\ldots$ which are i.i.d./f.i.d., centered, with variance $t>0$, we have, with $n\to\infty$, in moments,
$$\frac{1}{\sqrt{n}}\sum_{i=1}^nx_i\sim G_t/\Gamma_t$$
where $G_t/\Gamma_t$ are the complex normal and Voiculescu circular law of parameter $t$, given by:
$$G_t=law\left(\frac{1}{\sqrt{2}}(a+ib)\right)\quad,\quad 
\Gamma_t=law\left(\frac{1}{\sqrt{2}}(\alpha+i\beta)\right)$$
where $a,b/\alpha,\beta$ are independent/free, each following the law $g_t/\gamma_t$.
\end{theorem}

\begin{proof}
This follows indeed from the CLT, by taking real and imaginary parts.
\end{proof}

Finally, we have the following discrete version of the CLT:

\index{PLT}
\index{Poisson limit theorem}
\index{Poisson law}
\index{Marchenko-Pastur law}

\begin{theorem}[PLT]
The following Poisson limits converge, for any $t>0$,
$$p_t=\lim_{n\to\infty}\left(\left(1-\frac{t}{n}\right)\delta_0+\frac{t}{n}\delta_1\right)^{*n}\quad,\quad 
\pi_t=\lim_{n\to\infty}\left(\left(1-\frac{t}{n}\right)\delta_0+\frac{t}{n}\delta_1\right)^{\boxplus n}$$
the limiting measures being the Poisson law $p_t$, and the Marchenko-Pastur law $\pi_t$,
$$p_t=\frac{1}{e^t}\sum_{k=0}^\infty\frac{t^k\delta_k}{k!}\quad,\quad 
\pi_t=\max(1-t,0)\delta_0+\frac{\sqrt{4t-(x-1-t)^2}}{2\pi x}\,dx$$
with at $t=1$, the Marchenko-Pastur law being $\pi_1=\frac{1}{2\pi}\sqrt{4x^{-1}-1}\,dx$. 
\end{theorem}

\begin{proof}
This is again routine, by using the Fourier and $R$-transform.
\end{proof}

This was for the basic classical and free probability. In relation now with combinatorics, we have the following result, which reminds easiness, and is of interest for us:

\begin{theorem}
The moments of the various central limiting measures, namely
$$\xymatrix@R=45pt@C=45pt{
\pi_t\ar@{-}[r]\ar@{-}[d]&\gamma_t\ar@{-}[r]\ar@{-}[d]&\Gamma_t\ar@{-}[d]\\
p_t\ar@{-}[r]&g_t\ar@{-}[r]&G_t
}$$
are always given by the same formula, involving partitions, namely
$$M_k=\sum_{\pi\in D(k)}t^{|\pi|}$$
with the sets of partitions $D(k)$ in question being respectively
$$\xymatrix@R=50pt@C=50pt{
NC\ar[d]&NC_2\ar[d]\ar[l]&\mathcal{NC}_2\ar[l]\ar[d]\\
P&P_2\ar[l]&\mathcal P_2\ar[l]}$$
and with $|.|$ being the number of blocks. 
\end{theorem}

\begin{proof}
This follows indeed from the various computations leading to Theorems 2.31, 2.32, 2.33, and details can be found in any free probability book. See \cite{nsp}, \cite{vdn}.
\end{proof}

It is possible to say more on this, following Rota in the classical case, Speicher in the free case, and Bercovici-Pata for the classical/free correspondence. We first have:

\begin{definition}
The cumulants of a self-adjoint variable $a\in A$ are given by
$$\log F_a(\xi)=\sum_{n=1}^\infty k_n(a)\,\frac{(i\xi)^n}{n!}$$ 
and the free cumulants of the same variable $a\in A$ are given by:
$$R_a(\xi)=\sum_{n=1}^\infty\kappa_n(a)\xi^{n-1}$$
Moreover, we have extensions of these notions to the non-self-adjoint case.
\end{definition}

In what follows we will only discuss the self-adjoint case, which is simpler, and illustrating. Since the classical and free cumulants are by definition certain linear combinations of the moments, we should have conversion formulae. The result here is as follows:

\begin{theorem}
The moments can be recaptured out of cumulants via
$$M_n(a)=\sum_{\pi\in P(n)}k_\pi(a)\quad,\quad M_n(a)=\sum_{\pi\in NC(n)}\kappa_\pi(a)$$
with the convention that $k_\pi,\kappa_\pi$ are defined by multiplicativity over blocks. Also,
$$k_n(a)=\sum_{\nu\in P(n)}\mu_P(\nu,1_n)M_\nu(a)\quad,\quad \kappa_n(a)=\sum_{\nu\in NC(n)}\mu_{NC}(\nu,1_n)M_\nu(a)$$
where $\mu_P,\mu_{NC}$ are the M\"obius functions of $P(n),NC(n)$.
\end{theorem}

\begin{proof}
Here the first formulae follow from Definition 2.35, by doing some combinatorics, and the second formulae follow from them, via M\"obius inversion.
\end{proof}

In relation with the various laws that we are interested in, we have:

\begin{proposition}
The classical and free cumulants are as follows:
\begin{enumerate}
\item For $\mu=\delta_c$ both the classical and free cumulants are $c,0,0,\ldots$

\item For $\mu=g_t/\gamma_t$ the classical/free cumulants are $0,t,0,0,\ldots$

\item For $\mu=p_t/\pi_t$ the classical/free cumulants are $t,t,t,\ldots$
\end{enumerate}
\end{proposition}

\begin{proof}
Here (1) is something trivial, and (2,3) can be deduced either directly, starting from the definition of the various laws involved, or by using Theorem 2.34.
\end{proof}

Following now Bercovici-Pata \cite{bpa}, let us formulate the following definition:

\begin{definition}
If the classical cumulants of $\eta$ equal the free cumulants of  $\mu$,
$$k_n(\eta)=\kappa_n(\mu)$$
we say that $\eta$ is the classical version of $\mu$, and that $\mu$ is the free version of $\eta$.
\end{definition}

All this is quite interesting, and we have now a better understanding of Theorem 2.34, the point there being that on the vertical, we have measures in Bercovici-Pata bijection. Now back to quantum groups, we first have the following result, from \cite{bc1}:

\index{Gaussian law}
\index{Poisson law}
\index{Wigner law}
\index{Marchanko-Pastur law}

\begin{theorem}
The asymptotic laws of characters for the basic quantum groups,
$$\xymatrix@R=15mm@C=15mm{
S_N^+\ar[r]&O_N^+\ar[r]&U_N^+\\
S_N\ar[r]\ar[u]&O_N\ar[r]\ar[u]&U_N\ar[u]}$$
are precisely the main laws in classical and free probability at $t=1$.
\end{theorem}

\begin{proof}
This follows indeed from our various easiness considerations before, and from Theorem 2.34 applied at $t=1$, which gives $M_k=|D(k)|$ in this case.
\end{proof}

More generally, again following \cite{bc1}, let us discuss now the computation for the truncated characters. These are variables constructed as follows:

\index{truncated character}
\index{main character}

\begin{definition}
Associated to any Woronowicz algebra $(A,u)$ are the variables
$$\chi_t=\sum_{i=1}^{[tN]}u_{ii}$$
depending on a parameter $t\in(0,1]$, called truncations of the main character.
\end{definition}

In order to understand what these variables $\chi_t$ are about, let us first investigate the symmetric group $S_N$. We have here the following result:

\index{Stirling numbers}

\begin{theorem}
For the symmetric group $S_N\subset O_N$, the truncated character
$$\chi_t(g)=\sum_{i=1}^{[tN]}u_{ii}$$
becomes, with $N\to\infty$, a Poisson variable of parameter $t$. 
\end{theorem}

\begin{proof}
This can be deduced via inclusion-exclusion, as in the proof of Theorem 2.16, but let us prove this via an alternative method, which is instructive as well. Our first claim is that the integrals over $S_N$ are given by the following formula:
$$\int_{S_N}u_{i_1j_1}\ldots u_{i_kj_k}=\begin{cases}
\frac{(N-|\ker i|)!}{N!}&{\rm if}\ \ker i=\ker j\\
0&{\rm otherwise}
\end{cases}$$

Indeed, according to the definition of $u_{ij}$, the above integrals are given by:
$$\int_{S_N}u_{i_1j_1}\ldots u_{i_kj_k}=\frac{1}{N!}\#\left\{\sigma\in S_N\Big|\sigma(j_1)=i_1,\ldots,\sigma(j_k)=i_k\right\}$$

But this proves our claim. Now with the above formula in hand, with $S_{kb}$ being the Stirling numbers, counting the partitions in $P(k)$ having $b$ blocks, we have:
\begin{eqnarray*}
\int_{S_N}\chi_t^k
&=&\sum_{i_1\ldots i_k=1}^{[tN]}\int_{S_N}u_{i_1i_1}\ldots u_{i_ki_k}\\
&=&\sum_{\pi\in P(k)}\frac{[tN]!}{([tN]-|\pi|!)}\cdot\frac{(N-|\pi|!)}{N!}\\
&=&\sum_{b=1}^{[tN]}\frac{[tN]!}{([tN]-b)!}\cdot\frac{(N-b)!}{N!}\cdot S_{kb}
\end{eqnarray*}

Thus with $N\to\infty$ the moments are $M_k\simeq\sum_{b=1}^kS_{kb}t^b$, which gives the result.
\end{proof}

Summarizing, we have nice results about $S_N$. In general, however, and in particular for $O_N,U_N$ and $S_N^+,O_N^+,U_N^+$, there is no simple trick as for $S_N$, and we must use general integration methods, from \cite{bc1}, \cite{csn}. We have here the following formula:

\index{Weingarten formula}
\index{Tannakian category}
\index{easiness}
\index{Gram matrix}

\begin{theorem}
For an easy quantum group $G\subset_uU_N^+$, coming from a category of partitions $D=(D(k,l))$, we have the Weingarten integration formula
$$\int_Gu_{i_1j_1}\ldots u_{i_kj_k}=\sum_{\pi,\sigma\in D(k)}\delta_\pi(i)\delta_\sigma(j)W_{kN}(\pi,\sigma)$$
where $D(k)=D(\emptyset,k)$, $\delta$ are usual Kronecker symbols, and $W_{kN}=G_{kN}^{-1}$, with 
$$G_{kN}(\pi,\sigma)=N^{|\pi\vee\sigma|}$$
where $|.|$ is the number of blocks. 
\end{theorem}

\begin{proof}
This follows from the general Weingarten formula from Theorem 2.4. Indeed, in the easy case we can take $D_k=D(k,k)$, and the Kronecker symbols are given by:
$$\delta_{\xi_\pi}(i)
=<\xi_\pi,e_{i_1}\otimes\ldots\otimes e_{i_k}>
=\delta_\pi(i_1,\ldots,i_k)$$

The Gram matrix being as well the correct one, we obtain the result. See \cite{bc1}.
\end{proof}

With the above formula in hand, we can go back to the question of computing the laws of truncated characters. First, we have the following moment formula, from \cite{bc1}:

\index{truncated character}

\begin{proposition}
The moments of truncated characters are given by the formula
$$\int_G(u_{11}+\ldots +u_{ss})^k=Tr(W_{kN}G_{ks})$$
where $G_{kN}$ and $W_{kN}=G_{kN}^{-1}$ are the associated Gram and Weingarten matrices.
\end{proposition}

\begin{proof}
We have indeed the following computation:
\begin{eqnarray*}
\int_G(u_{11}+\ldots +u_{ss})^k
&=&\sum_{i_1=1}^{s}\ldots\sum_{i_k=1}^s\int u_{i_1i_1}\ldots u_{i_ki_k}\\
&=&\sum_{\pi,\sigma\in D(k)}W_{kN}(\pi,\sigma)\sum_{i_1=1}^{s}\ldots\sum_{i_k=1}^s\delta_\pi(i)\delta_\sigma(i)\\
&=&\sum_{\pi,\sigma\in D(k)}W_{kN}(\pi,\sigma)G_{ks}(\sigma,\pi)\\
&=&Tr(W_{kN}G_{ks})
\end{eqnarray*}

Thus, we have obtained the formula in the statement.
\end{proof}

In order to process now the above formula, things are quite technical, and won't work well in general. We must impose here a uniformity condition, as follows:

\index{uniformity}

\begin{theorem}
For an easy quantum group $G=(G_N)$, coming from a category of partitions $D\subset P$, the following conditions are equivalent:
\begin{enumerate}
\item $G_{N-1}=G_N\cap U_{N-1}^+$, via the embedding $U_{N-1}^+\subset U_N^+$ given by $u\to diag(u,1)$.

\item $G_{N-1}=G_N\cap U_{N-1}^+$, via the $N$ possible diagonal embeddings $U_{N-1}^+\subset U_N^+$.

\item $D$ is stable under the operation which consists in removing blocks.
\end{enumerate}
If these conditions are satisfied, we say that $G=(G_N)$ is uniform.
\end{theorem}

\begin{proof}
This is something very standard, the idea being as follows:

\medskip

$(1)\iff(2)$ This equivalence is elementary, coming from the inclusion $S_N\subset G_N$, which makes everything $S_N$-invariant. 

\medskip

$(1)\iff(3)$ Given a closed subgroup $K\subset U_{N-1}^+$, with fundamental corepresentation $u$, consider the following $N\times N$ matrix:
$$v=\begin{pmatrix}u&0\\0&1\end{pmatrix}$$

Then for any $\pi\in P(k)$ a standard computation shows that we have:
$$\xi_\pi\in Fix(v^{\otimes k})\iff\xi_{\pi'}\in Fix(v^{\otimes k'}),\,\forall\pi'\in P(k'),\pi'\subset\pi$$

Now with this in hand, the result follows from Tannakian duality.
\end{proof}
 
By getting back now to the truncated characters, we have the following result:

\begin{theorem}
For a uniform easy quantum group $G=(G_N)$, we have the formula
$$\lim_{N\to\infty}\int_{G_N}\chi_t^k=\sum_{\pi\in D(k)}t^{|\pi|}$$
with $D\subset P$ being the associated category of partitions.
\end{theorem}

\begin{proof}
In the uniform case the Gram matrix, and so the Weingarten matrix too, are asymptotically diagonal, so the asymptotic moments are given by:
$$\int_{G_N}\chi_t^k
=Tr(W_{kN}G_{k[tN]})
\simeq\sum_{\pi\in D(k)}N^{-|\pi|}[tN]^{|\pi|}
\simeq\sum_{\pi\in D(k)}t^{|\pi|}$$

Thus, we are led to the conclusion in the statement. See \cite{bc1}, \cite{bsp}.
\end{proof}

We can now improve our quantum group results, as follows:

\index{Gaussian law}
\index{Poisson law}
\index{Wigner law}
\index{Marchanko-Pastur law}

\begin{theorem}
The asymptotic laws of truncated characters for the quantum groups
$$\xymatrix@R=15mm@C=15mm{
S_N^+\ar[r]&O_N^+\ar[r]&U_N^+\\
S_N\ar[r]\ar[u]&O_N\ar[r]\ar[u]&U_N\ar[u]}$$
are precisely the main limiting laws in classical and free probability, namely:
$$\xymatrix@R=45pt@C=45pt{
\pi_t\ar@{-}[r]\ar@{-}[d]&\gamma_t\ar@{-}[r]\ar@{-}[d]&\Gamma_t\ar@{-}[d]\\
p_t\ar@{-}[r]&g_t\ar@{-}[r]&G_t
}$$
\end{theorem}

\begin{proof}
This follows indeed from easiness, Theorem 2.34 and Theorem 2.45.
\end{proof}

\section*{2e. Exercises} 

Generally speaking, as a best complement to the above material, we recommend some probability reading. Here is a first exercise, in relation with the above:

\begin{exercise}
Work out all details for the classical and free CLT, CCLT, PLT.
\end{exercise}

Here the tools were discussed in the above, and left to do are some computations.

\begin{exercise}
Prove that $S_N$ is easy, directly.
\end{exercise}

To be more precise, our proof for $S_N$ was based on the fact that $S_N^+$ is easy. The problem is that of finding a direct proof, with no reference to quantum groups.

\begin{exercise}
Try finding Gram determinant formulae for the groups $O_N,U_N$ and for the quantum groups $S_N^+,O_N^+,U_N^+$, complementing the Lindst\"om formula for $S_N$.
\end{exercise}

This is actually a quite difficult exercise, but we will be back to this.

\chapter{Representation theory}

\section*{3a. Rotation groups}

We have seen that the inclusion $S_N\subset S_N^+$, and its companion inclusions $O_N\subset O_N^+$ and $U_N\subset U_N^+$, are all liberations in the sense of easy quantum group theory, and that some representation theory consequences, in the $N\to\infty$ limit, can be derived from this. We discuss here the case where $N\in\mathbb N$ is fixed, which is more technical.

\bigskip

Let us first study the representations of $O_N^+$. We know that in the $N\to\infty$ limit we have $\chi\sim\gamma_1$, and as a first question, we would like to know how the irreducible representations of a ``formal quantum group'' should look like, when subject to the condition $\chi\sim\gamma_1$. And fortunately, the answer here is very simple, coming from $SU_2$:

\index{Catalan number}
\index{Wigner law}
\index{Clebsch-Gordan rules}
\index{rotation group}

\begin{theorem}
The group $SU_2$ is as follows:
\begin{enumerate}
\item The main character is real, its odd moments vanish, and its even moments are the Catalan numbers:
$$\int_{SU_2}\chi^{2k}=C_k$$

\item This main character follows the Wigner semicircle law, $\chi\sim\gamma_1$.

\item The irreducible representations can be labelled by positive integers, $r_k$ with $k\in\mathbb N$, and the fusion rules for these representations are:
$$r_k\otimes r_l=r_{|k-l|}+r_{|k-l|+2}+\ldots+r_{k+l}$$ 

\item The dimensions of these representations are $\dim r_k=k+1$.
\end{enumerate}
\end{theorem}

\begin{proof}
There are many possible proofs here, the idea being as follows:

\medskip

(1,2) These statements are equivalent, and in order to prove them, a simple argument is by using the well-known isomorphism $SU_2\simeq S^3_\mathbb R$, coming from:
$$SU_2=\left\{\begin{pmatrix}x+iy&z+it\\ -z+it&x-iy\end{pmatrix}\Big|x^2+y^2+z^2+t^2=1\right\}$$

Indeed, in this picture the moments of $\chi=2x$ can be computed via spherical coordinates and some calculus, and follow to be the Catalan numbers:
$$C_k=\frac{1}{k+1}\binom{2k}{k}$$

As for the formula $\chi\sim\gamma_1$, this follows from this, and is geometrically clear as well.

\medskip

(3,4) Our claim is that we can construct, by recurrence on $k\in\mathbb N$, a sequence $r_k$ of irreducible, self-adjoint and distinct representations of $SU_2$, satisfying:
$$r_0=1\quad,\quad
r_1=u\quad,\quad 
r_{k-1}\otimes r_1=r_{k-2}+r_k$$

Indeed, assume that $r_0,\ldots,r_{k-1}$ are constructed, and let us construct $r_k$. We have:
$$r_{k-2}\otimes r_1=r_{k-3}+r_{k-1}$$

Thus $r_{k-1}\subset r_{k-2}\otimes r_1$, and since $r_{k-2}$ is irreducible, by Frobenius we have:
$$r_{k-2}\subset r_{k-1}\otimes r_1$$

We conclude there exists a certain representation $r_k$ such that:
$$r_{k-1}\otimes r_1=r_{k-2}+r_k$$

By recurrence, $r_k$ is self-adjoint. Now observe that according to our recurrence formula, we can split $u^{\otimes k}$ as a sum of the following type, with positive coefficients:  
$$u^{\otimes k}=c_kr_k+c_{k-2}r_{k-2}+c_{k-4}r_{k-4}+\ldots$$

We conclude by Peter-Weyl that we have an inequality as follows, with equality precisely when $r_k$ is irreducible, and non-equivalent to the other summands $r_i$:
$$\sum_ic_i^2\leq\dim(End(u^{\otimes k}))$$

But by (1) the number on the right is $C_k$, and some straightforward combinatorics, based on the fusion rules, shows that the number on the left is $C_k$ as well. Thus:
$$C_k=\sum_ic_i^2\leq\dim(End(u^{\otimes k}))=\int_{SU_2}\chi^{2k}=C_k$$

We conclude that we have equality in our estimate, so our representation $r_k$ is irreducible, and non-equivalent to $r_{k-2},r_{k-4},\ldots$ Moreover, this representation $r_k$ is not equivalent to $r_{k-1},r_{k-3},\ldots$ either, with this coming from $r_p\subset u^{\otimes p}$, and from:
$$\dim(Fix(u^{\otimes 2s+1}))=\int_{SU_2}\chi^{2s+1}=0$$

Thus, we have proved our claim. Now since each irreducible representation of $SU_2$ must appear in some tensor power $u^{\otimes k}$, and we know how to decompose each $u^{\otimes k}$ into sums of representations $r_k$, these representations $r_k$ are all the irreducible representations of $SU_2$, and we are done with (3). As for the formula in (4), this is clear.
\end{proof}

There are of course many other proofs for the above result, which are all instructive, and we recommend here any good book on geometry and physics. In what concerns us, the above will do, and we will be back to this later, with some further comments.

\bigskip

Getting back now to $O_N^+$, we know that in the $N\to\infty$ limit we have $\chi\sim\gamma_1$, so by the above when formally setting $N=\infty$, the fusion rules are the same as for $SU_2$. Miraculously, however, this happens in fact at any $N\geq2$, the result being as follows:

\index{free orthogonal group}
\index{Catalan number}
\index{Wigner law}
\index{Clebsch-Gordan rules}

\begin{theorem}
The quantum groups $O_N^+$ with $N\geq2$ are as follows:
\begin{enumerate}
\item The odd moments of the main character vanish, and the even moments are:
$$\int_{O_N^+}\chi^{2k}=C_k$$

\item This main character follows the Wigner semicircle law, $\chi\sim\gamma_1$.

\item The fusion rules for irreducible representations are as for $SU_2$, namely:
$$r_k\otimes r_l=r_{|k-l|}+r_{|k-l|+2}+\ldots+r_{k+l}$$ 

\item We have $\dim r_k=(q^{k+1}-q^{-k-1})/(q-q^{-1})$, with $q+q^{-1}=N$.
\end{enumerate}
\end{theorem}

\begin{proof}
The idea is to skilfully recycle the proof of Theorem 3.1, as follows:

\medskip

(1,2) These assertions are equivalent, and since we cannot prove them directly, we will simply say that these follow from the combinatorics in (3) below.

\medskip

(3,4) As before, our claim is that we can construct, by recurrence on $k\in\mathbb N$, a sequence $r_0,r_1,r_2,\ldots$ of irreducible, self-adjoint and distinct representations of $O_N^+$, satisfying:
$$r_0=1\quad,\quad
r_1=u\quad,\quad 
r_{k-1}\otimes r_1=r_{k-2}+r_k$$

In order to do so, we can use as before $r_{k-2}\otimes r_1=r_{k-3}+r_{k-1}$ and Frobenius, and we conclude there exists a certain representation $r_k$ such that:
$$r_{k-1}\otimes r_1=r_{k-2}+r_k$$

As a first observation, $r_k$ is self-adjoint, because its character is a certain polynomial with integer coefficients in $\chi$, which is self-adjoint. In order to prove now that $r_k$ is irreducible, and non-equivalent to $r_0,\ldots,r_{k-1}$, let us split as before $u^{\otimes k}$, as follows:  
$$u^{\otimes k}=c_kr_k+c_{k-2}r_{k-2}+c_{k-4}r_{k-4}+\ldots$$

The point now is that we have the following equalities and inequalities:
$$C_k=\sum_ic_i^2\leq\dim(End(u^{\otimes k}))\leq |NC_2(k,k)|=C_k$$

Indeed, the equality at left is clear as before, then comes a standard inequality, then an inequality coming from easiness, then a standard equality. Thus, we have equality, so $r_k$ is irreducible, and non-equivalent to $r_{k-2},r_{k-4},\ldots$ Moreover, $r_k$ is not equivalent to $r_{k-1},r_{k-3},\ldots$ either, by using the same argument as for $SU_2$, and the end of the proof of (3) is exactly as for $SU_2$. As for (4), by recurrence we obtain, with $q+q^{-1}=N$:
$$\dim r_k=q^k+q^{k-2}+\ldots+q^{-k+2}+q^{-k}$$

But this gives the dimension formula in the statement, and we are done.
\end{proof}

The above result raises several interesting questions. For instance we would like to know if Theorem 3.1 can be unified with Theorem 3.2. Also, combinatorially speaking, we would like to have a better understanding of the ``miracle'' making Theorem 3.2 hold at any $N\geq2$, instead of $N=\infty$ only. These questions will be answered in due time.

\bigskip

Regarding now the quantum group $U_N^+$, a similar result holds here, which is also elementary, using only algebraic techniques, based on easiness. Let us start with:

\begin{theorem}
We have isomorphisms as follows,
$$U_N^+=\widetilde{O_N^+}\quad,\quad PO_N^+=PU_N^+$$
modulo the usual equivalence relation for compact quantum groups.
\end{theorem}

\begin{proof}
The above isomorphisms both come from easiness, as follows:

\medskip

(1) We have embeddings as follows, with the first one coming by using the counit, and with the second one coming from the universality property of $U_N^+$:
$$O_N^+
\subset\widetilde{O_N^+}
\subset U_N^+$$

We must prove that the embedding on the right is an isomorphism. In order to do so, let us denote by $v,zv,u$ the fundamental representations of the above quantum groups. At the level of the associated Hom spaces we obtain reverse inclusions, as follows:
$$Hom(v^{\otimes k},v^{\otimes l})
\supset Hom((zv)^{\otimes k},(zv)^{\otimes l})
\supset Hom(u^{\otimes k},u^{\otimes l})$$

But the spaces on the left and on the right are known from chapter 2, the easiness result there stating that these are as follows:
$$span\left(T_\pi\Big|\pi\in NC_2(k,l)\right)\supset span\left(T_\pi\Big|\pi\in\mathcal{NC}_2(k,l)\right)$$

Regarding the spaces in the middle, these are obtained from those on the left by coloring, and we obtain the same spaces as those on the right. Thus, by Tannakian duality, our embedding $\widetilde{O_N^+}\subset U_N^+$ is an isomorphism, modulo the usual equivalence relation.

\medskip

(2) Regarding now the projective versions, the result here follows from:
$$PU_N^+=P\widetilde{O_N^+}=PO_N^+$$

Alternatively, with the notations in the proof of (1), we have:
$$Hom\left((v\otimes v)^k,(v\otimes v)^l\right)=span\left(T_\pi\Big|\pi\in NC_2((\circ\bullet)^k,(\circ\bullet)^l)\right)$$
$$Hom\left((u\otimes\bar{u})^k,(u\otimes\bar{u})^l\right)=span\left(T_\pi\Big|\pi\in \mathcal{NC}_2((\circ\bullet)^k,(\circ\bullet)^l)\right)$$

The sets on the right being equal, we conclude that the inclusion $PO_N^+\subset PU_N^+$ preserves the corresponding Tannakian categories, and so must be an isomorphism.
\end{proof}

Getting now to the representations of $U_N^+$, the result here is as follows:

\index{free unitary group}
\index{Voiculescu law}
\index{Clebsch-Gordan rules}

\begin{theorem}
The quantum groups $U_N^+$ with $N\geq2$ are as follows:
\begin{enumerate}
\item The moments of the main character count the matching pairings:
$$\int_{U_N^+}\chi^k=|\mathcal NC_2(k)|$$

\item The main character follows the Voiculescu circular law of parameter $1$:
$$\chi\sim\Gamma_1$$

\item The irreducible representations are indexed by $\mathbb N*\mathbb N$, with as fusion rules:
$$r_k\otimes r_l=\sum_{k=xy,l=\bar{y}z}r_{xz}$$

\item The corresponding dimensions $\dim r_k$ can be computed by recurrence.
\end{enumerate}
\end{theorem}

\begin{proof}
There are several proofs here, the idea being as follows:

\medskip

(1) The original proof, explained for instance in \cite{ba4}, is by construcing the representations $r_k$ by recurrence, exactly as in the proof of Theorem 3.2, and then arguing, also as there, that the combinatorics found proves the first two assertions as well. In short, what we have is a ``complex remake'' of Theorem 3.2, which can be proved in a similar way.

\medskip

(2) An alternative argument, discussed as well in \cite{ba4}, is by using Theorem 3.3. Indeed, the fusion rules for $U_N^+=\widetilde{O_N^+}$ can be computed  by using those of $O_N^+$, and we end up with the above ``free complexification'' of the Clebsch-Gordan rules. As for the first two assertions, these follow too from $U_N^+=\widetilde{O_N^+}$, via standard free probability.
\end{proof}

As a conclusion, our results regarding $O_N^+,U_N^+$ show that the $N\to\infty$ convergence of the law of the main character to $\gamma_1,\Gamma_1$, known since chapter 2, is in fact stationary, starting with $N=2$. And this is quite a miracle, for instance because for $O_N,U_N$, some elementary computations show that the same $N\to\infty$ convergence, this time to the normal laws $g_1,G_1$, is far from being stationary. Thus, it is tempting to formulate:

\begin{conclusion}
The free world is simpler than the classical world.
\end{conclusion}

And please don't get me wrong, especially if you're new to the subject, having struggled with the free material explained so far in this book. What I'm saying here is that, once you're reasonably advanced, and familiar with freeness, and so you will be soon, a second look at what has been said so far in this book can only lead to the above conclusion.

\bigskip

More on this later, in connection with permutations and quantum permutations too. Finally, as an extra piece of evidence, we have the isomorphism $PO_N^+=PU_N^+$ from Theorem 3.3, which is something quite intruiguing too, suggesting that the ``free projective geometry is scalarless''. We will be back to this later, with the answer that yes, free projective geometry is indeed scalarless, simpler than classical projective geometry.

\section*{3b. Clebsch-Gordan rules}

We discuss now the representation theory of $S_N^+$ at $N\geq4$. Let us begin our study exactly as for $O_N^+$. We know that in the $N\to\infty$ limit we have $\chi\sim\pi_1$, and as a first question, we would like to know how the irreducible representations of a ``formal quantum group'' should look like, when subject to the condition $\chi\sim\pi_1$. And fortunately, the answer here is very simple, involving this time the group $SO_3$:

\index{Catalan numbers}
\index{quantum permutation}
\index{Clebsch-Gordan rules}
\index{fusion rules}
\index{Marchenko-Pastur law}

\begin{theorem}
The group $SO_3$ is as follows:
\begin{enumerate}
\item The moments of the main character are the Catalan numbers:
$$\int_{SO_3}\chi^k=C_k$$

\item The main character follows the Marchenko-Pastur law of parameter $1$:
$$\chi\sim\pi_1$$

\item The fusion rules for irreducible representations are as follows:
$$r_k\otimes r_l=r_{|k-l|}+r_{|k-l|+1}+\ldots+r_{k+l}$$ 
 
\item The dimensions of these representations are $\dim r_k=2k-1$.
\end{enumerate}
\end{theorem}

\begin{proof}
As before with $SU_2$, there are many possible proofs here, which are all instructive. Here is our take on the subject, in the spirit of our proof for $SU_2$:

\medskip

(1,2) These statements are equivalent, and in order to prove them, a simple argument is by using the $SU_2$ result, and the double cover map $SU_2\to SO_3$. Indeed, let us recall from the proof for $SU_2$ that we have an isomorphism $SU_2\simeq S^3_\mathbb R$, coming from:
$$SU_2=\left\{\begin{pmatrix}x+iy&z+it\\ -z+it&x-iy\end{pmatrix}\Big|x^2+y^2+z^2+t^2=1\right\}$$

The point now is that we have a double cover map $SU_2\to SO_3$, which gives the following formula for the generic elements of $SO_3$, called Euler-Rodrigues formula:
$$U=\begin{pmatrix}
x^2+y^2-z^2-t^2&2(yz-xt)&2(xz+yt)\\
2(xt+yz)&x^2+z^2-y^2-t^2&2(zt-xy)\\
2(yt-xz)&2(xy+zt)&x^2+t^2-y^2-z^2
\end{pmatrix}$$

It follows that the main character of $SO_3$ is given by the following formula:
\begin{eqnarray*}
\chi(U)
&=&Tr(U)+1\\
&=&3x^2-y^2-z^2-t^2+1\\
&=&4x^2
\end{eqnarray*}

On the other hand, we know from Theorem 3.1 and its proof that $2x\sim\gamma_1$. Now since we have $f\sim\gamma_1\implies f^2\sim\pi_1$, we obtain $\chi\sim\pi_1$, as desired.

\medskip

(3,4) Our claim is that we can construct, by recurrence on $k\in\mathbb N$, a sequence $r_k$ of irreducible, self-adjoint and distinct representations of $SO_3$, satisfying:
$$r_0=1\quad,\quad
r_1=u-1\quad,\quad 
r_{k-1}\otimes r_1=r_{k-2}+r_{k-1}+r_k$$

Indeed, assume that $r_0,\ldots,r_{k-1}$ are constructed, and let us construct $r_k$. The Frobenius trick from the proof for $SU_2$ will no longer work, as you can verify yourself, so we have to invoke (1). To be more precise, by integrating characters we obtain:
$$r_{k-1},r_{k-2}\subset r_{k-1}\otimes r_1$$

Thus, there exists a representation $r_k$ such that:
$$r_{k-1}\otimes r_1=r_{k-2}+r_{k-1}+r_k$$

Once again by integrating characters, we conclude that $r_k$ is irreducible, and non-equivalent to $r_1,\ldots,r_{k-1}$, and this proves our claim. Also, since any irreducible representation of $SO_3$ must appear in some tensor power of $u$, and we can decompose each $u^{\otimes k}$ into sums of representations $r_p$, we conclude that these representations $r_p$ are all the irreducible representations of $SO_3$. Finally, the dimension formula is clear.
\end{proof} 

Based on the above result, and on what we know about the relation between $SU_2$ and the quantum groups $O_N^+$ at $N\geq2$, we can safely conjecture that the fusion rules for $S_N^+$ at $N\geq4$ should be the same as for $SO_3$. However, a careful inspection of the proof of Theorem 3.6 shows that, when trying to extend it to $S_4^+$, a bit in the same way as the proof of Theorem 3.1 was extended to $O_N^+$, we run into a serious problem, namely: 

\begin{problem}
Regarding $S_N^+$ with $N\geq4$, we can't get away with the estimate
$$\int_{S_N^+}\chi^k\leq C_k$$
because the Frobenius trick won't work. We need equality in this estimate.
\end{problem}

To be more precise, the above estimate comes from easiness, and we have seen that for $O_N^+$ with $N\geq2$, a similar easiness estimate, when coupled with the Frobenius trick, does the job. However, the proof of Theorem 3.6 makes it clear that no Frobenius trick is available, and so we need equality in the above estimate, as indicated.

\bigskip

So, how to prove the equality? The original argument, from \cite{ba1}, is something quick and advanced, saying that modulo some standard identifications, we are in need of the fact that the trace on the Temperley-Lieb algebra $TL_N(k)=span(NC_2(k,k))$ is faithful at index values $N\geq4$, and with this being true by the results of Jones in \cite{jo1}. However, while very quick, this remains something advanced, because the paper \cite{jo1} itself is based on a good deal of von Neumann algebra theory, covering a whole book or so. And so, we don't want to get into this, at least at this stage of our presentation.

\bigskip

In short, we are a bit in trouble. But no worries, there should be a pedestrian way of solving our problem, because that is how reasonable mathematics is made, always available to pedestrians. Here is an idea for a solution, which is a no-brainer:

\begin{solution}
We can get the needed equality at $N\geq4$, namely
$$\int_{S_N^+}\chi^k=C_k$$
by proving that the vectors $\{\xi_\pi|\pi\in NC(k)\}$ are linearly independent.
\end{solution}

Indeed, this is something coming from easiness, and since this problem does not look that scary, let us try to solve it. As a starting point for our study, we have:

\begin{proposition}
The following are linearly independent, at any $N\geq2$:
\begin{enumerate}
\item The linear maps $\left\{T_\pi\big|\pi\in NC_2(k,l)\right\}$, with $k+l\in 2\mathbb N$.

\item The vectors $\left\{\xi_\pi\big|\pi\in NC_2(2k)\right\}$, with $k\in\mathbb N$.

\item The linear maps $\left\{T_\pi\big|\pi\in NC_2(k,k)\right\}$, with $k\in\mathbb N$.
\end{enumerate}
\end{proposition}

\begin{proof}
All this follows from the dimension equalities established in the proof of Theorem 3.2, because in all cases, the number of partitions is a Catalan number.
\end{proof}

In order to pass now to quantum permutations, we can use the following trick:

\index{fattening of partitions}
\index{shrinking of partitions}

\begin{proposition}
We have a bijection $NC(k)\simeq NC_2(2k)$, constructed by fattening and shrinking, as follows:
\begin{enumerate}
\item The application $NC(k)\to NC_2(2k)$ is the ``fattening'' one, obtained by doubling all the legs, and doubling all the strings too.

\item Its inverse $NC_2(2k)\to NC(k)$ is the ``shrinking'' application, obtained by collapsing pairs of consecutive neighbors.
\end{enumerate}
\end{proposition}

\begin{proof}
The fact that the above two operations are indeed inverse to each other is clear, by drawing pictures, and computing the corresponding compositions.
\end{proof}

At the level of the associated Gram matrices, the result is as follows:

\index{Gram matrix}

\begin{proposition}
The Gram matrices of $NC_2(2k)\simeq NC(k)$ are related by
$$G_{2k,n}(\pi,\sigma)=n^k(\Delta_{kn}^{-1}G_{k,n^2}\Delta_{kn}^{-1})(\pi',\sigma')$$
where $\pi\to\pi'$ is the shrinking operation, and $\Delta_{kn}$ is the diagonal of $G_{kn}$.
\end{proposition}

\begin{proof}
In the context of the bijection from Proposition 3.10, we have:
$$|\pi\vee\sigma|=k+2|\pi'\vee\sigma'|-|\pi'|-|\sigma'|$$

We therefore have the following formula, valid for any $n\in\mathbb N$:
$$n^{|\pi\vee\sigma|}=n^{k+2|\pi'\vee\sigma'|-|\pi'|-|\sigma'|}$$

Thus, we are led to the formula in the statement.
\end{proof}

We can now formulate a ``projective'' version of Proposition 3.9, as follows:

\begin{proposition}
The following are linearly independent, for $N=n^2$ with $n\geq2$:
\begin{enumerate}
\item The linear maps $\left\{T_\pi\big|\pi\in NC(k,l)\right\}$, with $k,l\in 2\mathbb N$.

\item The vectors $\left\{\xi_\pi\big|\pi\in NC(k)\right\}$, with $k\in\mathbb N$.

\item The linear maps $\left\{T_\pi\big|\pi\in NC(k,k)\right\}$, with $k\in\mathbb N$.
\end{enumerate}
\end{proposition}

\begin{proof}
This follows from the various linear independence results from Proposition 3.9, by using the Gram matrix formula from Proposition 3.11, along with the well-known fact that vectors are linearly independent when their Gram matrix is invertible.
\end{proof}

Good news, we can now discuss $S_N^+$ with $N=n^2$, $n\geq2$, as follows:

\index{Catalan numbers}
\index{quantum permutation}
\index{Clebsch-Gordan rules}
\index{fusion rules}
\index{Marchenko-Pastur law}

\begin{theorem}
The quantum groups $S_N^+$ with $N=n^2$, $n\geq2$ are as follows:
\begin{enumerate}
\item The moments of the main character are the Catalan numbers:
$$\int_{S_N^+}\chi^k=C_k$$

\item The main character follows the Marchenko-Pastur law, $\chi\sim\pi_1$.

\item The fusion rules for irreducible representations are as for $SO_3$, namely:
$$r_k\otimes r_l=r_{|k-l|}+r_{|k-l|+1}+\ldots+r_{k+l}$$ 
 
\item We have $\dim r_k=(q^{k+1}-q^{-k})/(q-1)$, with $q+q^{-1}=N-2$.
\end{enumerate}
\end{theorem}

\begin{proof}
This is quite similar to the proof of Theorem 3.2, by using the linear independence result from Proposition 3.12 as main ingredient, as follows:

\medskip

(1) We have the following computation, using Peter-Weyl, then the easiness property of $S_N^+$, then Proposition 3.12 (2), then Proposition 3.10, and the definition of $C_k$:
$$\int_{S_N^+}\chi^k
=|NC(k)|
=|NC_2(2k)|
=C_k$$

(2) This is a reformulation of (1), using standard free probability theory.

\medskip

(3) This is identical to the proof of Theorem 3.6 (3), based on (1). 

\medskip

(4) Finally, the dimension formula is clear by recurrence.
\end{proof}

All this is very nice, and although there is still some work, in order to reach to results for $S_N^+$ at any $N\geq4$, let us just enjoy what we have. As a consequence, we have:

\begin{theorem}
The free quantum groups are as follows:
\begin{enumerate}
\item $U_N^+$ is not coamenable at $N\geq2$.

\item $O_N^+$ is coamenable at $N=2$, and  not coamenable at $N\geq3$.

\item $S_N^+$ is coamenable at $N\leq4$, and not coamenable at $N\geq5$.
\end{enumerate}
\end{theorem}

\begin{proof}
The various non-coamenability assertions are all clear, due to various examples of non-coamenable group dual subgroups $\widehat{\Gamma}\subset G$, coming from the theory in chapter 1. As for the amenability assertions, regarding $O_2^+$ and $S_4^+$, these come from Theorem 3.2 and Theorem 3.13, which show that the support of the spectral measure of $\chi$ is:
$$supp(\gamma_1)=[-2,2]\quad,\quad 
supp(\pi_1)=[0,4]$$

Thus the Kesten criterion from chapter 1, telling us that $G\subset O_N^+$ is coamenable precisely when $N\in supp(law(\chi))$, applies in both cases, and gives the result. 
\end{proof}

\section*{3c. Meander determinants}

Let us discuss now the extension of Theorem 3.13, to all the quantum groups $S_N^+$ with $N\geq4$. For this purpose we need an extension of the linear independence results from Proposition 3.12. This is something non-trivial, and the first thought goes to:

\begin{speculation}
There should be a theory of deformed compact quantum groups, alowing us to talk about $O_n^+$ with $n\in[2,\infty)$, having the same fusion rules as $SU_2$, and therefore solving via partition shrinking our $S_N^+$ problems at any $N\geq4$.
\end{speculation}

This speculation is legit, and in what concerns the first part, generalities, that theory is indeed available, from the Woronowicz papers \cite{wo1}, \cite{wo2}. Is it also possible to talk about deformations of $O_N^+$ in this setting, as explained in Wang's paper \cite{wa1}, with the new parameter $n\in[2,\infty)$ being of course not the dimension of the fundamental representation, but rather its ``quantum dimension''. And with this understood, all the rest is quite standard, and worked out in the quantum group literature. We refer to \cite{ba4} for more about this, but we will not follow this path, which is too complicated.

\bigskip

As a second speculation now, which is something complicated too, but is far more conceptual, we have the idea, already mentioned before, of getting what we want via the trace on the Temperley-Lieb algebra $TL_N(k)=span(NC_2(k,k))$. We will not follow this path either, which is quite complicated too, but here is how this method works:

\index{Temperley-Lieb algebra}

\begin{theorem}
Consider the Temperley-Lieb algebra of index $N\geq4$, defined as 
$$TL_N(k)=span(NC_2(k,k))$$
with product given by the rule $\bigcirc=N$, when concatenating. 
\begin{enumerate}
\item We have a representation $i:TL_N(k)\to B((\mathbb C^N)^{\otimes k})$, given by $\pi\to T_\pi$.

\item $Tr(T_\pi)=N^{loops(<\pi>)}$, where $\pi\to<\pi>$ is the closing operation.

\item The linear form $\tau=Tr\circ i:TL_N(k)\to\mathbb C$ is a faithful positive trace.

\item The representation $i:TL_N(k)\to B((\mathbb C^N)^{\otimes k})$ is faithful.
\end{enumerate}
In particular, the vectors $\left\{\xi_\pi|\pi\in NC(k)\right\}\subset(\mathbb C^N)^{\otimes k}$ are linearly independent.
\end{theorem}

\begin{proof}
All this is quite standard, but advanced, the idea being as follows:

\medskip

(1) This is clear from the categorical properties of $\pi\to T_\pi$.

\medskip

(2) This follows indeed from the following computation:
\begin{eqnarray*}
Tr(T_\pi)
&=&\sum_{i_1\ldots i_k}\delta_\pi\binom{i_1\ldots i_k}{i_1\ldots i_k}\\
&=&\#\left\{i_1,\ldots,i_k\in\{1,\ldots,N\}\Big|\ker\binom{i_1\ldots i_k}{i_1\ldots i_k}\geq\pi\right\}\\
&=&N^{loops(<\pi>)}
\end{eqnarray*}

(3) The traciality of $\tau$ is clear from definitions. Regarding now the faithfulness, this is something well-known, and we refer here to Jones' paper \cite{jo1}.

\medskip

(4) This follows from (3) above, via a standard positivity argument. As for the last assertion, this follows from (4), by fattening the partitions.
\end{proof}

We will be back to this later, when talking subfactors and planar algebras, with a closer look into Jones' paper \cite{jo1}. In the meantime, however, Speculation 3.15 and Theorem 3.16 will not do, being too advanced, so we have to come up with something else, more pedestrian. And this can only be the computation of the Gram determinant.

\medskip

We already know, from chapter 2, that for the group $S_N$ the formula of the corresponding Gram matrix determinant, due to Lindst\"om \cite{lin}, is as follows:

\index{Gram determinant}
\index{Lindst\"om formula}

\begin{theorem}
The determinant of the Gram matrix of $S_N$ is given by
$$\det(G_{kN})=\prod_{\pi\in P(k)}\frac{N!}{(N-|\pi|)!}$$
with the convention that in the case $N<k$ we obtain $0$.
\end{theorem}

\begin{proof}
This is something that we know from chapter 2, the idea being that $G_{kN}$ decomposes as a product of an upper triangular and lower triangular matrix.
\end{proof}

Although we will not need this here, let us discuss as well, for the sake of completness, the case of the orthogonal group $O_N$. Here the combinatorics is that of the Young diagrams. We denote by $|.|$ the number of boxes, and we use quantity $f^\lambda$, which gives the number of standard Young tableaux of shape $\lambda$. The result is then as follows:

\index{Young tableaux}

\begin{theorem}
The determinant of the Gram matrix of $O_N$ is given by
$$\det(G_{kN})=\prod_{|\lambda|=k/2}f_N(\lambda)^{f^{2\lambda}}$$
where the quantities on the right are $f_N(\lambda)=\prod_{(i,j)\in\lambda}(N+2j-i-1)$.
\end{theorem}

\begin{proof}
This follows from the results of Zinn-Justin in \cite{zin}. Indeed, it is known from there that the Gram matrix is diagonalizable, as follows:
$$G_{kN}=\sum_{|\lambda|=k/2}f_N(\lambda)P_{2\lambda}$$

Here $1=\sum P_{2\lambda}$ is the standard partition of unity associated to the Young diagrams having $k/2$ boxes, and the coefficients $f_N(\lambda)$ are those in the statement. Now since we have $Tr(P_{2\lambda})=f^{2\lambda}$, this gives the result. See \cite{bcu}, \cite{zin}.
\end{proof}

For the free orthogonal and symmetric groups, the results, by Di Francesco \cite{dif}, are substantially more complicated. Let us begin with some examples. We first have:

\begin{proposition}
The first Gram matrices and determinants for $O_N^+$ are
$$\det\begin{pmatrix}N^2&N\\N&N^2\end{pmatrix}=N^2(N^2-1)$$
$$\det\begin{pmatrix}
N^3&N^2&N^2&N^2&N\\
N^2&N^3&N&N&N^2\\
N^2&N&N^3&N&N^2\\
N^2&N&N&N^3&N^2\\
N&N^2&N^2&N^2&N^3
\end{pmatrix}=N^5(N^2-1)^4(N^2-2)$$
with the matrices being written by using the lexicographic order on $NC_2(2k)$.
\end{proposition}

\begin{proof}
The formula at $k=2$, where $NC_2(4)=\{\sqcap\sqcap,\bigcap\hskip-4.9mm{\ }_\cap\,\}$, is clear. At $k=3$ however, things are tricky. We have $NC(3)=\{|||,\sqcap|,\sqcap\hskip-3.2mm{\ }_|\,,|\sqcap,\sqcap\hskip-0.7mm\sqcap\}$, and the corresponding Gram matrix and its determinant are, according to Theorem 3.17:
$$\det\begin{pmatrix}
N^3&N^2&N^2&N^2&N\\
N^2&N^2&N&N&N\\
N^2&N&N^2&N&N\\
N^2&N&N&N^2&N\\
N&N&N&N&N
\end{pmatrix}=N^5(N-1)^4(N-2)$$

By using Proposition 3.11, the Gram determinant of $NC_2(6)$ is given by:
\begin{eqnarray*}
\det(G_{6N})
&=&\frac{1}{N^2\sqrt{N}}\times N^{10}(N^2-1)^4(N^2-2)\times\frac{1}{N^2\sqrt{N}}\\
&=&N^5(N^2-1)^4(N^2-2)
\end{eqnarray*}

Thus, we have obtained the formula in the statement.
\end{proof}

In general, such tricks won't work, because $NC(k)$ is strictly smaller than $P(k)$ at $k\geq4$. However, following Di Francesco \cite{dif}, we have the following result:

\index{meander determinant}
\index{Gram determinant}

\begin{theorem}
The determinant of the Gram matrix for $O_N^+$ is given by
$$\det(G_{kN})=\prod_{r=1}^{[k/2]}P_r(N)^{d_{k/2,r}}$$
where $P_r$ are the Chebycheff polynomials, given by
$$P_0=1\quad,\quad 
P_1=X\quad,\quad 
P_{r+1}=XP_r-P_{r-1}$$
and $d_{kr}=f_{kr}-f_{k,r+1}$, with $f_{kr}$ being the following numbers, depending on $k,r\in\mathbb Z$,
$$f_{kr}=\binom{2k}{k-r}-\binom{2k}{k-r-1}$$
with the convention $f_{kr}=0$ for $k\notin\mathbb Z$. 
\end{theorem}

\begin{proof}
This is something quite technical, obtained by using a decomposition as follows of the Gram matrix $G_{kN}$, with the matrix $T_{kN}$ being lower triangular:
$$G_{kN}=T_{kN}T_{kN}^t$$

Thus, a bit as in the proof of the Lindst\"om formula, we obtain the result, but the problem lies however in the construction of $T_{kN}$, which is non-trivial. See \cite{dif}.
\end{proof}

We refer to \cite{bcu} for further details regarding the above result, including a short proof, based on the bipartite planar algebra combinatorics developed by Jones in \cite{jo4}. Let us also mention that the Chebycheff polynomials have something to do with all this due to the fact that these are the orthogonal polynomials for the Wigner law. See \cite{bcu}.

\bigskip

Moving ahead now, regarding $S_N^+$, we have here the following formula, which is quite similar, obtained via shrinking, also from Di Francesco \cite{dif}:

\index{meander determinant}
\index{Gram determinant}

\begin{theorem}
The determinant of the Gram matrix for $S_N^+$ is given by
$$\det(G_{kN})=(\sqrt{N})^{a_k}\prod_{r=1}^kP_r(\sqrt{N})^{d_{kr}}$$
where $P_r$ are the Chebycheff polynomials, given by
$$P_0=1\quad,\quad 
P_1=X\quad,\quad 
P_{r+1}=XP_r-P_{r-1}$$
and $d_{kr}=f_{kr}-f_{k,r+1}$, with $f_{kr}$ being the following numbers, depending on $k,r\in\mathbb Z$,
$$f_{kr}=\binom{2k}{k-r}-\binom{2k}{k-r-1}$$
with the convention $f_{kr}=0$ for $k\notin\mathbb Z$, and where $a_k=\sum_{\pi\in \mathcal P(k)}(2|\pi|-k)$.
\end{theorem}

\begin{proof}
This follows indeed from Theorem 3.20, by using Proposition 3.11.
\end{proof}

Getting back now to our quantum permutation group questions, by using the above results we can produce a key technical ingredient, as follows:

\begin{proposition}
The following are linearly independent, for any $N\geq4$:
\begin{enumerate}
\item The linear maps $\left\{T_\pi\big|\pi\in NC(k,l)\right\}$, with $k,l\in 2\mathbb N$.

\item The vectors $\left\{\xi_\pi\big|\pi\in NC(k)\right\}$, with $k\in\mathbb N$.

\item The linear maps $\left\{T_\pi\big|\pi\in NC(k,k)\right\}$, with $k\in\mathbb N$.
\end{enumerate}
\end{proposition}

\begin{proof}
The statement is identical to Proposition 3.12, with the assumption $N=n^2$ lifted. As for the proof, this comes from the formula in Theorem 3.21.
\end{proof}

With this in hand, we have the following extension of Theorem 3.13:

\index{quantum permutation}
\index{Catalan numbers}
\index{Marchenko-Pastur law}
\index{Clebsch-Gordan rules}
\index{main character}

\begin{theorem}
The quantum groups $S_N^+$ with $N\geq4$ are as follows:
\begin{enumerate}
\item The moments of the main character are the Catalan numbers:
$$\int_{S_N^+}\chi^k=C_k$$

\item The main character follows the Marchenko-Pastur law, $\chi\sim\pi_1$.

\item The fusion rules for irreducible representations are as for $SO_3$, namely:
$$r_k\otimes r_l=r_{|k-l|}+r_{|k-l|+1}+\ldots+r_{k+l}$$ 
 
\item We have $\dim r_k=(q^{k+1}-q^{-k})/(q-1)$, with $q+q^{-1}=N-2$.
\end{enumerate}
\end{theorem}

\begin{proof}
This is identical to the proof of Theorem 3.13, by using this time the linear independence result from Proposition 3.22 as technical ingredient.
\end{proof}

So long for representations of $S_N^+$. All the above might seem quite complicated, but we repeat, up to some standard algebra, everything comes down to Proposition 3.22. And with some solid modern mathematical knowledge, be that operator algebras a la Jones, or deformed quantum groups a la Woronowicz, or meander determinants a la Di Francesco, the result there is in fact trivial. You can check here \cite{ba1}, \cite{bc1}, both short papers.

\bigskip

In what concerns us, we will be back to the similarity between $S_N^+$ and $SO_3$ on several occasions, with a number of further results on the subject, refining Theorem 3.23.

\section*{3d. Planar algebras}

In the remainder of this chapter we keep developing some useful theory for $U_N^+,O_N^+,S_N^+$. We will present among others a result from \cite{ba3}, refining the Tannakian duality for the quantum permutation groups $G\subset S_N^+$, stating that these quantum groups are in correspondence with the subalgebras of Jones' spin planar algebra $P\subset\mathcal S_N$.

\bigskip

In order to get started, we need a lot of preliminaries, the lineup being von Neumann algebras, ${\rm II}_1$ factors, subfactors, and finally planar algebras. We already met von Neumann algebras, in chapter 1. The advanced general theory regarding them is as follows:

\index{von Neumann algebra}
\index{factor}
\index{continuous dimension}
\index{Murray-von Neumann factor R}
\index{Connes theorem}
\index{hyperfinite factor}

\begin{theorem}
The von Neumann algebras $A\subset B(H)$ are as follows:
\begin{enumerate}
\item Any such algebra decomposes as $A=\int_XA_xdx$, with $X$ being the spectrum of the center, $Z(A)=L^\infty(X)$, and with the fibers $A_x$ being factors, $Z(A_x)=\mathbb C$.

\item The factors can be fully classified in terms of ${\rm II}_1$ factors, which are those factors satisfying $\dim A=\infty$, and having a faithful trace $tr:A\to\mathbb C$.

\item The ${\rm II}_1$ factors enjoy the ``continuous dimension geometry'' property, in the sense that the traces of their projections can take any values in $[0,1]$.

\item Among the ${\rm II}_1$ factors, the smallest one is the Murray-von Neumann hyperfinite factor $R$, obtained as an inductive limit of matrix algebras.
\end{enumerate}
\end{theorem}

\begin{proof}
This is something heavy, the idea being as follows:

\medskip

(1) This is von Neumann's reduction theory theorem, which follows in finite dimensions from $A=M_{n_1}(\mathbb C)\oplus\ldots\oplus M_{n_k}(\mathbb C)$, and whose proof in general is quite technical.

\medskip

(2) This comes from results of Murray-von Neumann and Connes, the idea being that the other factors can be basically obtained via crossed product constructions.

\medskip

(3) This is subtle functional analysis, with the rational traces being relatively easy to obtain, and with the irrational ones coming from limiting arguments.

\medskip

(4) Once again, heavy results, by Murray-von Neumann and Connes, the idea being that any finite dimensional construction always leads to the same factor, called $R$.
\end{proof}

Let us discuss now subfactor theory, following Jones' fundamental paper \cite{jo1}. Jones looked at the inclusions of ${\rm II}_1$ factors $A\subset B$, called subfactors, which are quite natural objects in physics. Given such an inclusion, we can talk about its index:

\index{subfactor}
\index{index of subfactor}

\begin{definition}
The index of an inclusion of ${\rm II}_1$ factors $A\subset B$ is the quantity
$$[B:A]=\dim_AB\in[1,\infty]$$
constructed by using the Murray-von Neumann continuous dimension theory.
\end{definition}

In order to explain Jones' result in \cite{jo1}, it is better to relabel our subfactor as $A_0\subset A_1$. We can construct the orthogonal projection $e_1:A_1\to A_0$, and set: 
$$A_2=<A_1,e_1>$$

This remarkable procedure, called ``basic construction'', can be iterated, and we obtain in this way a whole tower of ${\rm II}_1$ factors, as follows:
$$A_0\subset_{e_1}A_1\subset_{e_2}A_2\subset_{e_3}A_3\subset\ldots\ldots$$

Quite surprisingly, this construction leads to a link with the Temperley-Lieb algebra $TL_N=span(NC_2)$. The results can be summarized as follows:

\index{basic construction}
\index{Jones projection}
\index{Jones index theorem}
\index{Temperley-Lieb algebra}

\begin{theorem}
Let $A_0\subset A_1$ be an inclusion of ${\rm II}_1$ factors.
\begin{enumerate}
\item The sequence of projections $e_1,e_2,e_3,\ldots\in B(H)$ produces a representation of the Temperley-Lieb algebra of index $N=[A_1,A_0]$, as follows:
$$TL_N\subset B(H)$$

\item The index $N=[A_1,A_0]$, which is a Murray-von Neumann continuous quantity $N\in[1,\infty]$, must satisfy the following condition:
$$N\in\left\{4\cos^2\left(\frac{\pi}{n}\right)\Big|n\in\mathbb N\right\}\cup[4,\infty]$$
\end{enumerate}
\end{theorem}

\begin{proof}
This result, from \cite{jo1}, is something tricky, the idea being as follows:

\medskip

(1) The idea here is that the functional analytic study of the basic construction leads to the conclusion that the sequence of projections $e_1,e_2,e_3,\ldots\in B(H)$ behaves algebrically, when rescaled, exactly as the sequence of diagrams $\varepsilon_1,\varepsilon_2,\varepsilon_3,\ldots\in TL_N$ given by:
$$\varepsilon_1={\ }^\cup_\cap\quad,\quad 
\varepsilon_2=|\!{\ }^\cup_\cap\quad,\quad 
\varepsilon_3=||\!{\ }^\cup_\cap\quad,\quad 
\ldots$$

But these diagrams generate $TL_N$, and so we have an embedding $TL_N\subset B(H)$, where $H$ is the Hilbert space where our subfactor $A_0\subset A_1$ lives, as claimed.

\medskip

(2) This is something quite surprising, which follows from (1), via some clever positivity considerations, involving the Perron-Frobenius theorem. In fact, the subfactors having index $N\in[1,4]$ can be classified by ADE diagrams, and the obstruction $N=4\cos^2(\frac{\pi}{n})$ comes from the fact that $N$ must be the squared norm of such a graph.
\end{proof}

Quite remarkably, Theorem 3.26 is just the tip of the iceberg. One can prove indeed that the planar algebra structure of $TL_N$, taken in an intuitive sense, extends to a planar algebra structure on the sequence of relative commutants $P_k=A_0'\cap A_k$. In order to discuss this key result, due as well to Jones, from \cite{jo3}, and that we will need too, in connection with our quantum group problems, let us start with:

\index{higher commutant}
\index{planar algebra}

\begin{definition}
The planar algebras are defined as follows:
\begin{enumerate}
\item A $k$-tangle, or $k$-box, is a rectangle in the plane, with $2k$ marked points on its boundary, containing $r$ small boxes, each having $2k_i$ marked points, and with the $2k+\sum 2k_i$ marked points being connected by noncrossing strings.

\item A planar algebra is a sequence of finite dimensional vector spaces $P=(P_k)$, together with linear maps $P_{k_1}\otimes\ldots\otimes P_{k_r}\to P_k$, one for each $k$-box, such that the gluing of boxes corresponds to the composition of linear maps.
\end{enumerate}
\end{definition}

As basic example of a planar algebra, we have the Temperley-Lieb algebra $TL_N$. Indeed, putting $TL_N(k_i)$ diagrams into the small $r$ boxes of a $k$-box clearly produces a $TL_N(k)$ diagram, so we have indeed a planar algebra, of somewhat ``trivial'' type. 

\bigskip

In general, the planar algebras are more complicated than this, and we will be back later with some explicit examples. However, the idea is very simple, namely that ``the elements of a planar algebra are not necessarily diagrams, but they behave like diagrams". In relation now with subfactors, the result, which extends Theorem 3.26 (1), and which was found by Jones in \cite{jo3}, almost 20 years after \cite{jo1}, is as follows:

\index{higher commutant}
\index{planar algebra}

\begin{theorem} 
Given a subfactor $A_0\subset A_1$, the collection $P=(P_k)$ of linear spaces 
$$P_k=A_0'\cap A_k$$
has a planar algebra structure, extending the planar algebra structure of $TL_N$.
\end{theorem}

\begin{proof}
As a first observation, since $e_1:A_1\to A_0$ commutes with $A_0$ we have $e_1\in P_2'$. By translation we obtain $e_1,\ldots,e_{k-1}\in P_k$ for any $k$, and so:
$$TL_N\subset P$$

The point now is that the planar algebra structure of $TL_N$, obtained by composing diagrams, can be shown to extend into an abstract planar algebra structure of $P$. This is something quite technical, and we will not get into details here. See \cite{jo3}.
\end{proof}

Getting back to quantum groups, all this machinery is interesting for us. We will need the construction of the tensor and spin planar algebras $\mathcal T_N,\mathcal S_N$. Let us start with:

\index{tensor planar algebra}
\index{Kronecker symbol}

\begin{definition}
The tensor planar algebra $\mathcal T_N$ is the sequence of vector spaces 
$$P_k=M_N(\mathbb C)^{\otimes k}$$
with the multilinear maps $T_\pi:P_{k_1}\otimes\ldots\otimes P_{k_r}\to P_k$
being given by the formula
$$T_\pi(e_{i_1}\otimes\ldots\otimes e_{i_r})=\sum_j\delta_\pi(i_1,\ldots,i_r:j)e_j$$
with the Kronecker symbols $\delta_\pi$ being $1$ if the indices fit, and being $0$ otherwise.
\end{definition}

In other words, we are using here a construction which is very similar to the construction $\pi\to T_\pi$ that we used for easy quantum groups. We put the indices of the basic tensors on the marked points of the small boxes, in the obvious way, and the coefficients of the output tensor are then given by Kronecker symbols, exactly as in the easy case.

\bigskip

The fact that we have indeed a planar algebra, in the sense that the gluing of tangles corresponds to the composition of linear maps, as required by Definition 3.27, is something elementary, in the same spirit as the verification of the functoriality properties of the correspondence $\pi\to T_\pi$, discussed in chapter 2, and we refer here to Jones \cite{jo3}. 

\bigskip

Let us discuss now a second planar algebra of the same type, which is important as well for various reasons, namely the spin planar algebra $\mathcal S_N$. This planar algebra appears somehow as the ``square root'' of the tensor planar algebra $\mathcal T_N$. Let us start with:

\begin{definition}
We write the standard basis of $(\mathbb C^N)^{\otimes k}$ in $2\times k$ matrix form,
$$e_{i_1\ldots i_k}=
\begin{pmatrix}i_1 & i_1 &i_2&i_2&i_3&\ldots&\ldots\\
i_k&i_k&i_{k-1}&\ldots&\ldots&\ldots&\ldots 
\end{pmatrix}$$
by duplicating the indices, and then writing them clockwise, starting from top left.
\end{definition}

Now with this convention in hand for the tensors, we can formulate the construction of the spin planar algebra $\mathcal S_N$, also from \cite{jo3}, as follows:

\index{spin planar algebra}
\index{Kronecker symbol}

\begin{definition}
The spin planar algebra $\mathcal S_N$ is the sequence of vector spaces 
$$P_k=(\mathbb C^N)^{\otimes k}$$
written as above, with the multiplinear maps $T_\pi:P_{k_1}\otimes\ldots\otimes P_{k_r}\to P_k$ 
being given by
$$T_\pi(e_{i_1}\otimes\ldots\otimes e_{i_r})=\sum_j\delta_\pi(i_1,\ldots,i_r:j)e_j$$
with the Kronecker symbols $\delta_\pi$ being $1$ if the indices fit, and being $0$ otherwise.
\end{definition}

Here are some illustrating examples for the spin planar algebra calculus:

\medskip

(1) The identity $1_k$ is the $(k,k)$-tangle having vertical strings only. The solutions of $\delta_{1_k}(x,y)=1$ being the pairs of the form $(x,x)$, this tangle $1_k$ acts by the identity:
$$1_k\begin{pmatrix}j_1 & \ldots & j_k\\ i_1 & \ldots & i_k\end{pmatrix}=\begin{pmatrix}j_1 & \ldots & j_k\\ i_1 & \ldots & i_k\end{pmatrix}$$

(2) The multiplication $M_k$ is the $(k,k,k)$-tangle having 2 input boxes, one on top of the other, and vertical strings only. It acts in the following way:
$$M_k\left( 
\begin{pmatrix}j_1 & \ldots & j_k\\ i_1 & \ldots & i_k\end{pmatrix}
\otimes\begin{pmatrix}l_1 & \ldots & l_k\\ m_1 & \ldots & m_k\end{pmatrix}
\right)=
\delta_{j_1m_1}\ldots \delta_{j_km_k}
\begin{pmatrix}l_1 & \ldots & l_k\\ i_1 & \ldots & i_k\end{pmatrix}$$

(3) The inclusion $I_k$ is the $(k,k+1)$-tangle which looks like $1_k$, but has one more vertical string, at right of the input box. Given $x$, the solutions of $\delta_{I_k}(x,y)=1$ are the elements $y$ obtained from $x$ by adding to the right a vector of the form $(^l_l)$, and so:
$${I_k}\begin{pmatrix}j_1 & \ldots & j_k\\ i_1 & \ldots & i_k\end{pmatrix}=
\sum_l\begin{pmatrix}j_1 & \ldots & j_k& l\\ i_1 & \ldots & i_k& l\end{pmatrix}$$

(4) The expectation $U_k$ is the $(k+1,k)$-tangle which looks like $1_k$, but has one more string, connecting the extra 2 input points, both at right of the input box:
$$U_k
\begin{pmatrix}j_1 & \ldots &j_k& j_{k+1}\\ i_1 & \ldots &i_k& i_{k+1}\end{pmatrix}=
\delta_{i_{k+1}j_{k+1}}
\begin{pmatrix}j_1 & \ldots & j_k\\ i_1 & \ldots & i_k\end{pmatrix}$$

(5) The Jones projection $E_k$ is a $(0,k+2)$-tangle, having no input box. There are $k$ vertical strings joining the first $k$ upper points to the first $k$ lower points, counting from left to right. The remaining upper 2 points are connected by a semicircle, and the remaining lower 2 points are also connected by a semicircle. We have:
$$E_k(1)=\sum_{ijl}\begin{pmatrix}i_1 & \ldots &i_k&j&j\\ i_1 & \ldots &i_k&l&l\end{pmatrix}$$

The elements $e_k=N^{-1}E_k(1)$ are then projections, and define a representation of the infinite Temperley-Lieb algebra of index $N$ inside the inductive limit algebra $\mathcal S_N$.

\medskip

(6) The rotation $R_k$ is the $(k,k)$-tangle which looks like $1_k$, but the first 2 input points are connected to the last 2 output points, and the same happens at right:
$$R_k=\begin{matrix}
\hskip 0.3mm\Cap \ |\ |\ |\ |\hskip -0.5mm |\cr
|\hskip -0.5mm |\hskip 10.3mm |\hskip -0.5mm |\cr
\hskip -0.3mm|\hskip -0.5mm |\ |\ |\ |\ \hskip -0.1mm\Cup
\end{matrix}$$

The action of $R_k$ on the standard basis is by rotation of the indices, as follows:
$$R_k(e_{i_1i_2\ldots i_k})=e_{i_2\ldots i_ki_1}$$

There are many other interesting examples of $k$-tangles, but in view of our present purposes, we can actually stop here, due to the following fact:

\index{rotation tangle}

\begin{theorem}
The multiplications, inclusions, expectations, Jones projections and rotations generate the set of all tangles, via the gluing operation.
\end{theorem}

\begin{proof}
This is something well-known and elementary, obtained by ``chopping'' the various planar tangles into small pieces, as in the above list. See \cite{jo3}.
\end{proof}

Finally, in order for our discussion to be complete, we must talk as well about the $*$-structure of the spin planar algebra. Once again this is constructed as in the easy quantum group calculus, by turning upside-down the diagrams, as follows:
$$\begin{pmatrix}j_1 & \ldots & j_k\\ i_1 & \ldots & i_k\end{pmatrix}^*
=\begin{pmatrix}i_1 & \ldots & i_k\\ j_1 & \ldots & j_k\end{pmatrix}$$

Getting back now to quantum groups, following \cite{ba3}, we have the following result:

\index{tensor power coaction}
\index{spin planar algebra}

\begin{theorem}
Given $G\subset S_N^+$, consider the tensor powers of the associated coaction map on $C(X)$, where $X=\{1,\ldots,N\}$, which are the folowing linear maps:
$$\Phi^k:C(X^k)\to C(X^k)\otimes C(G)$$
$$e_{i_1\ldots i_k}\to\sum_{j_1\ldots j_k}e_{j_1\ldots j_k}\otimes u_{j_1i_1}\ldots u_{j_ki_k}$$
The fixed point spaces of these coactions, which are by definition the spaces
$$P_k=\left\{ x\in C(X^k)\Big|\Phi^k(x)=1\otimes x\right\}$$
are given by $P_k=Fix(u^{\otimes k})$, and form a subalgebra of the spin planar algebra $\mathcal S_N$.
\end{theorem}

\begin{proof}
Since the map $\Phi$ is a coaction, its tensor powers $\Phi^k$ are coactions too, and at the level of fixed point algebras we have the following formula:
$$P_k=Fix(u^{\otimes k})$$

In order to prove now the planar algebra assertion, we will use Theorem 3.32. Consider the rotation $R_k$. Rotating, then applying $\Phi^k$, and rotating backwards by $R_k^{-1}$ is the same as applying $\Phi^k$, then rotating each $k$-fold product of coefficients of $\Phi$. Thus the elements obtained by rotating, then applying $\Phi^k$, or by applying $\Phi^k$, then rotating, differ by a sum of Dirac masses tensored with commutators in $A=C(G)$:
$$\Phi^kR_k(x)-(R_k\otimes id)\Phi^k(x)\in C(X^k)\otimes [A,A]$$

Now let $\int_A$ be the Haar functional of $A$, and consider the conditional expectation onto the fixed point algebra $P_k$, which is given by the following formula:
$$\phi_k=\left(id\otimes\int_A\right)\Phi^k$$

Since $\int_A$ is a trace, it vanishes on commutators. Thus $R_k$ commutes with $\phi_k$:
$$\phi_kR_k=R_k\phi_k$$

The commutation relation $\phi_kT=T\phi_l$ holds in fact for any $(l,k)$-tangle $T$. These tangles are called annular, and the proof is by verification on generators of the annular category. In particular we obtain, for any annular tangle $T$:
$$\phi_kT\phi_l=T\phi_l$$

We conclude from this that the annular category is contained in the suboperad $\mathcal P'\subset\mathcal P$ of the planar operad consisting of tangles $T$ satisfying the following condition, where $\phi =(\phi_k)$, and where $i(.)$ is the number of input boxes:
$$\phi T\phi^{\otimes i(T)}=T\phi^{\otimes i(T)}$$

On the other hand the multiplicativity of $\Phi^k$ gives $M_k\in\mathcal P'$. Now since the planar operad $\mathcal P$ is generated by multiplications and annular tangles, it follows that we have $\mathcal P'=P$. Thus for any tangle $T$ the corresponding multilinear map between spaces $P_k(X)$ restricts to a multilinear map between spaces $P_k$. In other words, the action of the planar operad $\mathcal P$ restricts to $P$, and makes it a subalgebra of $\mathcal S_N$, as claimed.
\end{proof}

As a second result now, also from \cite{ba3}, completing our study, we have:

\index{planar algebra}
\index{spin planar algebra}
\index{Tannakian duality}
\index{annular category}

\begin{theorem}
We have a bijection between quantum permutation groups and subalgebras of the spin planar algebra,
$$(G\subset S_N^+)\quad\longleftrightarrow\quad (Q\subset\mathcal S_N)$$
given in one sense by the construction in Theorem 3.33, and in the other sense by a suitable modification of Tannakian duality.
\end{theorem}

\begin{proof}
The idea is that this will follow by applying Tannakian duality to the annular category over $Q$. Let $n,m$ be positive integers. To any element $T_{n+m}\in Q_{n+m}$ we associate a linear map $L_{nm}(T_{n+m}):P_n(X)\to P_m(X)$ in the following way:
$$L_{nm}\left(\begin{matrix}|\ |\ |\\ T_{n+m}\\ |\ |\ |\end{matrix}\right):
\left(\begin{matrix}|\\ a_n\\ |\end{matrix}\right)
\to \left(\begin{matrix}
\hskip 1.5mm |\hskip 3.0mm |\hskip 3.0mm \cap\\
\ \ T_{n+m}\hskip 0.0mm  |\\
\hskip 1.9mm |\hskip 1.2mm |\hskip 3.2mm |\hskip2.2mm |\\
a_n|\hskip 3.2mm |\hskip 2.2mm |\\
\hskip 2.1mm\cup \hskip3.5mm |\hskip 2.2mm |
\end{matrix}\right)$$

That is, we consider the planar $(n,n+m,m)$-tangle having an small input $n$-box, a big input $n+m$-box and an output $m$-box, with strings as on the picture of the right. This defines a certain multilinear map, as follows:
$$P_n(X)\otimes P_{n+m}(X)\to P_m(X)$$

If we put the element $T_{n+m}$ in the big input box, we obtain in this way a certain linear map $P_n(X)\to P_m(X)$, that we call $L_{nm}$. With this convention, let us set:
$$Q_{nm}=\left\{ L_{nm}(T_{n+m}):P_n(X)\to P_m(X)\Big| T_{n+m}\in Q_{n+m}\right\}$$

These spaces form a Tannakian category, so by \cite{wo2} we obtain a Woronowicz algebra $(A,u)$, such that the following equalities hold, for any $m,n$:
$$Hom(u^{\otimes m},u^{\otimes n})=Q_{mn}$$

We prove now that $u$ is a magic unitary. We have $Hom(1,u^{\otimes 2})=Q_{02}=Q_2$, so the unit of $Q_2$ must be a fixed vector of $u^{\otimes 2}$. But $u^{\otimes 2}$ acts on the unit of $Q_2$ as follows:
\begin{eqnarray*}
u^{\otimes 2}(1)
&=&u^{\otimes 2}\left( \sum_i \begin{pmatrix}i&i\\ i&i\end{pmatrix}\right)\\
&=&\sum_{ikl}\begin{pmatrix}k&k\\ l&l\end{pmatrix}\otimes u_{ki}u_{li}\\
&=&\sum_{kl}\begin{pmatrix}k&k\\ l&l\end{pmatrix}\otimes (uu^t)_{kl}
\end{eqnarray*}

From $u^{\otimes 2}(1)=1\otimes 1$ ve get that $uu^t$ is the identity matrix. Together with the unitarity of $u$, this gives the following formulae:
$$u^t=u^*=u^{-1}$$

Consider the Jones projection $E_1\in Q_3$. After isotoping, $L_{21}(E_1)$ looks as follows:
$$L_{21}\left( \Bigl| \begin{matrix}\cup\\\cap\end{matrix}\right) :
\begin{pmatrix} \,|\ |\\ {\ }^i_j{\ }^i_j\\ \,|\ |\end{pmatrix}\,\to\,
\begin{pmatrix}\hskip -5.8mm |\\ {\ }^i_j{\ }^i_j\supset\\ \hskip -5.8mm |\end{pmatrix}
=\,\delta_{ij}\begin{pmatrix}\,|\\ {\ }^i_i\\ \,|\end{pmatrix}$$

In other words, the linear map $M=L_{21}(E_1)$ is the multiplication $\delta_i\otimes\delta_j\to\delta_{ij}\delta_i$:
$$M\begin{pmatrix}i&i\\ j&j\end{pmatrix}
=\delta_{ij}\begin{pmatrix}i\\ i\end{pmatrix}$$

In order to finish, consider the following element of $C(X)\otimes A$:
$$(M\otimes id)u^{\otimes 2}\left(\begin{pmatrix}i&i\\ j&j\end{pmatrix}\otimes 1\right)
=\sum_k\begin{pmatrix}k\\ k\end{pmatrix}\delta_k\otimes u_{ki}u_{kj}$$

Since $M\in Q_{21}=Hom(u^{\otimes 2},u)$, this equals the following element of $C(X)\otimes A$:
$$u(M\otimes id)\left(\begin{pmatrix}i&i\\ j&j\end{pmatrix}\otimes 1\right)
=\sum_k\begin{pmatrix}k\\ k\end{pmatrix}\delta_k\otimes\delta_{ij}u_{ki}$$

Thus we have $u_{ki}u_{kj}=\delta_{ij}u_{ki}$ for any $i,j,k$, which shows that $u$ is a magic unitary. Now if $P$ is the planar algebra associated to $u$, we have $Hom(1,v^{\otimes n})=P_n=Q_n$, as desired. As for the uniqueness, this is clear from the Peter-Weyl theory.
\end{proof}

All the above might seem a bit technical, but is worth learning, and for good reason, because it is extremely powerful. As an example of immediate application, if you agree with the bijection $G\leftrightarrow Q$ in Theorem 3.34, then $G=S_N^+$ itself, which is the biggest object on the left, must correspond to the smallest object on the right, namely $Q=TL_N$. Thus, more or less everything that we learned so far in this book is trivial.

\bigskip

Welcome to planar algebras. Try to master this technology. And once this understood, get to know some analysis too, which comes after. But it will be among our main purposes here to do so, getting you familiar with algebra, and with some analysis as well.

\bigskip

Back now to work, the results established above, regarding the subgroups $G\subset S_N^+$, have several generalizations, to the subgroups $G\subset O_N^+$ and $G\subset U_N^+$, as well as subfactor versions, going beyond the combinatorial level. At the algebraic level, we have:

\begin{theorem}
The following happen:
\begin{enumerate}
\item The closed subgroups $G\subset O_N^+$ produce planar algebras $P\subset\mathcal T_N$, via the following formula, and any subalgebra $P\subset\mathcal T_N$ appears in this way:
$$P_k=End(u^{\otimes k})$$

\item The closed subgroups $G\subset U_N^+$ produce planar algebras $P\subset\mathcal T_N$, via the following formula, and any subalgebra $P\subset\mathcal T_N$ appears in this way:
$$P_k=End(\underbrace{u\otimes\bar{u}\otimes u\otimes\ldots}_{k\ terms})$$

\item In fact, the closed subgroups $G\subset PO_N^+\simeq PU_N^+$ are in correspondence with the subalgebras $P\subset\mathcal T_N$, with $G\to P$ being given by $P_k=Fix(u^{\otimes k})$.
\end{enumerate}
\end{theorem}

\begin{proof}
There is a long story with this result, whose origins go back to papers of mine written before the 1999 papers \cite{ba1}, \cite{jo3}, using Popa's standard lattice formalism, instead of the planar algebra one, and then to a number of papers written in the early 2000s, proving results which are more general. For the whole story, and a modern treatment of the subject, we refer to Tarrago-Wahl \cite{twa}. As in what regards the proof:

\medskip

(1) This is similar to the proof of Theorem 3.33 and Theorem 3.34, ultimately coming from Woronowicz's Tannakian duality in \cite{wo2}. Note however that the correspondence is not bijective, because the spaces $P_k$ determine $PG\subset PO_N^+$, but not $G\subset O_N^+$ itself.

\medskip

(2) This is an extension of (1), and the same comments apply. With the extra comment that the fact that the subgroups $PG\subset PO_N^+$ produce the same planar algebras as the subgroups $PG\subset PU_N^+$ should not be surprising, due to $PO_N^+=PU_N^+$.

\medskip

(3) This is an extension of (2), and a further extension of (1), and is in fact the best result on the subject, due to the fact that we have there a true, bijective correspondence. As before, this ultimately comes from Woronowicz's Tannakian duality in \cite{wo2}. 

\medskip

(4) As a final comment, you might say that, now that we have (3) as ultimate result on the subject, why not saying a few words about the proof. In answer, (3) is in fact just the tip of the iceberg, so we prefer to discuss this later, once we'll see the whole iceberg.
\end{proof}

Finally, in relation with subfactors, the result here is as follows:

\begin{theorem}
The planar algebras coming the subgroups $G\subset S_N^+$ appear from fixed point subfactors, of the following type,
$$A^G\subset(\mathbb C^N\otimes A)^G$$
and the planar algebras coming from the subgroups $G\subset PO_N^+=PU_N^+$ appear as well from fixed point subfactors, of the following type,
$$A^G\subset(M_N(\mathbb C)\otimes A)^G$$
with the action $G\curvearrowright A$ being assumed to be minimal, $(A^G)'\cap A=\mathbb C$.
\end{theorem}

\begin{proof}
Again, there is a long story with this result, and besides needing some explanations, regarding the proof, all this is in need of some unification. We will be back to this in chapter 4, and in the meantime we refer to \cite{ba1}, \cite{twa} and related papers.
\end{proof}

Finally, let us mention that an important question, which is still open, is that of understanding whether the above subfactors can be taken to be hyperfinite, $A^G\simeq R$. This is related to the axiomatization of hyperfinite subfactors, another open question, which is of central importance in von Neumann algebras. We will be back to this.

\section*{3e. Exercises} 

Things have been quite technical in this chapter, and our exercises will be quite technical too. But before that, in relation with $SU_2,SO_3$, we have:

\begin{exercise}
Learn more about $SU_2,SO_3$, namely Euler-Rodrigues, various proofs for the Clebsch-Gordan rules, and ADE/McKay correspondence for subgroups too.
\end{exercise}

All this, and we insist, is very important. The more you know about $SU_2,SO_3$, and with such things being a pleasure to learn, the better your mathematics and physics will be, no matter what mathematics and physics you are interested in.

\begin{exercise}
Prove that the quantum group $O_2^+$ appears as a twist of $SU_2$, 
$$O_2^+\simeq SU_2^{-1}$$
and deduce the Clebsch-Gordan fusion rules for $O_2^+$ by using this.
\end{exercise}

As a bonus exercise, you can try to understand as well the relation between $O_N^+$ and $SU_2$ at $N\geq3$, with the indication here that all this is related to two technical topics, namely the ``FRT deformation'' procedure, and the ``free symplectic groups''.

\begin{exercise}
Prove that the quantum group $S_4^+$ appears as a twist of $SO_3$, 
$$S_4^+\simeq SO_3^{-1}$$
and deduce the Clebsch-Gordan fusion rules for $S_4^+$ by using this.
\end{exercise}

As before with the previous exercise, all this requires a good knowledge of the cocycle twisting procedure, and basically the same comments as there apply.

\begin{exercise}
Establish, with full details, the linear independence of the vectors
$$\left\{\xi_\pi\Big|\pi\in NC(k)\right\}$$
using meander determinants, subfactor theory, and deformed quantum groups.
\end{exercise}

Here the first question, in relation with meander determinants, is that of fully understanding the proof of Di Francesco's formula, following his paper, or one of the subsequent other proofs. In relation with the subfactor approach, the question here is that of understanding, following Jones, the positivity of the Temperley-Lieb algebra trace. Finally, in relation with quantum group deformations, check here Wang \cite{wa1}, and related papers. 

\chapter{Twisted permutations}

\section*{4a. Symmetry groups}

We investigate here, following \cite{ba1}, \cite{wa2} and subsequent papers, the quantum symmetry groups $S_Z^+$ of the finite quantum spaces $Z$, generalizing the quantum group $S_N^+$, coming from $Z=\{1,\ldots,N\}$. Besides providing a useful extension of our results regarding $S_N^+$, this will eventually explain the connection with $SO_3$, in an elegant way. As a bonus, we will obtain as well a conceptual result on the connection between $S_N^+$ and $O_N^+$. 

\bigskip

In order to get started, we must talk about finite quantum spaces. In view of the general $C^*$-algebra theory explained in chapter 1, we have the following definition:

\index{quantum space}
\index{finite quantum space}
\index{counting measure}

\begin{definition}
A finite quantum space $Z$ is the abstract dual of a finite dimensional $C^*$-algebra $B$, according to the following formula:
$$C(Z)=B$$
The formal number of elements of such a space is $|Z|=\dim B$. By decomposing the algebra $B$, we have a formula of the following type:
$$C(Z)=M_{n_1}(\mathbb C)\oplus\ldots\oplus M_{n_k}(\mathbb C)$$
With $n_1=\ldots=n_k=1$ we obtain in this way the space $Z=\{1,\ldots,k\}$. Also, when $k=1$ the equation is $C(Z)=M_n(\mathbb C)$, and the solution will be denoted $Z=M_n$.
\end{definition}

In order to talk now about the quantum symmetry group $S_Z^+$, we must use universal coactions. As in chapter 1 when defining $S_N^+$, when dealing with universal coactions on the space $Z=\{1,\ldots,N\}$, we must endow our space $Z$ with its counting measure:

\begin{definition}
We endow each finite quantum space $Z$ with its counting measure, corresponding as the algebraic level to the integration functional
$$tr:C(Z)\to B(l^2(Z))\to\mathbb C$$
obtained by applying the regular representation, and then the normalized matrix trace.
\end{definition}

To be more precise, consider the algebra $B=C(Z)$, which is by definition finite dimensional. We can make act $B$ on itself, by left multiplication:
$$\pi:B\to\mathcal L(B)\quad,\quad 
a\to(b\to ab)$$

The target of $\pi$ being a matrix algebra, $\mathcal L(B)\simeq M_N(\mathbb C)$ with $N=\dim B$, we can further compose with the normalized matrix trace, and we obtain $tr$:
$$tr=\frac{1}{N}\,Tr\circ\pi$$

As basic examples, for both $Z=\{1,\ldots,N\}$ and $Z=M_N$ we obtain the usual trace. In general, with $C(Z)=M_{n_1}(\mathbb C)\oplus\ldots\oplus M_{n_k}(\mathbb C)$, the weights of $tr$ are:
$$c_i=\frac{n_i^2}{\sum_in_i^2}$$

Let us study now the quantum group actions $G\curvearrowright Z$. If we denote by $\mu,\eta$ the multiplication and unit map of the algebra $C(Z)$, we have the following result:

\index{coaction}

\begin{proposition}
Consider a linear map $\Phi:C(Z)\to C(Z)\otimes C(G)$, written as
$$\Phi(e_i)=\sum_je_j\otimes u_{ji}$$
with $\{e_i\}$ being a linear space basis of $C(Z)$, chosen orthonormal with respect to $tr$.
\begin{enumerate}
\item $\Phi$ is a linear space coaction $\iff$ $u$ is a corepresentation.

\item $\Phi$ is multiplicative $\iff$ $\mu\in Hom(u^{\otimes 2},u)$.

\item $\Phi$ is unital $\iff$ $\eta\in Hom(1,u)$.

\item $\Phi$ leaves invariant $tr$ $\iff$ $\eta\in Hom(1,u^*)$.

\item If these conditions hold, $\Phi$ is involutive $\iff$ $u$ is unitary.
\end{enumerate}
\end{proposition}

\begin{proof}
This is similar to the proof for $S_N^+$ from chapter 1, as follows:

\medskip

(1) There are two axioms to be processed here, and we have indeed:
$$(id\otimes\Delta)\Phi=(\Phi\otimes id)\Phi
\iff\Delta(u_{ji})=\sum_ku_{jk}\otimes u_{ki}$$
$$(id\otimes\varepsilon)\Phi=id
\iff\varepsilon(u_{ji})=\delta_{ji}$$

(2) By using $\Phi(e_i)=u(e_i\otimes 1)$ we have the following identities, which give the result:
$$\Phi(e_ie_k)
=u(\mu\otimes id)(e_i\otimes e_k\otimes 1)$$
$$\Phi(e_i)\Phi(e_k)
=(\mu\otimes id)u^{\otimes 2}(e_i\otimes e_k\otimes 1)$$

(3) From $\Phi(e_i)=u(e_i\otimes1)$ we obtain by linearity, as desired:
$$\Phi(1)=u(1\otimes1)$$

(4) This follows from the following computation, by applying the involution:
\begin{eqnarray*}
(tr\otimes id)\Phi(e_i)=tr(e_i)1
&\iff&\sum_jtr(e_j)u_{ji}=tr(e_i)1\\
&\iff&\sum_ju_{ji}^*1_j=1_i\\
&\iff&(u^*1)_i=1_i\\
&\iff&u^*1=1
\end{eqnarray*}

(5) Assuming that (1-4) are satisfied, and that $\Phi$ is involutive, we have:
\begin{eqnarray*}
(u^*u)_{ik}
&=&\sum_lu_{li}^*u_{lk}\\
&=&\sum_{jl}tr(e_j^*e_l)u_{ji}^*u_{lk}\\
&=&(tr\otimes id)\sum_{jl}e_j^*e_l\otimes u_{ji}^*u_{lk}\\
&=&(tr\otimes id)(\Phi(e_i)^*\Phi(e_k))\\
&=&(tr\otimes id)\Phi(e_i^*e_k)\\
&=&tr(e_i^*e_k)1\\
&=&\delta_{ik}
\end{eqnarray*}

Thus $u^*u=1$, and since we know from (1) that $u$ is a corepresentation, it follows that $u$ is unitary. The proof of the converse is standard too, by using a similar computation.
\end{proof}

Following now \cite{ba1}, \cite{wa2}, we have the following result, extending the basic theory of $S_N^+$ from chapter 1 to the present finite quantum space setting:

\index{twisted quantum permutation}
\index{quantum automorphism group}
\index{finite quantum space}

\begin{theorem}
Given a finite quantum space $Z$, there is a universal compact quantum group $S_Z^+$ acting on $Z$, and leaving the counting measure invariant. We have
$$C(S_Z^+)=C(U_N^+)\Big/\Big<\mu\in Hom(u^{\otimes2},u),\eta\in Fix(u)\Big>$$
where $N=|Z|$, and where $\mu,\eta$ are the multiplication and unit maps of the algebra $C(Z)$. For the classical space $Z=\{1,\ldots,N\}$ we have $S_Z^+=S_N^+$. 
\end{theorem}

\begin{proof}
Here the first two assertions follow from Proposition 4.3, by using the standard fact that the complex conjugate of a corepresentation is a corepresentation too. As for the last assertion, regarding $S_N^+$, this follows from the results in chapter 1.
\end{proof}

The above result is quite conceptual, and we will see some applications in a moment. However, for many concrete questions, nothing beats multimatrix bases and indices. So, following the original paper of Wang \cite{wa2}, let us discuss this. We first have:

\begin{definition}
Given a finite quantum space $Z$, we let $\{e_i\}$ be the standard basis of $B=C(Z)$, so that the multiplication, involution and unit of $B$ are given by
$$e_ie_j=e_{ij}\quad,\quad 
e_i^*=e_{\bar{i}}\quad,\quad
1=\sum_{i=\bar{i}}e_i$$
where $(i,j)\to ij$ is the standard partially defined multiplication on the indices, with the convention $e_\emptyset=0$, and where $i\to\bar{i}$ is the standard involution on the indices.
\end{definition}

To be more precise, let $\{e_{ab}^r\}\subset B$ be the multimatrix basis. We set $i=(abr)$, and with this convention, the multiplication, coming from $e_{ab}^re_{cd}^p=\delta_{rp}\delta_{bc}e_{ad}^r$, is given by:
$$(abr)(cdp)=\begin{cases}
(adr)&{\rm if}\ b=c,\ r=p\\
\emptyset&{\rm otherwise}
\end{cases}$$

As for the involution, coming from $(e_{ab}^r)^*=e_{ba}^r$, this is given by:
$$\overline{(a,b,r)}=(b,a,r)$$

Finally, the unit formula comes from the following formula for the unit $1\in B$:
$$1=\sum_{ar}e_{aa}^r$$

Regarding now the generalized quantum permutation groups $S_Z^+$, the construction in Theorem 4.4 reformulates as follows, by using the above formalism:

\begin{proposition}
Given a finite quantum space $Z$, with basis $\{e_i\}\subset C(Z)$ as above, the algebra $C(S_Z^+)$ is generated by variables $u_{ij}$ with the following relations,
$$\sum_{ij=p}u_{ik}u_{jl}=u_{p,kl}\quad,\quad\sum_{kl=p}u_{ik}u_{jl}=u_{ij,p}$$
$$\sum_{i=\bar{i}}u_{ij}=\delta_{j\bar{j}}\quad,\quad\sum_{j=\bar{j}}u_{ij}=\delta_{i\bar{i}}$$
$$u_{ij}^*=u_{\bar{i}\hskip0.3mm\bar{j}}$$
with the fundamental corepresentation being the matrix $u=(u_{ij})$. We call a matrix $u=(u_{ij})$ satisfying the above relations ``generalized magic''.
\end{proposition}

\begin{proof}
We recall from Theorem 4.4 that the algebra $C(S_Z^+)$ appears as follows, where $N=|Z|$, and where $\mu,\eta$ are the multiplication and unit maps of $C(Z)$:
$$C(S_F^+)=C(U_N^+)\Big/\Big<\mu\in Hom(u^{\otimes2},u),\eta\in Fix(u)\Big>$$

But the relations $\mu\in Hom(u^{\otimes2},u)$ and $\eta\in Fix(u)$ produce the 1st and 4th relations in the statement, then the biunitarity of $u$ gives the 5th relation, and finally the 2nd and 3rd relations follow from the 1st and 4th relations, by using the antipode.
\end{proof}

As an illustration, consider the case $Z=\{1,\ldots,N\}$. Here the index multiplication is $ii=i$ and $ij=\emptyset$ for $i\neq j$, and the involution is $\bar{i}=i$. Thus, our relations are as follows, corresponding to the standard magic conditions on a matrix $u=(u_{ij})$:
$$u_{ik}u_{il}=\delta_{kl}u_{ik}\quad,\quad u_{ik}u_{jk}=\delta_{ij}u_{ik}$$
$$\sum_iu_{ij}=1\quad,\quad\sum_ju_{ij}=1$$
$$u_{ij}^*=u_{ij}$$

As a second illustration now, which is something new, we have:

\begin{theorem}
For the space $Z=M_2$, coming via $C(Z)=M_2(\mathbb C)$, we have 
$$S_Z^+=SO_3$$
with the action $SO_3\curvearrowright M_2(\mathbb C)$ being the standard one, coming from $SU_2\to SO_3$.
\end{theorem}

\begin{proof}
This is something quite tricky, the idea being as follows:

\medskip

(1) First, we have an action by conjugation $SU_2\curvearrowright M_2(\mathbb C)$, and this action produces, via the canonical quotient map $SU_2\to SO_3$, an action as follows:
$$SO_3\curvearrowright M_2(\mathbb C)$$

(2) Then, it is routine to check, by using computations like those from the proof of $S_N^+=S_N$ at $N\leq3$, from chapter 1, that any action $G\curvearrowright M_2(\mathbb C)$ must come from a classical group. Thus the action $SO_3\curvearrowright M_2(\mathbb C)$ is universal, as claimed.
\end{proof}

Let us develop now some basic theory for the quantum symmetry groups $S_Z^+$, and their closed subgroups $G\subset S_Z^+$. We have here the following key result:

\index{quantum automorphism group}
\index{quantum symmetry group}
\index{Temperley-Lieb algebra}
\index{Marchenko-Pastur law}
\index{fusion rules}

\begin{theorem}
The quantum groups $S_Z^+$ have the following properties: 
\begin{enumerate}
\item The associated Tannakian categories are $TL_N$, with $N=|Z|$.

\item The main character follows the Marchenko-Pastur law $\pi_1$, when $|Z|\geq4$.

\item The fusion rules for $S_Z^+$ with $|Z|\geq4$ are the same as for $SO_3$.
\end{enumerate}
\end{theorem}

\begin{proof}
This result is from \cite{ba1}, the idea being as follows:

\medskip

(1) Let us pick our orthogonal basis $\{e_i\}$ as in Definition 4.5, so that we have $e_i^*=e_{\bar{i}}$, for a certain involution $i\to\bar{i}$ on the index set. With this convention, we have:
\begin{eqnarray*}
\Phi(e_i)=\sum_je_j\otimes u_{ji}
&\implies&\Phi(e_i)^*=\sum_je_j^*\otimes u_{ji}^*\\
&\implies&\Phi(e_{\bar{i}})=\sum_je_{\bar{j}}\otimes u_{ji}^*\\
&\implies&\Phi(e_i)=\sum_je_j\otimes u_{\bar{i}\hskip0.3mm\bar{j}}^*
\end{eqnarray*}

Thus $u_{ji}^*=u_{\bar{i}\hskip0.3mm\bar{j}}$, so $u\sim\bar{u}$. Now with this result in hand, the proof goes as for the proof for $S_N^+$, from chapter 3. To be more precise, the result follows from the fact that the multiplication and unit of any complex algebra, and in particular of the algebra $C(Z)$ that we are interested in here, can be modeled by the following two diagrams:
$$m=|\cup|\qquad,\qquad u=\cap$$

Indeed, this is certainly true algebrically, and well-known, with as an illustration here, the associativity formula $m(m\otimes id)=(id\otimes m)m$ being checked as follows:
$$\begin{matrix}
|&\cup&|\ |&|\\
|&&\cup&|
\end{matrix}
\ \ =\ \ \begin{matrix}
|&|\ |&\cup&|\\
|&\cup&&|
\end{matrix}$$

As in what regards the $*$-structure, things here are fine too, because our choice for the trace from Definition 4.2 leads to the following formula regarding the adjoints, corresponding to $mm^*=N$, and so to the basic Temperley-Lieb calculus rule $\bigcirc=N$:
$$\mu\mu^*=N\cdot id$$

We conclude that the Tannakian category associated to $S_Z^+$ is, as claimed:
$$C=<\mu,\eta>=<m,u>=<|\cup|,\cap>=TL_N$$

(2) The proof here is exactly as for $S_N^+$, by using moments. To be more precise, according to (1) these moments are the Catalan numbers, which are the moments of $\pi_1$.

\medskip

(3) Once again same proof as for $S_N^+$, by using the fact that the moments of $\chi$ are the Catalan numbers, which naturally leads to the Clebsch-Gordan rules.
\end{proof}

We can merge and reformulate our main results so far in the following way:

\index{projective version}

\begin{theorem}
The quantun groups $S_Z^+$ have the following properties: 
\begin{enumerate}
\item For $Z=\{1,\ldots,N\}$ we have $S_Z^+=S_N^+$.

\item For the space $Z=M_N$ we have $S_Z^+=PO_N^+=PU_N^+$.

\item In particular, for the space $Z=M_2$ we have $S_Z^+=SO_3$.

\item The fusion rules for $S_Z^+$ with $|Z|\geq4$ are independent of $Z$.

\item Thus, the fusion rules for $S_Z^+$ with $|Z|\geq4$ are the same as for $SO_3$.
\end{enumerate}
\end{theorem}

\begin{proof}
This is basically a compact form of what has been said above, with a new result added, and with some technicalities left aside, the idea being as follows:

\medskip

(1) This is something that we know from Theorem 4.4.

\medskip

(2) We recall from chapter 3 that we have $PO_N^+=PU_N^+$. Consider the standard vector space action of the free unitary group $U_N^+$, and its adjoint action:
$$U_N^+\curvearrowright\mathbb C^N\quad,\quad 
PU_N^+\curvearrowright M_N(\mathbb C)$$

By universality of $S_{M_N}^+$, we must have inclusions as follows:
$$PO_N^+\subset PU_N^+\subset S_{M_N}^+$$

On the other hand, the main character of $O_N^+$ with $N\geq2$ being semicircular, the main character of $PO_N^+$ must be Marchenko-Pastur. Thus the inclusion $PO_N^+\subset S_{M_N}^+$ has the property that it keeps fixed the law of main character, and by Peter-Weyl we conclude that this inclusion must be an isomorphism, as desired.

\medskip

(3) This is something that we know from Theorem 4.7, and that can be deduced as well from (2), by using the formula $PO_2^+=SO_3$, which is something elementary. Alternatively, this follows without computations from (4) below, because the inclusion of quantum groups $SO_3\subset S_{M_2}^+$ has the property that it preserves the fusion rules.

\medskip

(4) This is something that we know from Theorem 4.8.

\medskip

(5) This follows from (3,4), as already pointed out in Theorem 4.8.
\end{proof}

As another application of our extended formalism, the Cayley theorem for the finite quantum groups holds in the $S_Z^+$ setting. We have indeed the following result:

\index{Cayley theorem}
\index{no-Cayley theorem}
\index{Cayley embedding}

\begin{theorem}
Any finite quantum group $G$ has a Cayley embedding, as follows:
$$G\subset S_G^+$$
However, there are finite quantum groups which are not quantum permutation groups.
\end{theorem}

\begin{proof}
There are two statements here, the idea being as follows:

\medskip

(1) We have an action $G\curvearrowright G$, which leaves invariant the Haar measure. Now since the counting measure is both left and right invariant, it is the Haar measure. We conclude that $G\curvearrowright G$ leaves invariant the counting measure, and so $G\subset S_G^+$, as claimed.

\medskip

(2) Regarding the second assertion, this is something non-trivial, from \cite{bbn}, the simplest counterexample being a certain quantum group $G$ appearing as a split abelian extension associated to the factorization $S_4=\mathbb Z_4S_3$, having cardinality $|G|=24$. 
\end{proof}

Finally, some interesting phenomena appear in the ``homogeneous'' case, where our quantum space is of the form $Z=M_K\times \{1\ldots,L\}$. Here we first have:

\index{wreath product}

\begin{proposition}
The classical symmetry group of $Z=M_K\times \{1\ldots,L\}$ is
$$S_Z=PU_K\wr S_L$$
with on the right a wreath product, equal by definition to $PU_K^L\rtimes S_L$.
\end{proposition}

\begin{proof}
The fact that we have an inclusion $PU_K\wr S_L\subset S_Z$ is standard, and this follows as well by taking the classical version of the inclusion $PU_K^+\wr_*S_L^+\subset S_Z^+$, established below. As for the fact that this inclusion $PU_K\wr S_L\subset S_Z$ is an isomorphism, this can be proved by picking an arbitrary element $g\in S_Z$, and decomposing it.
\end{proof}

In order to discuss the quantum analogue of the above result, we will need a notion of free wreath product. The basic theory here, from \cite{bi2}, is as follows:

\index{free wreath product}

\begin{proposition}
Given closed subgroups $G\subset U_N^+$, $H\subset S_k^+$, with fundamental corepresentations $u,v$, the following construction produces a closed subgroup of $U_{Nk}^+$:
$$C(G\wr_*H)=(C(G)^{*k}*C(H))/<[u_{ij}^{(a)},v_{ab}]=0>$$
In the case where $G,H$ are classical, the classical version of $G\wr_*H$ is the usual wreath product $G\wr H$. Also, when $G$ is a quantum permutation group, so is $G\wr_*H$.
\end{proposition}

\begin{proof}
Consider the matrix $w_{ia,jb}=u_{ij}^{(a)}v_{ab}$, over the quotient algebra in the statement. It is routine to check that $w$ is unitary, and in the case $G\subset S_N^+$, our claim is that this matrix is magic. Indeed, the entries are projections, and we have:
$$\sum_{jb}w_{ia,jb}
=\sum_{jb}u_{ij}^{(a)}v_{ab}
=\sum_bv_{ab}\sum_ju_{ij}^{(a)}
=1$$
$$\sum_{ia}w_{ia,jb}
=\sum_{ia}u_{ij}^{(a)}v_{ab}
=\sum_av_{ab}\sum_iu_{ij}^{(a)}
=1$$

Thus, $G\wr_*H$ is indeed a quantum permutation group, with fundamental corepresentation $w$. Finally, the assertion regarding classical versions is clear too. See \cite{bi2}.
\end{proof}

With the above notion in hand, we can now formulate the following result:

\index{free wreath product}

\begin{theorem}
The quantum symmetry group of $Z=M_K\times \{1\ldots,L\}$ satisfies:
$$PU_K^+\wr_*S_L^+\subset S_Z^+$$
However, this inclusion is not an isomorphism at $K,L\geq2$.
\end{theorem}

\begin{proof}
We have several assertions to be proved, the idea being as follows:

\medskip

(1) The fact that we have $PU_K^+\wr_*S_L^+\subset S_Z^+$ is well-known and routine, by checking the fact that the matrix $w_{ija,klb}=u_{ij,kl}^{(a)}v_{ab}$ is a generalized magic unitary.

\medskip

(2) The inclusion $PU_K^+\wr_*S_L^+\subset S_Z^+$ is not an isomorphism, for instance by using \cite{twa}, along with the fact that $\pi_1\boxtimes\pi_1\neq\pi_1$ where $\pi_1$ is the Marchenko-Pastur law.
\end{proof}

Finally, let us upgrade our previous planar algebra results from chapter 3, to the present quantum symmetry group setting. The result here is as follows:

\index{planar algebra}
\index{Tannakian duality}
\index{fixed point subfactor}

\begin{theorem}
The closed subgroups $G\subset S_Z^+$ are in correspondence with the subalgebras $P\subset P_Z$ of the planar algebra associated to $Z$, with $G\to P$ being given by:
$$P_k=Fix(u^{\otimes k})$$
Moreover, the associated subfactors can be chosen to be fixed point subfactors,
$$A^G\subset (C(Z)\otimes A)^G$$
with the action $G\curvearrowright A$ being assumed to be minimal, $(A^G)'\cap A=\mathbb C$.
\end{theorem}

\begin{proof}
As before with our previous planar algebra results from chapter 3, there is a long story with this result, the idea being as follows:

\medskip

(1) In what regards the statement, the planar algebra $P_Z$ associated to $Z$ is a joint generalization of the spin and tensor planar algebras $\mathcal S_N,\mathcal T_N$, which appear for $Z=\{1,\ldots,N\}$ and $Z=M_N$, and whose construction is via standard tensor calculus. 

\medskip

(2) From a modern perspective, $P_Z$ appears as the planar algebra associated to the bipartite graph of the inclusion $\mathbb C\subset C(Z)$, by using Jones' construction in \cite{jo4}. We will discuss this in detail in chapter 7, when dealing with more general such constructions.

\medskip

(3) In what regards the proof, this can be obtained along the lines of the proofs for $Z=\{1,\ldots,N\}$ and $Z=M_N$, from chapter 3. For the full story here, and generalizations, we refer to Tarrago-Wahl \cite{twa}. We will be back to this, in chapter 7 below.
\end{proof}

\section*{4b. Twisting results}

Let us go back to the case $N=4$. According to our various considerations above, the link between $S_4^+$ and $SO_3$ should come via some sort of twisting. To be more precise, since the classical space $\{1,2,3,4\}$ and the quantum space $M_2$ both have 4 elements, in the formal sense of Definition 4.1, we can expect to have a twisting result, as follows:
$$\{1,2,3,4\}\sim M_2$$

It is possible to be a bit more precise here, by developing some abstract algebra for this, but going ahead now towards what we are interested in, namely quantum permutation groups, this suggests that we should have a twisting relationship, as follows:
$$S_4^+=S^+_{\{1,2,3,4\}}\sim S^+_{M_2}=SO_3$$

We will see that this is indeed the case, with the subject being quite interesting. In order to discuss twisting, in general, let us start with the following construction:

\index{signature map}
\index{number of crossings}

\begin{proposition}
There is a signature map $\varepsilon:P_{even}\to\{-1,1\}$, given by 
$$\varepsilon(\tau)=(-1)^c$$
where $c$ is the number of switches needed to make $\tau$ noncrossing. In addition:
\begin{enumerate}
\item For $\tau\in S_k$, this is the usual signature.

\item For $\tau\in P_2$ we have $(-1)^c$, where $c$ is the number of crossings.

\item For $\tau\leq\pi\in NC_{even}$, the signature is $1$.
\end{enumerate}
\end{proposition}

\begin{proof}
The fact that $\varepsilon$ is indeed well-defined comes from the fact that the number $c$ in the statement is well-defined modulo 2, which is standard combinatorics. In order to prove now the remaining assertion, observe that any partition $\tau\in P(k,l)$ can be put in ``standard form'', by ordering its blocks according to the appearence of the first leg in each block, counting clockwise from top left, and then by performing the switches as for block 1 to be at left, then for block 2 to be at left, and so on:
$$\xymatrix@R=3mm@C=3mm{\circ\ar@/_/@{.}[drr]&\circ\ar@{-}[dddl]&\circ\ar@{-}[ddd]&\circ\\
&&\ar@/_/@{.}[ur]&\\
&&\ar@/^/@{.}[dr]&\\
\circ&\circ\ar@/^/@{.}[ur]&\circ&\circ}
\xymatrix@R=4mm@C=1mm{&\\\to\\&\\& }
\xymatrix@R=3mm@C=3mm{\circ\ar@/_/@{.}[dr]&\circ\ar@{-}[dddl]&\circ&\circ\ar@{-}[dddl]\\
&\ar@/_/@{.}[ur]&&\\
&&\ar@/^/@{.}[dr]&\\
\circ&\circ\ar@/^/@{.}[ur]&\circ&\circ}
\xymatrix@R=4mm@C=1mm{&\\\to\\&\\&}
\xymatrix@R=3mm@C=3mm{\circ\ar@/_/@{.}[r]&\circ&\circ\ar@{-}[dddll]&\circ\ar@{-}[dddl]\\
&&&\\
&&\ar@/^/@{.}[dr]&\\
\circ&\circ\ar@/^/@{.}[ur]&\circ&\circ}
\xymatrix@R=4mm@C=1mm{&\\\to\\&\\& }
\xymatrix@R=3mm@C=3mm{\circ\ar@/_/@{.}[r]&\circ&\circ\ar@{-}[dddll]&\circ\ar@{-}[dddll]\\
&&&\\
&&&\\
\circ&\circ&\circ\ar@/^/@{.}[r]&\circ}$$

With this convention, the proof of the remaining assertions is as follows:

\medskip

(1) For $\tau\in S_k$ the standard form is $\tau'=id$, and the passage $\tau\to id$ comes by composing with a number of transpositions, which gives the signature. 

\medskip

(2) For a general $\tau\in P_2$, the standard form is of type $\tau'=|\ldots|^{\cup\ldots\cup}_{\cap\ldots\cap}$, and the passage $\tau\to\tau'$ requires $c$ mod 2 switches, where $c$ is the number of crossings. 

\medskip

(3) Assuming that $\tau\in P_{even}$ comes from $\pi\in NC_{even}$ by merging a certain number of blocks, we can prove that the signature is 1 by proceeding by recurrence.
\end{proof}

With the above result in hand, we can now formulate:

\index{twisted linear maps}
\index{twisted Kronecker symbols}

\begin{definition}
Associated to a partition $\pi\in P_{even}(k,l)$ is the linear map
$$\bar{T}_\pi(e_{i_1}\otimes\ldots\otimes e_{i_k})=\sum_{j_1\ldots j_l}\bar{\delta}_\pi\begin{pmatrix}i_1&\ldots&i_k\\ j_1&\ldots&j_l\end{pmatrix}e_{j_1}\otimes\ldots\otimes e_{j_l}$$
where the signed Kronecker symbols
$$\bar{\delta}_\pi\in\{-1,0,1\}$$
are given by $\bar{\delta}_\pi=\varepsilon(\tau)$ if $\tau\geq\pi$, and $\bar{\delta}_\pi=0$ otherwise, with $\tau=\ker(^i_j)$.
\end{definition}

In other words, what we are doing here is to add signatures to the usual formula of $T_\pi$. Indeed, observe that the usual formula for $T_\pi$ can be written as folllows:
$$T_\pi(e_{i_1}\otimes\ldots\otimes e_{i_k})=\sum_{j:\ker(^i_j)\geq\pi}e_{j_1}\otimes\ldots\otimes e_{j_l}$$

Now by inserting signs, coming from the signature map $\varepsilon:P_{even}\to\{\pm1\}$, we are led to the following formula, which coincides with the one from Definition 4.16:
$$\bar{T}_\pi(e_{i_1}\otimes\ldots\otimes e_{i_k})=\sum_{\tau\geq\pi}\varepsilon(\tau)\sum_{j:\ker(^i_j)=\tau}e_{j_1}\otimes\ldots\otimes e_{j_l}$$

We have the following key categorical result:

\begin{proposition}
The assignement $\pi\to\bar{T}_\pi$ is categorical, in the sense that
$$\bar{T}_\pi\otimes\bar{T}_\sigma=\bar{T}_{[\pi\sigma]}\quad,\quad 
\bar{T}_\pi \bar{T}_\sigma=N^{c(\pi,\sigma)}\bar{T}_{[^\sigma_\pi]}\quad,\quad  
\bar{T}_\pi^*=\bar{T}_{\pi^*}$$
where $c(\pi,\sigma)$ are certain positive integers.
\end{proposition}

\begin{proof}
We have to check three conditions, as follows:

\medskip

\underline{1. Concatenation}. It is enough to check the following formula:
$$\varepsilon\left(\ker\begin{pmatrix}i_1\ldots i_p\\ j_1\ldots j_q\end{pmatrix}\right)
\varepsilon\left(\ker\begin{pmatrix}k_1\ldots k_r\\ l_1\ldots l_s\end{pmatrix}\right)=
\varepsilon\left(\ker\begin{pmatrix}i_1\ldots i_p&k_1\ldots k_r\\ j_1\ldots j_q&l_1\ldots l_s\end{pmatrix}\right)$$

Let us denote by $\tau,\nu$ the partitions on the left, so that the partition on the right is of the form $\rho\leq[\tau\nu]$. Now by switching to the noncrossing form, $\tau\to\tau'$ and $\nu\to\nu'$, the partition on the right transforms into $\rho\to\rho'\leq[\tau'\nu']$. Now since $[\tau'\nu']$ is noncrossing, we can use Proposition 4.15 (3), and we obtain the result.

\medskip

\underline{2. Composition}. Here we must establish the following formula:
$$\varepsilon\left(\ker\begin{pmatrix}i_1\ldots i_p\\ j_1\ldots j_q\end{pmatrix}\right)
\varepsilon\left(\ker\begin{pmatrix}j_1\ldots j_q\\ k_1\ldots k_r\end{pmatrix}\right)
=\varepsilon\left(\ker\begin{pmatrix}i_1\ldots i_p\\ k_1\ldots k_r\end{pmatrix}\right)$$

Let $\tau,\nu$ be the partitions on the left, so that the partition on the right is of the form $\rho\leq[^\tau_\nu]$. Our claim is that we can jointly switch $\tau,\nu$ to the noncrossing form. Indeed, we can first switch as for $\ker(j_1\ldots j_q)$ to become noncrossing, and then switch the upper legs of $\tau$, and the lower legs of $\nu$, as for both these partitions to become noncrossing. Now observe that when switching in this way to the noncrossing form, $\tau\to\tau'$ and $\nu\to\nu'$, the partition on the right transforms into $\rho\to\rho'\leq[^{\tau'}_{\nu'}]$. Now since $[^{\tau'}_{\nu'}]$ is noncrossing, we can apply Proposition 4.15 (3), and we obtain the result.

\medskip

\underline{3. Involution}. Here we must prove the following formula:
$$\bar{\delta}_\pi\begin{pmatrix}i_1\ldots i_p\\ j_1\ldots j_q\end{pmatrix}=\bar{\delta}_{\pi^*}\begin{pmatrix}j_1\ldots j_q\\ i_1\ldots i_p\end{pmatrix}$$

But this is clear from the definition of $\bar{\delta}_\pi$, and we are done.
\end{proof}

As a conclusion, our construction $\pi\to\bar{T}_\pi$ has all the needed properties for producing quantum groups, via Tannakian duality. So, we can now formulate:

\begin{theorem}
Given a category of partitions $D\subset P_{even}$, the construction
$$Hom(u^{\otimes k},u^{\otimes l})=span\left(\bar{T}_\pi\Big|\pi\in D(k,l)\right)$$
produces via Tannakian duality a quantum group $G_N^{-1}\subset U_N^+$, for any $N\in\mathbb N$.
\end{theorem}

\begin{proof}
This follows indeed from the Tannakian results from chapter 2, exactly as in the easy case, by using this time Proposition 4.17 as technical ingredient.
\end{proof}

We can unify the easy quantum groups, or at least the examples coming from categories $D\subset P_{even}$, with the quantum groups constructed above, as follows:

\index{quizzy quantum group}
\index{Schur-Weyl twist}

\begin{definition}
A closed subgroup $G\subset U_N^+$ is called $q$-easy, or quizzy, with deformation parameter $q=\pm1$, when its tensor category appears as
$$Hom(u^{\otimes k},u^{\otimes l})=span\left(\dot{T}_\pi\Big|\pi\in D(k,l)\right)$$
for a certain category of partitions $D\subset P_{even}$, where, for $q=-1,1$: 
$$\dot{T}=\bar{T},T$$
The Schur-Weyl twist of $G$ is the quizzy quantum group $G^{-1}\subset U_N^+$ obtained via $q\to-q$.
\end{definition}

In order to compute now the twists of the basic compact groups, we recall that the M\"obius function of any lattice, and in particular of $P_{even}$, is given by:
$$\mu(\sigma,\pi)=\begin{cases}
1&{\rm if}\ \sigma=\pi\\
-\sum_{\sigma\leq\tau<\pi}\mu(\sigma,\tau)&{\rm if}\ \sigma<\pi\\
0&{\rm if}\ \sigma\not\leq\pi
\end{cases}$$

With this notation, we have the following result:

\index{M\"obius inversion}

\begin{proposition}
For any partition $\pi\in P_{even}$ we have the formula
$$\bar{T}_\pi=\sum_{\tau\leq\pi}\alpha_\tau T_\tau$$
where $\alpha_\sigma=\sum_{\sigma\leq\tau\leq\pi}\varepsilon(\tau)\mu(\sigma,\tau)$, with $\mu$ being the M\"obius function of $P_{even}$.
\end{proposition}

\begin{proof}
The linear combinations $T=\sum_{\tau\leq\pi}\alpha_\tau T_\tau$ acts on tensors as follows:
\begin{eqnarray*}
T(e_{i_1}\otimes\ldots\otimes e_{i_k})
&=&\sum_{\tau\leq\pi}\alpha_\tau T_\tau(e_{i_1}\otimes\ldots\otimes e_{i_k})\\
&=&\sum_{\tau\leq\pi}\alpha_\tau\sum_{\sigma\leq\tau}\sum_{j:\ker(^i_j)=\sigma}e_{j_1}\otimes\ldots\otimes e_{j_l}\\
&=&\sum_{\sigma\leq\pi}\left(\sum_{\sigma\leq\tau\leq\pi}\alpha_\tau\right)\sum_{j:\ker(^i_j)=\sigma}e_{j_1}\otimes\ldots\otimes e_{j_l}
\end{eqnarray*}

Thus, in order to have $\bar{T}_\pi=\sum_{\tau\leq\pi}\alpha_\tau T_\tau$, we must have, for any $\sigma\leq\pi$:
$$\varepsilon(\sigma)=\sum_{\sigma\leq\tau\leq\pi}\alpha_\tau$$

But this problem can be solved by using the M\"obius inversion formula, and we obtain the numbers $\alpha_\sigma=\sum_{\sigma\leq\tau\leq\pi}\varepsilon(\tau)\mu(\sigma,\tau)$ in the statement.
\end{proof}

We can now twist $O_N$, and $U_N$ as well. The result here is as follows:

\index{twisted orthogonal group}
\index{twisted unitary group}

\begin{theorem}
The twists of $O_N,U_N$ are obtained by replacing the commutation relations $ab=ba$ between the coordinates $u_{ij}$ and their adjoints $u_{ij}^*$ with the relations
$$ab=\pm ba$$
with anticommutation on rows and columns, and commutation otherwise.
\end{theorem}

\begin{proof}
The basic crossing, $\ker\binom{ij}{ji}$ with $i\neq j$, comes from the transposition $\tau\in S_2$, so its signature is $-1$. As for its degenerated version $\ker\binom{ii}{ii}$, this is noncrossing, so here the signature is $1$. We conclude that the linear map associated to the basic crossing is:
$$\bar{T}_{\slash\!\!\!\backslash}(e_i\otimes e_j)
=\begin{cases}
-e_j\otimes e_i&{\rm for}\ i\neq j\\
e_j\otimes e_i&{\rm otherwise}
\end{cases}$$

We can proceed now as in the untwisted case, and since the intertwining relations coming from $\bar{T}_{\slash\!\!\!\backslash}$ correspond to the relations defining $O_N^{-1},U_N^{-1}$, we obtain the result.
\end{proof}

Now going towards $S_4^+$, let us start with the following definition, from \cite{bb2}:

\index{twisting}
\index{twisted determinant}
\index{anticommutation}

\begin{definition}
We let $SO_3^{-1}\subset O_3^{-1}$ be the subgroup coming from the relation
$$\sum_{\sigma\in S_3}u_{1\sigma(1)}u_{2\sigma(2)}u_{3\sigma(3)}=1$$
called twisted determinant one condition.
\end{definition}

Normally, we should prove here that $C(SO_3^{-1})$ is indeed a Woronowicz algebra. This is of course possible, by doing some computations, but we will not need to do these computations, because this follows from the following result, also from \cite{bb2}:

\index{Klein group}

\begin{theorem}
We have an isomorphism of compact quantum groups
$$S_4^+=SO_3^{-1}$$
given by the Fourier transform over the Klein group $K=\mathbb Z_2\times\mathbb Z_2$.
\end{theorem}

\begin{proof}
Consider the following matrix, coming from the action of $SO_3^{-1}$ on $\mathbb C^4$:
$$u^+=\begin{pmatrix}1&0\\0&u\end{pmatrix}$$

We apply to this matrix the Fourier transform over the Klein group $K=\mathbb Z_2\times\mathbb Z_2$: 
$$v=
\frac{1}{4}
\begin{pmatrix}
1&1&1&1\\
1&-1&-1&1\\
1&-1&1&-1\\
1&1&-1&-1
\end{pmatrix}
\begin{pmatrix}
1&0&0&0\\
0&u_{11}&u_{12}&u_{13}\\
0&u_{21}&u_{22}&u_{23}\\
0&u_{31}&u_{32}&u_{33}
\end{pmatrix}
\begin{pmatrix}
1&1&1&1\\
1&-1&-1&1\\
1&-1&1&-1\\
1&1&-1&-1
\end{pmatrix}$$

This matrix is then magic, and vice versa, so the Fourier transform over $K$ converts the relations in Definition 4.22 into the magic relations. But this gives the result.
\end{proof}

There are many more things that can be said here, and we have:

\begin{theorem}
The quantum group $S_4^+=SO_3^{-1}$ has the following properties:
\begin{enumerate}
\item It appears as a cocycle twist of $SO_3$.

\item Its fusion rules are the same as for $SO_3$.

\item Its subgroups are basically twists of the subgroups of $SO_3$.
\end{enumerate}
\end{theorem}

\begin{proof}
These are more advanced results, from \cite{bb2}, the idea being as follows:

\medskip

(1) This follows by suitably reformulating the definition of $SO_3^{-1}$ given above in purely algebraic terms, using cocycles, and for details here, we refer to \cite{bb2}. In what concerns us, we will actually discuss a generalization of this, right next, following \cite{bbs}.

\medskip

(2) This is something that we know well, via numerous proofs, and we can add to our trophy list one more proof, coming from (1), via standard cocycle twisting theory.

\medskip

(3) The idea here is that the closed subgroups $G\subset SO_3$ are subject to a well-known ADE/McKay classification result, and the subgroups of $SO_3^{-1}$ are basically twists of these, $G^{-1}\subset SO_3^{-1}$. We will discuss this, following \cite{bb2}, in chapter 10 below.
\end{proof}

\section*{4c. Cocycle twisting}

An interesting extension of the $S_4^+=SO_3^{-1}$ result comes by looking at the general case $N=n^2$, with $n\in\mathbb N$. We will prove that we have a twisting result, as follows: 
$$PO_n^+=(S_N^+)^\sigma$$

In order to explain this material, from \cite{bbs}, which is quite technical, requiring good algebraic knowledge, let us begin with some generalities. We first have:

\begin{proposition}
Given a finite group $G$, the algebra $C(S_{\widehat{G}}^+)$ is isomorphic to the abstract algebra presented by generators $x_{gh}$ with $g,h\in G$, with the following relations:
$$x_{1g}=x_{g1}=\delta_{1g}\quad,\quad
x_{s,gh}=\sum_{t\in G}x_{st^{-1},g}x_{th}\quad,\quad 
x_{gh,s}=\sum_{t\in G}x_{gt^{-1}}x_{h,ts}$$
The comultiplication, counit and antipode are given by the formulae
$$\Delta(x_{gh})=\sum_{s\in G}x_{gs}\otimes x_{sh}\quad,\quad 
\varepsilon(x_{gh})=\delta_{gh}\quad,\quad
S(x_{gh})=x_{h^{-1}g^{-1}}$$
on the standard generators $x_{gh}$.
\end{proposition}

\begin{proof}
This follows indeed from a direct verification, based either on Theorem 4.4, or on its equivalent formulation from Proposition 4.6. See \cite{bbs}.
\end{proof}

Let us discuss now the twisted version of the above result. Consider a 2-cocycle on $G$, which is by definition a map $\sigma:G\times G\to\mathbb C^*$ satisfying:
$$\sigma_{gh,s}\sigma_{gh}=\sigma_{g,hs}\sigma_{hs}\quad,\quad 
\sigma_{g1}=\sigma_{1g}=1$$

Given such a cocycle, we can construct the associated twisted group algebra $C(\widehat{G}_\sigma)$, as being the vector space $C(\widehat{G})=C^*(G)$, with product $e_ge_h=\sigma_{gh}e_{gh}$. We have:

\begin{proposition}
The algebra $C(S_{\widehat{G}_\sigma}^+)$ is isomorphic to the abstract algebra presented by generators $x_{gh}$ with $g,h\in G$, with the relations $x_{1g}=x_{g1}=\delta_{1g}$ and:
$$\sigma_{gh}x_{s,gh}=\sum_{t\in G}\sigma_{st^{-1},t}x_{st^{-1},g}x_{th}\quad,\quad
\sigma_{gh}^{-1}x_{gh,s}=\sum_{t\in G}\sigma_{t^{-1},ts}^{-1}x_{gt^{-1}}x_{h,ts}$$
The comultiplication, counit and antipode are given by the formulae
$$\Delta(x_{gh})=\sum_{s\in G}x_{gs}\otimes x_{sh}\quad,\quad 
\varepsilon(x_{gh})=\delta_{gh}\quad,\quad
S(x_{gh})=\sigma_{h^{-1}h}\sigma_{g^{-1}g}^{-1}x_{h^{-1}g^{-1}}$$
on the standard generators $x_{gh}$.
\end{proposition}

\begin{proof}
Once again, this follows from a direct verification, explained in \cite{bbs}.
\end{proof}

\index{twisting}
\index{cocycle twisting}

In what follows, we will prove that the quantum groups $S_{\widehat{G}}^+$ and $S_{\widehat{G}_\sigma}^+$ are related by a cocycle twisting operation. Let $A$ be a Hopf algebra. We recall that a left 2-cocycle is a convolution invertible linear map
$\sigma:A\otimes A\to\mathbb C$ satisfying:
$$\sigma_{x_1y_1}\sigma_{x_2y_2,z}=\sigma_{y_1z_1}\sigma_{x,y_2z_2}\quad,\quad 
\sigma_{x1}=\sigma_{1x}=\varepsilon(x)$$

Note that $\sigma$ is a left 2-cocycle if and only if $\sigma^{-1}$, the convolution inverse of $\sigma$, is a right 2-cocycle, in the sense that we have:
$$\sigma^{-1}_{x_1y_1,z}\sigma^{-1}_{x_1y_2}=\sigma^{-1}_{x,y_1z_1}\sigma^{-1}_{y_2z_2}\quad,\quad 
\sigma^{-1}_{x1}=\sigma^{-1}_{1x}=\varepsilon(x)$$

Given a left 2-cocycle $\sigma$ on $A$, one can form the 2-cocycle twist $A^\sigma$ as follows. As a coalgebra, $A^\sigma=A$, and an element $x\in A$, when considered in $A^\sigma$, is denoted $[x]$. The product in $A^\sigma$ is then defined, in Sweedler notation, by: 
$$[x][y]=\sum\sigma_{x_1y_1}\sigma^{-1}_{x_3y_3}[x_2y_2]$$

We can now state and prove a main theorem from \cite{bbs}, as follows:

\begin{theorem}
If $G$ is a finite group and $\sigma$ is a $2$-cocycle on $G$, the Hopf algebras
$$C(S_{\widehat{G}}^+)\quad,\quad C(S_{\widehat{G}_\sigma}^+)$$
are $2$-cocycle twists of each other, in the above sense.
\end{theorem}

\begin{proof}
In order to prove this result, we use the following Hopf algebra map: 
$$\pi:C(S_{\widehat{G}}^+)\to C(\widehat{G})\quad,\quad
x_{gh}\to\delta_{gh}e_g$$

Our 2-cocycle $\sigma:G\times G\to\mathbb C^*$ can be extended by linearity into a linear map as follows, which is a left and right 2-cocycle in the above sense:
$$\sigma:C(\widehat{G})\otimes C(\widehat{G})\to\mathbb C$$

Consider now the following composition:
$$\alpha=\sigma(\pi \otimes \pi):C(S_{\widehat{G}}^+)\otimes C(S_{\widehat{G}}^+)\to C(\widehat{G})\otimes C(\widehat{G})\to\mathbb C$$

Then $\alpha$ is a left and right 2-cocycle, because it is induced by a cocycle on a group algebra, and so is its convolution inverse $\alpha^{-1}$. Thus we can construct the twisted algebra $C(S_{\widehat{G}}^+)^{\alpha^{-1}}$, and inside this algebra we have the following computation:
$$[x_{gh}][x_{rs}]
=\alpha^{-1}(x_g,x_r)\alpha(x_h,x_s)[x_{gh}x_{rs}]
=\sigma_{gr}^{-1}\sigma_{hs}[x_{gh}x_{rs}]$$

By using this, we obtain the following formula:
$$\sum_{t\in G}\sigma_{st^{-1},t}[x_{st^{-1},g}][x_{th}]
=\sum_{t\in G}\sigma_{st^{-1},t}\sigma_{st^{-1},t}^{-1}\sigma_{gh}[x_{st^{-1},g}x_{th}]
=\sigma_{gh}[x_{s,gh}]$$

Similarly, we have the following formula:
$$\sum_{t\in G}\sigma_{t^{-1},ts}^{-1}[x_{g,t^{-1}}][x_{h,ts}]=\sigma_{gh}^{-1}[x_{gh,s}]$$

We deduce from this that there exists a Hopf algebra map, as follows:
$$\Phi:C(S_{\widehat{G}_\sigma}^+)\to C(S_{\widehat{G}}^+)^{\alpha^{-1}}\quad,\quad 
x_{gh}\to [x_{g,h}]$$

This map is clearly surjective, and is injective as well, by a standard fusion semiring argument, because both Hopf algebras have the same fusion semiring.
\end{proof}

Summarizing, we have proved our main twisting result. Our purpose in what follows will be that of working out versions and particular cases of it. We first have: 

\begin{proposition}
If $G$ is a finite group and $\sigma$ is a $2$-cocycle on $G$, then
$$\Phi(x_{g_1h_1}\ldots x_{g_mh_m})=\Omega(g_1,\ldots,g_m)^{-1}\Omega(h_1,\ldots,h_m)x_{g_1h_1}\ldots x_{g_mh_m}$$
with the coefficients on the right being given by the formula
$$\Omega(g_1,\ldots,g_m)=\prod_{k=1}^{m-1}\sigma_{g_1\ldots g_k,g_{k+1}}$$
is a coalgebra isomorphism $C(S_{\widehat{G}_\sigma}^+)\to C(S_{\widehat{G}}^+)$, commuting with the Haar integrals.
\end{proposition}

\begin{proof}
This is indeed just a technical reformulation of Theorem 4.27.
\end{proof}

Let us discuss now some concrete applications of the general results established above. Consider the group $G=\mathbb Z_n^2$, let $w=e^{2\pi i/n}$, and consider the following cocycle: 
$$\sigma:G\times G\to\mathbb C^*\quad,\quad 
\sigma_{(ij)(kl)}=w^{jk}$$ 

In order to understand what is the formula that we obtain, we must do some computations. Let $E_{ij}$ with $i,j \in\mathbb Z_n$ be the standard basis of $M_n(\mathbb C)$. We first have:

\begin{proposition}
The linear map given by
$$\psi(e_{(i,j)})=\sum_{k=0}^{n-1}{w}^{ki}E_{k,k+j}$$
defines an isomorphism of algebras $\psi:C(\widehat{G}_\sigma)\simeq M_n(\mathbb C)$. 
\end{proposition}

\begin{proof}
Consider indeed the following linear map:
$$\psi'(E_{ij})=\frac{1}{n}\sum_{k=0}^{n-1}{w}^{-ik}e_{(k,j-i)}$$
 
It is routine to check that both $\psi,\psi'$ are morphisms of algebras, and that these maps are inverse to each other. In particular, $\psi$ is an isomorphism of algebras, as stated.
\end{proof}

Next in line, we have the following result:

\begin{proposition}
The algebra map given by
$$\varphi(u_{ij}u_{kl}) = \frac{1}{n}\sum_{a,b=0}^{n-1}{w}^{ai-bj}x_{(a,k-i),(b,l-j)}$$
defines a Hopf algebra isomorphism $\varphi:C(S_{M_n}^+)\simeq C(S_{\widehat{G}_\sigma}^+)$.
\end{proposition}

\begin{proof}
Consider the universal coactions on the two algebras in the statement:
\begin{eqnarray*}
\alpha:M_n(\mathbb C)&\to&M_n({\mathbb C})\otimes C(S_{M_n}^+)\\
\beta:C(\widehat{G}_\sigma)&\to&C(\widehat{G}_\sigma)\otimes C(S_{\widehat{G}_\sigma}^+)
 \end{eqnarray*}
 
In terms of the standard bases, these coactions are given by:
\begin{eqnarray*}
\alpha(E_{ij})&=&\sum_{kl}E_{kl}\otimes u_{ki}u_{lj}\\
\beta(e_{(i,j)})&=&\sum_{kl} e_{(k,l)}\otimes x_{(k,l),(i,j)}
\end{eqnarray*}

We use now the identification $C(\widehat{G}_\sigma)\simeq M_n(\mathbb C)$ from Proposition 4.29. This identification produces a coaction map, as follows:
$$\gamma:M_n(\mathbb C)\to M_n(\mathbb C)\otimes C(S_{\widehat{G}_\sigma}^+)$$

Now observe that this map is given by the following formula:
$$\gamma(E_{ij})=\frac{1}{n}\sum_{ab}E_{ab}\otimes\sum_{kr}w^{ar-ik} x_{(r,b-a),(k,j-i)}$$

By comparing with the formula of $\alpha$, we obtain the isomorphism in the statement.
\end{proof}

We will need one more result of this type, as follows:

\begin{proposition}
The algebra map given by
$$\rho(x_{(a,b),(i,j)})=\frac{1}{n^2}\sum_{klrs}w^{ki+lj-ra-sb}p_{(r,s),(k,l)}$$
defines a Hopf algebra isomorphism $\rho:C(S_{\widehat{G}}^+)\simeq C(S_G^+)$.
\end{proposition}

\begin{proof}
We have a Fourier transform isomorphism, as follows:
$$C(\widehat{G})\simeq C(G)$$

Thus the algebras in the statement are indeed isomorphic.
\end{proof}

As a conclusion to all this, we have the following result, from \cite{bbs}:

\begin{theorem}
Let $n\geq 2$ and $w=e^{2\pi i/n}$. Then
$$\Theta(u_{ij}u_{kl})=\frac{1}{n}\sum_{ab=0}^{n-1}w^{-a(k-i)+b(l-j)}p_{ia,jb}$$
defines a coalgebra isomorphism $C(PO_n^+)\to C(S_{n^2}^+)$ commuting with the Haar integrals.
\end{theorem}
 
\begin{proof}
We know from Theorem 4.9 that we have identifications as follows, where the projective version of $(A,u)$ is the pair $(PA,v)$, with $PA=<v_{ij}>$ and $v=u\otimes\bar{u}$:
$$PO_n^+=PU_n^+=S_{M_n}^+$$

With this in hand, the result follows from Theorem 4.27 and Proposition 4.28, by combining them with the various isomorphisms established above.
\end{proof}

Summarizing, the twisting formula $S_4^+=SO_3^{-1}$ ultimately comes from something of type $X_4\simeq M_2$, where $X_4=\{1,2,3,4\}$ and $M_2=Spec(M_2(\mathbb C))$, and at $N\geq5$ there are some extensions of this, and notably when $N=n^2$ with $n\geq3$. For more on all this, advanced algebraic aspects of $S_N^+$, we refer to \cite{bb2}, \cite{bbs}, \cite{bne}, \cite{beh}.

\section*{4d. Hypergeometric laws}

Still following \cite{bbs}, let us discuss now some probabilistic consequences of the above. We first have the remarkable result, having no classical counterpart:

\begin{theorem}
The following families of variables have the same joint law,
\begin{enumerate}
\item $\{u_{ij}^2\}\in C(O_n^+)$,

\item $\{X_{ij}=\frac{1}{n}\sum_{ab}p_{ia,jb}\}\in C(S_{n^2}^+)$,
\end{enumerate}
where $u=(u_{ij})$ and $p=(p_{ia,jb})$ are the corresponding fundamental corepresentations.
\end{theorem}

\begin{proof}
This follows from Theorem 4.32, but we can recover this as well directly, by using the Weingarten formula for our two quantum groups, and the shrinking operation for partitions $\pi\to\pi'$. Indeed, we have the following moment formulae:
$$\int_{O_n^+}u_{ij}^{2k}=\sum_{\pi,\sigma\in NC_2(2k)}W_{2k,n}(\pi,\sigma)$$
$$\int_{S_{\!n^2}^+}X_{ij}^k=\sum_{\pi,\sigma\in NC_2(2k)}n^{|\pi'|+|\sigma'|-k}W_{k,n^2}(\pi',\sigma')$$

According to the fattening/shrinking results in chapter 3 the summands coincide, and so the moments are equal, as desired. The proof for joint moments is similar.
\end{proof}

In what follows we will be interested in single variables. We have here:

\index{hypergeometric law}
\index{free hypergeometric law}
\index{semicircle law}
\index{free Poisson law}

\begin{definition}
The noncommutative random variable
$$X(n,m,N)=\sum_{i=1}^n \sum_{j=1}^m u_{ij} \in C(S_N^+)$$
is called free hypergeometric, of parameters $(n,m,N)$.
\end{definition}

The terminology comes from the fact that the variable $X'(n,m,N)$, defined as above, but over the algebra $C(S_N)$, follows a hypergeometric law of parameters $(n,m,N)$. Following \cite{bbs}, here is an exploration of the basic properties of these laws:

\begin{theorem}
The free hypergeometric laws have the following properties:
\begin{enumerate}
\item Let $n,m,N\to\infty$, with $\frac{nm}{N}\to t\in(0,\infty)$. Then the law of 
$$X(n,m,N)$$ 
converges to Marchenko-Pastur law $\pi_t$.

\item Let $n,m,N\to\infty$, with $\frac{n}{N}\to\nu\in (0,1)$ and $\frac{m}{N}\to 0$. Then the law of 
$$S(n,m,N)=(X(n,m,N)-m\nu)/\sqrt{m\nu(1-\nu)}$$
converges to the semicircle law $\gamma_1$.
\end{enumerate}
\end{theorem}

\begin{proof}
This is standard, by using the Weingarten formula, as follows:

\medskip

(1) From the Weingarten formula, we have:
$$\int_{S_N^+}X(n,m,N)^k=\sum_{\pi,\sigma\in NC(k)}W_{kN}(\pi,\sigma)n^{|\pi|}m^{|\sigma|}$$

The point now is that we have the following estimate:
$$W_{kN}(\pi,\sigma)= 
\begin{cases}N^{-|\pi|}+O(N^{-|\pi|-1})&{\rm if}\ \pi=\sigma\\ 
O(N^{|\pi\vee\sigma|-|\pi|-|\sigma|})&{\rm if}\ \pi\neq\sigma
\end{cases}$$

It follows that we have the following estimate, which gives the result:
$$W_{kN}(\pi,\sigma) n^{|\pi|}m^{|\sigma|} \to 
\begin{cases}
t^{|\pi|}&{\rm if}\ \pi=\sigma\\ 
0&{\rm if}\ \pi\neq\sigma 
\end{cases}$$

(2) We need to show that the free cumulants satisfy:
$$\kappa^{(p)}[S(n,m,N),\ldots,S(n,m,N)]
\to\begin{cases}
1&{\rm if}\ p=2\\
0&{\rm if}\ p\neq2
\end{cases}$$

The case $p =1$ is trivial, so suppose $p\geq 2$.  We have:
$$\kappa^{(p)}[S(n,m,N),\ldots,S(n,m,N)]
=(m\nu(1-\nu))^{-p/2}\kappa^{(p)}[X(n,m,N),\ldots,X(n,m,N)]$$

On the other hand, from the Weingarten formula, we have:
\begin{eqnarray*}
&&\kappa^{(p)}[X(n,m,N),\dotsc,X(n,m,N)] \\
&=&\sum_{w\in NC(p)} \mu_p({w},1_p) \prod_{V \in {w}} \sum_{\pi_V,\sigma_V \in NC(V)} W_{NC(V),N}(\pi_V,\sigma_V)n^{|\pi_V|}m^{|\sigma_V|}\\
&=&\sum_{{w} \in NC(p)} \mu_p({w},1_p) \prod_{V \in {w}} \sum_{\pi_V,\sigma_V \in NC(V)} (N^{-|\pi_V|}\mu_{|V|}(\pi_V,\sigma_V) + O(N^{-|\pi_V|-1}))n^{|\pi_V|}m^{|\sigma_V|}\\
&=&\sum_{\substack{\pi,\sigma \in NC(p)\\ \pi \leq \sigma}} (N^{-|\pi|}\mu_p(\pi,\sigma) + O(N^{-|\pi|-1}))n^{|\pi|}m^{|\sigma|} \sum_{\substack{{w} \in NC(p)\\ \sigma \leq {w}}} \mu_p({w},1_p)
\end{eqnarray*}

We use now the following standard identity: 
$$\sum_{\substack{{w} \in NC(p)\\ \sigma \leq {w}}} \mu_p({w},1_p)=
\begin{cases}
1&{\rm if}\ \sigma=1_p\\
0&{\rm if}\ \sigma\neq1_p
\end{cases}$$

This gives the following formula for the cumulants:
$$\kappa^{(p)}[X(n,m,N),\ldots,X(n,m,N)]
=m\sum_{\pi \in NC(p)} (N^{-|\pi|}\mu_p(\pi,1_p) + O(N^{-|\pi|-1}))n^{|\pi|}$$

It follows that for $p \geq 3$ we have, as desired:
$$\kappa^{(p)}[S(n,m,N),\dotsc,S(n,m,N)]\to 0$$

As for the remaining case $p = 2$, here we have:
\begin{eqnarray*}
 \kappa^{(2)}[S(n,m,N),S(n,m,N)] 
&\to&\frac{1}{\nu(1-\nu)} \sum_{\pi \in NC(2)} \nu^{|\pi|}\mu_2(\pi,1_2)\\
&=&\frac{1}{\nu(1-\nu)}\bigl(\nu - \nu^2)=1
\end{eqnarray*}

Thus, we are led to the conclusion in the statement.
\end{proof}

Let us discuss as well some similar computations for the free hyperspherical laws, which are quite interesting. To start with, we have the following definition:

\index{free sphere}
\index{free real sphere}

\begin{definition}
The free real sphere is the quantum algebraic manifold given by:
$$C(S^{N-1}_{\mathbb R,+})=C^*\left(x_1,\ldots,x_N\Big|x_i=x_i^*,\sum_i x_i^2=1\right)$$
We endow this sphere with the uniform measure coming from the action $O_N^+\curvearrowright S^{N-1}_{\mathbb R,+}$.
\end{definition}

Here, as usual in the present book, the word ``quantum'' refers to the general framework of Gelfand duality. In the case of $S^{N-1}_{\mathbb R,+}$, however, things are much more concrete, because in analogy with what happens in the classical case, we have an embedding of algebras as follows, after GNS construction with respect to the Haar functional:
$$C(S^{N-1}_{\mathbb R,+})\subset C(O_N^+)\quad,\quad x_i=u_{1i}$$

Thus, the integration questions over $S^{N-1}_{\mathbb R,+}$ correspond to integration questions over $O_N^+$, and we have the following result, that we basically know from chapter 2:

\index{free hyperspherical law}

\begin{theorem}
For the free sphere $S^{N-1}_{\mathbb R,+}$, the rescaled coordinates 
$$y_i=\sqrt{N}x_i$$
become semicircular and free, in the $N\to\infty$ limit.
\end{theorem}

\begin{proof}
By using the above identification $x_i=u_{1i}$, the Weingarten formula for $O_N^+$ gives a Weingarten formula for $S^{N-1}_{\mathbb R,+}$, and together with the standard fact that the Gram matrix, and hence the Weingarten matrix too, is asymptotically diagonal, this gives:
$$\int_{S^{N-1}_{\mathbb R,+}}x_{i_1}\ldots x_{i_k}\,dx\simeq N^{-k/2}\sum_{\sigma\in NC_2(k)}\delta_\sigma(i_1,\ldots,i_k)$$

With this formula in hand, we can compute the asymptotic moments of each coordinate $x_i$. Indeed, by setting $i_1=\ldots=i_k=i$, all Kronecker symbols are 1, and we obtain:
$$\int_{S^{N-1}_{\mathbb R,+}}x_i^k\,dx\simeq N^{-k/2}|NC_2(k)|$$

Thus the rescaled coordinates $y_i=\sqrt{N}x_i$ become semicircular in the $N\to\infty$ limit, as claimed. As for the asymptotic freeness result, this follows as well from the above general joint moment estimate, via standard free probability theory. See \cite{ba4}.
\end{proof}

The problem now, which is highly non-trivial, is that of computing the moments of the coordinates of the free sphere at fixed values of $N\in\mathbb N$. The answer here, from \cite{bcz}, based on advanced quantum group and calculus techniques, is as follows:

\index{free hyperspherical law}
\index{special functions}
\index{deformed quantum groups}

\begin{theorem}
The moments of the free hyperspherical law are given by
$$\int_{S^{N-1}_{\mathbb R,+}}x_1^{2l}\,dx=\frac{1}{(N+1)^l}\cdot\frac{q+1}{q-1}\cdot\frac{1}{l+1}\sum_{r=-l-1}^{l+1}(-1)^r\begin{pmatrix}2l+2\cr l+r+1\end{pmatrix}\frac{r}{1+q^r}$$
where $q\in [-1,0)$ is such that $q+q^{-1}=-N$.
\end{theorem}

\begin{proof}
This is something quite technical, as follows:

\medskip

(1) To any $F\in GL_N(\mathbb R)$ satisfying $F^2=1$ we associate the following algebra:
$$C(O_F^+)=C^*\left((u_{ij})_{i,j=1,\ldots,N}\Big|u=F\bar{u}F={\rm unitary}\right)$$

Observe that we have $O_{I_N}^+=O_N^+$. In general, the above algebra satisfies Woronowicz's generalized axioms in \cite{wo1}, which do not include the strong antipode axiom $S^2=id$. 

\medskip

(2) At $N=2$, up to a trivial equivalence relation on the matrices $F$, and on the quantum groups $O_F^+$, we can assume that $F$ is as follows, with $q\in [-1,0)$:
$$F=\begin{pmatrix}0&\sqrt{-q}\\
1/\sqrt{-q}&0\end{pmatrix}$$

Our claim is that for this matrix $F$ we have an identification $O_F^+=SU^q_2$. Indeed, the relations $u=F\bar{u}F$ tell us that $u$ must be of the following form:
$$u=\begin{pmatrix}\alpha&-q\gamma^*\cr \gamma&\alpha^*\end{pmatrix}$$

Thus $C(O_F^+)$ is the universal algebra generated by two elements $\alpha,\gamma$, with the relations making the above matrix $u$ unitary. But these unitarity conditions are as follows:
$$\alpha\gamma=q\gamma\alpha\quad,\quad
\alpha\gamma^*=q\gamma^*\alpha\quad,\quad 
\gamma\gamma^*=\gamma^*\gamma$$
$$\alpha^*\alpha+\gamma^*\gamma=1\quad,\quad 
\alpha\alpha^*+q^2\gamma\gamma^*=1$$

We recognize here the relations in \cite{wo1} defining the algebra $C(SU^q_2)$, and it follows that we have an isomorphism of Hopf $C^*$-algebras, as follows:
$$C(O_F^+)\simeq C(SU^q_2)$$

(3) Now back to the general case, let us try to understand the integration over $O_F^+$. Given a noncrossing pairing $\pi\in NC_2(2k)$ and an index $i=(i_1,\ldots,i_{2k})$, we set:
$$\delta_\pi^F(i)=\prod_{s\in\pi}F_{i_{s_l}i_{s_r}}$$

Here the product is over all strings $s=\{s_l\curvearrowright s_r\}$ of $\pi$. Our claim is that the following family of vectors, with $\pi\in NC_2(2k)$, spans the space of fixed vectors of $u^{\otimes 2k}$:
$$\xi_\pi=\sum_i\delta_\pi^F(i)e_{i_1}\otimes\ldots\otimes e_{i_{2k}}$$ 

Indeed, having $\xi_\cap$ fixed by $u^{\otimes 2}$ is equivalent to assuming that $u=F\bar{u}F$ is unitary. By using now the above vectors, we obtain the following Weingarten formula:
$$\int_{O_F^+}u_{i_1j_1}\ldots u_{i_{2k}j_{2k}}=\sum_{\pi\sigma}\delta_\pi^F(i)\delta_\sigma^F(j)W_{kN}(\pi,\sigma)$$

(4) With these preliminaries in hand, let us start the computation. Let $q\in [-1,0)$ be such that $q+q^{-1}=-N$. Our claim is that we have the following formula:
$$\int_{O_N^+}\varphi(\sqrt{N+2}\,u_{ij})=\int_{SU^q_2}\varphi(\alpha+\alpha^*+\gamma-q\gamma^*)$$

Indeed, the moments of the variable on the left are given by:
$$\int_{O_N^+}u_{ij}^{2k}=\sum_{\pi\sigma}W_{kN}(\pi,\sigma)$$

On the other hand, the moments of the variable on the right, which in terms of the fundamental corepresentation $v=(v_{ij})$ is given by $w=\sum_{ij}v_{ij}$, are given by:
$$\int_{SU^q_2}w^{2k}=\sum_{ij}\sum_{\pi\sigma}\delta_\pi^F(i)\delta_\sigma^F(j)W_{kN}(\pi,\sigma)$$

We deduce that $w/\sqrt{N+2}$ has the same moments as $u_{ij}$, which proves our claim. 

\medskip

(5) In order to do now the computation over $SU^q_2$, we can use a matrix model due to Woronowicz \cite{wo1}, where the standard generators $\alpha,\gamma$ are mapped as follows:
$$\pi_u(\alpha)e_k=\sqrt{1-q^{2k}}e_{k-1}\quad,\quad 
\pi_u(\gamma)e_k=uq^k e_k$$

Here $u\in\mathbb T$ is a parameter, and $(e_k)$ is the standard basis of $l^2(\mathbb N)$. The point with this representation is that it allows the computation of the Haar functional. Indeed, if $D$ is the diagonal operator given by $D(e_k)=q^{2k}e_k$, then the formula is as follows:
$$\int _{SU^q_2}x=(1-q^2)\int_{\mathbb T}tr(D\pi_u(x))\frac{du}{2\pi iu}$$

(6) Thus, the law of the variable that we are interested in is of the following form:
$$\int_{SU^q_2}\varphi(\alpha+\alpha^*+\gamma-q\gamma^*)=(1-q^2)\int_{\mathbb T}tr(D\varphi(M))\frac{du}{2\pi iu}$$

To be more precise, this formula holds indeed, with:
$$M(e_k)=e_{k+1}+q^k(u-qu^{-1})e_k+(1-q^{2k})e_{k-1}$$

The point now is that the integral on the right can be computed, by using advanced calculus methods, and this gives the result. We refer here to \cite{bcz}. 
\end{proof}

For more regarding the above, including open problems, we refer to \cite{bbs}, \cite{bcz}.

\section*{4e. Exercises} 

Things have been quite technical in this chapter, especially towards the end, and our exercises here will focus on the first part, which is more elementary. First, we have: 

\begin{exercise}
Write down a complete, elementary proof of
$$S_{M_2}^+=SO_3$$
and then unify this with the other such result that we have, namely $S_3^+=S_3$.
\end{exercise}

This is something quite tricky, the key word here being ``elementary''. Indeed, as we have seen in the above, all this can be understood via representation theory.

\begin{exercise}
Come up with explicit constructions for the subgroups of
$$S_4^+=SO_3^{-1}$$
and explain what the associated $ADE$ diagrams should be.
\end{exercise}

This is again an instructive exercise, potentially leading you into a lot of interesting mathematics. You will have to learn here first, as a preliminary exercise, more about $SU_2,SO_3$ and their subgroups, with the theory here going under different names, with the keywords being ``ADE'' and ``McKay''. Then, get into the quantum group case, and with a look into subfactors too, where there are many interesting things to be learned. In what concerns us, we will be back to this, but only later in this book. 

\begin{exercise}
Perform a study of the group dual subgroups
$$\widehat{\Gamma}\subset S_Z^+$$
in the homogeneous space case, where $Z=M_K\times\{1,\ldots,L\}$.
\end{exercise}

We will be back to such subgroup questions later on, when analyzing the quantum subgroups $G\subset S_Z^+$ coming from the finite quantum graphs with vertex set $Z$.

\begin{exercise}
Meditate about the lack of relationship between $S_N,O_N$, as opposed to the relationship between $S_N^+,O_N^+$, explained above. Would you consider, in view of this, that the free world is simpler than the classical one? 
\end{exercise}

This looks more like a philosophy exercise, but philosophy is important too. You need good motivations for staying up late doing math computations, and personally, my belief that the answer to the above question is ``yes'' has always motivated me.

\part{Quantum reflections}

\ \vskip50mm

\begin{center}
{\em I'm gonna make me a good sharp axe

Shining steel tempered in the fire

I'll chop you down like an old dead tree

Dirty old town, dirty old town}
\end{center}

\chapter{Finite graphs}

\section*{5a. Finite graphs}

In the classical case, many interesting permutation groups $G\subset S_N$ appear as symmetry groups $G(X)$ of the graphs $X$ having $N$ vertices. In analogy with this, many interesting examples of quantum permutation groups $G\subset S_N^+$ appear as particular cases of the following general construction from \cite{ba3}, involving finite graphs:

\index{graph symmetry}
\index{graph automorphism}

\begin{proposition}
Given a finite graph $X$, with adjacency matrix $d\in M_N(0,1)$, the following construction produces a quantum permutation group, 
$$C(G^+(X))=C(S_N^+)\Big/\Big<du=ud\Big>$$
whose classical version $G(X)$ is the usual  automorphism group of $X$.
\end{proposition}

\begin{proof}
The fact that we have a quantum group comes from the fact that $du=ud$ reformulates as $d\in End(u)$, which makes it clear that we are dividing by a Hopf ideal. Regarding the second assertion, we must establish here the following equality:
$$C(G(X))=C(S_N)\Big/\Big<du=ud\Big>$$

For this purpose, recall that we have $u_{ij}(\sigma)=\delta_{\sigma(j)i}$. We therefore obtain:
$$(du)_{ij}(\sigma)
=\sum_kd_{ik}u_{kj}(\sigma)
=\sum_kd_{ik}\delta_{\sigma(j)k}
=d_{i\sigma(j)}$$
$$(ud)_{ij}(\sigma)
=\sum_ku_{ik}(\sigma)d_{kj}
=\sum_k\delta_{\sigma(k)i}d_{kj}
=d_{\sigma^{-1}(i)j}$$

Thus $du=ud$ reformulates as $d_{ij}=d_{\sigma(i)\sigma(j)}$, which gives the result.
\end{proof}

Let us work out some basic examples. We have the following result:

\index{simplex}
\index{square graph}
\index{disconnected union}

\begin{proposition}
The construction $X\to G^+(X)$ has the following properties:
\begin{enumerate}
\item For the $N$-point graph, having no edges at all, we obtain $S_N^+$.

\item For the $N$-simplex, having edges everywhere, we obtain as well $S_N^+$.

\item We have $G^+(X)=G^+(X^c)$, where $X^c$ is the complementary graph.

\item For a disconnected union, we have $G^+(X)\,\hat{*}\,G^+(Y)\subset G^+(X\sqcup Y)$.

\item For the square, we obtain a non-classical, proper subgroup of $S_4^+$.
\end{enumerate}
\end{proposition}

\begin{proof}
All these results are elementary, the proofs being as follows:

\medskip

(1) This follows from definitions, because here we have $d=0$.

\medskip

(2) Here $d=\mathbb I-1$, where $\mathbb I$ is the all-one matrix, and the magic condition gives $u\mathbb I=\mathbb Iu=N\mathbb I$. We conclude that $du=ud$ is automatic, and so $G^+(X)=S_N^+$.

\medskip

(3) The adjacency matrices of $X,X^c$ being related by the following formula: $$d_X+d_{X^c}=\mathbb I-1$$

By using now the above formula $u\mathbb I=\mathbb Iu=N\mathbb I$, we conclude that $d_Xu=ud_X$ is equivalent to $d_{X^c}u=ud_{X^c}$. Thus, we obtain, as claimed, $G^+(X)=G^+(X^c)$.

\medskip

(4) The adjacency matrix of a disconnected union is given by:
$$d_{X\sqcup Y}=diag(d_X,d_Y)$$

Now let $w=diag(u,v)$ be the fundamental corepresentation of $G^+(X)\,\hat{*}\,G^+(Y)$. Then $d_Xu=ud_X$ and $d_Yv=vd_Y$, and we obtain, as desired, $d_{X\sqcup Y}w=wd_{X\sqcup Y}$.

\medskip

(5) We know from (3) that we have $G^+(\square)=G^+(|\ |)$. We know as well from (4) that we have $\mathbb Z_2\,\hat{*}\,\mathbb Z_2\subset G^+(|\ |)$. It follows that $G^+(\square)$ is non-classical. Finally, the inclusion $G^+(\square)\subset S_4^+$ is indeed proper, because $S_4\subset S_4^+$ does not act on the square.
\end{proof}

In order to further advance, and to explicitely compute various quantum automorphism groups, we can use the spectral decomposition of $d$, as follows:

\index{spectral decomposition}

\begin{theorem}
A closed subgroup $G\subset S_N^+$ acts on a graph $X$ precisely when
$$P_\lambda u=uP_\lambda\quad,\quad\forall\lambda\in\mathbb R$$
where $d=\sum_\lambda\lambda\cdot P_\lambda$ is the spectral decomposition of the adjacency matrix of $X$.
\end{theorem}

\begin{proof}
Since $d\in M_N(0,1)$ is a symmetric matrix, we can consider indeed its spectral decomposition, $d=\sum_\lambda\lambda\cdot P_\lambda$, and we have the following formula:
$$<d>=span\left\{P_\lambda\Big|\lambda\in\mathbb R\right\}$$

Thus $d\in End(u)$ when $P_\lambda\in End(u)$ for all $\lambda\in\mathbb R$, which gives the result.
\end{proof}

In order to exploit Theorem 5.3, we will often combine it with the following fact:

\begin{proposition}
Given a closed subgroup $G\subset S_N^+$, with associated coaction
$$\Phi:\mathbb C^N\to \mathbb C^N\otimes C(G)\quad,\quad e_i\to\sum_je_j\otimes u_{ji}$$
and a linear subspace $V\subset\mathbb C^N$, the following are equivalent:
\begin{enumerate}
\item The magic matrix $u=(u_{ij})$ commutes with $P_V$.

\item $V$ is invariant, in the sense that $\Phi(V)\subset V\otimes C(G)$.
\end{enumerate}
\end{proposition}

\begin{proof}
Let $P=P_V$. For any $i\in\{1,\ldots,N\}$ we have the following formula:
$$\Phi(P(e_i))
=\Phi\left(\sum_kP_{ki}e_k\right) 
=\sum_{jk}P_{ki}e_j\otimes u_{jk}
=\sum_je_j\otimes (uP)_{ji}$$

On the other hand the linear map $(P\otimes id)\Phi$ is given by a similar formula:
$$(P\otimes id)(\Phi(e_i))
=\sum_kP(e_k)\otimes u_{ki}
=\sum_{jk}P_{jk}e_j\otimes u_{ki}
=\sum_je_j\otimes (Pu)_{ji}$$

Thus $uP=Pu$ is equivalent to $\Phi P=(P\otimes id)\Phi$, and the conclusion follows.
\end{proof}

As an application of the above results, we have the following computation:

\index{square graph}
\index{cycle graph}
\index{circulant graph}

\begin{theorem}
The quantum automorphism group of the $N$-cycle is, at $N\neq4$:
$$G^+(X)=D_N$$
However, at $N=4$ we have $D_4\subset G^+(X)\subset S_4^+$, with proper inclusions.
\end{theorem}

\begin{proof}
We know from Proposition 5.2, and from $S_N=S_N^+$ at $N\leq3$, that the various assertions hold indeed at $N\leq4$. So, assume $N\geq5$. Given a $N$-th root of unity, $w^N=1$, the vector $\xi=(w^i)$ is an eigenvector of $d$, with eigenvalue as follows:
$$\lambda=w+w^{N-1}$$

Now by taking $w=e^{2\pi i/N}$, it follows that the are eigenvectors of $d$ are:
$$1,f,f^2,\ldots ,f^{N-1}$$

More precisely, the invariant subspaces of $d$ are as follows, with the last subspace having dimension 1 or 2 depending on the parity of $N$:
$$\mathbb C 1,\, \mathbb C f\oplus\mathbb C f^{N-1},\, \mathbb C f^2\oplus\mathbb C f^{N-2},\ldots$$

Assuming $G\subset G^+(X)$, consider the coaction $\Phi:\mathbb C^N\to \mathbb C^N\otimes C(G)$, and write:
$$\Phi(f)=f\otimes a+f^{N-1}\otimes b$$

By taking the square of this equality we obtain:
$$\Phi(f^2)=f^2\otimes a^2+f^{N-2}\otimes b^2+1\otimes(ab+ba)$$

It follows that $ab=-ba$, and that $\Phi(f^2)$ is given by the following formula:
$$\Phi(f^2)=f^2\otimes a^2+f^{N-2}\otimes b^2$$

By multiplying this with $\Phi(f)$ we obtain the following formula:
$$\Phi(f^3)=f^3\otimes a^3+f^{N-3}\otimes b^3+f^{N-1}\otimes ab^2+f\otimes ba^2$$

Now since $N\geq 5$ implies that $1,N-1$ are different from $3,N-3$, we must have $ab^2=ba^2=0$. By using this and $ab=-ba$, we obtain by recurrence on $k$ that:
$$\Phi(f^k)=f^k\otimes a^k+f^{N-k}\otimes b^k$$

In particular at $k=N-1$ we obtain the following formula:
$$\Phi(f^{N-1})=f^{N-1}\otimes a^{N-1}+f\otimes b^{N-1}$$

On the other hand we have $f^*=f^{N-1}$, so by applying $*$ to $\Phi(f)$ we get:
$$\Phi(f^{N-1})=f^{N-1}\otimes a^*+f\otimes b^*$$

Thus $a^*=a^{N-1}$ and $b^*=b^{N-1}$. Together with $ab^2=0$ this gives:
$$(ab)(ab)^*
=abb^*a^*
=ab^Na^{N-1}
=(ab^2)b^{N-2}a^{N-1}
=0$$

From positivity we get from this $ab=0$, and together with $ab=-ba$, this shows that $a,b$ commute. On the other hand $C(G)$ is generated by the coefficients of $\Phi$, which are powers of $a,b$, and so $C(G)$ must be commutative, and we obtain the result.
\end{proof}

The above result is quite suprising, but we will be back to this, with a more conceptual explanation for the fact that the square $\square$ has quantum symmetry. Back to theory now, we have the following useful result from \cite{ba3}, complementary to Theorem 5.3:

\index{color decomposition}

\begin{theorem}
Given a matrix $p\in M_N(\mathbb C)$, consider its ``color'' decomposition
$$p=\sum_{c\in\mathbb C}c\cdot p_c$$
with the color components $p_c\in M_N(0,1)$ with $c\in\mathbb C$ being constructed as follows:
$$(p_c)_{ij}=\begin{cases}
1&{\rm if}\ p_{ij}=c\\
0&{\rm otherwise}
\end{cases}$$
Then a magic matrix $u=(u_{ij})$ commutes with $p$ iff it commutes with all matrices $p_c$.
\end{theorem}

\begin{proof}
Consider the multiplication and counit maps of the algebra $\mathbb C^N$:
$$M:e_i\otimes e_j\to e_ie_j\quad,\quad 
C:e_i\to e_i\otimes e_i$$

Since $M,C$ intertwine $u,u^{\otimes 2}$, their iterations $M^{(k)},C^{(k)}$ intertwine $u,u^{\otimes k}$, and so:
$$M^{(k)}p^{\otimes k}C^{(k)}
=\sum_{c\in\mathbb C}c^kp_c
\in End(u)$$

Now since this formula holds for any $k\in\mathbb N$, we obtain the result.
\end{proof}

The above general results can be combined, and we are led to the following statement:

\index{spectral-color decomposition}
\index{color-spectral decomposition}

\begin{theorem}
A closed subgroup $G\subset S_N^+$ acts on a graph $X$ precisely when 
$$u=(u_{ij})$$
commutes with all the matrices coming from the color-spectral decomposition of $d$.
\end{theorem}

\begin{proof}
This follows by combining Theorem 5.3 and Theorem 5.6, with the ``color-spectral'' decomposition in the statement referring to what comes out by succesively doing the color and spectral decomposition, until the process stabilizes.
\end{proof}

All this might seem in need of some further discussion. In answer to this, the point is that we are in fact doing planar algebras. Following \cite{ba3}, we have indeed:

\index{planar algebra}
\index{spin planar algebra}
\index{two-box}

\begin{theorem}
The planar algebra associated to $G^+(X)$ is equal to the planar algebra generated by $d$, viewed as a $2$-box in the spin planar algebra $\mathcal S_N$, with $N=|X|$.
\end{theorem}

\begin{proof}
We recall from chapter 3 that any quantum permutation group $G\subset S_N^+$ produces a subalgebra $P\subset\mathcal S_N$ of the spin planar algebra, given by:
$$P_k=Fix(u^{\otimes k})$$

In our case, the idea is that $G=G^+(X)$ comes via the relation $d\in End(u)$, but we can view this relation, via Frobenius duality, as a relation of the following type:
$$\xi_d\in Fix(u^{\otimes 2})$$

Indeed, let us view the adjacency matrix $d\in M_N(0,1)$ as a 2-box in $\mathcal S_N$, by using the canonical identification between $M_N(\mathbb C)$ and the algebra of 2-boxes $\mathcal S_N(2)$:
$$(d_{ij})\leftrightarrow \sum_{ij} d_{ij}\begin{pmatrix}i&i\\ j&j\end{pmatrix}$$

Let $P$ be the planar algebra associated to $G^+(X)$ and let $Q$ be the planar algebra generated by $d$. The action of $u^{\otimes 2}$ on $d$ viewed as a 2-box is given by:
$$u^{\otimes 2}\left(\sum_{ij} d_{ij}\begin{pmatrix}i&i\\ j&j\end{pmatrix}\right)
=\sum_{ijkl} d_{ij}\begin{pmatrix}k&k\\ l&l\end{pmatrix}\otimes u_{ki}u_{lj}
=\sum_{kl}\begin{pmatrix}k&k\\ l&l\end{pmatrix}\otimes (udu^t)_{kl}$$

Since $v$ is a magic unitary commuting with $d$ we have:
$$udu^t=duu^t=d$$

But this means that $d$, viewed as a 2-box, is in the algebra $P_2$ of fixed points of $u^{\otimes 2}$. Thus $Q\subset P$. As for $P\subset Q$, this follows from the duality found in chapter 3.
\end{proof}

As a conclusion to all this, the construction $X\to G^+(X)$ is something quite interesting. In analogy with the usual construction $X\to G(X)$, this often leads to the spectral theory of $X$. And also, there are some interesting planar algebra aspects.

\section*{5b. The hypercube}

Let us go back to the square $\square$. This is naturally part of the series of $N$-cycles, but Theorem 5.5 shows that, within this series, $\square$ is an exceptional object. With the reason for this coming somehow from the fact that the complement $\square^c=|\ |$ is not connected.

\bigskip

We will discuss later this phenomenon, with a systematic study of the graphs of type $|\ |$, which appear as products. Before that, however, let us try to compute $G^+(\square)$, with the tools that we have. Quite remarkably, the result here is as follows:

\index{hypercube}
\index{Cayley graph}

\begin{theorem}
The quantum symmetry group of the $N$-hypercube is
$$G^+(\square_N)=O_N^{-1}$$
with the corresponding coaction map on the vertex set being given by
$$\Phi:C^*(\mathbb Z_2^N)\to C^*(\mathbb Z_2^N)\otimes C(O_N^{-1})\quad,\quad 
g_i\to\sum_jg_j\otimes u_{ji}$$
via the standard identification $\square_N=\widehat{\mathbb Z_2^N}$. In particular we have $G^+(\square)=O_2^{-1}$.
\end{theorem} 

\begin{proof}
This result is from \cite{bbc}, with its $N=2$ particular case, corresponding to the last assertion, going back to \cite{bi1}, the idea being as follows:

\medskip

(1) Our first claim is that $\square_N$ is the Cayley graph of $\mathbb Z_2^N=<\tau_1,\ldots ,\tau_N>$. Indeed, the vertices of this latter Cayley graph are the products of the following form:
$$g=\tau_1^{i_1}\ldots\tau_N^{i_N}$$

The sequence of 0-1 exponents defining such an element determines a point of $\mathbb R^N$, which is a vertex of the cube. Thus the vertices of the Cayley graph are the vertices of the cube. Now regarding the edges, in the Cayley graph these are drawn between elements $g,h$ having the property $g=h\tau_i$ for some $i$. In terms of coordinates, the operation $h\to h\tau_i$ means to switch the sign of the $i$-th coordinate, and to keep the other coordinates fixed. In other words, we get in this way the edges of the cube, as desired.

\medskip

(2) Our second claim is that, when identifying the vector space spanned by the vertices of $\square_N$ with the algebra $C^*(\mathbb Z_2^N)$, the eigenvectors and eigenvalues of $\square_N$ are given by:
$$v_{i_1\ldots i_N}=\sum_{j_1\ldots j_N} (-1)^{i_1j_1
+\ldots+i_Nj_N}\tau_1^{j_1}\ldots\tau_N^{j_N}$$
$$\lambda_{i_1\ldots i_N}=(-1)^{i_1}+\ldots +(-1)^{i_N}$$

Indeed, let us recall that the action of $d$ on the functions on vertices is given by the following formula, with $q\sim p$ standing for the fact that $q,p$ are joined by an edge:
$$df(p)=\sum_{q\sim p}f(q)$$

Now by identifying the vertices with the the elements of $\mathbb Z_2^N$, hence the functions on the vertices with the elements of the algebra
$C^*(\mathbb Z_2^N)$, we get the following formula:
$$dv=\tau_1v+\ldots +\tau_Nv$$

With $v_{i_1\ldots i_N}$ as above, we have the following computation, which proves our claim:
\begin{eqnarray*}
dv_{i_1\ldots i_N}
&=&\sum_s\tau_s\sum_{j_1\ldots j_N}(-1)^{i_1j_1+\ldots+i_Nj_N}\tau_1^{j_1}\ldots\tau_N^{j_N}\\
&=&\sum_s\sum_{j_1\ldots j_N}(-1)^{i_1j_1+\ldots+i_Nj_N}\tau_1^{j_1}\ldots\tau_s^{j_s+1}\ldots\tau_N^{j_N}\\
&=&\sum_s\sum_{j_1\ldots j_N}(-1)^{i_s}(-1)^{i_1j_1+\ldots+i_Nj_N}\tau_1^{j_1}\ldots\tau_s^{j_s}\ldots\tau_N^{j_N}\\
&=&\sum_s(-1)^{i_s}\sum_{j_1\ldots j_N}(-1)^{i_s}(-1)^{i_1j_1+\ldots+i_Nj_N}\tau_1^{j_1}\ldots\tau_N^{j_N}\\
&=&\lambda_{i_1\ldots i_N}v_{i_1\ldots i_N}
\end{eqnarray*}

(3) We prove now that the quantum group $O_N^{-1}$ acts on the cube $\square_N$. For this purpose, observe first that we have a map as follows:
$$\Phi :C^*(\mathbb Z_2^N)\to C^*(\mathbb Z_2^N)\otimes C(O_N^{-1})
\quad,\quad\tau_i\to\sum_j\tau_j\otimes u_{ji}$$

It is routine to check that for $i_1\neq i_2\neq\ldots\neq i_l$ we have:
$$\Phi(\tau_{i_1}\ldots\tau_{i_l})=\sum_{j_1\neq\ldots\neq j_l}\tau_{j_1}\ldots\tau_{j_l} \otimes u_{j_1i_1}\ldots u_{j_li_l}$$

In terms of eigenspaces $E_s$ of the adjacency matrix, this gives:
$$\Phi(E_s)\subset E_s\otimes C(O_N^{-1})$$

Thus $\Phi$ preserves the adjacency matrix of $\square_N$, so is a coaction on $\square_N$, as claimed.

\medskip

(4) Conversely now, consider the universal coaction on the cube:
$$\Psi:C^*(\mathbb Z_2^N)\to C^*(\mathbb Z_2^N)\otimes C(G)\quad,\quad\tau_i\to\sum_j\tau_j\otimes u_{ji}$$

By applying $\Psi$ to the relation $\tau_i\tau_j=\tau_j\tau_i$ we get $u^tu=1$, so the matrix $u=(u_{ij})$ is orthogonal. By applying $\Psi$ to the relation $\tau_i^2=1$ we get:
$$1\otimes\sum_ku_{ki}^2+\sum_{k<l}\tau_k\tau_l\otimes(u_{ki}u_{li}+u_{li}u_{ki})=1\otimes 1$$

This gives $u_{ki}u_{li}=-u_{li}u_{ki}$ for $i\neq j$, $k\neq l$, and by using the antipode we get $u_{ik}u_{il}=-u_{il}u_{ik}$ for $k\neq l$. Also, by applying $\Psi$ to $\tau_i\tau_j=\tau_j\tau_i$ with $i\neq j$ we get:
$$\sum_{k<l}\tau_k\tau_l\otimes(u_{ki}u_{lj}+u_{li}u_{kj})=\sum_{k<l}\tau_k\tau_l\otimes
(u_{kj}u_{li}+u_{lj}u_{ki})$$

It follows that for $i\neq j$ and $k\neq l$, we have:
$$u_{ki}u_{lj}+u_{li}u_{kj}=u_{kj}u_{li}+u_{lj}u_{ki}$$

In other words, we have $[u_{ki},u_{lj}]=[u_{kj},u_{li}]$. By using the antipode we get $[u_{jl},u_{ik}]=[u_{il},u_{jk}]$. Now by combining these relations we get:
$$[u_{il},u_{jk}]=[u_{ik},u_{jl}]=[u_{jk},u_{il}]=-[u_{il},u_{jk}]$$

Thus $[u_{il},u_{jk}]=0$, so the elements $u_{ij}$ satisfy the relations for $C(O_N^{-1})$, as desired.
\end{proof}

Our purpose now is to understand which representation of $O_N$ produces by twisting the magic representation of $O_N^{-1}$. In order to solve this question, we will need:

\begin{proposition}
The Fourier transform over $\mathbb Z_2^N$ is the map
$$\alpha:C(\mathbb Z_2^N)\to C^*(\mathbb Z_2^N)\quad,\quad 
\delta_{g_1^{i_1}\ldots g_N^{i_N}}\to\frac{1}{2^N}\sum_{j_1\ldots j_N}(-1)^{<i,j>}g_1^{j_1}\ldots g_N^{j_N}$$
with the usual convention $<i,j>=\sum_ki_kj_k$, and its inverse is the map
$$\beta:C^*(\mathbb Z_2^N)\to C(\mathbb Z_2^N)\quad,\quad 
g_1^{i_1}\ldots g_N^{i_N}\to\sum_{j_1\ldots j_N}(-1)^{<i,j>}\delta_{g_1^{j_1}\ldots g_N^{j_N}}$$
with all the exponents being binary, $i_1,\ldots,i_N,j_1,\ldots,j_N\in\{0,1\}$.
\end{proposition}

\begin{proof}
Observe first that the group $\mathbb Z_2^N$ can be written as follows:
$$\mathbb Z_2^N=\left\{g_1^{i_1}\ldots g_N^{i_N}\Big|i_1,\ldots,i_N\in\{0,1\}\right\}$$

Thus both $\alpha,\beta$ are well-defined, and it is elementary to check that both are morphisms of algebras. Also, we have $\alpha\beta=\beta\alpha=id$, coming from the following standard formula:
$$\frac{1}{2^N}\sum_{j_1\ldots j_N}(-1)^{<i,j>}
=\prod_{k=1}^N\left(\frac{1}{2}\sum_{j_r}(-1)^{i_rj_r}\right)
=\delta_{i0}$$

Thus we have indeed a pair of inverse Fourier morphisms, as claimed.
\end{proof}

By using now these Fourier transforms, we obtain following result:

\begin{proposition}
The magic unitary for the embedding $O_N^{-1}\subset S_{2^N}^+$ is given by
$$w_{i_1\ldots i_N,k_1\ldots k_N}=\frac{1}{2^N}\sum_{j_1\ldots j_N}\sum_{b_1\ldots b_N}(-1)^{<i+k_b,j>}\left(\frac{1}{N}\right)^{\#(0\in j)}u_{1b_1}^{j_1}\ldots u_{Nb_N}^{j_N}$$
where $k_b=(k_{b_1},\ldots,k_{b_N})$, with respect to multi-indices $i,k\in\{0,1\}^N$ as above.
\end{proposition}

\begin{proof}
By composing the coaction map $\Phi$ from Theorem 5.9 with the above Fourier transform isomorphisms $\alpha,\beta$, we have a diagram as follows:
$$\xymatrix@R=20mm@C=25mm{
C^*(\mathbb Z_2^N)\ar[r]^\Phi&C^*(\mathbb Z_2^N)\otimes C(O_N^{-1}) \ar[d]^{\beta\otimes id}\\
C(\mathbb Z_2^N)\ar[u]^\alpha\ar@.[r]^\Psi&C(\mathbb Z_2^N)\otimes C(O_N^{-1})}$$

In order to compute the composition on the bottom $\Psi$, we first recall from Theorem 5.9 that the coaction map $\Phi$ is defined by the following formula:
$$\Phi(g_b)=\sum_ag_a\otimes u_{ab}$$

Now by making products of such quantities, we obtain the following global formula for $\Phi$, valid for any exponents $i_1,\ldots,i_N\in\{1,\ldots,N\}$:
$$\Phi(g_1^{i_1}\ldots g_N^{i_N})=\left(\frac{1}{N}\right)^{\#(0\in i)}\sum_{b_1\ldots b_N}g_{b_1}^{i_1}\ldots g_{b_N}^{i_N}\otimes u_{1b_1}^{i_1}\ldots u_{Nb_N}^{i_N}$$

The term on the right can be put in ``standard form'' as follows:
$$g_{b_1}^{i_1}\ldots g_{b_N}^{i_N}=g_1^{\sum_{b_x=1}i_x}\ldots g_N^{\sum_{b_x}i_x}$$

We therefore obtain the following formula for the coaction map $\Phi$:
$$\Phi(g_1^{i_1}\ldots g_N^{i_N})=\left(\frac{1}{N}\right)^{\#(0\in i)}
\sum_{b_1\ldots b_N}g_1^{\sum_{b_x=1}i_x}\ldots g_N^{\sum_{b_x=N}i_x}\otimes u_{1b_1}^{i_1}\ldots u_{Nb_N}^{i_N}$$

Now by applying the Fourier transforms, we obtain the following formula:
\begin{eqnarray*}
&&\Psi(\delta_{g_1^{i_1}\ldots g_N^{i_N}})\\
&=&(\beta\otimes id)\Phi\left(\frac{1}{2^N}\sum_{j_1\ldots j_N}(-1)^{<i,j>}g_1^{j_1}\ldots g_N^{j_N}\right)\\
&=&\frac{1}{2^N}\sum_{j_1\ldots j_N}\sum_{b_1\ldots b_N}(-1)^{<i,j>}
\left(\frac{1}{N}\right)^{\#(0\in j)}
\beta\left( g_1^{\sum_{b_x=1}j_x}\ldots g_N^{\sum_{b_x=N}j_x}\right)\otimes u_{1b_1}^{j_1}\ldots u_{Nb_N}^{j_N}
\end{eqnarray*}

By using now the formula of $\beta$ from Proposition 5.10, we obtain:
\begin{eqnarray*}
\Psi(\delta_{g_1^{i_1}\ldots g_N^{i_N}})
&=&\frac{1}{2^N}\sum_{j_1\ldots j_N}\sum_{b_1\ldots b_N}\sum_{k_1\ldots k_N}\left(\frac{1}{N}\right)^{\#(0\in j)}\\
&&(-1)^{<i,j>}(-1)^{<(\sum_{b_x=1}j_x,\ldots,\sum_{b_x=N}j_x),(k_1,\ldots,k_N)>}\\
&&\delta_{g_1^{k_1}\ldots g_N^{k_N}}\otimes u_{1b_1}^{j_1}\ldots u_{Nb_N}^{j_N}
\end{eqnarray*}

Now observe that, with the notation $k_b=(k_{b_1},\ldots,k_{b_N})$, we have:
$$\left<\left(\sum_{b_x=1}j_x,\ldots,\sum_{b_x=N}j_x\right),(k_1,\ldots,k_N)\right>=<j,k_b>$$

Thus, we obtain the following formula for our map $\Psi$:
$$\Psi(\delta_{g_1^{i_1}\ldots g_N^{i_N}})
=\frac{1}{2^N}\sum_{j_1\ldots j_N}\sum_{b_1\ldots b_N}\sum_{k_1\ldots k_N}(-1)^{<i+k_b,j>}\left(\frac{1}{N}\right)^{\#(0\in j)}\delta_{g_1^{k_1}\ldots g_N^{k_N}}\otimes u_{1b_1}^{j_1}\ldots u_{Nb_N}^{j_N}$$

But this gives the formula in the statement for the corresponding magic unitary, with respect to the basis $\{\delta_{g_1^{i_1}\ldots g_N^{i_N}}\}$ of the algebra $C(\mathbb Z_2^N)$, and we are done.
\end{proof}

We can now solve our original question, namely understanding where the magic representation of $O_N^{-1}$ really comes from, with the following final answer to it:

\index{antisymmetric representation}
\index{twisting}

\begin{theorem}
The magic representation of $O_N^{-1}$, coming from its action on the $N$-cube, corresponds to the antisymmetric representation of $O_N$, via twisting.
\end{theorem}

\begin{proof}
This follows from the formula of $w$ in Proposition 5.11, by computing the character, and then interpreting the result via twisting, as follows:

\medskip

(1) By applying the trace to the formula of $w$, we obtain:
$$\chi
=\sum_{j_1\ldots j_N}\sum_{b_1\ldots b_N}\left(\frac{1}{2^N}\sum_{i_1\ldots i_N}(-1)^{<i+i_b,j>}\right)\left(\frac{1}{N}\right)^{\#(0\in j)}u_{1b_1}^{j_1}\ldots u_{Nb_N}^{j_N}$$

(2) By computing the Fourier sum in the middle, we are led to the following formula, with binary indices $j_1,\ldots,j_N\in\{0,1\}$, and plain indices $b_1,\ldots,b_N\in\{1,\ldots,N\}$:
$$\chi=\sum_{j_1\ldots j_N}\sum_{b_1\ldots b_N}\left(\frac{1}{N}\right)^{\#(0\in j)}\delta_{j_1,\sum_{b_x=1}j_x}\ldots\delta_{j_N,\sum_{b_x=N}j_x}u_{1b_1}^{j_1}\ldots u_{Nb_N}^{j_N}$$

(3) With the notation $r=\#(1\in j)$ we obtain a decomposition as follows:
$$\chi=\sum_{r=0}^N\chi_r$$

To be more precise, the variables $\chi_r$ are as follows:
$$\chi_r=\frac{1}{N^{N-r}}\sum_{\#(1\in j)=r}\sum_{b_1\ldots b_N}\delta_{j_1,\sum_{b_x=1}j_x}\ldots\delta_{j_N,\sum_{b_x=N}j_x}u_{1b_1}^{j_1}\ldots u_{Nb_N}^{j_N}$$

(4) Consider now the set $A\subset\{1,\ldots,N\}$ given by:
$$A=\{a|j_a=1\}$$

The binary multi-indices $j\in\{0,1\}^N$ satisfying $\#(1\in j)=r$ being in bijection with such subsets $A$, satisfying $|A|=r$, we can replace the sum over $j$ with a sum over such subsets $A$. We obtain a formula as follows, where $j$ is the index corresponding to $A$:
$$\chi_r=\frac{1}{N^{N-r}}\sum_{|A|=r}\sum_{b_1\ldots b_N}\delta_{j_1,\sum_{b_x=1}j_x}\ldots\delta_{j_N,\sum_{b_x=N}j_x}\prod_{a\in A}u_{ab_a}$$

(5) Let us identify $b$ with the corresponding function $b:\{1,\ldots,N\}\to\{1,\ldots,N\}$, via $b(a)=b_a$. Then for any $p\in\{1,\ldots,N\}$ we have:
$$\delta_{j_p,\sum_{b_x=p}j_x}=1
\iff |b^{-1}(p)\cap A|=\chi_A(p)\ ({\rm mod}\ 2)$$

We conclude that the multi-indices $b\in\{1,\ldots,N\}^N$ which effectively contribute to the sum are those coming from the functions satisfying $b<A$. Thus, we have:
$$\chi_r=\frac{1}{N^{N-r}}\sum_{|A|=r}\sum_{b<A}\prod_{a\in A}u_{ab_a}$$

(6) We can further split each $\chi_r$ over the sets $A\subset\{1,\ldots,N\}$ satisfying $|A|=r$. The point is that for each of these sets we have:
$$\frac{1}{N^{N-r}}\sum_{b<A}\prod_{a\in A}u_{ab_a}=\sum_{\sigma\in S_N^A}\prod_{a\in A}u_{a\sigma(a)}$$

Thus, the magic character of $O_N^{-1}$ splits as $\chi=\sum_{r=0}^N\chi_r$, the components being:
$$\chi_r=\sum_{|A|=r}\sum_{\sigma\in S_N^A}\prod_{a\in A}u_{a\sigma(a)}$$

(7) The twisting operation $O_N\to O_N^{-1}$ makes correspond the following products:
$$\varepsilon(\sigma)\prod_{a\in A}u_{a\sigma(a)}\to\prod_{a\in A}u_{a\sigma(a)}$$

Now by summing over sets $A$ and permutations $\sigma$, we conclude that the twisting operation $O_N\to O_N^{-1}$ makes correspond the following quantities:
$$\sum_{|A|=r}\sum_{\sigma\in S_N^A}\varepsilon(\sigma)\prod_{a\in A}u_{a\sigma(a)}\to \sum_{|A|=r}\sum_{\sigma\in S_N^A}\prod_{a\in A}u_{a\sigma(a)}$$

Thus, we are led to the conclusion in the statement.
\end{proof}

\section*{5c. Product operations}

We discuss now the behavior of the operation $X\to G^+(X)$ under taking various products of graphs. We use the notation $X=(X,\sim)$, where the $X$ on the right is the set of vertices, and where we write $i\sim j$ when $i,j$ are connected by an edge. We have:

\begin{definition}
Let $X,Y$ be two finite graphs.
\begin{enumerate}
\item The direct product $X\times Y$ has vertex set $X\times Y$, and edges:
$$(i,\alpha)\sim(j,\beta)\Longleftrightarrow i\sim j,\, \alpha\sim\beta$$

\item The Cartesian product $X\,\square\,Y$ has vertex set $X\times Y$, and edges:
$$(i,\alpha)\sim(j,\beta)\Longleftrightarrow i=j,\, \alpha\sim\beta\mbox{ \rm{or} }i\sim j,\alpha=\beta$$

\item The lexicographic product $X\circ Y$ has vertex set $X\times Y$, and edges:
$$(i,\alpha)\sim(j,\beta)\Longleftrightarrow \alpha\sim\beta\mbox{ \rm{or} }\alpha=\beta,\,
 i\sim j$$
\end{enumerate}
\end{definition}

The direct product $X\times Y$ is the usual one in a categorical sense. The Cartesian product $X\,\square\,Y$ is quite natural too from a geometric viewpoint, for instance because a product by a segment gives a prism. As for the lexicographic product $X\circ Y$, this is something interesting too, obtained by putting a copy of $X$ at each vertex of $Y$.

\bigskip

The above products are all well-known in graph theory, and at the level of symmetry groups, we have some straightforward embeddings, as follows:
$$G(X)\times G(Y)\subset G(X \times Y)$$
$$G(X)\times G(Y)\subset G(X\,\square\,Y)$$
$$G(X)\wr G(Y)\subset G(X\circ Y)$$

Following \cite{bb1}, these embeddings have the following quantum analogues, using the various product operations constructed in chapter 1 and chapter 4:

\begin{proposition}
We have embeddings as follows,
$$G^+(X)\times G^+(Y)\subset G^+(X \times Y)$$
$$G^+(X)\times G^+(Y)\subset G^+(X\,\square\,Y)$$
$$G^+(X)\wr_*G^+(Y)\subset G^+(X\circ Y)$$
with the operation $\wr_*$ being a free wreath product. 
\end{proposition}

\begin{proof}
All the assertions are elementary, coming from definitions, as follows:

\medskip

(1) We use the following identification, given by $\delta_{(i,\alpha)}=\delta_i\otimes\delta_\alpha$:
$$C(X\times Y)=C(X)\otimes C(Y)$$

The adjacency matrix of the direct product is then given by:
$$d_{X\times Y}=d_X\otimes d_Y$$

Thus if $u$ commutes with $d_X$ and $v$ commutes with $d_Y$, then $u\otimes v=(u_{ij}v_{\alpha\beta})_{(i\alpha,j\beta)}$ is a magic unitary that commutes with $d_{X\times Y}$. But this gives a morphism as follows:
$$C(G^+(X\times Y))\to C(G^+(X)\times G^+(Y))$$ 

Finally, the surjectivity of this morphism follows by summing over $i$ and $\beta$.

\medskip

(2) The proof here is nearly identical to the one of (1). Indeed, with the identification there, the adjacency matrix of the Cartesian product is given by:
$$d_{X\,\square\,Y}=d_X\otimes1+1\otimes d_Y$$

Thus if $u$ commutes with $d_X$ and $v$ commutes with $d_Y$, then $u\otimes v=(u_{ij}v_{\alpha\beta})_{(i\alpha,j\beta)}$ is a magic unitary that commutes with $d_{X\,\square\,Y}$, and this gives the result.

\medskip

(3) With the same identification as before, namely $\delta_{(i,\alpha)}=\delta_i\otimes\delta_\alpha$, the adjacency matrix of the lexicographic product $X\circ Y$ is given by:
$$d_{X\circ Y}=d_X\otimes1+\mathbb I\otimes d_Y$$

Now let $u,v$ be the magic unitary matrices of $G^+(X),G^+(Y)$. The magic unitary matrix of $G^+(X)\wr_*G^+(Y)$ is then given by:
$$w_{ia,jb}= u_{ij}^{(a)}v_{ab}$$

Since $u,v$ commute with $d_X,d_Y$, we get that $w$ commutes with $d_{X\circ Y}$. But this gives the desired morphism, and the surjectivity follows by summing over $i$ and $b$. 
\end{proof}

The problem now is that of deciding when the embeddings in Proposition 5.14 are isomorphisms. This is something non-trivial, even at the level of the classical symmetry groups, and the results here will be quite technical. Following \cite{bb1}, we first have: 

\begin{theorem}
Let $X$ and $Y$ be finite connected regular graphs. If their spectra $\{\lambda\}$ and $\{\mu\}$ do not contain $0$ and satisfy
$$\{ \lambda_i/\lambda_j\} \cap \{\mu_k/\mu_l\}
= \{1\}$$
then $G^+(X \times Y)=G^+(X)\times G^+(Y)$. Also, if their spectra satisfy
$$\{\lambda_i - \lambda_j \} \cap \{\mu_k - \mu_l\}
= \{0\}$$
then $G^+(X\,\square\,Y)=G^+(X)\times G^+(Y)$.
\end{theorem}

\begin{proof}
Let $\lambda_1$ be the valence of $X$. Since $X$ is regular we have $\lambda_1 \in Sp(X)$, with $1$ as eigenvector, and since $X$ is connected $\lambda_1$ has multiplicity 1. Hence if $P_1$ is the orthogonal projection onto ${\mathbb C}1$, the spectral decomposition of $d_X$ is of the following form:
$$d_X = \lambda_1 P_1 + \sum_{i\not=1}\lambda_i P_i$$

We have a similar formula for the adjacency matrix $d_Y$, namely:
$$d_Y = \mu_1 Q_1 + \sum_{j\not=1}\mu_j Q_j$$

This gives the following formulae for products:
$$d_{X\times Y}=\sum_{ij}(\lambda_i\mu_j)P_{i}\otimes Q_{j}\quad,\quad 
d_{X\,\square\,Y} = \sum_{i,j}(\lambda_i + \mu_i)P_i \otimes Q_j$$

Here the projections form partitions of unity, and the scalar are distinct, so these are spectral decompositions. The coactions will commute with any of the spectral projections, and hence with both $P_1 \otimes 1$, $1 \otimes Q_1$. In both cases the universal coaction $v$ is the tensor product of its restrictions to the images of $P_1\otimes 1$, $1\otimes Q_1$, which gives the result.
\end{proof}

Regarding the lexicographic products, we have the following result, also from \cite{bb1}:

\begin{theorem}
Let $X,Y$ be regular graphs, with $X$ connected. If their spectra $\{\lambda_i\}$ and $\{\mu_j\}$ satisfy the condition
$$\{ \lambda_1-\lambda_i\mid i\neq 1 \} \cap \{-n\mu_j\} = \emptyset$$
where $n$ and $\lambda_1$ are the order and valence of $X$, then $G^+(X\circ Y)=G^+(X)\wr_*G^+(Y)$.   
\end{theorem}

\begin{proof}
We denote by $P_i,Q_j$ the spectral projections corresponding to $\lambda_i,\mu_j$. Since $X$ is connected we have $P_1=\mathbb I/n$, and we obtain:
\begin{eqnarray*}
d_{X\circ Y}
&=&d_X\otimes 1+{\mathbb I}\otimes d_Y\\
&=&\left(\sum_i\lambda_iP_i\right)\otimes\left(\sum_jQ_j\right)+\left(nP_1\right)\otimes \left(\sum_i\mu_jQ_j\right)\\
&=&\sum_j(\lambda_1+n\mu_j)(P_1 \otimes Q_j) + \sum_{i\not=1}\lambda_i (P_i\otimes 1)
\end{eqnarray*} 

In this formula the projections form a partition of unity and scalars are distinct, so this is the spectral decomposition of $d_{X\circ Y}$. Now let $W$ be the universal magic matrix for $X\circ Y$. Then $W$ must commute with all spectral projections, and in particular:
$$[W,P_1 \otimes Q_j]=0$$

Summing over $j$ gives $[W, P_1 \otimes 1]=0$, so $1\otimes C(Y)$ is invariant under the coaction. The corresponding restriction of $W$ gives a coaction of $G^+(X\circ Y)$ on $1 \otimes C(Y)$, say:
$$W(1 \otimes e_a) = \sum_b 1 \otimes e_b \otimes y_{ba}$$

Here $y$ is a magic unitary. On the other hand we can write:
$$W(e_i \otimes 1) = \sum_{jb} e_j \otimes e_b \otimes x_{ji}^b$$  

By multiplying by the previous relation we get:
$$W(e_i \otimes e_a)
=\sum_{jb} e_j \otimes e_b \otimes
y_{ba}x_{ji}^b
=\sum_{jb} e_j \otimes e_b \otimes x_{ji}^b y_{ba}$$

This shows that coefficients of $W$ are of the following form:
$$W_{jb,ia} = y_{ba} x_{ji}^b=x_{ji}^b y_{ba}$$

Consider now the matrix $x^b=(x_{ij}^b)$. Since $W$ is a morphism of algebras, each row of $x^b$ is a partition  of unity. Also using the antipode, we have
\begin{eqnarray*}
S\left(\sum_jx_{ji}^{b}\right)
&=&S\left(\sum_{ja}x_{ji}^{b}y_{ba}\right)
=S\left(\sum_{ja}W_{jb,ia}\right)\\
&=&\sum_{ja}W_{ia,jb}
=\sum_{ja}x_{ij}^ay_{ab}\\
&=&\sum_ay_{ab}
=1
\end{eqnarray*}

Thus $x^b$ is magic. We check now that $x^a,y$ commute with $d_X,d_Y$. We have:
$$(d_{X\circ Y})_{ia,jb} = (d_X)_{ij}\delta_{ab} + (d_Y)_{ab}$$

Thus the two products between $W$ and $d_{X\circ Y}$ are given by:
$$(Wd_{X\circ Y})_{ia,kc}=\sum_j W_{ia,jc} (d_X)_{jk} + \sum_{jb}W_{ia,jb}(d_Y)_{bc}$$
$$(d_{X\circ Y}W)_{ia,kc}=\sum_j (d_X)_{ij} W_{ja,kc} + \sum_{jb}(d_Y)_{ab}W_{jb,kc}$$

Now since $W$ commutes with $d_{X\circ Y}$, the terms on the right are equal, and by summing over $c$ we get:
$$\sum_j x_{ij}^a(d_X)_{jk} + \sum_{cb} y_{ab}(d_Y)_{bc}
= \sum_{j} (d_X)_{ij}x_{jk}^a + \sum_{cb} (d_Y)_{ab}y_{bc}$$

The second sums in both terms are equal to the valency of $Y$, so we get $[x^a,d_X]=0$. Now once again from the formula coming from commutation of $W$ with $d_{X\circ Y}$, we get $[y,d_Y] =0$. Summing up, the coefficients of $W$ are of the following form, where $x^b$ are magic unitaries commuting with $d_X$, and $y$ is a magic unitary commuting with $d_Y$: 
$$W_{jb,ia}=x_{ji}^by_{ba}$$

But this gives a morphism $C(G^+(X)\wr_*G^+(Y)\to G^+(X\circ Y)$ mapping $u_{ji}^{(b)}\to x_{ji}^b$ and $v_{ba}\to y_{ba}$, which is inverse to the morphism in Proposition 5.14, as desired.
\end{proof}

As an application, we have the following result:

\index{connected graph}

\begin{theorem}
Given a connected graph $X$, and $k\in\mathbb N$, we have the formulae
$$G(kX)=G(X)\wr S_k\quad,\quad 
G^+(kX)=G^+(X)\wr_*S_k^+$$
where $kX=X\sqcup\ldots\sqcup X$ is the $k$-fold disjoint union of $X$ with itself.
\end{theorem}

\begin{proof}
The first formula is something well-known, which follows as well from the second formula, by taking the classical version. Regarding now the second formula, it is elementary to check that we have an inclusion as follows, for any finite graph $X$:
$$G^+(X)\wr_*S_k^+\subset G^+(kX)$$

Regarding now the reverse inclusion, which requires $X$ to be connected, this follows by doing some matrix analysis, by using the commutation with $u$. To be more precise, let us denote by $w$ the fundamental corepresentation of $G^+(kX)$, and set:
$$u_{ij}^{(a)}=\sum_bw_{ia,jb}\quad,\quad 
v_{ab}=\sum_iv_{ab}$$

It is then routine to check, by using the fact that $X$ is indeed connected, that we have here magic unitaries, as in the definition of the free wreath products. Thus, we obtain:
$$G^+(kX)\subset G^+(X)\wr_*S_k^+$$

But this gives the result, as a consequence of Theorem 5.16. See \cite{bb1}.
\end{proof}

We are led in this way to the following result, from \cite{bbc}:

\index{hyperoctahedral group}
\index{free hyperoctahedral group}

\begin{theorem}
Consider the graph consisting of $N$ segments.
\begin{enumerate}
\item Its symmetry group is the hyperoctahedral group $H_N=\mathbb Z_2\wr S_N$.

\item Its quantum symmetry group is the quantum group $H_N^+=\mathbb Z_2\wr_*S_N^+$.
\end{enumerate}
\end{theorem}

\begin{proof}
Here the first assertion is clear from definitions, with the remark that the relation with the formula $H_N=G(\square_N)$ comes by viewing the $N$ segments as being the $[-1,1]$ segments on each of the $N$ coordinate axes of $\mathbb R^N$. Indeed, a symmetry of the $N$-cube is the same as a symmetry of the $N$ segments, and so, as desired: 
$$G(\square_N)=\mathbb Z_2\wr S_N$$

As for the second assertion, this follows from Theorem 5.17, applied to the segment graph. Observe also that (2) implies (1), by taking the classical version.
\end{proof}

We will be back to quantum reflections later, in chapter 6 below, with a systematic study of such quantum groups, and with some generalizations as well.

\bigskip

As an application of the above methods, as explained in \cite{bb1}, it is possible to compute the quantum isometry groups of all vertex-transitive graphs of order $\leq 11$, except for the Petersen graph. This latter graph $P_{10}$, which is a famous graph, is as follows:
$$\xymatrix@R=1pt@C=5pt{
&&&&\circ\ar@{-}[dddddrrrr]\ar@{-}[dddddllll]\\
\\
\\
\\
\\
\circ\ar@{-}[ddddddddr]&&&&\circ\ar@{-}[uuuuu]\ar@{-}[ddddddl]\ar@{-}[ddddddr]&&&&\circ\ar@{-}[ddddddddl]\\
\\
&&
\circ\ar@{-}[uull]\ar@{-}[ddddrrr]\ar@{-}[rrrr]&&&&\circ\ar@{-}[uurr]\ar@{-}[ddddlll]\\
\\
\\
\\
&&&\circ&&\circ\\
\\
&\circ\ar@{-}[rrrrrr]\ar@{-}[uurr]&&&&&&\circ\ar@{-}[uull]}
$$

In order to explain the computation for $P_{10}$, done by Schmidt in \cite{sc1}, we will need a number of preliminaries. Let us start with the following notion, from \cite{bi1}:

\begin{definition}
The reduced quantum automorphism group of $X$ is given by
$$C(G^*(X))=C(G^+(X))\Big/\left<u_{ij}u_{kl}=u_{kl}u_{ij}\Big|\forall i\sim k,j\sim l\right>$$
with $i\sim j$ standing as usual for the fact that $i,j$ are connected by an edge.
\end{definition}

As explained by Bichon in \cite{bi1}, the above construction produces indeed a quantum group $G^*(X)$, which sits as an intermediate subgroup, as follows:
$$G(X)\subset G^*(X)\subset G^+(X)$$

There are many things that can be said about this construction, but in what concerns us, we will rather use it as a technical tool. Following Schmidt \cite{sc1}, we have:

\begin{proposition}
Assume that a regular graph $X$ is strongly regular, with parameters $\lambda=0$ and $\mu=1$, in the sense that:
\begin{enumerate}
\item $i\sim j$ implies that $i,j$ have $\lambda$ common neighbors.

\item $i\not\sim j$ implies that $i,j$ have $\mu$ common neighbors.
\end{enumerate}
The quantum group inclusion $G^*(X)\subset G^+(X)$ is then an isomorphism.
\end{proposition}

\begin{proof}
This is something quite tricky, the idea being as follows:

\medskip

(1) First of all, regarding the statement, a graph is called regular, with valence $k$, when each vertex has exactly $k$ neighbors. Then we have the notion of strong regularity, given by the conditions (1,2) in the statement. And finally we have the notion of strong regularity with parameters $\lambda=0,\mu=1$, that the statement is about, and with as main example here $P_{10}$, which is 3-regular, and strongly regular with $\lambda=0,\mu=1$.

\medskip

(2) Regarding now the proof, we must prove that the following commutation relation holds, with $u$ being the magic unitary of the quantum group $G^+(X)$:
$$u_{ij}u_{kl}=u_{kl}u_{ij}\ ,\ \forall i\sim k,j\sim l$$

But for this purpose, we can use the $\lambda=0,\mu=1$ strong regularity of our graph, by inserting some neighbors into our computation. To be more precise, we have:
\begin{eqnarray*}
u_{ij}u_{kl}
&=&u_{ij}u_{kl}\sum_{s\sim l}u_{is}\\
&=&u_{ij}u_{kl}u_{ij}+\sum_{s\sim l,s\neq j}u_{ij}u_{kl}u_{is}\\
&=&u_{ij}u_{kl}u_{ij}+\sum_{s\sim l,s\neq j}u_{ij}\left(\sum_au_{ka}\right)u_{is}\\
&=&u_{ij}u_{kl}u_{ij}+\sum_{s\sim l,s\neq j}u_{ij}u_{is}\\
&=&u_{ij}u_{kl}u_{ij}
\end{eqnarray*}

(3) But this gives the result. Indeed, we conclude from this that $u_{ij}u_{kl}$ is self-adjoint, and so, by conjugating, that we have $u_{ij}u_{kl}=u_{kl}u_{ij}$, as desired.
\end{proof}

In the particular case of the Petersen graph $P_{10}$, which in addition is 3-regular, we can further build on the above result, and still following Schmidt \cite{sc1}, we have:

\begin{theorem}
The Petersen graph has no quantum symmetry,
$$G^+(P_{10})=G(P_{10})=S_5$$
with $S_5$ acting in the obvious way.
\end{theorem}

\begin{proof}
In view of Proposition 5.20, we must prove that the following commutation relation holds, with $u$ being the magic unitary of the quantum group $G^+(P_{10})$:
$$u_{ij}u_{kl}=u_{kl}u_{ij}\ ,\ \forall i\not\sim k,j\not\sim l$$

We can assume $i\neq k$, $j\neq l$. Now if we denote by $s,t$ the unique vertices having the property $i\sim s,k\sim s$ and $j\sim t,l\sim t$, a routine study shows that we have:
$$u_{ij}u_{kl}=u_{ij}u_{st}u_{kl}$$

With this in hand, if we denote by $q$ the third neighbor of $t$, we obtain:
\begin{eqnarray*}
u_{ij}u_{kl}
&=&u_{ij}u_{st}u_{kl}(u_{ij}+u_{il}+u_{iq})\\
&=&u_{ij}u_{st}u_{kl}u_{ij}+0+0\\
&=&u_{ij}u_{st}u_{kl}u_{ij}\\
&=&u_{ij}u_{kl}u_{ij}
\end{eqnarray*}

Thus $u_{ij}u_{kl}$ is self-adjoint, and so $u_{ij}u_{kl}=u_{kl}u_{ij}$, as desired.
\end{proof}

As an application of all this, we have the following classification table from \cite{bb1}, improved using \cite{sc1}, containing all vertex-transitive graphs of order $\leq 11$ modulo complementation, with their classical and quantum symmetry groups:
\begin{center}\begin{tabular}[t]{|l|l|l|l|}\hline
Order&Graph&Classical group&Quantum group\\ \hline\hline 2&$K_2$&$\mathbb Z_2$&$\mathbb  Z_2$\\
\hline\hline 3&$K_3$&$S_3$&$S_3$\\ \hline\hline
4&$2K_2$&$H_2$&$H_2^+$\\
\hline 4&$K_4$&$S_4$&$S_4^+$\\
\hline\hline 5&$C_5$&$D_5$&$D_5$\\ \hline
5&$K_5$&$S_5$&$S_5^+$\\ \hline\hline 6&$C_6$&$D_6$&$D_6$\\ \hline
6&$2K_3$&$S_3\wr\mathbb Z_2$&$S_3{\,\wr_*\,}\mathbb Z_2$\\
\hline 6&$3K_2$&$H_3$&$H_3^+$\\ \hline 6&$K_6$&$S_6$&$S_6^+$\\
\hline\hline 7&$C_7$&$D_7$&$D_7$ \\ \hline
7&$K_7$&$S_7$&$S_7^+$\\ \hline\hline 8&$C_8$, $C_8^+$&$D_8$&$D_8$\\
\hline 8&$P(C_4)$& $H_3$&$S_4^+\times Z_2$ \\ \hline 8&$2K_4$&$S_4\wr \mathbb Z_2$&$S_4^+{\,\wr_*\,}\mathbb Z_2$ \\
\hline 8&$2C_4$& $H_2\wr\mathbb Z_2$ & $H_2^+{\,\wr_*\,}\mathbb Z_2$
\\ \hline 8&$4K_2$&$H_4$&$H_4^+$ \\ \hline 8&$K_8$&$S_8$&$S_8^+$\\
\hline\hline 9&$C_9$, $C_9^3$&$D_9$&$D_9$\\ \hline 9 & $K_3\times K_3$&$S_3\wr\mathbb Z_2$&$S_3\wr\mathbb Z_2$
\\ \hline 9&$3K_3$&$S_3\wr S_3$&$S_3{\,\wr_*\,}S_3$ \\ \hline
9&$K_9$&$S_9$&$S_9^+$ \\ \hline \hline 10&$C_{10}$, $C_{10}^2$,
$C_{10}^+$, $P(C_5)$&$D_{10}$&$D_{10}$\\ \hline 10 &
$P(K_5)$&$S_5\times\mathbb  Z_2$&$S_5^+\times\mathbb  Z_2$\\ \hline
10&$C_{10}^4$&$\mathbb Z_2\wr D_5$&$\mathbb Z_2{\,\wr_*\,}D_5$\\ \hline
10&$2C_5$&$D_5\wr\mathbb  Z_2$&$D_5{\,\wr_*\,}\mathbb Z_2$\\ \hline
10&$2K_{5}$&$S_5\wr\mathbb  Z_2$&$S_5^+{\,\wr_*\,}\mathbb Z_2$\\ \hline
10&$5K_2$&$H_5$&$H_5^+$\\ \hline 
10&$K_{10}$&$S_{10}$&$S_{10}^+$\\ \hline
10&$P_{10}$&$S_5$&$S_5$\\
\hline\hline 11&$C_{11}$, $C_{11}^2$,
$C_{11}^3$&$D_{11}$&$D_{11}$\\ \hline
11&$K_{11}$&$S_{11}$&$S_{11}^+$\\ \hline
\end{tabular}\end{center}

Here $K$ denote the complete graphs, $C$ the cycles with chords, and $P$ stands for prisms. Moreover, by using more advanced techniques, the above table can be considerably extended. For more on all this, we refer to Schmidt's papers \cite{sc1}, \cite{sc2}, \cite{sc3}.

\section*{5d. Circulant graphs}

Following \cite{bbg}, let us discuss now the case of the circulant graphs, which is quite interesting, due to the fact that we can use the Fourier transform. We first have:

\index{type of circulant graph}

\begin{definition}
Associated to any circulant graph $X$ having $N$ vertices are:
\begin{enumerate}
\item The set $S\subset {\mathbb Z}_N$ given by $i\sim_X j \iff j-i \in S$.

\item The group $E\subset{\mathbb Z}_N^*$ consisting of elements $g$ such that $gS=S$.

\item The number $k=|E|$, called type of $X$.
\end{enumerate}
\end{definition}

In what follows $X$ will be a circulant graph having $p$ vertices, with $p$ prime. We denote by $\xi$ the column vector $(1,w,w^2,\ldots ,w^{p-1})$, where $w=e^{2\pi i/p}$. We have:

\begin{proposition}
The eigenspaces of $d$ are given by $V_0={\mathbb C}1$ and
$$V_x=\bigoplus_{a\in E}{\mathbb C}\,\xi^{xa}$$
with $x\in \mathbb Z_p^*$. Moreover, we have $V_x=V_y$ if and only if $xE=yE$.
\end{proposition}

\begin{proof}
Since $d$ is circulant, we have $d(\xi^x)=f(x)\xi^x$, with $f:{\mathbb Z}_p\to{\mathbb C}$ being:
$$f(x)=\sum_{t\in S}w^{xt}$$

Let $K={\mathbb Q}(w)$ and let $H$ be the Galois group of the Galois extension $\mathbb Q \subset K$. We have then a group isomorphism as follows:
$$\mathbb Z_p^*\simeq H\quad,\quad x\to s_x=[w\to w^x]$$

Also, we know from a theorem of Dedekind that the family $\{s_x\mid x\in{\mathbb Z}_p^*\}$ is free in ${\rm End}_{\mathbb Q}(K)$. Now for $x,y\in \mathbb Z_p^*$ consider the following operator:
$$L = \sum_{t \in S} s_{xt} - \sum_{t \in S} s_{yt} \in
End_{\mathbb Q}(K)$$

We have $L({w}) = f(x)-f(y)$, and since $L$ commutes with the action of $H$, we have:
$$L=0 \iff L({w}) =0 \iff f(x)=f(y)$$

By linear independence of the family $\{s_x\mid x\in \mathbb Z_p^*\}$ we get:
$$f(x) = f(y) \iff xS=yS \iff xE=yE$$

It follows that $d$ has precisely $1+(p-1)/k$ distinct eigenvalues, the corresponding
eigenspaces being those in the statement.
\end{proof}

Consider now a commutative ring $(R,+,\cdot)$. We denote by $R^*$ the group of invertibles, and we assume $2\in R^*$. A subgroup $G\subset R^*$ is called even if $-1\in G$. We have:

\begin{definition}
An even subgroup $G\subset R^*$ is called $2$-maximal if, inside $G$:
$$a-b=2(c-d)\implies a=\pm b$$
We call $a=b,c=d$ trivial solutions, and $a=-b=c-d$ hexagonal solutions. 
\end{definition}

To be more precise, in what regards our terminology, consider the group $G\subset{\mathbb C}$ formed by $k$-th roots of unity, with $k$ even. An equation of the form $a-b=2(c-d)$ with $a,b,c,d\in G$ says that the diagonals $a-b$ and $c-d$ must be parallel, and that the first one is twice as much as the second one. But this can happen only when $a,c,d,b$ are consecutive vertices of a regular hexagon, and here we have $a+b=0$.

\bigskip

The relation with our quantum symmetry considerations comes from:

\begin{proposition}
Assume that $R$ has the property $3\neq 0$, and consider a $2$-maximal subgroup $G\subset R^*$. Then, the following happen:
\begin{enumerate}
\item $2,3\not\in G$. 

\item $a+b=2c$ with $a,b,c\in G$ implies $a=b=c$. 

\item $a+2b=3c$ with $a,b,c\in G$ implies $a=b=c$.
\end{enumerate}
\end{proposition}

\begin{proof}
All these assertions are elementary, as follows:

\medskip

(1) This follows from the following formulae, which cannot hold in
$G$:
$$4-2=2(2-1)\quad,\quad 
3-(-1)=2(3-1)$$

Indeed, the first one would imply $4=\pm 2$, and the second one
would imply $3=\pm 1$. But from $2\in R^*$ and $3\neq 0$ we get
$2,4,6\neq 0$, contradiction.

\medskip

(2) We have $a-b=2(c-b)$. For a trivial solution we have
 $a=b=c$, and for a hexagonal
 solution we have $a+b=0$,
hence $c=0$, hence $0\in{G}$, contradiction.

\medskip

(3) We have $a-c=2(c-b)$. For a trivial solution we have
 $a=b=c$, and for a hexagonal
 solution we have $a+c=0$,
hence $b=-2a$, hence $2\in{G}$, contradiction.
\end{proof}

We can now formulate a general non quantum symmetry result, as follows:

\begin{theorem}
A circulant graph $X$ with a prime number $p\geq5$ of vertices, such that the corresponding group $E\subset {\mathbb Z}_p$ is $2$-maximal, has no quantum
symmetry.
\end{theorem}

\begin{proof}
We use Proposition 5.23, which ensures that $V_1,V_2,V_3$ are eigenspaces of $d$. By $2$-maximality of $E$, these three eigenspaces are different. From eigenspace preservation and 2-maximality we get formulae of the
following type, with $r_a\in A=C(G^+(X))$:
$$\alpha (\xi)=\sum_{a\in E}\xi^a\otimes r_a\quad,\quad 
\alpha(\xi^2)=\sum_{a\in E}\xi^{2a}\otimes r_a^2\quad,\quad 
\alpha(\xi^3)=\sum_{a\in E}\xi^{3a}\otimes r_a^3$$

Our claim now is that for $a\neq b$, we have the following key formula:
$$r_ar_b^3=0$$

Indeed, as explained in \cite{bbg}, this follows by examining the following equality:
$$\left(\sum_{a\in E}\xi^a\otimes r_a\right)
\left(\sum_{b\in E}\xi^{2b}\otimes r_b^2\right)=\sum_{c\in
E}\xi^{3c}\otimes r_c^3$$

By using this key formula, we get by recurrence on $s\geq 3$ that we have:
$$\alpha\left(\xi^{1+s}\right)
=\sum_{a\in E}\xi^{(1+s)a}\otimes r_a^{1+s}$$

But this gives $r_a^* = r_a^{p-1}$ for any $a$, and now by using the key formula, we get:
$$(r_ar_b)(r_ar_b)^*
=r_ar_br_b^*r_a^*
=r_ar_b^pr_a^*
=(r_ar_b^3)(r_b^{p-3}r_a^*)
=0$$

Thus $r_ar_b=r_br_a=0$. On the other hand, $A$ is generated by coefficients of $\alpha$, which are in turn powers of elements $r_a$. It follows that $A$ is commutative, and we are done.
\end{proof}

Still following \cite{bbg}, we can now formulate a main result, as follows:

\begin{theorem}
A type $k$ circulant graph having $p>>k$ vertices, with $p$ prime,
has no quantum symmetry.
\end{theorem}

\begin{proof}
This follows from Theorem 5.26 and some arithmetics, as follows:

\medskip

(1) Let $k$ be an even number, and consider the group of $k$-th roots of unity $G=\{1,\zeta,\ldots ,\zeta^{k-1}\}$, where $\zeta=e^{2\pi i/k}$. By some standard arithmetics, $G$ is $2$-maximal in ${\mathbb C}$.

\medskip

(2) As a continuation of this, again by some standard arithmetics, for $p>6^{\varphi(k)}$, with $\varphi$ being the Euler function, any subgroup $E\subset{\mathbb Z}_p^*$ of order $k$ is $2$-maximal.

\medskip

(3) But this proves our result. Indeed, by using (2), we can apply Theorem 5.26 provided that we have $p>6^{{\varphi(k)}}$, and our graph has no quantum symmetry, as desired.
\end{proof}

\section*{5e. Exercises} 

The constructions in this chapter produce many interesting examples of quantum permutation groups, and we will keep studying them. As a first exercise, we have:

\begin{exercise}
Learn more about quantum symmetries of finite graphs, from Lupini, Man\v cinska, Roberson, Schmidt et al., with a look into Musto-Reutter-Verdon too.
\end{exercise}

The idea here is that a lot of theory, going beyond what we can reasonably explain in this book, has been recently developed by the above people, and their collaborators. We will be back to this later, in this book, but only with some very basic results.

\begin{exercise}
Work out some theoretical generalizations of the no quantum symmetry result for the $N$-cycle, in the $N\to\infty$ limit.
\end{exercise}

This is something that we already discussed in the above, with the comment that some technology is available in the case where $N$ is prime, and that all this is in need of some further extensions, to the case where $N$ is not necessarily prime.

\chapter{Quantum reflections}

\section*{6a. Real reflections}

In this chapter and in the next one, which are central to the present book, we discuss the quantum reflection groups and their twisted analogues, following \cite{ba2}, \cite{bb+}, \cite{bbc}, \cite{bve} and related papers. The material will be quite technical, and we will be here at the core of modern quantum group theory. The point indeed is that, no matter what applications of your quantum groups you are looking for, it is most likely that these applications will come from a suitable quantum reflection group. Or at least, that is a common belief.

\bigskip

In order to get started, let us go back to the square graph problem. In order to present the correct, final solution to it, the idea will be that of looking at $G^+(|\ |)$, which is equal to $G^+(\square)$. We recall from chapter 5 that we have the following result, from \cite{bbc}:

\index{hyperoctahedral group}
\index{free hyperoctahedral group}

\begin{theorem}
Consider the graph consisting of $N$ segments.
\begin{enumerate}
\item Its symmetry group is the hyperoctahedral group $H_N=\mathbb Z_2\wr S_N$.

\item Its quantum symmetry group is the quantum group $H_N^+=\mathbb Z_2\wr_*S_N^+$.
\end{enumerate}
\end{theorem}

\begin{proof}
This is something that we know, the idea being as follows:

\medskip

(1) This is clear from definitions, and with the remark that $H_N$ appears as well as the symmetry group of the hypercube, $G(\square_N)=H_N$.

\medskip

(2) This is something that we know from the previous chapter, and with the remark that for the hypercube we obtain something different, $G(\square_N)=O_N^{-1}$.
\end{proof}

Now back to the square, we have $G^+(\square)=H_2^+$, and our claim is that this is the ``good'' and final formula. In order to prove this, we must work out the easiness theory for $H_N,H_N^+$, and find a compatibility there. We first have the following result:

\index{cubic matrix}
\index{sudoku matrix}

\begin{proposition}
The algebra $C(H_N^+)$ can be presented in two ways, as follows:
\begin{enumerate}
\item As the universal algebra generated by the entries of a $2N\times 2N$ magic unitary having the ``sudoku'' pattern $w=(^a_b{\ }^b_a)$, with $a,b$ being square matrices.

\item As the universal algebra generated by the entries of a $N\times N$ orthogonal matrix which is ``cubic'', in the sense that $u_{ij}u_{ik}=u_{ji}u_{ki}=0$, for any $j\neq k$.
\end{enumerate}
As for $C(H_N)$, this has similar presentations, among the commutative algebras.
\end{proposition}

\begin{proof}
We must prove that the algebras $A_s,A_c$ coming from (1,2) coincide. We can define a morphism $A_c\to A_s$ by the following formula:
$$\varphi(u_{ij})=a_{ij}-b_{ij}$$

We construct now the inverse morphism. Consider the following elements:
$$\alpha_{ij}=\frac{u_{ij}^2+u_{ij}}{2}\quad,\quad
\beta_{ij}=\frac{u_{ij}^2-u_{ij}}{2}$$

These are projections, and the following matrix is a sudoku unitary:
$$M=\begin{pmatrix}
(\alpha_{ij})&(\beta_{ij})\\
(\beta_{ij})&(\alpha_{ij})
\end{pmatrix}$$

Thus we can define a morphism $A_s\to A_c$ by the following formula:
$$\psi(a_{ij})=\frac{u_{ij}^2+u_{ij}}{2}\quad,\quad 
\psi(b_{ij})=\frac{u_{ij}^2-u_{ij}}{2}$$

We check now the fact that $\psi,\varphi$ are indeed inverse
morphisms:
$$\psi\varphi(u_{ij})
=\psi(a_{ij}-b_{ij})
=\frac{u_{ij}^2+u_{ij}}{2}-\frac{u_{ij}^2-u_{ij}}{2}
=u_{ij}$$

As for the other composition, we have the following computation:
$$\varphi\psi(a_{ij})
=\varphi\left(\frac{u_{ij}^2+u_{ij}}{2}\right)
=\frac{(a_{ij}-b_{ij})^2+(a_{ij}-b_{ij})}{2}
=a_{ij}$$

A similar computation gives $\varphi\psi(b_{ij})=b_{ij}$, as desired. As for the final assertion, regarding $C(H_N)$, this follows from the above results, by taking classical versions.
\end{proof}

We can now work out the easiness property of $H_N,H_N^+$, with respect to the cubic representations, and we are led to the following result, from \cite{bbc}:

\begin{theorem}
The quantum groups $H_N,H_N^+$ are both easy, as follows:
\begin{enumerate}
\item $H_N$ corresponds to the category $P_{even}$.

\item $H_N^+$ corresponds to the category $NC_{even}$.
\end{enumerate}
\end{theorem}

\begin{proof}
This is something quite routine, the idea being as follows:

\medskip

(1) We know that $H_N^+\subset O_N^+$ appears via the cubic relations, namely:
$$u_{ij}u_{ik}=u_{ji}u_{ki}=0\quad,\quad\forall j\neq k$$

Our claim is that, in Tannakian terms, these relations reformulate as follows, with $H\in P(2,2)$ being the 1-block partition, joining all 4 points:
$$T_H\in End(u^{\otimes 2})$$

(2) In order to prove our claim, observe first that we have, by definition of $T_H$:
$$T_H(e_i\otimes e_j)=\delta_{ij}e_i\otimes e_i$$

With this formula in hand, we have the following computation:
\begin{eqnarray*}
T_Hu^{\otimes 2}(e_i\otimes e_j\otimes1)
&=&T_H\left(\sum_{abij}e_{ai}\otimes e_{bj}\otimes u_{ai}u_{bj}\right)(e_i\otimes e_j\otimes1)\\
&=&T_H\sum_{ab}e_a\otimes e_b\otimes u_{ai}u_{bj}\\
&=&\sum_ae_a\otimes e_a\otimes u_{ai}u_{aj}
\end{eqnarray*}

On the other hand, we have as well the following computation:
\begin{eqnarray*}
u^{\otimes 2}T_H(e_i\otimes e_j\otimes1)
&=&\delta_{ij}u^{\otimes 2}(e_i\otimes e_j\otimes1)\\
&=&\delta_{ij}\left(\sum_{abij}e_{ai}\otimes e_{bj}\otimes u_{ai}u_{bj}\right)(e_i\otimes e_j\otimes1)\\
&=&\delta_{ij}\sum_{ab}e_a\otimes e_b\otimes u_{ai}u_{bi}
\end{eqnarray*}

We conclude that $T_Hu^{\otimes 2}=u^{\otimes 2}T_H$ means that $u$ is cubic, as desired.

\medskip

(3) With our claim proved, we can go back to $H_N^+$. Indeed, it follows from Tannakian duality that this quantum group is easy, coming from the following category:
$$D=<H>=NC_{even}$$

(4) But this proves as well the result for $H_N$. Indeed, since this group is the classical version of $H_N^+$, we have as desired easiness, the corresponding category being:
$$E=<NC_{even},\slash\hskip-2.1mm\backslash>=P_{even}$$

Thus, we are led to the conclusions in the statement.
\end{proof}

As an immediate consequence of the above result, we have:

\begin{theorem}
The operation $H_N\to H_N^+$ is a liberation in the sense of easy quantum groups, in the sense that the category of partitions for $H_N^+$ appears as
$$D^+=D\cap NC$$
with $D$ being the category of partitions for $H_N$.
\end{theorem}

\begin{proof}
We already know, from definitions, that $H_N\to H_N^+$ is a liberation, in the sense that the classical version of $H_N^+$ is $H_N$. However, by using Theorem 6.3, we can see that much more is true, in the sense that $H_N\to H_N^+$ is an easy quantum group liberation, as stated, and with this coming from $NC_{even}=P_{even}\cap NC$.
\end{proof}

We refer to \cite{bbc} for more regarding the above results. In what concerns us, we will be back to this in a moment, with some probabilistic consequences as well.

\section*{6b. Complex reflections}

The reflection groups $H_N$ and their liberations $H_N^+$ belong in fact to two remarkable series, depending on a parameter $s\in\mathbb N\cup\{\infty\}$, constructed as follows:
$$H_N^s=\mathbb Z_s\wr S_N\quad,\quad 
H_N^{s+}=\mathbb Z_s\wr_*S_N^+$$

We discuss now, following \cite{bb+}, \cite{bve}, the algebraic and analytic structure of these latter quantum groups. The main motivation comes from the cases $s=1,2,\infty$, where we recover respectively $S_N,S_N^+$ and $H_N,H_N^+$, and the full reflection groups $K_N,K_N^+$. Let us start with a brief discussion concerning the classical case. The result that we will need is:

\begin{proposition}
The group $H_N^s=\mathbb Z_s\wr S_N$ of $N\times N$ permutation-like matrices having as nonzero entries the $s$-th roots of unity is as follows:
\begin{enumerate}
\item $H_N^1=S_N$ is the symmetric group.

\item $H_N^2=H_N$ is the hyperoctahedral group.

\item $H_N^\infty=K_N$ is the group of unitary permutation-like matrices.
\end{enumerate}
\end{proposition}

\begin{proof}
Everything here is clear from definitions, and with permutation-like meaning of course that each row and column contain exactly one nonzero entry.
\end{proof}

The free analogues of the reflection groups $H_N^s$ can be constructed as follows:

\begin{definition}
The algebra $C(H_N^{s+})$ is the universal $C^*$-algebra generated by $N^2$ normal elements $u_{ij}$, subject to the following relations,
\begin{enumerate}
\item $u=(u_{ij})$ is unitary,

\item $u^t=(u_{ji})$ is unitary,

\item $p_{ij}=u_{ij}u_{ij}^*$ is a projection,

\item $u_{ij}^s=p_{ij}$,
\end{enumerate}
with Woronowicz algebra maps $\Delta,\varepsilon,S$ constructed by universality.
\end{definition}

Here we allow the value $s=\infty$, with the convention that the last axiom simply disappears in this case. Observe that at $s<\infty$ the normality condition is actually redundant. This is because a partial isometry $a$ subject to the relation $aa^*=a^s$ is normal. As a first result now, making the connection with $H_N^s$, we have:

\begin{theorem}
We have an inclusion of quantum groups
$$H_N^s\subset H_N^{s+}$$
which is a liberation, in the sense that the classical version of $H_N^{s+}$, obtained by dividing by the commutator ideal, is the group $H_N^s$.
\end{theorem}

\begin{proof}
This follows as before for $O_N\subset O_N^+$ or for $S_N\subset S_N^+$, by using the Gelfand theorem, applied to the quotient of $C(H_N^{s+})$ by its commutator ideal.
\end{proof}

In analogy with the results from the real case, we have the following result:

\begin{proposition}
The algebras $C(H_N^{s+})$ with $s=1,2,\infty$, and their presentation relations in terms of the entries of the matrix $u=(u_{ij})$, are as follows:
\begin{enumerate}
\item For $C(H_N^{1+})=C(S_N^+)$, the matrix $u$ is magic: all its entries are projections, summing up to $1$ on each row and column.

\item For $C(H_N^{2+})=C(H_N^+)$ the matrix $u$ is cubic: it is orthogonal, and the products of pairs of distinct entries on the same row or the same column vanish.

\item For $C(H_N^{\infty+})=C(K_N^+)$ the matrix $u$ is unitary, its transpose is unitary, and all its entries are normal partial isometries.
\end{enumerate}
\end{proposition}

\begin{proof}
This is something elementary, from \cite{bb+}, \cite{bve}, the idea being as follows:

\medskip

(1) This follows from definitions and from standard operator algebra tricks.

\medskip

(2) This follows as well from definitions and standard operator algebra tricks.

\medskip

(3) This is just a translation of the definition of $C(H_N^{s+})$, at $s=\infty$.
\end{proof}

Let us prove now that $H_N^{s+}$ with $s<\infty$ is a quantum permutation group. For this purpose, we must change the fundamental representation. Let us start with:

\index{sudoku matrix}

\begin{definition}
A $(s,N)$-sudoku matrix is a magic unitary of size $sN$, of the form
$$m=\begin{pmatrix}
a^0&a^1&\ldots&a^{s-1}\\
a^{s-1}&a^0&\ldots&a^{s-2}\\
\vdots&\vdots&&\vdots\\
a^1&a^2&\ldots&a^0
\end{pmatrix}$$
where $a^0,\ldots,a^{s-1}$ are $N\times N$ matrices.
\end{definition}

The basic examples of such matrices come from the group $H_n^s$. Indeed, with $w=e^{2\pi i/s}$, each of the $N^2$ matrix coordinates $u_{ij}:H_N^s\to\mathbb C$ takes values in the following set:
$$S=\{0\}\cup\{1,w,\ldots,w^{s-1}\}$$

Thus, this coordinate function $u_{ij}:H_N^s\to\mathbb C$ decomposes as follows:
$$u_{ij}=\sum_{r=0}^{s-1}w^ra^r_{ij}$$

Here each $a^r_{ij}$ is a function taking values in $\{0,1\}$, and so a projection in the $C^*$-algebra sense, and it follows from definitions that these projections form a sudoku matrix. With this notion in hand, we have the following result, from \cite{bve}:

\begin{theorem}
The following happen:
\begin{enumerate}
\item The algebra $C(H_N^s)$ is isomorphic to the universal commutative $C^*$-algebra generated by the entries of a $(s,N)$-sudoku matrix.

\item The algebra $C(H_N^{s+})$ is isomorphic to the universal $C^*$-algebra generated by the entries of a $(s,N)$-sudoku matrix.
\end{enumerate}
\end{theorem}

\begin{proof}
The first assertion follows from the second one, via Theorem 6.7. In order to prove the second assertion, consider the universal algebra in the statement, namely:
$$A=C^*\left(a_{ij}^p\ \Big\vert \left(a^{q-p}_{ij}\right)_{pi,qj}=(s,N)-\mbox{sudoku }\right)$$

Consider also the algebra $C(H_N^{s+})$. According to Definition 6.6, this is presented by certain relations $R$, that we will call here level $s$ cubic conditions:
$$C(H_N^{s+})=C^*\left(u_{ij}\ \Big\vert\  u=N\times N\mbox{ level $s$ cubic }\right)$$

We will construct a pair of inverse morphisms between these algebras.

\medskip

(1) Our first claim is that $U_{ij}=\sum_pw^{-p}a^p_{ij}$ is a level $s$ cubic unitary. Indeed, by using the sudoku condition, the verification of (1-4) in Definition 6.6 is routine.

\medskip

(2) Our second claim is that the elements $A^p_{ij}=\frac{1}{s}\sum_rw^{rp}u^r_{ij}$, with the convention $u_{ij}^0=p_{ij}$, form a level $s$ sudoku unitary. Once again, the proof here is routine.

\medskip

(3) According to the above, we can define a morphism $\Phi:C(H_N^{s+})\to A$ by the formula $\Phi(u_{ij})=U_{ij}$, and a morphism $\Psi:A\to C(H_N^{s+})$ by the formula $\Psi(a^p_{ij})=A^p_{ij}$.

\medskip

(4) We check now the fact that $\Phi,\Psi$ are indeed inverse morphisms:
\begin{eqnarray*}
\Psi\Phi(u_{ij})
&=&\sum_pw^{-p}A^p_{ij}\\
&=&\frac{1}{s}\sum_pw^{-p}\sum_rw^{rp}u_{ij}^r\\
&=&\frac{1}{s}\sum_{pr}w^{(r-1)p}u_{ij}^r\\
&=&u_{ij}
\end{eqnarray*}

As for the other composition, we have the following computation:
\begin{eqnarray*}
\Phi\Psi(a^p_{ij})
&=&\frac{1}{s}\sum_rw^{rp}U_{ij}^r\\
&=&\frac{1}{s}\sum_rw^{rp}\sum_qw^{-rq}a_{ij}^q\\
&=&\frac{1}{s}\sum_qa_{ij}^q\sum_rw^{r(p-q)}\\
&=&a^p_{ij}
\end{eqnarray*}

Thus we have an isomorphism $C(H_N^{s+})=A$, as claimed.
\end{proof}

We will need the following simple fact: 

\begin{proposition}
A $sN\times sN$ magic unitary commutes with the matrix
$$\Sigma=
\begin{pmatrix}
0&I_N&0&\ldots&0\\
0&0&I_N&\ldots&0\\
\vdots&\vdots&&\ddots&\\
0&0&0&\ldots&I_N\\
I_N&0&0&\ldots&0
\end{pmatrix}$$
if and only if it is a sudoku matrix in the sense of Definition 6.9.
\end{proposition}

\begin{proof}
This follows from the fact that commutation with $\Sigma$ means that the matrix is circulant. Thus, we obtain the sudoku relations from Definition 6.9.
\end{proof}

Now let $Z_s$ be the oriented cycle with $s$ vertices, and consider the graph $NZ_s$ consisting of $N$ disjoint copies of it. Observe that, with a suitable labeling of the vertices, the adjacency matrix of this graph is the above matrix $\Sigma$. We obtain from this:

\begin{theorem}
We have the following results:
\begin{enumerate}
\item $H_N^s$ is the symmetry group of $NZ_s$.

\item $H_N^{s+}$ is the quantum symmetry group of $NZ_s$.
\end{enumerate}
\end{theorem}

\begin{proof}
This is something elementary, the idea being as follows:

\medskip

(1) This follows from definitions.

\medskip

(2) This follows from Theorem 6.10 and Proposition 6.11, because the algebra $C(H_N^{s+})$ is the quotient of the algebra $C(S_{sN}^+)$ by the relations making the fundamental corepresentation commute with the adjacency matrix of $NZ_s$.
\end{proof}

Next in line, we must talk about wreath products. We have here:

\index{reflection group}
\index{quantum reflection group}

\begin{theorem}
We have the following results:
\begin{enumerate}
\item $H_N^s=\mathbb Z_s\wr S_N$.

\item $H_N^{s+}=\mathbb Z_s\wr_*S_N^+$.
\end{enumerate}
\end{theorem}

\begin{proof}
This follows from the following formulae, valid for any connected graph $X$, and explained before, in chapter 5, applied to the graph $Z_s$:
$$G(NX)=G(X)\wr S_N\quad,\quad 
G^+(NX)=G^+(X)\wr_*S_N^+$$

Alternatively, (1) follows from definitions, and (2) can be proved directly, by constructing a pair of inverse morphisms. For details here, we refer to \cite{bve}.
\end{proof}

Regarding now the easiness property of $H_N^s,H_N^{s+}$, we already know that this happens at $s=1,2$. The point is that this happens at $s=\infty$ too, the result being as follows:

\index{complex reflection group}

\begin{theorem}
The quantum groups $K_N,K_N^+$ are easy, the corresponding categories
$$\mathcal P_{even}\subset P\quad,\quad 
\mathcal{NC}_{even}\subset NC$$
consisting of the partitions satisfying $\#\circ=\#\bullet$, as a weighted equality, in each block.
\end{theorem}

\begin{proof}
This is something which is routine, along the lines of the proof of Theorem 6.3, and for full details here, we refer to the paper \cite{bb+}.
\end{proof}

More generally now, we have the following result, from \cite{bb+}:

\begin{theorem}
The quantum groups $H_N^s,H_N^{s+}$ are easy, the corresponding categories
$$P^s\subset P\quad,\quad 
NC^s\subset NC$$
consisting of partitions satisfying $\#\circ=\#\bullet(s)$, as a weighted sum, in each block.
\end{theorem}

\begin{proof}
Observe that the result holds at $s=1$, trivially, then at $s=2$ as well, where our condition is equivalent to $\#\circ=\#\bullet(2)$ in each block, as found in Theorem 6.3, and finally at $s=\infty$ too, as explained in Theorem 6.14. In general, this follows as in the case of $H_N,H_N^+$, by using the one-block partition in $P(s,s)$. See \cite{bb+}.
\end{proof}

The above proofs were of course quite brief, but we will be back to this, with details, in the next section, with full proofs for certain more general results.

\section*{6c. Fusion rules}

Let us discuss now, following \cite{bve}, the classification of the irreducible representations of $H_N^{s+}$, and the computation of their fusion rules. For this purpose, let us go back to the elements $u_{ij},p_{ij}$ in Definition 6.6. We recall that, as a consequence of Proposition 6.8, the matrix $p=(p_{ij})$ is a magic unitary. We first have the following result:

\begin{proposition}
The elements $u_{ij}$ and $p_{ij}$ satisfy:
\begin{enumerate}
\item $p_{ij}u_{ij}=u_{ij}$.

\item $u_{ij}^*=u_{ij}^{s-1}$.

\item $u_{ij}u_{ik}=0$ for $j\neq k$.
\end{enumerate}
\end{proposition}

\begin{proof}
We use the fact that in a $C^*$-algebra, $aa^*=0$ implies $a=0$.

\medskip

(1) This follows from the following computation, with $a=(p_{ij}-1)u_{ij}$:
$$aa^*
=(p_{ij}-1)p_{ij}(p_{ij}-1)
=0$$

(2) With $a=u_{ij}^*-u_{ij}^{s-1}$ we have $aa^*=0$, which gives the result.

\medskip

(3) With $a=u_{ij}u_{ik}$ we have $aa^*=0$, which gives the result.
\end{proof}

In what follows, we make the convention $u_{ij}^0=p_{ij}$. We have then:

\begin{proposition}
The algebra $C(H_N^{s+})$ has a family of $N$-dimensional corepresentations $\{u_k|k\in\mathbb Z\}$, satisfying the following conditions:
\begin{enumerate}
\item $u_k=(u_{ij}^k)$ for any $k\geq 0$.

\item $u_k=u_{k+s}$ for any $k\in\mathbb Z$.

\item $\bar{u}_k=u_{-k}$ for any $k\in\mathbb Z$.
\end{enumerate}
\end{proposition}

\begin{proof}
This is something elementary, the idea being as follows:

\medskip

(1) Let us set $u_k=(u_{ij}^k)$. By using Proposition 6.16 (3), we have:
$$\Delta(u_{ij}^k)
=\sum_{l_1\ldots l_k}u_{il_1}\ldots u_{il_k}\otimes u_{l_1j}\ldots u_{l_kj}
=\sum_lu_{il}^k\otimes u_{lj}^k$$

We have as well, trivially, the following two formulae:
$$\varepsilon(u_{ij}^k)=\delta_{ij}\quad,\quad 
S(u_{ij}^k)=u_{ji}^{*k}$$

(2) This follows once again from Proposition 6.16 (3), as follows:
$$u_{ij}^{k+s}
=u_{ij}^ku_{ij}^s
=u_{ij}^kp_{ij}
=u_{ij}^k$$

(3) This follows from Proposition 6.16 (2), and we are done.
\end{proof}

Let us compute now the intertwiners between the various tensor products between the above corepresentations $u_i$. For this purpose, we make the assumption $N\geq 4$, which brings linear independence. In order to simplify the notations, we will use:

\begin{definition}
For $i_1,\ldots,i_k\in\mathbb Z$ we use the notation
$$u_{i_1\ldots i_k}=u_{i_1}\otimes\ldots\otimes u_{i_k}$$
where $\{u_i|i\in\mathbb Z\}$ are the corepresentations in Proposition 6.17.
\end{definition}

Observe that in the particular case $i_1,\ldots,i_k\in\{\pm 1\}$, we obtain in this way all the possible tensor products between $u=u_1$ and $\bar{u}=u_{-1}$, known by \cite{wo1} to contain any irreducible corepresentation of $C(H_N^{s+})$. Here is now our main result:

\begin{theorem}
We have the following equality of linear spaces,
$$Hom(u_{i_1\ldots i_k},u_{j_1\ldots j_l})=span\left\{T_p\Big|p\in NC_s(i_1\ldots i_k,j_1\ldots j_l)\right\}$$
where the set on the right consists of elements of $NC(k,l)$ having the property that in each block, the sum of $i$ indices equals the sum of $j$ indices, modulo $s$.
\end{theorem}

\begin{proof}
This result is from \cite{bve}, the idea of the proof being as follows:

\medskip

(1) Our first claim is that, in order to prove $\supset$, we may restrict attention to the case $k=0$. This follow indeed from the Frobenius duality isomorphism.

\medskip

(2) Our second claim is that, in order to prove $\supset$ in the case $k=0$, we may restrict attention to the one-block partitions. Indeed, this follows once again from a standard trick. Consider the following disjoint union:
$$NC_s=\bigcup_{k=0}^\infty\bigcup_{i_1\ldots i_k} NC_s(0,i_1\ldots i_k)$$

This is a set of labeled partitions, having property that each $p\in NC_s$ is noncrossing, and that for $p\in NC_s$, any block of $p$ is in $NC_s$. But it is well-known that under these assumptions, the global algebraic properties of $NC_s$ can be checked on blocks.

\medskip

(3) Proof of $\supset$. According to the above considerations, we just have to prove that the vector associated to the one-block partition in $NC(l)$ is fixed by $u_{j_1\ldots j_l}$, when:
$$s|j_1+\ldots+j_l$$

Consider the standard generators $e_{ab}\in M_N(\mathbb C)$, acting on the basis vectors by:
$$e_{ab}(e_c)=\delta_{bc}e_a$$

The corepresentation $u_{j_1\ldots j_l}$ is given by the following formula:
$$u_{j_1\ldots j_l}=\sum_{a_1\ldots a_l}\sum_{b_1\ldots b_l}u_{a_1b_1}^{j_1}\ldots u_{a_lb_l}^{j_l}\otimes e_{a_1b_1}\otimes\ldots\otimes e_{a_lb_l}$$

As for the vector associated to the one-block partition, this is:
$$\xi_l=\sum_be_b^{\otimes l}$$

By using now several times the relations in Proposition 6.16, we obtain, as claimed: 
\begin{eqnarray*}
u_{j_1\ldots j_l}(1\otimes\xi_l)
&=&\sum_{a_1\ldots a_l}\sum_bu_{a_1b}^{j_1}\ldots u_{a_lb}^{j_l}\otimes e_{a_1}\otimes\ldots\otimes e_{a_l}\\
&=&\sum_{ab}u_{ab}^{j_1+\ldots+j_l}\otimes e_a^{\otimes l}\\
&=&1\otimes\xi_l
\end{eqnarray*}

(4) Proof of $\subset$. The spaces in the statement form a Tannakian category, so they correspond to a Woronowicz algebra $A$, coming with corepresentations $\{v_i\}$, such that:
$$Hom(v_{i_1\ldots i_k},v_{j_1\ldots j_l})={\rm span}\left\{T_p\Big|p\in NC_s(i_1\ldots i_k,j_1\ldots j_l)\right\}$$

On the other hand, the inclusion $\supset$ that we just proved shows that $C(H_N^{s+})$ is a model for the category. Thus we have a quotient map as follows:
$$A\to C(H_N^{s+})\quad,\quad 
v_i\to u_i$$

But this latter map can be shown to be an isomorphism, by suitably adapting the proof from the $s=1$ case, for the quantum permutation group $S_N^+$. See \cite{bb+}, \cite{bve}.
\end{proof}

As an illustration for the above result, we have the following statement:

\begin{proposition}
The basic corepresentations $u_0,\ldots,u_{s-1}$ are as follows:
\begin{enumerate}
\item $u_1,\ldots,u_{s-1}$ are irreducible.

\item $u_0=1+r_0$, with $r_0$ irreducible.

\item $r_0,u_1,\ldots,u_{s-1}$ are distinct.
\end{enumerate}
\end{proposition}

\begin{proof}
We apply Theorem 6.19 with $k=l=1$ and $i_1=i,j_1=j$. This gives:
$$\dim (Hom(u_i,u_j))=\# NC_s(i,j)$$

We have two candidates for the elements of $NC_s(i,j)$, namely the two partitions in $NC(1,1)$. So, consider these two partitions, with the points labeled by $i,j$:
$$p=\left\{\begin{matrix}i\cr\Big\vert\cr j\end{matrix}\right\}\qquad\ \qquad
q=\left\{\begin{matrix}i\cr|\cr\cr |\cr j\end{matrix}\right\}$$

We have to check for each of these partitions if the sum of $i$ indices equals or not the sum of $j$ indices, modulo $s$, in each block. The answer is as follows:
\begin{eqnarray*}
p\in NC_s(i,j)&\iff&i=j\\
q\in NC_s(i,j)&\iff&i=j=0
\end{eqnarray*}

By collecting together these two answers, we obtain:
$$\# NC_s(i,j)=
\begin{cases}
0&{\rm if\ }i\neq j\\
1&{\rm if\ }i=j\neq 0\\
2&{\rm if\ }i=j=0
\end{cases}$$

We can now prove the various assertions, as follows:

\medskip

(1) This follows from the second equality.

\medskip

(2) This follows from the third equality and from the fact that we have $1\in u_s$.

\medskip

(3) This follows from the first equality.
\end{proof}

We can now compute the fusion rules for $H_N^{s+}$. The result, from \cite{bve}, is as follows:

\index{complex reflection group}
\index{fusion rules}

\begin{theorem}
Let $F=<\mathbb Z_s>$ be the set of words over $\mathbb Z_s$, with involution given by $(i_1\ldots i_k)^-=(-i_k)\ldots(-i_1)$, and with fusion product given by:
$$(i_1\ldots i_k)\cdot (j_1\ldots j_l)=i_1\ldots i_{k-1}(i_k+j_1)j_2\ldots j_l$$
The irreducible representations of $H_N^{s+}$ can then be labeled $r_x$ with $x\in F$, such that the involution and fusion rules are $\bar{r}_x=r_{\bar{x}}$ and
$$r_x\otimes r_y=\sum_{x=vz,y=\bar{z}w}r_{vw}+r_{v\cdot w}$$
and such that we have $r_i=u_i-\delta_{i0}1$ for any $i\in\mathbb Z_s$.
\end{theorem}

\begin{proof}
This basically follows from Theorem 6.19, the idea being as follows:

\medskip

(1) Consider the monoid $A=\{a_x|x\in F\}$, with multiplication $a_xa_y=a_{xy}$. We denote by $\mathbb NA$ the set of linear combinations of elements in $A$, with coefficients in $\mathbb N$, and we endow it with fusion rules as in the statement:
$$a_x\otimes a_y=\sum_{x=vz,y=\bar{z}w}a_{vw}+a_{v\cdot w}$$

With these notations, $(\mathbb NA,+,\otimes)$ is a semiring. We will use as well the set $\mathbb ZA$, formed by the linear combinations of elements of $A$, with coefficients in $\mathbb Z$. The above tensor product operation extends to $\mathbb ZA$, and $(\mathbb ZA,+,\otimes)$ is a ring.

\medskip

(2) Our claim is that the fusion rules on $\mathbb ZA$ can be uniquely described by conversion formulae as follows, with $C$ being positive integers, and $D$ being integers:
$$a_{i_1}\otimes\ldots\otimes a_{i_k}=\sum_l\sum_{j_1\ldots j_l}C_{i_1\ldots i_k}^{j_1\ldots j_l}a_{j_1\ldots j_l}$$
$$a_{i_1\ldots i_k}=\sum_l\sum_{j_1\ldots j_l}D_{i_1\ldots i_k}^{j_1\ldots j_l}a_{j_1}\otimes\ldots\otimes a_{j_l}$$

Indeed, the existence and uniqueness of such decompositions follow from the definition of the tensor product operation, and by recurrence over $k$ for the $D$ coefficients.

\medskip

(3) Our claim is that there is a unique morphism of rings $\Phi:\mathbb ZA\to R$, such that $\Phi(a_i)=r_i$ for any $i$. Indeed, consider the following elements of $R$:
$$r_{i_1\ldots i_k}=\sum_l\sum_{j_1\ldots j_l}D_{i_1\ldots i_k}^{j_1\ldots j_l}r_{j_1}\otimes\ldots\otimes r_{j_l}$$

In case we have a morphism as claimed, we must have $\Phi(a_x)=r_x$ for any $x\in F$. Thus our morphism is uniquely determined on  $A$, so it is uniquely determined on $\mathbb ZA$. In order to prove now the existence, we can set $\Phi(a_x)=r_x$ for any $x\in F$, then extend $\Phi$ by linearity to the whole $\mathbb ZA$. Since $\Phi$ commutes with the above conversion formulae, which describe the fusion rules, it is indeed a morphism. 

\medskip

(4) Our claim is that $\Phi$ commutes with the linear forms $x\to\#(1\in x)$. Indeed, by linearity we just have to check the following equality:
$$\#(1\in a_{i_1}\otimes\ldots\otimes a_{i_k})=\#(1\in r_{i_1}\otimes\ldots\otimes r_{i_k})$$

Now remember that the elements $r_i$ are defined as $r_i=u_i-\delta_{i0}1$. So, consider the elements $c_i=a_i+\delta_{i0}1$. Since the operations $r_i\to u_i$ and $a_i\to c_i$ are of the same nature, by linearity the above formula is equivalent to:
$$\#(1\in c_{i_1}\otimes\ldots\otimes c_{i_k})=\#(1\in u_{i_1}\otimes\ldots\otimes u_{i_k})$$

Now by using Theorem 6.19, what we have to prove is:
$$\#(1\in c_{i_1}\otimes\ldots\otimes c_{i_k})=\#NC_s(i_1\ldots i_k)$$

In order to prove this formula, consider the product on the left:
$$P=(a_{i_1}+\delta_{i_10}1)\otimes(a_{i_2}+\delta_{i_20}1)\otimes\ldots\otimes (a_{i_k}+\delta_{i_k0}1)$$

This quantity can be computed by using the fusion rules on $A$. A recurrence on $k$ shows that the final components of type $a_x$ will come from the different ways of grouping and summing the consecutive terms of the sequence $(i_1,\ldots,i_k)$, and removing some of the sums which vanish modulo $s$, as to obtain the sequence $x$. But this can be encoded by families of noncrossing partitions, and in particular the 1 components will come from the partitions in $NC_s(i_1\ldots i_k)$. Thus $\#(1\in P)=\# NC_s(i_1\ldots i_k)$, as claimed.

\medskip

(5) Our claim now is that $\Phi$ is injective. Indeed, this follows from the result in the previous step, by using a standard positivity argument, namely:
\begin{eqnarray*}
\Phi(\alpha)=0
&\implies&\Phi(\alpha\alpha^*)=0\\
&\implies&\#(1\in \Phi(\alpha\alpha^*))=0\\
&\implies&\#(1\in \alpha\alpha^*)=0\\
&\implies&\alpha=0
\end{eqnarray*}

Here $\alpha$ is arbitrary in the domain of $\Phi$, we use the notation $a_x^*=a_{\bar{x}}$, where $a\to\#(1,a)$ is the unique linear extension of the operation consisting of counting the number of 1's.  Observe that this latter linear form is indeed positive definite, according to the identity $\#(1,a_xa_y^*)=\delta_{xy}$, which is clear from the definition of the product of $\mathbb ZA$.

\medskip

(6) Our claim is that $\Phi(A)\subset R_{irr}$. This is the same as saying that $r_x\in R_{irr}$ for any $x\in F$, and we will prove it by recurrence. Assume that the assertion is true for all the words of length $<k$, and consider an arbitrary length $k$ word, $x=i_1\ldots i_k$. We have:
$$a_{i_1}\otimes a_{i_2\ldots i_k}=a_x+a_{i_1+i_2,i_3\ldots i_k}+\delta_{i_1+i_2,0}a_{i_3\ldots i_k}$$

By applying $\Phi$ to this decomposition, we obtain:
$$r_{i_1}\otimes r_{i_2\ldots i_k}=r_x+r_{i_1+i_2,i_3\ldots i_k}+\delta_{i_1+i_2,0}r_{i_3\ldots i_k}$$

We have the following computation, which is valid for $y=i_1+i_2,i_3\ldots i_k$, as well as for $y=i_3\ldots i_k$ in the case $i_1+i_2=0$:
\begin{eqnarray*}
\#(r_y\in r_{i_1}\otimes r_{i_2\ldots i_k})
&=&\#(1,r_{\bar{y}}\otimes r_{i_1}\otimes r_{i_2\ldots i_k})\\
&=&\#(1,a_{\bar{y}}\otimes a_{i_1}\otimes a_{i_2\ldots i_k})\\
&=&\#(a_y\in a_{i_1}\otimes a_{i_2\ldots i_k})\\
&=&1  
\end{eqnarray*}

Moreover, we know from the previous step that we have $r_{i_1+i_2,i_3\ldots i_k}\neq r_{i_3\ldots i_k}$, so we conclude that the following formula defines an element of $R^+$:
$$\alpha=r_{i_1}\otimes r_{i_2\ldots i_k}-r_{i_1+i_2,i_3\ldots i_k}-\delta_{i_1+i_2,0}r_{i_3\ldots i_k}$$

On the other hand, we have $\alpha=r_x$, so we conclude that we have $r_x\in R^+$. Finally, the irreducibility of $r_x$ follows from the following computation:
\begin{eqnarray*}
\#(1\in r_x\otimes\bar{r}_x)
&=&\#(1\in r_x\otimes r_{\bar{x}})\\
&=&\#(1\in a_x\otimes a_{\bar{x}})\\
&=&\#(1\in a_x\otimes\bar{a}_x)\\
&=&1
\end{eqnarray*}

(7) Summarizing, we have constructed an injective ring morphism, as follows:
$$\Phi:\mathbb ZA\to R\quad,\quad 
\Phi(A)\subset R_{irr}$$

The remaining fact to be proved, namely that we have $\Phi(A)=R_{irr}$, is clear from the general results in \cite{wo1}. Indeed, since each element of $\mathbb NA$ is a sum of elements in $A$, by applying $\Phi$ we get that each element in $\Phi(\mathbb NA)$ is a sum of irreducible corepresentations in $\Phi(A)$. But since $\Phi(\mathbb NA)$ contains all the tensor powers between the fundamental corepresentation and its conjugate, we get $\Phi(A)=R_{irr}$, and we are done.
\end{proof}

Still following \cite{bve}, let us present now a useful formulation of Theorem 6.21. We begin with a slight modification of Theorem 6.21, as follows:

\begin{theorem}
Consider the free monoid $A=<a_i|i\in\mathbb Z_s>$ with the involution $a_i^*=a_{-i}$, and define inductively the following fusion rules on it:
$$pa_i\otimes a_jq=pa_ia_jq+pa_{i+j}q+\delta_{i+j,0}p\otimes q$$
Then the irreducible representations of $H_N^{s+}$ can be indexed by the elements of $A$, and the fusion rules and involution are the above ones.
\end{theorem}

\begin{proof}
Our claim is that this follows from Theorem 6.21, by performing the following relabeling of the irreducible corepresentations:
$$r_{i_1\ldots i_k}\to a_{i_1}\ldots a_{i_k}$$

Indeed, with the notations in Theorem 6.21 we have the following computation, valid for any two elements $i,j\in\mathbb Z_s$ and any two words $x,y\in F$:
\begin{eqnarray*}
r_{xi}\otimes r_{jy}
&=&\sum_{xi=vz,jy=\bar{z}w}r_{vw}+r_{v\cdot w}\\
&=&r_{xijy}+r_{x,i+j,y}+\delta_{i+j,0}\sum_{x=vz,y=\bar{z}w}r_{vw}+r_{v\cdot w}\\
&=&r_{xijy}+r_{x,i+j,y}+\delta_{i+j,0}r_x\otimes r_y
\end{eqnarray*}

With the above relabeling $r_{i_1\ldots i_k}\to a_{i_1}\ldots a_{i_k}$, this gives the formula in the statement (with $r_x\to p$ and $r_y\to q$), and we are done.
\end{proof}

Our alternative reformulation of Theorem 6.21 is based on the idea of embedding $R^+$ into a bigger fusion semiring. Given a fusion monoid $M$ and an element $b\in M$, we denote by $<b>$ the fusion monoid generated by $b$. That is, $<b>$ is the smallest subset of $M$ containing $b$, which is stable by composition, involution and fusion rules. We have:

\begin{theorem}
Consider the monoid $M=<a,z|z^s=1>$ with the involution $a^*=a,z^*=z^{-1}$, and define inductively the following fusion rules on it:
$$vaz^i\otimes z^jaw=vaz^{i+j}aw+\delta_{s|i+j}v\otimes w$$
Then the irreducible representations of $H_N^{s+}$ can be indexed by the elements of the monoid $N=<aza>$, and the fusion rules and involution are the above ones.
\end{theorem}

\begin{proof}
It is routine to check that the elements $az^ia$ with $i=1,\ldots,s$ are free inside $M$. In other words, the submonoid $N'=<az^ia>$ is free on $s$ generators, so it can be identified with the free monoid $A$ in Theorem 6.22, via $a_i=az^ia$. We have $(az^ia)^*=az^{-i}a$, so this identification is involution-preserving. Consider now two arbitrary elements $p,q\in N'$. By using twice the formula in the statement, we obtain:
\begin{eqnarray*}
pa_i\otimes a_jq
&=&paz^ia\otimes az^jaq\\
&=&paz^iaaz^jaq+paz^i\otimes z^jaq\\
&=&paz^iaaz^jaq+paz^{i+j}aq+\delta_{i+j,0}p\otimes q\\
&=&pa_ia_jq+pa_{i+j}q+\delta_{i+j,0}p\otimes q
\end{eqnarray*}

Thus our identification $N'\simeq A$ is fusion rule-preserving. In order to conclude, it remains to prove that the inclusion $N\subset N'$ is actually an equality. But this follows from the fact that $A$ is generated as a fusion monoid by $a_1$. Indeed, by using the identification $N'\simeq A$ this shows that $N'$ is generated as a fusion monoid by $aza$, and we are done. 
\end{proof}

\section*{6d. Bessel laws}

Let us discuss now the computation of the asymptotic laws of characters. We begin with a discussion for $H_N$, from \cite{bbc}, which has its own interest:

\begin{theorem}
The asymptotic law of $\chi_t$ for the group $H_N$ is given by
$$b_t=e^{-t}\sum_{k=-\infty}^\infty\delta_k
\sum_{p=0}^\infty \frac{(t/2)^{|k|+2p}}{(|k|+p)!p!}$$ where
$\delta_k$ is the Dirac mass at $k\in\mathbb Z$.
\end{theorem}

\begin{proof}
We regard the group $H_N$ as being the symmetry group of the graph $I_N=\{I^1,\ldots ,I^N\}$ formed by $N$ segments. The diagonal coefficients are then:
$$u_{ii}(g)=\begin{cases}
\ 0\ \mbox{ if $g$ moves $I^i$}\\
\ 1\ \mbox{ if $g$ fixes $I^i$}\\
-1\mbox{ if $g$ returns $I^i$}
\end{cases}$$

Let $s=[tN]$, and denote by $\uparrow g,\downarrow g$ the number of segments among $\{I^1,\ldots ,I^s\}$ which are fixed, respectively returned by an element $g\in H_N$. With this notation, we have:
$$u_{11}+\ldots+u_{ss}=\uparrow g-\downarrow g$$

We denote by $P_N$ probabilities computed over the group $H_N$. The density of the law of $u_{11}+\ldots+u_{ss}$ at a point $k\geq 0$ is given by the following formula:
\begin{eqnarray*}
D(k)
&=&P_N(\uparrow g-\downarrow g=k)\\
&=&\sum_{p=0}^\infty P_N(\uparrow g=k+p, \downarrow g=p)
\end{eqnarray*}

At $t=1$, since the probability for $\sigma\in S_N$ to have no fixed points is asymptotically $P_0=\frac{1}{e}$, the probability of $\sigma\in S_N$ to have $m$ fixed points is asymptotically:
$$P_m=\frac{1}{em!}$$

In terms of probabilities over $H_N$, we obtain, as desired:
\begin{eqnarray*}
\lim_{N\to\infty}D(k)
&=&\lim_{N\to\infty}\sum_{p=0}^\infty(1/2)^{k+2p}\begin{pmatrix}k+2p\\ k+p\end{pmatrix} P_N(\uparrow g+\downarrow g=k+2p)\\
&=&\sum_{p=0}^\infty(1/2)^{k+2p}\begin{pmatrix}k+2p\\ k+p\end{pmatrix}\frac{1}{e(k+2p)!}\\ 
&=&\frac{1}{e}\sum_{p=0}^\infty\frac{(1/2)^{k+2p}}{(k+p)!p!}
\end{eqnarray*}

The general case $t\in(0,1]$ follows by performing some modifications in the above computation. The asymptotic density is computed as follows:
\begin{eqnarray*}
\lim_{N\to\infty}D(k)
&=&\lim_{N\to\infty}\sum_{p=0}^\infty(1/2)^{k+2p}\begin{pmatrix}k+2p\\ k+p\end{pmatrix} P_N(\uparrow g+\downarrow g=k+2p)\\ \
&=&\sum_{p=0}^\infty(1/2)^{k+2p}\begin{pmatrix}k+2p\\ k+p\end{pmatrix}\frac{t^{k+2p}}{e^t(k+2p)!}\\
&=&e^{-t}\sum_{p=0}^\infty\frac{(t/2)^{k+2p}}{(k+p)!p!}
\end{eqnarray*}

On the other hand, we have $D(-k)=D(k)$, so we obtain the result.
\end{proof}

Next, we have the following result, once again from \cite{bbc}:

\index{Bessel law}

\begin{theorem}
The Bessel laws $b_t$ have the additivity property
$$b_s*b_t=b_{s+t}$$
so they form a truncated one-parameter semigroup with respect to convolution.
\end{theorem}

\begin{proof}
The Fourier transform of $b_t$ is given by:
$$Fb_t(y)=e^{-t}\sum_{k=-\infty}^\infty e^{ky}\,f_k(t/2)$$

We compute now the derivative with respect to $t$:
$$Fb_t(y)'=-Fb_t(y)+\frac{e^{-t}}{2}\sum_{k=-\infty}^\infty e^{ky}\,f_k'(t/2)$$

On the other hand, the derivative of $f_k$ with $k\geq 1$ is given by:
\begin{eqnarray*}
f_k'(t)
&=&\sum_{p=0}^\infty\frac{(k+2p)t^{k+2p-1}}{(k+p)!p!}\\
&=&\sum_{p=0}^\infty\frac{(k+p)t^{k+2p-1}}{(k+p)!p!}
+\sum_{p=0}^\infty\frac{p\,t^{k+2p-1}}{(k+p)!p!}\\
&=&\sum_{p=0}^\infty \frac{t^{k+2p-1}}{(k+p-1)!p!}
+\sum_{p=1}^\infty\frac{t^{k+2p-1}}{(k+p)!(p-1)!}\\
&=&\sum_{p=0}^\infty \frac{t^{(k-1)+2p}}{((k-1)+p)!p!}
+\sum_{p=1}^\infty\frac{t^{(k+1)+2(p-1)}}{((k+1)+(p-1))!
(p-1)!}\cr &=&f_{k-1}(t)+f_{k+1}(t)
\end{eqnarray*}

This computation works in fact for any $k$, so we get:
\begin{eqnarray*}
Fb_t(y)'
&=&-Fb_t(y)+\frac{e^{-t}}{2}
\sum_{k=-\infty}^\infty e^{ky} (f_{k-1}(t/2)+f_{k+1}(t/2))\\
&=&-Fb_t(y)+\frac{e^{-t}}{2} \sum_{k=-\infty}^\infty
e^{(k+1)y}f_{k}(t/2)+e^{(k-1)y}f_{k}(t/2)\\
&=&-Fb_t(y)+\frac{e^{y}+e^{-y}}{2}\,Fb_t(y)\\
&=&\left(\frac{e^{y}+e^{-y}}{2}-1\right)Fb_t(y)
\end{eqnarray*}

Thus the log of the Fourier transform is linear in $t$, and we get the assertion.
\end{proof}

In order to discuss now the free analogues $\beta_t$ of the above measures $b_t$, and then the $s$-analogues $b^s_t,\beta^s_t$ of the measures $b_t,\beta_t$, we need some free probability. We have the following notion, extending the Poisson limit theory from chapter 1:

\index{compound Poisson limit}
\index{compound PLT}
\index{CPLT}

\begin{definition}
Associated to any compactly supported positive measure $\rho$ are the probability measures
$$p_\rho=\lim_{n\to\infty}\left(\left(1-\frac{c}{n}\right)\delta_0+\frac{1}{n}\rho\right)^{*n}$$
$$\pi_\rho=\lim_{n\to\infty}\left(\left(1-\frac{c}{n}\right)\delta_0+\frac{1}{n}\rho\right)^{\boxplus n}$$
where $c=mass(\rho)$, called compound Poisson and compound free Poisson laws.
\end{definition}

In what follows we will be interested in the case where $\rho$ is discrete, as is for instance the case for $\rho=t\delta_1$ with $t>0$, which produces the Poisson and free Poisson laws. The following result allows us to detect compound Poisson/free Poisson laws:

\begin{proposition}
For $\rho=\sum_{i=1}^sc_i\delta_{z_i}$ with $c_i>0$ and $z_i\in\mathbb R$  we have
$$F_{p_\rho}(y)=\exp\left(\sum_{i=1}^sc_i(e^{iyz_i}-1)\right)$$
$$R_{\pi_\rho}(y)=\sum_{i=1}^s\frac{c_iz_i}{1-yz_i}$$
where $F,R$ denote respectively the Fourier transform, and Voiculescu's $R$-transform.
\end{proposition}

\begin{proof}
Let $\mu_n$ be the measure in Definition 6.26, under the convolution signs:
$$\mu_n=\left(1-\frac{c}{n}\right)\delta_0+\frac{1}{n}$$

In the classical case, we have the following computation:
\begin{eqnarray*}
F_{\mu_n}(y)=\left(1-\frac{c}{n}\right)+\frac{1}{n}\sum_{i=1}^sc_ie^{iyz_i}
&\implies&F_{\mu_n^{*n}}(y)=\left(\left(1-\frac{c}{n}\right)+\frac{1}{n}\sum_{i=1}^sc_ie^{iyz_i}\right)^n\\
&\implies&F_{p_\rho}(y)=\exp\left(\sum_{i=1}^sc_i(e^{iyz_i}-1)\right)
\end{eqnarray*}

In the free case now, we use a similar method. The Cauchy transform of $\mu_n$ is:
$$G_{\mu_n}(\xi)=\left(1-\frac{c}{n}\right)\frac{1}{\xi}+\frac{1}{n}\sum_{i=1}^s\frac{c_i}{\xi-z_i}$$

Consider now the $R$-transform of the measure $\mu_n^{\boxplus n}$, which is given by:
$$R_{\mu_n^{\boxplus n}}(y)=nR_{\mu_n}(y)$$

The above formula of $G_{\mu_n}$ shows that the equation for $R=R_{\mu_n^{\boxplus n}}$ is as follows:
\begin{eqnarray*}
&&\left(1-\frac{c}{n}\right)\frac{1}{y^{-1}+R/n}+\frac{1}{n}\sum_{i=1}^s\frac{c_i}{y^{-1}+R/n-z_i}=y\\
&\implies&\left(1-\frac{c}{n}\right)\frac{1}{1+yR/n}+\frac{1}{n}\sum_{i=1}^s\frac{c_i}{1+yR/n-yz_i}=1
\end{eqnarray*}

Now by multiplying by $n$, rearranging the terms, and letting $n\to\infty$, we get:
\begin{eqnarray*}
\frac{c+yR}{1+yR/n}=\sum_{i=1}^s\frac{c_i}{1+yR/n-yz_i}
&\implies&c+yR_{\pi_\rho}(y)=\sum_{i=1}^s\frac{c_i}{1-yz_i}\\
&\implies&R_{\pi_\rho}(y)=\sum_{i=1}^s\frac{c_iz_i}{1-yz_i}
\end{eqnarray*}

This finishes the proof in the free case, and we are done.
\end{proof}

We have the following result, providing an alternative to Definition 6.26, and which will be our formulation here of the Compound Poisson Limit Theorem (CPLT):

\begin{theorem}
For $\rho=\sum_{i=1}^sc_i\delta_{z_i}$ with $c_i>0$ and $z_i\in\mathbb R$ we have
$$p_\rho/\pi_\rho={\rm law}\left(\sum_{i=1}^sz_i\alpha_i\right)$$
where the variables $\alpha_i$ are Poisson/free Poisson$(c_i)$, independent/free.
\end{theorem}

\begin{proof}
Let $\alpha$ be the sum of Poisson/free Poisson variables in the statement. We will show that the Fourier/$R$-transform of $\alpha$ is given by the formulae in Proposition 6.27. Indeed, by using some well-known Fourier transform formulae, we have:
\begin{eqnarray*}
F_{\alpha_i}(y)=\exp(c_i(e^{iy}-1))
&\implies&F_{z_i\alpha_i}(y)=\exp(c_i(e^{iyz_i}-1))\\
&\implies&F_\alpha(y)=\exp\left(\sum_{i=1}^sc_i(e^{iyz_i}-1)\right)
\end{eqnarray*}

Also, by using some well-known $R$-transform formulae, we have:
\begin{eqnarray*}
R_{\alpha_i}(y)=\frac{c_i}{1-y}
&\implies&R_{z_i\alpha_i}(y)=\frac{c_iz_i}{1-yz_i}\\
&\implies&R_\alpha(y)=\sum_{i=1}^s\frac{c_iz_i}{1-yz_i}
\end{eqnarray*}

Thus we have indeed the same formulae as those in Proposition 6.27.
\end{proof}

Summarizing, we have now a full extension of the basic Poisson limit theory from chapter 1, in the classical and free cases. Back to quantum reflection groups, we have:

\begin{theorem}
The asymptotic laws of truncated characters are as follows, where $\varepsilon_s$ with $s\in\{1,2,\ldots,\infty\}$ is the uniform measure on the $s$-th roots of unity:
\begin{enumerate}
\item For $H_N^s$ we obtain the compound Poisson law $b_t^s=p_{t\varepsilon_s}$.

\item For $H_N^{s+}$ we obtain the compound free Poisson law $\beta_t^s=\pi_{t\varepsilon_s}$.
\end{enumerate}
These measures are in Bercovici-Pata bijection.
\end{theorem}

\begin{proof}
This follows from easiness, and from the Weingarten formula, exactly as for the classical and quantum permutation groups. For details here, we refer to \cite{bb+}.
\end{proof}

The Bessel and free Bessel laws have particularly interesting properties at the parameter values $s=2,\infty$. So, let us record the precise statement here:

\index{Bessel law}
\index{free Bessel law}
\index{Bercovici-Pata bijection}

\begin{theorem}
The asymptotic laws of truncated characters are as follows:
\begin{enumerate}
\item For $H_N$ we obtain the real Bessel law $b_t=p_{t\varepsilon_2}$.

\item For $K_N$ we obtain the complex Bessel law $B_t=p_{t\varepsilon_\infty}$.

\item For $H_N^+$ we obtain the free real Bessel law $\beta_t=\pi_{t\varepsilon_2}$.

\item For $K_N^+$ we obtain the free complex Bessel law $\mathfrak B_t=\pi_{t\varepsilon_\infty}$.
\end{enumerate}
\end{theorem}

\begin{proof}
This follows indeed from Theorem 6.29, at $s=2,\infty$.
\end{proof}

Let us discuss now, as a final topic, the computation of the moments of the free Bessel laws. The idea will be that of expressing these moments in terms of generalized binomial coefficients. We recall that the coefficient corresponding to $\alpha\in\mathbb R$, $k\in\mathbb N$ is:
$$\binom{\alpha}{k}=\frac{\alpha(\alpha-1)\ldots(\alpha-k+1)}{k!}$$

We denote by $m_1,m_2,m_3,\ldots$ the sequence of moments of a given probability measure. With this convention, we first have the following result, from \cite{bb+}:

\index{Fuss-Catalan numbers}

\begin{theorem}
The moments of $\beta^s_1$ with $s>0$ are
$$m_k=\frac{1}{sk+1}\binom{sk+k}{k}$$
which are the Fuss-Catalan numbers.
\end{theorem}

\begin{proof}
In the case $s\in\mathbb N$, we know that we have $m_k=\# NC_s(k)$. The formula in the statement follows then by counting such partitions. In the general case $s>0$, observe first that the Fuss-Catalan number in the statement is a polynomial in $s$:
$$\frac{1}{sk+1}\binom{sk+k}{k}=\frac{(sk+2)(sk+3)\ldots(sk+k)}{k!}$$

Thus, in order to pass from the case $s\in\mathbb N$ to the case $s>0$, it is enough to check that the $k$-th moment of $\pi_{s1}$ is analytic in $s$. But this is clear from the equation  $f=1+zf^{s+1}$ of the Stieltjes transform of $\pi_{s1}$, and this gives the result.
\end{proof}

We have as well the following result, which deals with the general case $t>0$:

\index{Fuss-Narayana numbers}

\begin{theorem}
The moments of $\beta^s_t$ with $s>0$ are
$$m_k=\sum_{b=1}^k\frac{1}{b}\binom{k-1}{b-1}\binom{sk}{b-1}t^b$$
which are the Fuss-Narayana numbers.
\end{theorem}

\begin{proof}
In the case $s\in\mathbb N$, we know from the above that we have the following formula, where $F_{kb}$ is the number of partitions in $NC_s(k)$ having $b$ blocks:
$$m_k=\sum_bF_{kb}t^b$$

With this observation in hand, the formula in the statement follows by counting such partitions, with this count being well-known. This result can be then extended to any parameter $s>0$, by using a standard complex variable argument, as before. See \cite{bb+}.
\end{proof}

In the case $s\notin\mathbb N$, the moments of $\beta^s_t$ can be further expressed in terms of gamma functions. In the case $s=1/2$, the result, also from \cite{bb+}, is as follows:

\begin{theorem}
The moments of $\beta^{1/2}_1$ are given by the following formulae:
$$m_{2p}=\frac{1}{p+1}\binom{3p}{p}\quad,\quad 
m_{2p-1}=\frac{2^{-4p+3}p}{(6p-1)(2p+1)}\cdot\frac{p!(6p)!}{(2p)!(2p)!(3p)!}$$
\end{theorem}

\begin{proof}
According to our various results above, the even moments of the free Bessel law $\beta^s_t$ with $s=n-1/2$, $n\in\mathbb N$, are given by:
$$m_{2p}
=\frac{1}{(n-1/2)(2p)+1}\binom{(n+1/2)2p}{2p}
=\frac{1}{(2n-1)p+1}\binom{(2n+1)p}{2p}$$

With $n=1$ we get the formula in the statement. Now for the odd moments, we can use here the following well-known identity:
$$\begin{pmatrix}m-1/2\cr k\end{pmatrix}=\frac{4^{-k}}{k!}\cdot\frac{(2m)!}{m!}\cdot\frac{(m-k)!}{(2m-2k)!}$$

With $m=2np+p-n$ and $k=2p-1$ we get:
\begin{eqnarray*}
m_{2p-1}
&=&\frac{1}{(n-1/2)(2p-1)+1}\binom{(n+1/2)(2p-1)}{2p-1}\\
&=&\frac{2}{(2n-1)(2p-1)+2}\binom{(2np+p-n)-1/2}{2p-1}\\
&=&\frac{2^{-4p+3}}{(2p-1)!}\cdot\frac{(4np+2p-2n)!}{(2np+p-n)!}\cdot\frac{(2np-p-n+1)!}{(4np-2p-2n+3)!}
\end{eqnarray*}

In particular with $n=1$ we obtain:
\begin{eqnarray*}
m_{2p-1}
&=&\frac{2^{-4p+3}}{(2p-1)!}\cdot\frac{(6p-2)!}{(3p-1)!}\cdot\frac{p!}{(2p+1)!}\\
&=&\frac{2^{-4p+3}(2p)}{(2p)!}\cdot\frac{(6p)!(3p)}{(3p)!(6p-1)6p}\cdot\frac{p!}{(2p)!(2p+1)}
\end{eqnarray*}

But this gives the formula in the statement.
\end{proof}

There are many other interesting things, of both combinatorial and complex analytic nature, that can be said about the free Bessel laws, their moments and their densities, and for a full discussion here, we refer to \cite{bb+}, and subsequent papers. 

\bigskip

Finally, we have the following result, which is of theoretical interest:

\begin{theorem}
The moments of the various central limiting measures, namely
$$\xymatrix@R=45pt@C=45pt{
\beta^s_t\ar@{-}[r]\ar@{-}[d]&\gamma_t\ar@{-}[r]\ar@{-}[d]&\Gamma_t\ar@{-}[d]\\
b^s_t\ar@{-}[r]&g_t\ar@{-}[r]&G_t
}$$
are always given by the same formula, involving partitions, namely
$$M_k=\sum_{\pi\in D(k)}t^{|\pi|}$$
where the sets of partitions $D(k)$ in question are respectively
$$\xymatrix@R=50pt@C=50pt{
NC^s\ar[d]&NC_2\ar[d]\ar[l]&\mathcal{NC}_2\ar[l]\ar[d]\\
P^s&P_2\ar[l]&\mathcal P_2\ar[l]}$$
and where $|.|$ is the number of blocks. 
\end{theorem}

\begin{proof}
This follows indeed by putting together the various moment results that we have, from this chapter and from the previous ones.
\end{proof}

As already mentioned, in what regards the Bessel and free Bessel laws $b^s_t,\beta^s_t$, the important particular cases are $s=1,2,\infty$. It is therefore tempting to leave one of these 3 cases aside, and fold the corresponding diagram into a cube. We obtain:

\begin{theorem}
The moments of the various central limiting measures,
$$\xymatrix@R=20pt@C=22pt{
&\mathfrak B_t\ar@{-}[rr]\ar@{-}[dd]&&\Gamma_t\ar@{-}[dd]\\
\beta_t\ar@{-}[rr]\ar@{-}[dd]\ar@{-}[ur]&&\gamma_t\ar@{-}[dd]\ar@{-}[ur]\\
&B_t\ar@{-}[rr]\ar@{-}[uu]&&G_t\ar@{-}[uu]\\
b_t\ar@{-}[uu]\ar@{-}[ur]\ar@{-}[rr]&&g_t\ar@{-}[uu]\ar@{-}[ur]
}$$
are always given by the same formula, involving partitions, namely
$$M_k=\sum_{\pi\in D(k)}t^{|\pi|}$$
where the sets of partitions $D(k)$ in question are respectively:
$$\xymatrix@R=20pt@C=4pt{
&\mathcal{NC}_{even}\ar[dl]\ar[dd]&&\ \ \ \mathcal{NC}_2\ \ \ \ar[ll]\ar[dd]\ar[dl]\\
NC_{even}\ar[dd]&&NC_2\ar[ll]\ar[dd]\\
&\mathcal P_{even}\ar[dl]&&\mathcal P_2\ar[ll]\ar[dl]\\
P_{even}&&P_2\ar[ll]
}$$
\end{theorem}

\begin{proof}
This follows indeed from the various moment results that we have.
\end{proof}

The above result is quite interesting, with the cubes in question, along with the cube formed by the corresponding quantum groups, being called ``standard cubes''. The point indeed is that all these cubes provide us with some useful 3D intuition, in relation with various aspects of quantum algebra, be them algebraic or probabilistic. So, saved.

\section*{6e. Exercises} 

As a first key exercise, in relation with the above, we have:

\begin{exercise}
Prove that $K_N,K_N^+$, regarded as intermediate quantum groups
$$S_N\subset K_N\subset U_N\quad,\quad 
S_N^+\subset K_N^+\subset U_N^+$$
are indeed easy, the corresponding categories of partitions being as follows,
$$P\supset \mathcal P_{even}\supset\mathcal P_2\quad,\quad 
NC\supset\mathcal{NC}_{even}\supset\mathcal{NC}_2$$
consisting of the partitions satisfying $\#\circ=\#\bullet$, as a weighted equality, in each block.
\end{exercise}

This is something that we briefly discussed in the above, and the problem now is that of working out all the details. Along the same lines, more generally, we have:

\begin{exercise}
Prove that $H_N^s,H_N^{s+}$, regarded as intermediate quantum groups
$$S_N\subset H_N^s\subset U_N\quad,\quad 
S_N^+\subset H_N^{s+}\subset U_N^+$$
are indeed easy, the corresponding categories of partitions being as follows,
$$P\supset P^s\supset\mathcal P_2\quad,\quad 
NC\supset NC^s\supset\mathcal{NC}_2$$
consisting of the partitions satisfying $\#\circ=\#\bullet(s)$, as a weighted equality, in each block.
\end{exercise}

Again, this is something that we discussed in the above, and fully proved at $s=1,2$, and the problem now is that of working out all the details. Assuming that you solved the previous exercise, corresponding to the case $s=\infty$, this should not be very hard.

\begin{exercise}
Reformulate the Schur-Weyl twisting theory from chapter 4 in terms of the intermediate easy quantum groups
$$H_N\subset G\subset U_N^+$$
as a duality between such intermediate quantum groups.
\end{exercise}

The point here is that, as explained in chapter 4, in order to perform the Schur-Weyl twisting operation we need a signature map for the partitions, and this signature map is only defined on $P_{even}$. Thus, we have a link here with $H_N$, and so the whole twisting material from chapter 4 must be now reviewed, by taking this into account.

\begin{exercise}
Prove that the quantum reflection groups $H_N,H_N^+$ equal their own Schur-Weyl twists, and then do the same for $H_N^s,H_N^{s+}$.
\end{exercise}

This is not something very difficult, normally coming from definitions, one subtlety however coming from the fact that for the second question we must assume $s\in2\mathbb N$, in order for the quantum groups $H_N^s,H_N^{s+}$ to be indeed twistable in our sense.

\begin{exercise}
Prove that the intermediate easy quantum groups $H_N\subset G\subset H_N^+$ are subject to a dichotomy, as follows,
$$H_N\subset G\subset H_N^{[\infty]}\qquad\big/\qquad 
H_N^{[\infty]}\subset G\subset H_N^+$$
where $H_N^{[\infty]}$ is a certain suitably chosen such intermediate quantum group.
\end{exercise}

As already mentioned, this is a rather difficult exercise.

\begin{exercise}
Fully classify the intermediate easy quantum groups
$$H_N\subset G\subset H_N^+$$
by solving the $2$ classification problems raised by the previous exercise.
\end{exercise}

As before, this is a difficult exercise, and the problem here is rather that of finding the relevant literature, reading it, and writing down a brief account of that.

\chapter{Twisted reflections}

\section*{7a. Quantum graphs}

We have seen in the previous two chapters that some general theory can be developed for the closed subgroups $G\subset S_N^+$, in particular with a notion of ``quantum reflection group''. We discuss here the twisted extension of these results, and in particular the twisted analogues of the quantum reflection groups, obtained by using generalized quantum permutation groups $S_Z^+$, with $Z$ being an arbitrary finite quantum space. 

\bigskip

In order to discuss this, let us go back to chapter 4, where the finite quantum spaces $Z$ and the quantum groups $S_Z^+$ were introduced. We recall from there that we have:

\begin{definition}
A finite quantum space $Z$ is the abstract dual of a finite dimensional $C^*$-algebra $B$, according to the following formula:
$$C(Z)=B$$
Alternatively, by decomposing the algebra $B$, we have the following formula:
$$C(Z)=M_{n_1}(\mathbb C)\oplus\ldots\oplus M_{n_k}(\mathbb C)$$
With the choice $n_1=\ldots=n_k=1$ we obtain the space $Z=\{1,\ldots,k\}$. Also, with $k=1$ we must have $C(Z)=M_n(\mathbb C)$, and we obtain a quantum space denoted $Z=M_n$.
\end{definition}

In order to do some mathematics on such spaces, the very first observation is that we can talk about the formal number of points of such a space, as follows:
$$|Z|=\dim B$$

Alternatively, by decomposing the algebra $B$ as a sum of matrix algebras, as in Definition 7.1, we have the following formula for the formal number of points:
$$|Z|=n_1^2+\ldots+n_k^2$$

Pictorially, this suggests representing $Z$ as a set of $|Z|$ points in the plane, arranged in squares having sides $n_1,\ldots,n_k$, coming from the matrix blocks of $B$, as follows:
$$\begin{matrix}
\circ&\circ&\circ\\
\circ&\circ&\circ&&\ldots&&\circ&\circ\\
\circ&\circ&\circ&&&&\circ&\circ
\end{matrix}$$

As a second piece of mathematics, we can talk about counting measures, as follows:

\begin{definition}
Given a finite quantum space $Z$, we construct the functional
$$tr:C(Z)\to B(l^2(Z))\to\mathbb C$$
obtained by applying the regular representation, and the normalized matrix trace, and we call it integration with respect to the normalized counting measure on $Z$.
\end{definition}

As an illustration here, for the space $Z=\{1,\ldots,k\}$ we obtain the integration with respect to the normalized counting measure, $\mu(\{i\})=1/k$ for any $i$. Also, for the space $Z=M_n$ we obtain the normalized trace of matrices, $tr=Tr/n$. In general, in terms of the matrix decomposition of $B=C(Z)$, as in Definition 7.1, the formula is:
$$tr(a_1,\ldots,a_k)=\frac{1}{|Z|}\sum_{i=1}^kn_i^2\,tr(a_i)$$

Pictorially, this suggests fine-tuning our previous picture of $Z$, by adding to each point the unnormalized trace of the corresponding element of $B$, as follows:
$$\begin{matrix}
\bullet_{n_1}&\circ_0&\circ_0\\
\circ_0&\bullet_{n_1}&\circ_0&&\ldots&&\bullet_{n_k}&\circ_0\\
\circ_0&\circ_0&\bullet_{n_1}&&&&\circ_0&\bullet_{n_k}
\end{matrix}$$

Here we have represented the points on the diagonals with solid circles, since they are of different nature from the off-diagonal ones, the attached numbers being nonzero. However, this picture is not complete either, and we can do better, as follows:

\begin{definition}
Given a finite quantum space $Z$, coming via a formula of type
$$C(Z)=M_{n_1}(\mathbb C)\oplus\ldots\oplus M_{n_k}(\mathbb C)$$
we use the following equivalent conventions for drawing $Z$:
\begin{enumerate}
\item Triple indices. We represent $Z$ as a set of $N=|Z|$ points, with each point being decorated with a triple index $ija$, coming from the standard basis $\{e_{ij}^a\}\subset B$. 

\item Double indices. As before, but by ignoring the index $a$, with the convention that $i,j$ belong to various indexing sets, one for each of the matrix blocks of $B$.

\item Single indices. As before, but with each point being now decorated with a single index, playing the role of the previous triple indices $ija$, or double indices $ij$.
\end{enumerate}
\end{definition}

All the above conventions are useful, and in practice, we will be mostly using the single index convention from (3). As an illustration, consider the space $Z=\{1,\ldots,k\}$. According to our single index convention, we can represent this space as a set of $k$ points, decorated by some indices, which must be chosen different. But the obvious choice for these $k$ different indices is $1,\ldots,k$, and we are led to the following picture:
$$\bullet_1\quad\bullet_2\quad\ldots\quad\bullet_k$$

As another illustration, consider the space $Z=M_n$. Here the picture is as follows, using double indices, which can be regarded as well as being single indices:
$$\begin{matrix}
\bullet_{11}&\circ_{12}&\circ_{13}\\
\circ_{21}&\bullet_{22}&\circ_{23}\\
\circ_{31}&\circ_{32}&\bullet_{33}
\end{matrix}$$

\smallskip

As yet another illustration, for the space $Z=M_3\sqcup M_2$, which appears by definition from the algebra $B=M_3(\mathbb C)\oplus M_2(\mathbb C)$, we are in need of triple indices, which can be of course regarded as single indices, in order to label all the points, and the picture is:
$$\begin{matrix}
\bullet_{111}&\circ_{121}&\circ_{131}\\
\circ_{211}&\bullet_{221}&\circ_{231}&&\ &&\bullet_{112}&\circ_{122}\\
\circ_{311}&\circ_{321}&\bullet_{331}&&&&\circ_{212}&\bullet_{222}
\end{matrix}$$

\smallskip

So long for finite quantum spaces $Z$ and their pictures, and we will stop here, because all this looks a bit like reinventing the wheel. As a last piece of preliminaries now, let us recall from chapter 4 that we have the following result, coming from \cite{wa2}:

\begin{theorem}
Associated to any finite quantum space $Z$ are its classical and quantum counting measure-preserving symmetry groups $S_Z\subset S_Z^+$, with $S_Z^+$ being given by
$$C(S_Z^+)=C(U_N^+)\Big/\Big<\mu\in Hom(u^{\otimes2},u),\eta\in Fix(u)\Big>$$
where $N=|Z|$, and $\mu,\eta$ are the multiplication and unit maps of $C(Z)$, and with $S_Z$ being its classical version. Alternatively, the relations defining $S_Z^+$ are
$$\sum_{ij=p}u_{ik}u_{jl}=u_{p,kl}\quad,\quad\sum_{kl=p}u_{ik}u_{jl}=u_{ij,p}$$
$$\sum_{i=\bar{i}}u_{ij}=\delta_{j\bar{j}}\quad,\quad\sum_{j=\bar{j}}u_{ij}=\delta_{i\bar{i}}$$
$$u_{ij}^*=u_{\bar{i}\hskip0.3mm\bar{j}}$$
in single index notation, with the usual multiplication and involution for the indices.
\end{theorem}

\begin{proof}
This is something that we know well from chapter 4, the idea being that the relations in the statement, be them in abstract $\mu,\eta$ formulation, or in their equivalent single index formulation, are what comes out of the coaction axioms. 
\end{proof}

Getting now to the point where we wanted to get, quantum graphs and their symmetries, let us start with the following straightforward extension of the usual notion of finite graph, from \cite{fpi}, obtained by using a finite quantum space as set of vertices:

\index{quantum graph}
\index{finite quantum graph}

\begin{definition}
We call ``finite quantum graph'' a pair of type
$$X=(Z,d)$$
with $Z$ being a finite quantum space, and with $d\in M_N(\mathbb C)$, where $N=|Z|$.
\end{definition}

This is of course something quite general. In the case $Z=\{1,\ldots,N\}$ for instance, what we have here is a directed graph, with the edges $i\to j$ colored by complex numbers $d_{ij}\in\mathbb C$, and with self-edges $i\to i$ allowed too, again colored by numbers $d_{ii}\in\mathbb C$. In the general case, however, where $Z$ is arbitrary, the need for extra conditions of type $d=d^*$, or $d_{ii}=0$, or $d\in M_N(\mathbb R)$, or $d\in M_N(0,1)$ and so on, is not very natural, as we will soon discover, and it is best to use Definition 7.5 as such, with no restrictions on $d$.

\bigskip

In general, a quantum graph can be represented as a colored oriented graph on $\{1,\ldots,N\}$, where $N=|Z|$, with the vertices being decorated by single indices $i$, and with the colors being complex numbers, namely the entries of $d$. This is similar to the formalism from chapter 5, but there is a discussion here in what regards the exact choice of the colors, which are usually irrelevant in connection with our symmetry problematics, and so can be true colors instead of complex numbers. More on this later.

\bigskip

With the above notion in hand, we have the following definition, also from \cite{fpi}:

\begin{definition}
The quantum automorphism group of $X=(Z,d)$ is the subgroup 
$$G^+(X)\subset S_Z^+$$
obtained via the relation $du=ud$, where $u=(u_{ij})$ is the fundamental corepresentation.
\end{definition}

Again, this is something very natural, coming as a continuation of the constructions in chapter 5. We refer to \cite{fpi}, \cite{twa} for more on this notion, and for a number of advanced computations, in relation with free wreath products. At an elementary level, which will be ours for the moment, a first problem is that of working out the basics of the correspondence $X\to G^+(X)$, following \cite{ba3}. And there are several things to be done here, namely simplices, complementation, color independence, multi-simplices, and reflections. 

\bigskip

Let us start with the simplices. As we will soon discover, things are quite tricky here, leading us in particular to the conclusion that the simplex based on an arbitrary finite quantum space $Z$ is not a usual graph, with $d\in M_N(0,1)$ where $N=|Z|$, but rather a sort of ``signed graph'', with $d\in M_N(-1,0,1)$. Let us start our study with:

\index{empty graph}
\index{simplex}

\begin{theorem}
Given a finite quantum space $Z$, we have
$$G^+(Z_{empty})=G^+(Z_{full})=S_Z^+$$
where $Z_{empty}$ is the empty graph on the vertex set $Z$, coming from the matrix $d=0$, and where $Z_{full}$ is the simplex on the vertex set $Z$, coming from the matrix 
$$d=NP_1-1_N$$
where $N=|Z|$, and where $P_1$ is the orthogonal projection on the unit $1\in C(Z)$.
\end{theorem}

\begin{proof}
This is something quite tricky, the idea being as follows:

\medskip

(1) First of all, the formula $G^+(Z_{empty})=S_Z^+$ is clear from definitions, because the commutation of $u$ with the matrix $d=0$ is automatic.

\medskip

(2) Regarding $G^+(Z_{full})=S_Z^+$, let us first discuss the classical case, $Z=\{1,\ldots,N\}$. Here the simplex $Z_{full}$ is the graph having having edges between any two vertices, whose adjacency matrix is $d=\mathbb I_N-1_N$, where $\mathbb I_N$ is the all-1 matrix. The commutation of $u$ with $1_N$ being automatic, and the commutation with $\mathbb I_N$ being automatic too, $u$ being bistochastic, we have $[u,d]=0$, and so $G^+(Z_{full})=S_Z^+$ in this case, as stated.

\medskip

(3) In general, we know from Theorem 7.4 that we have by definition $\eta\in Fix(u)$, with $\eta:\mathbb C\to C(Z)$ being the unit map. Thus we have $P_1\in End(u)$, and so the condition $[u,P_1]=0$ is automatic. Together with the fact that in the classical case we have $\mathbb I_N=NP_1$, this suggests to define the adjacency matrix of the simplex as being $d=NP_1-1_N$, and with this definition, we have indeed $G^+(Z_{full})=S_Z^+$, as claimed.
\end{proof}

Let us study now the simplices $Z_{full}$ found in Theorem 7.7. In the classical case, $Z=\{1,\ldots,N\}$, what we have is of course the usual simplex. However, in the general case things are more mysterious, the first result here being as follows:

\begin{proposition}
The adjacency matrix of the simplex $Z_{full}$, given by definition by $d=NP_1-1_N$, is a matrix $d\in M_N(-1,0,1)$, which can be computed as follows:
\begin{enumerate}
\item In single index notation, $d_{ij}=\delta_{i\bar{i}}\delta_{j\bar{j}}-\delta_{ij}$.

\item In double index notation, $d_{ab,cd}=\delta_{ab}\delta_{cd}-\delta_{ac}\delta_{bd}$.

\item In triple index notation, $d_{abp,cdq}=\delta_{ab}\delta_{cd}-\delta_{ac}\delta_{bd}\delta_{pq}$.
\end{enumerate}
\end{proposition}

\begin{proof}
According to our single index conventions, from Definition 7.3, the adjacency matrix of the simplex is the one in the statement, namely:
\begin{eqnarray*}
d_{ij}
&=&(NP_1-1_N)_{ij}\\
&=&\bar{1}_i1_j-\delta_{ij}\\
&=&\delta_{i\bar{i}}\delta_{j\bar{j}}-\delta_{ij}
\end{eqnarray*}

In double index notation now, with $i=(ab)$ and $j=(cd)$, and $a,b,c,d$ being usual matrix indices, each thought to be attached to the corresponding matrix block of $C(Z)$, the formula that we obtain in the second one in the statement, namely:
\begin{eqnarray*}
d_{ab,cd}
&=&\delta_{ab,ba}\delta_{cd,dc}-\delta_{ab,cd}\\
&=&\delta_{ab}\delta_{cd}-\delta_{ac}\delta_{bd}
\end{eqnarray*}

Finally, in standard triple index notation, $i=(abp)$ and $j=(cdq)$, with $a,b,c,d$ being now usual numeric matrix indices, ranging in $1,2,3,\ldots\,$, and with $p,q$ standing for corresponding blocks of the algebra $C(Z)$, the formula that we obtain is:
\begin{eqnarray*}
d_{abp,cdq}
&=&\delta_{abp,bap}\delta_{cdq,dcq}-\delta_{abp,cdq}\\
&=&\delta_{ab}\delta_{cd}-\delta_{ac}\delta_{bd}\delta_{pq}
\end{eqnarray*}

Thus, we are led to the conclusions in the statement.
\end{proof}

At the level of examples, for $Z=\{1,\ldots,N\}$ the best is to use the above formula (1). The involution on the index set is $\bar{i}=i$, and we obtain, as we should:
$$d_{ij}=1-\delta_{ij}$$

As a more interesting example now, for the quantum space $Z=M_n$, coming by definition via the formula $C(Z)=M_n(\mathbb C)$, the situation is as follows:

\begin{proposition}
The simplex $Z_{full}$ with $Z=M_n$ is as follows: 
\begin{enumerate}
\item The vertices are $n^2$ points in the plane, arranged in square form.

\item Usual edges, worth $1$, are drawn between distinct points on the diagonal.

\item In addition, each off-diagonal point comes with a self-edge, worth $-1$.
\end{enumerate}
\end{proposition}

\begin{proof}
Here the most convenient is to use the double index formula from Proposition 7.8 (2), which tells us that $d$ is as follows, with indices $a,b,c,d\in\{1,\ldots,n\}$:
$$d_{ab,cd}=\delta_{ab}\delta_{cd}-\delta_{ac}\delta_{bd}$$

This quantity can be $-1,0,1$, and the study goes as follows:

\medskip

-- Case $d_{ab,cd}=1$. This can only happen when $\delta_{ab}\delta_{cd}=1$ and $\delta_{ac}\delta_{bd}=0$, corresponding to a formula of type $d_{aa,cc}=0$, with $a\neq c$, and so to the edges in (2).

\medskip

-- Case $d_{ab,cd}=-1$. This can only happen when $\delta_{ab}\delta_{cd}=0$ and $\delta_{ac}\delta_{bd}=1$, corresponding to a formula of type $d_{ab,ab}=0$, with $a\neq b$, and so to the self-edges in (3).
\end{proof}

The above result is quite interesting, and as an illustration, here is the pictorial representation of the simplex $Z_{full}$ on the vertex set $Z=M_3$, with the convention that the solid arrows are worth $-1$, and the dashed arrows are worth 1:
$$\xymatrix@R=15pt@C=15pt
{\bullet\ar@{--}[dr]\ar@/^/@{--}[ddrr]&\circ^\circlearrowright&\circ^\circlearrowright\\
\circ^\circlearrowright&\bullet\ar@{--}[dr]&\circ^\circlearrowright\\
\circ^\circlearrowright&\circ^\circlearrowright&\bullet}$$

More generally, we can in fact compute $Z_{full}$ for any finite quantum space $Z$, with the result here, which will be our final saying on the subject, being as follows:

\begin{theorem}
Consider a finite quantum space $Z$, and write it as follows, according to the decomposition formula $C(Z)=M_{n_1}(\mathbb C)\oplus\ldots\oplus M_{n_k}(\mathbb C)$ for its function algebra:
$$Z=M_{n_1}\sqcup\ldots\sqcup M_{n_k}$$
The simplex $Z_{full}$ is then the classical simplex formed by the points lying on the diagonals of $M_{n_1},\ldots,M_{n_k}$, with self-edges added, each worth $-1$, at the non-diagonal points.
\end{theorem}

\begin{proof}
The study here is quite similar to the one from the proof of Proposition 7.9, but by using this time the triple index formula from Proposition 7.8 (3), namely:
$$d_{abp,cdq}=\delta_{ab}\delta_{cd}-\delta_{ac}\delta_{bd}\delta_{pq}$$

Indeed, this quantity can be $-1,0,1$, and the $1$ case appears precisely as follows, leading to the classical simplex mentioned in the statement:
$$d_{aap,ccq}=1\quad,\quad\forall ap\neq cq$$

As for the remaining $-1$ case, this appears precisely as follows, leading this time to the self-edges worth $-1$, also mentioned in the statement:
$$d_{abp,abp}=1\quad,\quad\forall a\neq b$$

Thus, we are led to the conclusion in the statement.
\end{proof}

As an illustration, here is the simplex on the vertex set $Z=M_3\sqcup M_2$, with again the convention that the solid arrows are worth $-1$, and the dashed arrows are worth 1:
$$\xymatrix@R=18pt@C=20pt
{\bullet\ar@{--}[dr]\ar@/^/@{--}[ddrr]
\ar@{--}[drrrr]\ar@/^/@{--}[ddrrrrr]
&\circ^\circlearrowright&\circ^\circlearrowright\\
\circ^\circlearrowright&\bullet\ar@{--}[dr]\ar@/_/@{--}[rrr]\ar@{--}[drrrr]
&\circ^\circlearrowright&&\bullet\ar@{--}[dr]&\circ^\circlearrowright\\
\circ^\circlearrowright&\circ^\circlearrowright&\bullet\ar@{--}[urr]\ar@/_/@{--}[rrr]
&&\circ^\circlearrowright&\bullet}$$

\medskip

Long story short, we know what the simplex $Z_{full}$ is, and we have the formula $G^+(Z_{empty})=G^+(Z_{full})=S_Z^+$, exactly as in the $Z=\{1,\ldots,N\}$ case. Now with the above results in hand, we can talk as well about complementation, as follows:

\begin{theorem}
For any finite quantum graph $X$ we have the formula
$$G^+(X)=G^+(X^c)$$
where $X\to X^c$ is the complementation operation, given by $d_X+d_{X^c}=d_{Z_{full}}$.
\end{theorem}

\begin{proof}
This follows from Theorem 7.7, and more specifically from the following commutation relation, which is automatic, as explained there:
$$[u,d_{Z_{full}}]=0$$

Let us mention too that, in what concerns the pictorial representation of $X^c$, this can be deduced from what we have Theorem 7.10, in the obvious way.
\end{proof}

Following now \cite{ba3}, let us discuss an important point, namely the ``independence on the colors'' question. The idea indeed is that given a classical graph $X$ with edges colored by complex numbers, or other types of colors, $G(X)$ does not change when changing the colors. This is obvious, and a quantum analogue of this fact, involving $G^+(X)$, can be shown to hold as well, as explained in \cite{ba3}, and in chapter 5. In the quantum graph setting things are more complicated. Let us start with the following technical definition:

\index{color independence}

\begin{definition}
We call a quantum graph $X=(Z,d)$ washable if, with
$$d=\sum_ccd_c$$
being the color decomposition of $d$, we have the equivalence
$$[u,d]=0\iff[u,d_c]=0,\forall c$$
valid for any magic unitary matrix $u$, having size $|Z|$.
\end{definition}

Obviously, this is something which is not very beautiful, but the point is that some quantum graphs are washable, and some other are not, and so we have to deal with the above definition, as stated. As a first observation, we have the following result:

\begin{proposition}
Assuming that $X=(Z,d)$ is washable, its quantum symmetry group $G^+(X)$ does not depend on the precise colors of $X$. That is, whenever we have another quantum graph $X'=(Z,d')$ with same color scheme, in the sense that 
$$d_{ij}=d_{kl}\iff d'_{ij}=d'_{kl}$$
we have $G^+(X)=G^+(X')$.
\end{proposition}

\begin{proof}
This is something which is clear from the definition of $G^+(X)$, namely:
$$C(G^+(X))=C(S_Z^+)\Big/\Big<[u,d]=0\Big>$$

Indeed, assuming that our graph is washable in the above sense, we have:
$$C(G^+(X))=C(S_Z^+)\Big/\Big<[u,d_c]=0,\forall c\Big>$$

Thus, we are led to the conclusion in the statement.
\end{proof}

As already mentioned, it was proved in \cite{ba3} that in the classical case, $Z=\{1,\ldots,N\}$, all graphs are washable. This is a key result, and this for several reasons: 

\medskip

(1) First, this gives some intuition on what is going on with respect to colors, in analogy with what happens for $G(X)$. Also, it allows the use of true colors, like black, blue, red and so on, when drawing colored graphs, instead of complex numbers.

\medskip

(2) Also, this can be combined with the fact that $G^+(X)$ is invariant as well via some similar changes in the spectral decomposition of $d$, at the level of eigenvalues, with all this leading to some powerful combinatorial methods for the computation of $G^+(X)$.

\medskip

All these things do not necessarily hold in general, and to start with, we have:

\begin{theorem}
There are quantum graphs, such as the simplex in the homogeneous quantum space case, where
$$Z=M_K\times\{1,\ldots,L\}$$
with $K,L\geq2$, which are not washable.
\end{theorem}

\begin{proof}
We know that the simplex, in the case $Z=M_K\times\{1,\ldots,L\}$, has as adjacency matrix a certain matrix $d\in M_N(-1,0,1)$, with $N=K^2L$. Moreover, assuming $K,L\geq2$ as above, entries of all types, $-1,0,1$, are possible. Thus, the color decomposition of the adjacency matrix is as follows, with all 3 components being nonzero:
$$d=-1\cdot d_{-1}+0\cdot d_0+1\cdot d_1$$

Now assume that our simplex $X=Z_{full}$ is washable, and let $u$ be the fundamental corepresentation of $G^+(X)$. We have then the following commutation relations:
$$d_{-1},d_0,d_1\in End(u)$$

Now since the matrices $d_{-1},d_0,d_1$ are all nonzero, we deduce from this that:
$$\dim(End(u))\geq3$$

On the other hand, we know from Theorem 7.7 that we have $G^+(X)=S_Z^+$. Also, we know from chapter 4 that the Tannakian category of $S_Z^+$ is the Temperley-Lieb category $TL_N$, with $N=K^2L$ as above. By putting these two results together, we obtain:
$$\dim(End(u))=\dim\Big(span\Big(\ |\ |\ ,{\ }^\cup_\cap\ \Big)\,\Big)\leq2$$

Thus, we have a contradiction, and so our simplex is not washable, as claimed.
\end{proof}

In order to come up with some positive results as well, the general idea will be that of using the method in \cite{ba3}. We have the following statement, coming from there:

\begin{theorem}
The following matrix belongs to $End(u)$, for any $n\in\mathbb N$:
$$d^{\times n}_{ij}=\sum_{i=k_1\ldots k_n}\sum_{j=l_1\ldots l_n}d_{k_1l_1}\ldots d_{k_nl_n}$$
In particular, in the classical case, $Z=\{1,\ldots,N\}$, all graphs are washable.
\end{theorem}

\begin{proof}
We have two assertions here, the idea being as follows:

\medskip

(1) Consider the multiplication and comultiplication maps of the algebra $C(Z)$, which in single index notation are given by:
$$\mu(e_i\otimes e_j)=e_{ij}\quad,\quad 
\gamma(e_i)=\sum_{i=jk}e_j\otimes e_k$$

Observe that we have $\mu^*=\gamma$, with the adjoint taken with respect to the scalar product coming from the canonical trace. We conclude that we have:
$$\mu\in Hom(u^{\otimes 2},u)\quad,\quad 
\gamma\in Hom(u,u^{\otimes 2})$$

The point now is that we can consider the iterations $\mu^{(n)},\gamma^{(n)}$ of $\mu,\gamma$, constructed in the obvious way, and we have then, for any $n\in\mathbb N$:
$$\mu^{(n)}\in Hom(u^{\otimes n},u)\quad,\quad 
\gamma^{(n)}\in Hom(u,u^{\otimes n})$$

Now if we assume that we have $d\in End(u)$, we have $d^{\otimes n}\in End(u^{\otimes n})$ for any $n\in\mathbb N$, and we conclude that we have the following formula:
$$\mu^{(n)}d^{\otimes n}\gamma^{(n)}\in End(u)$$

But, in single index notation, we have the following formula:
$$(\mu^{(n)}d^{\otimes n}\gamma^{(n)})_{ij}=\sum_{i=k_1\ldots k_n}\sum_{j=l_1\ldots l_n}d_{k_1l_1}\ldots d_{k_nl_n}$$

Thus, we are led to the conclusion in the statement.

\medskip

(2) Assuming that we are in the case $Z=\{1,\ldots,N\}$, the matrix $d^{\times n}$ in the statement is simply the componentwise $n$-th power of $d$, given by:
$$d^{\times n}_{ij}=d_{ij}^n$$ 

As explained in \cite{ba3} or in chapter 5, a simple analytic argument based on this, using $n\to\infty$ and then a recurrence on the number of colors, shows that we have washability indeed. Thus, we are led to the conclusions in the statement.
\end{proof}

In order now to further exploit the findings in Theorem 7.15, an idea would be that of assuming that we are in the homogeneous case, $Z=M_K\times\{1,\ldots,L\}$, and that the adjacency matrix is split, in the sense that one of the following happens:
$$d_{ab,cd}=e_{ab}f_{cd}\quad,\quad
d_{ab,cd}=e_{ac}f_{bd}\quad,\quad 
d_{ab,cd}=e_{ad}f_{bc}$$

Normally the graph should be washable in this case, but the computations are quite complex, and there is no clear result known in this sense. Thus, as a conclusion to all this, the basic theory of the quantum groups $G^+(X)$ from chapter 5 extends well to the present quantum graph setting, modulo some subtleties in connection with the colors. 

\section*{7b. Cayley graphs}

With the basic theory of quantum symmetry groups of quantum graphs $G^+(X)\subset S_Z^+$ developed, we can now start hunting for examples. We already know from chapter 5 that, even in the $Z=\{1,\ldots,N\}$ case, going beyond trivialities is quite tricky. In the present, general case, there are however at least 3 things that can be done, as follows:

\bigskip

(1) We can try to study the quantum graphs modelled on the simplest, non-trivial quantum space, namely $Z=M_2$. Here we have $S_Z^+=S_Z=SO_3$, and for the few graphs here, having $|M_2|=4$ vertices, we recover certain subgroups of $SO_3$. We will discuss this in chapter 10 below, when talking about the ADE classification of such subgroups.

\bigskip

(2) We can look for the $S_Z^+$ analogues of the quantum reflection groups. This is something quite fundamental, and potentially useful, and we defer the discussion here, which is a bit technical, to the end of the present chapter. We will see that such ``generalized quantum reflection groups'' exist indeed, and are quite interesting objects.

\bigskip

(3) Finally, we can talk about Cayley theorems and Frucht theorems, in the present quantum space and quantum graph setting. We already know from chapter 4 that the passage $\{1,\ldots,N\}\to Z$ fixes a number of issues in respect to Cayley theorems, and we will see here that the same happens in relation with Frucht theorems. 

\bigskip

In order to get started, we have the following result, that we basically know from chapter 4, and that we reproduce here for convenience:

\index{Cayley theorem}
\index{no-Cayley theorem}
\index{Cayley embedding}

\begin{theorem}
Any finite quantum group $G$ has a Cayley embedding, as follows:
$$G\subset S_G^+$$
However, there are finite quantum groups which are not quantum permutation groups.
\end{theorem}

\begin{proof}
There are two statements here, the idea being as follows:

\medskip

(1) We have an action $G\curvearrowright G$, which leaves invariant the Haar measure. Now since the counting measure is left and right invariant, so is the Haar measure, we conclude that $G\curvearrowright G$ leaves invariant the counting measure, and so we have $G\subset S_G^+$, as claimed.

\medskip

(2) Regarding the second assertion, this is something non-trivial, from \cite{bbn}, the simplest counterexample being a certain quantum group $G$ appearing as a split abelian extension associated to the factorization $S_4=\mathbb Z_4S_3$, having cardinality $|G|=24$. 
\end{proof}

The question now is whether we can do better, by viewing $G$ as the quantum symmetry group of a certain finite quantum graph $X_G$, modelled on the finite quantum space $G$. Such questions are quite interesting, and in order to discuss them, let us start with a study for the classical finite groups $G$. Here we first have the following result: 

\begin{proposition}
Any finite group $G$ appears as $G=G(X)$, with $X$ being the colored oriented graph having $G$ as set of vertices, and with the edges being colored by
$$d_{hk}=h^{-1}k$$
according to the usual colored graph conventions, with color set $C=G$.
\end{proposition}

\begin{proof}
Consider indeed the Cayley action of $G$ on itself, which is given by:
$$G\subset S_G\quad,\quad g\to[h\to gh]$$

We have $d_{gh,gk}=d_{hk}$, which gives an action $G\curvearrowright X$, and so an inclusion $G\subset G(X)$. Conversely now, pick an arbitrary permutation $\sigma\in S_G$. We have then:
\begin{eqnarray*}
\sigma\in G(X)
&\implies&d_{\sigma(h)\sigma(k)}=d_{hk}\\
&\implies&\sigma(h)^{-1}\sigma(k)=h^{-1}k\\
&\implies&\sigma(1)^{-1}\sigma(k)=k\\
&\implies&\sigma(k)=\sigma(1)k\\
&\implies&\sigma\in G
\end{eqnarray*}

Thus, the inclusion $G\subset G(X)$ is an equality, as desired.
\end{proof}

We have the following improvement of the above result, in the quantum setting:

\begin{theorem}
Any finite group $G$ appears as $G=G^+(X)$, with $X$ being the colored oriented graph from Proposition 7.17, with vertex set $G$ and adjacency matrix
$$d_{hk}=h^{-1}k$$
with the conventions there, namely that the set of colors is $C=G$ itself.
\end{theorem}

\begin{proof}
We know from chapter 5 that the magic matrix $u=(u_{ij})$ of the quantum group $G^+(X)$ must commute with all the color components of $d$. But these color components, $d^c\in M_G(0,1)$ with $c\in G$, are given by the following formula:
$$d^c_{hk}
=\begin{cases}
1&{\rm if}\ h^{-1}k=c\\
0&{\rm otherwise}
\end{cases}$$

With this formula in hand, the commutation computation goes as follows:
\begin{eqnarray*}
d^cu=ud^c
&\iff&(d^cu)_{hk}=(ud^c)_{hk}\\
&\iff&\sum_gd^c_{hg}u_{gk}=\sum_gu_{hg}d^c_{gk}\\
&\iff&\sum_g\delta_{h^{-1}g,c}u_{gk}=\sum_gu_{hg}\delta_{g^{-1}k,c}\\
&\iff&u_{hc,k}=u_{h,kc^{-1}}\\
&\iff&u_{hc,dc}=u_{hd}
\end{eqnarray*}

In particular, with $h=1$ we obtain the following formula:
$$u_{c,dc}=u_{1d}$$

Thus all the rows of our magic matrix $u=(u_{ij})$ appear as permutations of the 1st row. Now since the entries in the 1st row commute, as being pairwise orthogonal projections, it follows that all the entries of $u$ commute. Thus our graph has no quantum symmetry, $G^+(X)=G(X)$, and the result follows from Proposition 7.17.
\end{proof}

Our next goal will be that of improving Proposition 7.17 and Theorem 7.18, as to get rid of orientation, or colors. We first have the following extension of Proposition 7.17:

\begin{proposition}
Given a finite group $G$, written as $G=<S>$, define the associated Cayley graph $X$ by saying that the edges $g\to h$ exist when 
$$g^{-1}h\in S$$
and are colored by $g^{-1}h\in S$ precisely. Then we have $G=G(X)$.
\end{proposition}

\begin{proof}
We use the same method as in the proofs of Proposition 7.17 and Theorem 7.18. The adjacency matrix of the graph $X$ in the statement is given by:
$$d_{hk}
=\begin{cases}
h^{-1}k&{\rm if}\ h^{-1}k\in S\\
0&{\rm otherwise}
\end{cases}$$

Thus we have an action $G\curvearrowright X$, and so on inclusion $G\subset G(X)$. Conversely now, pick an arbitrary permutation $\sigma\in S_G$. We know that $\sigma$ must preserve all the color components of $d$, which are the following matrices, depending on a color $c\in S$:
$$d_{hk}^c
=\begin{cases}
1&{\rm if}\ h^{-1}k=c\\
0&{\rm otherwise}
\end{cases}$$

In other words, we have the following equivalences:
\begin{eqnarray*}
\sigma\in G(X)
&\iff&d^c_{\sigma(h)\sigma(k)}=d_{hk},\forall c\in S\\
&\iff&\sigma(h)^{-1}\sigma(k)=h^{-1}k,\forall h^{-1}k\in S
\end{eqnarray*}

Now observe that with $h=1$ we obtain from this that we have:
\begin{eqnarray*}
k\in S
&\implies&\sigma(1)^{-1}\sigma(k)=k\\
&\implies&\sigma(k)=\sigma(1)k
\end{eqnarray*}

Next, by taking $h\in S$, we obtain from the above formula that we have:
\begin{eqnarray*}
k\in hS
&\implies&\sigma(h)^{-1}\sigma(k)=h^{-1}k\\
&\implies&\sigma(k)=\sigma(h)h^{-1}k\\
&\implies&\sigma(k)=(\sigma(1)h)h^{-1}k\\
&\implies&\sigma(k)=\sigma(1)k
\end{eqnarray*}

Thus with $g=\sigma(1)$ we have the following formula, for any $k\in S$:
$$\sigma(k)=gk$$

But  the same method shows that this formula holds as well for any $k\in S^2$, then for any $k\in S^3$, any $k\in S^4$, and so on. Thus the above formula $\sigma(k)=gk$ holds for any $k\in G$, and so the inclusion $G\subset G(X)$ is an equality, as desired.
\end{proof}

We have the following joint extension of Theorem 7.18 and Proposition 7.19:

\begin{theorem}
Given a finite group $G$, written as $G=<S>$, define the associated Cayley graph $X$ by saying that the edges $g\to h$ exist when 
$$g^{-1}h\in S$$
and are colored by $g^{-1}h\in S$ precisely. Then we have $G=G^+(X)$.
\end{theorem}

\begin{proof}
We use the same method as in the proof of Theorem 7.18. The computation there applies, but this time under the assumption $c\in S$, and gives:
$$d^cu=ud^c\iff u_{hc,dc}=u_{hd}$$

Now observe that with $h=1$ we obtain from this that we have:
$$u_{c,dc}=u_{1d}$$

But this shows that the entries of the various $c$-rows of $u=(u_{ij})$, with $c\in\{1\}\cup S$, pairwise commute. Next, by taking $h=c'\in S$, we obtain from the above formula:
$$h_{c'c,d}=u_{c'd}$$

We conclude from this that the entries of the various $c$-rows of $u=(u_{ij})$, indexed by group elements $c\in\{1\}\cup S\cup S^2$, pairwise commute. By continuing the procedure we obtain in this way that all the entries of $u$ commute. Thus our graph has no quantum symmetry, $G^+(X)=G(X)$, and the result follows from Proposition 7.19.
\end{proof}

The above result is not the end of the story, but rather the beginning of it. Indeed, at a more advanced level, we have the following classical result, due to Frucht:

\begin{theorem}[Frucht]
Any finite group $G$ appears as the symmetry group
$$G=G(X)$$
of a certain uncolored, unoriented graph $X$.
\end{theorem}

\begin{proof}
The idea is to start with the oriented graph in Proposition 7.19, and make a unoriented graph out of it, by replacing each edge with a copy of the following graph, with the height being in a chosen bijection with the corresponding element $g^{-1}h\in S$: 
$$\xymatrix@R=20pt@C=30pt{
&&\circ\\
&\circ&\circ\ar@{-}[u]\\
\\
\circ_g\ar@{-}[r]&\circ\ar@{.}[uu]\ar@{-}[r]&\circ\ar@{.}[uu]\ar@{-}[r]&\circ_h}$$

With these replacements made, under suitable assumptions on the generating $S$, namely $1\notin S$, plus the fact that $S\cap S^{-1}$ must consist only of involutions, one can prove that $G$ appears indeed of the symmetry group of this graph $X$. 
\end{proof}

So long for classical groups $G$, viewed as symmetry groups. Further questions, which are more technical, include getting rid of orientation in Theorem 7.20, as for $G$ to appear as classical or quantum symmetry group of its Cayley graph, viewed as metric space, and also the study of the quantum symmetry group of the graph from Theorem 7.21.

\bigskip

In the quantum group case, similar problems make sense, with Theorem 7.16 as a starting point. For more on all this, we refer to \cite{fpi}, \cite{gro} and related papers.

\section*{7c. Twisted reflections}

With the above technology in hand, we can talk about twisted quantum reflections. The idea here, from \cite{ba2}, will be that the twisted analogues of the quantum reflection groups $H_N^{s+}\subset S_{sN}^+$ from chapter 6 will be the quantum automorphism groups $S_{Z\to Y}^+$ of the fibrations of finite quantum spaces $Z\to Y$, which correspond by definition to the Markov inclusions of finite dimensional $C^*$-algebras $C(Y)\subset C(Z)$.

\bigskip

In order to discuss this material, let us start with the following definition:

\index{Markov inclusion}
\index{Markov fibration}

\begin{definition}
A fibration of finite quantum spaces $Z\to Y$ corresponds to an inclusion of finite dimensional $C^*$-algebras
$$C(Y)\subset C(Z)$$
which is Markov, in the sense that it commutes with the canonical traces.
\end{definition}

Here the commutation condition with the canonical traces means that the composition $C(Y)\subset C(Z)\to\mathbb C$ should equal the canonical trace $C(Y)\to\mathbb C$. At the level of the corresponding quantum spaces, this means that the quotient map $Z\to Y$ must commute with the corresponding counting measures, and this is where our term ``fibration'' comes from. As for the term ``Markov'', this is something standard in subfactor theory.

\bigskip

In order to talk about the quantum symmetry groups $S_{Z\to Y}^+$, we will need:

\begin{proposition}
Given a fibration $Z\to Y$, a closed subgroup $G\subset S_Z^+$ leaves invariant $Y$ precisely when its magic unitary $u=(u_{ij})$ satisfies the condition
$$e\in End(u)$$
where $e:C(Z)\to C(Z)$ is the Jones projection, onto the subalgebra $C(Y)\subset C(Z)$.
\end{proposition}

\begin{proof}
This is something that we know well from chapter 5, in the commutative case, where $Z$ is a usual finite set, and the proof in general is similar.
\end{proof}

We can now talk about twisted quantum reflection groups, as follows:

\begin{theorem}
Any fibration of finite quantum spaces $Z\to Y$ has a quantum symmetry group, which is the biggest acting on $Z$ by leaving $Y$ invariant:
$$S_{Z\to Y}^+\subset S_Z^+$$
At the level of algebras of functions, this quantum group $S_{Z\to Y}^+$ is obtained as follows, with $e:C(Z)\to C(Y)$ being the Jones projection:
$$C(S_{Z\to Y}^+)=C(S_Z^+)\Big/\Big<e\in End(u)\Big>$$
We call these quantum groups $S_{Z\to Y}^+$ twisted quantum reflection groups. 
\end{theorem}

\begin{proof}
This follows indeed from Proposition 7.23.
\end{proof}

As a basic example, let us discuss the commutative case. Here we have:

\begin{proposition}
In the commutative case, the fibration $Z\to Y$ must be of the following special form, with $N,s$ being certain integers,
$$\{1,\ldots,N\}\times\{1,\ldots,s\}\to\{1,\ldots,N\}\quad,\quad (i,a)\to i$$
and we obtain the quantum reflection groups studied in chapter 6,
$$(S_{Z\to Y}^+\subset S_Z^+)\ = \ (H_N^{s+}\subset S_{sN}^+)$$
via some standard identifications.
\end{proposition}

\begin{proof}
In the commutative case our fibration must be a usual fibration of finite spaces, $\{1,\ldots,M\}\to\{1,\ldots,N\}$, commuting with the counting measures. But this shows that our fibration must be of the following special form, with $N,s\in\mathbb N$:
$$\{1,\ldots,N\}\times\{1,\ldots,s\}\to\{1,\ldots,N\}\quad,\quad (i,a)\to i$$

Regarding now the quantum symmetry group, we have the following formula for it, with $e:\mathbb C^N\otimes\mathbb C^s\to\mathbb C^N$ being the Jones projection for the inclusion $\mathbb C^N\subset\mathbb C^N\otimes\mathbb C^s$:
$$C(S_{Z\to Y}^+)=C(S_{sN}^+)\Big/\Big<e\in End(u)\Big>$$

On the other hand, recall from chapter 6 that the quantum reflection group $H_N^{s+}\subset S_{sN}^+$ appears via the condition that the corresponding magic matrix must be sudoku:
$$u=\begin{pmatrix}
a^0&a^1&\ldots&a^{s-1}\\
a^{s-1}&a^0&\ldots&a^{s-2}\\
\vdots&\vdots&&\vdots\\
a^1&a^2&\ldots&a^0
\end{pmatrix}$$

But, as explained in \cite{ba2}, \cite{ba3}, this is the same as saying that the quantum group $H_N^{s+}\subset S_{sN}^+$ appears as the symmetry group of the multi-simplex associated to the fibration $\{1,\ldots,N\}\times\{1,\ldots,s\}\to\{1,\ldots,N\}$, so we have an identification as follows:
$$(S_{Z\to Y}^+\subset S_Z^+)\ = \ (H_N^{s+}\subset S_{sN}^+)$$

Thus, we are led to the conclusions in the statement.
\end{proof}

Observe that in Proposition 7.25 the fibration $Z\to Y$ is ``trivial'', in the sense that it is of the following special form:
$$Y\times T\to Y\quad,\quad (i,a)\to i$$

However, in the general quantum case, there are many interesting fibrations $Z\to Y$ which are not trivial, and in what follows we will not make any assumption on our fibrations, and use Definition 7.22 and Theorem 7.24 as stated.

\section*{7d. Subfactor results}

Following \cite{ba2}, we will prove now that the Tannakian category of $S_{Z\to Y}^+$, which is by definition a generalization of $S_Z^+$, is the Fuss-Catalan category, which is a generalization of the Temperley-Lieb category, introduced by Bisch and Jones in \cite{bjo}.

\bigskip

In order to do so, let us first reformulate Theorem 7.24 in a more convenient way, in purely functional analytic terms, and also as a self-contained statement, as follows:

\begin{theorem}
Any Markov inclusion of finite dimensional algebras $D\subset B$ has a quantum symmetry group $S_{D\subset B}^+$. The corresponding Woronowicz algebra is generated by the coefficients of a biunitary matrix $v=(v_{ij})$ subject to the conditions
$$m\in Hom(v^{\otimes 2},v)\quad,\quad 
u\in Hom(1,v)\quad,\quad 
e\in End(v)$$
where $m:B\otimes B\to B$ is the multiplication, $u:\mathbb C\to B$ is the unit and $e:B\to B$ is the projection onto $D$, with respect to the scalar product $<x,y>=tr(xy^*)$.
\end{theorem}

\begin{proof}
This is a reformulation of Theorem 7.24, with several modifications made. Indeed, by using the algebras $D=C(Y)$, $B=C(Z)$ instead of the quantum spaces $Y,Z$ used there, and also by calling the fundamental corepresentation $v=(v_{ij})$, in order to avoid confusion with the unit $u:\mathbb C\to B$, the formula in Theorem 7.24 reads:
$$C(S_{D\subset B}^+)=C(S_B^+)\Big/\Big<e\in End(v)\Big>$$

Also, we know from Theorem 7.4 that we have the following formula, again by using $B$ instead of $Z$, and by calling the fundamental corepresentation $v=(v_{ij})$:
$$C(S_B^+)=C(U_N^+)\Big/\Big<m\in Hom(v^{\otimes2},v),u\in Fix(v)\Big>$$

Thus, we are led to the conclusion in the statement.
\end{proof}

As already mentioned, our goal will be that of proving that the Tannakian category of $S_{D\subset B}^+$ is the Fuss-Catalan category, introduced by Bisch and Jones in \cite{bjo}. This will be something quite technical, from \cite{ba2}, jointly generalizing our previous Temperley-Lieb computations for $S_Z^+$, and our previous results regrading quantum reflection groups.

\bigskip

The Fuss-Catalan category, to be introduced in a moment, is a tensor $C^*$-category having $(\mathbb N,+)$ as monoid of objects. In what follows, we will call such a tensor category a $\mathbb N$-algebra. If $C$ is a $\mathbb N$-algebra we use the following notations:
$$C(m,n)=Hom_C(m,n)\quad,\quad 
C(m)=End_C(m)$$

Let us first discuss in detail the Temperley-Lieb algebra, as a continuation of the material from chapters 1-4. In the present context, we have the following definition:

\index{Fuss-Catalan algebra}

\begin{definition}
The $\mathbb N$-algebra $TL^2$ of index $\delta >0$ is defined as follows:
\begin{enumerate}
\item The space $TL^2(m,n)$ consists of linear combinations
of noncrossing pairings between $2m$ points and $2n$ points:
$${{TL^2}}(m,n)=
\left\{\sum\,\alpha\,\,
\begin{matrix}\cdots\cdots&\leftarrow&2m\,\,\ {\rm points}\\
\mathfrak{W}&\leftarrow&m+n\ \ {\rm strings}\\
\cdot\,\cdot&\leftarrow&2n\,\,{\rm points}
\end{matrix}\right\}$$

\item The operations $\circ$, $\otimes$, $*$ are induced by the vertical and horizontal concatenation and the upside-down turning of diagrams:
$$A\circ B=\binom{B}{A}\quad,\quad A\otimes B=AB\quad,\quad A^*=\forall$$

\item With the rule $\bigcirc=\delta$, erasing a circle is the same as multiplying by $\delta$.
\end{enumerate}
\end{definition}

Our first task will be that of finding a suitable presentation for this algebra. Consider the following two elements $u\in TL^2(0,1)$ and $m\in TL^2(2,1)$:
$$u=\delta^{-\frac{1}{2}}\,\cap
\quad,\quad
m=\delta^{\frac{1}{2}}\,|\cup|$$ 

With this convention, we have the following result:

\begin{theorem}
The following relations are a presentation of $TL^2$ by the above rescaled diagrams $u\in TL^2(0,1)$ and $m\in TL^2(2,1)$:
\begin{enumerate}
\item $mm^*=\delta^2$.

\item $u^*u=1$.

\item $m(m\otimes1)=m(1\otimes m)$.

\item $m(1\otimes u)=m(u\otimes1)=1$.

\item $(m\otimes1)(1\otimes m^*)=(1\otimes m)(m^*\otimes1)=m^*m$.
\end{enumerate}
\end{theorem}

\begin{proof}
This is something well-known, and elementary, by drawing diagrams, and for details here, we refer for instance to \cite{bjo}.
\end{proof}

In more concrete terms, the above result says that $u,m$ satisfy the above relations, which is something clear, and that if $C$ is a $\mathbb N$-algebra and $v\in C(0,1)$ and $n\in C(2,1)$ satisfy the same relations then there exists a $\mathbb N$-algebra morphism as follows:
$$TL^2\to C\quad,\quad 
u\to v\quad,\quad 
m\to n$$

Now let $B$ be a finite dimensional $C^*$-algebra, with its canonical trace. We have a scalar product $<x,y>=tr(xy^*)$ on $B$, so $B$ is an object in the category of finite dimensional Hilbert spaces. Consider the unit $u$ and the multiplication $m$ of $B$:
$$u\in\mathbb NB(0,1)\quad,\quad 
m\in\mathbb NB(2,1)$$

The relations in Theorem 7.28 are then satisfied, and one can deduce from this that in this case, the category of representations of $S_B^+$ is the completion of $TL^2$, as we already know from chapter 4. Getting now to Fuss-Catalan algebras, we have here:

\begin{definition}
A Fuss-Catalan diagram is a planar diagram formed by an upper row of $4m$ points, a lower row of $4n$ points, both colored
$$\circ\bullet\bullet\circ\circ\bullet\bullet\ldots$$
and by $2m+2n$ noncrossing strings joining these $4m+4n$ points, with the rule that the points which are joined must have the same color.
\end{definition}

Fix $\beta>0$ and $\omega>0$. The $\mathbb N$-algebra $FC$ is defined as follows. The spaces $FC(m,n)$ consist of linear combinations of Fuss-Catalan diagrams:
$$FC(m,n)=\left\{\sum\,\alpha\,\,
\begin{matrix}\circ\bullet\bullet\circ\circ\bullet\bullet\circ\ldots\ldots&\leftarrow&4m\,\,\ {\rm colored\ points}\\
\ &\ &m+n {\rm\ black\ strings}\\
\mathfrak{W}\mbox{\ \ }&\leftarrow&{\rm and}\\
\ & \ & m+n{\rm\ white\ strings}\\
\circ\bullet\bullet\circ\circ\bullet\bullet\circ\ldots\ldots&\leftarrow&4n\,\,{\rm colored\ points}\end{matrix}\right\}$$

As before with the Temperley-Lieb algebra, the operations $\circ$, $\otimes$, $*$ are induced by vertical and horizontal concatenation and upside-down turning of diagrams, but this time with the rule that erasing a black/white circle is the same as multiplying by $\beta$/$\omega$:
$$A\circ B=\binom{B}{A}\quad,\quad 
A\otimes B=AB\quad,\quad 
A^*=\forall$$
$$\mbox{\tiny{black}}\rightarrow\bigcirc=\beta
\qquad,\qquad
\mbox{\tiny{white}}\rightarrow\bigcirc=\omega$$

Let $\delta=\beta\omega$. We have the following bicolored analogues of the elements $u,m$:
$$u=\delta^{-\frac{1}{2}}\,\,\bigcap\;\!\!\!\!\!\!\!\!\cap
\quad,\quad
m=\delta^{\frac{1}{2}}\,\,||\bigcup\;\!\!\!\!\!\!\!\!\cup\ ||$$

Consider also the black and white Jones projections, namely:
$$e=\omega^{-1}\,\,|\ \begin{matrix}\cup\cr\cap\end{matrix}\ |
\quad,\quad 
f=\beta^{-1}\,\,|||\ \begin{matrix}\cup\cr\cap\end{matrix}\ |||$$

For simplifying writing we identify $x$ and $x\otimes 1$. We have the following result:

\begin{theorem}
The following relations, with $f=\beta^{-2}(1\otimes me)m^*$, are a presentation of $FC$ by $m\in FC(2,1)$, $u\in FC(0,1)$ and $e\in FC(1)$:
\begin{enumerate}
\item The relations in Theorem 7.28, with $\delta =\beta\omega$.

\item $e=e^2=e^*$, $f=f^*$ and $(1\otimes f)f=f(1\otimes f)$.

\item $eu=u$.

\item $mem^*=m(1\otimes e)m^*=\beta^2$.

\item $mm(e\otimes e\otimes e)=emm(e\otimes1\otimes e)$.
\end{enumerate}
\end{theorem}

\begin{proof}
As for any presentation result, we have to prove two assertions:

\medskip

(I) The elements $m,u,e$ satisfy the relations (1-5) and generate the
$\mathbb N$-algebra $FC$.

\medskip

(II) If  $M$, $U$ and $E$ in a $\mathbb N$-algebra $C$ satisfy the relations
(1-5), then there exists a morphism of $\mathbb N$-algebras $FC\to C$
sending $m\to M$, $u\to U$, $e\to E$.

\medskip

(I) First, the relations (1-5) are easily verified by drawing pictures. Let us show now that the $\mathbb N$-subalgebra $C=<m,u,e>$ of $FC$ is equal to $FC$. First, $C$ contains the infinite sequence of black and white Jones projections:
$$p_1=e=\omega^{-1}\,\,|\ \begin{matrix}\cup\\\cap\end{matrix}\ |$$
$$p_2=f=\beta^{-1}\,\,|||\ \begin{matrix}\cup\\\cap\end{matrix}\ |||$$
$$p_3=1\otimes e=\omega^{-1}\,\,|||||\ \begin{matrix}\cup\\\cap\end{matrix}\ |$$
$$p_4=1\otimes f=\beta^{-1}\,\,|||||||\ \begin{matrix}\cup\cr\cap\end{matrix}\ |||$$
$$\vdots$$

The algebra $C$ contains as well the infinite sequence of bicolored Jones projections:
$$e_1=uu^*=\delta^{-1}\,\,\begin{matrix}\bigcup\hskip -3.9mm\cup\\ \bigcap\hskip -3.9mm\cap\end{matrix}$$
$$e_2=\delta^{-2}m^*m=\delta^{-1}\,\,||\ \begin{matrix}\bigcup\hskip -3.86mm\cup\\\bigcap\hskip -3.86mm\cap\end{matrix}\ ||$$
$$e_3=1\otimes uu^*=\delta^{-1}\,\,||||\ \begin{matrix}\bigcup\hskip -3.86mm\cup\\ \bigcap\hskip -3.86mm\cap\end{matrix}$$
$$e_4=\delta^{-2}(1\otimes m^*m)=\delta^{-1}\,\,||||||\ \begin{matrix}\bigcup\hskip -3.86mm\cup\\\bigcap\hskip -3.86mm\cap\end{matrix}\ ||$$
$$\vdots$$

By the results of Bisch and Jones in \cite{bjo}, these latter projections generate the diagonal $\mathbb N$-algebra $\Delta FC$. Thus we have inclusions as follows:
$$\Delta FC\subset C\subset FC$$

By definition of $C$, we have as well the following equality:
$$\Delta FC=\Delta C$$

Also, the existence of semicircles shows that the objects of $C$ and $FC$ are self-dual, and by Frobenius reciprocity we obtain that for $m+n$ even, we have:
$$\dim(C(m,n))
=\dim(FC(m,n))$$

By tensoring with $u$ and $u^*$ we get embeddings as follows:
$$C(m,n)\subset C(m,n+1)\quad,\quad 
FC(m,n)\subset FC(m,n+1)$$

But this shows that the above dimension equalities hold for any $m$ and
$n$. Together with $\Delta FC\subset C\subset FC$, this shows that $C=FC$.

\medskip

(II) Assume that $M,U,E$ in a $\mathbb N$-algebra $C$ satisfy the
relations (1-5). We have to construct a morphism $FC\to C$ sending:
$$m\to M\quad,\quad 
u\to U\quad,\quad 
e\to E$$

As a first task, we would like to construct a morphism $\Delta FC\to\Delta C$ sending:
$$m^*m\to M^*M\quad,\quad
uu^*\to UU^*\quad,\quad 
e\to E$$

By constructing the corresponding Jones projections $E_i$ and $P_i$, we must send:
$$e_i\to E_i\quad,\quad
p_i\to P_i$$

In order to construct these maps, we use now the fact, from \cite{bjo}, that the following relations are a presentation of $\Delta FC$: 

\bigskip

(a) $e_i^2=e_i$, $e_ie_j=e_je_i$ if $|i-j|\geq2$ and $e_ie_{i\pm1}e_i=\delta^{-2}e_i$.

\medskip

(b) $p_i^2=p_i$ and $p_ip_j=p_jp_i$.

\medskip

(c) $e_ip_i=p_ie_i=e_i$ and $p_ie_j=e_jp_i$ if $|i-j|\geq2$.

\medskip

(d) $e_{2i\pm1}p_{2i}e_{2i\pm1}=\beta^{-2}e_{2i\pm1}$ and $e_{2i}p_{2i\pm1}e_{2i}=\omega^{-2}e_{2i}$.

\medskip

(e) $p_{2i}e_{2i\pm1}p_{2i}=\beta^{-2}p_{2i\pm 1}p_{2i}$ and $p_{2i\pm1}e_{2i}p_{2i\pm1}=\omega^{-2}p_{2i}p_{2i\pm1}$.

\bigskip

Thus, it remains to verify that we have the following implication, where $m,u,e$ are now abstract objects, and we are no longer allowed to draw pictures:
$$(1-5)\implies (a-e)$$

But these relations are all easy to verify. The conclusion is that we constructed a certain $\mathbb N$-algebra morphism, as follows:
$$\Delta J:\Delta FC\to\Delta C$$

We have to extend now this morphism into a morphism $J:FC\to C$ sending $u\to U$ and $m\to M$. We will use a standard argument. For $w\geq k,l$ we define:
$$\phi:FC(l,k)\to FC(w)$$
$$x\to(u^{\otimes(w-k)}\otimes1_k)\,x\,((u^*)^{\otimes(w-l)}\otimes1_l)$$ 

We can define as well a morphism as follows:
$$\theta:FC(w)\to FC(l,k)$$
$$x\to((u^*)^{\otimes (w-k)}\otimes1_k)\,x\,(u^{\otimes (w-l)}\otimes1_l)$$

Here $1_k=1^{\otimes k}$, and the convention $x=x\otimes1$ is no longer used. We define $\Phi$ and $\Theta$ in $C$ by similar formulae. We have $\theta\phi=\Theta\Phi=Id$. We define a map $J$ by:
$$\xymatrix@R=40pt@C=40pt{
FC(l,k)\ar[r]^J\ar[d]^\phi&C(l,k)\\
FC(w)\ar[r]^{\Delta J}&C(w)\ar[u]^\Theta}$$

Since $J(a)$ does now depend on the choice of $w$, these $J$ maps are the components of a global map $J:FC\to C$, which sends $u\to U$ and $m\to M$, as desired.
\end{proof}

Getting back now to the inclusions $D\subset B$, we have the following result:

\begin{theorem}
Given a Markov inclusion $D\subset B$, we have
$$<m,u,e>=FC$$
as an equality of $\mathbb N$-algebras.
\end{theorem}

\begin{proof}
It is routine to check that the linear maps $m,u,e$ associated to an inclusion $D\subset B$ as in the statement satisfy the relations (1-5) in Theorem 7.30. Thus, we obtain a certain $\mathbb N$-algebra surjective morphism, as follows:
$$J:FC\to<m,u,e>$$

It remains to prove that this morphism $J$ is faithful. For this purpose, consider the following map, where $v=m^*u\in FC(0,2)$:
$$\phi_n:FC(n)\to FC(n-1)$$
$$x\to(1^{\otimes(n-1)}\otimes v^*)(x\otimes1)(1^{\otimes(n-1)}\otimes v)$$

Consider as well the following map, where $v=m^*u\in FC(0,2)$ is as above:
$$\psi_n:C(n)\to C(n-1)$$ 
$$x\to(1^{\otimes(n-1)}\otimes J(v)^*)(x\otimes1)(1^{\otimes(n-1)}\otimes J(v))$$

These maps make then the following diagram commutative:
$$\xymatrix@R=40pt@C=40pt{
FC(n)\ar[r]^J\ar[d]^{\phi_n}&C(n)\ar[d]^{\psi_n}\\
FC(n-1)\ar[r]^J&C(n-1)}$$
  
By gluing such diagrams we get a factorization by $J$ of the composition on the left of conditional expectations, which is the Markov trace. By positivity $J$ is faithful on $\Delta FC$, then by Frobenius reciprocity faithfulness has to hold on the whole $FC$.
\end{proof}

Getting back now to quantum groups, we have:

\index{planar algebra}
\index{Fuss-Catalan algebra}
\index{twisted reflection group}

\begin{theorem}
Given a Markov inclusion $D\subset B$, the category of representations of its quantum symmetry group $S_{D\subset B}^+$ is the completion of $FC$.
\end{theorem}

\begin{proof}
Since $S_{D\subset B}^+$ comes by definition from the relations corresponding to $m,u,e$, its tensor category of corepresentations is the completion of the tensor category $<m,u,e>$. Thus Theorem 7.31 applies, and gives an isomorphism $<m,u,e>\simeq FC$.
\end{proof}

As already mentioned in the above, in terms of finite quantum spaces and quantum graphs, the conclusion of all this is that the quantum automorphism groups $S_{Z\to Y}^+$ of the Markov fibrations $Z\to Y$, which can be thought of as being the ``twisted versions'' of the quantum reflection groups $H_N^{s+}$ studied in chapter 6, correspond to the Fuss-Catalan algebras. We refer to \cite{ba2} and related papers for more on these topics.

\bigskip

Let us also mention that the various planar algebra results formulated so far throughout this book have extensions to the present setting. We refer here to \cite{twa}.

\section*{7e. Exercises} 

Things have been quite technical in this chapter, and as an exercise on all this, which is quite technical, but is quite instructive, we have:

\begin{exercise}
Make the connection between the easiness results from chapter 6, regarding the quantum groups $H_N^{s+}$, and the representation theory results from here, regarding the quantum automorphism groups $S_{Z\to Y}^+$ of the Markov fibrations $Z\to Y$.
\end{exercise}

This is something a bit tricky, in the sense that no new computations are needed, with the work instead consisting of making lots of identifications, in order to view the main results from chapter 6 as particular cases of those discussed here.

\chapter{Partial permutations}

\section*{8a. Partial permutations}

We discuss in this chapter an extension of some of the results that we have seen so far, both of algebraic and analytic nature, from the case of the quantum permutations, to the case of quantum partial permutations. Our motivations are as follows:

\medskip

(1) In the quantum case, several things are potentially simpler for partial permutations than for permutations, such as the diagonal subalgebras $<u_{ii}>\subset C(G)$.

\medskip

(2) Along the same lines, the quantum partial permutations of a graph $X$, finite or not, are sometimes simpler to understand than the quantum permutations.

\medskip

(3) We will see later in this book that the combinatorics of a partial Hadamard matrix $H\in M_{M\times N}(\mathbb T)$ is encoded by a certain quantum partial permutation semigroup. 

\medskip

Summarizing, we have good motivations, and potential applications in mind, and all that will follow will be far from being anecdotical. Getting started now, we have:

\index{partial permutation}

\begin{definition}
$\widetilde{S}_N$ is the semigroup of partial permutations of $\{1\,\ldots,N\}$,
$$\widetilde{S}_N=\left\{\sigma:X\simeq Y\Big|X,Y\subset\{1,\ldots,N\}\right\}$$
with the usual composition operation for such partial permutations, namely
$$\sigma'\sigma:\sigma^{-1}(X'\cap Y)\simeq\sigma'(X'\cap Y)$$
being the composition of $\sigma':X'\simeq Y'$ and $\sigma:X\simeq Y$.
\end{definition}

Observe that $\widetilde{S}_N$ is not simplifiable, because the null permutation $\emptyset\in\widetilde{S}_N$, having the empty set as domain/range, satisfies the following formula, , for any $\sigma\in\widetilde{S}_N$:
$$\emptyset\sigma=\sigma\emptyset=\emptyset$$

Observe also that our semigroup $\widetilde{S}_N$ has a kind of ``subinverse'' map, which is not a true inverse in the semigroup sense, sending a partial permutation $\sigma:X\to Y$ to its usual inverse $\sigma^{-1}:Y\to X$. As a first interesting result now about $\widetilde{S}_N$, which shows that we are dealing with some non-trivial combinatorics, we have:

\begin{proposition}
The number of partial permutations is given by
$$|\widetilde{S}_N|=\sum_{k=0}^Nk!\binom{N}{k}^2$$
that is, $1,2,7,34,209,\ldots\,$, and we have the formula
$$|\widetilde{S}_N|\simeq N!\sqrt{\frac{\exp(4\sqrt{N}-1)}{4\pi\sqrt{N}}}$$
in the $N\to\infty$ limit.
\end{proposition}

\begin{proof}
The first assertion is clear, because in order to construct a partial permutation $\sigma:X\to Y$ we must choose an integer $k=|X|=|Y|$, then we must pick two subsets $X,Y\subset\{1,\ldots,N\}$ having cardinality $k$, and there are $\binom{N}{k}$ choices for each, and finally we must construct a bijection $\sigma:X\to Y$, and there are $k!$ choices here. As for the estimate, which is non-trivial, this is something standard, and well-known.
\end{proof}

Another result, which is trivial, but quite fundamental, is as follows:

\begin{proposition}
We have a semigroup embedding $u:\widetilde{S}_N\subset M_N(0,1)$, given by 
$$u_{ij}(\sigma)=
\begin{cases}
1&{\rm if}\ \sigma(j)=i\\
0&{\rm otherwise}
\end{cases}$$
whose image are the matrices having at most one nonzero entry, on each row and column.
\end{proposition}

\begin{proof}
This is trivial from definitions, with $u:\widetilde{S}_N\subset M_N(0,1)$ extending the standard embedding $u:S_N\subset M_N(0,1)$, that we have been heavily using, so far.
\end{proof}

Let us discuss now some probabilistic aspects, related to the Poisson law computations from chapter 2. We denote by $\kappa:\widetilde{S}_N\to\mathbb N$ the cardinality of the domain/range, and by $\chi^l:\widetilde{S}_N\to\mathbb N$ the number of fixed points among $\{1,\ldots,l\}$. In terms of the standard coordinates $u_{ij}$ from Proposition 8.3, these variables are given by:
$$\kappa=\sum_{ij}u_{ij}\quad,\quad 
\chi_l=\sum_{i=1}^lu_{ii}$$

Generally speaking, we are interested in computing the joint law of $(\chi_l,\kappa)$. There are many interesting questions here, and as a main result on the subject, we have:

\index{Poisson law}

\begin{theorem}
The measures $\mu_k^l=law\left(\chi_l\big|\kappa=k\right)$ are given by
$$\mu_k^l=\sum_{q\geq0}\binom{k}{q}\binom{l}{q}\binom{N}{q}^{-2}\frac{(\delta_1-\delta_0)^{*q}}{q!}$$
and become Poisson $(st)$ in the $k=sN,l=tN,N\to\infty$ limit.
\end{theorem}

\begin{proof}
We can use the same method as for $S_N$, from chapter 2. Let us set:
$$\widetilde{S}_N^{(k)}=\left\{\sigma\in\widetilde{S}_N\Big|\kappa(\sigma)=k\right\}$$

By using the inclusion-exclusion principle, we obtain the following formula:
\begin{eqnarray*}
P\left(\chi_l=p\Big|\kappa=k\right)
&=&\frac{1}{|\widetilde{S}_N^{(k)}|}\binom{l}{p}\#\left\{\sigma\in\widetilde{S}_{N-p}^{(k-p)}\Big|\sigma(i)\neq i,\forall i\leq l-p\right\}\\
&=&\frac{1}{|\widetilde{S}_N^{(k)}|}\binom{l}{p}\sum_{r\geq0}(-1)^r\binom{l-p}{r}\left|\widetilde{S}_{N-p-r}^{(k-p-r)}\right|
\end{eqnarray*}

Here the index $r$, which counts the fixed points among $\{1,\ldots,l-p\}$, runs a priori up to $\min(k,l)-p$. However, since the binomial coefficient or the cardinality of the set on the right vanishes by definition at $r>\min(k,l)-p$, we can sum over $r\geq0$. We have:
$$|\widetilde{S}_N^{(k)}|=k!\binom{N}{k}^2$$

By using this and then cancelling various factorials, and grouping back into binomial coefficients, we obtain the following formula:
\begin{eqnarray*}
P\left(\chi_l=p\Big|\kappa=k\right)
&=&\frac{1}{k!\binom{N}{k}^2}\binom{l}{p}\sum_{r\geq0}(-1)^r\binom{l-p}{r}(k-p-r)!\binom{N-p-r}{k-p-r}^2\\
&=&\sum_{r\geq0}\frac{(-1)^r}{p!r!}\binom{k}{p+r}\binom{l}{p+r}\binom{N}{p+r}^{-2}
\end{eqnarray*}

We can now compute the measure itself. With $p=q-r$, we obtain:
\begin{eqnarray*}
law\left(\chi_l\Big|\kappa=k\right)
&=&\sum_{p\geq0}\sum_{r\geq0}\frac{(-1)^r}{p!r!}\binom{k}{p+r}\binom{l}{p+r}\binom{N}{p+r}^{-2}\,\delta_p\\
&=&\sum_{q\geq0}\sum_{r\geq0}\frac{(-1)^r}{(q-r)!r!}\binom{k}{q}\binom{l}{q}\binom{N}{q}^{-2}\,\delta_{q-r}\\
&=&\sum_{q\geq0}\binom{k}{q}\binom{l}{q}\binom{N}{q}^{-2}\cdot\frac{1}{q!}\sum_{r\geq0}(-1)^r\binom{q}{r}\delta_{q-r}
\end{eqnarray*}

The sum on the right being $(\delta_1-\delta_0)^{*q}$, this gives the formula in the statement. Regarding now the asymptotics, in the regime $k=sN,l=tN,N\to\infty$ from the statement, the coefficient of $(\delta_1-\delta_0)^{*q}/q!$ in the formula of $\mu_k^l$ is:
\begin{eqnarray*}
c_q
&=&\binom{k}{q}\binom{l}{q}\binom{N}{q}^{-2}
=\frac{\binom{k}{q}}{\binom{N}{q}}\cdot\frac{\binom{l}{q}}{\binom{N}{q}}\\
&\simeq&\left(\frac{k}{N}\right)^q\left(\frac{l}{N}\right)^q
=(st)^q
\end{eqnarray*}

We deduce that the Fourier transform of $\mu_k^l$ is given by:
$$F_{\mu_k^l}(y)
\simeq\sum_{q\geq0}(st)^q\frac{(e^y-1)^q}{q!}
=e^{st(e^y-1)}$$

But this is the Fourier transform of the law Poisson ($st$), and we are done.
\end{proof}

Following \cite{bsk}, let us discuss now the construction of $\widetilde{S}_N^+$, in analogy with the above considerations. Let us go back to the embedding in Proposition 8.3, namely:
$$u:\widetilde{S}_N\subset M_N(0,1)\quad,\quad 
u_{ij}(\sigma)=\delta_{i\sigma(j)}$$

The image of this embedding being formed by the matrices in $M_N(0,1)$ having at most one nonzero entry on each row and column, the matrix $u=(u_{ij})$ is ``submagic'', in the sense that its entries are projections, which are pairwise orthogonal on each row and column. In fact, Gelfand duality shows that we have the following formula:
$$C(\widetilde{S}_N)=C^*_{comm}\left((u_{ij})_{i,j=1,\ldots,N}\Big|u={\rm submagic}\right)$$

This suggests the following definition, given in \cite{bsk}:

\index{submagic matrix}
\index{quantum partial permutation}
\index{quantum semigroup}

\begin{definition}
$C(\widetilde{S}_N^+)$ is the universal $C^*$-algebra generated by the entries of a $N\times N$ submagic matrix $u$, with comultiplication and counit maps given by
$$\Delta(u_{ij})=\sum_ku_{ik}\otimes u_{kj}\quad,\quad 
\varepsilon(u_{ij})=\delta_{ij}$$
where ``submagic'' means formed of projections, which are pairwise orthogonal on rows and columns. We call $\widetilde{S}_N^+$ semigroup of quantum partial permutations of $\{1,\ldots,N\}$.
\end{definition}

Observe that the morphisms $\Delta,\varepsilon$ constructed above satisfy the usual axioms for a comultiplication and an antipode, in the bialgebra setting, namely:
$$(\Delta\otimes id)\Delta=(id\otimes \Delta)\Delta$$
$$(\varepsilon\otimes id)\Delta=(id\otimes\varepsilon)\Delta=id$$

Observe also that our bialgebra $C(\widetilde{S}_N^+)$ has a kind of ``subantipode'' map, defined by $S(u_{ij})=u_{ji}$. The axiom satisfied by this subantipode map $S$ is as follows, where $m^{(3)}$ is the triple multiplication, and $\Delta^{(2)}$ is the double comultiplication:
$$m^{(3)}(S\otimes id\otimes S)\Delta^{(2)}=S$$

As a conclusion to this discussion, the basic properties of the quantum semigroup $\widetilde{S}_N^+$ that we constructed in Definition 8.5 can be summarized as follows:

\begin{proposition}
We have maps as follows,
$$\begin{matrix}
C(S_N^+)&\leftarrow&C(\widetilde{S}_N^+)\\
\\
\downarrow&&\downarrow\\
\\
C(S_N)&\leftarrow&C(\widetilde{S}_N)
\end{matrix}
\quad \quad \quad:\quad \quad\quad
\begin{matrix}
S_N^+&\subset&\widetilde{S}_N^+\\
\\
\cup&&\cup\\
\\
S_N&\subset&\widetilde{S}_N
\end{matrix}$$
with the bialgebras at left corresponding to the quantum semigroups at right.
\end{proposition}

\begin{proof}
This is clear from the above discussion, and from the well-known fact that projections which sum up to $1$ are pairwise orthogonal. 
\end{proof}

We recall from chapter 1 that we have $S_N^+\neq S_N$ starting from $N=4$. At the semigroup level things get interesting starting from $N=2$, where, following \cite{bsk}, we have:

\begin{proposition}
We have an isomorphism as follows,
$$C(\widetilde{S}_2^+)\simeq\left\{(x,y)\in C^*(D_\infty)\oplus C^*(D_\infty)\Big|\varepsilon(x)=\varepsilon(y)\right\}$$
where $\varepsilon:C^*(D_\infty)\to\mathbb C1$ is the counit, given by the formula
$$u=\begin{pmatrix}p\oplus 0&0\oplus r\\0\oplus s&q\oplus 0\end{pmatrix}$$ where $p,q$ and $r,s$ are the standard generators of the two copies of $C^*(D_\infty)$.
\end{proposition}

\begin{proof}
Consider an arbitrary $2\times 2$ matrix formed by projections:
$$u=\begin{pmatrix}P&R\\S&Q\end{pmatrix}$$

This matrix is submagic when the following conditions are satisfied:
$$PR=PS=QR=QS=0$$

Now observe that these conditions tell us that the non-unital algebras $X=<P,Q>$ and $Y=<R,S>$ must commute, and must satisfy $xy=0$, for any $x\in X,y\in Y$. Thus, if we denote by $Z$ the universal algebra generated by two projections, we have:
$$C(\widetilde{S}_2^+)\simeq\mathbb C1\oplus Z\oplus Z$$

Now since we have $C^*(D_\infty)=\mathbb C1\oplus Z$, we obtain an isomorphism as follows:
$$C(\widetilde{S}_2^+)\simeq\left\{(\lambda+a,\lambda+b)\Big|\lambda\in\mathbb C, a,b\in Z\right\}$$

Thus, we are led to the conclusion in the statement.
\end{proof}

Let us extend now Theorem 8.4 to the free setting. We first have:

\begin{proposition}
The following two elements of $C(\widetilde{S}_N^+)$ are self-adjoint,
$$\chi=\sum_iu_{ii},\qquad\kappa=\sum_{ij}u_{ij}$$
satisfy $0\leq\chi,\kappa\leq N$, and coincide with the usual $\chi,\kappa$ on the quotient $C(\widetilde{S}_N)$.
\end{proposition}

\begin{proof}
All the above assertions are clear from definitions, and with the inequalities $0\leq\chi,\kappa\leq N$ being taken of course in an operator-theoretic sense.
\end{proof}

More generally, we can talk about truncations of the variable $\chi$ constructed above, with respect to a parameter $l\in\{1,\ldots,N\}$, which are constructed as follows:
$$\chi_l=\sum_{i=1}^lu_{ii}$$

Let us look now at Theorem 8.4. Since the algebra $C(\widetilde{S}_N^+)$ has no integration functional, we cannot talk about the joint law of $(\chi,\kappa)$. Thus, we need an alternative approach to $\mu_k^l$. For this purpose, in the classical case, we use the following simple fact:

\begin{proposition}
Any partial permutation $\sigma:X\simeq Y$ can be factorized as
$$\xymatrix@R=40pt@C=40pt
{X\ar[r]^{\sigma}\ar[d]_\gamma&Y\\\{1,\ldots,k\}\ar[r]_\beta&\{1,\ldots,k\}\ar[u]_\alpha}$$
with $\alpha,\beta,\gamma\in S_k$ being certain non-unique permutations, where $k=\kappa(\sigma)$.
\end{proposition}

\begin{proof}
We can choose any two bijections $X\simeq\{1,\ldots,k\}$ and $\{1,\ldots,k\}\simeq Y$, and then complete them up to permutations $\gamma,\alpha\in S_N$. The remaining permutation $\beta\in S_k$ is then uniquely determined by the formula $\sigma=\alpha\beta\gamma$.
\end{proof}

We can now formulate an alternative definition for the measures $\mu_k^l$. We fix $k\leq N$, and we denote by $p,q,r$ the magic matrices for $S_N,S_k,S_N$. We have:

\begin{proposition}
Consider the map $\varphi:S_N\times S_k\times S_N\to\widetilde{S}_N$, sending $(\alpha,\beta,\gamma)$ to the partial permutation $\sigma:\gamma^{-1}\{1,\ldots,k\}\simeq\alpha\{1,\ldots,k\}$ given by:
$$\sigma(\gamma^{-1}(t))=\alpha(\beta(t))$$
\begin{enumerate}
\item The image of $\varphi$ is the set $\widetilde{S}_N^{(k)}=\{\sigma\in\widetilde{S}_N|\kappa(\sigma)=k\}$.

\item The transpose of $\varphi$ is given by $\varphi^*(u_{ij})=\sum_{s,t\leq k}p_{is}\otimes q_{st}\otimes r_{tj}$.

\item $\mu_k^l$ equals the law of the variable $\varphi^*(\chi_l)\in C(S_N\times S_k\times S_N)$.
\end{enumerate}
\end{proposition}

\begin{proof}
This is an elementary statement, whose proof goes as follows:

\medskip

(1) Since $\alpha,\gamma\in S_N$, the domain and range of the associated element $\sigma\in\widetilde{S}_N$ have indeed cardinality $k$. The surjectivity follows from Proposition 8.9.

\medskip

(2) For the element $\sigma\in\widetilde{S}_N$ in the statement, we have:
\begin{eqnarray*}
u_{ij}(\sigma)=1
&\iff&\sigma(j)=i\\
&\iff&\exists t\leq k,\gamma^{-1}(t)=j,\alpha(\beta(t))=i\\
&\iff&\exists s,t\leq k,\gamma^{-1}(t)=j,\beta(t)=s,\alpha(s)=i\\
&\iff&\exists s,t\leq k,r_{tj}(\gamma)=1,q_{st}(\beta)=1,p_{is}(\alpha)=1\\
&\iff&\exists s,t\leq k,(p_{is}\otimes q_{st}\otimes r_{tj})(\alpha,\beta,\gamma)=1
\end{eqnarray*}

Now since the numbers $s,t\leq k$ are uniquely determined by $\alpha,\beta,\gamma,i,j$, if they exist, we conclude that we have the following formula:
$$u_{ij}(\sigma)=\sum_{s,t\leq k}(p_{is}\otimes q_{st}\otimes r_{tj})(\alpha,\beta,\gamma)$$

But this gives the formula in the statement, and we are done.

\medskip

(3) This comes from the fact that the map $\varphi_k:S_N\times S_k\times S_N\to\widetilde{S}_N^{(k)}$ obtained by restricting the target of $\varphi$ commutes with the normalized counting measures. At $k=N$ this follows from the well-known fact that given $(\alpha,\beta,\gamma)\in S_N\times S_N\times S_N$ random, the product $\alpha\beta\gamma\in S_N$ is random, and the general case is clear as well.
\end{proof}

The point now is that we can use the same trick, ``$\sigma=\alpha\beta\gamma$'', in the free case. The precise preliminary statement that we will need is as follows:

\begin{proposition}
Let $p,q,r$ be the magic matrices for $S_N^+,S_k^+,S_N^+$.
\begin{enumerate}
\item The matrix $U_{ij}=\sum_{s,t\leq k}p_{is}\otimes q_{st}\otimes r_{tj}$ is submagic.

\item We have a representation $\pi:C(\widetilde{S}_N^+)\to C(S_N^+\times S_k^+\times S_N^+)$, $\pi(u_{ij})=U_{ij}$.

\item $\pi$ factorizes through the algebra $C(\widetilde{S}_N^{+(k)})=C(\widetilde{S}_N^+)/<\kappa=k>$.

\item At $k=N$, this factorization $\pi_k$ commutes with the Haar functionals.
\end{enumerate}
\end{proposition}

\begin{proof}
Once again, this is an elementary statement, whose proof goes at follows:

\medskip

(1) By using the fact that $p,q,r$ are magic, we obtain:
\begin{eqnarray*}
U_{ij}U_{il}
&=&\sum_{s,t\leq k}\sum_{v,w\leq k}p_{is}p_{iv}\otimes q_{st}q_{vw}\otimes r_{tj}r_{wl}\\
&=&\sum_{s,t\leq k}\sum_{w\leq k}p_{is}\otimes q_{st}q_{sw}\otimes r_{tj}r_{wl}\\
&=&\sum_{s,t\leq k}p_{is}\otimes q_{st}\otimes r_{tj}r_{tl}\\
&=&\delta_{jl}U_{ij}
\end{eqnarray*}

The proof of $U_{ij}U_{lj}=\delta_{il}U_{ij}$ is similar, and we conclude that $U$ is submagic.

\medskip

(2) This follows from (1), and from the definition of $C(\widetilde{S}_N^+)$.

\medskip

(3) By using the fact that $p,q,r$ are magic, we obtain indeed:
\begin{eqnarray*}
\sum_{ij}U_{ij}
&=&\sum_{ij}\sum_{s,t\leq k}p_{is}\otimes q_{st}\otimes r_{tj}\\
&=&\sum_{s,t\leq k}1\otimes q_{st}\otimes 1\\
&=&k
\end{eqnarray*}

Thus the representation $\pi$ factorizes indeed through the algebra in the statement.

\medskip

(4) This is a well-known analogue of the fact that ``the product of random permutations is a random permutation'', that we already used before. Here is a representation theory proof, using Peter-Weyl theory. With $P=Proj(Fix(u^{\otimes n}))$, we have:
\begin{eqnarray*}
\int_{S_N^+\times S_N^+\times S_N^+}U_{i_1j_1}\ldots U_{i_nj_n}
&=&\sum_{st}\int_{S_N^+}p_{i_1s_1}\ldots p_{i_ns_n}\int_{S_N^+}q_{s_1t_1}\ldots q_{s_nt_n}\int_{S_N^+}r_{t_1j_1}\ldots r_{t_nj_n}\\
&=&\sum_{st}P_{i_1\ldots i_n,s_1\ldots s_n}P_{s_1\ldots s_n,t_1\ldots t_n}P_{t_1\ldots t_n,j_1\ldots j_n}\\
&=&(P^3)_{i_1\ldots i_n,j_1\ldots j_n}\\
&=&P_{i_1\ldots i_n,j_1\ldots j_n}\\
&=&\int_{S_N^+}u_{i_1j_1}\ldots u_{i_nj_n}
\end{eqnarray*}

Thus $\pi_N$ commutes indeed with the Haar functionals, and we are done.
\end{proof}

Observe that, since the variable $\kappa$ is now continuous, $0\leq\kappa\leq N$, the various algebras $C(\widetilde{S}_N^{+(k)})$ constructed above don't sum any longer up to the algebra $C(\widetilde{S}_N^+)$ itself. Thus, in a certain sense, the above measures $\mu_k^l$ encode only a part of the ``probabilistic theory'' of $\widetilde{S}_N^+$. We can however formulate a free analogue of Theorem 8.4, as follows:

\index{free Poisson law}
\index{partial permutation}

\begin{theorem}
The measures $\mu_k^l=law(\pi_k(\chi_l))$, where $\pi_k$ is defined as
$$\pi_k:C(\widetilde{S}_N^+)\to C(S_N^+\times S_k^+\times S_N^+)$$
$$u_{ij}\to\sum_{s,t\leq k}p_{is}\otimes q_{st}\otimes r_{tj}$$
become free Poisson $(st)$ in the $k=sN,l=tN,N\to\infty$ limit.
\end{theorem}

\begin{proof}
The variable that we are interested in, $\chi_k^l=\pi_k(\chi_l)$, is given by:
$$\chi_k^l=\sum_{i\leq l}\sum_{s,t\leq k}p_{is}\otimes q_{st}\otimes r_{ti}$$

By raising to the power $n$ and integrating, we obtain the following formula:
$$\int_{S_N^+\times S_k^+\times S_N^+}(\chi_k^l)^n
=\sum_{i_a\leq l}\sum_{s_a,t_a\leq k}\int_{S_N^+}p_{i_1s_1}\ldots p_{i_ns_n}\int_{S_k^+}q_{s_1t_1}\ldots q_{s_nt_n}\int_{S_N^+}r_{t_1i_1}\ldots r_{t_ni_n}$$

By using now the Weingarten formula, the above moment is:
\begin{eqnarray*}
c_n
&=&\sum_{ist}\sum_{\alpha\ldots\rho\in NC(n)}\delta_\alpha(i)\delta_\beta(s)W_{nN}(\alpha,\beta)\cdot\delta_\gamma(s)\delta_\delta(t)W_{nk}(\gamma,\delta)\cdot\delta_\varepsilon(t)\delta_\rho(i)W_{nN}(\varepsilon,\rho)\\
&=&\sum_{\alpha\ldots\rho\in NC(n)}W_{nN}(\alpha,\beta)W_{nk}(\gamma,\delta)W_{nN}(\varepsilon,\rho)\sum_{ist}\delta_\alpha(i)\delta_\beta(s)\delta_\gamma(s)\delta_\delta(t)\delta_\varepsilon(t)\delta_\rho(i)\\
&=&\sum_{\alpha\ldots\rho\in NC(n)}W_{nN}(\alpha,\beta)W_{nk}(\gamma,\delta)W_{nN}(\varepsilon,\rho)\sum_{ist}\delta_{\alpha\vee\rho}(i)\delta_{\beta\vee\gamma}(s)\delta_{\delta\vee\varepsilon}(t)\\
&=&\sum_{\alpha\ldots\rho\in NC(n)}W_{nN}(\alpha,\beta)W_{nk}(\gamma,\delta)W_{nN}(\varepsilon,\rho)\cdot l^{|\alpha\vee\rho|}k^{|\beta\vee\gamma|}k^{|\delta\vee\varepsilon|}
\end{eqnarray*}

Since in the $N\to\infty$ limit the Gram and Weingarten matrices are concentrated on the diagonal, we obtain, in the regime in the statement, as in \cite{bc2}:
\begin{eqnarray*}
c_n
&\simeq&\sum_{\alpha,\gamma,\varepsilon\in NC(n)}N^{-|\alpha|}k^{-|\gamma|}N^{-|\varepsilon|}\cdot l^{|\alpha\vee\varepsilon|}k^{|\alpha\vee\gamma|}k^{|\gamma\vee\varepsilon|}\\
&\simeq&\sum_{\alpha,\gamma,\varepsilon\in NC(n)}N^{-|\alpha|-|\gamma|-|\varepsilon|+|\alpha\vee\varepsilon|+|\alpha\vee\gamma|+|\gamma\vee\varepsilon|}\cdot s^{-|\gamma|+|\alpha\vee\gamma|+|\gamma\vee\varepsilon|}\cdot t^{|\alpha\vee\varepsilon|}\\
&\simeq&\sum_{\alpha\in NC(n)}(st)^{|\alpha|}
\end{eqnarray*}

Thus, we are led to the conclusion in the statement.
\end{proof}

As a conclusion, the operation $\widetilde{S}_N\to\widetilde{S}_N^+$ is indeed a ``correct'' liberation, agreeing with the standard liberation operation from free probability theory \cite{bep}, \cite{vdn}. 

\section*{8b. Graph symmetries}

Some interesting questions, in relation with the classical and quantum partial permutations, appear in relation with graph symmetries. We first have:

\begin{proposition}
Given a graph $X$ with $N$ vertices, and adjacency matrix $d\in M_N(0,1)$, consider its partial automorphism semigroup, given by:
$$\widetilde{G}(X)=\left\{\sigma\in\widetilde{S}_N\Big|d_{ij}=d_{\sigma(i)\sigma(j)},\ \forall i,j\in Dom(\sigma)\right\}$$
We have then the following formula, with $R=diag(R_i)$, $C=diag(C_j)$, with $R_i,C_j$ being the row and column sums of the associated submagic matrix $u$:
$$C(\widetilde{G}(X))=C(\widetilde{S}_N)\Big/\Big<R(du-ud)C=0\Big>$$
Moreover, when using the relation $du=ud$ instead of the above one, we obtain a certain semigroup $\bar{G}(X)\subset\widetilde{G}(X)$, which can be strictly smaller.
\end{proposition}

\begin{proof}
The definition of $\widetilde{G}(X)$ in the statement reformulates as follows:
$$\widetilde{G}(X)=\left\{\sigma\in\widetilde{S}_N\Big|i\sim j,\exists\,\sigma(i),\exists\,\sigma(j)\implies \sigma(i)\sim\sigma(j)\right\}$$

We have the following computations:
$$(du)_{ij}(\sigma)=\sum_kd_{ik}u_{kj}(\sigma)=\sum_{k\sim i}u_{kj}(\sigma)
=\begin{cases}
1&{\rm if}\ \sigma(j)\sim i\\
0&{\rm otherwise}
\end{cases}$$
$$(ud)_{ij}(\sigma)=\sum_ku_{ik}d_{kj}(\sigma)=\sum_{k\sim j}u_{ik}(\sigma)
=\begin{cases}
1&{\rm if}\ \sigma^{-1}(i)\sim j\\
0&{\rm otherwise}
\end{cases}$$

Here the ``otherwise'' cases include by definition the cases where $\sigma(j)$, respectively $\sigma^{-1}(i)$, is undefined. We have as well the following formulae:
$$R_i(\sigma)=\sum_ju_{ij}(\sigma)
=\begin{cases}
1&{\rm if}\ \exists\,\sigma^{-1}(i)\\
0&{\rm otherwise}
\end{cases}$$
$$C_j(\sigma)=\sum_iu_{ij}(\sigma)
=\begin{cases}
1&{\rm if}\ \exists\,\sigma(j)\\
0&{\rm otherwise}
\end{cases}$$

Now by multiplying the above formulae, we obtain the following formulae:
$$(R_i(du)_{ij}C_j)(\sigma)
=\begin{cases}
1&{\rm if}\ \sigma(j)\sim i\ {\rm and}\ \exists\,\sigma^{-1}(i)\ {\rm and}\ \exists\,\sigma(j)\\
0&{\rm otherwise}
\end{cases}$$
$$(R_i(ud)_{ij}C_j)(\sigma)
=\begin{cases}
1&{\rm if}\ \sigma^{-1}(i)\sim j\ {\rm and}\ \exists\,\sigma^{-1}(i)\ {\rm and}\ \exists\,\sigma(j)\\
0&{\rm otherwise}
\end{cases}$$

We conclude that the relations in the statement, which read $R_i(du)_{ij}C_j=R_i(ud)_{ij}C_j$, when applied to a given $\sigma\in\widetilde{S}_N$, correspond to the following condition:
$$\exists\,\sigma^{-1}(i),\ \exists\,\sigma(j)\implies [\sigma(j)\sim i\iff\sigma^{-1}(i)\sim j]$$

But with $i=\sigma(k)$, this latter condition reformulates as follows:
$$\exists\,\sigma(k),\ \exists\,\sigma(j)\implies [\sigma(j)\sim\sigma(k)\iff k\sim j]$$

Thus we must have $\sigma\in\widetilde{G}(X)$, and we obtain the presentation result for $\widetilde{G}(X)$. Regarding now the second assertion, the simplest counterexample here is simplex $X_N$, having $N$ vertices and edges everywhere. Indeed, the adjacency matrix of this simplex is $d=\mathbb I_N-1_N$, with $\mathbb I_N$ being the all-1 matrix, and so the commutation of this matrix with $u$ corresponds to the fact that $u$ must be bistochastic. Thus, $u$ must be in fact magic, and we obtain $\bar{G}(X_N)=S_N$, which is smaller than $\widetilde{G}(X_N)=\widetilde{S}_N$.
\end{proof}

With the above result in hand, we are led to the following statement:

\begin{theorem}
The following construction, with $R,C$ being the diagonal matrices formed by the row and column sums of $u$, produces a subsemigroup $\widetilde{G}^+(X)\subset\widetilde{S}_N^+$,
$$C(\widetilde{G}^+(X))=C(\widetilde{S}_N^+)\Big/\Big<R(du-ud)C=0\Big>$$
called semigroup of quantum partial automorphisms of $X$, whose classical version is $\widetilde{G}(X)$. When using $du=ud$, we obtain a semigroup $\bar{G}^+(X)\subset\widetilde{G}^+(X)$ which can be smaller.
\end{theorem}

\begin{proof}
All this is elementary, the idea being as follows:

\medskip

(1) In order to construct the comultiplication $\Delta$, consider the following elements:
$$U_{ij}=\sum_ku_{ik}\otimes u_{kj}$$

By using the fact that $u$ is submagic, we deduce that we have:
$$R^U_i(dU)_{ij}C^U_j=\Delta(R_i(du)_{ij}C_j)$$
$$R^U_i(Ud)_{ij}C^U_j=\Delta(R_i(ud)_{ij}C_j)$$

Thus we can define $\Delta$ by mapping $u_{ij}\to U_{ij}$, as desired.

\medskip

(2) Regarding now $\varepsilon$, the algebra in the statement has indeed a morphism $\varepsilon$ defined by $u_{ij}\to\delta_{ij}$, because the following relations are trivially satisfied:
$$R_i(d1_N)_{ij}C_j=R_i(1_Nd)_{ij}C_j$$

(3) Regarding now $S$, we must prove that we have a morphism $S$ given by $u_{ij}\to u_{ji}$. For this purpose, we know that with $R=diag(R_i)$ and $C=diag(C_j)$, we have:
$$R(du-ud)C=0$$

Now when transposing this formula, we obtain:
$$C^t(u^td-du^t)R^t=0$$

Since $C^t,R^t$ are respectively the diagonal matrices formed by the row sums and column sums of $u^t$, we conclude that the relations $R(du-ud)C=0$ are satisfied by the transpose matrix $u^t$, and this gives the existence of the subantipode map $S$.

\medskip

(4) The fact that we have $\widetilde{G}^+(X)_{class}=\widetilde{G}(X)$ follows from $(S_N^+)_{class}=S_N$.

\medskip

(5) Finally, the last assertion follows from the last assertion in Proposition 8.13, by taking classical versions, the simplest counterexample being the simplex. 
\end{proof}

As a first result now regarding the correspondence $X\to\widetilde{G}^+(X)$, we have:

\begin{proposition}
For any finite graph $X$ we have
$$\widetilde{G}^+(X)=\widetilde{G}^+(X^c)$$
where $X^c$ is the complementary graph.
\end{proposition}

\begin{proof}
The adjacency matrices of a graph $X$ and of its complement $X^c$ are related by the following formula, where $\mathbb I_N$ is the all-1 matrix:
$$d_X+d_{X^c}=\mathbb I_N-1_N$$

Thus, in order to establish the formula in the statement, we must prove that:
$$R_i(\mathbb I_Nu)_{ij}C_j=R_i(u\mathbb I_N)_{ij}C_j$$

For this purpose, let us recall that, the matrix $u$ being submagic, its row sums and column sums $R_i,C_j$ are projections. By using this fact, we have:
$$R_i(\mathbb I_Nu)_{ij}C_j=R_iC_jC_j=R_iC_j$$
$$R_i(u\mathbb I_N)_{ij}C_j=R_iR_iC_j=R_iC_j$$

Thus we have proved our equality, and the conclusion follows.
\end{proof}

In order to discuss now various aspects of the correspondence $X\to\widetilde{G}^+(X)$, it is technically convenient to slightly enlarge our formalism, as follows:

\begin{definition}
Associated to any complex-colored oriented graph $X$, with adjacency matrix $d\in M_N(\mathbb C)$, is its semigroup of partial automorphisms, given by
$$\widetilde{G}(X)=\left\{\sigma\in\widetilde{S}_N\Big|d_{ij}=d_{\sigma(i)\sigma(j)},\ \forall i,j\in Dom(\sigma)\right\}$$
as well as its quantum semigroup of quantum partial automorphisms, given by 
$$C(\widetilde{G}^+(X))=C(\widetilde{S}_N^+)\Big/\Big<R(du-ud)C=0\Big>$$
where $R=diag(R_i)$, $C=diag(C_j)$, with $R_i,C_j$ being the row and column sums of $u$.
\end{definition}

With this notion in hand, following the material in chapter 5, let us discuss now the color independence. Let $m,\gamma$ be the multiplication and comultiplication of $\mathbb C^N$:
$$m(e_i\otimes e_j)=\delta_{ij}e_i\quad,\quad 
\gamma(e_i)=e_i\otimes e_i$$ 

We denote by $m^{(p)},\gamma^{(p)}$ their iterations, given by the following formulae:
$$m^{(p)}(e_{i_1}\otimes\ldots\otimes e_{i_1})=\delta_{i_1\ldots i_p}e_{i_1}\quad,\quad 
\gamma^{(p)}(e_i)=e_i\otimes\ldots\otimes e_i$$ 

Our goal is to use these iterations in the semigroup case, exactly as we did in chapter 5, in  the quantum group case. We will need some technical results. Let us start with:

\begin{proposition}
We have the following formulae,
$$m^{(p)}u^{\otimes p}=um^{(p)}\quad,\quad 
u^{\otimes p}\gamma^{(p)}=\gamma^{(p)}u$$
valid for any submagic matrix $u$.
\end{proposition}

\begin{proof}
We have the following computations, which prove the first formula:
$$m^{(p)}u^{\otimes p}(e_{i_1}\otimes\ldots\otimes e_{i_p})
=\sum_je_j\otimes u_{ji_1}\ldots u_{ji_p}
=\delta_{i_1\ldots i_p}\sum_je_j\otimes u_{ji_1}$$
$$um^{(p)}(e_{i_1}\otimes\ldots\otimes e_{i_p})
=\delta_{i_1\ldots i_p}u(e_{i_1})
=\delta_{i_1\ldots i_p}\sum_je_j\otimes u_{ji_1}$$

We have as well the following computations, which prove the second formula:
$$u^{\otimes p}\gamma^{(p)}(e_i)
=u^{\otimes p}(e_i\otimes\ldots\otimes e_i)
=\sum_je_j\otimes\ldots\otimes e_j\otimes u_{ji}$$
$$\gamma^{(p)}u(e_i)
=\gamma^{(p)}\left(\sum_je_j\otimes u_{ji}\right)
=\sum_je_j\otimes\ldots\otimes e_j\otimes u_{ji}$$

Summarizing, we have proved both formulae in the statement.
\end{proof}

We will need as well a second technical result, as follows:

\begin{proposition}
We have the following formulae, with $u,m,\gamma$ being as before,
$$m^{(p)}R^{\otimes p}d^{\otimes p}\gamma^{(p)}=Rd^{\times p}\quad,\quad
m^{(p)}d^{\otimes p}C^{\otimes p}\gamma^{(p)}=d^{\times p}C$$
and with $\times$ being the componentwise, or Hadamard, product of matrices.
\end{proposition}

\begin{proof}
We have the following computations, which prove the first formula:
$$m^{(p)}R^{\otimes p}d^{\otimes p}\gamma^{(p)}(e_i)
=m^{(p)}R^{\otimes p}d^{\otimes p}(e_i\otimes\ldots\otimes e_i)
=\sum_je_j\otimes R_jd_{ji}^p$$
$$Rd^{\times p}(e_i)
=R\left(\sum_je_j\otimes d_{ji}^p\right)
=\sum_je_j\otimes R_jd_{ji}^p$$

We have as well the following computations, which prove the second formula:
$$m^{(p)}d^{\otimes p}C^{\otimes p}\gamma^{(p)}(e_i)
=m^{(p)}d^{\otimes p}(e_i\otimes\ldots\otimes e_i\otimes C_i)
=\sum_je_j\otimes d_{ji}^pC_i$$
$$d^{\times p}C(e_i)
=d^{\times p}(e_i\otimes C_i)
=\sum_je_j\otimes d_{ji}^pC_i$$

Thus, we have proved both formulae in the statement.
\end{proof}

We can now prove a key result, as follows:

\begin{proposition}
We have the following formulae, with $u,m,\gamma$ being as before,
$$m^{(p)}(Rdu)^{\otimes p}\gamma^{(p)}=Rd^{\times p}u\quad,\quad 
m^{(p)}(udC)^{\otimes p}\gamma^{(p)}=ud^{\times p}C$$
and with $\times$ being the componentwise product of matrices.
\end{proposition}

\begin{proof}
By using the formulae in Proposition 8.17 and Proposition 8.18, we get:
\begin{eqnarray*}
m^{(p)}(Rdu)^{\otimes p}\gamma^{(p)}
&=&m^{(p)}R^{\otimes p}d^{\otimes p}u^{\otimes p}\gamma^{(p)}\\
&=&m^{(p)}R^{\otimes p}d^{\otimes p}\gamma^{(p)}u\\
&=&Rd^{\times p}u
\end{eqnarray*}

Once again by using Proposition 8.17 and Proposition 8.18, we have:
\begin{eqnarray*}
m^{(p)}(udC)^{\otimes p}\gamma^{(p)}
&=&m^{(p)}u^{\otimes p}d^{\otimes p}C^{\otimes p}\gamma^{(p)}\\
&=&um^{(p)}d^{\otimes p}C^{\otimes p}\gamma^{(p)}\\
&=&ud^{\times p}C
\end{eqnarray*}

Thus, we have proved both formulae in the statement.
\end{proof}

We can now prove the color independence result, as follows:

\begin{theorem}
The quantum semigroups of quantum partial isomorphisms of finite graphs are subject to the ``independence on the colors'' formula
$$\Big[d_{ij}=d_{kl}\iff d'_{ij}=d'_{kl}\Big]\implies\widetilde{G}^+(X)=\widetilde{G}^+(X')$$
valid for any graphs $X,X'$, having adjacency matrices $d,d'$.
\end{theorem}

\begin{proof}
Given a matrix $d\in M_N(\mathbb C)$, consider its color decomposition, which is as follows, with the color components $d_c$ being by definition 0-1 matrices:
$$d=\sum_{c\in\mathbb C}c\cdot d_c$$

We want to prove that a given quantum semigroup $G$ acts on $(X,d)$ if and only if it acts on $(X,d_c)$, for any $c\in\mathbb C$. For this purpose, consider the following linear space:
$$E_u=\left\{f\in M_N(\mathbb C)\Big|Rfu=ufC\right\}$$

In terms of this space, we want to prove that we have:
$$d\in E_u\implies d_c\in E_u,\forall c\in\mathbb C$$

For this purpose, observe that we have the following implication, as a consequence of the formulae established in Proposition 8.19:
$$Rdu=udC\implies Rd^{\times p}u=ud^{\times p}C$$

We conclude that we have the following implication:
$$d\in E_u\implies d^{\times p}\in E_u,\forall p\in\mathbb N$$

But this gives the result, exactly as in \cite{ba3}, via the standard linear algebra fact that the color components $d_c$ can be obtained from the componentwise powers $d^{\times p}$.
\end{proof}

In contrast with what happens for the groups or quantum groups, in the semigroup setting we do not have a spectral decomposition result as well. To be more precise, consider as before the following linear space, associated to a submagic matrix $u$:
$$E_u=\left\{d\in M_N(\mathbb C)\Big|Rdu=udC\right\}$$

It is clear that $E_u$ is a linear space, containing 1, and stable under the adjoint operation $*$ too. We also know from Theorem 8.20 that $E_u$ is stable under color decomposition. However, $E_u$ is not stable under taking products, and so is not an algebra, in general.

\bigskip

In general, the computation of $\widetilde{G}^+(X)$ remains a very interesting question. Interesting as well is the question of generalizing all this to the infinite graph case, $|X|=\infty$, with the key remark that this might be simpler than talking about $G^+(X)$ with $|X|=\infty$.

\section*{8c. Configuration spaces}

We discuss now a number of further questions in relation with quantum partial permutations, and more specifically, a number of potential applications of $\widetilde{S}_N^+$ to questions regarding $S_N^+$ itself. Let us start with the following key notion, that will be central:

\index{diagonal algebra}

\begin{definition}
The diagonal algebra of $C(S_N^\times)$ is defined as
$$D(S_N^\times)=<u_{11},\ldots,u_{NN}>\subset C(S_N^\times)$$
with $u_{ij}$ being as usual the standard generators of the algebra $C(S_N^\times)$.
\end{definition}

In this definition $S_N^\times$ can be the usual symmetric group $S_N$, or its free version $S_N^+$, that we are mainly interested in. However, we can perform this construction as well for $\widetilde{S}_N$, or for its free version $\widetilde{S}_N^+$, and we will see that this leads to interesting conclusions.

\bigskip

Before getting into the study of diagonal algebras, let us complement Definition 8.21, which is something algebraic, with its quantum space counterpart, as follows:

\index{configuration space}

\begin{definition}
The configuration space of $S_N^\times$ is the quotient $S_N^\times\to E_N^\times$ given 
by
$$C(E_N^\times)=D(S_N^\times)$$
where $D(S_N^\times)=<u_{11},\ldots,u_{NN}>\subset C(S_N^\times)$ is the diagonal algebra of $C(S_N^\times)$.
\end{definition}

According to the above definitions, we have diagrams as follows, with the quantum groups and semigroups on the left producing the configuration spaces on the right:
$$\xymatrix@R=50pt@C=50pt{
S_N^+\ar[r]&\widetilde{S}_N^+\\
S_N\ar[u]\ar[r]&\widetilde{S}_N\ar[u]}\qquad
\xymatrix@R=28pt@C=50pt{\\ :}
\qquad
\xymatrix@R=51pt@C=50pt{
E_N^+\ar[r]&\widetilde{E}_N^+\\
E_N\ar[u]\ar[r]&\widetilde{E}_N\ar[u]}$$

As mentioned, we are mainly interested in the diagonal algebra $D(S_N^+)\subset C(S_N^+)$, and in the corresponding configuration space $S_N^+\to E_N^+$. Indeed, we know that at $N\geq4$ the irreducible representations of $S_N^+$ have the same fusion rules as $SO_3$, and in particular their characters are polynomials in the main character $\chi=\sum_iu_{ii}$, so we have:
$$C(S_N^+)_{central}\subset D(S_N^+)$$

Thus, the diagonal algebra $D(S_N^+)\subset C(S_N^+)$ and the corresponding configuration space $S_N^+\to E_N^+$ are very interesting objects, encapsulating key information regarding the quantum group $S_N^+$. However, and here comes our point, enlarging the attention to quantum semigroups is technically a good idea, as shown by the following result:

\begin{theorem}
The configuration spaces for the main quantum permutation groups and semigroups are contained into certain tori, as follows, 
$$\xymatrix@R=49pt@C=50pt{
E_N^+\ar[r]&\widetilde{E}_N^+\\
E_N\ar[u]\ar[r]&\widetilde{E}_N\ar[u]}\qquad
\xymatrix@R=28pt@C=50pt{\\ \subset}
\qquad
\xymatrix@R=50pt@C=50pt{
T_N^+\ar[r]&T_N^+\\
T_N^\circ\ar[u]\ar[r]&T_N\ar[u]}$$
with $T_N=\mathbb Z_2^N$ being the real torus, $T_N^+=\mathbb Z_2^{*N}$, and $T_N^\circ\subset T_N$ being a certain subspace, obtained by removing $N$ points. The inclusions on the left are both isomorphisms.
\end{theorem}

\begin{proof}
There are several things going on here, the idea being as follows:

\medskip

(1) Let us first look at the diagonal algebra and configuration space for $S_N$. According to our general construction from Definition 8.21, the diagonal algebra $D(S_N)$ is generated by the characteristic functions $u_{ii}$ of the following subsets of $S_N$:
$$S_N^i=\left\{\sigma\in S_N\Big|\sigma(i)=i\right\}$$

(2) Our first claim is that we have the following dimension formula:
$$\dim\left(D(S_N)\right)=2^N-N$$

Indeed, the sets $S_N^i$ are in ``almost'' generic position, leading to a dimension of $2^N$, but up to the following constraint, which lowers the dimension by $N$:
$$S_N^{i_1}\cap\ldots\cap S_N^{i_{N-1}}=S_N^1\cap\ldots\cap S_N^N$$ 

(3) Our second claim is that, with $C(E_N)=D(S_N)$, we have an embedding as follows:
$$E_N\subset T_N$$

Indeed, since the diagonal coordinates $u_{ii}$ are pairwise commuting projections, we have a certain quotient map $C^*(\mathbb Z_2^N)\to D(S_N)$, obtained via Fourier transform, which at the dual level corresponds to an embedding of finite spaces $E_N\subset T_N$, as above.

\medskip

(4) By combining now (2) and (3), we conclude that the configuration space for $S_N$ is a certain subspace $E_N=T_N^\circ\subset T_N$, obtained by removing $N$ points.

\medskip

(5) With the case of $S_N$ done, let us turn now to $\widetilde{S}_N$. Here the diagonal algebra $D(\widetilde{S}_N)$ is generated by the characteristic functions $u_{ii}$ of the following sets:
$$\widetilde{S}_N^i=\left\{\sigma\in\widetilde{S}_N\Big|\sigma(i)=i\right\}$$

Since the sets $\widetilde{S}_N^i$ are now in generic position, we have the following formula:
$$\dim\big(D(\widetilde{S}_N)\big)=2^N$$

Thus, things are simpler for $\widetilde{S}_N$ than for $S_N$. As before, we have a certain quotient map $C^*(\mathbb Z_2^N)\to D(\widetilde{S}_N)$, obtained via Fourier transform, which at the dual level corresponds to an embedding of spaces $\widetilde{E}_N\subset T_N$, which this time must be an equality.

\medskip

(6) Finally, for the quantum group $S_N^+$ and its semigroup version $\widetilde{S}_N^+$, there is nothing much elementary to be done. Since the standard generators $u_{ii}$ are still  projections, in both cases, we can only say that the configuration spaces are as follows:
$$E_N^+\subset\widetilde{E}_N^+\subset T_N^+$$

Thus, we are led to the conclusions in the statement.
\end{proof}

As a more advanced result now, complementing Theorem 8.23, we have:

\begin{theorem}
The following happen:
\begin{enumerate}
\item The inclusion $E_N^+\subset T_N^+$ is an equality at $N\geq4$.

\item The inclusion $\widetilde{E}_N^+\subset T_N^+$ is an equality at any $N\geq2$.
\end{enumerate}
\end{theorem}

\begin{proof}
This is something quite tricky, the idea being as follows:

\medskip

(1) The problem here is that of computing the following functional:
$$\varphi:C(T_N^+)\to C(E_N^+)\subset C(S_N^+)\to\mathbb C$$

But a routine study here leads to the conclusion that this functional $\varphi$ is strictly positive, giving the faithfulness of the quotient map $C(T_N^+)\to C(E_N^+)$.

\medskip

(2) We know from the above that the result holds at $N\geq4$. Thus, we are led to a study at $N=2,3$, which is something elementary, and we obtain the result.
\end{proof}

All this is quite interesting, and we will be back to such questions, regarding diagonal algebras and configuration spaces, later in this book, in chapters 9-12 below.

\section*{8d. Partial isometries}

We have seen so far that the passage from $S_N\subset S_N^+$ to $\widetilde{S}_N\subset\widetilde{S}_N^+$ leads to some interesting combinatorics, and to some potential applications as well. More on this later, on various occasions, in chapters 9-12 below, and in chapters 13-16 as well. 

\bigskip

In the remainder of this chapter we discuss one more piece of related theory, namely the twisted version of all this, involving this time isometries and quantum isometries, which are very basic noncommutative geometry objects. Our starting point will be:

\index{partial isometry}

\begin{definition}
$\widetilde{O}_N$ is the semigroup of partial linear isometries of $\mathbb R^N$,
$$\widetilde{O}_N=\left\{T:A\to B\ {\rm isometry}\Big|A,B\subset\mathbb R^N\right\}$$
with the usual composition operation for such maps, namely
$$T'T:T^{-1}(A'\cap B)\to T'(A'\cap B)$$
being the composition of $T':A'\to B'$ with $T:A\to B$.
\end{definition} 

As a first remark, $\widetilde{O}_N$ is indeed a semigroup, with respect to the operation in the statement, and this is best seen in the matrix model picture, as follows:

\begin{proposition}
We have an embedding $\widetilde{O}_N\subset M_N(\mathbb R)$, obtained by completing maps $T:A\to B$ into linear maps $U:\mathbb R^N\to\mathbb R^N$, by setting $U_{|A^\perp}=0$. Moreover:
\begin{enumerate}
\item This embedding makes $\widetilde{O}_N$ correspond to the set of matrix-theoretic partial isometries, i.e. to the matrices $U\in M_N(\mathbb R)$ satisfying:
$$UU^tU=U$$

\item The semigroup operation on $\widetilde{O}_N$ corresponds in this way to the semigroup operation for matrix-theoretic partial isometries, namely:
$$U\circ V=U(U^tU\wedge VV^t)V$$
\end{enumerate}
\end{proposition}

\begin{proof}
All these assertions are well-known, and elementary. For a vector space $C=A,B$ let indeed $I_C:C\subset\mathbb R^N$ be the inclusion, and $P_C:\mathbb R^N\to C$ be the projection. The correspondence $T\leftrightarrow U$ in the statement is then given by:
$$\xymatrix@R=15mm@C=25mm{A\ar[r]^T&B\ar[d]^{I_B}\\\mathbb R^N\ar[u]^{P_A}\ar[r]_U&\mathbb R^N}\qquad\qquad
\xymatrix@R=15mm@C=25mm{A\ar[r]^T\ar[d]_{I_A}&B\\\mathbb R^N\ar[r]_U&\mathbb R^N\ar[u]_{P_B}}$$

The fact that the composition $U\circ V$ is indeed a partial isometry comes from the fact that the projections $U^tU$ and $VV^t$ are absorbed when performing the product:
$$U(U^tU\wedge VV^t)V\cdot V^t(U^tU\wedge VV^t)U^t\cdot U(U^tU\wedge VV^t)V=U(U^tU\wedge VV^t)V$$

Thus, we are led to the conclusions in the statement.
\end{proof}

In general, the multiplication formula $U\circ V=U(U^tU\wedge VV^t)V$ in Proposition 8.26 (2), while being quite complicated, is quite unavoidable. In view of some future liberation purposes, we would need a functional analytic interpretation of it. We have here:

\begin{proposition}
$C(\widetilde{O}_N)$ is the universal commutative $C^*$-algebra generated by the entries of a $N\times N$ matrix $u=(u_{ij})$ satisfying the relations
$$u=\bar{u}\quad,\quad uu^tu=u$$
with comultiplication given by the formula
$$(id\otimes\Delta)u=u_{12}(p_{13}\wedge q_{12})u_{13}=\lim_{n\to\infty}\underbrace{UU^t\ldots U^tU}_{2n+1\ {\rm terms}}$$
where $p=uu^t,q=u^tu$ and $U_{ij}=\sum_ku_{ik}\otimes u_{kj}$.
\end{proposition}

\begin{proof}
The presentation assertion is standard, by using the Gelfand and Stone-Weierstrass theorems. Let us find now the comultiplication map of $C(\widetilde{O}_N)$. For this purpose, consider the following canonical isomorphism:
$$\Phi:C(\widetilde{O}_N)\otimes C(\widetilde{O}_N)\to C(\widetilde{O}_N\times \widetilde{O}_N)$$

Consider as well the following map:
$$L_{ij}(U,V)=(U\circ V)_{ij}$$

With these conventions, the comultiplication map of $C(\widetilde{O}_N)$ is given by:
$$\Delta(u_{ij})=\Phi^{-1}(L_{ij})$$

In order to write now the map $L_{ij}$ in tensor product form, we can use:
$$P\wedge Q=\lim_{n\to\infty}(PQ)^n$$

More precisely, with $P=VV^t$ and $Q=U^tU$, we obtain the following formula:
$$(U\circ V)_{ij}
=\sum_{kl}U_{ik}(P\wedge Q)_{kl}V_{lj}
=\lim_{n\to\infty}\sum_{kl}U_{kl}(PQ)^n_{kl}V_{lj}$$

With $a_0=k,a_{2n}=l$, and by expanding the product, we obtain:
\begin{eqnarray*}
(U\circ V)_{ij}
&=&\lim_{n\to\infty}\sum_{a_0\ldots a_{2n}}U_{ia_0}P_{a_0a_1}Q_{a_1a_2}\ldots P_{a_{2n-2}a_{2n-1}}Q_{a_{2n-1}a_{2n}}V_{a_{2n}j}\\
&=&\lim_{n\to\infty}\sum_{a_0\ldots a_{2n}}U_{ia_0}Q_{a_1a_2}\ldots Q_{a_{2n-1}a_{2n}}\cdot P_{a_0a_1}\ldots P_{a_{2n-2}a_{2n-1}}V_{a_{2n}j}
\end{eqnarray*}

Now by getting back to $\Delta(u_{ij})=\Phi^{-1}(L_{ij})$, with $L_{ij}(U,V)=(U\circ V)_{ij}$, we conclude that we have the following formula, with $p=uu^t$ and $q=u^tu$:
$$\Delta(u_{ij})=\lim_{n\to\infty}\sum_{a_0\ldots a_{2n}}u_{ia_0}q_{a_1a_2}\ldots q_{a_{2n-1}a_{2n}}\otimes p_{a_0a_1}\ldots p_{a_{2n-2}a_{2n-1}}u_{a_{2n}j}$$

Let us expand now both matrix products $p=uu^t$ and $q=u^tu$. In terms of the element $U_{ij}=\sum_ku_{ik}\otimes u_{kj}$ in the statement, the sum on the right, say $S_{ij}^{(n)}$, becomes:
\begin{eqnarray*}
S_{ij}^{(n)}
&=&\sum_{a_s}u_{ia_0}(u^tu)_{a_1a_2}\ldots(u^tu)_{a_{2n-1}a_{2n}}\otimes(uu^t)_{a_0a_1}\ldots(uu^t)_{a_{2n-2}a_{2n-1}}u_{a_{2n}j}\\
&=&\sum_{a_sb_sc_s}u_{ia_0}u_{b_1a_1}u_{b_1a_2}\ldots u_{b_na_{2n-1}}u_{b_na_{2n}}\otimes u_{a_0c_1}u_{a_1c_1}\ldots u_{a_{2n-2}c_n}u_{a_{2n-1}c_n}u_{a_{2n}j}\\
&=&\sum_{b_sc_s}U_{ic_1}U_{b_1c_1}U_{b_1c_2}\ldots U_{b_nc_n}U_{b_nj}\\
&=&(\underbrace{UU^t\ldots U^tU}_{2n+1\ {\rm terms}})_{ij}
\end{eqnarray*}

Thus we have obtained the second formula in the statement. Regarding now the first formula, observe that we have $U=u_{12}u_{13}$. This gives:
\begin{eqnarray*}
\underbrace{UU^t\ldots U^tU}_{2n+1\ {\rm terms}}
&=&(u_{12}u_{13})(u_{13}^tu_{12}^t)\ldots(u_{13}^tu_{12}^t)(u_{12}u_{13})\\
&=&u_{12}(u_{13}u_{13}^t)(u_{12}^tu_{12})\ldots(u_{13}u_{13}^t)(u_{12}^tu_{12})u_{13}\\
&=&u_{12}p_{13}q_{12}\ldots p_{13}q_{12}u_{13}
\end{eqnarray*}

But this gives the first formula in the statement, and we are done.
\end{proof}

Let us construct now the liberations. We have here the following definition:

\index{quantum partial isometry}

\begin{definition}
To any $N\in\mathbb N$ we associate the following algebra,
$$C(\widetilde{O}_N^+)=C^*\left((u_{ij})_{i,j=1,\ldots,N}\Big|u_{ij}=u_{ij}^*,uu^t=p={\rm projection}\right)$$
and we call the underlying object $\widetilde{O}_N^+$ space of quantum partial isometries.
\end{definition}

As a first result regarding these liberations, we have:

\begin{proposition}
We have embeddings of compact quantum spaces
$$\xymatrix@R=40pt@C=40pt
{O_N^+\ar[r]&\widetilde{O}_N^+\\
O_N\ar[u]\ar[r]&\widetilde{O}_N\ar[u]}$$
and the spaces on the right produce the compact those on the left by dividing by the relations $p=p'=q=q'=1$, where $q=u^*u$ and $q'=u^t\bar{u}$, as in Definition 8.28.
\end{proposition}

\begin{proof}
It follows from definitions that we have embeddings as above. Regarding now the second assertion, in the case of $\widetilde{O}_N^+$, the relations $p=p'=q=q'=1$ read:
$$uu^t=\bar{u}u^t=u^tu=u^t\bar{u}=1$$

We deduce that both $u,u^t$ are unitaries, and so when dividing by these relations we obtain the quantum group $O_N^+$. As for the result regarding the classical versions, this is clear too, by dividing by the commutation relations $ab=ba$.
\end{proof}

Let us discuss now the multiplicative structure. We have here:

\begin{proposition}
$\widetilde{O}_N^+$ has a non-associative multiplication given by
$$(id\otimes\Delta)u=u_{12}(p_{13}\wedge q_{12})u_{13}=\lim_{n\to\infty}\underbrace{UU^t\ldots U^tU}_{2n+1\ {\rm terms}}$$
where $p=uu^t,q=u^tu$ and $U_{ij}=\sum_ku_{ik}\otimes u_{kj}$, compatible with that of $\widetilde{O}_N,O_N^+$.
\end{proposition}

\begin{proof}
First of all, the equality between the two matrices on the right in the statement follows as in the proof of Proposition 8.27. Let us call $W=(W_{ij})$ this matrix. In order to check that $\Delta(u_{ij})=W_{ij}$ defines indeed a morphism, we must verify that $W=(W_{ij})$ satisfies the conditions in Definition 8.28. We have:
\begin{eqnarray*}
WW^tW
&=&u_{12}(p_{13}\wedge q_{12})u_{13}\cdot u_{13}^t(p_{13}\wedge q_{12})u_{12}^t\cdot u_{12}(p_{13}\wedge q_{12})u_{13}\\
&=&u_{12}(p_{13}\wedge q_{12})p_{13}(p_{13}\wedge q_{12})p_{12}(p_{13}\wedge q_{12})u_{13}\\
&=&u_{12}(p_{13}\wedge q_{12})u_{13}=W
\end{eqnarray*}

Regarding now the last assertion, this is clear from definitions.
\end{proof}

Finally, let us discuss probabilistic aspects. We use the same method as for $\widetilde{S}_N^+$. So, pick an exponent $\circ\in\{\emptyset,+\}$, set $\kappa=\sum_{ij}u_{ij}u_{ij}^t$, and consider the following algebra:
$$C(\widetilde{O}_N^{\circ(k)})=C(\widetilde{O}_N^{\circ})/<\kappa=k>$$

With this convention, we have the following result:

\begin{proposition}
For any $\circ\in\{\emptyset,+\}$ we have a representation
$$\pi_k:C(\widetilde{O}_N^{\circ(k)})\to C(O_N^\circ\times O_k^\circ\times O_N^\circ)$$
$$\pi_k(u_{ij})=\sum_{s,t\leq k}p_{is}\otimes q_{st}\otimes r_{tj}$$
which commutes with the Haar functionals at $k=N$.
\end{proposition}

\begin{proof}
In the classical case, $\circ=\emptyset$, the first observation is that any partial isometry $T:A\to B$, with the linear spaces $A,B\subset\mathbb R^N$ having dimension $\dim(A)=\dim(B)=k$, decomposes as $T=UVW$, with $U,W\in O_N$ and $V\in O_k$:
$$\xymatrix@R=15mm@C=25mm
{A\ar[r]^T\ar[d]_W&B\\\mathbb R^k\ar[r]_V&\mathbb R^k\ar[u]_U}$$

We conclude that we have a surjection $\varphi:O_N\times O_k\times O_N\to\widetilde{O}_N^{(k)}$ mapping $(U,V,W)$ to the partial isometry $T:W^{-1}(\mathbb R^k)\to U(\mathbb R^k)$ given by:
$$T(W^{-1}x)=U(Vx)$$

By proceeding now as in the proof of Proposition 8.10 (2), we see that the transpose map $\pi=\varphi^*$ is the representation in the statement, and we are done with the classical case. In the free case, this is a routine extension of Proposition 8.10. We have indeed:
\begin{eqnarray*}
(UU^tU)_{ij}
&=&\sum_{kl}U_{ik}U_{lk}U_{lj}\\
&=&\sum_{kl}\sum_{s,t\leq k}\sum_{v,w\leq k}\sum_{y,z\leq k}p_{is}p_{lv}p_{ly}\otimes q_{st}q_{vw}q_{yz}\otimes r_{tk}r_{wk}r_{zj}\\
&=&\sum_{s,t\leq k}\sum_{y,z\leq k}p_{is}\otimes q_{st}q_{yt}q_{yz}\otimes r_{zj}\\
&=&\sum_{s,z\leq k}p_{is}\otimes q_{sz}\otimes r_{zj}\\
&=&U_{ij}
\end{eqnarray*}

Since $u_{ij}=u_{ij}^*$ implies $U_{ij}=U_{ij}^*$, this proves the partial isometry condition. Let us ckeck now that this representation vanishes on the ideal $<\kappa=k>$. We have:
\begin{eqnarray*}
\sum_{ij}U_{ij}U_{ij}^t
&=&\sum_{ij}\sum_{s,t\leq k}\sum_{v,w\leq k}p_{is}p_{iv}\otimes q_{st}q_{vw}\otimes r_{tj}r_{wj}\\
&=&\sum_{s,t\leq k}1\otimes q_{st}q_{st}\otimes 1\\
&=&k
\end{eqnarray*}

Thus we have a representation $\pi_k$ as in the statement. Finally, the last assertion is already known, from the proof of Proposition 8.10 (3).
\end{proof}

With the above result in hand, we can construct variables $\chi_k^l$ and then real probability measures $\mu_k^l$ exactly as in the discrete case, in the following way:
$$\chi_k^l=\pi_k(\chi_l)\quad,\quad 
\mu_k^l=law(\chi_k^l)$$

With these conventions, we have the following result, which is similar to the analytic liberation result obtained in the above for the liberation operation $\widetilde{S}_N\to\widetilde{S}_N^+$:

\index{liberation}
\index{quantum partial isometry}

\begin{theorem}
The operation $\widetilde{O}_N\to\widetilde{O}_N^+$ is an analytic liberation, in the sense that we have the Bercovici-Pata bijection for 
$$\mu_k^l=\pi_k(\chi_l)$$
in the $k=sN,l=tN,N\to\infty$ limit.
\end{theorem}

\begin{proof}
This follows by using standard integration technology, from \cite{bc1}, \cite{bsp}, \cite{csn}. More precisely, the Weingarten computation in the proof of Theorem 8.12 gives the following formula, in the $k=sN$, $l=tN$, $N\to\infty$ limit, where $D(n)\subset P(n)$ denotes the set of partitions associated to the quantum group $O_N^\circ$ under consideration:
$$\lim_{N\to\infty}\int_{O_N^\circ\times O_k^\circ\times O_N^\circ}(\chi_k^l)^n=\sum_{\alpha\in D(n)}(st)^{|\alpha|}$$

On the other hand, we know from \cite{bc1}, \cite{bsp}, \cite{csn} that the law of the truncated character $\chi_l$ is given by the following formula, in the $l=tN$, $N\to\infty$ limit:
$$\lim_{N\to\infty}\int_{O_N^\circ}(\chi_l)^n=\sum_{\alpha\in D(n)}t^{|\alpha|}$$

We conclude that in the $k=sN$, $l=tN$, $N\to\infty$ limit, we have:
$$\lim_{N\to\infty}\mu_k^l=\lim_{N\to\infty}\mu_N^{sl}$$

Thus, we are led to the conclusion in the statement.
\end{proof}

Summarizing, we have basic results regarding $\widetilde{O}_N\subset\widetilde{O}_N^+$, in analogy with what we know about $\widetilde{S}_N\subset\widetilde{S}_N^+$, and these results can be though of as being a ``twisted'' extension of what we know about $\widetilde{S}_N\subset\widetilde{S}_N^+$, due to the isomorphism $S_{M_N}^+=PO_N^+$. There are far more things that can be said on this subject, and generally speaking, the whole subject belongs to the relatively new and exciting ``free geometry'' area. And here, the number of things that can be done, inspired by classical geometry, is potentially infinite.
 
\section*{8e. Exercises} 

Things in this chapter have often been related to various research questions, and our exercises here will be research-level as well. First, we have:

\begin{exercise}
Construct an explicit embedding of type
$$\widetilde{S}_N\subset S_{2N}$$
and then variables over $S_{2N}$ extending our variables $\chi,\kappa$.
\end{exercise}

Here there are of course many possible answers to the first question, but assuming that the answer to this first question is found in the most straightforward way, our variables $\chi,\kappa:\widetilde{S}_N\to\mathbb N$ should then correspond to the variables $\chi_l,\chi_r:S_{2N}\to\mathbb N$ counting respectively the number of fixed points in $\{1,\ldots,N\}$, and in $\{N+1,\ldots,2N\}$.

\begin{exercise}
Reformulate all the probabilistic computations for
$$\xymatrix@R=14.5mm@C=15mm{
\widetilde{S}_N^+\ar[r]&\widetilde{O}_N^+\\
\widetilde{S}_N\ar[r]\ar[u]&\widetilde{O}_N\ar[u]}$$
in terms of suitable homogeneous spaces over the corresponding quantum groups.
\end{exercise}

The answer can be actually found in the literature, so the question is that of finding that literature, and making a brief account of it, in the cases that we are interested in.

\begin{exercise}
Develop some theory for the diagonal algebra $D(S_\infty^+)$.
\end{exercise}

To be more precise here, $S_\infty^+$ itself is not defined, and the problem is that of talking however about $D(S_\infty^+)$, based on the computations did in the above.

\part{Transitive subgroups}

\ \vskip50mm

\begin{center}
{\em I'm coming up only to hold you under

I'm coming up only to show you wrong

And to know you is hard, we wonder

To know you all wrong, we warn}
\end{center}

\chapter{Orbits, orbitals}

\section*{9a. Orbits, orbitals}

We have seen so far that several classes of subgroups $G\subset S_N$ liberate into clossed subgroups $G^+\subset S_N^+$. Our goal here, in this third part of the present book, will be to have a systematic look at the liberation operation $G\to G^+$, going beyond easiness. We will be particularly interested in the transitive subgroups $G\subset S_N$, and how they liberate into subgroups $G^+\subset S_N^+$, this being perhaps the most important case.

\bigskip

Getting started now, a useful tool for the study of the subgroups $G\subset S_N$ are the orbits of the action $G\curvearrowright\{1,\ldots,N\}$, and also the orbitals and higher orbitals of this action. In the quantum case, $G\subset S_N^+$, the theory goes back to Bichon's paper \cite{bi3}, and then to the paper of Lupini, Man\v cinska, Roberson \cite{lmr}. Following \cite{bi3}, we first have:

\index{orbits}

\begin{theorem}
Given a closed subgroup $G\subset S_N^+$, with standard coordinates denoted $u_{ij}\in C(G)$, the following defines an equivalence relation on $\{1,\ldots,N\}$,
$$i\sim j\iff u_{ij}\neq0$$
that we call orbit decomposition associated to the action $G\curvearrowright\{1,\ldots,N\}$. In the classical case, $G\subset S_N$, this is the usual orbit equivalence.
\end{theorem}

\begin{proof}
We first check the fact that we have indeed an equivalence relation. The reflexivity axiom $i\sim i$ follows by using the counit, as follows:
$$\varepsilon(u_{ii})=1
\implies u_{ii}\neq0$$

The symmetry axiom $i\sim j\implies j\sim i$ follows by using the antipode:
$$S(u_{ji})=u_{ij}\implies[u_{ij}\neq0\implies u_{ji}\neq0]$$

As for the transitivity axiom $i\sim k,k\sim j\implies i\sim j$, this follows by using the comultiplication, and positivity. Consider indeed the following formula:
$$\Delta(u_{ij})=\sum_ku_{ik}\otimes u_{kj}$$

On the right we have a sum of projections, and we obtain from this, as desired:
\begin{eqnarray*}
u_{ik}\neq0,u_{kj}\neq0
&\implies&u_{ik}\otimes u_{kj}>0\\
&\implies&\Delta(u_{ij})>0\\
&\implies&u_{ij}\neq0
\end{eqnarray*}

Finally, in the classical case, where $G\subset S_N$, the standard coordinates are:
$$u_{ij}=\chi\left(\sigma\in G\Big|\sigma(j)=i\right)$$

Thus $u_{ij}\neq0$ means that $i,j$ must be in the same orbit, as claimed.
\end{proof}

Generally speaking, the theory from the classical case extends well to the quantum group setting, and we have in particular the following result, also from \cite{bi3}:

\index{fixed points}

\begin{theorem}
Given a closed subgroup $G\subset S_N^+$, with magic matrix $u=(u_{ij})$, consider the associated coaction map, on the space $X=\{1,\ldots,N\}$:
$$\Phi:C(X)\to C(X)\otimes C(G)\quad,\quad e_i\to\sum_je_j\otimes u_{ji}$$
The following three subalgebras of $C(X)$ are then equal,
$$Fix(u)=\left\{\xi\in C(X)\Big|u\xi=\xi\right\}$$
$$Fix(\Phi)=\left\{\xi\in C(X)\Big|\Phi(\xi)=\xi\otimes1\right\}$$
$$Fix(\sim)=\left\{\xi\in C(X)\Big|i\sim j\implies \xi_i=\xi_j\right\}$$
where $\sim$ is the orbit equivalence relation constructed in Theorem 9.1.
\end{theorem}

\begin{proof}
The fact that we have $Fix(u)=Fix(\Phi)$ is standard, with this being valid for any corepresentation $u=(u_{ij})$. Indeed, we first have the following computation:
\begin{eqnarray*}
\xi\in Fix(u)
&\iff&u\xi=\xi\\
&\iff&(u\xi)_j=\xi_j,\forall j\\
&\iff&\sum_iu_{ji}\xi_i=\xi_j,\forall j
\end{eqnarray*}

On the other hand, we have as well the following computation:
\begin{eqnarray*}
\xi\in Fix(\Phi)
&\iff&\Phi(\xi)=\xi\otimes1\\
&\iff&\sum_i\Phi(e_i)\xi_i=\xi\otimes1\\
&\iff&\sum_{ij}e_j\otimes u_{ji}\xi_i=\sum_je_j\otimes\xi_j\\
&\iff&\sum_iu_{ji}\xi_i=\xi_j,\forall j
\end{eqnarray*}

Thus we have $Fix(u)=Fix(\Phi)$, as claimed. Regarding now the equality of this algebra with $Fix(\sim)$, observe first that given a vector $\xi\in Fix(\sim)$, we have:
\begin{eqnarray*}
\sum_iu_{ji}\xi_i
&=&\sum_{i\sim j}u_{ji}\xi_i\\
&=&\sum_{i\sim j}u_{ji}\xi_j\\
&=&\sum_iu_{ji}\xi_j\\
&=&\xi_j
\end{eqnarray*}

Thus $\xi\in Fix(u)=Fix(\Phi)$. Finally, for the reverse inclusion, we know from Theorem 9.1 that the magic unitary $u=(u_{ij})$ is block-diagonal, with respect to the orbit decomposition there. But this shows that the algebra $Fix(u)=Fix(\Phi)$ decomposes as well with respect to the orbit decomposition, so in order to prove the result, we are left with a study in the transitive case. More specifically we must prove that if the action is transitive, then $u$ is irreducible, and this being clear, we obtain the result. See \cite{bi3}.
\end{proof}

We have as well a useful analytic result, as follows:

\begin{theorem}
Given a closed subgroup $G\subset S_N^+$, the matrix
$$P_{ij}=\int_Gu_{ij}$$
is the orthogonal projection onto $Fix(\sim)$, and determines the orbits of $G\curvearrowright\{1,\ldots,N\}$.
\end{theorem}

\begin{proof}
This follows from Theorem 9.2, and from the standard fact, coming from Peter-Weyl theory, that $P$ is the orthogonal projection onto $Fix(u)$. 
\end{proof}

As a comment on the above result, we can see there an interesting relation between the orbit problematics and the Weingarten function problematics. We will be back to this key phenomenon on various occasions, in what follows.

\bigskip

Following now Lupini, Man\v cinska, Roberson \cite{lmr}, let us discuss the higher orbitals. Things are quite tricky here, and we have the following result, to start with:

\index{higher orbitals}

\begin{theorem}
Let $G\subset S_N^+$ be a closed subgroup, with magic unitary $u=(u_{ij})$, and let $k\in\mathbb N$. The relation 
$$(i_1,\ldots,i_k)\sim(j_1,\ldots,j_k)\iff u_{i_1j_1}\ldots u_{i_kj_k}\neq0$$
is then reflexive and symmetric, and is transitive at $k=1,2$. In the classical case, $G\subset S_N$, this relation is transitive at any $k\in\mathbb N$, and is the usual $k$-orbital equivalence.
\end{theorem}

\begin{proof}
This is known from \cite{lmr}, the proof being as follows:

\medskip

(1) The reflexivity of $\sim$ follows by using the counit, as follows:
\begin{eqnarray*}
\varepsilon(u_{i_ri_r})=1,\forall r
&\implies&\varepsilon(u_{i_1i_1}\ldots u_{i_ki_k})=1\\
&\implies&u_{i_1i_1}\ldots u_{i_ki_k}\neq0\\
&\implies&(i_1,\ldots,i_k)\sim(i_1,\ldots,i_k)
\end{eqnarray*}

(2) The symmetry follows by applying the antipode, and then the involution:
\begin{eqnarray*}
(i_1,\ldots,i_k)\sim(j_1,\ldots,j_k)
&\implies&u_{i_1j_1}\ldots u_{i_kj_k}\neq0\\
&\implies&u_{j_ki_k}\ldots u_{j_1i_1}\neq0\\
&\implies&u_{j_1i_1}\ldots u_{j_ki_k}\neq0\\
&\implies&(j_1,\ldots,j_k)\sim(i_1,\ldots,i_k)
\end{eqnarray*}

(3) The transitivity at $k=1,2$ is more tricky. Here we need to prove that:
$$u_{i_1j_1}\ldots u_{i_kj_k}\neq0\ ,\ u_{j_1l_1}\ldots u_{j_kl_k}\neq0\implies u_{i_1l_1}\ldots u_{i_kl_k}\neq0$$

In order to do so, we use the following formula:
$$\Delta(u_{i_1l_1}\ldots u_{i_kl_k})=\sum_{s_1\ldots s_k}u_{i_1s_1}\ldots u_{i_ks_k}\otimes u_{s_1l_1}\ldots u_{s_kl_k}$$

To be more precise, at $k=1$ the result is clear from this by positivity, and known since Theorem 9.1. At $k=2$ now, we can use the following trick, from \cite{lmr}:
\begin{eqnarray*}
(u_{i_1j_1}\otimes u_{j_1l_1})\Delta(u_{i_1l_1}u_{i_2l_2})(u_{i_2j_2}\otimes u_{j_2l_2})
&=&\sum_{s_1s_2}u_{i_1j_1}u_{i_1s_1}u_{i_2s_2}u_{i_2j_2}\otimes u_{j_1l_1}u_{s_1l_1}u_{s_2l_2}u_{j_2l_2}\\
&=&u_{i_1j_1}u_{i_2j_2}\otimes u_{j_1l_1}u_{j_2l_2}
\end{eqnarray*}

Indeed, we obtain from this the following implication, as desired:
$$u_{i_1j_1}u_{i_2j_2}\neq0,u_{j_1l_1}u_{j_2l_2}\neq0\implies u_{i_1l_1}u_{i_2l_2}\neq0$$

(4) Finally, assume that we are in the classical case, $G\subset S_N$. We have:
$$u_{ij}=\chi\left(\sigma\in G\Big|\sigma(j)=i\right)$$

But this formula shows that we have the following equivalence:
$$u_{i_1j_1}\ldots u_{i_kj_k}\neq0\iff\exists \sigma\in G,\ \sigma(i_1)=j_1,\ldots,\sigma(i_k)=j_k$$

In other words, $(i_1,\ldots,i_k)\sim(j_1,\ldots,j_k)$ happens precisely when $(i_1,\ldots,i_k)$ and $(j_1,\ldots,j_k)$ are in the same $k$-orbital of $G$, and this gives the last assertion.
\end{proof}

The above result raises the question about what exactly happens at $k=3$, in relation with transitivity, and the answer here, due to McCarthy \cite{mc1}, is as follows:

\begin{theorem}
There are closed subgroups $G\subset S_N^+$, as for instance the Kac-Paljutkin quantum group $G\subset S_4^+$, for which $\sim$ is not transitive at $k=3$.
\end{theorem}

\begin{proof}
This is something quite technical, from \cite{mc1}, and we will come back to this, with details, later on, when talking about $S_4^+$ and its subgroups.
\end{proof}

In view of the above results, and as a conclusion to all this, what we have on the higher orbitals, we can only formulate a modest definition, as follows:

\index{algebraic orbitals}

\begin{definition}
Given a closed subgroup $G\subset S_N^+$, consider the relation defined by: 
$$(i_1,\ldots,i_k)\sim_k(j_1,\ldots,j_k)\iff u_{i_1j_1}\ldots u_{i_kj_k}\neq0$$
\begin{enumerate}
\item The equivalence classes with respect to $\sim_1$ are called orbits of $G$.

\item The equivalence classes with respect to $\sim_2$ are called orbitals of $G$.

\item If $\sim_k$ with $k\geq3$ is transitive, we call its equivalence classes $k$-orbitals of $G$.
\end{enumerate}
\end{definition}

We should mention here that having this definition, which might look a bit deceiving, should not be a source of worries, because in practice, orbits and orbitals are what we need. As for the higher orbitals, existing or not, we will be back to them later.

\section*{9b. Group duals}

As an application of the above general theory, still following \cite{bi3}, let us discuss now the group duals $\widehat{\Gamma}\subset S_N^+$. We first have the following result:

\begin{theorem}
Given a quotient group $\mathbb Z_{N_1}*\ldots*\mathbb Z_{N_k}\to\Gamma$, we have an embedding $\widehat{\Gamma}\subset S_N^+$, with $N=N_1+\ldots+N_k$, having the following properties:
\begin{enumerate}
\item This embedding appears by diagonally joining the embeddings $\widehat{\mathbb Z_{N_k}}\subset S_{N_k}^+$, and the corresponding magic matrix has blocks of sizes $N_1,\ldots,N_k$.

\item The equivalence relation on $X=\{1,\ldots,N\}$ coming from the orbits of the action $\widehat{\Gamma}\curvearrowright X$ appears by refining the partition $N=N_1+\ldots+N_k$.
\end{enumerate}
\end{theorem}

\begin{proof}
This is something elementary, the idea being as follows:

\medskip

(1) Given a quotient group $\mathbb Z_{N_1}*\ldots*\mathbb Z_{N_k}\to\Gamma$, we have indeed a standard embedding as follows, with $N=N_1+\ldots+N_k$, that we actually know well since chapter 1:
\begin{eqnarray*}
\widehat{\Gamma}
&\subset&\widehat{\mathbb Z_{N_1}*\ldots*\mathbb Z_{N_k}}
=\widehat{\mathbb Z_{N_1}}\,\hat{*}\,\ldots\,\hat{*}\,\widehat{\mathbb Z_{N_k}}\\
&\simeq&\mathbb Z_{N_1}\,\hat{*}\,\ldots\,\hat{*}\,\mathbb Z_{N_k}
\subset S_{N_1}\,\hat{*}\,\ldots\,\hat{*}\,S_{N_k}\\
&\subset&S_{N_1}^+\,\hat{*}\,\ldots\,\hat{*}\,S_{N_k}^+
\subset S_N^+
\end{eqnarray*}

(2) Regarding the magic matrix, our claim is that this is as follows, $F_N=\frac{1}{\sqrt{N}}(w_N^{ij})$ with $w_N=e^{2\pi i/N}$ being Fourier matrices, and $g_l$ being the standard generator of $\mathbb Z_{N_l}$:
$$u=\begin{pmatrix}
F_{N_1}I_1F_{N_1}^*\\
&\ddots\\
&&F_{N_k}I_kF_{N_k}^*
\end{pmatrix}
\quad,\quad 
I_l=\begin{pmatrix}
1\\
&g_l\\
&&\ddots\\
&&&g_l^{N_l-1}
\end{pmatrix}$$

(3) Indeed, let us recall that the magic matrix for $\mathbb Z_N\subset S_N\subset S_N^+$ is given by:
$$v_{ij}
=\chi\left(\sigma\in\mathbb Z_N\Big|\sigma(j)=i\right)
=\delta_{i-j}$$

Let us apply now the Fourier transform. According to our Pontrjagin duality conventions from chapter 1, we have a pair of inverse isomorphisms, as follows:
$$\Phi:C(\mathbb Z_N)\to C^*(\mathbb Z_N)\quad,\quad\delta_i\to\frac{1}{N}\sum_kw^{ik}g^k$$
$$\Psi:C^*(\mathbb Z_N)\to C(\mathbb Z_N)\quad,\quad g^i\to\sum_kw^{-ik}\delta_k$$

Here $w=e^{2\pi i/N}$, and we use the standard Fourier analysis convention that the indices are $0,1,\ldots,N-1$. With $F=\frac{1}{\sqrt{N}}(w^{ij})$ and $I=diag(g^i)$ as above, we have:
\begin{eqnarray*}
u_{ij}
&=&\Phi(v_{ij})\\
&=&\frac{1}{N}\sum_kw^{(i-j)k}g^k\\
&=&\frac{1}{N}\sum_kw^{ik}g^kw^{-jk}\\
&=&(FIF^*)_{ij}
\end{eqnarray*}

Thus, the magic matrix that we are looking for is $u=FIF^*$, as claimed.

\medskip

(4) Finally, the second assertion in the statement is clear from the fact that $u$ is block-diagonal, with blocks corresponding to the partition $N=N_1+\ldots+N_k$.
\end{proof}

As a first comment on the above result, not all group dual subgroups $\widehat{\Gamma}\subset S_N^+$ appear exactly as above, a well-known counterexample here being the Klein group:
$$K=\mathbb Z_2\times\mathbb Z_2\subset S_4\subset S_4^+$$

Indeed, with $K=\{1,a,b,c\}$, where $c=ab$, consider the embedding $K\subset S_4$ given by $a=(12)(34)$, $b=(13)(24)$, $c=(14)(23)$. The corresponding magic matrix is:
$$u=\begin{pmatrix}
\delta_1&\delta_a&\delta_b&\delta_c\\
\delta_a&\delta_1&\delta_c&\delta_b\\
\delta_b&\delta_c&\delta_1&\delta_a\\
\delta_c&\delta_b&\delta_a&\delta_1
\end{pmatrix}\in M_4(C(K))$$

Now since this matrix is not block-diagonal, the only choice for $K=\widehat{K}$ to appear as in Theorem 9.7 would be via a quotient map $\mathbb Z_4\to K$, which is impossible. As a second comment now on Theorem 9.7, in the second assertion there we really have a possible refining operation, as shown by the example provided by the trivial group, namely:
$$\mathbb Z_{N_1}*\ldots*\mathbb Z_{N_k}\to\{1\}$$

Summarizing, what we have in Theorem 9.7 is rather a beginning of something, and some further study is needed, in order to clarify the structure of the arbitrary group dual subgroups $\widehat{\Gamma}\subset S_N^+$, and the structure of their orbits, in connection with Theorem 9.1. As a complementary result here, regarding the Klein group, we have:

\begin{proposition}
The magic unitary for $\widehat{K}\subset S_4$ diagonalizes as
$$FuF=\begin{pmatrix}
\delta_1&0&0&0\\
0&\delta_b&0&0\\
0&0&\delta_a&0\\
0&0&0&\delta_c
\end{pmatrix}\in M_4(C^*(K))$$
where $F=F_2\otimes F_2$ is the Fourier matrix of $K$.
\end{proposition}

\begin{proof}
Consider indeed the Fourier matrix of $K=\mathbb Z_2\times\mathbb Z_2$, which is:
$$F=F_2\otimes F_2
=\frac{1}{\sqrt{2}}\begin{pmatrix}1&1\\1&-1\end{pmatrix}
\otimes\frac{1}{\sqrt{2}}\begin{pmatrix}1&1\\1&-1\end{pmatrix}
=\frac{1}{2}\begin{pmatrix}
1&1&1&1\\
1&-1&1&-1\\
1&1&-1&-1\\
1&-1&-1&1
\end{pmatrix}$$

By conjugating the magic matrix $u$ found above by this matrix $F$, we obtain:
$$FuF=\begin{pmatrix}
\delta_1+\delta_a+\delta_b+\delta_c&0&0&0\\
0&\delta_1-\delta_a+\delta_b-\delta_c&0&0\\
0&0&\delta_1+\delta_a-\delta_b-\delta_c&0\\
0&0&0&\delta_1-\delta_a-\delta_b+\delta_c
\end{pmatrix}$$

Thus, we are led via $K\simeq\widehat{K}$ to the conclusion in the statement.
\end{proof}

In order to further discuss all this, how the group dual subgroups $\widehat{\Gamma}\subset S_N^+$ and the corresponding magic matrices appear, and how these magic matrices diagonalize, plainly or over blocks, let us first enlarge the attention to the group dual subgroups $\widehat{\Gamma}\subset G$ of an arbitrary closed subgroup $G\subset U_N^+$. The theory here, coming from \cite{wo1}, is as follows: 

\begin{proposition}
Given a closed subgroup $G\subset U_N^+$ and a matrix $Q\in U_N$, we let $T_Q\subset G$ be the diagonal torus of $G$, with fundamental representation spinned by $Q$:
$$C(T_Q)=C(G)\Big/\left<(QuQ^*)_{ij}=0\Big|\forall i\neq j\right>$$
This torus is then a group dual, $T_Q=\widehat{\Lambda}_Q$, where $\Lambda_Q=<g_1,\ldots,g_N>$ is the discrete group generated by the elements $g_i=(QuQ^*)_{ii}$, which are unitaries inside $C(T_Q)$.
\end{proposition}

\begin{proof}
Since $v=QuQ^*$ is a unitary corepresentation, its diagonal entries $g_i=v_{ii}$, when regarded inside $C(T_Q)$, are unitaries, and satisfy:
$$\Delta(g_i)=g_i\otimes g_i$$

Thus $C(T_Q)$ is a group algebra, and more specifically we have $C(T_Q)=C^*(\Lambda_Q)$, where $\Lambda_Q=<g_1,\ldots,g_N>$ is the group in the statement, and this gives the result.
\end{proof}

Summarizing, associated to any closed subgroup $G\subset U_N^+$ is a whole family of tori, indexed by the unitaries $U\in U_N$. As a first result regarding these tori, we have:

\begin{proposition}
Any torus $T\subset G$ appears as follows, for a certain $Q\in U_N$:
$$T\subset T_Q\subset G$$
In other words, any torus appears inside a standard torus.
\end{proposition}

\begin{proof}
Given a torus $T\subset G$, we have an inclusion as follows:
$$T\subset G\subset U_N^+$$

On the other hand, we know that each torus $T=\widehat{\Lambda}\subset U_N^+$, coming from a group $\Lambda=<g_1,\ldots,g_N>$, has a fundamental corepresentation as follows, with $Q\in U_N$:
$$u=Q\begin{pmatrix}
g_1\\
&\ddots\\
&&g_N
\end{pmatrix}Q^*$$

But this shows that we have $T\subset T_Q$, and this gives the result.
\end{proof}

Generally speaking, the family $\{T_Q|Q\in U_N\}$ constructed above can be thought of as being a kind of ``maximal torus'' for $G\subset U_N^+$, the idea being that various algebraic or analytic properties of $G$ can be read on the tori $T_Q$. We refer here to \cite{bpa}.

\bigskip

In view of the above theory and results, a main problem that we are interested in, when discussing the group dual subgroups $\widehat{\Gamma}\subset S_N^+$, is the computation of the standard tori of $S_N^+$. And the result here, complementing Theorem 9.7, is as follows:

\begin{theorem}
For the quantum permutation group $S_N^+$, the discrete group quotient $F_N\to\Lambda_Q$ with $Q\in U_N$ comes from the following relations:
$$\begin{cases}
g_i=1&{\rm if}\ \sum_lQ_{il}\neq 0\\
g_ig_j=1&{\rm if}\ \sum_lQ_{il}Q_{jl}\neq 0\\ 
g_ig_jg_k=1&{\rm if}\ \sum_lQ_{il}Q_{jl}Q_{kl}\neq 0
\end{cases}$$
Also, given a decomposition $N=N_1+\ldots+N_k$, for the matrix $Q=diag(F_{N_1},\ldots,F_{N_k})$, where $F_N=\frac{1}{\sqrt{N}}(\xi^{ij})_{ij}$ with $\xi=e^{2\pi i/N}$ is the Fourier matrix, we obtain
$$\Lambda_Q=\mathbb Z_{N_1}*\ldots*\mathbb Z_{N_k}$$
with dual embedded into $S_N^+$ in a standard way, as in Theorem 9.7.
\end{theorem}

\begin{proof}
This can be proved by a direct computation, as follows:

\medskip

(1) Fix a unitary matrix $Q\in U_N$, and consider the following quantities:
$$\begin{cases}
c_i=\sum_lQ_{il}\\
c_{ij}=\sum_lQ_{il}Q_{jl}\\
d_{ijk}=\sum_l\bar{Q}_{il}\bar{Q}_{jl}Q_{kl}
\end{cases}$$

We write $w=QvQ^*$, where $v$ is the fundamental corepresentation of $C(S_N^+)$. Assume $X\simeq\{1,\ldots,N\}$, and let $\alpha$ be the coaction of $C(S_N^+)$ on $C(X)$. Let us set:
$$\varphi_i=\sum_l\bar{Q}_{il}\delta_l\in C(X)$$

Also, let $g_i=(QvQ^*)_{ii}\in C^*(\Lambda_Q)$. If $\beta$ is the restriction of $\alpha$ to $C^*(\Lambda_Q)$, then:
$$\beta(\varphi_i)=\varphi_i\otimes g_i$$

(2) Now recall that $C(X)$ is the universal $C^*$-algebra generated by elements $\delta_1,\ldots,\delta_N$ which are pairwise orthogonal projections. Writing these conditions in terms of the linearly independent elements $\varphi_i$ by means of the formulae $\delta_i=\sum_lQ_{il}\varphi_l$, we find that the universal relations for $C(X)$ in terms of the elements $\varphi_i$ are as follows:
$$\begin{cases}
\sum_ic_i\varphi_i=1\\
\varphi_i^*=\sum_jc_{ij}\varphi_j\\
\varphi_i\varphi_j=\sum_kd_{ijk}\varphi_k
\end{cases}$$

(3) Let $\widetilde{\Lambda}_Q$ be the group in the statement. Since $\beta$ preserves these relations, we get:
$$\begin{cases}
c_i(g_i-1)=0\\
c_{ij}(g_ig_j-1)=0\\
d_{ijk}(g_ig_j-g_k)=0
\end{cases}$$

We conclude from this that $\Lambda_Q$ is a quotient of $\widetilde{\Lambda}_Q$. On the other hand, it is immediate that we have a coaction map as follows:
$$C(X)\to C(X)\otimes C^*(\widetilde{\Lambda}_Q)$$

Thus $C(\widetilde{\Lambda}_Q)$ is a quotient of $C(S_N^+)$. Since $w$ is the fundamental corepresentation of $S_N^+$ with respect to the basis $\{\varphi_i\}$, it follows that the generator $w_{ii}$ is sent to $\widetilde{g}_i\in\widetilde{\Lambda}_Q$, while $w_{ij}$ is sent to zero. We conclude that $\widetilde{\Lambda}_Q$ is a quotient of $\Lambda_Q$. Since the above quotient maps send generators on generators, we conclude that $\Lambda_Q=\widetilde{\Lambda}_Q$, as desired.

\medskip

(4) We apply the result found in (3), with the $N$-element set $X$ there being:
$$X=\mathbb Z_{N_1}\sqcup\ldots\sqcup\mathbb Z_{N_k}$$

With this choice, we have $c_i=\delta_{i0}$ for any $i$. Also, we have $c_{ij}=0$, unless $i,j,k$ belong to the same block to $Q$, in which case $c_{ij}=\delta_{i+j,0}$, and also $d_{ijk} =0$, unless $i,j,k$ belong to the same block of $Q$, in which case $d_{ijk}=\delta_{i+j,k}$. We conclude from this that $\Lambda_Q$ is the free product of $k$ groups which have generating relations as follows:
$$g_ig_j=g_{i+j}\quad,\quad g_i^{-1}=g_{-i}$$

But this shows that our group is $\Lambda_Q=\mathbb Z_{N_1}*\ldots*\mathbb Z_{N_k}$, as stated.
\end{proof}

The above result does of course not close the discussion, and the Klein group embedding, from Proposition 9.8, remains for instance to be unified with the general embeddings from Theorem 9.7. This is something quite technical, and we will get back to this question in chapter 10 below, after studying more in detail the transitive subgroups $G\subset S_N^+$.

\section*{9c. Higher orbitals}

Let us get away now from group dual subgroups, and go back to the general theory of orbits and orbitals for the arbitary closed subgroups $G\subset S_N^+$. As a complement to the results that we have, whose conclusions were summarized in Definition 9.6, we will discuss now an analytic approach to the higher orbitals, which is particularly useful when $\sim_k$ is not transitive. Let us begin with the following standard result:

\begin{proposition}
For a subgroup $G\subset S_N$, which fundamental corepresentation denoted $u=(u_{ij})$, the following numbers are equal:
\begin{enumerate}
\item The number of $k$-orbitals.

\item The dimension of space $Fix(u^{\otimes k})$.

\item The number $\int_G\chi^k$, where $\chi=\sum_iu_{ii}$.
\end{enumerate}
\end{proposition}

\begin{proof}
This is well-known, the proof being as follows:

\medskip

$(1)=(2)$ Given $\sigma\in G$ and vector $\xi=\sum_{i_1\ldots i_k}\alpha_{i_1\ldots i_k}e_{i_1}\otimes\ldots\otimes e_{i_k}$, we have:
\begin{eqnarray*}
\sigma^{\otimes k}\xi&=&\sum_{i_1\ldots i_k}\alpha_{i_1\ldots i_k}e_{\sigma(i_1)}\otimes\ldots\otimes e_{\sigma(i_k)}\\
\xi&=&\sum_{i_1\ldots i_k}\alpha_{\sigma(i_1)\ldots\sigma(i_k)}e_{\sigma(i_1)}\otimes\ldots\otimes e_{\sigma(i_k)}
\end{eqnarray*}

Thus $\sigma^{\otimes k}\xi=\xi$ holds for any $\sigma\in G$ precisely when $\alpha$ is constant on the $k$-orbitals of $G$, and this gives the equality between the numbers in (1) and (2).

\medskip

$(2)=(3)$ This follows from the Peter-Weyl theory, because $\chi=\sum_iu_{ii}$ is the character of the fundamental corepresentation $u$.
\end{proof}

In the quantum case now, where we have a closed subgroup $G\subset S_N^+$, we have:

\begin{proposition}
Given a closed subgroup $G\subset S_N^+$, and a number $k\in\mathbb N$, consider the following linear space:
$$F_k=\left\{\xi\in(\mathbb C^N)^{\otimes k}\Big|\xi_{i_1\ldots i_k}=\xi_{j_1\ldots j_k},\forall(i_1,\ldots,i_k)\sim(j_1,\ldots,j_k)\right\}$$
\begin{enumerate}
\item We have $F_k\subset Fix(u^{\otimes k})$.

\item At $k=1,2$ we have $F_k=Fix(u^{\otimes k})$.

\item In the classical case, we have $F_k=Fix(u^{\otimes k})$.

\item For $G=S_N^+$ with $N\geq4$ we have $F_3\neq Fix(u^{\otimes 3})$.
\end{enumerate}
\end{proposition}

\begin{proof}
The tensor power $u^{\otimes k}$ being the corepresentation $(u_{i_1,\ldots i_k,j_1\ldots j_k})_{i_1\ldots i_k,j_1\ldots j_k}$, the corresponding fixed point space $Fix(u^{\otimes k})$ consists of the vectors $\xi$ satisfying:
$$\sum_{j_1\ldots j_k}u_{i_1j_1}\ldots u_{i_kj_k}\xi_{j_1\ldots j_k}=\xi_{i_1\ldots i_k}\quad,\quad\forall i_1,\ldots,i_k$$

With this formula in hand, the proof goes as follows:

\medskip

(1) Assuming $\xi\in F_k$, the above fixed point formula holds indeed, because:
$$\sum_{j_1\ldots j_k}u_{i_1j_1}\ldots u_{i_kj_k}\xi_{j_1\ldots j_k}
=\sum_{j_1\ldots j_k}u_{i_1j_1}\ldots u_{i_kj_k}\xi_{i_1\ldots i_k}
=\xi_{i_1\ldots i_k}$$

(2) This is something more tricky, coming from the following formulae:
$$u_{ik}\left(\sum_ju_{ij}\xi_j-\xi_i\right)=u_{ik}(\xi_k-\xi_i)$$
$$u_{i_1k_1}\left(\sum_{j_1j_2}u_{i_1j_1}u_{i_2j_2}\xi_{j_1j_2}-\xi_{i_1i_2}\right)u_{i_2k_2}=u_{i_1k_1}u_{i_2k_2}(\xi_{k_1k_2}-\xi_{i_1i_2})$$

(3) This follows indeed from Proposition 9.12.

\medskip

(4) This follows from the representation theory of $S_N^+$ with $N\geq4$, and from some elementary computations, the dimensions of the two spaces involved being $4<5$. To be more precise, let us start with the symmetric group $S_N$. It follows from definitions that the $k$-orbitals are indexed by the partitions $\pi\in P(k)$, as follows:
$$C_\pi=\left\{(i_1,\ldots,i_k)\Big|\ker i=\pi\right\}$$

In particular at $k=3$ we have $5$ such orbitals, corresponding to:
$$\sqcap\hskip-0.7mm\sqcap\quad,\quad\sqcap\,|\quad,\quad|\,\sqcap\quad,\quad\sqcap\hskip-3.2mm{\ }_|\quad,\quad |\,|\,|$$

Regarding now $S_N^+$, the $3$-orbitals are exactly as for $S_N$, except for the fact that the $\sqcap\hskip-3.2mm{\ }_|$\ \ and $|\,|\,|$ 3-orbitals get merged. Thus, we have $4$ such orbitals, corresponding to:
$$\sqcap\hskip-0.7mm\sqcap\quad,\quad\sqcap\,|\quad,\quad|\,\sqcap\quad,\quad{\ }_|\!{\ }_|\!{\ }_|\hskip-5.75mm{\ }^{.....}$$

As for the number of analytic orbitals, this is the same as for $S_N$, namely 5.
\end{proof}

The above considerations suggest formulating the following definition:

\index{analytic orbitals}

\begin{definition}
Given a closed subgroup $G\subset U_N^+$, the integer
$$\dim Fix(u^{\otimes k})=\int_G\chi^k$$
is called number of analytic $k$-orbitals.
\end{definition}

To be more precise, in the classical case the situation is of course well understood, and this is the number of $k$-orbitals. The same goes for the general case, with $k=1,2$, where this is the number of $k$-orbitals. At $k=3$ and higher, however, even in the case where the algebraic $3$-orbitals are well-defined, their number is not necessarily the above one. In the particular case $k=3$, we have as well the following result:

\begin{proposition}
For a closed subgroup $G\subset S_N^+$, and an integer $k\leq3$, the following conditions are equivalent:
\begin{enumerate}
\item $G$ is $k$-transitive, in the sense that $Fix(u^{\otimes k})$ has dimension $1,2,5$.

\item The $k$-th moment of the main character is $\int_G\chi^k=1,2,5$.

\item We have the integration formula
$$\int_Gu_{i_1j_1}\ldots u_{i_kj_k}=\frac{(N-k)!}{N!}$$
for distinct indices $i_r$ and distinct indices $j_r$.

\item An arbitrary polynomial integral
$$\int_Gu_{i_1j_1}\ldots u_{i_kj_k}$$
equals $\frac{(N-|\ker i|)!}{N!}$ when $\ker i=\ker j$, and equals $0$, otherwise.
\end{enumerate}
\end{proposition}

\begin{proof}
Most of these implications are well-known, the idea being as follows:

\medskip

$(1)\iff(2)$ This follows from the Peter-Weyl type theory from \cite{wo1}, because the $k$-th moment of the character counts the number of fixed points of $u^{\otimes k}$.

\medskip

$(2)\iff(3)$ This follows from the Schur-Weyl duality results for $S_N,S_N^+$ and from $P(k)=NC(k)$ at $k\leq3$. 

\medskip

$(3)\iff(4)$ Once again this follows from $P(k)=NC(k)$ at $k\leq3$, and from a standard integration result for $S_N$, stating that we have:
$$\int_{S_N}u_{i_1j_1}\ldots u_{i_kj_k}=\delta_{\ker i,\ker j}\frac{(N-|\ker i|)!}{N!}$$

Thus, we are led to the conclusions in the statement.
\end{proof}

As a conclusion to all these considerations, we have:

\begin{theorem}
For a closed subgroup $G\subset S_N^+$, and an integer $k\in\mathbb N$, the number $\dim(Fix(u^{\otimes k}))=\int_G\chi^k$ of ``analytic $k$-orbitals'' has the following properties:
\begin{enumerate}
\item In the classical case, this is the number of $k$-orbitals.

\item In general, at $k=1,2$, this is the number of $k$-orbitals.

\item At $k=3$, when this number is minimal, $G$ is $3$-transitive in the above sense.
\end{enumerate}
\end{theorem}

\begin{proof}
This follows indeed from the above considerations.
\end{proof}

We will be back to all this later, with a number of more powerful tools.

\section*{9d. Finite subgroups}

We discuss now an alternative take on the above questions, in the finite quantum group case. Let us start with the following standard definition:

\index{finite quantum group}

\begin{definition}
Associated to any finite quantum group $F$ is its dual finite quantum group $G=\widehat{F}$, given by $C(G)=C(F)^*$, with Hopf $C^*$-algebra structure as follows:
\begin{enumerate}
\item Multiplication $(\varphi\psi)a=(\varphi\otimes\psi)\Delta(a)$.

\item Unit $1=\varepsilon$.

\item Involution $\varphi^*(a)=\overline{\varphi(S(a)^*)}$.

\item Comultiplication $(\Delta \varphi)(a\otimes b)=\varphi(ab)$.

\item Counit $\varepsilon(\varphi)=\varphi(1)$.

\item Antipode $(S\varphi)a=\varphi(S(a))$.
\end{enumerate}
\end{definition}

Our aim will be that of reformulating in terms of $G=\widehat{F}$ the quantum permutation group condition $F\subset S_ N^+$. We will see later how this can potentially help, by dropping the assumption that $F,G$ are finite, in connection with various quantum permutation group questions. In order to get started, we have the following well-known fact:

\begin{proposition}
Given $F$ and $G=\widehat{F}$ as in Definition 9.17, the formula
$$\pi:C(G)\to M_N(\mathbb C)\quad,\quad 
\varphi\to[\varphi(u_{ij})]_{ij}$$
defines a $*$-algebra representation precisely when $u$ is a corepresentation.
\end{proposition}

\begin{proof}
In one sense, the fact that $\pi$ is multiplicative follows from the fact that $u$ is comultiplicative, in the sense that $\Delta(u_{ij})=\sum_ku_{ik}\otimes u_{kj}$, as follows:
\begin{eqnarray*}
\pi(\varphi\psi)
&=&[(\varphi\psi)u_{ij}]_{ij}\\
&=&[(\varphi\otimes\psi)\Delta(u_{ij})]_{ij}\\
&=&\Big[\sum_k\varphi(u_{ik})\psi(u_{kj})\Big]_{ij}\\
&=&[\varphi(u_{ij})]_{ij}[\psi(u_{ij})]_{ij}\\
&=&\pi(\varphi)\pi(\psi)
\end{eqnarray*}

The fact that the morphism $\pi$ constructed above is unital is clear, coming from the fact that $u$ is counital, in the sense that $\varepsilon(u_{ij})=\delta_{ij}$, as follows:
$$\pi(\varepsilon)
=[\varepsilon(u_{ij})]_{ij}
=1$$

Regarding now the fact that $\pi$ is involutive, observe first that we have:
$$\varphi^*(u_{ij})
=\overline{\varphi(S(u_{ij})^*)}
=\overline{\varphi(u_{ji})}$$

Thus, we can prove that $\pi$ is indeed involutive, as follows, using the fact that $u$ is coinvolutive, in the sense that $S(u_{ij})=u_{ji}^*$, as follows:
$$\pi(\varphi^*)
=[\varphi^*(u_{ij})]_{ij}
=[\overline{\varphi(u_{ji})}]_{ij}
=\Big[[\varphi(u_{ij})]_{ij}\Big]^*
=\pi(\varphi)^*$$

Finally, the proof in the other sense follows from exactly the same computations.
\end{proof}

In order to reach now to the condition $F\subset S_N^+$, we must impose several conditions on the matrix $u=(u_{ij})$. Let us start with the bistochasticity condition. We have here:

\begin{proposition}
Given $F$ and $G=\widehat{F}$ as in Proposition 9.18, the matrix $u=(u_{ij})$ is bistochastic, in the sense that all the row and column sums are $1$, precisely when the associated $*$-algebra representation $\pi:C(G)\to M_N(\mathbb C)$ satisfies the conditions
$$\pi(\varphi)\xi=\varphi(1)\xi\quad,\quad 
\pi(\varphi)^t\xi=\varphi(1)\xi$$
where $\xi\in\mathbb C^N$ is the all-one vector.
\end{proposition}

\begin{proof}
We want the following two conditions to be verified:
$$\sum_ju_{ij}=1\quad,\quad 
\sum_iu_{ij}=1$$ 

In what regards the condition $\sum_ju_{ij}=1$, observe that in terms of $\pi$, we have:
$$\sum_j\pi(\varphi)_{ij}
=\sum_j\varphi(u_{ij})
=\varphi\left(\sum_ju_{ij}\right)$$

Thus, we want this quantity to be $\varphi(1)$, for any $i$, and this leads to the condition $\pi(\varphi)\xi=\varphi(1)\xi$ in the statement. As for the second condition, namely $\sum_iu_{ij}=1$, this leads to the second condition in the statement, namely $\pi(\varphi)^t\xi=\varphi(1)\xi$.
\end{proof}

Independently of the above result, we must impose the condition that the coordinates $u_{ij}$ are self-adjoint. The result here is as follows:

\begin{proposition}
Given $F$ and $G=\widehat{F}$ as in Proposition 9.18, we have $u_{ij}=u_{ij}^*$ precisely when the associated $*$-algebra representation $\pi:C(G)\to M_N(\mathbb C)$ satisfies:
$$\pi S(\varphi)=\pi(\varphi)^t$$
\end{proposition}

\begin{proof}
According to formula $(S\varphi)a=\varphi(S(a))$ from Definition 9.17, we have:
$$\pi S(\varphi)
=[S\varphi(u_{ij})]_{ij}
=[\varphi(u_{ji}^*)]_{ij}$$

With this formula in hand, we see that the condition $u_{ij}=u_{ij}^*$ means that this latter matrix should be $[\varphi(u_{ji})]_{ij}=\pi(\varphi)^t$, as claimed.
\end{proof}

Let us put now what we have together. We are led to the following statement:

\begin{proposition}
Given $F$ and $G=\widehat{F}$ as in Proposition 9.18, $u=(u_{ij})$ is bistochastic, with self-adjoint entries, precisely when associated $*$-algebra representation
$$\pi:C(G)\to M_N(\mathbb C)\quad,\quad 
\varphi\to[\varphi(u_{ij})]_{ij}$$
satisfying the following conditions,
$$\pi(\varphi)\xi=\varphi(1)\xi\quad,\quad 
\pi(\varphi)^t\xi=\varphi(1)\xi\quad,\quad 
\pi S(\varphi)=\pi(\varphi)^t$$
with $\xi\in\mathbb C^N$ being the all-one vector.
\end{proposition}

\begin{proof}
This follows indeed from Proposition 9.19 and Proposition 9.20.
\end{proof}

In order to reach now to $F\subset S_N^+$, we must impose one final condition, stating that the entries of $u=(u_{ij})$ are idempotents, $u_{ij}^2=u_{ij}$. This is something more technical:

\begin{proposition}
Given $F$ and $G=\widehat{F}$ as in Proposition 9.21, we have $u_{ij}^2=u_{ij}$ precisely when the associated $*$-algebra representation $\pi:C(G)\to M_N(\mathbb C)$ satisfies
$$m(\pi\otimes\pi)\Delta(\varphi)=\pi(\varphi)m$$
as an equality of maps $\mathbb C^N\otimes\mathbb C^N\to\mathbb C^N$, where $m$ is the multiplication of $\mathbb C^N$.
\end{proposition}

\begin{proof}
We have the following computation, by using the Sweedler notation:
\begin{eqnarray*}
m(\pi\otimes\pi)\Delta(\varphi)(e_i\otimes e_j)
&=&m(\pi\otimes\pi)\left(\sum\varphi_1\otimes\varphi_2\right)(e_i\otimes e_j)\\
&=&m\left(\sum\sum_{kl}\varphi_1(u_{ki})e_k\otimes\varphi_2(u_{lj})e_l\right)\\
&=&\sum\sum_k\varphi_1(u_{ki})\varphi_2(u_{kj})e_k\\
&=&\sum\sum_k(\varphi_1\otimes\varphi_2)(u_{ki}\otimes u_{kj})e_k\\
&=&\sum_k\Delta(\varphi)(u_{ki}\otimes u_{kj})e_k\\
&=&\sum_k\varphi(u_{ki}u_{kj})e_k
\end{eqnarray*}

On the other hand, we have as well the following computation:
\begin{eqnarray*}
\pi(\varphi)m(e_i\otimes e_j)
&=&\pi(\varphi)\delta_{ij}e_i\\
&=&\big[\varphi(u_{ij})\big]_{ij}\delta_{ij}e_i\\
&=&\delta_{ij}\sum_k\varphi(u_{ki})e_k
\end{eqnarray*}

Thus, the condition in the statement simply reads $u_{ki}u_{kj}=\delta_{ij}u_{ki}$, for any $i,j,k$. In particular with $i=j$ we obtain, as desired, the idempotent condition:
$$u_{ki}^2=u_{ki}$$

Conversely now, if this idempotent condition is satisfied, then $u=(u_{ij})$ follows to be a matrix of projections, which is bistochastic. Thus this matrix is magic, and so we have $u_{ki}u_{kj}=\delta_{ij}u_{ki}$ for any $i,j,k$, and this leads to the formula in the statement.
\end{proof}

Let us put now what we have together. We are led to the following statement:

\begin{theorem}
Given $F$ and $G=\widehat{F}$ as in Definition 9.17, we have $F\subset S_N^+$, with associated magic matrix $u=(u_{ij})$, precisely when we have a $*$-algebra representation
$$\pi:C(G)\to M_N(\mathbb C)\quad,\quad 
\varphi\to[\varphi(u_{ij})]_{ij}$$
satisfying the following conditions,
$$\pi(\varphi)\xi=\varphi(1)\xi\quad,\quad 
\pi(\varphi)^t\xi=\varphi(1)\xi$$
$$\pi S(\varphi)=\pi(\varphi)^t\quad,\quad 
m(\pi\otimes\pi)\Delta(\varphi)=\pi(\varphi)m$$
where $\xi\in\mathbb C^N$ is the all-one vector, and $m$ is the multiplication of $\mathbb C^N$.
\end{theorem}

\begin{proof}
This follows indeed from Proposition 9.21 and Proposition 9.22, and from the well-known fact, already mentioned in the proof of Proposition 9.22, that a magic matrix $u=(u_{ij})$ is the same as a matrix of projections which is bistochastic.
\end{proof}

As a first illustration, in the classical case, we have:

\begin{proposition}
Given a closed subgroup $F\subset U_N$, the associated $*$-algebra representation constructed in Theorem 9.23 is given by
$$\pi:C^*(F)\to M_N(\mathbb C)\quad,\quad 
\sum_g\lambda_gg\to\sum_g\lambda_gg$$
and we have $F\subset S_N$ precisely when the conditions in Theorem 9.18 are satisfied.
\end{proposition}

\begin{proof}
Here the first assertion is clear from definitions. As for the second assertion, this is something that we know from Theorem 9.23, but here is a direct check as well:

\medskip

(1) For $\varphi\in C^*(F)$ given by $\varphi=\sum_g\lambda_gg$ we have $\pi(\varphi)=\sum_g\lambda_gg$, and also $\varphi(1)=\sum_g\lambda_g$ via $C^*(F)\simeq C(F)^*$ so the bistochasticity condition $F\subset C_N$ corresponds indeed to the conditions $\pi(\varphi)\xi=\varphi(1)\xi$ and $\pi(\varphi)^t\xi=\varphi(1)\xi$ from Theorem 9.23.

\medskip

(2) Once again with $\varphi=\sum_g\lambda_gg$, we have the following formulae:
$$\pi S\varphi=\pi\left(\sum_g\lambda_gg^{-1}\right)=\sum_g\lambda_gg^{-1}$$
$$\pi(\varphi)^t=\left(\sum_g\lambda_gg\right)^t=\sum_g\lambda_gg^t$$

Thus $F\subset O_N$, which is the same as saying that $g^{-1}=g^t$, for any $g\in F$, is indeed equivalent to the condition $\pi S(\varphi)=\pi(\varphi)^t$ from Theorem 9.23.

\medskip

(3) As before with $\varphi=\sum_g\lambda_gg$, assuming $F\subset S_N$, we have the following formula:
\begin{eqnarray*}
m(\pi\otimes\pi)\Delta(\varphi)(e_i\otimes e_j)
&=&m\left(\sum_g\lambda_gg\otimes g\right)(e_i\otimes e_j)\\
&=&m\left(\sum_g\lambda_g e_{g(i)}\otimes e_{g(j)}\right)\\
&=&\delta_{ij}\sum_{g\in G}\lambda_ge_{g(i)}
\end{eqnarray*}

On the other hand, we have as well the following formula:
\begin{eqnarray*}
\pi(\varphi)m(e_i\otimes e_j)
&=&\left(\sum_g\lambda_gg\right)m(e_i\otimes e_j)\\
&=&\left(\sum_g\lambda_gg\right)(\delta_{ij}e_i)\\
&=&\delta_{ij}\sum_g\delta_ge_{g(i)}
\end{eqnarray*}

Thus the condition $m(\pi\otimes\pi)\Delta(\varphi)=\pi(\varphi)m$ in Theorem 9.23 must be indeed satisfied, and the proof of the converse is similar, using the same computations.
\end{proof}

In the group dual case now, the result is as follows:

\index{group dual}

\begin{theorem}
Given a finite group $G$, and setting $F=\widehat{G}$, the associated $*$-algebra representation constructed in Theorem 9.23 appears as follows, for a certain family of generators $g_1,\ldots,g_N\in H$, and for a certain unitary $U\in U_N$,
$$\pi:C(G)\to M_N(\mathbb C)\quad,\quad 
\varphi\to U
\begin{pmatrix}
g_1\\
&\ddots\\
&&g_N
\end{pmatrix}U^*$$
and we have $F\subset S_N^+$ precisely when the conditions in Theorem 9.23 are satisfied, which in turn mean that the representation $\pi$ appears as in Theorem 9.7.
\end{theorem}

\begin{proof}
Here the first assertion is standard, coming from Woronowicz's Peter-Weyl type theory from \cite{wo1}. As for the second assertion, since the algebra $C(G)$ is commutative, its matrix representation $\pi$ must appear diagonally, spinned by a unitary, as claimed.
\end{proof}

There are many things that can be done with finite quantum permutation groups, in the present dual formalism. In order to discuss this, let us start with:

\begin{proposition}
Given $G=\widehat{F}$ as in Theorem 9.23, the $*$-algebra representation
$$\pi:C(G)\to M_N(\mathbb C)$$ 
gives rise to a family of $*$-algebra representations as follows, for any $k\in\mathbb N$,
$$\pi^k:C(G)\to M_N(\mathbb C)^{\otimes k}\quad,\quad 
\pi^k=\pi^{\otimes k}\Delta^{(k)}$$
that we will still denote by $\pi$, when there is no confusion, which are given by
$$\pi_{i_1\ldots i_k,j_1\ldots j_k}(\varphi)=\varphi(u_{i_1j_1}\ldots u_{i_kj_k})$$
in standard multi-index notation for the elements of $M_N(\mathbb C)^{\otimes k}$.
\end{proposition}

\begin{proof}
We have the following computation, in Sweedler notation:
\begin{eqnarray*}
\pi_{i_1\ldots i_k,j_1\ldots j_k}(\varphi)
&=&<\pi^{\otimes k}\Delta^{(k)}(\varphi)(e_{j_1}\otimes\ldots\otimes e_{j_k}),e_{i_1}\otimes\ldots\otimes e_{i_k}>\\
&=&\left<\pi^{\otimes k}\left(\sum\varphi_1\otimes\ldots\otimes\varphi_k\right)(e_{j_1}\otimes\ldots\otimes e_{j_k}),e_{i_1}\otimes\ldots\otimes e_{i_k}\right>\\
&=&\sum<\pi(\varphi_1)e_{j_1},e_{i_1}>\ldots<\pi(\varphi_k)e_{j_k},e_{i_k}>\\
&=&\sum\varphi_1(u_{i_1j_1})\ldots\varphi_k(u_{i_kj_k})\\
&=&\Delta^{(k)}(\varphi)(u_{i_1j_1}\ldots u_{i_kj_k})\\
&=&\varphi(u_{i_1j_1}\ldots u_{i_kj_k})
\end{eqnarray*}

Thus, we are led to the conclusion in the statement.
\end{proof}

Following \cite{bc1} and related papers, and also \cite{mc1}, let us discuss as well integration questions, in the present dual setting. We have here the following result:

\begin{theorem}
The polynomial integrals over $G$ are given by
$$\left[\int u_{i_1j_1}\ldots u_{i_kj_k}\right]_{i_1\ldots i_k,j_1\ldots j_k}=\pi^{\otimes k}\Delta^{(k)}(\smallint)$$
and the moments of the main character $\chi=\sum_iu_{ii}$ are given by
$$\int \chi^k=Tr(\pi^{\otimes k}\Delta^{(k)}(\smallint))$$
where $\smallint\in C(G)$ is the Haar integration functional.
\end{theorem}

\begin{proof}
The first formula is clear from Proposition 9.26. Regarding now the moments of the main character, observe that we have the following general formula:
$$Tr(\pi^{\otimes k}\Delta^{(k)}(\varphi))
=\sum_{i_1\ldots i_k}\varphi(u_{i_1i_1}\ldots u_{i_ki_k})
=\varphi(\chi^k)$$

In particular, with $\varphi=\int$, we are led to the formula in the statement.
\end{proof}

As a second topic, let us discuss the orbit and orbital theory. We first have:

\begin{theorem}
The orbits of $F\subset S_N^+$ can be defined dually by $i\sim j$ when 
$$\pi_{ij}(\varphi)>0$$
for a certain positive linear form $\varphi>0$.
\end{theorem}

\begin{proof}
We know from \cite{bi3} that $i\sim j$ when $u_{ij}\neq0$ is an equivalence relation on $\{1,\ldots,N\}$. Here is a proof of this fact, using our present, dual formalism:

\medskip

(1) The reflexivity of $\sim$ as defined in the statement is clear, coming from:
$$\pi(1)=1\implies\pi_{ii}\neq0$$

(2) The symmetry is clear too, coming from $\pi S(\varphi)=\pi(\varphi)^t$. Alternatively:
$$\pi(\varphi^*)=\pi(\varphi)^*\implies\pi_{ij}(\varphi^*)=\overline{\pi_{ij}(\varphi)}$$

(3) Regarding now the transitivity, things are a bit more tricky. We have:
$$\pi_{ij}(\varphi\psi)=\sum_k\pi_{ik}(\varphi)\pi_{kj}(\psi)$$

Now since $\varphi\geq0$ implies $\varphi(u_{ij})\geq0$ for any $i,j$, we obtain the result.
\end{proof}

Regarding the orbitals, following \cite{lmr}, we have:

\begin{proposition}
The relation on $\{1,\ldots,N\}^k$ given by $i\sim j$ when
$$\pi_{i_1\ldots i_k,j_1\ldots j_k}(\varphi)>0$$
for a certain positive linear form $\varphi>0$, is reflexive and symmetric.
\end{proposition}

\begin{proof}
The reflexivity is clear exactly as at $k=1$, coming from:
$$\pi(1)=1\implies\pi_{i_1\ldots i_k,i_1\ldots i_k}\neq0$$

The symmetry is clear too, coming from $\pi S(\varphi)=\pi(\varphi)^t$. Alternatively:
$$\pi(\varphi^*)=\pi(\varphi)^*\implies\pi_{i_1\ldots i_k,j_1\ldots j_k}(\varphi^*)=\overline{\pi_{i_1\ldots i_k,j_1\ldots j_k}(\varphi)}$$

Thus, we are led to the conclusion in the statement.
\end{proof}

Finally, let us discuss the actions $F\curvearrowright X$ of the finite quantum groups on the finite graphs, formulated in our dual setting, in terms of $G=\widehat{F}$. Let us start with:

\begin{definition}
We say that a compact quantum group $G$ acts dually on $X$ when 
$$[Im(\pi),d]=0$$
where $d\in M_N(\mathbb C)$ is the adjacency matrix of $X$. Also, when the condition
$$Im(\pi)=\{d\}'$$
is satisfied, we say that $G$ is the dual quantum automorphism group of $X$.
\end{definition}

Observe also that we have assumed here $d\in M_N(\mathbb C)$, which means that our finite graph $X$ can be oriented, with the edges colored by complex numbers, as in \cite{ba3} and in chapter 5. Regarding now the general theory, in the dual setting, we first have:

\begin{proposition}
The following happen, in regards with the usual graphs:
\begin{enumerate}
\item The empty graph is always invariant.

\item The same goes for the complete graph.

\item $X$ and its complement $X^c$ have the same quantum automorphism group.
\end{enumerate}
\end{proposition}

\begin{proof}
The empty graph, having adjacency matrix $d=0_N$, is obviously invariant. The same goes for the complete graph, having adjacency matrix $d=\mathbb I_N-1_N$, due to our bistochasticity conditions from our main set of axioms, namely: 
$$\pi(\varphi)\xi=\varphi(1)\xi\quad,\quad 
\pi(\varphi)^t\xi=\varphi(1)\xi$$

More generally, we obtain from axioms, and more specifically from the above bistochasticity conditions, that we have invariance under complementation, as claimed.
\end{proof}

Among the other basic results from \cite{ba3} is the fact that the quantum actions on graphs are stable under spectral and color decomposition, at the level of the associated adjacency matrices. We can easily recover these key results with our dual formalism, as follows:

\begin{theorem}
The quantum actions on graphs are stable under the following operations, at the level of the associated adjacency matrices:
\begin{enumerate}
\item Spectral decomposition.

\item Color decomposition.
\end{enumerate}
\end{theorem}

\begin{proof}
We have two results to be proved, the idea being as follows:

\medskip

(1) The spectral decomposition result is clear by definition, because the spectral projections of the adjacency matrix $d$ belong to the algebra $<d>$, that the algebra $Im(\pi)$ must commute with, according to Definition 9.30. 

\medskip

(2) Regarding now the color decomposition, following \cite{ba3}, we want to prove that we have a commutation relation as follows, for any color $c\in\mathbb C$:
$$[Im(\pi),d_c]=0$$

By a standard analytic argument, this is equivalent to the following fact, which must be valid for any $n\in\mathbb N$, with $\times$ being the componentwise, of Hadamard product:
$$[Im(\pi),d^{\times n}]=0$$

In order to establish this formula, observe that we have:
\begin{eqnarray*}
[\pi(\varphi),d]=0
&\implies&[\pi(\varphi_1)\otimes\ldots\otimes\pi(\varphi_n),d^{\otimes n}]=0\\
&\implies&[\pi^{\otimes n}\Delta^{(n)}(\varphi),d^{\otimes n}]=0
\end{eqnarray*}

We have the following formula:
\begin{eqnarray*}
[\pi^n(\varphi)\cdot d^{\otimes n}]_{i_1\ldots i_n,j_1\ldots j_n}
&=&\sum_{k_1\ldots k_n}(\pi^n(\varphi))_{i_1\ldots i_n,k_1\ldots k_n}(d^{\otimes n})_{k_1\ldots k_n,j_1\ldots j_n}\\
&=&\sum_{k_1\ldots k_n}\varphi(u_{i_1k_1}\ldots u_{i_nk_n})d_{k_1j_1}\ldots d_{k_nj_n}
\end{eqnarray*}

We have as well the following formula:
\begin{eqnarray*}
[d^{\otimes n}\cdot\pi^n(\varphi)]_{i_1\ldots i_n,j_1\ldots j_n}
&=&\sum_{k_1\ldots k_n}(d^{\otimes n})_{i_1\ldots i_n,k_1\ldots k_n}(\pi^n(\varphi))_{k_1\ldots k_n,j_1\ldots j_n}\\
&=&\sum_{k_1\ldots k_n}d_{i_1k_1}\ldots d_{i_nk_n}\varphi(u_{k_1j_1}\ldots u_{k_nj_n})
\end{eqnarray*}

Now observe that with $i_1=\ldots=i_k=i$ in the first formula, we obtain:
\begin{eqnarray*}
[\pi^n(\varphi)\cdot d^{\otimes n}]_{i\ldots i,j\ldots j}
&=&\sum_{k_1\ldots k_n}\varphi(u_{ik_1}\ldots u_{ik_n})d_{k_1j}\ldots d_{k_nj}\\
&=&\sum_k\varphi(u_{ik})d^n_{kj}\\
&=&[\pi(\varphi)\cdot d^{\times n}]_{ij}
\end{eqnarray*}

Also, with $j_1=\ldots=j_k=i$ in the first formula, we obtain:
\begin{eqnarray*}
[d^{\otimes n}\cdot\pi^n(\varphi)]_{i\ldots i,j\ldots j}
&=&\sum_{k_1\ldots k_n}d_{ik_1}\ldots d_{ik_n}\varphi(u_{k_1j}\ldots u_{k_nj})\\
&=&\sum_kd^n_{ik}\varphi(u_{kj})\\
&=&[d^{\times n}\cdot\pi(\varphi)]_{ij}
\end{eqnarray*}

Thus, we obtain $[\pi(\varphi),d^{\times n}]=0$, and by arguing as in \cite{ba3}, we obtain the result.
\end{proof}

\section*{9e. Exercises} 

We have seen in this chapter how to deal with the arbitrary subgroups $G\subset S_N^+$, by using orbitals, group duals, and related techniques. As a first exercise, we have:

\begin{exercise}
Formulate and study a notion of $k$-transitivity for the subgroups
$$G\subset S_N^+$$
by requiring that the action $G\curvearrowright\{1,\ldots,N\}$ has a minimal number of $k$-orbitals.
\end{exercise}

As before, there are many things that can be done here, in analogy with the theory of the higher transitive subgroups $G\subset S_N$, but in the quantum case some of the things will work only at $k=1$, or $k\leq2$, or $k\leq 3$. We will be back to all this.

\begin{exercise}
Extend the theory of the duals $G=\widehat{F}$ of the finite subgroups
$$F\subset S_N^+$$
developed here, notably by solving the orbital problem formulated above.
\end{exercise}

As with the previous exercises, this is something quite tricky, partly going into unexplored territory, and the more the results, the better.

\chapter{Transitive subgroups}

\section*{10a. Transitivity}

In this chapter we restrict the attention to the transitive subgroups $G\subset S_N^+$. Let us first review the basic theory here, that we will need in what follows. The notion of transitivity, which goes back to Bichon's paper \cite{bi3}, can be introduced as follows:

\index{orbits}
\index{transitivity}

\begin{definition}
Let $G\subset S_N^+$ be a closed subgroup, with magic unitary $u=(u_{ij})$, and consider the equivalence relation on $\{1,\ldots,N\}$ given by $i\sim j\iff u_{ij}\neq0$.
\begin{enumerate}
\item The equivalence classes under $\sim$ are called orbits of $G$.

\item $G$ is called transitive when the action has a single orbit. 
\end{enumerate}
In other words, we call a subgroup $G\subset S_N^+$ transitive when $u_{ij}\neq0$, for any $i,j$.
\end{definition}

This notion of transitivity is standard, coming in a straightforward way from the orbit theory from chapter 9. In the classical case, we obtain of course the usual notion of transitivity. We have the following result, once again coming from \cite{bi3}:

\begin{proposition}
For a closed subgroup $G\subset S_N^+$, the following are equivalent:
\begin{enumerate}
\item $G$ is transitive.

\item $Fix(u)=\mathbb C\xi$, where $\xi$ is the all-one vector.

\item $\int_Gu_{ij}=\frac{1}{N}$, for any $i,j$.
\end{enumerate}
\end{proposition}

\begin{proof}
This is well-known in the classical case. In general, the proof is as follows:

\medskip

$(1)\iff(2)$ We use the standard fact that the fixed point space of a corepresentation coincides with the fixed point space of the associated coaction:
$$Fix(u)=Fix(\Phi)$$

As explained in chapter 9, the fixed point space of the magic corepresentation $u=(u_{ij})$ has the following interpretation, in terms of orbits:
$$Fix(u)=\left\{\xi\in C(X)\Big|i\sim j\implies \xi(i)=\xi(j)\right\}$$

In particular, the transitivity condition corresponds to $Fix(u)=\mathbb C\xi$, as stated.

\medskip

$(2)\iff(3)$ This is clear from the general properties of the Haar integration, and more precisely from the fact that $(\int_Gu_{ij})_{ij}$ is the projection onto $Fix(u)$.
\end{proof}

Here is now a list of basic examples of transitive quantum groups, that we already know, coming from the various constructions from the previous chapters:

\begin{proposition}
The following are transitive subgroups $G\subset S_N^+$:
\begin{enumerate}
\item The quantum permutation group $S_N^+$ itself.

\item The transitive subgroups $G\subset S_N$. These are the classical examples.

\item The quantum automorphism groups of transitive graphs $G^+(X)$, with $|X|=N$.

\item In particular, the hyperoctahedral quantum group $H_n^+\subset S_N^+$, with $N=2n$.

\item We have as well the twisted orthogonal group $O_n^{-1}\subset S_N^+$, with $N=2^n$.
\end{enumerate}
In addition, the class of transitive quantum permutation groups $\{G\subset S_N^+|N\in\mathbb N\}$ is stable under direct products $\times$, wreath products $\wr$ and free wreath products $\wr_*$.
\end{proposition}

\begin{proof}
All these assertions are elementary, the idea being as follows:

\medskip

(1) This comes from the fact that we have an inclusion $S_N\subset S_N^+$. Indeed, since $S_N$ is transitive, so must be $S_N^+$, because its coordinates $u_{ij}$ map to those of $S_N$.

\medskip

(2) This is again something trivial. Indeed, for a classical group $G\subset S_N$, the variables $u_{ij}=\chi(\sigma\in S_N|\sigma(j)=i)$ are all nonzero precisely when $G$ is transitive.

\medskip

(3) This is trivial, because $X$ being transitive means that $G(X)\curvearrowright X$ is transitive, and by definition of $G^+(X)$, we have $G(X)\subset G^+(X)$.

\medskip

(4) This comes from the result from \cite{bbc}, stating that we have $H_n^+=G^+(I_n)$, where $I_n$ is the graph formed by $n$ segments, having $N=2n$ vertices.

\medskip

(5) Once again this comes from a result from \cite{bbc}, stating that we have $O_n^{-1}=G^+(\square_n)$, where $\square_n$ is the $n$-dimensional hypercube, having $N=2^n$ vertices.

\medskip

(6) Finally, the stability assertion is clear from the definition of the various products involved, from \cite{bi2}, \cite{wa1}. This is well-known, and we will be back later to this.
\end{proof}

Let us study now the transitive subgroups $G\subset S_N^+$. As a first result here, in the classical case the situation is very simple, as follows:

\index{transitive action}
\index{quotient space}

\begin{proposition}
Let $G$ be a finite group.
\begin{enumerate}
\item Assuming that we have a transitive action $G\curvearrowright\{1,\ldots,N\}$, by setting $H=\left\{\sigma\in G|\sigma(1)=1\right\}$, we have an identification $G/H=\{1,\ldots,N\}$.

\item Conversely, any $H\subset G$ produces an action $G\curvearrowright G/H$ given by $g(hH)=(gh)H$, so a morphism $G\to S_N$, with $N=[G:H]$.

\item This latter morphism is injective when the subgroup $H\subset G$ satisfies the condition $hgh^{-1}\in H,\forall h\in G\implies g=1$.

\item In the case where $G\subset S_N$ is abelian and transitive, the subgroup $H\subset G$ is trivial, $H=\{1\}$, and so we have $|G|=N$.
\end{enumerate}
\end{proposition}

\begin{proof}
All the above assertions are well-known and standard, coming from the definition of the quotient space $G/H$, as being the space of cosets $gH$.
\end{proof}

In the quantum case it is quite unclear how to generalize the above result, which fails as stated. However, we can at least try to extend the obvious fact that $G=N|H|$ must be a multiple of $N$. And here, we have the following result, from \cite{bch}:

\begin{theorem}
If $G\subset S_N^+$ is finite and transitive, then $N$ divides $|G|$. Moreover:
\begin{enumerate}
\item The case $|G|=N$ comes from the classical finite groups, of order $N$, acting on themselves.

\item The case $|G|=2N$ is possible, in the non-classical setting, an example here being the Kac-Paljutkin quantum group, at $N=4$.
\end{enumerate} 
\end{theorem}

\begin{proof}
We use the standard coaction of $C(G)$ on $\mathbb C^N=C(1,\ldots,N)$, given by:
$$\Phi:\mathbb C^N\to \mathbb C^N\otimes C(G)\quad,\quad 
e_i\to\sum_je_j\otimes u_{ji}$$

For $a\in\{1,\ldots,N\}$ consider the evaluation map $ev_a:\mathbb C^N\to\mathbb C$ at $a$. By composing $\Phi$ with $ev_a\otimes id$ we obtain a $C(G)$-comodule map, as follows:
$$I_a:\mathbb C^N\to C(G)\quad,\quad 
e_i\to u_{ia}$$

Our transitivity assumption on $G$ ensures that $I_a$ is injective. Thus, we have realized $\mathbb C^N$ as a coideal subalgebra of $C(G)$. We recall now that a finite dimensional Hopf algebra is free as a module over a coideal subalgebra $A$ provided that the latter is Frobenius, in the sense that there exists a non-degenerate bilinear form $b:A\otimes A\to \mathbb{C}$ satisfying:
$$b(xy,z)=b(x,yz)$$

We can apply this result to the coideal subalgebra $I_a(\mathbb{C}^N)\subset C(G)$, with the remark that $\mathbb{C}^N$ is indeed Frobenius, with bilinear form as follows:
$$b(fg)=\frac{1}{N}\sum_if(i)g(i)$$

Thus $C(G)$ is a free module over the $N$-dimensional algebra $\mathbb{C}^N$, and this gives the result. Regarding now the remaining assertions, the proof here goes as follows:

\medskip

(1) Since $C(G)=<u_{ij}>$ is of dimension $N$, and its commutative subalgebra $<u_{i1}>$ is of dimension $N$ already, $C(G)$ must be commutative. Thus $G$ must be classical, and by transitivity, the inclusion $G\subset S_N$ must come from the action of $G$ on itself.

\medskip

(2) The closed subgroups $G\subset S_4^+$ were classified in \cite{bb2}, and among them we have indeed the Kac-Paljutkin quantum group, which satisfies $|G|=8$, and is transitive.
\end{proof}

There are many interesting questions in relation with the above considerations, with a main one being the classification of the transitive group duals $\widehat{\Gamma}\subset S_N^+$. As a main example here, which is quite illustrating, we have the Klein group $K\subset S_4\subset S_4^+$.

\bigskip

Following \cite{bfr}, let us discuss now a useful of extension of the notion of transitivity, that we will need later, in relation with matrix modeling questions, as follows:

\index{quasi-transitivity}

\begin{definition}
A quantum permutation group $G\subset S_N^+$ is called quasi-transitive when all its orbits have the same size. In other words:
\begin{enumerate}
\item The usual orbit equivalence relation $\sim$, given by $i\sim j\iff u_{ij}\neq0$, must have its equivalence classes of the same size.

\item Equivalently, the binary matrix $\varepsilon\in M_N(0,1)$ given by $\varepsilon_{ij}=\delta_{u_{ij},0}$ must be block-diagonal, with flat matrices of same size as blocks.
\end{enumerate}
\end{definition}

As a first example, if $G$ is transitive then it is of course quasi-transitive, trivially. In general now, if we denote by $K\in\mathbb N$ the common size of the blocks, and by $M\in\mathbb N$ their multiplicity, then we must have $N=KM$. We have the following result:

\begin{proposition}
Assuming that $G\subset S_N^+$ is quasi-transitive, we must have
$$G\subset\underbrace{S_K^+\,\hat{*}\,\ldots\,\hat{*}\,S_K^+}_{M\ terms}$$
where $K\in\mathbb N$ is the common size of the orbits, and $M\in\mathbb N$ is their number.
\end{proposition}

\begin{proof}
This follows indeed from definitions, because for a quasi-transitive subgroup $G\subset S_N^+$, with orbits having common size $K|N$, the corresponding magic matrix is by definition block-diagonal, with the common size of the blocks being $K$.
\end{proof}

Let us discuss now the examples. Assume that we are in the following situation:
$$G\subset S_K^+\,\hat{*}\,\ldots\,\hat{*}\,S_K^+$$ 

If $u,v$ are the fundamental corepresentations of $C(S_N^+),C(S_K^+)$, consider the quotient map $\pi_i:C(S_N^+)\to C(S_K^+)$ constructed as follows: 
$$u\to diag(1_K,\ldots,1_K,\underbrace{v}_{i-th\ term},1_K,\ldots,1_K)$$

We can then set $C(G_{i}) = \pi_{i}(C(G))$, and we have the following result:

\begin{proposition}
If $G_{i}$ is transitive for all $i$, then $G$ is quasi-transitive.
\end{proposition}

\begin{proof}
We know that we have embeddings as follows:
$$G_1\times\ldots\times G_M\subset G\subset\underbrace{S_K^+\,\hat{*}\,\ldots\,\hat{*}\,S_K^+}_{M\ terms}$$

It follows that the size of any orbit of $G$ is at least $K$, because it contains $G_1\times\ldots\times G_M$, and at most $K$, because it is contained in $S_K^+\,\hat{*}\,\ldots\,\hat{*}\,S_K^+$. Thus, $G$ is quasi-transitive.
\end{proof}

We call the quasi-transitive subgroups appearing as above ``of product type''. There are quasi-transitive groups which are not of product type, as for instance:
$$G=S_2\subset S_2\times S_2\subset S_4\quad,\quad
\sigma\to(\sigma,\sigma)$$

Indeed, the quasi-transitivity is clear, say by letting $G$ act on the vertices of a square. On the other hand, since we have $G_1=G_2=\{1\}$, this group is not of product type. In general, we can construct examples by using various product operations:

\begin{proposition}
Given transitive subgroups $G_1,\ldots,G_M\subset S_K^+$, the following constructions produce quasi-transitive subgroups as follows, of product type:
$$G\subset\underbrace{S_K^+\,\hat{*}\,\ldots\,\hat{*}\,S_K^+}_{M\ terms}$$
\begin{enumerate}
\item The usual product: $G=G_1\times\ldots\times G_M$.

\item The dual free product: $G=G_1\,\hat{*}\,\ldots\,\hat{*}\, G_M$.
\end{enumerate}
\end{proposition}

\begin{proof}
All these assertions are clear from definitions, because in each case, the quantum groups $G_i\subset S_K^+$ constructed before are those in the statement.
\end{proof}

In the group dual case, we have the following result:

\begin{proposition}
The group duals which are of product type 
$$\widehat{\Gamma}\subset\underbrace{S_K^+\,\hat{*}\,\ldots\,\hat{*}\,S_K^+}_{M\ terms}$$
are precisely those appearing from intermediate groups of the following type:
$$\underbrace{\mathbb Z_K*\ldots*\mathbb Z_K}_{M\ terms}\to\Gamma\to\underbrace{\mathbb Z_K\times\ldots\times\mathbb Z_K}_{M\ terms}$$
\end{proposition}

\begin{proof}
In one sense, this is clear. Conversely, consider a group dual $\widehat{\Gamma}\subset S_N^+$, coming from a quotient group $\mathbb Z_K^{*M}\to\Gamma$. The subgroups $G_i\subset\widehat{\Gamma}$ constructed above must be group duals as well, $G_i=\widehat{\Gamma}_i$, for certain quotient groups $\Gamma\to\Gamma_i$. Now if $\widehat{\Gamma}$ is of product type, $\widehat{\Gamma}_i\subset S_K^+$ must be transitive, and hence equal to $\widehat{\mathbb Z}_K$.  
Thus we have $\Gamma\to\mathbb Z_K^M$.
\end{proof}

In order to construct some other classes of examples, we can use the notion of normality for compact quantum groups. This notion, from \cite{dpr}, is introduced as follows:

\index{normal subgroup}

\begin{definition}
Given a quantum subgroup $H\subset G$, coming from a quotient map $\pi:C(G)\to C(H)$, the following conditions are equivalent:
\begin{enumerate}
\item The following algebra satisfies $\Delta(A)\subset A\otimes A$:
$$A=\left\{a\in C(G)\Big|(id\otimes\pi)\Delta(a)=a\otimes1\right\}$$
 
\item The following algebra satisfies $\Delta(B)\subset B\otimes B$:
$$B=\left\{a\in C(G)\Big|(\pi\otimes id)\Delta(a)=1\otimes a\right\}$$

\item We have $A=B$, as subalgebras of $C(G)$.
\end{enumerate}
If these conditions are satisfied, we say that $H\subset G$ is a normal subgroup.
\end{definition}

Now with this notion in hand, we have, following \cite{bfr}:

\begin{theorem}
Assuming that $G\subset S_N^+$ is transitive, and that $H\subset G$ is normal, $H\subset S_N^+$ follows to be quasi-transitive.
\end{theorem}

\begin{proof}
Consider the quotient map $\pi:C(G)\to C(H)$, given at the level of standard coordinates by $u_{ij}\to v_{ij}$. Consider two orbits $O_1,O_2$ of $H$ and set:
$$x_i=\sum_{j\in O_1}u_{ij}\quad,\quad
y_i=\sum_{j\in O_2}u_{ij}$$

These two elements are orthogonal projections in $C(G)$ and they are nonzero, because they are sums of nonzero projections by transitivity of $G$. We have:
\begin{eqnarray*}
(id\otimes\pi)\Delta(x_i)
&=&\sum_{k\in O_1}\sum_{j\in O_1}u_{ik}\otimes v_{kj}\\
&=&\sum_{k\in O_1}u_{ik}\otimes 1\\
&=&x_i\otimes 1
\end{eqnarray*}

Thus by normality of $H$ we have the following formula: 
$$(\pi\otimes id)\Delta(x_i)=1\otimes x_i$$

On the other hand, assuming that we have $i\in O_2$, we obtain:
$$(\pi\otimes id)\Delta(x_i)
=\sum_k\sum_{j\in O_1}v_{ik}\otimes u_{kj}
=\sum_{k\in O_2}v_{ik}\otimes x_k$$

Multiplying this by $v_{ik}\otimes 1$ with $k\in O_2$ yields $v_{ik}\otimes x_k=v_{ik}\otimes x_i$, that is to say, $x_k=x_i$. In other words, $x_i$ only depends on the orbit of $i$. The same is of course true for $y_i$. By using this observation, we can compute the following element:
$$z
=\sum_{k\in O_2}\sum_{j\in O_1}u_{kj}
=\sum_{k\in O_2}x_k
=\vert O_2\vert x_i$$

On the other hand, by applying the antipode, we have as well:
$$S(z)
=\sum_{k\in O_2}\sum_{j\in O_1}u_{jk}
=\sum_{j\in O_1}y_j
=\vert O_1\vert y_j$$

We therefore obtain the following formula:
$$S(x_i)=\frac{\vert O_1\vert}{\vert O_2\vert}y_j$$

Now since both $x_i$ and $y_j$ have norm one, we conclude that the two orbits have the same size, and this finishes the proof.
\end{proof}

We will be back to all this later, when talking about matrix models.

\section*{10b. Higher transitivity}

Let us discuss now the notion of double transitivity. Following \cite{lmr}, and the theory of orbitals developed in chapter 9, we have the following definition:

\begin{definition}
Let $G\subset S_N^+$ be a closed subgroup, with magic unitary $u=(u_{ij})$, and consider as before the equivalence relation on $\{1,\ldots,N\}^2$ given by:
$$(i,k)\sim (j,l)\iff u_{ij}u_{kl}\neq0$$ 
\begin{enumerate}
\item The equivalence classes under $\sim$ are called orbitals of $G$.

\item $G$ is called doubly transitive when the action has two orbitals. 
\end{enumerate}
In other words, we call $G\subset S_N^+$ doubly transitive when $u_{ij}u_{kl}\neq0$, for any $i\neq k,j\neq l$.
\end{definition}

To be more precise, it is clear from definitions that the diagonal $D\subset\{1,\ldots,N\}^2$ is an orbital, and it follows that its complement $D^c$ must be a union of orbitals. With this remark in hand, the meaning of (2) is that the orbitals must be $D,D^c$.

\bigskip

In order to study the above notion, we will need the following fact, also from \cite{lmr}:

\begin{theorem}
Given a closed subgroup $G\subset S_N^+$, with magic matrix $u=(u_{ij})$, consider the following vector space coaction map, where $X=\{1,\ldots,N\}$:
$$\Phi:C(X\times X)\to C(X\times X)\otimes C(G)\quad,\quad e_{ik}\to\sum_{jl}e_{jl}\otimes u_{ji}u_{lk}$$
The following three algebras are then isomorphic,
$$End(u)=\left\{d\in M_N(\mathbb C)\Big|du=ud\right\}$$
$$Fix(\Phi)=\left\{\xi\in C(X\times X)\Big|\Phi(\xi)=\xi\otimes1\right\}$$
$$Fix(\sim)=\left\{\xi\in C(X\times X)\Big|(i,k)\sim(j,l)\implies \xi_{ik}=\xi_{jl}\right\}$$
where $\sim$ is the orbital equivalence relation from Definition 10.13.
\end{theorem}

\begin{proof}
This follows indeed by doing some computations, as those from chapter 9, for the similar result there regarding the orbits, and we refer to \cite{lmr} for details.
\end{proof}

Before going further, let us point out that the above result makes a useful connection with the graph problematics, the precise statement here being as follows:

\begin{theorem}
In order for a quantum permutation group $G\subset S_N^+$ to act on a graph $X$, having $N$ vertices, the adjacency matrix $d\in M_N(0,1)$ of the graph must be, when regarded as function on the set $\{1,\ldots,N\}^2$, constant on the orbitals of $G$. 
\end{theorem}

\begin{proof}
This follows indeed from the following isomorphism, from Theorem 10.14:
$$End(u)\simeq Fix(\sim)$$

For more on all this, details, examples, and applications too, we refer to \cite{lmr}.
\end{proof}

Now back to our notion of double transitivity, from Definition 10.13, in more analytic terms, we have the following result, also from \cite{lmr}:

\begin{theorem}
For a doubly transitive subgroup $G\subset S_N^+$, we have:
$$\int_Gu_{ij}u_{kl}=\begin{cases}
\frac{1}{N}&{\rm if}\ i=k,j=l\\
0&{\rm if}\ i=k,j\neq l\ {\rm or}\ i\neq k,j=l\\
\frac{1}{N(N-1)}&{\rm if}\ i\neq k,j\neq l
\end{cases}$$
Moreover, this formula characterizes the double transitivity.
\end{theorem}

\begin{proof}
We use the standard fact, from \cite{wo1}, that the integrals in the statement form the projection onto $Fix(u^{\otimes 2})$. Now if we assume that $G$ is doubly transitive, $Fix(u^{\otimes 2})$ has dimension 2, and therefore coincides with $Fix(u^{\otimes 2})$ for the usual symmetric group $S_N$. Thus the integrals in the statement coincide with those for the symmetric group $S_N$, which are given by the above formula. Finally, the converse is clear as well.
\end{proof}

We refer to \cite{lmr} and subsequent papers for more on all this.

\bigskip

Let us discuss the notion of $k$-transitivity, at $k\in\mathbb N$. We begin our study by recalling a few standard facts regarding the symmetric group $S_N$, and its subgroups $G\subset S_N$, from a representation theory/probabilistic viewpoint. We first have the following result:

\begin{proposition}
Consider the symmetric group $S_N$, together with its standard matrix coordinates $u_{ij}=\chi(\sigma\in S_N|\sigma(j)=i)$. We have the formula
$$\int_{S_N}u_{i_1j_1}\ldots u_{i_kj_k}=\begin{cases}
\frac{(N-|\ker i|)!}{N!}&{\rm if}\ \ker i=\ker j\\
0&{\rm otherwise}
\end{cases}$$
where $\ker i$ denotes as usual the partition of $\{1,\ldots,k\}$ whose blocks collect the equal indices of $i$, and where $|.|$ denotes the number of blocks.
\end{proposition}

\begin{proof}
According to the definition of $u_{ij}$, the integrals in the statement are given by:
$$\int_{S_N}u_{i_1j_1}\ldots u_{i_kj_k}=\frac{1}{N!}\#\left\{\sigma\in S_N\Big|\sigma(j_1)=i_1,\ldots,\sigma(j_k)=i_k\right\}$$

Since the existence of $\sigma\in S_N$ as above requires $i_m=i_n\iff j_m=j_n$, this integral vanishes when $\ker i\neq\ker j$. As for the case $\ker i=\ker j$, if we denote by $b\in\{1,\ldots,k\}$ the number of blocks of this partition, we have $N-b$ points to be sent bijectively to $N-b$ points, and so $(N-b)!$ solutions, and the integral is $\frac{(N-b)!}{N!}$, as claimed.
\end{proof}

We recall now that each action $G\curvearrowright\{1,\ldots,N\}$ produces an action $G\curvearrowright\{1,\ldots,N\}^k$ for any $k\in\mathbb N$, and by restriction, $G$ acts on the following set:
$$I_N^k=\left\{(i_1,\ldots,i_k)\in\{1,\ldots,N\}^k\Big|i_1,\ldots,i_k\ {\rm distinct}\right\}$$

We have the following well-known result:

\begin{theorem}
Given a subgroup $G\subset S_N$, with standard matrix coordinates denoted $u_{ij}=\chi(\sigma|\sigma(j)=i)$, and a number $k\leq N$, the following are equivalent:
\begin{enumerate}
\item $G$ is $k$-transitive, in the sense that $G\curvearrowright I_N^k$ is transitive.

\item $Fix(u^{\otimes k})$ is minimal, i.e. is the same as for $G=S_N$.

\item $\dim Fix(u^{\otimes k})=B_k$, where $B_k$ is the $k$-th Bell number.

\item $\int_Gu_{i_1j_1}\ldots u_{i_kj_k}=\frac{(N-k)!}{N!}$, for any $i,j\in I_N^k$.

\item $\int_Gu_{i_1j_1}\ldots u_{i_kj_k}\neq0$, for any $i,j\in I_N^k$.

\item $u_{i_1j_1}\ldots u_{i_kj_k}\neq0$, for any $i,j\in I_N^k$.
\end{enumerate}
\end{theorem}

\begin{proof}
All this is well-known, the idea being as follows:

\medskip

$(1)\implies(2)$ This follows from the fact that $u^{\otimes k}$ comes by summing certain actions $G\curvearrowright I_N^r$ with $r\leq k$, and the transitivity at $k$ implies the transitivity at any $r\leq k$.

\medskip

$(2)\implies(3)$ This comes from the well-known fact that for the symmetric group $S_N$, the multiplicity $\#(1\in u^{\otimes k})$ equals the $k$-th Bell number $B_k$, for any $k\leq N$.

\medskip

$(3)\implies(4)$ We can use the fact that $P_{i_1\ldots i_k,j_1\ldots j_k}=\int_Gu_{i_1j_1}\ldots u_{i_kj_k}$ is the orthogonal projection onto $Fix(u^{\otimes k})$. Thus we can assume $G=S_N$, and here we have:
$$\int_{S_N}u_{i_1j_1}\ldots u_{i_kj_k}
=\int_{S_N}\chi\left(\sigma\Big|\sigma(j_1)=i_1,\ldots,\sigma(j_k)=i_k\right)
=\frac{(N-k)!}{N!}$$

$(4)\implies(5)$ This is trivial.

\medskip

$(5)\implies(6)$ This is trivial too.

\medskip

$(6)\implies(1)$ This is clear, because if $u_{i_1j_1}\ldots u_{i_kj_k}=\chi(\sigma|\sigma(j_1)=i_1,\ldots,\sigma(j_k)=i_k)$ is nonzero, we can find an element $\sigma\in G$ such that $\sigma(j_s)=i_s$, $\forall s$.
\end{proof}

In the quantum case now, each magic unitary matrix $u=(u_{ij})$ produces a corepresentation $u^{\otimes k}=(u_{i_1j_1}\ldots u_{i_kj_k})$, and so a coaction map, constructed as follows:
$$\Phi:(\mathbb C^N)^{\otimes k}\to (\mathbb C^N)^{\otimes k}\otimes C(G)$$
$$\Phi(e_{i_1\ldots i_k})=\sum_{j_1\ldots j_k}e_{j_1\ldots j_k}\otimes u_{j_1i_1}\ldots u_{j_ki_k}$$

The problem is that $span(I_N^k)$ is no longer invariant, due to the fact that the variables $u_{ij}$ no longer commute. We can only say that $span(J_N^k)$ is invariant, where:
$$J_N^k=\left\{(i_1,\ldots,i_k)\in\{1,\ldots,N\}^k\Big|i_1\neq i_2\neq\ldots\neq i_k\right\}$$

Indeed, by using the fact, coming from the magic condition on $u$, that  $a\neq c,b=d$ implies $u_{ab}u_{cd}=0$, we obtain that for $i\in J_N^k$ we have, as desired:
$$\Phi(e_{i_1\ldots i_k})=\sum_{j_1\neq j_2\neq\ldots\neq j_k}e_{j_1\ldots j_k}\otimes u_{j_1i_1}\ldots u_{j_ki_k}$$

We can study the transitivity properties of this coaction, as follows:

\begin{proposition}
Given a closed subgroup $G\subset S_N^+$, consider the associated coaction map $\Phi:span(J_N^k)\to span(J_N^k)\otimes C(G)$. The following are then equivalent:
\begin{enumerate}
\item The fixed point space $Fix(\Phi)=\left\{\xi|\Phi(\xi)=\xi\otimes1\right\}$ is $1$-dimensional.

\item We have $Fix(\Phi)=\mathbb C\eta$, with $\eta=\sum_{i\in J_N^k}e_{i_1\ldots i_k}$.

\item We have the formula $\sum_{i\in J_N^k}\int_Gu_{i_1i_1}\ldots u_{i_ki_k}=1$, for any multi-index $j$.
\end{enumerate}
If these conditions are satisfied, we say that the coaction $\Phi$ is transitive.
\end{proposition}

\begin{proof}
The equivalences are elementary, the idea being as follows:

\medskip

$(1)\iff(2)$. Here we just have to check that we have indeed $\Phi(\eta)=\eta\otimes1$, with $\eta$ being as in the statement. By definition of $\Phi$, we have:
$$\Phi(\eta)=\sum_{i_1\neq i_2\neq\ldots\neq i_k}\sum_{j_1\neq j_2\neq\ldots\neq j_k}e_{j_1\ldots j_k}\otimes u_{j_1i_1}\ldots u_{j_ki_k}$$

Let us compute the middle sum $S$. When summing over indices $i_1\neq i_2$ we obtain:
$$(1-u_{j_2i_1})u_{j_2i_2}\ldots u_{j_ki_k}=u_{j_2i_2}\ldots u_{j_ki_k}$$

Then when summing over indices $i_2\neq i_3$ we obtain:
$$(1-u_{j_3i_2})u_{j_3i_3}\ldots u_{j_ki_k}=u_{j_3i_3}\ldots u_{j_ki_k}$$

And so on, up to obtaining in the end the following formula:
$$\sum_{i_k}u_{i_ki_k}=1$$

Thus we have $S=1$, and so the condition $\Phi(\eta)=\eta\otimes1$ is satisfied indeed.

\medskip

$(1)\iff(3)$ This comes from the following general formula, where $\chi$ is the character of the corepresentation associated to $\Phi$:
$$\dim Fix(\Phi)=\int_G\chi$$

Indeed, in the standard basis $\{e_i|i\in J_N^k\}$ we have:
$$\chi=\sum_{i\in J_N^k}u_{i_1i_1}\ldots u_{i_ki_k}$$

But this gives the result, by integrating. 
\end{proof}

We have the following partial analogue of Theorem 10.18:

\begin{proposition}
Given a closed subgroup $G\subset S_N^+$, with $N\geq4$, with matrix coordinates denoted $u_{ij}$, and a number $k\in\mathbb N$, the following conditions are equivalent:
\begin{enumerate}
\item The action $G\curvearrowright span(J_N^k)$ is transitive.

\item $Fix(u^{\otimes k})$ is minimal, i.e. is the same as for $G=S_N^+$.

\item $\dim Fix(u^{\otimes k})=C_k$, where $C_k$ is the $k$-th Catalan number.

\item $\int_Gu_{i_1j_1}\ldots u_{i_kj_k}$ is the same as for $G=S_N^+$, for any $i,j\in J_N^k$.
\end{enumerate}
\end{proposition}

\begin{proof}
This follows as in the first part of the proof of Theorem 10.18, by performing changes where needed, and by using the general theory from \cite{bc1}, as an input:

\medskip

$(1)\implies(2)$ This follows from the fact that $u^{\otimes k}$ comes by summing certain actions $G\curvearrowright J_N^r$ with $r\leq k$, and the transitivity at $k$ implies the transitivity at any $r\leq k$.

\medskip

$(2)\implies(3)$ This comes from the well-known fact that for the quantum group $S_N^+$ with $N\geq4$, the multiplicity $\#(1\in u^{\otimes k})$ equals the $k$-th Catalan number $C_k$.

\medskip

$(3)\implies(4)$ This comes from the well-known fact that $P_{i_1\ldots i_k,j_1\ldots j_k}=\int_Gu_{i_1j_1}\ldots u_{i_kj_k}$ is the orthogonal projection onto $Fix(u^{\otimes k})$, coming from \cite{wo1}.

\medskip

$(4)\implies(1)$ This follows by taking $i=j$ and then summing over this index, by using the transitivity criterion for $G\curvearrowright span(J_N^k)$ from Proposition 10.19 (3).
\end{proof}

Now let us compare our main results, Theorem 10.18 and Proposition 10.20. We conclude that the notion of $k$-transitivity for the subgroups $G\subset S_N$ extends to the quantum group case, $G\subset S_N^+$, depending on the value of $k$, as follows:

\medskip

(1) At $k=1,2$ everything extends well, due to the results in \cite{bi2}, \cite{lmr}.

\medskip

(2) At $k=3$ we have a good phenomenon, $P_3=NC_3$, and a bad one, $I_N^3\neq J_N^3$.

\medskip

(3) At $k\geq4$ we have two bad phenomena, namely $P_k\neq NC_k$ and $I_N^k\neq J_N^k$.

\medskip

Summarizing, our study suggests that things basically stop at $k=3$. So, as a conclusion, let us record the definition and main properties of the $3$-transitivity:

\index{higher transitivity}
\index{triple transitivity}

\begin{theorem}
A closed subgroup $G\subset S_N^+$ is $3$-transitive, in the sense that we have $\dim(Fix(u^{\otimes 3}))=5$, if and only if, for any $i,k,p$ distinct and any $j,l,q$ distinct:
$$\int_Gu_{ij}u_{kl}u_{pq}=\frac{1}{N(N-1)(N-2)}$$
In addition, in the classical case, we recover in this way the usual notion of $3$-transitivity.
\end{theorem}

\begin{proof}
We know from Proposition 10.20 that the 3-transitivity condition is equivalent to the fact that the integrals of type $\int_Gu_{ij}u_{kl}u_{pq}$ with $i\neq k\neq p$ and $j\neq l\neq q$ have the same values as those for $S_N^+$. But these values are computed by Proposition 10.17 and the Weingarten formula, and the 3-transitivity condition follows to be equivalent to:
$$\int_Gu_{ij}u_{kl}u_{pq}=
\begin{cases}
\frac{1}{N(N-1)(N-2)}&{\rm if}\ \ker(ikp)=\ker(jlq)=|||\\
\frac{1}{N(N-1)}&{\rm if}\ \ker(ikp)=\ker(jlq)=\sqcap\hskip-4.55mm{\ }_{{\ }^|}\\
0&{\rm if}\ \{\ker(ikp),\ker(jlq)\}=\{|||,\sqcap\hskip-4.55mm{\ }_{{\ }^|}\,\}
\end{cases}$$

Now observe that the last formula is automatic, by using the traciality of the integral and the magic assumption on $u$, and that the middle formula follows from the first one, by summing over $i,j$. Thus we have are left with the first formula, as stated. Finally, the last assertion follows from Theorem 10.18, applied at $k=3$.
\end{proof}

\section*{10c. Classification results}

Let us discuss now classification results for the transitive subgroups $G\subset S_N^+$. At $N=4$, we will need the following result from \cite{bb2}, that we know from chapter 4:

\index{twisting}

\begin{theorem}
We have an isomorphism of compact quantum groups
$$S_4^+=SO_3^{-1}$$
given by the Fourier transform over the Klein group $K=\mathbb Z_2\times\mathbb Z_2$.
\end{theorem}

\begin{proof}
This is something that we know from chapter 4. Consider indeed the following matrix, corresponding to the standard vector space action of $SO_3^{-1}$ on $\mathbb C^4$:
$$u^+=\begin{pmatrix}1&0\\ 0&u\end{pmatrix}$$

We apply to this matrix the Fourier transform over the Klein group $K=\mathbb Z_2\times\mathbb Z_2$: 
$$v=
\frac{1}{4}
\begin{pmatrix}
1&1&1&1\\
1&-1&-1&1\\
1&-1&1&-1\\
1&1&-1&-1
\end{pmatrix}
\begin{pmatrix}
1&0&0&0\\
0&u_{11}&u_{12}&u_{13}\\
0&u_{21}&u_{22}&u_{23}\\
0&u_{31}&u_{32}&u_{33}
\end{pmatrix}
\begin{pmatrix}
1&1&1&1\\
1&-1&-1&1\\
1&-1&1&-1\\
1&1&-1&-1
\end{pmatrix}$$

Then $v$ is magic, and a converse of this holds too, and this gives the result.
\end{proof}

We have the following classification result, also from \cite{bb2}:

\index{ADE classification}

\begin{theorem}
The closed subgroups of $S_4^+=SO_3^{-1}$ are as follows:
\begin{enumerate}
\item Infinite quantum groups: $S_4^+$, $O_2^{-1}$, $\widehat{D}_\infty$.

\item Finite groups: $S_4$, and its subgroups.

\item Finite group twists: $S_4^{-1}$, $A_5^{-1}$.

\item Series of twists: $D_{2n}^{-1}$ $(n\geq 3)$, $DC^{-1}_{2n}$ $(n\geq 2)$.

\item A group dual series: $\widehat{D}_n$, with $n\geq 3$.
\end{enumerate}
Moreover, these quantum groups are subject to an ADE classification result.
\end{theorem}

\begin{proof}
This is something quite technical. Regarding the precise statement, the idea is that, with the convention that prime stands for twists, which all unique in the cases below, and that double prime denotes pseudo-twists, the classification is as follows:

\medskip

(A) $\mathbb Z_1$, $\mathbb Z_2$, $\mathbb Z_3$, $K$, $\widehat{D}_n$ ($n=2,3,\ldots,\infty$), $S_4^+$.

\medskip

(D) $\mathbb Z_4$, $D_{2n}'$, $D_{2n}''$ ($n=2,3,\ldots$), $H_2^+$, $D_1$, $S_3$.

\medskip

(E) $A_4$, $S_4$, $S_4'$, $A_5'$.

\medskip

There are many comments to be made here, regarding our various conventions, and the construction of some of the above quantum groups, as follows:

\medskip

-- To start with, the 2-element group $\mathbb Z_2=\{1,\tau\}$ can act in 2 ways on 4 points: either with the transposition $\tau$ acting without fixed point, and we use here the notation $\mathbb Z_2$, or with $\tau$ acting with 2 fixed points, and we use here the notation $D_1$.

\medskip

-- Similarly, the Klein group $K=\mathbb Z_2\times\mathbb Z_2$ can act in 2 ways on 4 points: either with 2 non-trivial elements having 2 fixed points each, and we use here the notation $K$, or with all non-trivial elements having no fixed points, and we use here the notation $D_2=\widehat{D}_2$.

\medskip

-- We have $D_4'=D_4$, and $D_4''=G_0$, the Kac-Paljutkin quantum group. Besides being a pseudo-twist of $D_{2n}$, the quantum group $D_{2n}''$ with $n\geq2$ is known to be as well a pseudo-twist of the dicyclic, or binary cyclic group $DC_{2n}$.

\medskip

-- Finally, the definition of $D_{2n}'$, $D_{2n}''$ can be extended at $n=1,\infty$, and we formally have $D_2'=D_2''=K$, and $D_\infty'=D_\infty''=H_2^+$, but these conventions are not very useful. Also, as explained in \cite{bb2}, the groups $D_1,S_3$ are a bit special at (D).

\medskip

This was for the idea, and we refer to \cite{bb2} for the construction and properties of the various twists and pseudo-twists in the above ADE list, and for the classification.
\end{proof}

By restricting now the attention to the transitive case, we obtain:

\begin{theorem}
The small order transitive quantum groups are as follows:
\begin{enumerate}
\item At $N=1,2,3$ we have $\{1\}$, $\mathbb Z_2$, $\mathbb Z_3$, $S_3$.

\item At $N=4$ we have $\mathbb Z_2\times\mathbb Z_2,\mathbb Z_4,D_4,A_4,S_4,O_2^{-1},S_4^+$ and $S_4^{-1}$, $A_5^{-1}$.
\end{enumerate}
\end{theorem}

\begin{proof}
This follows from the above result, the idea being as follows:

\medskip

(1) This follows from the fact that we have $S_N=S_N^+$ at $N\leq3$, from \cite{wa2}.

\medskip

(2) This follows indeed from the above ADE classification of the subgroups $G\subset S_4^+$, from \cite{bb2}, with all the twists appearing in the statement being standard twists.
\end{proof}

As an interesting consequence of the above result, we have:

\begin{proposition}
The inclusion of compact quantum groups
$$S_4\subset S_4^+$$
is maximal, in the sense that there is no quantum group in between.
\end{proposition}

\begin{proof} 
This follows indeed from the above classification result.
\end{proof}

Let us study now the subgroups $G\subset S_5^+$. This is something substantially more complicated, which will require the use of some advanced results from subfactor theory. We first have the following elementary observations, regarding such subgroups:

\begin{proposition}
We have the following examples of subgroups $G\subset S_5^+$:
\begin{enumerate}
\item The classical subgroups, $G\subset S_5$. There are $16$ such subgroups, having order:
$$1,2,3,4,4,5,6,6,8,10,12,12,20,24,60,120$$

\item The group duals, $G=\widehat{\Gamma}\subset S_5^+$. These appear, via a Fourier transform construction, from the various quotients $\Gamma$ of the following groups:
$$\mathbb Z_4,\mathbb Z_2*\mathbb Z_2,\mathbb Z_2*\mathbb Z_3$$
\end{enumerate}
In addition, we have as well all the ADE quantum groups $G\subset S_4^+\subset S_5^+$ from Theorem 10.23, embedded via the $5$ standard embeddings $S_4^+\subset S_5^+$.
\end{proposition}

\begin{proof}
These results are well-known, the proof being as follows:

\medskip

(1) This is a classical result, with the groups which appear being respectively:

\medskip

-- The cyclic groups $\{1\},\mathbb Z_2,\mathbb Z_3,\mathbb Z_4$.

\medskip

-- The Klein group $K=\mathbb Z_2\times\mathbb Z_2$.

\medskip

-- The groups $\mathbb Z_5,\mathbb Z_6,S_3,D_4,D_5,A_4$.

\medskip

-- A copy of $S_3\rtimes\mathbb Z_2$.

\medskip

-- The general affine group $GA_1(5)=\mathbb Z_5\rtimes\mathbb Z_4$.

\medskip

-- And finally $S_4,A_5,S_5$.

\medskip

(2) This follows from Bichon's result in \cite{bi3}, stating that the group dual subgroups $G=\widehat{\Gamma}\subset S_N^+$ appear from the various quotients $\mathbb Z_{N_1}*\ldots*\mathbb Z_{N_k}\to\Gamma$, with $N_1+\ldots+N_k=N$. At $N=5$ the partitions are $5=1+4,1+2+2,2+3$, and this gives the result.
\end{proof}

Summarizing, the classification of the subgroups $G\subset S_5^+$ is a particularly difficult task, the situation here being definitely more complicated than at $N=4$.  Let us restrict now the attention to the transitive subgroups. We first have the following result:

\begin{proposition}
We have the following examples of transitive subgroups $G\subset S_5^+$:
\begin{enumerate}
\item The classical transitive subgroups $G\subset S_5$. There are only $5$ such subgroups, namely $\mathbb Z_5,D_5,GA_1(5),A_5,S_5$.

\item The transitive group duals, $G=\widehat{\Gamma}\subset S_5^+$. There is only one example here, namely the dual of $\Gamma=\mathbb Z_5$, which is $\mathbb Z_5$, already appearing above.
\end{enumerate}
In addition, all the ADE quantum groups $G\subset S_4^+\subset S_5^+$ are not transitive.
\end{proposition}

\begin{proof}
This follows indeed by examining the lists in Proposition 10.26:

\medskip

(1) The result here is well-known, and elementary. Observe that $GA_1(5)=\mathbb Z_5\rtimes\mathbb Z_4$, which is by definition the general affine group of $\mathbb F_5$, is indeed transitive.

\medskip

(2) This follows from the results in \cite{bi3}, because with $\mathbb Z_{N_1}*\ldots*\mathbb Z_{N_k}\to\Gamma$ as in the proof of Proposition 10.26 (2), the orbit decomposition is precisely $N=N_1+\ldots+N_k$.

\medskip

(3) Finally, the last assertion is clear, because the embedding $S_4^+\subset S_5^+$ is obtained precisely by fixing a point. Thus $S_4^+$ and its subgroups are not transitive, as claimed.
\end{proof}

In order to prove the uniqueness result, we will use the recent progress in subfactor theory \cite{jms}, concerning the classification of the small index subfactors. For our purposes, the most convenient formulation of the result in \cite{jms} is:

\begin{theorem}
The principal graphs of the irreducible index $5$ subfactors are:
\begin{enumerate}
\item $A_\infty$, and a non-extremal perturbation of $A_\infty^{(1)}$.

\item The McKay graphs of $\mathbb Z_5,D_5,GA_1(5),A_5,S_5$.

\item The twists of the McKay graphs of $A_5,S_5$.
\end{enumerate}
\end{theorem}

\begin{proof}
This is a heavy result, and we refer to \cite{jms} for the whole story. The above formulation is the one from \cite{jms}, with the subgroup subfactors there replaced by fixed point subfactors \cite{ba2}, and with the cyclic groups denoted as usual by $\mathbb Z_N$.
\end{proof}

In the quantum permutation group setting, this result becomes:

\begin{theorem}
The set of principal graphs of the transitive subgroups $G\subset S_5^+$ coincide with the set of principal graphs of the following subgroups:
$$\mathbb Z_5,D_5,GA_1(5),A_5,S_5,S_5^+$$
\end{theorem}

\begin{proof}
We must take the list of graphs in Theorem 10.28, and exclude some of the graphs, on the grounds that the graph cannot be realized by a transitive subgroup $G\subset S_5^+$. We have 3 cases here to be studied, as follows:

\medskip

(1) The graph $A_\infty$ corresponds to $S_5^+$ itself. As for the perturbation of $A_\infty^{(1)}$, this dissapears, because our notion of transitivity requires the subfactor extremality.

\medskip

(2) For the McKay graphs of $\mathbb Z_5,D_5,GA_1(5),A_5,S_5$ there is nothing to be done, all these graphs being solutions to our problem.

\medskip

(3)  The possible twists of $A_5,S_5$, coming from the graphs in Theorem 10.28 (3), cannot contain $S_5$, because their cardinalities are smaller or equal than $|S_5|=120$.
\end{proof}

In connection now with our maximality questions, we have:

\begin{theorem}
The inclusion $S_5\subset S_5^+$ is maximal.
\end{theorem}

\begin{proof}
This follows indeed from Theorem 10.29, with the remark that $S_5$ being transitive, so must be any intermediate subgroup $S_5\subset G\subset S_5^+$.
\end{proof}

As a comment here, with a little more work, the above considerations give the full list of  transitive subgroups $G\subset S_5^+$. To be more precise, we have here the various subgroups appearing in Theorem 10.29, plus some possible twists of $A_5,S_5$, which remain to be investigated. Finally, let us also mention that the classification of arbitrary subfactors at $N=6$ is not known, but some results are available under some strong transitivity assumptions, in the subfactor sense. We refer here to \cite{jms} and related papers.

\section*{10d. Maximality questions}

In general, the maximality of $S_N\subset S_N^+$ is a difficult question. We discuss here the ``standard approach'' to the maximality conjecture, via representation theory and diagrams. We first have the following result, coming from Tannakian duality:

\begin{proposition}
Consider a quantum group $S_N\subset G\subset S_N^+$, with fundamental corepresentation denoted $v$. We have then inclusions as follows, for any $k\in\mathbb N$, 
$$span\left(\xi_\pi\Big|\pi\in P(k)\right)\supset Fix(v^{\otimes k})\supset span\left(\xi_\pi\Big|\pi\in NC(k)\right)$$
and equality on the left or on the right, for any $k\in\mathbb N$, is equivalent to having equality on the left or on the right in the inclusions $S_N\subset G\subset S_N^+$.
\end{proposition}

\begin{proof}
Consider a quantum group $S_N\subset G\subset S_N^+$, and let $w,v,u$ be the fundamental corepresentations of these quantum groups. We have then inclusions as follows:
$$Fix(w^{\otimes k})\supset Fix(v^{\otimes k})\supset Fix(u^{\otimes k})$$

Moreover, by Peter-Weyl, equality on the left or on the right, for any $k\in\mathbb N$, is equivalent to having equality on the left or on the right in the inclusions $S_N\subset G\subset S_N^+$. Now by using the easiness property of $S_N,S_N^+$, this gives the result.
\end{proof}

The above result is good news, because what we have there is a purely combinatorial reformulation of the maximality conjecture, in terms of partitions, noncrossing partitions, and the associated vectors. To be more precise, we have the following statement:

\begin{theorem}
The following conditions are equivalent:
\begin{enumerate}
\item There is no intermediate quantum group, as follows:
$$S_N\subset G\subset S_N^+$$

\item Any linear combination of vectors of type
$$\xi\in span\left(\xi_\pi\Big|\pi\in P(k)\right)-span\left(\xi_\pi\Big|\pi\in NC(k)\right)$$
produces via Tannakian operations the flip map, $\Sigma(a\otimes b)=b\otimes a$.
\end{enumerate}
\end{theorem}

\begin{proof}
According to Proposition 10.31, the non-existence of the quantum groups $S_N\subset G\subset S_N^+$ is equivalent to the non-existence of Tannakian categories as follows:
$$span\left(\xi_\pi\Big|\pi\in P(k)\right)\supset C_k\supset span\left(\xi_\pi\Big|\pi\in NC(k)\right)$$

But this means that whenever we pick an element $\xi$ which is on the left, but not on the right, the Tannakian category that it generates should be the one on the left:
$$<\xi>=span\left(\xi_\pi\Big|\pi\in P(k)\right)$$

Now since the category of all partitions $P=(P(k))$ is generated by the basic crossing $\slash\hskip-2.1mm\backslash$, this amounts in saying that the Tannakian category generated by $\xi$ should contain the vector associated to this basic crossing, which is $\xi_{\slash\hskip-1.5mm\backslash}=\Sigma$, as desired.
\end{proof}

The above result might look quite encouraging, and the first thought goes into inventing some kind of tricky ``averaging operation'', perhaps probability-inspired, made up of Tannakian operations, which in practice means made of basic planar operations, which converts the crossing partitions $\pi\in P(k)-NC(k)$ into the basic crossing $\slash\hskip-2.1mm\backslash$. However, this is something difficult, and in fact such questions are almost always difficult.

\bigskip

Of course, we are not saying here that such things are hopeless, but rather that they require considerable work. In connection with the above-mentioned mysterious ``averaging operation'', our feeling is that this cannot be found with bare hands, and that a heavy use of a computer, in order to understand what is going on, is required. To our knowledge, no one has ever invested much time in all this, and so things here remain open. Getting back to Earth now, here are some concrete results, obtained in this way:

\begin{theorem}
The following happen:
\begin{enumerate}
\item There is no intermediate easy quantum group $S_N\subset G\subset S_N^+$.

\item A generalization of this fact holds, at easiness level $2$, instead of $1$.
\end{enumerate}
\end{theorem}

\begin{proof}
The idea here is that everything follows from Theorem 10.32, with suitable definitions for the various easiness notions involved, and by doing some combinatorics:

\medskip

(1) Here what happens is that any $\pi\in P-NC$ has the following property:
$$<\pi>=P$$

Indeed, the idea is to cap $\pi$ with semicircles, as to preserve one crossing, chosen in advance, and to end up, by a recurrence procedure, with the standard crossing:
$$\slash\hskip-2.1mm\backslash\in<\pi>$$

Now in terms of the notions in Theorem 10.32, the conclusion is that the criterion (2) there holds for the linear combinations $\xi$ having lenght 1, and this gives the result. Indeed, according to \cite{bsp}, the easy quantum groups are by definition those having Tannakian categories as follows, with $D=(D(k))$ being a certain category of partitions:
$$Fix(v^{\otimes k})=span\left(\xi_\pi\Big|\pi\in D(k)\right)$$

Thus, the generation formula $<\pi>=P$ established above does the job, and proves that an intermediate easy quantum group $S_N\subset G\subset S_N^+$ cannot exist. See \cite{ez1}.

\medskip

(2) This is a generalization of (1), the idea being that of looking at the combinations having length 2, of type $\xi=\alpha\xi_\pi+\beta\xi_\sigma$. Our first claim is that, assuming that $G\subset H$ comes from an inclusion of categories $D\subset E$, the maximality at order $2$ is equivalent to the following condition, for any $\pi,\sigma\in E$, not both in $D$, and for any $\alpha,\beta\neq0$:
$$<span(D),\alpha T_\pi+\beta T_\sigma>=span(E)$$

Consider indeed a category $span(D)\subset C\subset span(E)$, corresponding to a quantum group $G\subset K\subset H$ having order 2. The order 2 condition means that we have $C=<C\cap span_2(P)>$, where $span_2$ denotes the space of linear combinations having 2 components. Since we have $span(E)\cap span_2(P)=span_2(E)$, the order 2 formula reads:
$$C=<C\cap span_2(E)>$$

Now observe that the category on the right is generated by the categories $C_{\pi\sigma}^{\alpha\beta}$ constructed in the statement. Thus, the order 2 condition reads:
$$C=\left<C_{\pi\sigma}^{\alpha\beta}\Big|\pi,\sigma\in E,\alpha,\beta\in\mathbb C\right>$$

Now since the maximality at order 2 of the inclusion $G\subset H$ means that we have $C\in\{span(D),span(E)\}$, for any such $C$, we are led to the following condition:
$$C_{\pi\sigma}^{\alpha\beta}\in\{span(D),span(E)\}\quad,\quad\forall\pi,\sigma\in E,\alpha,\beta\in\mathbb C$$

Thus, we have proved our claim. In order to prove now that $S_N\subset S_N^+$ is maximal at order $2$, we can use semicircle capping. The statement that we have to prove is as follows: ``for $\pi\in P-NC,\sigma\in P$ and $\alpha,\beta\neq0$ we have $<\alpha T_\pi+\beta T_\sigma>=span(P)$''. 

\medskip

In order to do this, our claim is that the same method as at level 1 applies, after some suitable modifications. We have indeed two cases, as follows:

\medskip

-- Assuming that $\pi,\sigma$ have at least one different crossing, we can cap the partition $\pi$ as to end up with the basic crossing, and $\sigma$ becomes in this way an element of $P(2,2)$ different from this basic crossing, and so a noncrossing partition, from $NC(2,2)$. Now by substracting this noncrossing partition, which belongs to $C_{S_N^+}=span(NC)$, we obtain that the standard crossing belongs to $<\alpha T_\pi+\beta T_\sigma>$, and we are done.

\medskip

-- In the case where $\pi,\sigma$ have exactly the same crossings, we can start our descent procedure by selecting one common crossing, and then two strings of $\pi,\sigma$ which are different, and then joining the crossing to these two strings. We obtain in this way a certain linear combination $\alpha' T_{\pi'}+\beta'T_{\sigma'}\in <\alpha T_\pi+\beta T_\sigma>$ which satisfies the conditions in the first case discussed above, and we can continuate as indicated there.
\end{proof}

As a last topic now, we would like to discuss a ``continuous version'' of the maximality conjecture. The point indeed is that the quantum groups $S_N^+$, at square values of the parameter, $N=n^2$, are related to the quantum groups $O_n^+$, via a twisting procedure.

\bigskip

In practice, the phenomenon is best understood at the Tannakian level, where the Brauer theorem for $O_n^+$ is as follows, with $NC_2$ being the noncrossing pairings:
$$Fix(v^{\otimes k})=span\left(\xi_\pi\Big|\pi\in NC_2(k)\right)$$

Observe the similarity with the Brauer result for $S_N^+$, which holds as above, but with $NC_2(k)$ replaced by $NC(k)$. The point now is that we have a bijection as follows, obtained in one sense by fattening the partitions, and in the other sense, by shrinking them:
$$NC(k)\simeq NC_2(2k)$$

Thus, there is a link between $S_N^+$ and $O_n^+$, and the result here, from \cite{bbs}, that we know from chapter 4, states that $PO_n^+$ appears as a twist of $S_N^+$, with $N=n^2$.

\bigskip

All this suggests that the maximality conjectures for $S_N^+$ should have a continuous counterpart, and this is indeed the case, the results being as follows:

\begin{theorem}
The following happen:
\begin{enumerate}
\item The only intermediate easy quantum group $O_n\subset G\subset O_n^+$ is the half-classical orthogonal group $O_n^*$, obtained via the relations $abc=cba$.

\item There is a generalization of this result, stating that the uniqueness of $O_n^*$ as intermediate object appears at easiness level $2$.

\item Regarding the inclusion $O_n\subset O_n^*$, this is maximal, in the sense that there is no compact quantum group in between, of any kind.
\end{enumerate}
\end{theorem}

\begin{proof}
This is a collection of trivial and non-trivial results, as follows:

\medskip

(1) We must compute here the categories of pairings $NC_2\subset D\subset P_2$, and this can be done via some standard combinatorics, in three steps, as follows:

\medskip

(a) Let $\pi\in P_2-NC_2$, having $s\geq 4$ strings. Our claim is that:

\smallskip

-- If $\pi\in P_2-P_2^*$, there exists a semicircle capping $\pi'\in P_2-P_2^*$.

\smallskip

-- If $\pi\in P_2^*-NC_2$, there exists a semicircle capping $\pi'\in P_2^*-NC_2$.

\smallskip

Indeed, both these assertions can be easily proved, by drawing pictures.

\medskip

(b) Consider now a partition $\pi\in P_2(k,l)-NC_2(k,l)$. Our claim is that:

\smallskip

-- If $\pi\in P_2(k, l)-P_2^*(k,l)$ then $<\pi>=P_2$.

\smallskip

-- If $\pi\in P_2^*(k,l)-NC_2(k,l)$ then $<\pi>=P_2^*$.

\smallskip

This can be indeed proved by recurrence on the number of strings, $s=(k+l)/2$, by using (a), which provides us with a descent procedure $s\to s-1$, at any $s\geq4$.

\medskip

(c) Finally, assume that we are given an easy quantum group $O_N\subset G\subset O_N^+$, coming from certain sets of pairings $D(k,l)\subset P_2(k,l)$. We have three cases:

\smallskip

-- If $D\not\subset P_2^*$, we obtain $G=O_N$.

\smallskip

-- If $D\subset P_2,D\not\subset NC_2$, we obtain $G=O_N^*$.

\smallskip

-- If $D\subset NC_2$, we obtain $G=O_N^+$.

\smallskip

Thus, we are led to the conclusion in the statement.

\medskip

(2) Here the study is quite similar to the study that we did for the inclusion $S_N\subset S_N^+$, in the proof of Theorem 10.33 (2), by using standard semicircle capping.

\medskip

(3) This is something non-trivial, from \cite{bbs}, that can be done as follows:

\medskip

-- Our first claim is that the inclusion $\mathbb TO_N\subset U_N$ is maximal in the category of compact groups. But this is something standard, using some Lie theory.

\medskip

-- Our second claim is that the inclusion $PO_N\subset PU_N$ is maximal in the category of compact groups. Indeed, assuming that $PO_N\subset G \subset PU_N$ is a proper intermediate subgroup, then its preimage under the quotient map $U_N\to PU_N$ would be a proper intermediate subgroup of $\mathbb TO_N\subset U_N$, which is a contradiction.

\medskip

-- Finally, our claim is that the inclusion $O_N\subset O_n^*$ is maximal in the category of compact compact groups. But this follows from the maximality of the inclusion of projective versions $PO_N\subset PU_N$, by lifting. For further details on this, we refer to \cite{bbs}.
\end{proof}

Summarizing, there are some difficult questions going on here. Moving forward, still in relation with the maximality of the inclusion $S_N\subset S_N^+$, let us formulate:

\begin{question}
What can be said, constructively, about the ``exotic'' quantum permutation groups, $S_N\subset G\subset S_N^+$?
\end{question}

In answer, given such a quantum group $S_N\subset G\subset S_N^+$, we can expect for instance all the combinatorics of $G$ to be invariant under $S_N$, due to the embedding $S_N\subset G$. In order to establish some concrete results here, we can use the following formula, from \cite{mc1}:

\index{evaluation map}
\index{convolution}

\begin{proposition}
Assuming $S_N\subset G\subset S_N^+$, consider the quotient map 
$$\pi:C(G)\to C(S_N)$$
and set $ev_\sigma(a)=\pi(a)(\sigma)$, for any $a\in C(G)$. Then with 
$$\varphi_{\sigma\tau}=ev_{\sigma^{-1}}*\varphi*ev_\tau$$
we have the following formula,
$$\varphi_{\sigma\tau}(u_{i_1j_1}\ldots u_{i_pj_p})=\varphi(u_{\sigma(i_1)\tau(j_1)}\ldots u_{\sigma(i_p)\tau(j_p)})$$
valid for any $p\in\mathbb N$, and any indices $i_1,\ldots,i_p$ and $j_1,\ldots,j_p$.
\end{proposition}

\begin{proof}
We have indeed the following computation:
\begin{eqnarray*}
&&\varphi_{\sigma\tau}(u_{i_1j_1}\ldots u_{i_pj_p})\\
&=&(ev_{\sigma^{-1}}\otimes\varphi\otimes ev_\tau)\Delta^{(2)}(u_{i_1j_1}\ldots u_{i_pj_p})\\
&=&\sum_{k_1\ldots k_p}\sum_{l_1\ldots l_p}ev_{\sigma^{-1}}(u_{i_1k_1}\ldots u_{i_pk_p})\varphi(u_{k_1l_1}\ldots u_{k_pl_p})ev_\tau(u_{l_1j_1}\ldots u_{l_pj_p})\\
&=&\sum_{k_1\ldots k_p}\sum_{l_1\ldots l_p}\delta_{\sigma(i_1)k_1}\ldots\delta_{\sigma(i_p)k_p}
\varphi(u_{k_1l_1}\ldots u_{k_pl_p})\delta_{\tau(j_1)l_1}\ldots\delta_{\tau(j_p)l_p}\\
&=&\varphi(u_{\sigma(i_1)\tau(j_1)}\ldots u_{\sigma(i_p)\tau(j_p)})
\end{eqnarray*}

Thus, we obtain the formula in the statement. See \cite{mc1}.
\end{proof}

As a first application, we can apply this to the existence problem for the  algebraic $k$-orbitals. We obtain somehow ``half'' of the proof of that, as follows:

\begin{proposition}
Assuming $S_N\subset G\subset S_N^+$, the relation
$$(i_1,\ldots i_k)\sim(j_1,\ldots,j_k)\iff u_{i_1j_1}\ldots u_{i_kj_k}\neq0$$
depends only on $\ker i,\ker j\in P(k)$, and not of the particular multi-indices $i,j$.
\end{proposition}

\begin{proof}
By using Proposition 10.36 and Hahn-Banach we obtain, for any two permutations $\sigma,\tau\in S_N$, the following implication:
$$u_{i_1j_1}\ldots u_{i_kj_k}\neq0
\implies u_{\sigma(i_1)\tau(j_1)}\ldots u_{\sigma(i_k)\tau(j_k)}\neq0$$

Thus, we are led to the conclusion in the statement.
\end{proof}

In short, the conjectural algebraic $k$-orbitals are invariant under $S_N$, and the problem now is if the above relation $\sim$ is transitive on the partitions of $P(k)$. 

\bigskip

As another application, we have the following result:

\begin{proposition}
Given a quantum permutation group $S_N\subset G\subset S_N^+$, we have 
$$[u_{ij},u_{kl}]\neq0$$
for any indices $i\neq k$ and $j\neq l$. More generally, the following algebra
$$<u_{ij},u_{kl}>\subset C(G)$$
does not depend on the choice of $i\neq k$ and $j\neq l$. Even more generally, the algebra
$$<u_{i_1j_1},\ldots,u_{i_pj_p}>\subset C(G)$$
depends only on the relative position of the indices $(i_r,j_r)$ inside the standard square $\{1,\ldots,N\}^2$, and not on the precise value of these indices.
\end{proposition}

\begin{proof}
This is similar to the proof of Proposition 10.37, by using Proposition 10.36 applied to commutators, or more general quantities, and then using Hahn-Banach.
\end{proof}

As an interesting consequence of Proposition 10.38, we have:

\index{exotic quantum permutation}
\index{diagonal algebra}
\index{configuration space}

\begin{proposition}
Given an exotic quantum permutation group $S_N\subset G\subset S_N^+$, the corresponding diagonal algebra, which is by definition given by
$$D(G)=<u_{11},\ldots,u_{NN}>\subset C(G)$$
when regarded as a quotient of $C^*(\mathbb Z_2^{*N})$ is invariant under the action of $S_N$.
\end{proposition}

\begin{proof}
According to Proposition 10.38, we have an isomorphism as follows, between subalgebras of $D(G)$, for any indices $i\neq j$ and any $k\neq l$:
$$<u_{ii},u_{jj}>\simeq<u_{kk},u_{ll}>$$

More generally, still according to Proposition 10.38, we have such isomorphisms in higher length $3,4,5,\ldots$\,, and at length $N$ we obtain the result.
\end{proof}

The above result has some potentially interesting consequences. Indeed, since we have $S_N\subset G\subset S_N^+$, at the level of diagonal algebras we have quotient maps as follows:
$$D(S_N^+)\to D(G)\to D(S_N)$$

Thus, Proposition 10.39 suggests that $D(G)$ should come from a certain uniform intermediate quotient $\mathbb Z_2^{*N}\to\Gamma\to\mathbb Z_2^N$. This would be very interesting, first in order to understand the structure of $C(G)$ itself, which is not that much bigger than $D(G)$, and also for making the link with the quantum reflection groups.

\bigskip

We discuss now a number of freeness questions regarding the algebra $C(S_N^+)$, or rather its dense subalgebra $\mathcal C(S_N^+)\subset C(S_N^+)$ generated by the entries of $u=(u_{ij})$. The problem is that of understanding which are the polynomial relations relating the variables $u_{ij}$:
$$P(\{u_{ij}\})=0$$

That is, we would like to understand what is the kernel of the following map:
$$\mathbb C<\{X_{ij}\}>\to\mathcal C(S_N^+)\quad,\quad 
X_{ij}\to u_{ij}$$

We know that this kernel is the ideal generated by the magic relations, and the conjecture would be that at $N\geq4$ the relations $P(\{u_{ij}\})=0$ could only appear for ``trivial reasons''. However, it seems difficult to formulate a precise conjecture in this sense. Let us look instead at the monomial case. Here we would simply like to understand which are the monomials in the variables $u_{ij}$ which vanish, and we have:

\begin{conjecture}
For the quantum group $S_N^+$ with $N\geq4$, an equality of type 
$$u_{i_1j_1}\ldots u_{i_kj_k}=0$$
can only appear ``for trivial reasons'', due to an occurrence inside this relation of a cancellation formula of type $pq=0$, with $p,q\in\{u_{ij}\}$. 
\end{conjecture}

Here is an equivalent form of Conjecture 10.40, using the orbital formalism:

\begin{conjecture}
The relation $i\sim j\iff u_{i_1j_1}\ldots u_{i_kj_k}\neq 0$ is given by: 
$$(i_1,\ldots,i_k)\sim(j_1,\ldots,j_k)\iff
\begin{cases}
i_1=i_2\iff j_1=j_2&\\
i_2=i_3\iff j_2=j_3&\\
\ldots&\\
i_{k-1}=i_k\iff j_{k-1}=j_k
\end{cases}$$
In particular $\sim$ is transitive, and we have $2^{k-1}$ equivalence classes, or orbitals.
\end{conjecture}

Here the first assertion is a reformulation of Conjecture 10.40, using the definition of the equivalence relation $\sim$, and the transitivity assumption and the orbital count are trivial consequences of it. As before, the cases $k=3,4$ are of particular interest, but this time not in view of their potential doability, but rather in view of their potential applications. In relation with this latter conjecture, we have the following useful result:

\begin{proposition}
Given a quantum permutation group $G\subset S_N^+$, the relation
$$(i_1,\ldots,i_k)\sim(j_1,\ldots,j_k)\iff u_{i_1j_1}\ldots u_{i_kj_k}\neq 0$$
is in fact a relation between the $k$-orbitals of $G_{class}\subset S_N$. In particular, with $G=S_N^+$, what we have is a relation on the set $P(k)$ of partitions of $\{1,\ldots,k\}$.
\end{proposition}

\begin{proof}
This is explained above, the idea being that any linear form $\varphi\in C(G)^*$ can be suitably modified by permutations $\sigma,\tau\in G_{class}$ as to take the same values on the quantities $u_{i_1j_1}\ldots u_{i_kj_k}$, under the action of $G_{class}$ on the indices. As explained there, together with Hahn-Banach, this gives the result. Finally, in the case of $G=S_N^+$ we have $G_{class}=S_N$, and the orbits here are $\{1,\ldots,N\}^k/\sim\ =P(k)$, as claimed.
\end{proof}

Finally, getting back now to the general freeness questions formulated in the beginning, an interesting question, which is complementary to the above ``orbital'' ones, regards the diagonal coefficients $u_{ii}$. We have here the following conjecture:

\index{higher transitivity}
\index{freeness conjecture}

\begin{conjecture}
For the quantum group $S_N^+$ with $N\geq4$, the variables 
$$u_{ii}\in C(S_N^+)$$
are algebrically free. In particular, the diagonal algebra
$$D(S_N^+)=<u_{11},\ldots,u_{NN}>\subset C(S_N^+)$$
is isomorphic to the group algebra $C^*(\mathbb Z_2^{*N})$.
\end{conjecture}

Observe that this is indeed complementary to the above orbital questions, because we are dealing here with linear products of type $u_{i_1i_1}\ldots u_{i_ki_k}$, and their linear combinations, which are ``trivial'' from the point of view of the orbital theory.

\bigskip

As before, it is very unclear on how to prove such things. In relation now with the ``exotic'' case, $S_N\subset G\subset S_N^+$, we would like such things to be understood, because the orbitals for $G$, if well-defined indeed, should lie between those for $S_N$ and for $S_N^+$, and the same should happen for the diagonal algebra. Thus, the exotic quantum permutation groups  $S_N\subset G\subset S_N^+$ would lead to ``exotic'' relations on $P(k)$, and also to ``exotic'' quotients of $\mathbb Z_2^{*N}$, and with a bit of luck, such exotic things could not really exist.

\bigskip

Finally, let us mention that Conjecture 10.40 was recently proved by McCarthy \cite{mc2}. For more on all this, and related topics, we refer to his papers \cite{mc1}, \cite{mc2}.

\section*{10e. Exercises} 

Things have been fairly technical in this chapter, and as a unique exercise on all this, in connection with the questions at the end, we have:

\begin{exercise}
Prove that the diagonal algebra of the quantum group $S_N^+$,
$$D(S_N^+)=<u_{11},\ldots,u_{NN}>\subset C(S_N^+)$$
is isomorphic to the group algebra $C^*(\mathbb Z_2^{*N})$.
\end{exercise}

This is something that we already discussed before, in relation with the quantum partial permutations and their associated diagonal algebras, and with the related functional analysis problematics of constructing a quantum group of type $S_\infty^+$. In view of the above considerations, making the link with the exotic quantum groups $S_N\subset G\subset S_N^+$, this question appears to be of great importance, related to many things.

\chapter{Liberation theory}

\section*{11a. Reflection groups}

In this chapter we discuss the liberation theory for the closed subgroups $G\subset S_N$, with the aim of finding their free analogues $G^+\subset S_N^+$. We will have in mind as well the twisted case, involving groups $G\subset S_Z$ and quantum groups $G^+\subset S_Z^+$, with the study here for $Z=M_N$, relying on the isomomorphism $S_{M_N}^+=PO_N^+=PU_N^+$, leading us, up to taking projective versions, into the liberation question for the closed subgroups $G\subset U_N$.

\bigskip

Long story short, we are interested in the general liberation question for the closed subgroups $G\subset U_N$, with particular attention to the permutation group case, $G\subset S_N$. There are many things that can be said here, and we have seen quite a number of non-trivial results on this subject, so far in this book. Our goal will be that of further building on this material, with some advanced results on the subject. The presentation will be quite technical, for the most mixing difficult theorems and conjectural statements.

\bigskip

In order to get started, the considerations from the previous chapter show that for $G=S_N$ we are into a no-go problem, the precise conjecture being as follows:

\begin{conjecture}
There is no proper liberation $S_N\subset G\subset S_N^+$.
\end{conjecture}

We refer to the previous chapter for more on this conjecture, including evidence, and difficulties involved. We will be actually back to this in the next chapter too, with some probabilistic considerations, in relation with De Finetti theorems and other, relating Conjecture 11.1 to the well-known, folklore saying in probability theory, that ``there is nothing much interesting between classical independence, and Voiculescu freeness''.

\bigskip

Another topic to be discussed, before anything more advanced, is the notion of quantum group freeness. We already have some knowledge of that, from easiness, the idea being that the easy groups $S_N\subset G\subset U_N$ are in correspondence, although not always one-to-one, with the easy quantum groups $S_N^+\subset G^+\subset U_N^+$, via the following operations, connecting the associated categories of partitions, $D$ and $D^+$:
$$D^+=D\cap NC\quad,\quad D=<D^+,\slash\hskip-2.05mm\backslash>$$

But this shows in particular that the easy quantum groups $S_N\subset G\subset U_N^+$ which are ``free'', in whatever natural sense of freeness that you would prefer, are simply those satisfying $S_N^+\subset G$. So, forgetting now about easiness, we are led into the question of studying the closed subgroups $G\subset U_N^+$ satisfying $S_N^+\subset G$, and here we have:

\begin{conjecture}
The free quantum groups, $S_N^+\subset G\subset U_N^+$, are all easy.
\end{conjecture}

To be more precise, the above easiness discussion suggests calling ``free'' the intermediate quantum groups $S_N^+\subset G\subset U_N^+$. However, in practice, since there are no other known examples besides the easy ones, we are led to the above conjecture.

\bigskip

As before with Conjecture 11.1, this should be regarded as being a difficult combinatorial question. Indeed, Tannakian duality allows us to reformulate everything in terms of partitions, and the resulting partition question is well-known to be difficult. Finally, let us mention that the easy quantum groups $S_N^+\subset G\subset U_N^+$ can be fully classified, as done by Tarrago-Weber in \cite{twa}, and to be explained later, and so Conjecture 11.2 would eventually provide us with the precise list of free quantum groups.

\bigskip

Summarizing, the very first thoughts about the liberation of subgroups $G\subset S_N$, or more generally subgroups $G\subset U_N$, lead us into two difficult conjectures. However, there is solution to anything, and on the same topic as Conjectures 11.1 and 11.2, we have:

\begin{question}
What are the liberations of the complex reflection groups $G\subset K_N$? Also, what about the half-classical quantum groups, $G\subset U_N^*$? 
\end{question}

Here the first question is obviously related to Conjecture 11.1, and the point is that, while this still remains something quite undoable, at least the degree of flexibility brought by looking at arbitrary reflection groups will allow us to do some interesting, non-trivial work on the subject. As for the second question, this is a sort of downgrade of Conjecture 11.2, the idea being that, in the lack of ideas so far on how to deal with the free quantum groups, we can try instead to look at the half-classical ones, and do some work there.

\bigskip

Let us first discuss the first question. We have seen in the chapters 5-8 that the basic reflection groups $H_N^s=\mathbb Z_s\wr S_N$, which are all easy, have free analogues $H_N^{s+}=\mathbb Z_s\wr_*S_N^+$, and that the theory of these quantum groups, both classical and free, is very interesting, algebrically and analytically speaking. However, the world of quantum reflection groups is in fact much wider than this. In the classical case already, the classification theorem for the complex reflection groups, a celebrated result by Shephard and Todd \cite{sto}, from the 1950s, is something non-trivial, which can be briefly stated as follows:

\index{complex reflection group}
\index{Shephard-Todd}

\begin{theorem}
The irreducible complex reflection groups are
$$H_N^{sd}=\left\{U\in H_N^s\Big|(\det U)^d=1\right\}$$
along with $34$ exceptional examples.
\end{theorem}

\begin{proof}
This is something quite advanced, and we refer here to the paper of Shephard and Todd \cite{sto}, and to the subsequent literature on the subject.
\end{proof}

In the general quantum case now, the axiomatization and classification of the quantum reflection groups is a key problem, which is not understood yet. However, at least from an intuitive point of view, we can say that such a quantum group should be by definition a liberation of a complex reflection group $G\subset U_N$. And with this notion in hand, getting now into the case where $G\subset U_N$ is in addition assumed to be irreducible, and leaving aside the above 34 exceptional examples, we are led into the question of understanding the liberations of the groups $H_N^{sd}$. And here, the conjecture is as follows:

\begin{conjecture}
Among the complex reflection groups $H_N^{sd}$, only the easy ones, namely $H_N^s$, admit liberations. In addition, these latter liberations should be easy.
\end{conjecture}

This might look a bit complicated, and here are some explanations. In what regards the first assertion, this comes from the widespread belief that the usual determinant has no ``liberated analogues'', and so the constructions involving it, including that of the groups $H_N^{sd}$ with $d<s$, have no quantum analogues. An obvious question here, to start with, dealing with the simplest possible case, regards the lack of liberations of the alternating group $A_N$, but this is not solved yet. However, we will be back to this later, with some concrete results on the subject, regarding the half-liberation operation.

\bigskip

In what regards now the second assertion in Conjecture 11.5, this is something quite subtle, related to both Conjecture 11.1 and Conjecture 11.2. Indeed, at $s=1$, where $H_N^s=S_N$, this would tell us that the liberations $S_N\subset G\subset S_N^+$ are automatically easy, solving Conjecture 11.2 in the case $G_{class}=S_N$, and solving as well Conjecture 11.1, via the elementary fact, explained in chapter 10, that there is no proper easy liberation $S_N\subset G\subset S_N^+$. Let us also mention that, by using Theorem 11.4, it is possible to further complicate Conjecture 11.5, as to fully cover both Conjecture 11.5 and Conjecture 11.6. But too many conjectures anyway, and so probably time to stop here.

\bigskip

To summarize now, our hunt for complex reflection groups, or even for a doable problem regarding them, has been so far unsuccessful. As a matter of not giving up, however, and having some bush meat on the grill for tonight, let us restrict the attention to the easy case. Here our questions become far more reasonable, basically reducing to:

\begin{question}
What are the easy quantum groups $H_N^s\subset G\subset H_N^{s+}$?
\end{question}

This latter question looks fully doable and reasonable, and we already know the answer to it in the case $s=1$, form chapter 10, the answer being that there is no intermediate easy group $S_N\subset G\subset S_N^+$, with proof coming via some simple combinatorics.

\bigskip

In order to upgrade this result that we have, at $s=1$, let us study first the other cases of main interest, namely $s=2,\infty$. Here the groups $H_N^s$ are respectively $H_N,K_N$, and we will be interested in classifying the intermediate easy quantum groups, as follows:
$$H_N\subset G\subset H_N^+\quad,\quad K_N\subset G\subset K_N^+$$

The problem, however, is that there are many such quantum groups, and we will have to develop our study slowly. Let us start with the following result, from \cite{rwe}:

\index{intermediate liberation}

\begin{theorem}
We have intermediate liberations $H_N^{[\infty]},K_N^{[\infty]}$ as follows, constructed by using the relations $\alpha\beta\gamma=0$, for any $a\neq c$ on the same row or column of $u$,
$$\xymatrix@R=15mm@C=17mm{
K_N\ar[r]&K_N^*\ar[r]&K_N^{[\infty]}\ar[r]&K_N^+\\
H_N\ar[r]\ar[u]&H_N^*\ar[r]\ar[u]&H_N^{[\infty]}\ar[r]\ar[u]&H_N^+\ar[u]}$$
with the convention $\alpha=a,a^*$, and so on. These quantum groups are easy, the corresponding categories $P_{even}^{[\infty]}\subset P_{even}$ and $\mathcal P_{even}^{[\infty]}\subset\mathcal P_{even}$ being generated by $\eta=\ker(^{iij}_{jii})$.
\end{theorem}

\begin{proof}
This is routine, by using the fact that the relations $\alpha\beta\gamma=0$ in the statement are equivalent to the following condition, with $|k|=3$:
$$\eta\in End(u^{\otimes k})$$ 

For further details on these quantum groups, we refer to \cite{ez1}, \cite{rwe}.
\end{proof}

Leaving aside for a moment classification questions, we would like to present now, as a somewhat unexpected application of the above constructions, the solution to some twisting problems, left open before. We will need the following technical result:

\begin{proposition}
We have the following equalities,
\begin{eqnarray*}
P_{even}^*&=&\left\{\pi\in P_{even}\Big|\varepsilon(\tau)=1,\forall\tau\leq\pi,|\tau|=2\right\}
\\
P_{even}^{[\infty]}&=&\left\{\pi\in P_{even}\Big|\sigma\in P_{even}^*,\forall\sigma\subset\pi\right\}\\
P_{even}^{[\infty]}&=&\left\{\pi\in P_{even}\Big|\varepsilon(\tau)=1,\forall\tau\leq\pi\right\}
\end{eqnarray*}
where $\varepsilon:P_{even}\to\{\pm1\}$ is the signature of even permutations.
\end{proposition}

\begin{proof}
This is routine combinatorics, the idea being as follows:

\medskip

(1) Given $\pi\in P_{even}$, we have $\tau\leq\pi,|\tau|=2$ precisely when $\tau=\pi^\beta$ is the partition obtained from $\pi$ by merging all the legs of a certain subpartition $\beta\subset\pi$, and by merging as well all the other blocks. Now observe that $\pi^\beta$ does not depend on $\pi$, but only on $\beta$, and that the number of switches required for making $\pi^\beta$ noncrossing is $c=N_\bullet-N_\circ$ modulo 2, where $N_\bullet/N_\circ$ is the number of black/white legs of $\beta$, when labelling the legs of $\pi$ counterclockwise $\circ\bullet\circ\bullet\ldots$ Thus $\varepsilon(\pi^\beta)=1$ holds precisely when $\beta\in\pi$ has the same number of black and white legs, and this gives the result.

\medskip

(2) This simply follows from the equality $P_{even}^{[\infty]}=<\eta>$ coming from Theorem 11.7, by computing $<\eta>$, and for the complete proof here we refer to Raum-Weber \cite{rwe}.

\medskip

(3) We use here the fact, also from \cite{rwe}, that the relations $g_ig_ig_j=g_jg_ig_i$ are trivially satisfied for real reflections. This leads to the following conclusion:
$$P_{even}^{[\infty]}(k,l)=\left\{\ker\begin{pmatrix}i_1&\ldots&i_k\\ j_1&\ldots&j_l\end{pmatrix}\Big|g_{i_1}\ldots g_{i_k}=g_{j_1}\ldots g_{j_l}\ {\rm inside}\ \mathbb Z_2^{*N}\right\}$$

In other words, the partitions in $P_{even}^{[\infty]}$ are those describing the relations between free variables, subject to the conditions $g_i^2=1$. We conclude that $P_{even}^{[\infty]}$ appears from $NC_{even}$ by ``inflating blocks'', in the sense that each $\pi\in P_{even}^{[\infty]}$ can be transformed into a partition $\pi'\in NC_{even}$ by deleting pairs of consecutive legs, belonging to the same block. Now since this inflation operation leaves invariant modulo 2 the number $c\in\mathbb N$ of switches in the definition of the signature, it leaves invariant the signature $\varepsilon=(-1)^c$ itself, and we obtain in this way the inclusion ``$\subset$'' in the statement. Conversely now, given $\pi\in P_{even}$ satisfying $\varepsilon(\tau)=1$, $\forall\tau\leq\pi$, our claim is that:
$$\rho\leq\sigma\subset\pi,|\rho|=2\implies\varepsilon(\rho)=1$$

Indeed, let us denote by $\alpha,\beta$ the two blocks of $\rho$, and by $\gamma$ the remaining blocks of $\pi$, merged altogether. We know that the partitions $\tau_1=(\alpha\wedge\gamma,\beta)$, $\tau_2=(\beta\wedge\gamma,\alpha)$, $\tau_3=(\alpha,\beta,\gamma)$ are all even. On the other hand, putting these partitions in noncrossing form requires respectively $s+t,s'+t,s+s'+t$ switches, where $t$ is the number of switches needed for putting $\rho=(\alpha,\beta)$ in noncrossing form. Thus $t$ is even, and we are done. With the above claim in hand, we conclude, by using the second equality in the statement, that we have $\sigma\in P_{even}^*$. Thus $\pi\in P_{even}^{[\infty]}$, which ends the proof of ``$\supset$''.
\end{proof}

With the above result in hand, we can now prove:

\index{twisting}

\begin{theorem}
The basic quantum reflection groups, namely 
$$\xymatrix@R=15mm@C=17mm{
K_N\ar[r]&K_N^*\ar[r]&K_N^{[\infty]}\ar[r]&K_N^+\\
H_N\ar[r]\ar[u]&H_N^*\ar[r]\ar[u]&H_N^{[\infty]}\ar[r]\ar[u]&H_N^+\ar[u]}$$
are equal to their own Schur-Weyl twists.
\end{theorem}

\begin{proof}
This basically comes from the results that we have, as follows:

\medskip

(1) In the real case, the verifications are as follows:

\medskip

-- $H_N^+$. We know from chapter 4 that for $\pi\in NC_{even}$ we have $\bar{T}_\pi=T_\pi$, and since we are in the situation $D\subset NC_{even}$, the definitions of $G,\bar{G}$ coincide.

\medskip

-- $H_N^{[\infty]}$. Here we can use the same argument as in (1), based this time on the description of $P_{even}^{[\infty]}$ involving the signature found in Proposition 11.8.

\medskip

-- $H_N^*$. We have $H_N^*=H_N^{[\infty]}\cap O_N^*$, so $\bar{H}_N^*\subset H_N^{[\infty]}$ is the subgroup obtained via the defining relations for $\bar{O}_N^*$. But all the $abc=-cba$ relations defining $\bar{H}_N^*$ are automatic, of type $0=0$, and it follows that $\bar{H}_N^*\subset H_N^{[\infty]}$ is the subgroup obtained via the relations $abc=cba$, for any $a,b,c\in\{u_{ij}\}$. Thus we have $\bar{H}_N^*=H_N^{[\infty]}\cap O_N^*=H_N^*$, as claimed.

\medskip

-- $H_N$. We have $H_N=H_N^*\cap O_N$, and by functoriality, $\bar{H}_N=\bar{H}_N^*\cap\bar{O}_N=H_N^*\cap\bar{O}_N$. But this latter intersection is easily seen to be equal to $H_N$, as claimed.

\medskip

(2) In the complex case the proof is similar, by using the same arguments.
\end{proof}

Getting back now to classification, and to Question 11.6, let us first discuss the classification of the easy quantum groups $H_N\subset G\subset H_N^+$. Following Raum-Weber \cite{rwe}, it is actually convenient to discuss, more generally, the classification of the easy quantum groups $H_N\subset G\subset O_N^+$. And here, we first have the folowing result:

\index{easy reflection group}

\begin{proposition}
The easy quantum groups $H_N\subset G\subset O_N^+$ are as follows,
$$\xymatrix@R=7mm@C=35mm{
H_N^+\ar[r]&O_N^+\\
H_N^{[\infty]}\ar@.[u]&O_N^*\ar[u]\\
H_N\ar@.[u]\ar[r]&O_N\ar[u]}$$
with the dotted arrows indicating that we have intermediate quantum groups there.
\end{proposition}

\begin{proof}
This is something quite technical, the idea being as follows:

\medskip

(1) We have a first dichotomy concerning the quantum groups in the statement, namely $H_N\subset G\subset O_N^+$, which must fall into one of the following two classes:
$$O_N\subset G\subset O_N^+\quad,\quad 
H_N\subset G\subset H_N^+$$

(2) Moreover, the early classification results from \cite{ez1} solve as well the first problem, namely $O_N\subset G\subset O_N^+$, with $G=O_N^*$ being the unique non-trivial solution.

\medskip

(3) We have then a second dichotomy, concerning the quantum groups which are left, namely $H_N\subset G\subset H_N^+$, which must fall into one of the following two classes:
$$H_N\subset G\subset H_N^{[\infty]}\quad,\quad 
H_N^{[\infty]}\subset G\subset H_N^+$$

All this comes indeed from various papers, and mainly from the final classification paper of Raum and Weber \cite{rwe}, where the quantum groups $S_N\subset G\subset H_N^+$ with $G\not\subset H_N^{[\infty]}$ were classified, and shown to contain $H_N^{[\infty]}$. For full details here, we refer to \cite{rwe}.
\end{proof}

Summarizing, we are left with classifying the following easy quantum groups:
$$H_N\subset G\subset H_N^{[\infty]}\quad,\quad 
H_N^{[\infty]}\subset G\subset H_N^+$$

Regarding the second case, namely $H_N^{[\infty]}\subset G\subset H_N^+$, the result here, also from \cite{rwe}, which is quite technical, but has a simple formulation, is as follows:

\begin{proposition}
Let $H_N^{[r]}\subset H_N^+$ be the easy quantum group coming from:
$$\pi_r=\ker\begin{pmatrix}1&\ldots&r&r&\ldots&1\\1&\ldots&r&r&\ldots&1\end{pmatrix}$$
We have then inclusions of quantum groups as follows,
$$H_N^+=H_N^{[1]}\supset H_N^{[2]}\supset H_N^{[3]}\supset\ldots\ldots\supset H_N^{[\infty]}$$
and we obtain all easy quantum groups $H_N^{[\infty]}\subset G\subset H_N^+$, satisfying $G\neq H_N^{[\infty]}$.
\end{proposition}

\begin{proof}
Once again, this is something technical, and we refer here to \cite{rwe}.
\end{proof}

It remains to discuss the easy quantum groups $H_N\subset G\subset H_N^{[\infty]}$, with the endpoints $G=H_N,H_N^{[\infty]}$ included. Once again, we follow here \cite{rwe}. First, we have:

\begin{definition}
A discrete group generated by real reflections, $g_i^2=1$,
$$\Gamma=<g_1,\ldots,g_N>$$
is called uniform if each $\sigma\in S_N$ produces a group automorphism, $g_i\to g_{\sigma(i)}$.
\end{definition}

Consider now a uniform reflection group, $\mathbb Z_2^{*N}\to\Gamma\to\mathbb Z_2^N$. We can associate to it a family of subsets $D(k,l)\subset P(k,l)$, which form a category of partitions, as follows:
$$D(k,l)=\left\{\pi\in P(k,l)\Big|\ker\binom{i}{j}\leq\pi\implies g_{i_1}\ldots g_{i_k}=g_{j_1}\ldots g_{j_l}\right\}$$

Observe that we have inclusions $P_{even}^{[\infty]}\subset D\subset P_{even}$. Conversely now, given such a category, we can associate to it a uniform reflection group $\mathbb Z_2^{*N}\to\Gamma\to\mathbb Z_2^N$, as follows:
$$\Gamma=\left\langle g_1,\ldots g_N\Big|g_{i_1}\ldots g_{i_k}=g_{j_1}\ldots g_{j_l},\forall i,j,k,l,\ker\binom{i}{j}\in D(k,l)\right\rangle$$

As explained in \cite{rwe}, the above correspondences $\Gamma\to D$ and $D\to\Gamma$ are bijective, and inverse to each other, at $N=\infty$. We have in fact the following result, from \cite{rwe}:

\begin{proposition}
We have correspondences between:
\begin{enumerate}
\item Uniform reflection groups $\mathbb Z_2^{*\infty}\to\Gamma\to\mathbb Z_2^\infty$.

\item Categories of partitions $P_{even}^{[\infty]}\subset D\subset P_{even}$.

\item Easy quantum groups $G=(G_N)$, with $H_N^{[\infty]}\supset G_N\supset H_N$.
\end{enumerate}
\end{proposition}

\begin{proof}
This is something quite technical, which follows along the lines of the above discussion. As an illustration, if we denote by $\mathbb Z_2^{\circ N}$ the quotient of $\mathbb Z_2^{*N}$ by the relations of type $abc=cba$ between the generators, we have the following correspondences:
$$\xymatrix@R=15mm@C=15mm{
\mathbb Z_2^N\ar@{~}[d]&\mathbb Z_2^{\circ N}\ar[l]\ar@{~}[d]&\mathbb Z_2^{*N}\ar[l]\ar@{~}[d]\\
H_N\ar[r]&H_N^*\ar[r]&H_N^{[\infty]}}$$

More generally, for any $s\in\{2,4,\ldots,\infty\}$, the quantum groups $H_N^{(s)}\subset H_N^{[s]}$ constructed in \cite{ez1} come from the quotients of $\mathbb Z_2^{\circ N}\leftarrow\mathbb Z_2^{*N}$ by the relations $(ab)^s=1$. See \cite{rwe}.
\end{proof}

We can now formulate a final classification result, as follows:

\begin{theorem}
The easy quantum groups $H_N\subset G\subset O_N^+$ are as follows,
$$\xymatrix@R=3mm@C=50mm{
H_N^+\ar[r]&O_N^+\\
H_N^{[r]}\ar[u]\\
H_N^{[\infty]}\ar[u]&O_N^*\ar[uu]\\
H_N^\Gamma\ar[u]\\
H_N\ar[u]\ar[r]&O_N\ar[uu]}$$
with the family $H_N^\Gamma$ covering $H_N,H_N^{[\infty]}$, and with the series $H_N^{[r]}$ covering $H_N^+$.
\end{theorem}

\begin{proof}
This follows from the above, and we refer here to \cite{rwe}. Let us also mention that the above diagram would look better flipped, but it is better to leave it like this, because this diagram as drawn is a face of the ``standard cube''. More on this later.
\end{proof}

\section*{11b. Soft liberation}

Getting back now to the non-easy case, we will be interested in what follows in the ``twistable'' case, where the theory is more advanced. Let us start with:

\index{homogeneous quantum group}
\index{twistable quantum group}
\index{half-homogeneous quantum group}

\begin{definition}
A closed subgroup $G\subset U_N^+$ is called:
\begin{enumerate}
\item Half-homogeneous, when it contains the alternating group, $A_N\subset G$.

\item Homogeneous, when it contains the symmetric group, $S_N\subset G$.

\item Twistable, when it contains the hyperoctahedral group, $H_N\subset G$.
\end{enumerate}
\end{definition}

These notions are mostly motivated by the easy case. Here we have by definition $S_N\subset G\subset U_N^+$, so our quantum group is automatically homogeneous. The point now is that the twistability assumption, as formulated above, corresponds to the following condition, at the level of the associated category of partitions $D\subset P$:
$$D\subset P_{even}$$

We recognize here the condition which is needed for performing the Schur-Weyl twisting operation, explained in chapter 4, and based on the signature map:
$$\varepsilon:P_{even}\to\{\pm1\}$$

As a conclusion, in the easy case our notion of twistability is the correct one. In general, there are of course more general twisting methods, usually requiring $\mathbb Z_2^N\subset G$ only. But in the half-homogeneous case, the condition $\mathbb Z_2^N\subset G$ is equivalent to $H_N\subset G$. With this discussion done, let us formulate now the following definition:

\index{twistable quantum reflection group}

\begin{definition}
A twistable quantum reflection group is an intermediate subgroup
$$H_N\subset K\subset K_N^+$$
between the group $H_N=\mathbb Z_2\wr S_N$, and the quantum group $K_N^+=\mathbb T\wr_*S_N^+$.
\end{definition}

Here is now another definition, which is important for general compact quantum group purposes, and which provides motivations for our formalism from Definition 11.15:

\index{reflection subgroup}

\begin{definition}
Given a closed subgroup $G\subset U_N^+$ which is twistable, in the sense that we have $H_N\subset G$, we define its associated reflection subgroup to be
$$K=G\cap K_N^+$$
with the intersection taken inside $U_N^+$. We say that $G$ appears as a soft liberation of its classical version $G_{class}=G\cap U_N$ when $G=<G_{class},K>$.
\end{definition}

These notions are important in the classification theory of compact quantum groups, and in connection with certain noncommutative geometry questions as well. As a first observation, with $K$ being as above, we have an intersection diagram, as follows:
$$\xymatrix@R=15mm@C=15mm{
K\ar[r]&G\\
K_{class}\ar[r]\ar[u]&G_{class}\ar[u]}$$

The soft liberation condition states that this diagram must be a generation diagram. We will be back to this in a moment, with some further theoretical comments. Let us first work out some examples. As a basic result, we have:

\begin{theorem}
The reflection subgroups of the basic unitary quantum groups
$$\xymatrix@R=15mm@C=15mm{
U_N\ar[r]&U_N^*\ar[r]&U_N^+\\
O_N\ar[r]\ar[u]&O_N^*\ar[r]\ar[u]&O_N^+\ar[u]}$$
are as follows,
$$\xymatrix@R=15mm@C=15mm{
K_N\ar[r]&K_N^*\ar[r]&K_N^+\\
H_N\ar[r]\ar[u]&H_N^*\ar[r]\ar[u]&H_N^+\ar[u]}$$
and these unitary quantum groups all appear via soft liberation.
\end{theorem}

\begin{proof}
The fact that the reflection subgroups of the quantum groups on the left are those on the right is clear in all cases, with the middle objects being by definition:
$$H_N^*=H_N\cap O_N^*\quad,\quad 
K_N^*=K_N\cap U_N^*$$

Regarding the second assertion, things are quite tricky here, as follows:

\medskip

(1) In the classical case there is nothing to prove, because any classical group is by definition a soft liberation of itself.

\medskip

(2) In the half-classical case the results are non-trivial, but can be proved by using the technology developed by Bichon and Dubois-Violette in \cite{bdu}.

\medskip

(3) In the free case the results are highly non-trivial, and the only known proof so far uses the recurrence methods developed by Chirvasitu in \cite{chi}. 
\end{proof}

Let us discuss now a number of more specialized classification results, for the twistable easy quantum groups, and for more general intermediate quantum groups as follows:
$$H_N\subset G\subset U_N^+$$

The idea will be that of viewing $G$ as sitting inside the standard cube:
$$\xymatrix@R=20pt@C=20pt{
&K_N^+\ar[rr]&&U_N^+\\
H_N^+\ar[rr]\ar[ur]&&O_N^+\ar[ur]\\
&K_N\ar[rr]\ar[uu]&&U_N\ar[uu]\\
H_N\ar[uu]\ar[ur]\ar[rr]&&O_N\ar[uu]\ar[ur]
}$$

To be more precise, beyond easiness, let us start with the following definition:

\begin{definition}
Associated to any closed subgroup $G_N\subset U_N^+$ are its classical, discrete and real versions, given by
$$G_N^c=G_N\cap U_N\quad,\quad 
G_N^d=G_N\cap K_N^+\quad,\quad
G_N^r=G_N\cap O_N^+$$
as well as its free, smooth and unitary versions, given by
$$G_N^f=<G_N,H_N^+>\quad,\quad 
G_N^s=<G_N,O_N>\quad,\quad 
G_N^u=<G_N,K_N>$$
where $<\,,>$ is the topological generation operation.
\end{definition}

We will need as well a second definition, as follows:

\begin{definition}
Associated to any closed subgroup $G_N\subset U_N^+$ are the mixes of its classical, discrete and real versions, given by
$$G_N^{cd}=G_N\cap K_N\quad,\quad 
G_N^{cr}=G_N\cap O_N^+\quad,\quad
G_N^{dr}=G_N\cap H_N^+$$
as well as the mixes of its free, smooth and unitary versions, given by
$$G_N^{fs}=<G_N,O_N^+>\quad,\quad 
G_N^{fu}=<G_N,K_N^+>\quad,\quad 
G_N^{us}=<G_N,U_N>$$
where $<\,,>$ is the topological generation operation.
\end{definition}

With the above notions in hand, we can formulate:

\index{orientability}

\begin{definition}
A closed subgroup $G_N\subset U_N^+$ is called ``oriented'' if
$$G_N=<G_N^{cd},G_N^{cr},G_N^{dr}>\quad,\quad 
G_N=G_N^{fs}\cap G_N^{fu}\cap G_N^{su}$$
and ``weakly oriented'' if the following conditions hold,
$$G_N=<G_N^c,G_N^d,G_N^r>\quad,\quad 
G_N=G_N^f\cap G_N^s\cap G_N^u$$
where the various versions are those in Definition 11.19 and Definition 11.20.
\end{definition}

We refer to the book \cite{ba4} for more on the above orientability notions, and what can be done with them, the idea being that we notably have a ``Ground Zero'' theorem, stating that when imposing very strong combinatorial axioms, including orientabilty, there will be only 8 quantum groups surviving, namely the vertices of the standard cube.

\bigskip

As a continuation of this, our claim is that more classification results are possible:

\bigskip

(1) In the classical case, we believe that the uniform, half-homogeneous, oriented groups should be the obvious examples of such groups. This is of course something quite heavy, well beyond easiness, with the potential tools available for proving such things coming from advanced finite group theory, and from Lie algebra theory.

\bigskip

(2) In the free case, under similar assumptions, we believe that the solutions should be the obvious examples of such quantum groups. This is something heavy, too, related to a well-known freeness conjecture, namely $<G_N,S_N^+>=\{G_N',S_N^+\}$, with the prime standing for the easy envelope. Indeed, assuming that we would have such a formula, along with some further ingredients, we can in principle work out our way inside the cube, from the edge and face projections to $G_N$ itself, and in this process $G_N$ would become easy.

\bigskip

(3) In the group dual case, the orientability axiom simplifies, because the group duals are discrete in our sense. We believe that the uniform, twistable, oriented group duals should appear as combinations of certain abelian groups, which appear in the classical case, with duals of varieties of real reflection groups, which appear in the real case. This is probably the easiest question in the present series, and the most reasonable one, to start with. However, there are no concrete results so far, in this direction.

\section*{11c. Toral subgroups}

Getting back now to Theorem 11.18 as it is, we are here into recent and interesting quantum group theory. We will discuss a bit later the concrete applications of Theorem 11.18. There is a connection here as well with the notion of diagonal torus, introduced in chapter 1. We can indeed refine Definition 11.17, in the following way:

\index{diagonal subgroup}
\index{soft liberation}
\index{hard liberation}

\begin{definition}
Given $H_N\subset G\subset U_N^+$, the diagonal tori $T=G\cap\mathbb T_N^+$ and reflection subgroups $K=G\cap K_N^+$ for $G$ and for $G_{class}=G\cap U_N$ form a diagram as follows:
$$\xymatrix@R=15mm@C=15mm{
T\ar[r]&K\ar[r]&G\\
T_{class}\ar[r]\ar[u]&K_{class}\ar[r]\ar[u]&G_{class}\ar[u]}$$
We say that $G$ appears as a soft/hard liberation when it is generated by $G_{class}$ and by $K/T$, which means that the right square/whole rectangle should be generation diagrams.
\end{definition}

All this is quite technical, and as a concrete result in connection with the above notion of hard liberation, we have the following statement, improving Theorem 11.18:

\begin{theorem}
The diagonal tori of the basic unitary quantum groups
$$\xymatrix@R=15mm@C=15mm{
U_N\ar[r]&U_N^*\ar[r]&U_N^+\\
O_N\ar[r]\ar[u]&O_N^*\ar[r]\ar[u]&O_N^+\ar[u]}$$
are as follows,
$$\xymatrix@R=15mm@C=15mm{
\mathbb T_N\ar[r]&\mathbb T_N^*\ar[r]&\mathbb T_N^+\\
T_N\ar[r]\ar[u]&T_N^*\ar[r]\ar[u]&T_N^+\ar[u]}$$
and these unitary quantum groups all appear via hard liberation.
\end{theorem}

\begin{proof}
The first assertion is something that we already know, from chapter 1. As for the second assertion, this can be proved by carefully examining the proof of Theorem 11.18, and performing some suitable modifications, where needed.
\end{proof}

As an interesting remark, some subtleties appear in the following way:

\begin{proposition}
The diagonal tori of the basic quantum reflection groups
$$\xymatrix@R=15mm@C=15mm{
K_N\ar[r]&K_N^*\ar[r]&K_N^+\\
H_N\ar[r]\ar[u]&H_N^*\ar[r]\ar[u]&H_N^+\ar[u]}$$
are as follows,
$$\xymatrix@R=15mm@C=15mm{
\mathbb T_N\ar[r]&\mathbb T_N^*\ar[r]&\mathbb T_N^+\\
T_N\ar[r]\ar[u]&T_N^*\ar[r]\ar[u]&T_N^+\ar[u]}$$
and these quantum reflection groups do not all appear via hard liberation.
\end{proposition}

\begin{proof}
The first assertion is clear, as a consequence of Theorem 11.23, because the diagonal torus is the same for a quantum group, and for its reflection subgroup:
$$G\cap\mathbb T_N^+=(G\cap K_N^+)\cap\mathbb T_N^+$$

Regarding the second assertion, things are quite tricky here, as follows:

\medskip

(1) In the classical case the hard liberation property definitely holds, because any classical group is by definition a hard liberation of itself.

\medskip

(2) In the half-classical case the answer is again positive, and this can be proved by using the technology developed by Bichon and Dubois-Violette in \cite{bdu}.

\medskip

(3) In the free case the hard liberation property fails, due to the intermediate quantum groups $H_N^{[\infty]}$, $K_N^{[\infty]}$, where ``hard liberation stops''. We will be back to this.
\end{proof}

As a conjectural solution to these latter difficulties, we have the notion of Fourier liberation, that we will discuss now. Let us first study the group dual subgroups of the arbitrary compact quantum groups $G\subset U_N^+$. To start with, we have:

\begin{proposition}
Let $G\subset U_N^+$ be a compact quantum group, and consider the group dual subgroups $\widehat{\Lambda}\subset G$, also called toral subgroups, or simply ``tori''.
\begin{enumerate}
\item In the classical case, where $G\subset U_N$ is a compact Lie group, these are the usual tori, where by torus we mean here closed abelian subgroup.

\item In the group dual case, $G=\widehat{\Gamma}$ with $\Gamma=<g_1,\ldots,g_N>$ being a discrete group, these are the duals of the various quotients $\Gamma\to\Lambda$.
\end{enumerate}
\end{proposition}

\begin{proof}
Both these assertions are elementary, as follows:

\medskip

(1) This follows indeed from the fact that a closed subgroup $H\subset U_N^+$ is at the same time classical, and a group dual, precisely when it is classical and abelian.

\medskip

(2) This follows from the general propreties of the Pontrjagin duality, and more precisely from the fact that the subgroups $\widehat{\Lambda}\subset\widehat{\Gamma}$ correspond to the quotients $\Gamma\to\Lambda$.
\end{proof}

At a more concrete level now, most of the tori that we met appear as diagonal tori. However, for certain quantum groups like the bistochastic ones, or the quantum permutation group ones, this torus collapses to $\{1\}$, and so it cannot be of use in the study of $G$. In order to deal with this issue, the idea, from \cite{bpa}, will be that of using:

\index{toral subgroup}
\index{standard tori}

\begin{proposition}
Given a closed subgroup $G\subset U_N^+$ and a matrix $Q\in U_N$, we let $T_Q\subset G$ be the diagonal torus of $G$, with fundamental representation spinned by $Q$:
$$C(T_Q)=C(G)\Big/\left<(QuQ^*)_{ij}=0\Big|\forall i\neq j\right>$$
This torus is then a group dual, $T_Q=\widehat{\Lambda}_Q$, where $\Lambda_Q=<g_1,\ldots,g_N>$ is the discrete group generated by the elements $g_i=(QuQ^*)_{ii}$, which are unitaries inside $C(T_Q)$.
\end{proposition}

\begin{proof}
This follows indeed from our results, because, as said in the statement, $T_Q$ is by definition a diagonal torus. Equivalently, since $v=QuQ^*$ is a unitary corepresentation, its diagonal entries $g_i=v_{ii}$, when regarded inside $C(T_Q)$, are unitaries, and satisfy:
$$\Delta(g_i)=g_i\otimes g_i$$

Thus $C(T_Q)$ is a group algebra, and more specifically we have $C(T_Q)=C^*(\Lambda_Q)$, where $\Lambda_Q=<g_1,\ldots,g_N>$ is the group in the statement, and this gives the result.
\end{proof}

Summarizing, associated to any closed subgroup $G\subset U_N^+$ is a whole family of tori, indexed by the unitaries $U\in U_N$. As a first result regarding these tori, we have:

\begin{theorem}
Any torus $T\subset G$ appears as follows, for a certain $Q\in U_N$:
$$T\subset T_Q\subset G$$
In other words, any torus appears inside a standard torus.
\end{theorem}

\begin{proof}
Given a torus $T\subset G$, we have an inclusion $T\subset G\subset U_N^+$. On the other hand, we know that each torus $T=\widehat{\Lambda}\subset U_N^+$, coming from a discrete group $\Lambda=<g_1,\ldots,g_N>$, has a fundamental corepresentation as follows, with $Q\in U_N$:
$$u=Q
\begin{pmatrix}
g_1\\
&\ddots\\
&&g_N
\end{pmatrix}
Q^*$$

But this shows that we have $T\subset T_Q$, and this gives the result.
\end{proof}

Let us do now some computations. In the classical case, the result is as follows:

\begin{proposition}
For a closed subgroup $G\subset U_N$ we have
$$T_Q=G\cap(Q^*\mathbb T^NQ)$$
where $\mathbb T^N\subset U_N$ is the group of diagonal unitary matrices.
\end{proposition}

\begin{proof}
This is indeed clear at $Q=1$, where $\Gamma_1$ appears by definition as the dual of the compact abelian group $G\cap\mathbb T^N$. In general, this follows by conjugating by $Q$.
\end{proof}

In the group dual case now, we have the following result, from \cite{bpa}:

\begin{proposition}
Given a discrete group $\Gamma=<g_1,\ldots,g_N>$, consider its dual compact quantum group $G=\widehat{\Gamma}$, diagonally embedded into $U_N^+$. We have then
$$\Lambda_Q=\Gamma\Big/\left<g_i=g_j\Big|\exists k,Q_{ki}\neq0,Q_{kj}\neq0\right>$$
with the embedding $T_Q\subset G=\widehat{\Gamma}$ coming from the quotient map $\Gamma\to\Lambda_Q$.
\end{proposition}

\begin{proof}
Assume indeed that $\Gamma=<g_1,\ldots,g_N>$ is a discrete group, with $\widehat{\Gamma}\subset U_N^+$ coming via $u=diag(g_1,\ldots,g_N)$. With $v=QuQ^*$, we have:
\begin{eqnarray*}
\sum_s\bar{Q}_{si}v_{sk}
&=&\sum_{st}\bar{Q}_{si}Q_{st}\bar{Q}_{kt}g_t\\
&=&\sum_t\delta_{it}\bar{Q}_{kt}g_t\\
&=&\bar{Q}_{ki}g_i
\end{eqnarray*}

Thus $v_{ij}=0$ for $i\neq j$ gives $\bar{Q}_{ki}v_{kk}=\bar{Q}_{ki}g_i$, which is the same as saying that $Q_{ki}\neq0$ implies $g_i=v_{kk}$. But this latter equality reads:
$$g_i=\sum_j|Q_{kj}|^2g_j$$

We conclude from this that $Q_{ki}\neq0,Q_{kj}\neq0$ implies $g_i=g_j$, as desired. As for the converse, this is elementary to establish as well.
\end{proof}

In view of the above, we can expect the collection $\{T_Q|Q\in U_N\}$ to encode various algebraic and analytic properties of $G$. We have the following result, from \cite{bpa}:

\index{toral conjectures}

\begin{theorem}
The following results hold, both for the compact Lie groups, and for the duals of the finitely generated discrete groups:
\begin{enumerate}
\item Generation: any closed quantum subgroup $G\subset U_N^+$ has the generation property $G=<T_Q|Q\in U_N>$. In other words, $G$ is generated by its tori.

\item Characters: if $G$ is connected, for any nonzero $P\in C(G)_{central}$ there exists $Q\in U_N$ such that $P$ becomes nonzero, when mapped into $C(T_Q)$.

\item Amenability: a closed subgroup $G\subset U_N^+$ is coamenable if and only if each of the tori $T_Q$ is coamenable, in the usual discrete group sense.

\item Growth: assuming $G\subset U_N^+$, the discrete quantum group $\widehat{G}$ has polynomial growth if and only if each the discrete groups $\widehat{T_Q}$ has polynomial growth.
\end{enumerate}
\end{theorem}

\begin{proof}
In the classical case, where $G\subset U_N$, the proof is elementary, based on standard facts from linear algebra, and goes as follows:

\medskip

(1) Generation. We use the following formula, established above: 
$$T_Q=G\cap Q^*\mathbb T^NQ$$

Since any group element $U\in G$ is diagonalizable, $U=Q^*DQ$ with $Q\in U_N,D\in\mathbb T^N$, we have $U\in T_Q$ for this value of $Q\in U_N$, and this gives the result.

\medskip

(2) Characters. We can take here $Q\in U_N$ to be such that $QTQ^*\subset\mathbb T^N$, where $T\subset U_N$ is a maximal torus for $G$, and this gives the result.

\medskip

(3) Amenability. This conjecture holds trivially in the classical case, $G\subset U_N$, due to the fact that these latter quantum groups are all coamenable.

\medskip

(4) Growth. This is something non-trivial, which is well-known from the theory of compact Lie groups, and we refer here to \cite{dpr} and related papers.

\medskip

Regarding now the group duals, here everything is trivial. Indeed, when the group duals are diagonally embedded we can take $Q=1$, and when the group duals are embedded by using a spinning matrix $Q\in U_N$, we can use precisely this matrix $Q$.
\end{proof}

The various statements above are conjectured to hold for any compact quantum group, and for a number of verifications, we refer to \cite{bpa} and to subsequent papers. In relation with our questions, let us focus now on the generation property. We will need:

\begin{proposition}
Given a closed subgroup $G\subset U_N^+$ and a matrix $Q\in U_N$, the corresponding standard torus and its Tannakian category are given by
$$T_Q=G\cap\mathbb T_Q\quad,\quad 
C_{T_Q}=<C_G,C_{\mathbb T_Q}>$$
where $\mathbb T_Q\subset U_N^+$ is the dual of the free group $F_N=<g_1,\ldots,g_N>$, with the fundamental corepresentation of $C(\mathbb T_Q)$ being the matrix $u=Qdiag(g_1,\ldots,g_N)Q^*$.
\end{proposition}

\begin{proof}
The first assertion comes from the well-known, and elementary fact that given two closed subgroups $G,H\subset U_N^+$, the corresponding quotient algebra $C(U_N^+)\to C(G\cap H)$ appears by dividing by the kernels of both the following quotient maps:
$$C(U_N^+)\to C(G)\quad,\quad 
C(U_N^+)\to C(H)$$

Indeed, the construction of $T_Q$ amounts in performing this operation, with $H=\mathbb T_Q$, and so we obtain $T_Q=G\cap\mathbb T_Q$, as claimed. As for the Tannakian category formula, this follows from this, and from the general duality formula $C_{G\cap H}=<C_G,C_H>$.
\end{proof}

We have the following Tannakian reformulation of the generation property:

\index{generation property}

\begin{theorem}
Given a closed subgroup $G\subset U_N^+$, the subgroup 
$$G'=<T_Q|Q\in U_N>$$
generated by its standard tori has the following Tannakian category:
$$C_{G'}=\bigcap_{Q\in U_N}<C_G,C_{\mathbb T_Q}>$$
In particular we have $G=G'$ when this intersection reduces to $C_G$.
\end{theorem}

\begin{proof}
Consider indeed the subgroup $G'\subset G$ in the statement. We have:
$$C_{G'}=\bigcap_{Q\in U_N}C_{T_Q}$$

Together with the formula in Proposition 11.31, this gives the result.
\end{proof}

The above result can be used for investigating the toral generation conjecture, but the combinatorics is quite difficult, and there are no results yet, along these lines. Let us further discuss now the toral generation property, with some modest results, regarding its behaviour with respect to product operations. We first have:

\begin{proposition}
Given two closed subgroups $G,H\subset U_N^+$, and $Q\in U_N$, we have:
$$<T_Q(G),T_Q(H)>\subset T_Q(<G,H>)$$
Also, the toral generation property is stable under the operation $<\,,>$.
\end{proposition}

\begin{proof}
The first assertion can be proved either by using Theorem 11.32, or directly. For the direct proof, which is perhaps the simplest, we have:
$$T_Q(G)
=G\cap\mathbb T_Q\subset<G,H>\cap\mathbb T_Q
=T_Q(<G,H>)$$

We have as well the following computation:
$$T_Q(H)
=H\cap\mathbb T_Q\subset<G,H>\cap\mathbb T_Q
=T_Q(<G,H>)$$

Now since $A,B\subset C$ implies $<A,B>\subset C$, this gives the result. Regarding now the second assertion, we have the following computation:
\begin{eqnarray*}
<G,H>
&=&<<T_Q(G)|Q\in U_N>,<T_Q(H)|Q\in U_N>>\\
&=&<T_Q(G),T_Q(H)|Q\in U_N>\\
&=&<<T_Q(G),T_Q(H)>|Q\in U_N>\\
&\subset&<T_Q(<G,H>)|Q\in U_N>
\end{eqnarray*}

Thus the quantum group $<G,H>$ is generated by its tori, as claimed.
\end{proof}

We have as well the following result:

\begin{proposition}
We have the following formula, for any $G,H$ and $R,S$:
$$T_{R\otimes S}(G\times H)=T_R(G)\times T_S(H)$$
Also, the toral generation property is stable under usual products $\times$.
\end{proposition}

\begin{proof}
The product formula is clear. Regarding now the second assertion, we have:
\begin{eqnarray*}
<T_Q(G\times H)|Q\in U_{MN}>
&\supset&<T_{R\otimes S}(G\times H)|R\in U_M,S\in U_N>\\
&=&<T_R(G)\times T_S(H)|R\in U_M,S\in U_N>\\
&=&<T_R(G)\times\{1\},\{1\}\times T_S(H)|R\in U_M,S\in U_N>\\
&=&<T_R(G)|R\in U_M>\times<T_G(H)|H\in U_N>\\
&=&G\times H
\end{eqnarray*}

Thus the quantum group $G\times H$ is generated by its tori, as claimed.
\end{proof}

\section*{11d. Fourier liberation}

Let us go back now to the quantum permutation groups. In relation with the tori,  let us start with the following basic fact, that we know well from chapter 9:

\begin{proposition}
Consider a discrete group generated by elements of finite order, written as a quotient group, as follows:
$$\mathbb Z_{N_1}*\ldots*\mathbb Z_{N_k}\to\Gamma$$
We have then an embedding $\widehat{\Gamma}\subset S_N^+$, where $N=N_1+\ldots+N_k$, with magic unitary
$$u=\begin{pmatrix}
F_{N_1}I_1F_{N_1}^*\\
&\ddots\\
&&F_{N_k}I_kF_{N_k}^*
\end{pmatrix}
\quad,\quad 
I_l=\begin{pmatrix}
1\\
&g_l\\
&&\ddots\\
&&&g_l^{N_l-1}
\end{pmatrix}$$
where $F_N=\frac{1}{\sqrt{N}}(w_N^{ij})$ with $w_N=e^{2\pi i/N}$, and $g_l$ is the standard generator of $\mathbb Z_{N_l}$.
\end{proposition}

\begin{proof}
We have indeed a sequence of embeddings and isomorphisms as follows:
\begin{eqnarray*}
\widehat{\Gamma}
&\subset&\widehat{\mathbb Z_{N_1}*\ldots*\mathbb Z_{N_k}}
=\widehat{\mathbb Z_{N_1}}\,\hat{*}\,\ldots\,\hat{*}\,\widehat{\mathbb Z_{N_k}}\\
&\simeq&\mathbb Z_{N_1}\,\hat{*}\,\ldots\,\hat{*}\,\mathbb Z_{N_k}
\subset S_{N_1}\,\hat{*}\,\ldots\,\hat{*}\,S_{N_k}\\
&\subset&S_{N_1}^+\,\hat{*}\,\ldots\,\hat{*}\,S_{N_k}^+
\subset S_N^+
\end{eqnarray*}

Thus, we are led to the conclusion in the statement.
\end{proof}

We have as well the following more specialized result, also from chapter 9:

\begin{theorem}
For the quantum permutation group $S_N^+$, the discrete group quotient $F_N\to\Lambda_Q$ with $Q\in U_N$ comes from the following relations:
$$\begin{cases}
g_i=1&{\rm if}\ \sum_lQ_{il}\neq 0\\
g_ig_j=1&{\rm if}\ \sum_lQ_{il}Q_{jl}\neq 0\\ 
g_ig_jg_k=1&{\rm if}\ \sum_lQ_{il}Q_{jl}Q_{kl}\neq 0
\end{cases}$$
Also, given a decomposition $N=N_1+\ldots+N_k$, for the matrix $Q=diag(F_{N_1},\ldots,F_{N_k})$, where $F_N=\frac{1}{\sqrt{N}}(\xi^{ij})_{ij}$ with $\xi=e^{2\pi i/N}$ is the Fourier matrix, we obtain
$$\Lambda_Q=\mathbb Z_{N_1}*\ldots*\mathbb Z_{N_k}$$
with dual embedded into $S_N^+$ in a standard way, as in Proposition 11.35.
\end{theorem}

\begin{proof}
This follows indeed from a direct computation, based on the definition of the diagonal tori $T_Q=\widehat{\Lambda_Q}$, explained in chapter 9.
\end{proof}

In connection with our liberation questions for the subgroups $G\subset S_N$, all this is quite interesting, and suggests formulating the following definition:

\index{Fourier torus}
\index{Fourier liberation}

\begin{definition}
Consider a closed subgroup $G\subset U_N^+$.
\begin{enumerate}
\item Its standard tori $T_F$, with $F=F_{N_1}\otimes\ldots\otimes F_{N_k}$, and $N=N_1+\ldots+N_k$ being regarded as a partition, are called Fourier tori.

\item In the case where we have $G_N=<G_N^c,(T_F)_F>$, we say that $G_N$ appears as a Fourier liberation of its classical version $G_N^c$.
\end{enumerate}
\end{definition}

The conjecture is then that the easy quantum groups should appear as Fourier liberations. With respect to the basic examples, the situation in the free case is as follows:

\medskip

(1) $O_N^+,U_N^+$ are diagonal liberations, so they are Fourier liberations as well.

\medskip

(2) $B_N^+,C_N^+$ are Fourier liberations too, with this being standard.

\medskip

(3) $S_N^+$ is a Fourier liberation too, being generated by its tori \cite{bcf}, \cite{chi}.

\medskip

(4) $H_N^+,K_N^+$ remain to be investigated, by using the general theory in \cite{rwe}.

\medskip

As a word of warning here, observe that an arbitrary classical group $G_N\subset U_N$ is not necessarily generated by its Fourier tori, and nor is an arbitrary discrete group dual, with spinned embedding. Thus, the Fourier tori, and the related notion of Fourier liberation, remain something quite technical, in connection with the easy case.

\bigskip

As an application of all this, let us go back to quantum permutation groups, and more specifically to the quantum symmetry groups of finite graphs, from chapter 9. One interesting question is whether $G^+(X)$ appears as a Fourier liberation of $G(X)$. Generally speaking, this is something quite difficult, because for the empty graph itself we are in need of the above-mentioned technical results from \cite{bcf}, \cite{chi}.

\bigskip

In order to discuss this, let us begin with the following elementary statement:

\begin{proposition}
In order for a closed subgroup $G\subset U_K^+$ to appear as $G=G^+(X)$, for a certain graph $X$ having $N$ vertices, the following must happen:
\begin{enumerate}
\item We must have a representation $G\subset U_N^+$.

\item This representation must be magic, $G\subset S_N^+$.

\item We must have a graph $X$ having $N$ vertices, such that $d\in End(u)$.

\item $X$ must be in fact such that the Tannakian category of $G$ is precisely $<d>$.
\end{enumerate}
\end{proposition}

\begin{proof}
This is more of an empty statement, coming from the definition of the quantum automorphism group $G^+(X)$, as formulated in chapter 5.
\end{proof}

In the group dual case, forgetting about Fourier transforms, and imagining that we are at step (1) in the general strategy outlined in Proposition 11.38, we must  compute the Tannakian category of $\widehat{\Gamma}\subset U_N^+$, diagonally embedded, for the needs of (3,4). We have:

\begin{proposition}
Given a discrete group $\Gamma=<g_1,\ldots,g_N>$, embed diagonally $\widehat{\Gamma}\subset U_N^+$, via the unitary matrix $u=diag(g_1,\ldots,g_N)$. We have then the formula
$$Hom(u^{\otimes k},u^{\otimes l})=\left\{T=(T_{j_1\ldots j_l,i_1\ldots i_k})\Big|g_{i_1}\ldots g_{i_k}\neq g_{j_1}\ldots g_{j_l}\implies T_{j_1\ldots j_l,i_1\ldots i_k}=0\right\}$$
and in particular, with $k=l=1$, we have the formula
$$End(u)=\left\{T=(T_{ji})\Big|g_i\neq g_j\implies T_{ji}=0\right\}$$
with the linear maps being identified with the corresponding scalar matrices.
\end{proposition}

\begin{proof}
This is well-known, and elementary, with the first assertion coming from:
\begin{eqnarray*}
T\in Hom(u^{\otimes k},u^{\otimes l})
&\iff&Tu^{\otimes k}=u^{\otimes l}T\\
&\iff&(Tu^{\otimes k})_{j_1\ldots j_l,i_1\ldots i_k}=(u^{\otimes l}T)_{j_1\ldots j_l,i_1\ldots i_k}\\
&\iff&T_{j_1\ldots j_l,i_1\ldots i_k}g_{i_1}\ldots g_{i_k}=g_{j_1}\ldots g_{j_l}T_{j_1\ldots j_l,i_1\ldots i_k}\\
&\iff&T_{j_1\ldots j_l,i_1\ldots i_k}(g_{i_1}\ldots g_{i_k}-g_{j_1}\ldots g_{j_l})=0
\end{eqnarray*}

As for the second assertion, this follows from the first one.
\end{proof}

Let us go ahead now, with respect to the general strategy outlined in Proposition 11.38, and apply \cite{bi3} in order to solve (2), and then reformulate (3,4), by using Proposition 11.39, and by choosing to put the multi-Fourier transform on the graph part. We are led in this way into the following refinement of Proposition 11.38, in the group dual setting:

\begin{theorem}
In order for a group dual $\widehat{\Gamma}$ as in Proposition 11.35 to appear as $G=G^+(X)$, for a certain graph $X$ having $N$ vertices, the following must happen:
\begin{enumerate}
\item First, we need a quotient map $\mathbb Z_{N_1}*\ldots*\mathbb Z_{N_k}\to\Gamma$.

\item Let $u=diag(I_1,\ldots,I_k)$, with $I_l=diag(\mathbb Z_{N_l})$, for any $l$.

\item Consider also the matrix $F=diag(F_{N_1},\ldots,F_{N_k})$.

\item We must then have a graph $X$ having $N$ vertices.

\item This graph must be such that $F^*dF\neq0\implies I_i=I_j$.

\item In fact, $<F^*dF>$ must be the category in Proposition 11.39.
\end{enumerate}
\end{theorem}

\begin{proof}
This is something rather informal, the idea being as follows:

\medskip

(1) This is the assumption of Proposition 11.35, explained above, with the remark that we can add to this a unitary base change, as piece of data.

\medskip

(2) This is just a notation, with $I_l=diag(\mathbb Z_{N_l})$ meaning that $I_l$ is the diagonal matrix formed by $1,g,g^2,\ldots,g^{N_l-1}$, with $g\in\mathbb Z_{N_l}$ being the standard generator.

\medskip

(3) This is another notation, with each Fourier matrix $F_{N_l}$ being the standard one, namely $F_{N_l}=\frac{1}{\sqrt{N_l}}(w^{ij})$, with $w=e^{2\pi i/N_l}$, and with indices $0,1,\ldots,N_l-1$.

\medskip

(4) This is a just an informal statement, with the precise graph formalism to be clarified later on, in view of the fact that $X$ will get Fourier-transformed anyway.

\medskip

(5) This is an actual result, our claim being that the condition $d\in End(u)$ from Proposition 11.38 (3) is equivalent to the condition $F^*dF\neq0\implies I_i=I_j$ in the statement. Indeed, we know that with $F,I$ being as in the statement, we have $u=FIF^*$. Now with this formula in hand, we have the following equivalences:
\begin{eqnarray*}
\widehat{\Gamma}\curvearrowright X
&\iff&du=ud\\
&\iff&dFIF^*=FIF^*d\\
&\iff&[F^*dF,I]=0
\end{eqnarray*}

Also, since the matrix $I$ is diagonal, with $M=F^*dF$ have:
\begin{eqnarray*}
MI=IM
&\iff&(MI)_{ij}=(IM)_{ij}\\
&\iff&M_{ij}I_j=I_iM_{ij}\\
&\iff&[M_{ij}\neq 0\implies I_i=I_j]
\end{eqnarray*}

We therefore conclude that we have, as desired:
$$\widehat{\Gamma}\curvearrowright X\iff[F^*dF\neq0\implies I_i=I_j]$$

(6) This is the Tannakian condition in Proposition 11.38 (4), with reference to the explicit formula for the Tannakian category of $G=\widehat{\Gamma}$ given in Proposition 11.39.
\end{proof}

Going ahead now, in connection with the Fourier tori, we have:

\index{Fourier torus}
\index{Fourier liberation}

\begin{proposition}
The Fourier tori of $G^+(X)$ are the biggest quotients
$$\mathbb Z_{N_1}*\ldots*\mathbb Z_{N_k}\to\Gamma$$
whose duals act on the graph, $\widehat{\Gamma}\curvearrowright X$.
\end{proposition}

\begin{proof}
We have indeed the following computation, at $F=1$:
\begin{eqnarray*}
C(T_1(G^+(X)))
&=&C(G^+(X))/<u_{ij}=0,\forall i\neq j>\\
&=&[C(S_N^+)/<[d,u]=0>]/<u_{ij}=0,\forall i\neq j>\\
&=&[C(S_N^+)/<u_{ij}=0,\forall i\neq j>]/<[d,u]=0>\\
&=&C(T_1(S_N^+))/<[d,u]=0>
\end{eqnarray*}

Thus, we obtain the result, with the remark that the quotient that we are interested in appears via relations of type $d_{ij}=1\implies g_i=g_j$. The proof in general is similar.
\end{proof}

In connection now with the above-mentioned questions, we have:

\begin{theorem}
Consider the following conditions:
\begin{enumerate}
\item We have $G(X)=G^+(X)$.

\item $G(X)\subset G^+(X)$ is a Fourier liberation.

\item $\widehat{\Gamma}\curvearrowright X$ implies that $\Gamma$ is abelian.
\end{enumerate}
We have then the equivalence $(1)\iff(2)+(3)$.
\end{theorem}

\begin{proof}
This is something elementary, the proof being as follows:

\medskip

$(1)\implies(2,3)$ Here both the implications are trivial.

\medskip

$(2,3)\implies(1)$ Assuming $G(X)\neq G^+(X)$, from (2) we know that $G^+(X)$ has at least one non-classical Fourier torus, and this contradicts (3).
\end{proof}

With this in hand, the question is whether $(3)\implies(1)$ holds. This is a good question, which in practice would make connections between the various conjectures that can be made about a given graph $X$, and its quantum symmetry group $G^+(X)$.

\bigskip

As an illustration for the potential interest of such considerations, it is known from \cite{lmr} that the random graphs have no quantum symmetries, with this being something advanced. Our point now is that, assuming that one day the general compact quantum Lie group theory will solve its Weyl-type questions in relation with the tori, and in particular include, as a theorem, the fact that any $G^+(X)$ appears as a Fourier liberation of $G(X)$, this deep graph result from \cite{lmr} would become accessible as well via its particular case for the group dual subgroups, which is something elementary, as follows:

\begin{theorem}
For a finite graph $X$, the probability for having an action
$$\widehat{\Gamma}\curvearrowright X$$
with $\Gamma$ being a non-abelian group goes to $0$ with $|X|\to\infty$.
\end{theorem}

\begin{proof}
This is something quite elementary, the idea being as follows:

\medskip

(1) First of all, the graphs $X$ having a fixed number $N\in\mathbb N$ of vertices correspond to the matrices $d\in M_N(0,1)$ which are symmetric, and have 0 on the diagonal. The probability mentioned in the statement is the uniform one on such 0-1 matrices.

\medskip

(2) Regarding now the proof, our claim is that this should come in a quite elementary way, from the $du=ud$ condition, as reformulated before. Indeed, observe first that in the cyclic case, where $F=F_N$ is a usual Fourier matrix, associated to a cyclic group $\mathbb Z_N$, we have the following formula, with $w=e^{2\pi i/N}$:
\begin{eqnarray*}
(F^*dF)_{ij}
&=&\sum_{kl}(F^*)_{ik}d_{kl}F_{lj}\\
&=&\sum_{kl}w^{lj-ik}d_{kl}\\
&=&\sum_{k\sim l}w^{lj-ik}
\end{eqnarray*}

(3) In the general setting now, where we have a quotient map $\mathbb Z_{N_1}*\ldots*\mathbb Z_{N_k}\to\Gamma$, with $N_1+\ldots+N_k=N$, the computation is similar, as follows, with $w_i=e^{2\pi i/N_i}$:
\begin{eqnarray*}
(F^*dF)_{ij}
&=&\sum_{kl}(F^*)_{ik}d_{kl}F_{lj}\\
&=&\sum_{k\sim l}(F^*)_{ik}F_{lj}\\
&=&\sum_{k:i,l:j,k\sim l}(w_{N_i})^{-ik}(w_{N_j})^{lj}
\end{eqnarray*}

Here the conditions $k:i$ and $l:j$ refer to the fact that $k,l$ must belong respectively to the same matrix blocks as $i,j$, with respect to the partition $N_1+\ldots+N_k=N$, and $k\sim l$ means as usual that there is an edge between $k,l$, in the graph $X$.

\medskip

(4) The point now is that with the partition $N_1+\ldots+N_k=N$ fixed, and with $d\in M_N(0,1)$ being random, we have $(F^*dF)_{ij}\neq 0$ almost everywhere in the $N\to\infty$ limit, so we have $I_i=I_j$ almost everywhere, and so abelianity of $\Gamma$, with $N\to\infty$.
\end{proof}

\section*{11e. Exercises} 

Things have been quite technical in this chapter, and as a unique exercise, summarizing the main problems that we have been talking about here, we have:

\begin{exercise}
Find an abstract framework for the ``quantum reflection groups'', as intermediate quantum groups of type
$$A_N\subset G\subset K_N^+$$
covering at the same time all the intermediate easy quantum groups of type
$$S_N\subset G\subset K_N^+$$
and all the classical, non-exceptional complex reflection groups, namely
$$H_N^{sd}=\left\{U\in H_N^s\Big|(\det U)^d=1\right\}$$
and then start classifying such beasts. 
\end{exercise}

The point here is that in the classical case, there is a definition and then classification result for the complex reflection groups, with the classification stating that we have as exemples the above groups $H_N^{sd}$, plus a number of exceptional examples, which can be classified as well. The problem is that of finding the correct quantum extension of this.

\chapter{Analytic aspects}

\section*{12a. Weingarten estimates}

We discuss here a number of advanced probabilistic aspects, with improvements of some of the character results from chapter 2, and with De Finetti theorems as well, following \cite{ez1}, \cite{ez2}. We will insist on our main examples of quantum groups, namely:
$$\xymatrix@R=15mm@C=15mm{
S_N^+\ar[r]&O_N^+\\
S_N\ar[r]\ar[u]&O_N\ar[u]
}$$

Here the quantum groups $S_N,S_N^+$ need no presentation, being those that we are mainly interested in. As for $O_N,O_N^+$, we will study them too, their projective versions being basic examples of quantum groups $S_Z,S_Z^+$, due to the isomorphism $S_{M_N}^+=PO_N^+$.

\bigskip

We already know, since chapter 2, that for everything probability we will need at some point Weingarten function estimates. So, let us begin with this topic, which is of independent interest as well. Let us first recall the Weingarten formula:

\begin{theorem}
The integrals over an easy quantum group $G\subset_uO_N^+$, coming from a category of partitions $D=(D(k,l))$, are given by the formula
$$\int_Gu_{i_1j_1}\ldots u_{i_kj_k}=\sum_{\pi,\nu\in D(k)}\delta_\pi(i)\delta_\nu(j)W_{kN}(\pi,\nu)$$
where $D(k)=D(0,k)$, the $\delta$ symbols are Kronecker type symbols, $\delta\in\{0,1\}$, and where $W_{kN}=G_{kN}^{-1}$ with $G_{kn}(\pi,\nu)=N^{|\pi\vee\nu|}$ is the Weingarten matrix.
\end{theorem}

\begin{proof}
This is something that we know from chapter 2, coming from the fact that the integrals in the statement, with $k\in\mathbb N$ fixed and with the multi-indices $i,j$ varying, form a $N^k\times N^k$ matrix which is the orthogonal projection onto $Fix(u^{\otimes k})$.
\end{proof}

In general, the computation of the Weingarten matrix is something quite delicate. However, in the case of the symmetric group $S_N$, the situation is in fact very simple, because we can explicitly compute and estimate this matrix, as follows:

\index{Weingarten function}

\begin{theorem}
For the group $S_N$ the Weingarten function is given by
$$W_{kN}(\pi,\nu)=\sum_{\tau\leq\pi\wedge\nu}\mu(\tau,\pi)\mu(\tau,\nu)\frac{(N-|\tau|)!}{N!}$$
and satisfies the folowing estimate,
$$W_{kN}(\pi,\nu)=N^{-|\pi\wedge\nu|}(
\mu(\pi\wedge\nu,\pi)\mu(\pi\wedge\nu,\nu)+O(N^{-1}))$$
with $\mu$ being the M\"obius function of $P(k)$.
\end{theorem}

\begin{proof}
The first assertion follows from the usual Weingarten formula, from Theorem 12.1. Indeed, in that formula the integrals on the left are in fact known, from the explicit integration formula over $S_N$ that we established in chapter 2, namely:
$$\int_{S_N}g_{i_1j_1}\ldots g_{i_kj_k}=\begin{cases}
\frac{(N-|\ker i|)!}{N!}&{\rm if}\ \ker i=\ker j\\
0&{\rm otherwise}
\end{cases}$$

But this allows the computation of the right term, via the M\"obius inversion formula, and we get the result. As for the second assertion, this follows from the first one.
\end{proof}

Regarding now the quantum group $S_N^+$, that we are particularly interested in here, let us begin with some explicit computations. We first have the following simple and final result at $k=2,3$, directly in terms of the quantum group integrals:

\begin{proposition}
At $k=2,3$ we have the following estimate:
$$\int_{S_N^+}u_{i_1j_1}\ldots u_{i_kj_k}=\begin{cases}
0&(\ker i\neq\ker j)\\
\simeq N^{-|\ker i|}&(\ker i=\ker j)
\end{cases}$$
\end{proposition}

\begin{proof}
Since at $k\leq3$ we have $NC(k)=P(k)$, the Weingarten integration formulae for $S_N$ and $S_N^+$ coincide, and we obtain, by using the above formula for $S_N$:
$$\int_{S_N^+}u_{i_1j_1}\ldots u_{i_kj_k}
=\delta_{\ker i,\ker j}\frac{(N-|\ker i|)!}{N!}$$

Thus, we obtain the formula in the statement.
\end{proof}

In general now, the idea will be that of working out a ``master estimate'' for the Weingarten function, as above. Before starting, let us record the formulae at $k=2,3$, which will be useful later, as illustrations. At $k=2$, with indices $||,\sqcap$ as usual, and with the convention that $\approx$ means componentwise dominant term, we have:
$$W_{2N}\approx\begin{pmatrix}N^{-2}&-N^{-2}\\-N^{-2}&N^{-1}\end{pmatrix}$$

At $k=3$ now, with indices $|||,|\sqcap,\sqcap|,\sqcap\hskip-3.2mm{\ }_|\,,\sqcap\hskip-0.8mm\sqcap$, and same meaning for $\approx$, we have:
$$W_{3N}\approx\begin{pmatrix}
N^{-3}&-N^{-3}&-N^{-3}&-N^{-3}&2N^{-3}\\
-N^{-3}&N^{-2}&N^{-3}&N^{-3}&-N^{-2}\\
-N^{-3}&N^{-3}&N^{-2}&N^{-3}&-N^{-2}\\
-N^{-3}&N^{-3}&N^{-3}&N^{-2}&-N^{-2}\\
2N^{-3}&-N^{-2}&-N^{-2}&-N^{-2}&N^{-1}
\end{pmatrix}$$

In order to deal now with the general case, let us start with some standard facts:

\begin{proposition}
The following happen, regarding the partitions in $P(k)$:
\begin{enumerate}
\item $|\pi|+|\nu|\leq|\pi\vee\nu|+|\pi\wedge\nu|$.

\item $|\pi\vee\tau|+|\tau\vee\nu|\leq|\pi\vee\nu|+|\tau|$.

\item $d(\pi,\nu)=\frac{|\pi|+|\nu|}{2}-|\pi\vee\nu|$ is a distance.
\end{enumerate}
\end{proposition}

\begin{proof}
All this is well-known, the idea being as follows:

\medskip

(1) This comes from the fact that $P(k)$ is a semi-modular lattice.

\medskip

(2) This follows from (1), as explained for instance in \cite{ez1}.

\medskip

(3) This follows from (2), which says that the following holds:
$$\frac{|\pi|+|\tau|}{2}-d(\pi,\tau)+\frac{|\tau|+|\nu|}{2}-d(\tau,\nu)
\leq\frac{|\pi|+|\nu|}{2}-d(\pi,\nu)+|\tau|$$

Thus, we obtain the triangle inequality, and the other axioms are all clear.
\end{proof}

Actually in what follows we will only need (3) in the above statement. For more on this, and on the geometry and combinatorics of partitions, we refer to \cite{nsp}. As a main result now regarding the Weingarten functions, we have:

\index{series expansion}
\index{Weingarten function}
\index{geodesicity defect}

\begin{theorem}
The Weingarten matrix $W_{kN}$ has a series expansion in $N^{-1}$,
$$W_{kN}(\pi,\nu)=N^{|\pi\vee\nu|-|\pi|-|\nu|}\sum_{g=0}^\infty K_g(\pi,\nu)N^{-g}$$
where the various objects on the right are defined as follows:
\begin{enumerate}
\item A path from $\pi$ to $\nu$ is a sequence as follows:
$$p=[\pi=\tau_0\neq\tau_1\neq\ldots\neq\tau_r=\nu]$$

\item The signature of such a path is $+$ when $r$ is even, and $-$ when $r$ is odd.

\item The geodesicity defect of such a path is:
$$g(p)=\sum_{i=1}^rd(\tau_{i-1},\tau_i)-d(\pi,\nu)$$

\item $K_g$ counts the signed paths from $\pi$ to $\nu$, with geodesicity defect $g$.
\end{enumerate} 
\end{theorem}

\begin{proof}
We recall that the Weingarten matrix $W_{kN}$ appears as the inverse of the Gram matrix $G_{kN}$, which is given by the following formula:
$$G_{kN}(\pi,\nu)=N^{|\pi\vee\nu|}$$

Now observe that the Gram matrix can be written in the following way:
\begin{eqnarray*}
G_{kN}(\pi,\nu)
&=&N^{|\pi\vee\nu|}\\
&=&N^{\frac{|\pi|}{2}}N^{|\pi\vee\nu|-\frac{|\pi|+|\nu|}{2}}N^{\frac{|\nu|}{2}}\\
&=&N^{\frac{|\pi|}{2}}N^{-d(\pi,\nu)}N^{\frac{|\nu|}{2}}
\end{eqnarray*}

This suggests considering the following diagonal matrix:
$$\Delta=diag(N^{\frac{|\pi|}{2}})$$

So, let us do this, and consider as well the following matrix:
$$H(\pi,\nu)=\begin{cases}
0&(\pi=\nu)\\
N^{-d(\pi,\nu)}&(\pi\neq\nu)
\end{cases}$$

In terms of these two matrices, the above formula for $G_{kN}$ simply reads:
$$G_{kN}=\Delta(1+H)\Delta$$

Thus, the Weingarten matrix $W_{kN}$ is given by the following formula:
$$W_{kN}=\Delta^{-1}(1+H)^{-1}\Delta^{-1}$$

In order to compute now the inverse of $1+H$, we will use the following formula:
$$(1+H)^{-1}=1-H+H^2-H^3+\ldots$$

Consider indeed the set $P_r(\pi,\nu)$ of length $r$ paths between $\pi$ and $\nu$. We have:
\begin{eqnarray*}
H^r(\pi,\nu)
&=&\sum_{p\in P_r(\pi,\nu)}H(\tau_0,\tau_1)\ldots H(\tau_{r-1},\tau_r)\\
&=&\sum_{p\in P_r(\pi,\nu)}N^{-d(\pi,\nu)-g(p)}
\end{eqnarray*}

Thus by using $(1+H)^{-1}=1-H+H^2-H^3+\ldots$ we obtain:
\begin{eqnarray*}
(1+H)^{-1}(\pi,\nu)
&=&\sum_{r=0}^\infty(-1)^rH^r(\pi,\nu)\\
&=&N^{-d(\pi,\nu)}\sum_{r=0}^\infty\sum_{p\in P_r(\pi,\nu)}(-1)^rN^{-g(p)}
\end{eqnarray*}

It follows that the Weingarten matrix is given by the following formula:
\begin{eqnarray*}
W_{kN}(\pi,\nu)
&=&\Delta^{-1}(\pi)(1+H)^{-1}(\pi,\nu)\Delta^{-1}(\nu)\\
&=&N^{-\frac{|\pi|}{2}-\frac{|\nu|}{2}-d(\pi,\nu)}\sum_{r=0}^\infty\sum_{p\in P_r(\pi,\nu)}(-1)^rN^{-g(p)}\\
&=&N^{|\pi\vee\nu|-|\pi|-|\nu|}\sum_{r=0}^\infty\sum_{p\in P_r(\pi,\nu)}(-1)^rN^{-g(p)}
\end{eqnarray*}

Now by rearranging the various terms in the above double sum according to their geodesicity defect $g=g(p)$, this gives the following formula:
$$W_{kN}(\pi,\nu)=N^{|\pi\vee\nu|-|\pi|-|\nu|}\sum_{g=0}^\infty K_g(\pi,\nu)N^{-g}$$

Thus, we are led to the conclusion in the statement.
\end{proof}

As an illustration for all this, we have the following explicit estimates:

\begin{theorem}
Consider an easy quantum group $G=(G_N)$, coming from a category of partitions $D=(D(k))$. For any $\pi\leq\nu$ we have the estimate
$$W_{kN}(\pi,\nu)=N^{-|\pi|}(\mu(\pi,\nu)+O(N^{-1}))$$
and for $\pi,\nu$ arbitrary we have
$$W_{kN}(\pi,\nu)=O(N^{|\pi\vee\nu|-|\pi|-|\nu|})$$
with $\mu$ being the M\"obius function of $D(k)$.
\end{theorem}

\begin{proof}
We have two assertions here, the idea being as follows:

\medskip

(1) The first estimate is clear from the formula in Theorem 12.5, namely: 
$$W_{kN}(\pi,\nu)=N^{|\pi\vee\nu|-|\pi|-|\nu|}\sum_{g=0}^\infty K_g(\pi,\nu)N^{-g}$$

(2) In the case $\pi\leq\nu$ it is known that $K_0$ coincides with the M\"obius function of $NC(k)$, as explained for instance in \cite{ez1}, so we obtain once again from Theorem 12.5 the fine estimate in the statement as well, namely:
$$W_{kN}(\pi,\nu)=N^{-|\pi|}(\mu(\pi,\nu)+O(N^{-1}))\qquad\forall\pi\leq\nu$$

Observe that, by symmetry of $W_{kN}$, we obtain as well that we have:
$$W_{kN}(\pi,\nu)=N^{-|\nu|}(\mu(\nu,\pi)+O(N^{-1}))\qquad\forall\pi\geq\nu$$

Thus, we are led to the conclusions in the statement.
\end{proof}

When $\pi,\nu$ are not comparable by $\leq$, things are quite unclear. The simplest example appears at $k=3$, where we have the following formula, which is elementary:
$$W_{3N}(|\sqcap,\sqcap|)\simeq N^{-3}$$

As for the corresponding coefficient, $K_0(|\sqcap,\sqcap|)=1$, this is definitely not the M\"obius function, which vanishes for partitions which are not comparable by $\leq$. According to Theorem 12.5, this is rather the number of signed geodesic paths from $|\sqcap$ to $\sqcap|$.

\section*{12b. Traces of powers}

Let us discuss now, following \cite{ez1}, the computation of the asymptotic laws of the following variables, depending on an integer $k\in\mathbb N$:
$$\chi^{(k)}=Tr(u^k)$$

\index{Diaconis-Shahshahani variables}

These variables, called Diaconis-Shahshahani variables \cite{dsh}, generalize the usual characters, which appear at $k=1$. Let us start with the following standard definition:

\begin{definition}
Associated to any integers $k_1,\ldots,k_s\in\mathbb N$ is the ``trace'' permutation $\gamma\in S_k$, having cycles as follows, with $k=\sum k_i$: 
$$(1,\ldots,k_1)\quad,\quad 
(k_1+1,\ldots,k_1+k_2)\quad,\quad 
\ldots\quad,\quad 
(k-k_s+1,\ldots,k)$$
We also denote by $\gamma(\sigma)$ the partition given by $i\sim_\sigma j\iff \gamma(i)\sim_{\gamma(\sigma)}\gamma(j)$. 
\end{definition}

These conventions are standard in free probability, and we refer to \cite{nsp} for more on all this. Now with these conventions, we have the following result:

\begin{proposition}
Given an easy quantum group $G$, we have:
$$\int_GTr(u^{k_1})\ldots Tr(u^{k_s})\,du=\#\left\{\pi\in D_k\Big|\pi=\gamma(\pi)\right\} +O(N^{-1})$$
If $G$ is classical, this estimate is exact, without any lower order corrections.
\end{proposition}

\begin{proof}
We have two assertions to be proved, the idea being as follows:

\medskip

(1) Let $I$ be the integral to be computed. According to the definition of $\gamma$, we have:
\begin{eqnarray*}
I
&=&\int_GTr(u^{k_1})\ldots Tr(u^{k_s})\,du\\
&=&\sum_{i_1\ldots i_k}\int_G(u_{i_1i_2}\ldots u_{i_ki_1})\ldots\ldots (u_{i_{k-k_s+1}i_{k-k_s+2}}\ldots u_{i_ki_{k-k_s+1}})\\
&=&\sum_{i_1\ldots i_k}\int_Gu_{i_1i_{\gamma(1)}}\ldots u_{i_ki_{\gamma(k)}}
\end{eqnarray*}

We use now the Weingarten formula. We obtain the following formula:
\begin{eqnarray*}
I
&=&\sum_{i_1\ldots i_k}\sum_{\pi\leq\ker i}\sum_{\sigma\leq\ker i\gamma}W_{kN}(\pi,\sigma)\\
&=&\sum_{i_1\ldots i_k}\sum_{\pi\leq\ker i}\sum_{\gamma(\sigma)\leq\ker i}W_{kN}(\pi,\sigma)\\
&=&\sum_{\pi,\sigma\in D_k}N^{|\pi\vee\gamma(\sigma)|}W_{kN}(\pi,\sigma)\\
&=&\sum_{\pi,\sigma\in D_k}N^{|\pi\vee\gamma(\sigma)|}N^{|\pi\vee\sigma|-|\pi|-|\sigma|}(1+O(N^{-1}))
\end{eqnarray*}

Let us look now at the power of $N$ in the above, namely:
$$N^{|\pi\vee\gamma(\sigma)|+|\pi\vee\sigma|-|\pi|-|\sigma|}$$

The leading order is $N^0$, which is achieved if and only if $\sigma\geq\pi$ and $\pi\geq\gamma(\sigma)$, or equivalently when $\pi=\sigma=\gamma(\sigma)$. But this gives the formula in the statement.

\medskip

(2) In the classical case, instead of using the approximation for $W_{Nk}(\pi,\sigma)$, we can write $N^{|\pi\vee\gamma(\sigma)|}=G_{kN}(\gamma(\sigma),\pi)$, and we can continue as follows:
\begin{eqnarray*}
I
&=&\sum_{\pi,\sigma\in D_k}G_{kN}(\gamma(\sigma),\pi)W_{kN}(\pi,\sigma)\\
&=&\sum_{\sigma\in D_k}\delta(\gamma(\sigma),\sigma)\\
&=&\#\{\sigma\in D_k|\sigma=\gamma(\pi)\}
\end{eqnarray*}

Thus, we are led to the conclusion in the statement.
\end{proof}

If $c$ is a cycle we use the notation $c^1=c$, and $c^*$= cycle opposite to $c$. We have the following definition, generalizing Definition 12.7:

\index{trace permutation}

\begin{definition}
Associated to any $k_1,\ldots,k_s\in\mathbb N$ and any $e_1,\ldots,e_s\in\{1,*\}$ is the trace permutation $\gamma\in S_k$, with $k=\sum k_i$, having as cycles
$$(1,\ldots,k_1)^{e_1}\quad,\quad 
(k_1+1,\ldots,k_1+k_2)^{e_2}\quad,\quad 
\ldots\quad,\quad 
(k-k_s+1,\ldots,k)^{e_s}$$
called trace permutation associated to $k_1,\ldots,k_s\in\mathbb N$ and $e_1,\ldots,e_s\in\{1,*\}$.
\end{definition}

Again, this convention is standard in free probability, and we refer to \cite{nsp} for more on all this. With this convention, Proposition 12.8 extends as follows:

\begin{proposition}
Given an easy quantum group $G$, we have:
$$\int_GTr(u^{k_1})^{e_1}\ldots Tr(u^{k_s})^{e_s}\,du=\#\left\{\pi\in D_k\Big|\pi=\gamma(\pi)\right\} +O(N^{-1})$$
If $G$ is classical, this estimate is exact, without any lower order corrections.
\end{proposition}

\begin{proof}
This is similar to the proof of Proposition 12.8, by performing modifications where needed, and with the computations being the same as before. See \cite{ez1}.
\end{proof}

In terms of cumulants, we have the following result, also from \cite{ez1}:

\begin{proposition}
For $G=O_N,S_N$ we have the following cumulant formula:
$$c_r(Tr(u^{k_1})^{e_1},\ldots,Tr(u^{k_r})^{e_r})
=\#\left\{\pi\in D_k\Big|\pi\vee\gamma=1_k,\,\pi=\gamma(\pi)\right\}$$
Also, for $G=O_N^+,S_N^+$ we have the following free cumulant formula:
$$\kappa_r(Tr(u^{k_1})^{e_1},\ldots,Tr(u^{k_r})^{e_r})
=\#\left\{\pi\in D_k\Big|\pi\vee\gamma=1_k,\,\pi=\gamma(\pi)\right\}+O(N^{-1})$$
\end{proposition}

\begin{proof}
We have two assertions to be proved, the idea being as follows:

\medskip

(1) Let $c_r$ be the considered cumulant. We write, for those partitions $\pi\in P_k$ such that the restriction of $\pi$ to a block of $\sigma$
is an element in the corresponding set $D_{|v|}$:
$$D_\sigma=\left\{\pi\in P_k\Big|p_{|v}\in D_{|v|}\, \forall v\in\sigma\right\}$$

We have then the following equivalent formula:
$$D_\sigma=\left\{\pi\in D_k\Big|\pi\leq\sigma^\gamma\right\}$$

Then, by the definition of the classical cumulants via M\"obius inversion, we get:
\begin{eqnarray*}
c_r
&=&\sum_{\sigma\in P(r)}\mu(\sigma,1_r)\cdot\#\{\pi\in D_\sigma|\pi=\gamma(\pi)\}\\
&=&\sum_{\sigma\in P(r)}\mu(\sigma,1_r)\cdot\#\{\pi\in D_k|\pi\leq\sigma^\gamma,\pi=\gamma(\pi)\}\\
&=&\sum_{\sigma\in P(r)}\mu(\sigma,1_r)\sum_{\pi\leq\sigma^\gamma,\pi=\gamma(\pi)}1
\end{eqnarray*}

In order to exchange the two summations, we first have to replace the summation over $\sigma\in P(r)$ by a summation over $\tau=\sigma^\gamma\in P(k)$. Note that the condition on the latter is exactly $\tau\geq\gamma$ and that we have $\mu(\sigma,1_r)=\mu(\sigma^\gamma,1_k)$. Thus:
$$c_r
=\sum_{\tau\geq\gamma}\mu(\tau,1_k)\sum_{\pi\leq\tau,\pi=\gamma(\pi)}1
=\sum_{\pi=\gamma(\pi)}\sum_{\pi\vee\gamma\leq\tau}\mu(\tau,1_k)$$

The definition of the M\"obius function gives for the second summation:
$$\sum_{\pi\vee\gamma\leq\tau}\mu(\tau,1_k)
=\begin{cases}
1&{\rm if}\ \pi\vee\gamma=1_k\\
0&{\rm otherwise}
\end{cases}$$

With this formula in hand, the assertion follows.

\medskip

(2) In the free case, the proof runs in the same way, by using free cumulants and the corresponding M\"obius function on noncrossing partitions. Note that we have the analogue of our equation in this case only for noncrossing $\sigma$.
\end{proof}

As a first application of the above cumulant formula, we can now recover the theorem of Diaconis and Shahshahani in \cite{dsh} in the orthogonal case, as follows:

\begin{theorem}
The variables $u_k=\lim_{N\to\infty}Tr(u^k)$ are as follows:
\begin{enumerate}
\item For $O_N$, the $u_k$ are real Gaussian variables, with variance $k$
and mean $0$ or $1$, depending on whether $k$ is odd or even. The $u_k$  are independent.

\item For $O_N^+$, at $k=1,2$ we get semicircular variables of  variance $1$
and mean $0$ for $u_1$ and mean $1$ for $u_2$, and at $k\geq 3$ we get circular variables of mean $0$ and covariance $1$. The $u_k$ are $*$-free.
\end{enumerate}
\end{theorem}

\begin{proof}
This follows by using the formula in Proposition 12.11, as follows:

\medskip

(1) In this case $D_k$ consists of all pairings of $k$ elements. We have to count all pairings $\pi$ with the properties that $\pi\vee\gamma=1_k$ and $\pi=\gamma(\pi)$.

\medskip

Note that if $\pi$ connects two different cycles of $\gamma$, say $c_i$ and $c_j$, then the property $\pi=\gamma(p)$ implies that each element from $c_i$ must be paired with an element from $c_j$. Thus those cycles cannot be connected to other cycles and they must contain the same number of elements. This means that for $s\geq 3$ there is no such $\pi$. Thus all cumulants of order 3 and higher vanish asymptotically and all traces are asymptotically Gaussian.

\medskip

Since in the case $s=2$ we only have permissible pairings if the two cycles have the same number of elements, we also see that the covariance between traces of different powers vanishes and thus different powers are asymptotically independent. The variance of $u_k$ is given by the number of matchings between $\{1,\ldots,k\}$ and $\{k+1,\ldots,2k\}$ which are invariant under rotations. 

\medskip

Now since such a matching is determined by the partner of the first element 1, for which we have $k$ possibilities, the variance of $u_k$ is $k$. For the mean, if $k$ is odd there is clearly no pairing at all, and if $k=2p$ is even then the only pairing of $\{1,\ldots,2p\}$ which is invariant under rotations is $(1,p+1),(2,p+2),\ldots,(p,2p)$. Thus the mean of $u_k$ is zero if $k$ is odd and 1 if $k$ is even.

\medskip

(2) In the quantum case $D_k$ consists of noncrossing pairings. We can essentially repeat the arguments from above but have to take care that only noncrossing pairings are counted. We also have to realize that for $k\geq 3$, the $u_k$ are not selfadjoint any longer, thus we have to consider also $u_k^*$ in these cases. This means that in our arguments we have to allow cycles which are rotated ``backwards'' under $\gamma$.

\medskip

By the same reasoning as before we see that free cumulants of order three and higher vanish. The pairing which gave mean 1 in the classical case is only in the case $k=2$ a noncrossing one, thus the mean of $u_2$ is 1, all other means are zero. 

\medskip

For the variances, by using the same argument, one has again that different powers allow no pairings at all and are asymptotically $*$-free. For the matchings between $\{1,\ldots,k\}$ and $\{k+1,\ldots,2k\}$ one has to observe that there is only one non-crossing possibility, namely $(1,2k),(2,2k-1),\ldots,(k,k+1)$ and this satisfies $\pi=\gamma(\pi)$ only if $\gamma$ rotates both cycles in different directions.

\medskip

For $k=1$ and $k=2$ there is no difference between both directions, but for $k\geq 3$ this implies that we get only a non-vanishing covariance between $u_k$ and $u_k^*$, with value 1. This shows that $u_1$ and $u_2$ are semicircular, whereas the higher $u_k$ are circular.
\end{proof}

In order to discuss now permutations and quantum permutations, let us start with the following technical result, further building on the formula in Proposition 12.11:

\begin{proposition}
The cumulants of the variables $u_k=\lim_{N\to\infty}Tr(u^k)$ are as follows, in the quantum permutation group case:
\begin{enumerate}
\item For $S_N$, the classical cumulants are given by:
$$c_r(u_{k_1},\dots ,u_{k_r})=\sum_{q\mid k_i \forall i=1,\dots,r} q^{r-1}$$

\item For $S_N^+$, the free cumulants are given by:
$$c_r(u_{k_1}^{e_1},\dots ,u_{k_r}^{e_r})
=\begin{cases}
2&{\rm if}\ r=1,\, k_1\geq 2\\
2&{\rm if}\ r=2,\, k_1=k_2, \, e_1=e_2^*\\
2&{\rm if}\ r=2,\, k_1=k_2=2\\
1&{\rm otherwise}
\end{cases}$$
\end{enumerate}
\end{proposition}

\begin{proof}
We have two assertions to be proved, the idea being as follows:

\medskip

(1) Here $D_k$ consists of all partitions. We have to count partitions $\pi$ which have the properties that $\pi\vee\gamma=1_k$ and $\pi=\gamma(\pi)$.

\medskip

Consider a partition $\pi$ which connects different cycles of $\gamma$. Consider the restriction of $\pi$ to one cycle. Let $k$ be the number of elements in this cycle and $t$ be the number of the points in the restriction. Then the orbit of those $t$ points under $\gamma$ must give a partition of that cycle, which means that $t$ is a divisor of $k$ and that the $t$ points are equally spaced. The same must be true for all cycles of $\gamma$ which are connected via $\pi$, and the ratio between $t$ and $k$ is the same for all those cycles. 

\medskip

But this means that if one block of $\pi$ connects some cycles then the orbit under $\gamma$ of this block connects exactly those cycles and exhausts all points of those cycles. So if we want to connect all cycles of $\gamma$ then this can only happen in the way that we have one block of $\pi$ intersecting each of the cycles of $\gamma$. 

\medskip

To be more precise, let us consider $c_r(u_{k_1},\dots ,u_{k_r})$. We have then to look for a common divisor $q$ of all $k_1,\dots,k_r$, and a contributing $\pi$ is then one the blocks of which are of the following form: $k_1/q$ points in the first cycle, equally spaced, and so on up to $k_r/q$ points in the last cycle, equally spaced. We can specify this by saying to which points in the other cycles the first point in the first cycle is connected. There are $q^{r-1}$ possibilities for such choices, and this gives the formula in the statement.

\medskip

(2) For $S_N^+$ we have to consider noncrossing partitions instead of all partitions. Most of the partitions from the classical case are crossing, so do not count for the quantum case. Actually, whenever a restriction of a block to one cycle has two or more elements then the corresponding partition is crossing, unless the restriction exhausts the whole group. 

\medskip

This is the case $q=1$ from the considerations above, corresponding to the partition which has only one block, giving a contribution 1 to each cumulant $c_r(u_{k_1},\dots,u_{k_r})$. 

\medskip

For cumulants of order 3 or higher there are no other contributions. For cumulants of second order one might also have contributions coming from pairings, where each restriction of a block to a cycle has one element. 

\medskip

But this is the same problem as in the $O_N^+$ case, and we only get an additional contribution for the second order cumulants $c_2(u_k,u_k^*)$. For first order cumulants, singletons can also appear and make an additional contribution. And this gives the result.
\end{proof}

With the above ingredients in hand, still following \cite{ez1}, we can now formulate a result about permutations and quantum permutations, as follows:

\begin{theorem}
The variables $u_k=\lim_{N\to\infty}Tr(u^k)$ are as follows:
\begin{enumerate}
\item For $S_N$ we have a decomposition of type
$$u_k=\sum_{l\mid k}lC_l$$
with the variables $C_k$ being Poisson of parameter $1/k$, and independent.

\item For $S_N^+$ we have a decomposition of the type
$$u_1=C_1\quad,\quad u_k=C_1+C_k \quad (k\geq 2)$$
where the variables $C_l$ are $*$-free, $C_1$ is free Poisson, $C_2$ is semicircular, and $C_k$ with $k\geq 3$ are circular.
\end{enumerate}
\end{theorem}

\begin{proof}
We have several assertions to be proved, the idea being as follows:

\medskip

(1) Let $C_k$ be the number of cycles of length $k$. We have $u_k=\sum_{l\mid k} l C_l$. We are claiming now that the $C_k$ are independent and each
is a Poisson variable of parameter $1/k$, i.e., that $c_r(C_{l_1},\dots,C_{l_r})$ is zero unless all the $l_i$ are the same, say $=l$, in which case it is $1/l$, independently of $r$. This is compatible with the cumulants for the $u_k$, according to the following computation:
\begin{eqnarray*}
c_r(u_{k_1},\dots ,u_{k_r})
&=&\sum_{l_1|k_1}\ldots\sum_{l_r|k_r}l_1\ldots l_r c_r(C_{l_1},\ldots,C_{l_r})\\
&=&\sum_{l|k_i\forall i}l^r\,\frac{1}{l}
\end{eqnarray*}

Since the $C_k$ are uniquely determined by the $u_k$, via some kind of M\"obius inversion, this shows that the $C_k$ are independent, and that $C_k$ is Poisson with parameter $1/k$.

\medskip

(2) In the classical case the random variable $C_l$ can be defined by:
$$C_l=\frac 1l\sum_{i_1 \dots i_l\ {\rm distinct}}u_{i_1i_2}u_{i_2i_3}\ldots
u_{i_li_1}$$

Note that we divide by $l$ because each term appears actually $l$ times, in cyclically permuted versions, which are all the same because our variables commute.

\medskip

Note that, by using commutativity and the monomial condition, in general the expression $u_{i_1i_2}u_{i_2i_3}\ldots u_{i_ki_1}$ has to be zero unless the indices $(i_1,\ldots,i_k)$ are of the form $(i_1,\ldots,i_l,i_1,\ldots,i_l,\ldots)$, where $l$ divides $k$ and $i_1,\dots,i_l$ are distinct. This yields then the following relation, which we used before to define $C_l$:
\begin{eqnarray*} 
Tr(u^k)
&=&\sum_{i_1\dots i_l}u_{i_1i_2}u_{i_2i_3}\ldots u_{i_li_1}\\
&=&\sum_{l|k}\sum_{i_1\dots i_l\ {\rm distinct}}(u_{i_1i_2}u_{i_2i_3}\ldots u_{i_li_1})^{k/l}\\
&=&\sum_{l|k}lC_l
\end{eqnarray*}

This explicit form of $C_l$ in terms of $u_{ij}$ can be used to give a direct proof, by using the Weingarten formula, of the fact that the $C_l$ are independent and Poisson.

\medskip

(3) In the free case we define the ``cycle'' $C_l$ by requiring neighboring indices to be different, as follows:
$$C_l=\sum_{i_1\not=i_2\not=\ldots\not=i_l\not=i_1}u_{i_1i_2}u_{i_2i_3}\ldots u_{i_li_1}$$

Note that if two adjacent indices are the same in $u_{i_1i_2}u_{i_2i_3}\ldots u_{i_li_1}$ then, because of the relation $u_{ij}u_{ik}=0$ for $j\not=k$, all must be the same or the term vanishes. For the case where all indices are the same we have:
$$\sum_iu_{ii}u_{ii}\ldots u_{ii}
=\sum_iu_{ii}
=C_1$$

But this gives then the following relation:
$$Tr(u^k)=C_k+C_1$$

Again, the $C_l$ are uniquely determined by the $Tr(u^k)$ and thus our calculations also show that the $C_l$ defined by our equation are $*$-free and have the distributions as stated. Thus, we are led to the conclusion in the statement.
\end{proof}

For more on the above, we refer to \cite{ez1} and related papers. As a comment, however, the story is not over here, because the usual characters $\chi=\sum_iu_{ii}$ are generalized by both the truncated characters $\chi_t$ studied in chapter 2, and the Diaconis-Shahshahani $Tr(u^k)$ variables studied here. Thus, there is certainly room for more general results, covering both these computations. Also, in addition to this, we have many interesting questions regarding the diagonal algebras $D(G)\subset C(G)$, previously discussed in this book.

\section*{12c. Invariance questions}

Following \cite{ez2}, let is discuss now probabilistic invariance questions with respect to the basic quantum permutation and rotation groups, namely:
$$\xymatrix@R=15mm@C=15mm{
S_N^+\ar[r]&O_N^+\\
S_N\ar[r]\ar[u]&O_N\ar[u]
}$$

We start by fixing some notations. We use as usual the formalism of the orthogonal quantum groups, which best covers the main quantum groups that we are interested in. In relation wih invariance questions, we first have the following definition:

\begin{definition}
Given a closed subgroup $G\subset O_N^+$, we denote by
$$\alpha:\mathbb C<t_1,\ldots,t_N>\to\mathbb C<t_1,\ldots,t_N>\otimes\,C(G)$$
$$t_i\to\sum_jt_j\otimes v_{ji}$$
the standard coaction of $C(G)$ on the free complex algebra on $N$ variables.
\end{definition}

Observe that the map $\alpha$ constructed above is indeed a coaction, in the sense that it satisfies the following standard coassociativity and counitality conditions:
$$(id\otimes\Delta)\alpha=(\alpha\otimes id)\alpha\quad,\quad 
(id\otimes\varepsilon)\alpha=id$$

With the above notion of coaction in hand, we can now talk about invariant sequences of classical or noncommutative random variables, in the following way:

\begin{definition}
Let $(B,tr)$ be a $C^*$-algebra with a trace, and $x_1,\ldots,x_N\in B$. We say that $x=(x_1,\ldots,x_N)$ is invariant under $G\subset O_N^+$ if the distribution functional 
$$\mu_x:\mathbb C<t_1,\ldots,t_N>\to\mathbb C\quad,\quad 
P\to tr(P(x_1,\ldots,x_N))$$
is invariant under the coaction $\alpha$, in the sense that we have
$$(\mu_x\otimes id)\alpha(P)=\mu_x(P)$$
for any noncommuting polynomial $P\in\mathbb C<t_1,\ldots,t_N>$.
\end{definition}

In the classical case, where $G\subset O_N$ is a usual group, we recover in this way the usual invariance notion from classical probability. In the general case, where $G\subset O_N^+$ is arbitrary, what we have is a natural generalization of this. We have the following equivalent formulation of the above invariance condition: 

\index{invariant sequence}

\begin{proposition}
Let $(B,tr)$ be a $C^*$-algebra with a trace, and $x_1,\ldots,x_N\in B$. Then $x=(x_1,\ldots,x_N)$ is invariant under $G\subset O_N^+$ precisely when
$$tr(x_{i_1}\ldots x_{i_k})=\sum_{j_1\ldots j_k}tr(x_{j_1}\ldots x_{j_k})v_{j_1i_1}\ldots v_{j_ki_k}$$
as an equality in $C(G)$, for any $k\in\mathbb N$, and any $i_1,\ldots,i_k\in\{1,\ldots,N\}$.  
\end{proposition}

\begin{proof}
By linearity, in order for a sequence $x=(x_1,\ldots,x_N)$ to be $G$-invariant in the sense of Definition 12.16, the formula there must be satisfied for any noncommuting monomial $P\in\mathbb C<t_1,\ldots,t_N>$. But an arbitrary such monomial can be written as follows, for a certain $k\in\mathbb N$, and certain indices $i_1,\ldots,i_k\in\{1,\ldots,N\}$:
$$P=t_{i_1}\ldots t_{i_k}$$

Now with this formula for $P$ in hand, we have the following computation:
\begin{eqnarray*}
(\mu_x\otimes id)\alpha(P)
&=&(\mu_x\otimes id)\sum_{j_1,\ldots,j_k}t_{j_1}\ldots t_{j_k}\otimes v_{j_1i_1}\ldots v_{j_ki_k}\\
&=&\sum_{j_1,\ldots,j_k}\mu_x(t_{j_1}\ldots t_{j_k})v_{j_1i_1}\ldots v_{j_ki_k}\\
&=&\sum_{j_1\ldots j_k}tr(x_{j_1}\ldots x_{j_k})v_{j_1i_1}\ldots v_{j_ki_k}
\end{eqnarray*}

On the other hand, by definition of the distribution $\mu_x$, we have:
$$\mu_x(P)
=\mu_x(t_{i_1}\ldots t_{i_k})
=tr(x_{i_1}\ldots x_{i_k})$$

Thus, we are led to the conclusion in the statement.
\end{proof}

As already mentioned after Definition 12.16, in the classical case, where $G\subset O_N$ is a usual compact group, our notion of $G$-invariance coincides with the usual $G$-invariance notion from classical probability. We have in fact the following result:

\begin{proposition}
In the classical group case, $G\subset O_N$, a sequence $(x_1,\ldots,x_N)$ is $G$-invariant in the above sense if and only if
$$tr(x_{i_1}\ldots x_{i_k})=\sum_{j_1\ldots j_k}g_{j_1i_1}\ldots g_{j_ki_k}tr(x_{j_1}\ldots x_{j_k})$$
for any $k\in\mathbb N$, any $i_1,\ldots,i_k\in\{1,\ldots,N\}$, and any $g=(g_{ij})\in G$, and this coincides with the usual notion of $G$-invariance for a sequence of classical random variables.
\end{proposition}

\begin{proof}
According to Proposition 12.17, the invariance property happens precisely when we have the following equality, for any $k\in\mathbb N$, and any $i_1,\ldots,i_k\in\{1,\ldots,N\}$:
$$tr(x_{i_1}\ldots x_{i_k})=\sum_{j_1\ldots j_k}tr(x_{j_1}\ldots x_{j_k})v_{j_1i_1}\ldots v_{j_ki_k}$$

Now by evaluating both sides of this equation at a given $g\in G$, we obtain:
$$tr(x_{i_1}\ldots x_{i_k})=\sum_{j_1\ldots j_k}g_{j_1i_1}\ldots g_{j_ki_k}tr(x_{j_1}\ldots x_{j_k})$$

Thus, we are led to the conclusion in the statement.
\end{proof}

With the above ingredients in hand, we can now investigate invariance questions for the sequences of classical or noncommutative random variables, with respect to the main quantum permutation and rotation groups that we are interested in here. To be more precise, we first have a reverse De Finetti theorem, from \cite{ez2}, as follows:

\index{reverse De Finetti}

\begin{theorem}
Let $(x_1,\ldots,x_N)$ be a sequence in $A$.
\begin{enumerate}
\item If $x_1,\ldots,x_N$ are freely independent and identically distributed with amalgamation over $B$, then the sequence is $S_N^+$-invariant.

\item If $x_1,\ldots,x_N$ are freely independent and identically distributed with amalgamation over $B$, and have centered semicircular distributions with respect to $E$, then the sequence is $O_N^+$-invariant.

\item If $<B,x_1,\ldots,x_N>$ is commutative and $x_1,\ldots,x_N$ are conditionally independent and identically distributed given $B$, then the sequence is $S_N$-invariant.

\item If $<x_1,\ldots,x_N>$ is commutative and $x_1,\ldots,x_N$ are conditionally independent and identically distributed given $B$, and have centered Gaussian distributions with respect to $E$, then the sequence is $O_N$-invariant.
\end{enumerate}
\end{theorem}

\begin{proof}
Assume that the joint distribution of $(x_1,\ldots,x_N)$ satisfies one of the conditions in the statement, and let $D$ be the category of partitions associated to the corresponding easy quantum group.  We have then the following computation:
\begin{eqnarray*}
\sum_{j_1\ldots j_k}tr(x_{j_1}\ldots x_{j_k})v_{j_1i_1}\ldots v_{j_ki_k}
&=&\sum_{j_1\ldots j_k}tr(E(x_{j_1}\ldots x_{j_k}))v_{j_1i_1}\ldots v_{j_ki_k}\\
&=&\sum_{j_1\ldots j_k}\sum_{\pi\leq\ker j}tr(\xi^{(\pi)}_E(x_1,\ldots,x_1))v_{j_1i_1}\ldots v_{j_ki_k}\\
&=&\sum_{\pi\in D(k)}tr(\xi^{(\pi)}_E(x_1,\ldots,x_1))\sum_{\ker j\geq\pi}v_{j_1i_1}\ldots v_{j_ki_k}
\end{eqnarray*}

Here $\xi$ denotes the free and classical cumulants in the cases (1,2) and (3,4) respectively. On the other hand, it follows from a direct computation that if $\pi\in D(k)$ then we have the following formula, in each of the 4 cases in the statement:
$$\sum_{\ker j\geq\pi}v_{j_1i_1}\ldots v_{j_ki_k}=
\begin{cases}1&{\rm if}\ \pi\leq\ker i\\ 
0&{\rm otherwise}
\end{cases}$$

By using this formula, we can finish our computation, in the following way:
\begin{eqnarray*}
\sum_{j_1\ldots j_k}tr(x_{j_1}\ldots x_{j_k})v_{j_1i_1}\ldots v_{j_ki_k}
&=&\sum_{\pi\in D(k)}tr(\xi^{(\pi)}_E(x_1,\ldots,x_1))\delta_{\pi\leq\ker i}\\
&=&\sum_{\pi\leq\ker i}tr(\xi_E^{(\pi)}(x_1,\ldots,x_1))\\
&=&tr(x_{i_1}\ldots x_{i_k})
\end{eqnarray*}

Thus, we are led to the conclusions in the statement.
\end{proof}

\section*{12d. De Finetti theorems}

Our goal here, still following \cite{ez2}, is to establish a converse to Theorem 12.19. Let us begin with some technical results. We will use the following standard fact:

\begin{proposition}
Assume that a sequence $(x_1,\ldots,x_N)$ is $G$-invariant. Then there is a coaction 
$$\widetilde{\alpha}:M_N(\mathbb C)\to M_N(\mathbb C)\otimes C(G)$$
determined by the following formula:
$$\widetilde{\alpha}(p(x))=(ev_x\otimes\pi_N)\alpha(p)$$
Moreover, the fixed point algebra of $\widetilde{\alpha}$ is the $G$-invariant subalgebra $B_N$.
\end{proposition}

\begin{proof}
This follows indeed after identifying the GNS representation of the algebra $\mathbb C<t_1,\ldots,t_N>$ for the state $\mu_x$ with the morphism $ev_x:\mathbb C<t_1,\ldots,t_N>\to M_N(\mathbb C)$.
\end{proof}

In order to further advance, we use the fact that there is a natural conditional expectation given by integrating the coaction $\widetilde{\alpha}$ with respect to the Haar state, as follows:
$$E_N:M_N(\mathbb C)\to B_N\quad,\quad 
E_N(m)=\left(id\otimes\int_G\right)\widetilde{\alpha}(m)$$

The point now is that by using the Weingarten formula, we can give a simple combinatorial formula for the moment functionals with respect to $E_N$, in the case where $G$ is one of the easy quantum groups under consideration. Following \cite{ez2}, we have:

\begin{theorem}
Assume that $(x_1,\ldots,x_N)$ is $G$-invariant, and that either we have $G=O_N^+,S_N^+$, or that $G=O_N,S_N$ and $(x_1,\ldots,x_N)$ commute. We have then
$$E_N^{(\pi)}(b_0x_1b_1,\ldots,x_1b_k)=\frac{1}{N^{|\pi|}}\sum_{\pi\leq\ker i} b_0x_{i_1}\ldots b x_{i_k}b_k$$
for any $\pi$ in the partition category $D(k)$ for $G$, and any $b_0,\ldots,b_k\in B_N$.
\end{theorem}

\begin{proof}
We prove this result by recurrence on the number of blocks of $\pi$. First suppose that $\pi=1_k$ is the partition with only one block. Then: 
$$E_N^{(1_k)}(b_0x_1b_1,\ldots,x_1b_k)
=E_N(b_0x_1\ldots x_1b_k)
=\sum_{i_1 \ldots i_k}b_0x_{i_1}\ldots x_{i_k}b_k\int_Gv_{i_11}\ldots v_{i_k1}$$

Here we have used the fact that the elements $b_0,\dotsc,b_k$ are fixed by the coaction $\widetilde{\alpha}$.  Applying now the Weingarten integration formula, we have:
\begin{eqnarray*}
E_N(b_0x_1\ldots x_1b_k)
&=&\sum_{i_1\ldots i_k}b_0x_{i_1}\ldots x_{i_k}b_k\sum_{\pi\leq\ker i}\sum_\sigma W_{kN}(\pi,\sigma)\\
&=&\sum_{\pi\in D(k)}\left(\sum_{\sigma\in D(k)}W_{kN}(\pi,\sigma)\right) \sum_{\pi\leq\ker i}b_0x_{i_1}\ldots x_{i_k}b_k
\end{eqnarray*}

Now observe that for any $\sigma\in D(k)$ we have the following formula:
$$G_{kN}(\sigma,1_k)=N^{|\sigma\vee1_k|}=N$$

It follows that for any partition $\pi \in D(k)$, we have:
$$N\sum_{\sigma\in D(k)}W_{kN}(\pi,\sigma)
=\sum_{\sigma \in D(k)}W_{kN}(\pi,\sigma)G_{kN}(\sigma,1_k)
=\delta_{\pi1_k}$$

Applying this in the above context, we find, as desired:
\begin{eqnarray*}
E_N(b_0x_1\ldots x_1b_k)
&=&\sum_{\pi\in D(k)}\frac{1}{N}\,\delta_{\pi1_k}\sum_{\pi\leq\ker i}b_0x_{i_1}\ldots x_{i_k}b_k\\
&=&\frac{1}{N}\sum_{i=1}^Nb_0x_i\ldots x_ib_k
\end{eqnarray*}

If the condition (3) or (4) is satisfied, then the general case follows from:
$$E_N^{(\pi)}(b_0x_1b_1,\ldots,x_1b_k)=b_1\ldots b_k\prod_{V\in\pi}E_N(V)(x_1,\ldots,x_1)$$

Indeed, the one thing that we must check here is that if $\pi\in D(k)$ and $V$ is a block of $\pi$ with $s$ elements, then $1_s\in D(s)$.  But this is easily verified, in each case.

\medskip

Assume now that the condition (1) or (2) is satisfied.  Let $\pi\in D(k)$. Since $\pi$ is noncrossing, $\pi$ contains an interval $V=\{l+1,\ldots,l+s+1\}$, and we have:
$$E_N^{(\pi)}(b_0x_1b_1,\ldots,x_1b_k)
=E_N^{(\pi-V)}(b_0x_1b_1,\ldots,E_N(x_1b_{l+1}\ldots x_1b_{l+s})x_1,\ldots,x_1b_k)$$

To apply induction, we must check that we have $\pi-V\in D(k-s)$ and $1_s\in D(s)$.  Indeed, this is easily verified for $NC,NC_2$. Applying induction, we have:
\begin{eqnarray*}
&&E_N^{(\pi)}(b_0x_1b_1,\ldots,x_1b_k)\\
&=&\frac{1}{N^{|\pi|-1}}\sum_{\pi-V\leq\ker i}b_0x_{i_1}\ldots b_l\left(E_n(x_1b_{l+1}\ldots x_1b_{l+s})\right)x_{i_{l+s}}\ldots x_{i_k}b_k\\
&=&\frac{1}{N^{|\pi|-1}}\sum_{\pi-V\leq\ker i}b_0x_{i_1}\ldots b_l\left(\frac{1}{N}\sum_{i=1}^Nx_ib_{l+1}\ldots b x_ib_{l+s}\right)x_{i_{l+s}}\ldots x_{i_k}b_k\\
&=&\frac{1}{N^{|\pi|}}\sum_{\pi\leq\ker i}b_0x_{i_1}\ldots x_{i_k}b_k
\end{eqnarray*}

Thus, we are led to the conclusion in the statement.
\end{proof}

As a continuation, still following \cite{ez2}, we have the following result:

\begin{theorem}
Suppose that $(x_1,\ldots,x_N)$ is $G_N$-invariant, and that $G_N=O_N^+,S_N^+$, or that $G_N=O_N,S_N$ and $(x_1,\ldots,x_N)$ commute. Let $(y_1,\ldots,y_N)$ be a sequence of $B_N$-valued random variables with $B_N$-valued joint distribution determined as follows:
\begin{enumerate}
\item $G=O^+$:  free semicircular, centered with same variance as $x_1$.

\item $G=S^+$: freely independent, $y_i$ has same distribution as $x_1$.

\item $G=O$: independent Gaussian, centered with same variance as $x_1$.

\item $G=S$: independent, $y_i$ has same distribution as $x_1$.
\end{enumerate}
Then if $1\leq j_1,\ldots,j_k \leq N$ and $b_0,\ldots,b_k\in B_N$, we have the following estimate,
$$\left|\left|E_N(b_0x_{j_1}\ldots x_{j_k}b_k)-E(b_0y_{j_1}\ldots y_{j_k}b_k)\right|\right|\leq\frac{C_k(G)}{N}||x_1||^k||b_0||\ldots||b_k||$$
with $C_k(G)$ being a constant depending only on $k$ and $G$.
\end{theorem}

\begin{proof}
First we note that it suffices to prove the result for $N$ large enough. We will assume that $N$ is sufficiently large, as for the Gram matrix $G_{kN}$ to be invertible. Let $1\leq j_1,\ldots,j_k\leq N$ and $b_0,\ldots,b_k\in B_N$. We have then:
\begin{eqnarray*}
E_N(b_0x_{j_1}\ldots x_{j_k}b_k)
&=&\sum_{i_1\ldots i_k}b_0x_{i_1}\ldots x_{i_k}b_k\int v_{i_1j_1}\ldots v_{i_kj_k}\\
&=&\sum_{i_1\ldots i_k}b_0x_{i_1}\ldots x_{i_k}b_k\sum_{\pi\leq\ker i}\sum_{\sigma\leq\ker j}W_{kN}(\pi,\sigma)\\
&=&\sum_{\sigma\leq\ker j}\sum_\pi W_{kN}(\pi,\sigma)\sum_{\pi\leq\ker i} b_0x_{i_1}\ldots x_{i_k}b_k
\end{eqnarray*}

On the other hand, it follows from our assumptions on $(y_1,\ldots,y_N)$, and from the various moment-cumulant formulae given before, that we have:
$$E(b_0y_{j_1}\ldots y_{j_k}b_k)=\sum_{\sigma\leq\ker j}\xi_{E_N}^{(\sigma)}(b_0x_1b_1,\ldots,x_1b_k)$$

Here, and in what follows, $\xi$ denote the relevant free or classical cumulants. The right hand side can be expanded, via the M\"obius inversion formula, in terms of expectation functionals of the following type, with $\pi$ being a partition in $NC,P$ according to the cases (1,2) or (3,4) in the statement, and with $\pi\leq\sigma$ for some $\sigma\in D(k)$: 
$$E_N^{(\pi)}(b_0x_1b_1,\ldots,x_1b_k)$$

Now if $\pi\notin D(k)$, we claim that this expectation functional is zero. Indeed this is only possible if $D= NC_2,P_2$, and if $\pi$ has a block with an odd number of legs. But it is easy to see that in these cases $x_1$ has an even distribution with respect to $E_N$, and therefore we have, as claimed, the following formula:
$$E_N^{(\pi)}(b_0x_1b_1,\ldots,x_1b_k)=0$$

Now this observation allows to to rewrite the above equation as follows:
$$E(b_0y_{j_1}\ldots y_{j_k}b_k)=\sum_{\sigma\leq\ker j}\sum_{\pi\leq \sigma} \mu_{D(k)}(\pi,\sigma)E_N^{(\pi)}(b_0x_1b_1,\ldots,x_1b_k)$$

We therefore obtain the following formula:
$$E(b_0y_{j_1}\ldots y_{j_k}b_k)=\sum_{\sigma\leq\ker j}\sum_{\pi\leq\sigma} \mu_{D(k)}(\pi,\sigma)N^{-|\pi|}\sum_{\pi\leq\ker i}b_0x_{i_1}\ldots x_{i_k}b_k$$

Comparing the above two equations, we find that:
\begin{eqnarray*}
&&E_N(b_0x_{j_1}\ldots x_{j_k}b_k)-E(b_0y_{j_1}\ldots y_{j_k}b_k)\\
&=&\sum_{\sigma\leq\ker j}\sum_\pi\left(W_{kN}(\pi,\sigma)-\mu_{D(k)}(\pi,\sigma)N^{-|\pi|}\right)\sum_{\pi\leq\ker i}b_0x_{i_1}\ldots x_{i_k}b_k
\end{eqnarray*}

Now since $x_1,\ldots,x_N$ are identically distributed with respect to the faithful state $\varphi$, it follows that these variables have the same norm.  Thus, for any $\pi \in D(k)$:
$$\left|\left|\sum_{\pi\leq\ker i}b_0x_{i_1}\ldots x_{i_k}b_k\right|\right|\leq N^{|\pi|}||x_1||^k||b_0||\ldots||b_k||$$

Combining this with the former equation, we obtain the following estimate:
\begin{eqnarray*}
&&\left|\left|E_N(b_0x_{j_1}\ldots x_{j_k}b_k)-E(b_0y_{j_1}\ldots y_{j_k}b_k)\right|\right|\\
&\leq&\sum_{\sigma\leq\ker j}\sum_\pi\left|W_{kN}(\pi,\sigma)N^{|\pi|}-\mu_{D(k)}(\pi,\sigma)\right|||x_1||^k||b_0||\ldots||b_k||
\end{eqnarray*}

Let us set now, according to the above:
$$C_k(G)=\sup_{N\in\mathbb N}\left(N\times\sum_{\sigma,\pi \in D(k)}\left|W_{kN}(\pi,\sigma)N^{|\pi|}-\mu_{D(k)}(\pi,\sigma)\right|\right)$$

But this number is finite by our main estimate, which completes the proof.
\end{proof}

We will use in what follows the inclusions $G_N\subset G_M$ for $N<M$, which correspond to the Hopf algebra morphisms $\omega_{N,M}:C(G_M)\to C(G_N)$ given by:
$$\omega_{N,M}(u_{ij})=
\begin{cases} 
u_{ij}&{\rm if}\ 1\leq i,j\leq N\\
\delta_{ij}&{\rm if}\ \max(i,j)>N\
\end{cases}$$

Still following \cite{ez2}, we begin by extending the notion of $G_N$-invariance to the infinite sequences of variables, in the following way:

\index{invariant sequence}

\begin{definition}
Let $(x_i)_{i\in\mathbb N}$ be a sequence in a noncommutative probability space $(A,\varphi)$. We say that $(x_i)_{i\in\mathbb N}$ is $G$-invariant if 
$$(x_1,\ldots,x_N)$$
is $G_N$-invariant for each $N\in\mathbb N$.
\end{definition}

In other words, the condition is that the joint distribution of $(x_1,\ldots,x_N)$ should be invariant under the following coaction map, for each $N\in\mathbb N$:
$$\alpha_N:\mathbb C<t_1,\ldots,t_N>\to\mathbb C<t_1,\ldots,t_N>\otimes\,C(G_N)$$

It is convenient to extend these coactions to a coaction on the algebra of noncommutative polynomials on an infinite number of variables, in the following way:
$$\beta_N:\mathbb C<t_i|i\in\mathbb N>\to\mathbb C<t_i|i\in\mathbb N>\otimes\,C(G_N)$$

Indeed, we can define $\beta_N$ to be the unique unital morphism satisfying:
$$\beta_N(t_j)=
\begin{cases}
\sum_{i=1}^Nt_i\otimes v_{ij}&{\rm if}\ 1\leq j\leq N\\
t_j\otimes 1&{\rm if}\ j>N
\end{cases}$$

It is clear that $\beta_N$ as constructed above is a coaction of $G_N$. Also, we have the following relations, where $\iota_N:\mathbb C<t_1,\ldots,t_N>\to\mathbb C<t_i|i\in\mathbb N>$ is the natural inclusion:
$$(id\otimes \omega_{N,M})\beta_M=\beta_N\quad,\quad 
(\iota_N\otimes id)\alpha_N=\beta_N\iota_N$$

By using these compatibility relations, we obtain the following result:

\begin{proposition}
An infinite sequence of random variables $(x_i)_{i\in\mathbb N}$ is $G$-invariant if and only if the joint distribution functional 
$$\mu_x:\mathbb C<t_i|i\in\mathbb N>\to\mathbb C\quad,\quad 
P\to tr(P(x))$$
is invariant under the coaction $\beta_N$, for each $N\in\mathbb N$.
\end{proposition}

\begin{proof}
This is clear indeed from the above discussion.
\end{proof}

In what follows $(x_i)_{i\in\mathbb N}$ will be a sequence of self-adjoint random variables in a von Neumann algebra $(M,tr)$. We will assume that $M$ is generated by $(x_i)_{i\in\mathbb N}$. We denote by $L^2(M,tr)$ the corresponding GNS Hilbert space, with inner product as follows:
$$<m_1,m_2>=tr(m_1m_2^*)$$

Also, the strong topology on $M$, that we will use in what follows, will be taken by definition with respect to the faithful representation on the space $L^2(M,tr)$. 

\bigskip

We let $P_N$ be the fixed point algebra of the action $\beta_N$, and we set: 
$$B_N=\left\{p(x)\Big|p\in P_N\right\}''$$

We have then an inclusion $B_{N+1}\subset B_N$, for any $N\geq1$, and we can then define the $G$-invariant subalgebra as the common intersection of these algebras:
$$B=\bigcap_{N\geq1}B_N$$

With these conventions, we have the following result, from \cite{ez2}:

\begin{proposition}
If an infinite sequence of random variables $(x_i)_{i\in\mathbb N}$ is $G$-invariant, then for each $N\in\mathbb N$ there is a coaction 
$$\widetilde{\beta}_N:M\to M\otimes L^\infty(G_N)$$
determined by the following formula, for any $p\in\mathcal P_\infty$:
$$\widetilde{\beta}_N(p(x))=(ev_x\otimes\pi_N)\beta_N(p)$$
The fixed point algebra of $\widetilde{\beta}_N$ is then $B_N$. 
\end{proposition}

\begin{proof}
This is indeed clear from definitions, and from the various compatibility formulae given above, between the coactions $\alpha_N$ and $\beta_N$.
\end{proof}

We have as well the following result, which is clear as well:

\begin{proposition}
In the above context, that of an infinite sequence of random variables belonging to an arbitrary von Neumann algebra $M$ with a trace
$$(x_i)_{i\in\mathbb N}$$
which is $G$-invariant, for each $N\in\mathbb N$ there is a trace-preserving conditional expectation $E_N:M\to B_N$ given by integrating the action $\widetilde{\beta}_N$: 
$$E_N(m)=\left(id\otimes\int_G\right)\widetilde{\beta}_N(m)$$
By taking the limit of these expectations as $N\to \infty$, we obtain a trace-preserving conditional expectation onto the $G$-invariant subalgebra.
\end{proposition}

\begin{proof}
Once again, this is clear from definitions, and from the various compatibility formulae given above, between the coactions $\alpha_N$ and $\beta_N$.
\end{proof}

We are now prepared to state and prove the main theorem, from \cite{ez2}, which comes as a complement to the reverse De Finetti theorem that we already established:

\index{De Finetti theorem}

\begin{theorem}
Let $(x_i)_{i\in\mathbb N}$ be a $G$-invariant sequence of self-adjoint random variables in $(M,tr)$, and assume that $M=<(x_i)_{i\in\mathbb N}>$. Then there exists a subalgebra $B\subset M$ and a trace-preserving conditional expectation $E:M\to B$ such that:
\begin{enumerate}
\item If $G=(S_N)$, then $(x_i)_{i\in\mathbb N}$ are conditionally independent and identically distributed given $B$.

\item If $G=(S_N^+)$, then $(x_i)_{i\in\mathbb N}$ are freely independent and identically distributed with amalgamation over $B$.

\item If $G=(O_N)$, then $(x_i)_{i\in\mathbb N}$ are conditionally independent, and have Gaussian distributions with mean zero and common variance, given $B$.

\item If $G=(O_N^+)$, then $(x_i)_{i\in\mathbb N}$ form a $B$-valued free semicircular family with mean zero and common variance.
\end{enumerate}
\end{theorem}

\begin{proof}
We use the various partial results and formulae established above. Let $j_1,\ldots,j_k \in \mathbb N$ and $b_0,\ldots,b_k\in B$. We have then the following computation:
\begin{eqnarray*}
E(b_0x_{j_1}\ldots x_{j_k}b_k)
&=&\lim_{N\to\infty}E_N(b_0x_{j_1}\ldots x_{j_k}b_k)\\
&=&\lim_{N\to\infty}\sum_{\sigma\leq\ker j}\sum_\pi W_{kN}(\pi,\sigma) \sum_{\pi\leq\ker i}b_0x_{i_1}\ldots x_{i_k}b_k\\
&=&\lim_{N\to\infty}\sum_{\sigma\leq\ker j}\sum_{\pi\leq\sigma}\mu_{D(k)}(\pi,\sigma)N^{-|\pi|}\sum_{\pi\leq\ker i}b_0x_{i_1}\ldots x_{i_k}b_k
\end{eqnarray*}

Let us recall now from the above that we have the following compatibility formula, where $\widetilde{\iota}_N:W^*(x_1,\ldots,x_N)\to M$ is the canonical inclusion, and $\widetilde{\alpha}_N$ is as before:
$$(\widetilde{\iota}_N\otimes id)\widetilde{\alpha}_N=\widetilde{\beta}_N\widetilde{\iota}_N$$

By using this formula, and the above cumulant results, we have:
$$E(b_0x_{j_1}\ldots x_{j_k}b_k)=\lim_{N\to\infty}\sum_{\sigma\leq\ker j} \sum_{\pi\leq\sigma}\mu_{D(k)}(\pi,\sigma)E_N^{(\pi)}(b_0x_1b_1,\ldots,x_1b_k)$$

We therefore obtain the following formula:
$$E(b_0x_{j_1}\ldots x_{j_k}b_k)=\sum_{\sigma\leq\ker j}\sum_{\pi\leq\sigma} \mu_{D(k)}(\pi,\sigma)E^{(\pi)}(b_0x_1b_1,\ldots,x_1b_k)$$

We can replace the sum of expectation functionals by cumulants, as to obtain:
$$E(b_0x_{j_1}\ldots x_{j_k}b_k)=\sum_{\sigma\leq\ker j}\xi_E^{(\sigma)}(b_0x_1b_1,\ldots,x_1b_k)$$

Here and in what follows $\xi$ denotes as usual the relevant free or classical cumulants, depending on the quantum group that we are dealing with, free or classical. Now since the cumulants are determined by the moment-cumulant formulae, we conclude that we have the following formula:
$$\xi_E^{(\sigma)}(b_0x_{j_1}b_1,\ldots,x_{j_k}b_k)
=\begin{cases}
\xi_E^{(\sigma)}(b_0x_1b_1,\ldots,x_1b_k)&{\rm if}\ \sigma\in D(k)\ {\rm and}\ \sigma\leq\ker j\\
0&{\rm otherwise}
\end{cases}$$

With this formula in hand, the result then follows from the characterizations of these joint distributions in terms of cumulants, and we are done.
\end{proof}

Summarizing, we are done with our first and main objective, namely establishing De Finetti theorems for the main quantum permutation and rotation groups, namely:
$$\xymatrix@R=15mm@C=15mm{
S_N^+\ar[r]&O_N^+\\
S_N\ar[r]\ar[u]&O_N\ar[u]
}$$

The story is of course not over here, and there are many related interesting questions left, which are more technical, in relation with the invariance questions with respect to these quantum groups. We refer here to \cite{ez1}, \cite{ez2}, \cite{cu1}, \cite{csp}, \cite{ksp} and related papers.

\section*{12e. Exercises} 

Things have been quite technical in this chapter, dealing with advanced probability theory, and so will be our exercises here. As a first exercise, we have:

\begin{exercise}
Formulate and prove the classical De Finetti theorem, concerning sequences which are invariant under $S_\infty$, without using representation theory methods.
\end{exercise}

This is something very standard, and is a must-do exercise, the point being that all the Weingarten technology used in this chapter, which is something quite heavy, was motivated by the fact that we want to deal with several quantum groups at the same time, in a ``uniform'' way. In the case of the symmetric group itself things are in fact much simpler, and the exercise is about understanding how this works.

\begin{exercise}
Formulate and prove the free De Finetti theorem, concerning sequences which are invariant under $(S_N^+)$, without using representation theory methods.
\end{exercise}

The same comments as for the previous exercise apply, the idea being that, once again, the Weingarten function machinery can be avoided in this case.

\begin{exercise}
Work out the full proof of the explicit formula for the Weingarten function for $S_N$, namely
$$W_{kN}(\pi,\nu)=\sum_{\tau\leq\pi\wedge\nu}\mu(\tau,\pi)\mu(\tau,\nu)\frac{(N-|\tau|)!}{N!}$$
then of the main estimate for this function, namely
$$W_{kN}(\pi,\nu)=N^{-|\pi\wedge\nu|}(
\mu(\pi\wedge\nu,\pi)\mu(\pi\wedge\nu,\nu)+O(N^{-1}))$$
where $\mu$ is the M\"obius function of $P(k)$.
\end{exercise}

This was something that was already discussed in the above, the idea being that all this comes from the explicit knowledge of the integrals over $S_N$, via the M\"obius inversion formula, and the problem now is that of working out all the details.

\begin{exercise}
Work out estimates for the integrals of type
$$\int_{S_N^+}v_{i_1j_1}v_{i_2j_2}v_{i_3j_3}v_{i_4j_4}$$
and then for the Weingarten function of $S_N^+$ at $k=4$.
\end{exercise}

Once again, this was something partly discussed in the above, with the comment that things are clear at $k=2,3$, due to the formula $P(k)=NC(k)$ valid here. The problem now is that of working out what happens at $k=4$, where things are non-trivial.

\begin{exercise}
Prove directly that the function
$$d(\pi,\nu)=\frac{|\pi|+|\nu|}{2}-|\pi\vee\nu|$$
is a distance on $P(k)$.
\end{exercise}

To be more precise here, this is something that we talked about in the above, with the idea being that this follows from a number of well-known facts regarding the partitions in $P(k)$. The problem now is that of proving directly this result.

\part{Matrix models}

\ \vskip50mm

\begin{center}
{\em And I wait praying to the Northern Star

I'm afraid it won't lead you anywhere

He's so cold raining on the world tonight

All the angels kneeling to the Northern Lights}
\end{center}

\chapter{Matrix models}

\section*{13a. Matrix models}

One potentially interesting method for the study of the closed subgroups $G\subset S_N^+$, that we have not tried yet, consists in modeling the standard coordinates $u_{ij}\in C(G)$ by concrete variables over some familiar $C^*$-algebra, $U_{ij}\in B$. Indeed, assuming that the model is faithful in some suitable sense, and that the variables $U_{ij}$ are not too complicated, all questions about $G$ would correspond in this way to routine questions inside $B$. We will discuss here such questions, which are quite interesting, first for the arbitrary closed subgroups $G\subset U_N^+$, and then for the quantum permutation groups $G\subset S_N^+$.

\bigskip

All this sounds good, mathematically speaking, and we will soon see that there are some potentially interesting connections with physics as well. Getting started now, we have a good idea, but we must first solve the following philosophical question:

\begin{question}
What type of target algebras $B$ shall we use for our matrix models $\pi:C(G)\to B$? We would like these to be simple enough, as for the computations inside them to be doable, but also general enough, as to model well our quantum groups.
\end{question}

In answer, a good idea would be probably that of using random matrix algebras, $B=M_K(C(T))$, with $K\geq1$ being an integer, and $T$ being a compact space. Indeed, these algebras generalize the most familiar algebras that we know, namely the matrix ones $M_K(\mathbb C)$, and the commutative ones $C(T)$, so they are definitely simple enough. As for their potential modeling power, my cat who knows some physics says okay. 

\bigskip

In short, time to start our study, with the following definition:

\index{matrix model}
\index{random matrix model}

\begin{definition}
A matrix model for $G\subset U_N^+$ is a morphism of $C^*$-algebras
$$\pi:C(G)\to M_K(C(T))$$
where $K\geq1$ is an integer, and $T$ is a compact space.
\end{definition}

As a first comment, focusing on such models might look a bit restrictive, but we will soon discover that, with some know-how, we can do many things with such models. For the moment, let us develop some general theory. The main question to be solved is that of understanding the suitable faithfulness assumptions needed on $\pi$, as for the model to ``remind'' the quantum group. As we will see, this is something quite tricky.

\bigskip

The simplest situation is when $\pi$ is faithful in the usual sense. Here $\pi$ obviously reminds $G$. However, this is something quite restrictive, because in this case the algebra $C(G)$ must be quite small, admitting an embedding as follows:
$$\pi:C(G)\subset M_K(C(T))$$

Technically, this means that $C(G)$ must be of type I, as an operator algebra, and we will discuss this in a moment, with the comment that this is indeed something quite restrictive. However, there are many interesting examples here, and all this is worth a detailed look. First, we have the following result, providing us with basic examples:

\begin{proposition}
The following closed subgroups $G\subset U_N^+$ have faithful models:
\begin{enumerate}
\item The compact Lie groups $G\subset U_N$.

\item The finite quantum groups $G\subset U_N^+$.
\end{enumerate}
In both cases, we can arrange for $\int_G$ to be restriction of the random matrix trace.
\end{proposition}

\begin{proof}
These assertions are all elementary, the proofs being as follows:

\medskip

(1) This is clear, because we can simply use here the identity map:
$$id:C(G)\to M_1(C(G))$$

(2) Here we can use the left regular representation $\lambda:C(G)\to M_{|G|}(\mathbb C)$. Indeed, let us endow the linear space $H=C(G)$ with the scalar product $<a,b>=\int_Gab^*$. We have then a representation of $*$-algebras, as follows:
$$\lambda:C(G)\to B(H)\quad,\quad 
a\to[b\to ab]$$

Now since we have $H\simeq\mathbb C^{|G|}$, we can view $\lambda$ as a matrix model map, as above. 

\medskip

(3) Finally, our claim is that we can choose our model as for the following formula to hold, where $\int_T$ is the integration with respect to a given probability measure on $T$:
$$\int_G=\left(tr\otimes\int_T\right)\pi$$

But this is clear for the model in (1), by definition, and is clear as well for the model in (2), by using the basic properties of the left regular representation.
\end{proof}

In the above result, the last assertion is quite interesting, and suggests formulating the following definition, somewhat independently on the notion of faithfulness:

\index{stationary model}

\begin{definition}
A matrix model $\pi:C(G)\to M_K(C(T))$ is called stationary when
$$\int_G=\left(tr\otimes\int_T\right)\pi$$
where $\int_T$ is the integration with respect to a given probability measure on $T$.
\end{definition}

Here the term ``stationary'' comes from a functional analytic interpretation of all this, with a certain Ces\`aro limit needed to be stationary, and this will be explained later. Yet another explanation comes from a certain relation with the lattice models, but this is something rather folklore, not axiomatized yet. We will be back to this.

\bigskip

We will see in a moment that stationarity implies faithfulness, so that stationarity can be regarded as being a useful, pragmatic version of faithfulness. But let us first discuss the examples. Besides those in Proposition 13.3, we can look at group duals. So, consider a discrete group $\Gamma$, and a model for the corresponding group algebra, as follows:
$$\pi:C^*(\Gamma)\to M_K(C(T))$$

Since a representation of a group algebra must come from a unitary representation of the group, such a matrix model must come from a representation as follows:
$$\rho:\Gamma\to C(T,U_K)$$

With this identification made, we have the following result:

\begin{proposition}
An matrix model $\rho:\Gamma\subset C(T,U_K)$ is stationary when:
$$\int_Ttr(g^x)dx=0,\forall g\neq1$$
Moreover, the examples include all abelian groups, and all finite groups.
\end{proposition}

\begin{proof}
Consider indeed a group embedding $\rho:\Gamma\subset C(T,U_K)$, which produces by linearity a matrix model, as follows:
$$\pi:C^*(\Gamma)\to M_K(C(T))$$

It is enough to formulate the stationarity condition on the group elements $g\in C^*(\Gamma)$. Let us set $\rho(g)=(x\to g^x)$. With this notation, the stationarity condition reads:
$$\int_Ttr(g^x)dx=\delta_{g,1}$$

Since this equality is trivially satisfied at $g=1$, where by unitality of our representation we must have $g^x=1$ for any $x\in T$, we are led to the condition in the statement. Regarding now the examples, these are both clear. More precisely:

\medskip

(1) When $\Gamma$ is abelian we can use the following trivial embedding:
$$\Gamma\subset C(\widehat{\Gamma},U_1)\quad,\quad 
g\to[\chi\to\chi(g)]$$

(2) When $\Gamma$ is finite we can use the left regular representation:
$$\Gamma\subset\mathcal L(\mathbb C\Gamma)\quad,\quad 
g\to[h\to gh]$$

Indeed, in both cases, the stationarity condition is trivially satisfied.
\end{proof}

In order to discuss now certain analytic aspects of the matrix models, let us go back to the von Neumann algebras, discussed in chapter 1, and in chapter 3. We recall from there that we have the following result, due to Murray-von Neumann and Connes:

\index{reduction theory}
\index{factor}
\index{Connes classification}

\begin{theorem}
Given a von Neumann algebra $A\subset B(H)$, if we write its center as 
$$Z(A)=L^\infty(X)$$
then we have a decomposition as follows, with the fibers $A_x$ having trivial center:
$$A=\int_XA_x\,dx$$
Moreover, the factors, $Z(A)=\mathbb C$, can be basically classified in terms of the ${\rm II}_1$ factors, which are those satisfying $\dim A=\infty$, and having a faithful trace $tr:A\to\mathbb C$.
\end{theorem}

\begin{proof}
This is something which is clear in finite dimensions, and in the commutative case too. In general, this is something heavy, the idea being as follows:

\medskip

(1) The first assertion, regarding the decomposition into factors, is von Neumann's reduction theory main result, which is actually one of the heaviest results in fundamental mathematics, and whose proof uses advanced functional analysis techniques. 

\medskip

(2) The classification of factors, due to Murray-von Neumann and Connes, is again something heavy, the idea being that the ${\rm II}_1$ factors are the ``building blocks'', with the other factors basically appearing from them via crossed product type constructions.
\end{proof}

Back now to matrix models, as a first general result, which is something which is not exactly trivial, and whose proof requires some functional analysis, we have:

\index{coamenability}
\index{type I algebra}
\index{stationarity}

\begin{theorem}
Assuming that a closed subgroup $G\subset U_N^+$ has a stationary model
$$\pi:C(G)\to M_K(C(T))$$
it follows that $G$ must be coamenable, and that the model is faithful. Moreover, $\pi$ extends into an embedding of von Neumann algebras, as follows,
$$L^\infty(G)\subset M_K(L^\infty(T))$$
which commutes with the canonical integration functionals.
\end{theorem}

\begin{proof}
Assume that we have a stationary model, as in the statement. By performing the GNS construction with respect to $\int_G$, we obtain a factorization as follows, which commutes with the respective canonical integration functionals:
$$\pi:C(G)\to C(G)_{red}\subset M_K(C(T))$$

Thus, in what regards the coamenability question, we can assume that $\pi$ is faithful. With this assumption made, we have an embedding as follows:
$$C(G)\subset M_K(C(T))$$

By performing the GNS construction we obtain a better embedding, as follows:
$$L^\infty(G)\subset M_K(L^\infty(T))$$

Now since the von Neumann algebra on the right is of type I, so must be its subalgebra $A=L^\infty(G)$. But this means that, when writing the center of this latter algebra as  $Z(A)=L^\infty(X)$, the whole algebra decomposes over $X$, as an integral of type I factors:
$$L^\infty(G)=\int_XM_{K_x}(\mathbb C)\,dx$$

In particular, we can see from this that $C(G)\subset L^\infty(G)$ has a unique $C^*$-norm, and so $G$ is coamenable. Thus we have proved our first assertion, and the second assertion follows as well, because our factorization of $\pi$ consists of the identity, and of an inclusion.
\end{proof}

In relation with the above, we have the following well-known result of Thoma:

\begin{theorem}
For a discrete group $\Gamma$, the following are equivalent:
\begin{enumerate}
\item $C^*(\Gamma)$ is of type I, so that we have an embedding $\pi:C^*(\Gamma)\subset M_K(C(X))$, with $X$ being a compact space.

\item $C^*(\Gamma)$ has a stationary model of type $\pi:C^*(\Gamma)\to M_F(C(L))$, with $F$ being a finite group, and $L$ being a compact abelian group.

\item $\Gamma$ is virtually abelian, in the sense that we have an abelian subgroup $\Lambda\triangleleft\Gamma$ such that the quotient group $F=\Gamma/\Lambda$ is finite. 

\item $\Gamma$ has an abelian subgroup $\Lambda\subset\Gamma$ whose index $K=[\Gamma:\Lambda]$ is finite.
\end{enumerate}
\end{theorem}

\begin{proof}
There are several proofs for this fact, the idea being as follows:

\medskip

$(1)\implies(4)$ This is the non-trivial implication, that we will not prove here. We refer instead to the literature, either Thoma's orignal paper, or books like those of Dixmier, mixing advanced group theory and advanced operator algebra theory.

\medskip

$(4)\implies(3)$ We choose coset representatives $g_i\in\Gamma$, and we set:
$$\Lambda'=\bigcap_ig_i\Gamma g_i^{-1}$$

Then $\Lambda'\subset\Lambda$ has finite index, and we have $\Lambda'\triangleleft\Gamma$, as desired.

\medskip

$(3)\implies(2)$ This follows by using the theory of induced representations. We can define a model $\pi:C^*(\Gamma)\to M_F(C(\widehat{\Lambda}))$ by setting:
$$\pi(g)(\chi)=Ind_\Lambda^\Gamma(\chi)(g)$$

Indeed, any character $\chi\in\widehat{\Lambda}$ is a 1-dimensional representation of $\Lambda$, and we can therefore consider the induced representation $Ind_\Lambda^\Gamma(\chi)$ of the group $\Gamma$. This representation is $|F|$-dimensional, and so maps the group elements $g\in\Gamma$ into order $|F|$ matrices $Ind_\Lambda^\Gamma(\chi)(g)$. Thus the above map $\pi$ is well-defined, and the fact that it is a representation is clear as well. In order to check now the stationarity property of this representation, we can use the following well-known character formula, due to Frobenius:
$$Tr\left(Ind_\Lambda^\Gamma(\chi)(g)\right)=\sum_{x\in F}\delta_{x^{-1}gx\in\Lambda}\chi(x^{-1}gx)$$

By integrating with respect to $\chi\in\widehat{\Lambda}$, we deduce from this that we have:
\begin{eqnarray*}
\left(Tr\otimes \int_{\widehat{\Lambda}}\right)\pi(g)
&=&\sum_{x\in F}\delta_{x^{-1}gx\in \Lambda}\int_{\widehat{\Lambda}}\chi(x^{-1}gx)d\chi\\
&=&\sum_{x\in F}\delta_{x^{-1}gx\in \Lambda}\delta_{g,1}\\
&=&|F|\cdot\delta_{g,1}
\end{eqnarray*}

Now by dividing by $|F|$ we conclude that the model is stationary, as claimed.

\medskip

$(2)\implies(1)$ This is the trivial implication, with the faithfulness of $\pi$ following from the abstract functional analysis arguments from the proof of Theorem 13.7.
\end{proof}

We refer to \cite{bch}, \cite{bfr} and related papers for more on all this, including for some partial extensions of Thoma's theorem, to the case of the discrete quantum groups.

\bigskip

Getting back now to Definition 13.2, more generally, we can model in that way the standard coordinates $x_i\in C(X)$ of various algebraic manifolds $X\subset S^{N-1}_{\mathbb C,+}$. Indeed, these manifolds generalize the compact matrix quantum groups, which appear as:
$$G\subset U_N^+\subset S^{N^2-1}_{\mathbb C,+}$$

Thus, we have many other interesting examples of such manifolds, such as the homogeneous spaces discussed in chapter 8. However, at this level of generality, not much general theory is available. It is elementary to show that, under the technical assumption $X^{class}\neq\emptyset$, there exists a universal $K\times K$ model for the algebra $C(X)$, which factorizes as follows, with $X^{(K)}\subset X$ being a certain algebraic submanifold: 
$$\pi_K:C(X)\to C(X^{(K)})\subset M_K(C(T_K))$$

To be more precise, the universal $K\times K$ model space $T_K$ appears by imposing to the complex $K\times K$ matrices the relations defining $X$, and the algebra $C(X^{(K)})$ is then by definition the image of $\pi_K$. In relation with this, we can set as well:
$$X^{(\infty)}=\bigcup_{K\in\mathbb N}X^{(K)}$$

We are led in this way to a filtration of $X$, as follows:
$$X^{class}= X^{(1)}\subset X^{(2)}\subset X^{(3)}\subset\ldots\ldots\subset X^{(\infty)}\subset X$$

It is possible to say a few non-trivial things about these manifolds $X^{(K)}$. In the compact quantum group case, however, that we are mainly interested in here, the matrix truncations $G^{(K)}\subset G$ are generically not quantum subgroups at $K\geq2$, and so this theory is a priori not very useful, at least in its basic form presented above.

\section*{13b. Inner faithfulness}

Let us discuss now the general, non-coamenable case, with the aim of finding a weaker notion of faithfulness, which still does the job, namely that of ``reminding'' the quantum group. The idea comes by looking at the group duals $G=\widehat{\Gamma}$. Consider indeed a general model for the associated group algebra, which can be written as follows:
$$\pi:C^*(\Gamma)\to M_K(C(T))$$

The point is that such a representation of the group algebra must come by linearization from a unitary group representation, as follows: 
$$\rho:\Gamma\to C(T,U_K)$$

Now observe that when this group representation $\rho$ is faithful, the representation $\pi$ is in general not faithful, for instance because when $T=\{.\}$ its target algebra is finite dimensional. On the other hand, this representation ``reminds'' $\Gamma$, so can be used in order to fully understand $\Gamma$. Thus, we have an idea here, basically saying that, for practical purposes, the faithfuless property can be replaced with something much weaker. 

\bigskip

This weaker notion, which will be of great interest for us, is called ``inner faithfulness''. The general theory here, from \cite{bb3}, starts with the following definition:

\index{Hopf image}
\index{inner faithfulness}

\begin{definition}
Let $\pi:C(G)\to M_K(C(T))$ be a matrix model. 
\begin{enumerate}
\item The Hopf image of $\pi$ is the smallest quotient Hopf $C^*$-algebra $C(G)\to C(H)$ producing a factorization as follows:
$$\pi:C(G)\to C(H)\to M_K(C(T))$$

\item When the inclusion $H\subset G$ is an isomorphism, i.e. when there is no non-trivial factorization as above, we say that $\pi$ is inner faithful.
\end{enumerate}
\end{definition}

The above notions are quite tricky, and having them well understood will take us some time. As a first example, motivated by the above discussion, in the case where $G=\widehat{\Gamma}$ is a group dual, $\pi$ must come from a group representation, as follows:
$$\rho:\Gamma\to C(T,U_K)$$

Thus the minimal factorization in (1) is obtained by taking the image:
$$\rho:\Gamma\to\Lambda\subset C(T,U_K)$$

Thus, as a conclusion, in this case $\pi$ is inner faithful precisely when we have:
$$\Gamma\subset C(T,U_K)$$

Dually now, given a compact Lie group $G$, and elements $g_1,\ldots,g_K\in G$, we have a diagonal representation $\pi:C(G)\to M_K(\mathbb C)$, appearing as follows:
$$f\to\begin{pmatrix}
f(g_1)\\
&\ddots\\
&&f(g_K)
\end{pmatrix}$$

The minimal factorization of this representation $\pi$, as in Definition 13.9 (1), is then via the algebra $C(H)$, with $H$ being the following closed subgroup of $G$:
$$H=\overline{<g_1,\ldots,g_K>}$$

Thus, as a conclusion, $\pi$ is inner faithful precisely when we have:
$$G=H$$

There are many other examples of inner faithful representations, which are however substantially more technically advanced, and we will discuss them later. 

\bigskip

Back to general theory now, in the framework of Definition 13.9, the existence and uniqueness of the Hopf image come by dividing $C(G)$ by a suitable ideal, with this being something standard. Alternatively, in Tannakian terms, as explained in \cite{bb3}, we have:

\index{Tannakian duality}

\begin{theorem}
Assuming $G\subset U_N^+$, with fundamental corepresentation $u=(u_{ij})$, the Hopf image of a model $\pi:C(G)\to M_K(C(T))$ comes from the Tannakian category
$$C_{kl}=Hom(U^{\otimes k},U^{\otimes l})$$
where $U_{ij}=\pi(u_{ij})$, and where the spaces on the right are taken in a formal sense.
\end{theorem}

\begin{proof}
Since the morphisms increase the intertwining spaces, when defined either in a representation theory sense, or just formally, we have inclusions as follows:
$$Hom(u^{\otimes k},u^{\otimes l})\subset Hom(U^{\otimes k},U^{\otimes l})$$

More generally, we have such inclusions when replacing $(G,u)$ with any pair producing a factorization of $\pi$. Thus, by Tannakian duality, the Hopf image must be given by the fact that the intertwining spaces must be the biggest, subject to the above inclusions. On the other hand, since $u$ is biunitary, so is $U$, and it follows that the spaces on the right form a Tannakian category. Thus, we have a quantum group $(H,v)$ given by:
$$Hom(v^{\otimes k},v^{\otimes l})=Hom(U^{\otimes k},U^{\otimes l})$$

By the above discussion, $C(H)$ follows to be the Hopf image of $\pi$, as claimed.
\end{proof}

Regarding now the study of the inner faithful models, a key problem is that of computing the Haar integration functional. The result here, from \cite{wa3}, is as follows:

\index{truncated integrals}

\begin{theorem}
Given an inner faithful model $\pi:C(G)\to M_K(C(T))$, we have
$$\int_G=\lim_{k\to\infty}\frac{1}{k}\sum_{r=1}^k\int_G^r$$
with the truncations of the integration on the right being given by
$$\int_G^r=(\varphi\circ\pi)^{*r}$$
with $\phi*\psi=(\phi\otimes\psi)\Delta$, and with $\varphi=tr\otimes\int_T$ being the random matrix trace.
\end{theorem}

\begin{proof}
This is something quite tricky, the idea being as follows:

\medskip

(1) As a first observation, there is an obvious similarity here with the Woronowicz construction of the Haar measure, explained in chapter 1. In fact, the above result holds more generally for any model $\pi:C(G)\to B$, with $\varphi\in B^*$ being a faithful trace.

\medskip

(2) In order to prove now the result, we can proceed as in chapter 1. If we denote by $\int_G'$ the limit in the statement, we must prove that this limit converges, and that:
$$\int_G'=\int_G$$

It is enough to check this on the coefficients of the Peter-Weyl corepresentations, and if we let $v=u^{\otimes k}$ be one of these corepresentations, we must prove that we have:
$$\left(id\otimes\int_G'\right)v=\left(id\otimes\int_G\right)v$$

(3) In order to prove this, we already know, from the Haar measure theory from chapter 1, that the matrix on the right is the orthogonal projection onto $Fix(v)$:
$$\left(id\otimes\int_G\right)v=Proj\Big[Fix(v)\Big]$$

Regarding now the matrix on the left, the trick in \cite{wo1} applied to the linear form $\varphi\pi$ tells us that this is the orthogonal projection onto the $1$-eigenspace of $(id\otimes\varphi\pi)v$:
$$\left(id\otimes\int_G'\right)v=Proj\Big[1\in (id\otimes\varphi\pi)v\Big]$$

(4) Now observe that, if we set $V_{ij}=\pi(v_{ij})$, we have the following formula:
$$(id\otimes\varphi\pi)v=(id\otimes\varphi)V$$

Thus, we can apply the trick in \cite{wo1}, and we conclude that the $1$-eigenspace that we are interested in equals $Fix(V)$. But, according to Theorem 13.10, we have:
$$Fix(V)=Fix(v)$$

Thus, we have proved that we have $\int_G'=\int_G$, as desired.
\end{proof}

In practice, Theorem 13.11 is something quite powerful. As an illustration, regarding the law of the main character, we obtain here the following result:

\begin{proposition}
Assume that $\pi:C(G)\to M_K(C(T))$ is inner faithful, let
$$\mu=law(\chi)$$
and let $\mu^r$ be the law of $\chi$ with respect to $\int_G^r=(\varphi\circ\pi)^{*r}$, where $\varphi=tr\otimes\int_T$.
\begin{enumerate}
\item We have the following convergence formula, in moments: 
$$\mu=\lim_{k\to\infty}\frac{1}{k}\sum_{r=0}^k\mu^r$$

\item The moments of $\mu^r$ are the numbers $c_\varepsilon^r=Tr(T_\varepsilon^r)$, where:
$$(T_\varepsilon)_{i_1\ldots i_p,j_1\ldots j_p}=\left(tr\otimes\int_T\right)(U_{i_1j_1}^{\varepsilon_1}\ldots U_{i_pj_p}^{\varepsilon_p})$$
\end{enumerate}
\end{proposition}

\begin{proof}
These formulae are both elementary, by using the convergence result established in Theorem 13.11, the proof being as follows:

\medskip

(1) This follows from the limiting formula in Theorem 13.11, by applying the linear forms there to the main character $\chi$.

\medskip

(2) This follows from the definitions of the measure $\mu^r$ and of the matrix $T_e$, by summing the entries of $T_e$ over equal indices, $i_r=j_r$.
\end{proof}

Interestingly, the above results regarding inner faithfulness have applications as well to the notion of stationarity introduced before, clarifying among others the use of the word ``stationary''. To be more precise, in order to detect the stationary models, we have the following useful criterion, mixing linear algebra and analysis, from \cite{bb3}:

\index{idempotent state}
\index{stationary on its image}

\begin{theorem}
For a model $\pi:C(G)\to M_K(C(T))$, the following are equivalent:
\begin{enumerate}
\item $Im(\pi)$ is a Hopf algebra, and the Haar integration on it is:
$$\psi=\left(tr\otimes\int_T\right)\pi$$

\item The linear form $\psi=(tr\otimes\int_T)\pi$ satisfies the idempotent state property:
$$\psi*\psi=\psi$$

\item We have $T_e^2=T_e$, $\forall p\in\mathbb N$, $\forall e\in\{1,*\}^p$, where:
$$(T_e)_{i_1\ldots i_p,j_1\ldots j_p}=\left(tr\otimes\int_T\right)(U_{i_1j_1}^{e_1}\ldots U_{i_pj_p}^{e_p})$$
\end{enumerate}
If these conditions are satisfied, we say that $\pi$ is stationary on its image.
\end{theorem}

\begin{proof}
Given a matrix model $\pi:C(G)\to M_K(C(T))$ as in the statement, we can factorize it via its Hopf image, as in Definition 13.9:
$$\pi:C(G)\to C(H)\to M_K(C(T))$$

Now observe that (1,2,3) above depend only on the factorized representation:
$$\nu:C(H)\to M_K(C(T))$$

Thus, we can assume in practice that we have $G=H$, which means that we can assume that $\pi$ is inner faithful. With this assumption made, the formula in Theorem 13.11 applies to our situation, and the proof of the equivalences goes as follows:

\medskip

$(1)\implies(2)$ This is clear from definitions, because the Haar integration on any compact quantum group satisfies the idempotent state equation:
$$\psi*\psi=\psi$$

$(2)\implies(1)$ Assuming $\psi*\psi=\psi$, we have $\psi^{*r}=\psi$ for any $r\in\mathbb N$, and Theorem 13.11 gives $\int_G=\psi$. By using now Theorem 13.7, we obtain the result.

\medskip

In order to establish now $(2)\Longleftrightarrow(3)$, we use the following elementary formula, which comes from the definition of the convolution operation:
$$\psi^{*r}(u_{i_1j_1}^{e_1}\ldots u_{i_pj_p}^{e_p})=(T_e^r)_{i_1\ldots i_p,j_1\ldots j_p}$$

\medskip

$(2)\implies(3)$ Assuming $\psi*\psi=\psi$, by using the above formula at $r=1,2$ we obtain that the matrices $T_e$ and $T_e^2$ have the same coefficients, and so they are equal.

\medskip

$(3)\implies(2)$ Assuming $T_e^2=T_e$, by using the above formula at $r=1,2$ we obtain that the linear forms $\psi$ and $\psi*\psi$ coincide on any product of coefficients $u_{i_1j_1}^{e_1}\ldots u_{i_pj_p}^{e_p}$. Now since these coefficients span a dense subalgebra of $C(G)$, this gives the result.
\end{proof}

\section*{13c. Half-liberation}

As a first illustration, we can apply the above criterion to certain models for $O_N^*,U_N^*$. We first have the following result, coming from the work in \cite{bb3}, \cite{ez2}, \cite{bdu}:

\index{half-liberation}

\begin{proposition}
We have a matrix model as follows, 
$$C(O_N^*)\to M_2(C(U_N))\quad,\quad 
u_{ij}\to\begin{pmatrix}0&v_{ij}\\ \bar{v}_{ij}&0\end{pmatrix}$$
where $v$ is the fundamental corepresentation of $C(U_N)$, as well as a model as follows,
$$C(U_N^*)\to M_2(C(U_N\times U_N))\quad,\quad 
u_{ij}\to\begin{pmatrix}0&v_{ij}\\ w_{ij}&0\end{pmatrix}$$
where $v,w$ are the fundamental corepresentations of the two copies of $C(U_N)$.
\end{proposition}

\begin{proof}
It is routine to check that the matrices on the right are indeed biunitaries, and since the first matrix is also self-adjoint, we obtain in this way models as follows:
$$C(O_N^+)\to M_2(C(U_N))\quad,\quad
C(U_N^+)\to M_2(C(U_N\times U_N))$$

Regarding now the half-commutation relations, this comes from something general, regarding the antidiagonal $2\times2$ matrices. Consider indeed matrices as follows:
$$X_i=\begin{pmatrix}0&x_i\\ y_i&0\end{pmatrix}$$

We have then the following computation:
$$X_iX_jX_k
=\begin{pmatrix}0&x_i\\ y_i&0\end{pmatrix}\begin{pmatrix}0&x_j\\ y_j&0\end{pmatrix}\begin{pmatrix}0&x_k\\ y_k&0\end{pmatrix}
=\begin{pmatrix}0&x_iy_jx_k\\ y_ix_jy_k&0\end{pmatrix}$$

Since this quantity is symmetric in $i,k$, we obtain from this:
$$X_iX_jX_k=X_kX_jX_i$$

Thus, the antidiagonal $2\times2$ matrices half-commute, and we conclude that our models for $C(O_N^+)$ and $C(U_N^+)$ constructed above factorize as in the statement.
\end{proof}

We can now formulate our first concrete modeling theorem, as follows:

\index{antidiagonal model}
\index{stationarity}

\begin{theorem}
The above antidiagonal models, namely
$$C(O_N^*)\to M_2(C(U_N))\quad,\quad 
C(U_N^*)\to M_2(C(U_N\times U_N))$$
are both stationary, and in particular they are faithful.
\end{theorem}

\begin{proof}
Let us first discuss the case of $O_N^*$. We will use Theorem 13.13 (3). Since the fundamental representation is self-adjoint, the various matrices $T_e$ with $e\in\{1,*\}^p$ are all equal. We denote this common matrix by $T_p$. We have, by definition:
$$(T_p)_{i_1\ldots i_p,j_1\ldots j_p}
=\left(tr\otimes\int_H\right)\left[\begin{pmatrix}0&v_{i_1j_1}\\\bar{v}_{i_1j_1}&0\end{pmatrix}\ldots\ldots\begin{pmatrix}0&v_{i_pj_p}\\\bar{v}_{i_pj_p}&0\end{pmatrix}\right]$$

Since when multipliying an odd number of antidiagonal matrices we obtain an atidiagonal matrix, we have $T_p=0$ for $p$ odd. Also, when $p$ is even, we have:
\begin{eqnarray*}
(T_p)_{i_1\ldots i_p,j_1\ldots j_p}
&=&\left(tr\otimes\int_H\right)\begin{pmatrix}v_{i_1j_1}\ldots\bar{v}_{i_pj_p}&0\\0&\bar{v}_{i_1j_1}\ldots v_{i_pj_p}\end{pmatrix}\\
&=&\frac{1}{2}\left(\int_Hv_{i_1j_1}\ldots\bar{v}_{i_pj_p}+\int_H\bar{v}_{i_1j_1}\ldots v_{i_pj_p}\right)\\
&=&\int_HRe(v_{i_1j_1}\ldots\bar{v}_{i_pj_p})
\end{eqnarray*}

We have $T_p^2=T_p=0$ when $p$ is odd, so we are left with proving that for $p$ even we have $T_p^2=T_p$. For this purpose, we use the following formula:
$$Re(x)Re(y)=\frac{1}{2}\left(Re(xy)+Re(x\bar{y})\right)$$

By using this identity for each of the terms which appear in the product, and multi-index notations in order to simplify the writing, we obtain:
\begin{eqnarray*}
(T_p^2)_{ij}
&=&\sum_{k_1\ldots k_p}(T_p)_{i_1\ldots i_p,k_1\ldots k_p}(T_p)_{k_1\ldots k_p,j_1\ldots j_p}\\
&=&\int_H\int_H\sum_{k_1\ldots k_p}Re(v_{i_1k_1}\ldots\bar{v}_{i_pk_p})Re(w_{k_1j_1}\ldots\bar{w}_{k_pj_p})dvdw\\
&=&\frac{1}{2}\int_H\int_H\sum_{k_1\ldots k_p}Re(v_{i_1k_1}w_{k_1j_1}\ldots\bar{v}_{i_pk_p}\bar{w}_{k_pj_p})+Re(v_{i_1k_1}\bar{w}_{k_1j_1}\ldots\bar{v}_{i_pk_p}w_{k_pj_p})dvdw\\
&=&\frac{1}{2}\int_H\int_HRe((vw)_{i_1j_1}\ldots(\bar{v}\bar{w})_{i_pj_p})+Re((v\bar{w})_{i_1j_1}\ldots(\bar{v}w)_{i_pj_p})dvdw
\end{eqnarray*}

Now since $vw\in H$ is uniformly distributed when $v,w\in H$ are uniformly distributed, the quantity on the left integrates up to $(T_p)_{ij}$. Also, since $H$ is conjugation-stable, $\bar{w}\in H$ is uniformly distributed when $w\in H$ is uniformly distributed, so the quantity on the right integrates up to the same quantity, namely $(T_p)_{ij}$. Thus, we have:
$$(T_p^2)_{ij}
=\frac{1}{2}\Big((T_p)_{ij}+(T_p)_{ij}\Big)
=(T_p)_{ij}$$

Summarizing, we have obtained that for any $p$, we have $T_p^2=T_p$. Thus Theorem 13.13 applies, and shows that our model is stationary, as claimed. As for the proof of the stationarity for the model for $U_N^*$, this is similar. See \cite{bne}.
\end{proof}

As a second illustration, regarding $H_N^*,K_N^*$, we have:

\begin{theorem}
We have a stationary matrix model as follows, 
$$C(H_N^*)\to M_2(C(K_N))\quad,\quad 
u_{ij}\to\begin{pmatrix}0&v_{ij}\\ \bar{v}_{ij}&0\end{pmatrix}$$
where $v$ is the fundamental corepresentation of $C(K_N)$, as well as a stationary model
$$C(K_N^*)\to M_2(C(K_N\times K_N))\quad,\quad 
u_{ij}\to\begin{pmatrix}0&v_{ij}\\ w_{ij}&0\end{pmatrix}$$
where $v,w$ are the fundamental corepresentations of the two copies of $C(K_N)$.
\end{theorem}

\begin{proof}
This follows by adapting the proof of Proposition 13.14 and Theorem 13.15, by adding there the $H_N^+,K_N^+$ relations. All this is in fact part of a more general phenomenon, concerning half-liberation in general, and we refer here to \cite{bb3}, \cite{bdu}.
\end{proof}

As a consequence of this, we can now work out the discrete group case:

\begin{proposition}
Any reflection group $\Gamma=<g_1,\ldots,g_N>$ which is half-abelian, in the sense that its standard generators half-commute, 
$$g_ig_jg_k=g_kg_jg_i$$
has an algebraic stationary model, with $K=2$.
\end{proposition}

\begin{proof}
This follows from Theorem 13.15. To be more precise, in the non-abelian case, the results in \cite{bdu} show that $\widehat{\Gamma}\subset O_N^*$ must come from a group dual $\widehat{\Lambda}\subset U_N$, via the construction there, and with $\Lambda=<h_1,\ldots,h_N>$, the corresponding model is:
$$\Gamma\subset C(\widehat{\Lambda},U_2)\quad,\quad g_i\to \left[\chi\to\begin{pmatrix}0&\chi(h_i)\\ \bar{\chi}(h_i)&0\end{pmatrix}\right]$$

As for the abelian case, the result here follows from Proposition 13.5. 
\end{proof}

More generally now, we have the following result, from \cite{bch}:

\begin{proposition}
If $L$ is a compact group, having a $N$-dimensional unitary corepresentation $v$, and an order $K$ automorphism $\sigma:L\to L$, we have a matrix model
$$\pi:C(U_N^*)\to M_K(C(L))\quad,\quad u_{ij}\to\tau[v_{ij}^{(1)},\ldots,v_{ij}^{(K)}]$$
where $v^{(i)}(g)=v(\sigma^i(g))$, and where $\tau[x_1,\ldots,x_K]$ is obtained by filling the standard $K$-cycle $\tau\in M_K(0,1)$ with the elements $x_1,\ldots,x_K$. We call such models ``cyclic''.
\end{proposition}

\begin{proof}
The matrices $U_{ij}=\tau[v_{ij}^{(1)},\ldots,v_{ij}^{(K)}]$ in the statement appear by definition as follows, with the convention that all the blank spaces denote 0 entries:
$$U_{ij}=\begin{pmatrix}
&&&v_{ij}^{(1)}\\
v_{ij}^{(2)}\\
&\ddots\\
&&v_{ij}^{(K)}
\end{pmatrix}$$

The matrix $U=(U_{ij})$ is then unitary, and so is $\bar{U}=(U_{ij}^*)$. Thus, if we denote by $w=(w_{ij})$ the fundamental corepresentation of $C(U_N^+)$, we have a model as follows:
$$\rho:C(U_N^+)\to M_K(C(L))\quad,\quad w_{ij}\to U_{ij}$$

Now observe that the matrices $U_{ij}U_{kl}^*,U_{ij}^*U_{kl}$ are all diagonal, so in particular, they commute. Thus the above morphism $\rho$ factorizes through $C(U_N^*)$, as claimed.
\end{proof}

In relation with the above models, we have the following result, also from \cite{bch}:

\begin{theorem}
Any cyclic model in the above sense,
$$\pi:C(U_N^*)\to M_K(C(L))$$
is stationary on its image, with the corresponding closed subgroup $[L]\subset U_N^*$, given by 
$$Im(\pi)=C([L])$$
being the quotient $L\rtimes\mathbb Z_K\to[L]$ having as coordinates the variables $u_{ij}=v_{ij}\otimes\tau$.
\end{theorem}

\begin{proof}
Assuming that $(L,\sigma)$ are as in Proposition 13.18, we have an action $\mathbb Z_K\curvearrowright L$, and we can therefore consider the following short exact sequence:
$$1\to\mathbb Z_K\to L\rtimes\mathbb Z_K\to L\to1$$

By performing a Thoma type construction we obtain a model as follows, where $x^{(i)}=\tilde{\sigma}^i(x)$, with $\tilde{\sigma}:C(L)\to C(L)$ being the automorphism induced by $\sigma:L\to L$:
$$\rho:C(L\rtimes\mathbb Z_K)\subset M_K(C(L))\quad,\quad x\otimes\tau^i\to\tau^i[x^{(1)},\ldots,x^{(K)}]$$

Consider now the quotient quantum group $L\rtimes\mathbb Z_K\to[L]$ having as coordinates the variables $u_{ij}=v_{ij}\otimes\tau$. We have then a injective morphism, as follows:
$$\nu:C([L])\subset C(L\rtimes\mathbb Z_K)\quad,\quad u_{ij}\to v_{ij}\otimes\tau$$

By composing the above two embeddings, we obtain an embedding as follows:
$$\rho\nu:C([L])\subset M_K(C(L))\quad,\quad u_{ij}\to\tau[v_{ij}^{(1)},\ldots,v_{ij}^{(K)}]$$

Now since $\rho$ is stationary, and since $\nu$ commutes with the Haar funtionals as well, it follows that this morphism $\rho\nu$ is stationary, and this finishes the proof.
\end{proof}

As an illustration, we can now recover the following result, from \cite{bdu}:

\begin{proposition}
For any non-classical $G\subset O_N^*$ we have a stationary model
$$\pi:C(G)\to M_2(C(L))\quad,\quad u_{ij}=\begin{pmatrix}0&v_{ij}\\ \bar{v}_{ij}&0\end{pmatrix}$$
where $L\subset U_N$, with coordinates denoted $v_{ij}$, is the lift of $PG\subset PO_N^*=PU_N$.
\end{proposition}

\begin{proof}
Assume first that $L\subset U_N$ is self-conjugate, in the sense that $g\in L\implies\bar{g}\in L$. If we consider the order 2 automorphism of $C(L)$ induced by $g_{ij}\to\bar{g}_{ij}$, we can apply Theorem 13.19, and we obtain a stationary model, as follows:
$$\pi:C([L])\subset M_2(C(L))\quad,\quad u_{ij}\otimes1=\begin{pmatrix}0&v_{ij}\\ \bar{v}_{ij}&0\end{pmatrix}$$

The point now is that, as explained in \cite{bdu}, any non-classical subgroup $G\subset O_N^*$ must appear as $G=[L]$, for a certain self-conjugate subgroup $L\subset U_N$. Moreover, since we have $PG=P[L]$, it follows that $L\subset U_N$ is the lift of $PG\subset PO_N^*=PU_N$, as claimed. 
\end{proof}

In the unitary case now, and with the matrix size $K\in\mathbb N$ being arbitrary, we recall from \cite{bb3}, \cite{mwe} and related papers that $U_N^*$ has a certain ``arithmetic version'' $U_{N,K}^*\subset U_N^*$, obtained by imposing some natural length $2K$ relations on the standard coordinates. As basic examples, at $K=1$ we have $U_{N,1}^*=U_N$, the defining relations being $ab=ba$ with $a,b\in\{u_{ij},u_{ij}^*\}$, and at $K=2$ we have $U_{N,2}^*=U_N^{**}$, with the latter quantum group appearing via the relations $ab\cdot cd=cd\cdot ab$, for any $a,b,c,d\in\{u_{ij},u_{ij}^*\}$. 

\bigskip

With this convention, we have the following result, also from \cite{bch}:

\begin{theorem}
For any subgroup $G\subset U_{N,K}^*$ which is $K$-symmetric, in the sense that $u_{ij}\to e^{2\pi i/K}u_{ij}$ defines an automorphism of $C(G)$, we have a stationary model
$$\pi:C(G)\to M_K(C(L))\quad,\quad u_{ij}\to\tau[v_{ij}^{(1)},\ldots,v_{ij}^{(K)}]$$
with $L\subset U_N^K$ being a closed subgroup which is symmetric, in the sense that it is stable under the cyclic action $\mathbb Z_K\curvearrowright U_N^K$.
\end{theorem}

\begin{proof}
This follows from what we have, as follows:

\medskip

(1) Assuming that $L\subset U_N^K$ is symmetric in the above sense, we have representations $v^{(i)}:L\subset U_N^K\to U_N^{(i)}$ for any $i$, and the cyclic action $\mathbb Z_K\curvearrowright U_N^K$ restricts into an order $K$ automorphism $\sigma:L\to L$. Thus we can apply Theorem 13.19, and we obtain a certain closed subgroup $[L]\subset U_{N,K}^*$, having a stationary model as in the statement.

\medskip

(2) Conversely now, assuming that $G\subset U_{N,K}^*$ is $K$-symmetric, the main result in \cite{bb3} applies, and shows that we must have $C(G)\subset C(L)\rtimes\mathbb Z_K$, for a certain closed subgroup $L\subset U_N^K$ which is symmetric. But this shows that we have $G=[L]$, and we are done.
\end{proof}

We refer to \cite{bb3}, \cite{bch}, \cite{bdu}, \cite{mwe} for more on the above.

\section*{13d. Group duals}

Let us discuss now the group dual case, where we have a closed subgroup $\widehat{\Gamma}\subset S_N^+$, with $\Gamma$ being a discrete group. Following \cite{bch}, \cite{bfr}, we use the following construction:

\begin{proposition}
The following happen:
\begin{enumerate}
\item Given integers $K_1,\ldots,K_M$ satisfying $K_1+\ldots+K_M=N$, the dual of any quotient group $\mathbb Z_{K_1}*\ldots*\mathbb Z_{K_M}\to\Gamma$ appears as a closed subgroup $\widehat{\Gamma}\subset S_N^+$. 

\item By refining if necessary the partition $N=K_1+\ldots+K_M$, we can always assume that the $M$ morphisms $\mathbb Z_{K_i}\to\Gamma$ are all injective.

\item Assuming that the partition $N=K_1+\ldots+K_M$ is refined, as above, this partition is precisely the one describing the orbit structure of $\widehat{\Gamma}\subset S_N^+$.
\end{enumerate}
\end{proposition}

\begin{proof}
The idea for (1) is that we have embeddings $\widehat{\mathbb Z}_{K_i}\simeq\mathbb Z_{K_i}\subset S_{K_i}\subset S_{K_i}^+$, and by performing a free product construction, we obtain an embedding as follows:
$$\widehat{\Gamma}\subset\widehat{\mathbb Z_{K_1}*\ldots*\mathbb Z_{K_M}}\subset S_N^+$$

To be more precise, the magic unitary that we get is as follows, where $F_i=\frac{1}{\sqrt{K_i}}(w_i^{ab})_{ab}$ with $w_i=e^{2\pi i/K_i}$, and  $V_i=(g_i^a)_a$, with $g_i$ being the standard generator of $\mathbb Z_{K_i}$:
$$u=diag(u_i)\quad,\quad u_i=\frac{1}{\sqrt{K_i}}\begin{pmatrix}(F_iV_i)_0&\ldots&(F_iV_i)_{K_i-1}\\ (F_iV_i)_{K_i-1}&\ldots&(F_iV_i)_{K_i-2}\\ \vdots&\vdots&\vdots\\ 
(F_iV_i)_1&\ldots&(F_iV_i)_0\end{pmatrix}$$

Regarding (2,3), the idea here is that the orbit structure of any $\widehat{\Gamma}\subset S_N^+$ produces a partition $N=K_1+\ldots+K_M$, and then a quotient map $\mathbb Z_{K_1}*\ldots*\mathbb Z_{K_M}\to\Gamma$.
\end{proof}

Following the material from chapter 10, we will be mainly interested in what follows in the quasi-transitive case. Let us start with the following definition:

\begin{definition}
Given a subgroup $G\subset S_N^+$, a random matrix model of type
$$\pi:C(G)\to M_K(C(T))$$
is called quasi-flat when the fibers $P_{ij}^x=\pi(u_{ij})(x)$ all have rank $\leq1$. 
\end{definition}

We will explore more in detail this notion in chapter 11. Now with this convention made, and getting back to the group duals, we have the following result, from \cite{bch}:

\begin{proposition}
The quasi-transitive group duals $\widehat{\Gamma}\subset S_N^+$, with orbits having $K$ elements, appearing as above, have the following properties:
\begin{enumerate}
\item These come from the quotients $\mathbb Z_K^{*M}\to\Gamma$, having the property that the corresponding $M$ morphisms $\mathbb Z_K^{(i)}\subset\mathbb Z_K^{*M}\to\Gamma$ are all injective.

\item For such a quotient, a matrix model $\pi:C^*(\Gamma)\to M_K(\mathbb C)$ is quasi-flat if and only if it is stationary on each subalgebra $C^*(\mathbb Z_K^{(i)})\subset C^*(\Gamma)$.
\end{enumerate}
\end{proposition}

\begin{proof}
The first assertion follows from Proposition 13.23. Regarding the second assertion, consider a matrix model $\pi:C^*(\Gamma)\to M_K(\mathbb C)$, mapping $g_i\to U_i$, where $g_i$ is the standard generator of $\mathbb Z_K^{(i)}$. With notations from the proof of Proposition 13.23, the images of the nonzero standard coordinates on $\widehat{\Gamma}\subset S_N^+$ are mapped as follows:
$$\pi:\frac{1}{\sqrt{K}}(FV_i)_c\to\frac{1}{\sqrt{K}}(FW_i)_c$$

Here $V_i=(g_i^a)_a$, $W_i=(U_i^a)_a$, and $F=\frac{1}{\sqrt{K}}(w^{ab})_{ab}$ with $w=e^{2\pi i/K}$. With this formula in hand, the flatness condition on $\pi$ simply states that we must have:
$$Tr((FW_i)_c)=\sqrt{K}\quad,\quad\forall i,\forall c$$

In terms of the trace vectors $T_i=(Tr(U_i^a))_a$ this condition becomes $FT_i=\sqrt{K}\xi$, where $\xi\in\mathbb C^K$ is the all-one vector. Thus we must have $T_i=\sqrt{K}F^*\xi$, which reads:
$$\begin{pmatrix}Tr(1)\\ Tr(U_i)\\\vdots\\ Tr(U_i^{K-1})\end{pmatrix}=\sqrt{K}F^*\begin{pmatrix}1\\ 1\\\vdots\\1\end{pmatrix}=\begin{pmatrix}K\\ 0\\\vdots\\0\end{pmatrix}\quad,\quad\forall i$$
 
In other words, we have reached to the conclusion that $\pi$ is flat precisely when its restrictions to each subalgebra $C^*(\mathbb Z_K^{(i)})\subset C^*(\Gamma)$ are stationary, as claimed.
\end{proof}

We would like to end our study with a purely group-theoretical formulation of these results, and of some related questions, that we believe of interest. Let us start with:

\begin{definition}
A discrete group $\Gamma$ is called uniform when:
\begin{enumerate}
\item $\Gamma$ is finitely generated, $\Gamma=<g_1,\ldots,g_M>$.

\item The generators $g_1,\ldots,g_M$ have common order $K<\infty$.

\item $\Gamma$ appears as an intermediate quotient $\mathbb Z_K^{*M}\to\Gamma\to\mathbb Z_K^M$.

\item We have an action $S_M\curvearrowright\Gamma$, given by $\sigma(g_i)=g_{\sigma(i)}$.
\end{enumerate}
\end{definition}

Here the conditions (1-3) basically come from \cite{bi3}, via Proposition 13.24 (1), and together with some extra considerations from \cite{bfr}, which prevent us from using groups of type $\Gamma=(\mathbb Z_K*\mathbb Z_K)\times\mathbb Z_K$, we are led to the condition (4) as well.

\bigskip

Observe that some of the above conditions are technically redundant, with (4) implying that the generators $g_1,\ldots,g_M$ have common order, as stated in (2), and also with (3) implying that the group is finitely generated, with the generators having finite order. We have as well the following notion, which is once again group-theoretical:

\begin{definition}
If a discrete group $\Gamma$ is uniform, as above, a unitary representation $\rho:\Gamma\to U_K$ is called quasi-flat when the eigenvalues of each 
$$U_i=\rho(g_i)\in U_K$$
are uniformly distributed.
\end{definition}

To be more precise, assuming that $\Gamma=<g_1,\ldots,g_M>$ with $ord(g_i)=K$ is as in Definition 13.25, any unitary representation $\rho:\Gamma\to U_K$ is uniquely determined by the images $U_i=\rho(g_i)\in U_K$ of the standard generators. Now since each of these unitaries satisfies $U_i^K=1$, its eigenvalues must be among the $K$-th roots of unity, and our quasi-flatness assumption states that each eigenvalue must appear with multiplicity $1$.

\bigskip

With these notions in hand, we have the following result:

\begin{theorem}
If $\Gamma=<g_1,\ldots,g_M>$ is uniform, with $ord(g_i)=K$, a model 
$$\pi:C^*(\Gamma)\to M_K(C(X))$$
is quasi-flat precisely when the associated unitary representation 
$$\rho:\Gamma\to C(X,U_K)$$
has quasi-flat fibers, in the sense of Definition 13.26.
\end{theorem} 

\begin{proof}
According to Proposition 13.24 (2), the model is quasi-flat precisely when the following compositions are all stationary:
$$\pi_i:C^*(\mathbb Z_K^{(i)})\subset C^*(\Gamma)\to M_K(C(X))$$

On the other hand, as already observed in the proof of Proposition 13.24, a matrix model $\rho:C^*(\mathbb Z_K)\to M_K(C(X))$ is stationary precisely when the unitary $U=\rho(g)$, where $g$ is the standard generator of $\mathbb Z_K$, satisfies the following condition:
$$\begin{pmatrix}tr(1)\\ tr(U)\\\vdots\\ tr(U^{K-1})\end{pmatrix}=\begin{pmatrix}1\\ 0\\\vdots\\0\end{pmatrix}$$

Thus, such a model is stationary precisely when the eigenvalues of $U$ are uniformly distributed, over the $K$-th roots of unity. We conclude that $\pi$ is quasi-flat precisely when the eigenvalues of each $U_i=\rho(g_i)$ are uniformly distributed, as in Definition 13.26.
\end{proof}

We are interested now in the matrix models for the discrete group algebras, which are stationary. We use a lift of the quasi-flat models, in the following sense:

\begin{proposition}
The affine lift of the universal quasi-flat model for $C^*(\mathbb Z_K^{\ast M})$,
$$\pi:C^*(\mathbb Z_K^{\ast M})\to M_K(C(U_K^M))$$
is given on the canonical generator $g_{i}$ of the $i$-th factor by
$$\pi(g_{i})(U^{1}, \dots, U^{M})=\sum_jw^jP_{U_j^i}$$
where $U_j^i$ is the $j$-th column of $U^i$ and $P_{\xi}$ denotes the orthogonal projection onto $\mathbb C\xi$.
\end{proposition}

\begin{proof}
There is indeed a canonical quotient map $U_K\to E_K$, obtained by parametrizing the orthonormal bases of $\mathbb C^K$ by the unitary group $U_K$, and this gives the result.
\end{proof}

We know that the maximal group dual subgroups $\widehat{\Gamma}\subset S_N^+$ are the free products of type $\mathbb Z_{K_1}*\ldots*\mathbb Z_{K_M}$ with $K_1+\ldots+K_M=N$. In the quasi-transitive case, where by definition $K_1=\ldots=K_M=K$ with $K|N$, we have the following result, from \cite{bfr}:

\begin{theorem}
The universal quasi-flat model for the group
$$\Gamma=\mathbb Z_K^{*M}$$
is inner faithful.
\end{theorem}

\begin{proof}
It is enough to prove that the affine lift of the universal model in the statement is inner faithful. For this purpose, let us consider a reduced word $\gamma\in\mathbb Z_K^{*M}$, and write it as follows, with indices $i_t\neq i_{t+1}$, and with exponents $1\leq k_t\leq K-1$:
$$\gamma=g_{i_1}^{k_1}\ldots g_{i_n}^{k_n}$$

With this convention, we have then the following computation:
\begin{eqnarray*}
\pi(\gamma)(U^1,\ldots ,U^M) 
&=&\sum_{j_1\ldots j_n=1}^Kw^{k_1j_1}P_{U^{i_1}_{j_1}}\ldots w^{k_nj_n}P_{U^{i_n}_{j_n}}\\
&=&\sum_{j_1\ldots j_n=1}^Kw^{<k,j>}P_{U^{i_1}_{j_1}}\ldots P_{U^{i_n}_{j_n}}
\end{eqnarray*}

Our aim is to prove that there is at least one sequence $(U^1,\ldots ,U^M)$ for which the above matrix is not the identity. For this purpose, we use the following formula:
$$P_{\xi_1}\ldots P_{\xi_l}(x)=<x,\xi_l><\xi_l,\xi_{l-1}>\ldots\ldots<\xi_2, \xi_1>\xi_1$$

To compute the trace of this operator, we can consider any orthonormal basis containing $\xi_l$, yielding $<\xi_1,\xi_l><\xi_l,\xi_{l-1}>\ldots\ldots<\xi_2, \xi_1>$. Applying this to $\pi(\gamma)$ and using the equality $<U_i,V_j>=\sum_lU_{ki}\bar{V}_{kj}=(V^*U)_{ji}$, we get:
\begin{eqnarray*}
tr\circ\pi(\gamma)
&=&\frac{1}{K}\sum_{j_1\ldots j_n=1}^Kw^{<k,j>}<U^{i_1}_{j_1}, U^{i_n}_{j_n}><U^{i_n}_{j_n},U^{i_{n-1}}_{j_{n-1}}>\ldots <U^{i_2}_{j_2},U^{i_1}_{j_1}>\\
&=&\frac{1}{K}\sum_{j_1\ldots j_n=1}^Kw^{<k,j>}(U^{i_n*}U^{i_1})_{j_nj_1}(U^{i_{n-1}*}U^{i_n})_{j_{n-1}j_n}\ldots (U^{i_1*}U^{i_2})_{j_1j_2}
\end{eqnarray*}

Denoting by $W$ the diagonal matrix given by $W_{ij}=\delta_{ij}w^i$, we have:
$$\sum_{j_1}w^{k_1j_1}U_{j_nj_{1}}U^*_{j_1j_2}=\sum_{j_1l}U_{j_nj_1}W^{k_1}_{j_1l}U^*_{lj_2}=(UW^{k_1}U^*)_{j_nj_2}$$

Applying this $n$ times in the above formula for $tr\circ\pi(\gamma)$ yields:
\begin{eqnarray*}
tr\circ\pi(\gamma)
&=&tr\left(U^{i_n*}U^{i_1}W^{k_1}U^{i_1*}U^{i_2}W^{k_2}\ldots W^{k_{n-1}}U^{i_{n-1}*}U^{i_n}W^{k_n}\right)\\
&=&tr\left(U^{i_1}W^{k_1}U^{i_1*}\ldots U^{i_n}W^{k_n}U^{i_n*}\right)
\end{eqnarray*}

Assume now that $\pi(\gamma)(U^1,\ldots,U^M) = Id$ for all sequences of unitary matrices. The trace of a unitary matrix can only be equal to $1$ if it is the identity, hence:
$$\prod_{p=1}^nU^{i_p}W^{k_p}U^{i_p*} = Id$$

In other words, the following noncommutative polynomial vanishes on $U_K^M$:
$$P = \prod_{p=1}^nX^{i_p}W^{k_p}X^{i_p*} - 1$$

But this is clearly impossible if $k_t\neq 0(K)$ for all $t$, hence $\pi(\gamma)$ is not always the identity, and so our representation $\pi$ is inner faithful, as desired.
\end{proof}

More generally now, we have the following result, also from \cite{bfr}:

\begin{theorem}
The universal quasi-flat model for the group
$$\Gamma=\mathbb Z_K^{M_1}*\ldots*\mathbb Z_K^{M_n}$$
is inner faithful.
\end{theorem}

\begin{proof}
We can consider the space $U_K\times S_K^M$ as the affine lift of our model. If $g_1,\ldots, g_M$ are the canonical generators of the direct product, their action is then given by:
$$\pi(g_i)(U, \sigma_1, \ldots, \sigma_M) = \sum_{j=1}^Kw^jP_{U_{\sigma_i^{-1}(j)}} = \sum_{j=1}^Kw^{\sigma_i(j)}P_{U_j}$$

This gives the following formula for a general element:
$$\pi(g_1^{k_1}\ldots g_M^{k_M})(U, \sigma_1, \ldots, \sigma_M) = \sum_{j=1}^Kw^{k_1\sigma_1(j) + \ldots + k_M\sigma_M(j)}P_{U_j}$$

Let $g_1(i), \ldots, g_{M_i}(i)$ be the generators of $\mathbb{Z}_K^{M_i}$, and let consider a reduced word:
$$\gamma = \left(g_1(i_1)^{k_1(1)}\ldots g_{M_{i_1}}(i_1)^{k_{M_{i_1}}(1)}\right)\ldots \left(g_1(i_n)^{k_1(n)}\ldots g_{M_{i_n}}(i_n)^{k_{M_{i_n}}(n)}\right)$$

The computation of $tr\circ\pi(\gamma)$ is then similar to the one in the proof of Theorem 13.29, until the introduction of the matrices $W$. Here we have to replace $W^{k_t}$ by $\prod_{s=1}^{M_t}W_{\sigma_s^t}^{k_s(t)}$, where $(U^t, \sigma^t_1, \ldots, \sigma^t_{M_t})_{1 \leq t\leq n}$ is the element to which we are applying $\pi(\gamma)$, and:
$$(W_{\sigma})_{ij} = \delta_{ij}w^{\sigma(i)}$$

Assuming that $\pi(\gamma) = 1$, we can apply the same strategy as before. Indeed, we have a polynomial which must vanish on all sequences of unitary matrices, and this is impossible unless all the matrices appearing in the polynomial are the identity. Thus, we have:
$$\prod_{s=1}^{M_t}W_{\sigma_s^t}^{k_s(t)} = Id\quad,\quad\forall t$$

But this condition translates into the following condition:
$$\prod_{s=1}^{M_t}w^{\sigma_s^t(i)k_s(t)} = 1$$

We can sum the above equation over all permutations, and we obtain:
$$\frac{1}{(K!)^M}\sum_{\sigma_1,\ldots,\sigma_{M_i}\in S_K}w^{k_1(t)\sigma_1(i) + \ldots + k_{M_t}(t)\sigma_{M_t}(i)} = \frac{1}{(K!)^M}\prod_{s=1}^{M_t}\left(\sum_{\sigma_s\in S_K}w^{k_s(t)\sigma_s(i)}\right)$$

For any $i'$, there are $(K-1)!$ permutations $\sigma$ such that $\sigma_s(i) = i'$. This leads to:
$$\sum_{\sigma_s\in S_K}w^{k_s(t)\sigma_s(i)}
=\sum_{i'=1}^K(K-1)!w^{k_s(t)i'}
=K!\delta_{k_s(t), 0}$$

Now by putting everything together, we obtain the following formula:
$$\prod_{s=1}^{M_i}\delta_{k_s(t), 0} = 1$$

Thus $k_s(t) = 0$ for all $t$ and all $s$, which is a contradiction, as desired.
\end{proof}

\section*{13e. Exercises} 

The matrix modeling problematics is something quite exciting, and we have several exercises here. In relation with the notion of stationarity, we have:

\begin{exercise}
Work out the details for the fact that the stationarity of a model
$$\pi:C(G)\to M_K(C(T))$$
implies its faithfulness.
\end{exercise}

This is something that we already discussed in the above, but with some standard functional analysis details missing. The problem is that of working out these details.

\begin{exercise}
Find an example of an inner faithful model 
$$\pi:C(G)\to M_K(C(T))$$
which is not faithful, not coming from a classical group, or a group dual.
\end{exercise}

This is something quite tricky, and it is of course possible to cheat a bit here, by using product operations. The exercise asks for a high-quality counterexample.

\chapter{Flat models}

\section*{14a. Quasi-flatness}

We have seen in the previous chapter how to use matrix model technology for the closed subgroups $G\subset U_N^+$, and we discussed some basic examples of models as well. In this chapter we restrict the attention to the quantum permutation groups, $G\subset S_N^+$, and we further develop the general theory that we have, for such quantum groups.

\bigskip

According to our general matrix model philosophy, we must look for concrete realizations of the standard coordinates $u_{ij}\in C(G)$. In the quantum permutation group case, $G\subset S_N^+$, these standard coordinates are projections, so a natural idea is that of looking at the rank of these projections, in the model. We are led to the following notion:

\begin{definition}
Given a subgroup $G\subset S_N^+$, a random matrix model of type
$$\pi:C(G)\to M_K(C(T))$$
is called flat when the fibers $P_{ij}^x=\pi(u_{ij})(x)$ all have rank $1$. 
\end{definition}

We will see in what follows that there are many interesting examples of such models. Observe however, as an obstruction, that in order for such a model to exist, the subgroup $G\subset S_N^+$ must be transitive, due to the following trivial fact:
$$P_{ij}\neq0\implies u_{ij}\neq0$$

In order to include as well non-transitive subgroups $G\subset S_N^+$ in our discussion, as for instance certain group duals $\widehat{\Gamma}\subset S_N^+$, it is convenient to enlarge our formalism, and we have here the following notion, that we already met in chapter 13:

\begin{definition}
Given a subgroup $G\subset S_N^+$, a random matrix model of type
$$\pi:C(G)\to M_K(C(T))$$
is called quasi-flat when the fibers $P_{ij}^x=\pi(u_{ij})(x)$ all have rank $\leq1$. 
\end{definition}

As a first observation, the functions $x\to r_{ij}^x=rank(P_{ij}^x)$ are locally constant over $T$, so they are constant over the connected components of $T$. Thus, when $T$ is connected, our assumption is that we have $r_{ij}^x=r_{ij}\in\{0,1\}$, for any $x\in T$, and any $i,j$. Observe also that when $K=N$ these questions disappear, because in this case we must have $r_{ij}^x=1$ for any $i,j$ and any $x\in t$, and the model is flat. As a first result now, we have:

\begin{proposition}
Assume that we have a quasi-flat matrix model,
$$\pi:C(G)\to M_K(C(T))\quad,\quad u_{ij}\to P_{ij}$$
and consider the matrix $r_{ij}=rank(P_{ij})$. The following happen:
\begin{enumerate}
\item $r$ is bistochastic, with row and column sums $K$.

\item $r_{ij}\leq\varepsilon_{ij}$, where $\varepsilon\in M_N(0,1)$ is given by $\varepsilon_{ij}=\delta_{u_{ij},0}$.

\item If $G$ is quasi-transitive, with orbits of size $K$, then $r=\varepsilon$.

\item If $\pi$ is assumed to be flat, $G$ must be transitive.
\end{enumerate}
\end{proposition}

\begin{proof}
These results are all elementary, the proofs being as follows:

\medskip

(1) This is clear from the fact that each $P^x=(P_{ij}^x)$ is bistochastic, with sums $1$.

\medskip

(2) This simply comes from $u_{ij}=0\implies P_{ij}=0$.

\medskip

(3) The matrices $r=(r_{ij})$ and $\varepsilon=(\varepsilon_{ij})$ are both bistochastic, with sums $K$, and they satisfy $r_{ij}\leq\varepsilon_{ij}$, for any $i,j$. Thus, these matrices must be equal, as stated.

\medskip

(4) This is clear, because $rank(P_{ij})=1$ implies $u_{ij}\neq0$, for any $i,j$.
\end{proof}

It is convenient to identify the rank one projections in $M_N(\mathbb C)$ with the elements of the complex projective space $P^{N-1}_\mathbb C$. We have then the following result, from \cite{bne}:

\begin{proposition}
The algebra $C(S_N^+)$ has a universal flat model, given by 
$$\pi_N:C(S_N^+)\to M_N(C(T_N))\quad,\quad\pi_N(u_{ij})=[P\to P_{ij}]$$
where $T_N$ is the set of matrices $P\in M_N(P^{N-1}_\mathbb C)$ which are bistochastic with sums $1$.
\end{proposition}

\begin{proof}
This is clear from definitions, because any flat model $C(S_N^+)\to M_N(\mathbb C)$ must map the magic corepresentation $u=(u_{ij})$ into a matrix $P=(P_{ij})$ belonging to $T_N$.
\end{proof}

We will be back to the above universal models later. Regarding now the general quasi-transitive case, we have here the following result, from \cite{bfr}:

\begin{theorem}
Given a quasi-transitive subgroup $G\subset S_N^+$, with orbits of size $K$, we have a universal quasi-flat model $\pi:C(G)\to M_K(C(T))$, constructed as follows:
\begin{enumerate}
\item For $G=\underbrace{S_K^+\,\hat{*}\,\ldots\,\hat{*}\,S_K^+}_{M\ terms}$ with $N=KM$, the model space is $T_{N,K}=\underbrace{T_K\times\ldots\times T_K}_{M\ terms}$, and with $u=diag(u^1,\ldots,u^M)$ the modelling map is:
$$\pi_{N,K}(u^r_{ij})=[(P^1,\ldots,P^M)\to P^r_{ij}]$$

\item In general, the model space is the submanifold $T_G\subset T_{N,K}$ obtained via the Tannakian relations defining $G$.
\end{enumerate}
\end{theorem}

\begin{proof}
The idea is to use Proposition 14.3 and Proposition 14.4, as follows:

\medskip

(1) This follows from Proposition 14.4, by using Proposition 14.3 (3), which tells us that the 0 entries of the model must appear exactly where $u=(u_{ij})$ has 0 entries.

\medskip

(2) Assume that $G\subset S_N^+$ is quasi-transitive, with orbits of size $K$. We have then an inclusion $G\subset\underbrace{S_K^+\,\hat{*}\,\ldots\,\hat{*}\,S_K^+}_{M\ terms}$, and in order to construct the universal quasi-flat model for $C(G)$, we need a universal solution to the following factorization problem:
$$\begin{matrix}
C(\underbrace{S_K^+\,\hat{*}\,\ldots\,\hat{*}\,S_K^+}_{M\ terms})&\to&M_K(C(T_{N,K}))\\
\\
\downarrow&&\downarrow\\
\\
C(G)&\to&M_K(C(T_G))
\end{matrix}$$

But, the solution to this latter question is given by the following construction, with the Hom-spaces at left being taken as usual in a formal sense:
$$C(T_G)=C(T_{N,K})\Big/\Big(T\in Hom(P^{\otimes k},P^{\otimes l}),\forall k,l\in\mathbb N,\forall T\in Hom(u^{\otimes k},u^{\otimes l})\Big)$$

With this result in hand, the Gelfand spectrum of the algebra on the left is then an algebraic submanifold $T_G\subset T_{N,K}$, having the desired universality property.
\end{proof}

Let us discuss now the classical case, $G\subset S_N$. We use the following notion:

\begin{definition}
A ``sparse Latin square'' is a square matrix 
$$L\in M_N(*,1,\ldots,K)$$
having the property that each of its rows and columns consists of a permutation of the numbers $1,\ldots,K$, completed with $*$ entries.
\end{definition}

In the case $K=N$, where there are no $*$ symbols, we recover the usual Latin squares. In general, however, the combinatorics of these matrices is more complicated than that of the usual Latin squares. Here are a few examples of such matrices:
$$\begin{pmatrix}1&2&*\\2&*&1\\ *&1&2\end{pmatrix}\quad,\quad
\begin{pmatrix}1&2&*&*\\ 2&*&1&*\\ *&1&*&2\\ *&*&2&1\end{pmatrix}\quad,\quad
\begin{pmatrix}1&2&*&*\\ 2&1&*&*\\ *&*&1&2\\ *&*&2&1\end{pmatrix}\quad,\quad
\begin{pmatrix}1&2&*&*&*\\ 2&*&1&*&*\\ *&1&2&*&*\\ *&*&*&1&2\\ *&*&*&2&1
\end{pmatrix}$$

With this notion in hand, the result that we need is as follows:

\begin{proposition}
The quasi-flat representations $\pi:C(S_N)\to M_K(\mathbb C)$ appear as
$$u_{ij}\to P_{L_{ij}}$$ 
where $P_1,\ldots,P_K\in M_K(\mathbb C)$ are rank $1$ projections, summing up to $1$, and where $L\in M_N(*,1,\ldots,K)$ is a sparse Latin square, with the convention $P_*=0$. 
\end{proposition}

\begin{proof}
This is something elementary, as follows:

\medskip

(1) Assuming that $\pi:C(S_N)\to M_K(\mathbb C)$ is quasi-flat, the elements $P_{ij}=\pi(u_{ij})$ are projections of rank $\leq1$, which pairwise commute, and form a magic unitary. 

\medskip

(2) Let $P_1,\ldots,P_K\in M_K(\mathbb C)$ be the rank one projections appearing in the first row of $P=(P_{ij})$. Since these projections form a partition of unity with rank one projections, any rank one projection $Q\in M_K(\mathbb C)$ commuting with all of them satisfies:
$$Q\in\{P_1,\ldots,P_K\}$$

In particular we have $P_{ij}\in\{P_1,\ldots,P_K\}$ for any $i,j$ such that $P_{ij}\neq 0$. Thus we can write $u_{ij}\to P_{L_{ij}}$, for a certain matrix $L\in M_N(*,1,\ldots,K)$, with the convention $P_*=0$.

\medskip

(3) The remark now is that $u_{ij}\to P_{L_{ij}}$ defines a representation $\pi:C(S_N)\to M_K(\mathbb C)$ precisely when the matrix $P=(P_{L_{ij}})_{ij}$ is magic. But this condition tells us precisely that $L$ must be a sparse Latin square, in the sense of Definition 14.6.
\end{proof}

Our task now is to compute the associated Hopf image. Following \cite{bfr}, we have:

\begin{theorem}
Given a sparse Latin square $L\in M_N(*,1,\ldots,K)$, consider the permutations $\sigma_1,\ldots,\sigma_K\in S_N$ given by:
$$\sigma_x(j)=i\iff L_{ij}=x$$
The Hopf image associated to a representation as above,
$$\pi:C(S_N)\to M_K(\mathbb C)\quad,\quad u_{ij}\to P_{L_{ij}}$$
is then the algebra $C(G_L)$, where $G_L=<\sigma_1,\ldots,\sigma_K>\subset S_N$.
\end{theorem}

\begin{proof}
The image of $\pi$ being generated by $P_1,\ldots,P_K$, we have an isomorphism of algebras $\alpha:Im(\pi)\simeq C(1,\ldots,K)$ given by $P_i\to\delta_i$. Consider the following diagram:
$$\xymatrix@R=10mm@C=10mm{
C(S_N)\ar[r]^\pi\ar[dr]_\varphi&Im(\pi)\ar[r]\ar[d]^\alpha&M_K(\mathbb C)\\
&C(1,\ldots,K)
}$$

Here the map on the right is the canonical inclusion and $\varphi=\alpha\pi$. Since the Hopf image of $\pi$ coincides with the one of $\varphi$, it is enough to compute the latter. We know that $\varphi$ is given by $\varphi(u_{ij})=\delta_{L_{ij}}$, with the convention $\delta_*=0$. By Gelfand duality, $\varphi$ must come from a certain map $\sigma:\{1,\ldots,K\}\to S_N$, via the following transposition formula:
$$\varphi(f)(x)=f(\sigma_x)$$

With the choice $f=u_{ij}$, we obtain the following formula:
$$\delta_{L_{ij}}(x)=u_{ij}(\sigma_x)$$

Now observe that these two quantities are by definition given by:
$$\delta_{L_{ij}}(x)=\begin{cases}
1&{\rm if}\ L_{ij}=x\\
0&{\rm otherwise}
\end{cases}\qquad,\qquad 
u_{ij}(\sigma_x)=\begin{cases}
1&{\rm if}\ \sigma_x(j)=i\\
0&{\rm otherwise}
\end{cases}
$$

We conclude that $\sigma_x$ is the permutation in the statement. Summarizing, we have shown that $\varphi$ comes by transposing the map $x\to\sigma_x$, with $\sigma_x$ being as in the statement. Thus the Hopf image of $\varphi$ is the algebra $C(G_L)$, with $G_L=<\sigma_1,\ldots,\sigma_K>$, as desired.
\end{proof}

Let us discuss now the construction of the universal model space. We agree as before to identify the rank one projections in $M_K(\mathbb C)$ with the corresponding elements of the projective space $P^{K-1}_\mathbb C$. With these conventions, we first have the following result:

\begin{proposition}
Assuming that $G\subset S_N$ is quasi-transitive, with orbits of size $K$, the universal quasi-flat model space for $G$ is given by $X_G=E_K\times L_{N,K}^G$, where:
$$E_K=\left\{(P_1,\ldots,P_K)\in (P^{K-1}_\mathbb C)^K\Big|P_i\perp P_j,\forall i\neq j\right\}$$
$$L_{N,K}^G=\left\{L\in M_N(*,1,\ldots,K)
\text{ sparse Latin square }
\Big|G_L\subset G\right\}$$ 
In particular, $X_G$ has a canonical probability measure, obtained as the homogeneous space measure on $E_K$ times the the normalized counting measure on $L_{N,K}^G$.
\end{proposition}

\begin{proof}
The first assertion follows by combining Proposition 14.7 and Theorem 14.8, and the second assertion is clear from definitions.
\end{proof}

Note that the model above is empty if there is no sparse Latin square $L$ such that $G_L\subset G$. The existence of such a sparse Latin square is a strong condition on $G$, and here is an intrinsic characterization of the groups satisfying this condition:

\begin{proposition}
Let $G\subset S_N$ be a classical group. The following are equivalent:
\begin{enumerate}
\item $L_{N,K}^G\neq \emptyset$.

\item There exist $K$ elements $\sigma_1, \dots, \sigma_K\in G$ such that $\sigma_1(i), \dots, \sigma_K(i)$ are pairwise distinct for all $1\leq i\leq N$.
\end{enumerate}
\end{proposition}

\begin{proof}
This is something elementary, as follows:

\medskip

$(1)\implies (2)$ If $\sigma_1, \dots, \sigma_K$ are the permutations associated to $L$, then $\sigma_x(i) = \sigma_y(i)$ means by definition that $x = L_{i\sigma_x(i)} = L_{i\sigma_y(i)} = y$.

\medskip

$(2)\implies (1)$ Consider such permutations $\sigma_1,\ldots,\sigma_K$. If $i$ is fixed then for each $j$ there is at most one index $x$ such that:
$$\sigma_x(i) = j$$

We set $L_{ij} = x$ in that case and $L_{ij} = *$ otherwise. Then, $L$ is a sparse Latin square and the associated permutations are $\sigma_1^{-1},\ldots,\sigma_K^{-1}$, which belong to $G$.
\end{proof}

It turns out that as soon as $L_{N,K}^G\neq \emptyset$, the universal quasi-flat model is stationary. In order to prove this, let us first formulate another basic observation, as follows:

\begin{proposition}
Given $G\subset S_N$, we have an action $G\curvearrowright L_{N,K}^G$, given by
$$(L^\tau)_{ij}=L_{\tau^{-1}(i)j}$$
\end{proposition}

\begin{proof}
Given a sparse Latin square $L\in M_N(*,1,\ldots,K)$, consider the matrix in the statement, $L^\tau\in M_N(*,1,\ldots,K)$. This matrix is a sparse Latin square too and $\tau\to (L\mapsto L^\tau)$ is a group morphism, so it remains to check that:
$$G_L\subset G\implies G_{L^\tau}\subset G$$

For this purpose, let us write our groups as follows:
$$G_L=<\sigma_1,\ldots,\sigma_K>\quad,\quad G_{L^\tau}=<\sigma_1',\ldots,\sigma_K'>$$

We have then the following computation:
\begin{eqnarray*}
\sigma'_x(j)=i
&\iff&(L^\tau)_{ij}=x\\
&\iff&L_{\tau^{-1}(i)j}=x\\
&\iff&\sigma_x(j)=\tau^{-1}(i)\\
&\iff&\tau\sigma_x(j)=i
\end{eqnarray*}

Thus, we have $\sigma_x'=\tau\sigma_x$, and it follows that we have:
\begin{eqnarray*}
G_L\subset G
&\iff&\sigma_1,\ldots,\sigma_K\in G\\
&\iff&\tau\sigma_1,\ldots,\tau\sigma_K\in G\\
&\iff&\sigma_1',\ldots,\sigma_K'\in G\\
&\iff&G_{L^\tau}\subset G
\end{eqnarray*}

Thus $G_L\subset G$ implies $G_{L^\tau}\subset G$, and this finishes the proof.
\end{proof}

As a consequence of the above study, we have the following result, from \cite{bfr}:

\begin{theorem}
If the space $L_{N,K}^G$ is not empty, then the universal flat model for a quasi-transitive subgroup $G\subset S_N$,
$$\pi:C(G)\to M_K(C(T_G))$$
is stationary with respect to $\nu\otimes m$, where $m$ is the normalized counting measure on $L_{N,K}^G$, and $\nu$ is any probability measure on $E_K$.
\end{theorem}

\begin{proof}
We use the following standard identification:
$$M_K(C(T_G)) = M_{K}(\mathbb C)\otimes C(E_K)\otimes C(L_{N,K}^G)$$

We consider, for $P\in E_K$, the following map:
$$\pi^P = (id\otimes ev_P\otimes id) : C(G)\to M_K(C(L_{N,K}^G))$$

Recall that for $\tau\in G$ and $f\in C(G)$, $\tau.f$ denotes the map $h\mapsto f(\tau^{-1}h)$. This corresponds to the regular action of $G$ on itself. Moreover:
$$\pi^P(u_{ij})(\tau^{-1}L) = \pi^P(u_{\tau(i)j}) = \pi^P(\tau.u_{ij})$$

Since the normalized counting measure $m$ on $L_{N,K}^G$ is $G$-invariant, it follows, writing again $m$ for the integration with respect to $m$ on $C(L_{N, K}^G)$, that $(tr\otimes m)\pi^P$ is a $G$-invariant state on $C(G)$, so it is equal to $\int_G$. Summarizing, we have proved that $(tr\otimes id\otimes m)\pi$ is the constant function equal to $\int_G$, and the result follows.
\end{proof}

Let us discuss as well the group dual case. We start with $\Gamma=\mathbb Z_N$, where the class of $1$ will be denoted by $g$, and called the canonical generator. We have:

\begin{proposition}
We have an isomorphism of Hopf algebras
$$C^*(\mathbb Z_N)=C(S_N^+)\Big/\left<u_{ij}=u_{kl}\Big|j-i=l-k\ ({\rm mod}\ N)\right>$$
which is such that $g=\sum_{j=1}^Nw^{j-i}u_{ij}$, for any $i$, where $w=e^{2\pi i/N}$.
\end{proposition}

\begin{proof}
The quotient algebra $A$ in the statement being generated by the entries of the first row of $u=(u_{ij})$, it is commutative. If we identify the elements of $\widehat{\mathbb Z}_N$ with the powers of $w$, and the elements of $\mathbb Z_N$ with the functions $w\to w^k$, then $\pi(u_{ij})=\delta_{w^{j-i}}$ defines an isomorphism between $A$ and the algebra $C^*(\mathbb Z_N) = C(\widehat{\mathbb Z}_N)$. Moreover:
$$(\pi\otimes \pi)\Delta(u_{ij})=\sum_k\delta_{w^{k-i}}\otimes \delta_{w^{j-k}} = \Delta_{\widehat{\mathbb Z}_K}(\delta_{w^{j-i}})$$

Thus we have an isomorphism of Hopf algebras, as stated. The formula for $g$ in the statement then follows from the formula for $\pi$.
\end{proof}

This translates into a convenient description of the universal flat matrix model:

\begin{proposition}
The universal flat model space for $C^*(\mathbb Z_N)$ is the space
$$E_N=\left\{(P_1,\ldots,P_N)\in (P^{N-1}_\mathbb C)^N\Big|P_i\perp P_j,\forall i\neq j\right\}$$
with the model map given by $\pi(g)(P_1,\ldots,P_N)=\sum_{j=1}^Nw^j P_j$.
\end{proposition}

\begin{proof}
According to Proposition 14.13, the following is a flat matrix model map, which gives the formula in the statement for $\pi(g)$:
$$\pi(u_{ij})(P_1,\ldots,P_N)=P_{j-i}$$

Let $T$ be a flat matrix model space for $C^*(\mathbb Z_N)$, with model map $\pi_T$. For each $x\in T$ the matrix $\pi_T(g)(x)$ is circulant, so this matrix is completely determined by its first row, which is an element of $E_N$. Let us set:
$$P_j^x=P_{Nj}^x$$

Let also $P^x$ be the circulant matrix with first row $(P_1^x,\ldots,P_N^x)$, and consider the following map:
$$\Phi:T\to E_N\quad,\quad \Phi(x)=P^x$$

This map being a continuous embedding, we get a surjective morphism as follows, through which $\pi_T$ factors by construction:
$$\Psi:M_N(C(E_N))\to M_N(C(T))$$

Thus, we obtain the universality of $E_N$, as claimed.
\end{proof}

Regarding now the general group dual case, we have here the following result:

\begin{theorem}
Given a quotient group $\mathbb Z_K^{* M}\to\Gamma$, the universal quasi-flat model $\pi:C^*(\Gamma)\to M_K(C(Y_\Gamma))$ appears as follows:
\begin{enumerate}

\item For $\Gamma=\mathbb Z_K^{* M}$, the model space is $Y_{\Gamma} = Y_{N,K}=\underbrace{E_K\times\ldots\times E_K}_{M\ terms}$, and the model map is given by the formula in Proposition 14.14, on each of the components.

\item In general, the model space is the subspace $Y_{\Gamma}\subset Y_{N,K}$ consisting of the elements $x$ such that $\pi(\gamma)(x)=Id$ for any $\gamma\in\ker(\mathbb Z_K^{*M}\to\Gamma)$.
\end{enumerate}
\end{theorem}

\begin{proof}
The first assertion is clear from Proposition 14.14. Regarding now the second assertion, let us go back to Theorem 14.5 (2). The statement there tells us that the model space for $\widehat{\Gamma}$ appears from the model space for $\underbrace{S_K^+\,\hat{*}\,\ldots\,\hat{*}\,S_K^+}_{M\ terms}$ via the Tannakian conditions defining $\widehat{\Gamma}$. But these Tannakian conditions, when expressed in terms of the generators of $\mathbb Z_K^{* M}$ are precisely the relations $\gamma=1$ with $\gamma\in\ker(\mathbb Z_K^{*M}\to\Gamma)$, as stated.
\end{proof}

As a final topic now, let us discuss the following special types of models:
$$\pi:C(S_K^+)\to M_2(\mathbb C)$$

In order to to this, let us that a model $\pi:C(G)\to M_N(\mathbb C)$ for a quantum permutation group $G\subset S_K^+$ is called flat if $N=K$, and the projections $p_{ij}=\pi(u_{ij})$ are all of rank 1, and quasi-flat is the projections $p_{ij}=\pi(u_{ij})$ are all of rank 0 or 1. Up to a simple manipulation, we are in this latter situation, of quasi-flat models:

\begin{proposition}
Up to substracting an identity diagonal block, which will not change the image of the model, any model of type
$$\pi:C(S_K^+)\to M_2(\mathbb C)$$
can be assumed to be quasi-flat, in the sense that the projections $p_{ij}=\pi(u_{ij})$, belonging to the matrix algebra $M_2(\mathbb C)$, are all of rank $0$ or $1$.
\end{proposition} 

\begin{proof}
We know that the matrix $p=(p_{ij})$ formed by the projections $p_{ij}=\pi(u_{ij})$ is a $K\times K$ magic matrix having as entries projections $r\in M_2(\mathbb C)$. The point now is that these latter projections  $r\in M_2(\mathbb C)$ can be of rank $0,1,2$, and up to a permutation of the indices, we can isolate the rank 2 projections in a matrix block, so that our magic matrix $p=(p_{ij})$ looks as follows, with only rank $0,1$ projections in the lower block $q$:
$$p=\begin{pmatrix}1&0\\0&q\end{pmatrix}$$

Thus, we are led to the conclusion in the statement.
\end{proof}

In view of the above result, we can assume that our model is quasi-flat. In what regards now the models of $S_K^+$ that we are interested in, we have:

\begin{theorem}
The matrix models of the following type, assumed to be quasi-flat, up to substracting an identity matrix block,
$$\pi:C(S_K^+)\to M_2(\mathbb C)$$
appear, via permutations of the rows and columns, in the following way,
$$\pi:C(S_K^+)\to C(S_{K_1}\,\hat{*}\,\ldots\,\hat{*}S_{K_n})\to M_2(\mathbb C)$$
by making a direct sum of sparse Latin square models.
\end{theorem} 

\begin{proof}
We know that the matrix $p=(p_{ij})$ formed by the projections $p_{ij}=\pi(u_{ij})$ is a $K\times K$ magic matrix having as entries projections $r\in M_2(\mathbb C)$, having rank 0 or 1. We must prove that, up to a permutation of the rows and columns, this matrix becomes block-diagonal, with commuting entries in each block. But this can be done as follows:

\medskip

(1) Up to a permutation of the indices, we can assume that $r=p_{11}$ has rank 1. But then, in order for the sum on the first row to be 1, this first row must consist of $r$, or $s=1-r$, and of $0$ projections. The same goes for the first column. Now by permuting the rows and columns, in follows that our matrix must be as follows:
$$p=\begin{pmatrix}
r&s&0&\ldots\\
s&*\\
0\\
\vdots
\end{pmatrix}$$

If the $*$ entry is $r$, we have our first block, and we can proceed by recurrence.

\medskip

(2) If the $*$ entry is 0, we need a copy of $r$ both in the second row and second column, and so by permuting the rows and columns, our matrix must be as follows:
$$p=\begin{pmatrix}
r&s&0&0&\ldots\\
s&0&r&0&\ldots\\
0&r&*\\
0&0\\
\vdots&\vdots
\end{pmatrix}$$

Now if the $*$ entry is $s$, we have our first block, and we can proceed by recurrence.

\medskip

(3) And so on, and we end up with a block decomposition of our representation, as in the statement. As for the precise quantum group formulation of the factorization that we obtain, as in the statement, this follows from the various results above.
\end{proof}

Observe that the various manipulations in the above proof classify as well the possible sparse Latin squares $L\in M_K(*,1,2)$ that we are interested in, in relation with our present purposes. To be more precise, there is one such square for each $K\in\mathbb N$, as follows:
$$\begin{pmatrix}1&2\\2&1\end{pmatrix}\quad,\quad
\begin{pmatrix}1&2&*\\2&*&1\\ *&1&2\end{pmatrix}$$
$$\begin{pmatrix}1&2&*&*\\ 2&*&1&*\\ *&1&*&2\\ *&*&2&1\end{pmatrix}\quad,\quad
\begin{pmatrix}
1&2&*&*&*\\ 
2&*&1&*&*\\ 
*&1&*&2&*\\ 
*&*&2&*&1\\ 
*&*&*&1&2
\end{pmatrix}$$

In order to compute the image of the representation in Theorem 14.17, we can use Theorem 14.8. We are led in this way to the following conclusion:

\begin{proposition}
The matrix models of the following type, assumed to be quasi-flat, up to substracting an identity matrix block,
$$\pi:C(S_K^+)\to M_2(\mathbb C)$$
appear, via permutations of the rows and columns, by making a direct sum of sparse Latin square models, and factorize as follows,
$$\pi:C(S_K^+)\to C(G_{K_1}\,\hat{*}\,\ldots\,\hat{*}G_{K_n})\to M_2(\mathbb C)$$
with the groups $G_i\subset S_{K_i}$ being constructed from the sparse Latin squares, and with the algebra in the middle being the Hopf image.
\end{proposition} 

\begin{proof}
This follows indeed from the above results.
\end{proof}

Now recall from the above that the sparse Latin squares that we are interested in are classified, with one such square, that we will denote $L_K$, at any $K\geq2$. With the convention $L_1=(1)$, we can convert Proposition 14.18 into a final result, as follows:

\begin{theorem}
The matrix models of the following type,
$$\pi:C(S_K^+)\to M_2(\mathbb C)$$
appear, via permutations of the rows and columns, by making a direct sum of sparse Latin square models, and factorize as follows,
$$\pi:C(S_K^+)\to C(D_{K_1}\,\hat{*}\,\ldots\,\hat{*}D_{K_n})\to M_2(\mathbb C)$$
with $D_{K_i}$ being dihedral groups, and with the algebra in the middle being the Hopf image.
\end{theorem} 

\begin{proof}
We must prove that the group associated to the unique sparse Latin square $L_K\in M_K(*,1,2)$ is the dihedral group $D_N\subset S_N$. For this purpose, observe that the permutations $\sigma=\sigma_1$ and $\tau=\sigma_2$ are the following two involutions of $S_K$:
$$\sigma=(1)(23)(45)(67)\ldots\quad,\quad 
\tau=(12)(34)(56)(78)\ldots$$

Since $G_N=<\sigma,\tau>$ is generated by two involutions, we have:
$$G_N=<1,\sigma,\tau,\sigma\tau,\tau\sigma,\sigma\tau\sigma,\tau\sigma\tau,\ldots>$$

On the other hand since the permutation $\sigma\tau\in S_K$ is generically given by $i\to i\pm2$, with the sign depending on the parity of $i$, we have:
$$(\sigma\tau)^K=1$$

We conclude that the above sequence defining our group $G_N$ by enumerating its elements stops after $2N=|D_N|$ steps. With a bit more care, by using the quotient map $D_\infty=\mathbb Z_2*\mathbb Z_2\to G_N$, we conclude that we have $G_N=D_N$, as claimed.
\end{proof}

To summarize this study, we have a good understanding of the following models:
$$\pi:C(S_K^+)\to M_2(\mathbb C)$$

As a basic example, we have the ``sudoku'' construction. Given rank 1 projections satisfying $r+s=1$, $r'+s'=1$, consider the following matrix, which is sudoku magic:
$$p=\begin{pmatrix}
r&0&s&0\\
0&r'&0&s'\\
s&0&r&0\\
0&s'&0&r'
\end{pmatrix}$$

By permuting rows and columns, we obtain the following matrix, which is still magic, and which appears as above, from two $2\times2$ Latin squares:
$$\tilde{p}=\begin{pmatrix}
r&s&0&0\\
s&r&0&0\\
0&0&r'&s'\\
0&0&s'&r'
\end{pmatrix}$$

There are many interesting questions here, and with a bit more work, it is probably possible to deduce from this a classification of the models of type $\pi:C(S_K^+)\to M_2(\mathbb C)$.

\section*{14b. Flat models}

Following \cite{bne}, let us discuss now some more specialized questions, in the flat case. Assume that $B$ is a $C^*$-algebra, of finite dimension $\dim B=N<\infty$. We can endow $B$ with its canonical trace, $tr:B\subset\mathcal L(B)\to\mathbb C$, and use the following scalar product:
$$<a,b>=tr(ab^*)$$

We recall that, in terms of the decomposition $B=\oplus_sM_{n_s}(\mathbb C)$, we have $N=\sum_sn_s^2$, and the weights of the canonical trace are  $tr(I_s)=n_s^2/N$. We have:

\begin{definition}
A magic unitary $u\in M_N(\mathcal L(B))$ is called:
\begin{enumerate}
\item Flat, if each $u_{ij}\in\mathcal L(B)$ is a rank $1$ projection.

\item Split, if $u_{ij}=Proj(e_if_j^*)$, for certain sets $\{e_i\},\{f_i\}\subset U_B$.

\item Fully split, if $u_{ij}=Proj(g_ixg_j^*)$, with $\{g_i\}\subset U_B$, and $x\in U_B$.
\end{enumerate}
\end{definition}

We will see in what follows many interesting examples of such models. Let us first clarify the relation with the previous results for classical groups. We first have:

\begin{proposition}
The split magic unitaries which produce Latin squares are those of the form 
$$u_{ij}=Proj(g_ig_j^*)$$
with $\{g_1,\ldots,g_N\}\subset U_B$ being pairwise orthogonal, and forming a group $G\subset PU_B$. For such a magic unitary, the associated Latin square is $L_G$.
\end{proposition}

\begin{proof}
Assume indeed that $u_{ij}=Proj(e_if_j^*)$ produces a Latin square.

\medskip

(1) Our first claim is that we can assume $e_1=f_1=1$. Indeed, given $x,y\in U_B$ the following matrix is still magic:
$$u_{ij}'=Proj(xe_if_j^*y)$$ 

In the case where $u$ comes from a Latin square, $u_{ij}=Proj(\xi_{L_{ij}})$, we have:
$$u_{ij}'=Proj(\xi_{L_{ij}}')\quad,\quad \xi'_{ab}=x\xi_{ab}y$$

Thus $u'$ comes from $L$ as well, so with $x=e_1^*,y=f_1$, we can assume $e_1=f_1=1$.

\medskip

(2) Our second claim is that we can assume $u_{ij}=Proj(e_ie_j^*)$. Indeed, since $u$ is magic, the first row of vectors $\{1,f_2^*,\ldots,f_N^*\}\subset PU_B$ must appear as a permutation of the first column of vectors $\{1,e_2,\ldots,e_N\}\subset PU_B$. Thus, up to a permutation of columns, and a rescaling of the columns by elements in $Z(U_B)$, we can assume $f_i=e_i$, and we obtain:
$$u_{ij}=Proj(e_ie_j^*)$$

Observe that this permutation/rescaling of the columns won't change the fact that the associated Latin square $L$ comes or not from a group.

\medskip

(3) Let us construct now $G$. The Latin square condition shows that for any $i,j$ there is a unique $k$ such that $e_ie_j=e_k$ inside $PU_B$, and our claim is that the operation $(i,j)\to k$ gives a group structure on the set of indices. Indeed, all the group axioms are clear from definitions, and we obtain in this way a subgroup $G\subset PU_B$, having order $N$. 

\medskip

(4) With the group $G$ being constructed as above, we have:
$$u_{ij}=Proj(e_ie_j^*)=Proj(e_{ij^{-1}})$$

Thus we have $u_{ij}'=Proj(\xi_{L_{ij}})$ with $\xi_k=e_k$ and $L_{ij}=ij^{-1}$, and we are done.
\end{proof}

In order to further process the above result, we will need:

\begin{definition}
A $2$-cocycle on a group $G$ is a function $\sigma:G\times G\to\mathbb T$ satisfying:
$$\sigma(gh,k)\sigma(g,h)=\sigma(g,hk)\sigma(h,k)$$
$$\sigma(g,1)=\sigma(1,g)=1$$
The algebra $C^*(G)$, with multiplication $g\cdot h=\sigma(g,h)gh$, is denoted $C^*_\sigma(G)$.
\end{definition}

Observe that $g\cdot h=\sigma(g,h)gh$ is associative, and that we have $g\cdot 1=1\cdot g=g$, due to the 2-cocycle condition. Thus $C^*_\sigma(G)$ is an associative algebra with unit 1. In fact, $C^*_\sigma(G)$ is a $C^*$-algebra, with the involution making the canonical generators $g\in C^*_\sigma(G)$ unitaries. Also, the canonical trace on $C^*_\sigma(G)$ coincides with that of $C^*(G)$. 

\bigskip

With these conventions, we have the following result, from \cite{bne}:

\begin{proposition}
The split magic unitaries which produce Latin squares are precisely those of the form 
$$u_{ij}=Proj(g_ig_j^*)$$
with $\{g_1,\ldots,g_N\}$ being the standard basis of a twisted group algebra $C^*_\sigma(G)$. In this case, the associated Latin square is $L_G$.
\end{proposition}

\begin{proof}
We use Proposition 14.21. With the notations there, $\{g_1,\ldots,g_N\}\subset U_B$ must form a group $G\subset PU_B$, and so there are scalars $\sigma(i,j)\in\mathbb T$ such that:
$$g_ig_j=\sigma(i,j)g_{ij}$$

It follows from definitions that $\sigma$ is a 2-cocycle, and our claim now is that we have $B=C^*_\sigma(G)$. Indeed, this is clear when $\sigma=1$, because by linear independence we can define a linear space isomorphism $B\simeq C^*(G)$, which follows to be a $C^*$-algebra isomorphism. In the general case, where $\sigma$ is arbitrary, the proof is similar.
\end{proof}

We have now all the needed ingredients for proving a key result, as follows:

\begin{theorem}
Given a $2$-cocycle $\sigma:G\times G\to\mathbb T$, we have a representation
$$\pi:C(S_N^+)\to C(U_B,\mathcal L(B))\quad:\quad w_{ij}\to[x\to Proj(g_ixg_j^*)]$$
where $\{g_1,\ldots,g_N\}\subset U_B$ is the standard basis of the algebra $B=C^*_\sigma(G)$.
\end{theorem}

\begin{proof}
This follows indeed from the above results.
\end{proof}

There are many interesting examples of models as above, and we will study them later. Still at the general level now, we have the following result, also from \cite{bne}:

\begin{proposition}
For a representation coming from a split matrix
$$u_{ij}=Proj(e_if_j^*)$$
the truncated measure $\mu^r$ is the law of the Gram matrix of the vectors
$$\xi_{i_1\ldots i_r}=e_{i_1}f_{i_2}^*\otimes e_{i_2}f_{i_3}^*\otimes\ldots\ldots\otimes e_{i_r}f_{i_1}^*$$
with respect to the normalized trace of the $N^r\times N^r$ matrices.
\end{proposition}

\begin{proof}
According to the formulae in chapter 13, the moments of $\mu^r$ are given by:
\begin{eqnarray*}
c_p^r
&=&\sum_{i_1^1\ldots i_p^r}(T_p)_{i_1^1\ldots i_p^1,i_1^2\ldots i_p^2}(T_p)_{i_1^2\ldots i_p^2,i_1^3\ldots i_p^3}\ldots\ldots(T_p)_{i_1^r\ldots i_p^r,i_1^1\ldots i_p^1}\\
&=&\sum_{i_1^1\ldots i_p^r}tr(u_{i_1^1i_1^2}\ldots u_{i_p^1i_p^2})tr(u_{i_1^2i_1^3}\ldots u_{i_p^2i_p^3})\ldots\ldots tr(u_{i_1^ri_1^1}\ldots u_{i_p^ri_p^1})
\end{eqnarray*}

In the case of a split magic unitary, $u_{ij}=Proj(e_if_j^*)$, since the vectors $e_if_j^*$ are all of norm 1, with respect to the canonical scalar product, we therefore obtain:
\begin{eqnarray*}
c_p^r
&=&\frac{1}{N^r}\sum_{i_1^1\ldots i_p^r}<e_{i_1^1}f_{i_1^2}^*,e_{i_2^1}f_{i_2^2}^*>\ldots<e_{i_p^1}f_{i_p^2}^*,e_{i_1^1}f_{i_1^2}^*>\\
&&\hskip15.6mm<e_{i_1^2}f_{i_1^3}^*,e_{i_2^2}f_{i_2^3}^*>\ldots<e_{i_p^2}f_{i_p^3}^*,e_{i_1^2}f_{i_1^3}^*>\\
&&\hskip43mm\ldots\ldots\\
&&\hskip15.6mm<e_{i_1^r}f_{i_1^1}^*,e_{i_2^r}f_{i_2^1}^*>\ldots<e_{i_p^r}f_{i_p^1}^*,e_{i_1^r}f_{i_1^1}^*>
\end{eqnarray*}

Now by changing the order of the terms in the product, this gives:
\begin{eqnarray*}
c_p^r
&=&\frac{1}{N^r}\sum_{i_1^1\ldots i_p^r}<e_{i_1^1}f_{i_1^2}^*,e_{i_2^1}f_{i_2^2}^*><e_{i_1^2}f_{i_1^3}^*,e_{i_2^2}f_{i_2^3}^*>\ldots<e_{i_1^r}f_{i_1^1}^*,e_{i_2^r}f_{i_2^1}^*>\\
&&\hskip43mm\ldots\ldots\\
&&\hskip15.6mm<e_{i_p^1}f_{i_p^2}^*,e_{i_1^1}f_{i_1^2}^*><e_{i_p^2}f_{i_p^3}^*,e_{i_1^2}f_{i_1^3}^*>\ldots<e_{i_p^r}f_{i_p^1}^*,e_{i_1^r}f_{i_1^1}^*>
\end{eqnarray*}

In terms of the vectors $\xi_{i_1\ldots i_r}=e_{i_1}f_{i_2}^*\otimes\ldots\otimes e_{i_r}f_{i_1}^*$ in the statement, and then of their Gram matrix $G_{i_1\ldots i_r,j_1\ldots j_r}=<\xi_{i_1\ldots i_r},\xi_{j_1\ldots j_r}>$, we obtain the following formula:
\begin{eqnarray*}
c_p^r&=&\frac{1}{N^r}\sum_{i_1^1\ldots i_p^r}<\xi_{i_1^1\ldots i_1^r},\xi_{i_2^1\ldots i_2^r}>\ldots\ldots<\xi_{i_p^1\ldots i_p^r},\xi_{i_1^1\ldots i_1^r}>\\
&=&\frac{1}{N^r}\sum_{i_1^1\ldots i_p^r}G_{i_1^1\ldots i_1^r,i_2^1\ldots i_2^r}\ldots\ldots G_{i_p^1\ldots i_p^r,i_1^1\ldots i_1^r}\\
&=&\frac{1}{N^r}Tr(G^p)=tr(G^p)
\end{eqnarray*}

But this gives the formula in the statement, and we are done.
\end{proof}

In the fully split case now, we have the following result:

\begin{theorem}
For a model coming from a fully split matrix,
$$u_{ij}=Proj(g_ixg_j^*)$$
the truncated measure $\mu^r$ is the law of the Gram matrix of the vectors
$$\xi_{i_1\ldots i_r}^{x_1\ldots x_r}=g_{i_1}x_1g_{i_2}^*\otimes g_{i_2}x_2g_{i_3}^*\otimes\ldots\ldots\otimes g_{i_r}x_rg_{i_1}^*$$
with respect to the usual integration over $M_{N^r}(C(U_B^r))$.
\end{theorem}

\begin{proof}
The idea is that the computations in the proof of Proposition 14.25 apply, with $e_i=g_ix$ and $f_i=g_i$, and with an integral $\int_{U_B^r}$ added. To be more precise, we can start with the same formula as there, stating that the moments of $\mu^r$ are given by:
$$c_p^r=\sum_{i_1^1\ldots i_p^r}tr(u_{i_1^1i_1^2}\ldots u_{i_p^1i_p^2})\ldots\ldots tr(u_{i_1^ri_1^1}\ldots u_{i_p^ri_p^1})$$

In the case of a fully split matrix, $u_{ij}=Proj(g_ixg_j^*)$, since the vectors $g_ixg_j^*$ are all of norm 1, we therefore obtain:
\begin{eqnarray*}
c_p^r
&=&\frac{1}{N^r}\sum_{i_1^1\ldots i_p^r}\int_{U_B}<g_{i_1^1}x_1g_{i_1^2}^*,g_{i_2^1}x_1g_{i_2^2}^*>\ldots<g_{i_p^1}x_1g_{i_p^2}^*,g_{i_1^1}x_1g_{i_1^2}^*>dx_1\\
&&\hskip58mm\ldots\ldots\\
&&\hskip15.6mm\int_{U_B}<g_{i_1^r}x_rg_{i_1^1}^*,g_{i_2^r}x_rg_{i_2^1}^*>\ldots<g_{i_p^r}x_rg_{i_p^1}^*,g_{i_1^r}x_rg_{i_1^1}^*>dx_r
\end{eqnarray*}

Now by changing the order of the terms in the product, this gives:
\begin{eqnarray*}
c_p^r
&=&\frac{1}{N^r}\sum_{i_1^1\ldots i_p^r}\int_{U_B^r}<g_{i_1^1}x_1g_{i_1^2}^*,g_{i_2^1}x_1g_{i_2^2}^*>\ldots<g_{i_1^r}x_rg_{i_1^1}^*,g_{i_2^r}x_rg_{i_2^1}^*>\\
&&\hskip58mm\ldots\ldots\\
&&\hskip23mm<g_{i_p^1}x_1g_{i_p^2}^*,g_{i_1^1}x_1g_{i_1^2}^*>\ldots<g_{i_p^r}x_rg_{i_p^1}^*,g_{i_1^r}x_rg_{i_1^1}^*>dx
\end{eqnarray*}

In terms of the vectors $\xi_{i_1\ldots i_r}^{x_1\ldots x_r}=g_{i_1}x_1g_{i_2}^*\otimes\ldots\ldots\otimes g_{i_r}x_rg_{i_1}^*$  in the statement, and then of their Gram matrix $G_{i_1\ldots i_r,j_1\ldots j_r}^{x_1\ldots x_r}=<\xi_{i_1\ldots i_r}^{x_1\ldots x_r},\xi_{j_1\ldots j_r}^{x_1\ldots x_r}>$, we therefore obtain:
\begin{eqnarray*}
c_p^r&=&\frac{1}{N^r}\int_{U_B^r}\sum_{i_1^1\ldots i_p^r}<\xi_{i_1^1\ldots i_1^r}^{x_1\ldots x_r},\xi_{i_2^1\ldots i_2^r}^{x_1\ldots x_r}>\ldots\ldots<\xi_{i_p^1\ldots i_p^r}^{x_1\ldots x_r},\xi_{i_1^1\ldots i_1^r}^{x_1\ldots x_r}>dx\\
&=&\frac{1}{N^r}\int_{U_B^r}\sum_{i_1^1\ldots i_p^r}G_{i_1^1\ldots i_1^r,i_2^1\ldots i_2^r}^{x_1\ldots x_r}\ldots\ldots G_{i_p^1\ldots i_p^r,i_1^1\ldots i_1^r}^{x_1\ldots x_r}dx\\
&=&\frac{1}{N^r}\int_{U_B^r}Tr((G^{x_1\ldots x_r})^p)dx\\
&=&\int_{U_B^r}tr((G^{x_1\ldots x_r})^p)dx
\end{eqnarray*}

But this gives the formula in the statement, and we are done.
\end{proof}

We will see later explicit applications of the above formula. Observe also, independently of the mathematics in this book, the obvious similarity with the computations with transfer matrices from statistical mechanics. For more on this, we refer to \cite{ba5}.

\section*{14c. Gram matrices}

Following now \cite{bch}, \cite{bne} and related papers, we have the following result, regarding the matrices $T_p$ introduced in chapter 13, in relation with stationarity:

\begin{proposition}
Assuming that a model $\pi:C(G)\to M_K(C(X))$ is inner faithful and quasi-flat, mapping $u_{ij}\to Proj(\xi_{ij}^x)$, with $||\xi_{ij}^x||\in\{0,1\}$, we have
$$T_p=\int_XT_p(\xi^x)dx$$
where the matrix $T_p(\xi)\in M_{N^p}(\mathbb C)$, associated to an array $\xi\in M_N(\mathbb C^K)$ is given by: 
$$T_p(\xi)_{i_1\ldots i_p,j_1\ldots j_p}=\frac{1}{K}<\xi_{i_1j_1},\xi_{i_pj_p}><\xi_{i_pj_p},\xi_{i_{p-1}j_{p-1}}>\ldots\ldots<\xi_{i_2j_2},\xi_{i_1j_1}>$$
\end{proposition}

\begin{proof}
We have the following well-known computation, valid for any vectors $\xi_1,\ldots,\xi_p$ having norms $||\xi_i||\in\{0,1\}$, with the scalar product being as usual linear at left:
\begin{eqnarray*}
&&Proj(\xi_i)x=<x,\xi_i>\xi_i,\forall i\\
&\implies&
Proj(\xi_1)\ldots Proj(\xi_p)(x)=<x,\xi_p><\xi_p,\xi_{p-1}>\ldots\ldots<\xi_2,\xi_1>\xi_1\\
&\implies&Tr(Proj(\xi_1)\ldots Proj(\xi_p))=<\xi_1,\xi_p><\xi_p,\xi_{p-1}>\ldots\ldots
<\xi_2,\xi_1>
\end{eqnarray*}

Thus, the matrices $T_p$ can be computed as follows:
\begin{eqnarray*}
(T_p)_{i_1\ldots i_p,j_1\ldots j_p}
&=&\int_Xtr\left(Proj(\xi_{i_1j_1}^x)Proj(\xi_{i_2j_2}^x)\ldots Proj(\xi_{i_pj_p}^x)\right)dx\\
&=&\frac{1}{K}\int_X<\xi_{i_1j_1}^x,\xi_{i_pj_p}^x><\xi_{i_pj_p}^x,\xi_{i_{p-1}j_{p-1}}^x>\ldots\ldots<\xi_{i_2j_2}^x,\xi_{i_1j_1}^x>dx\\
&=&\int_X(T_p(\xi^x))_{i_1\ldots i_p,j_1\ldots j_p}dx
\end{eqnarray*}

We therefore obtain the formula in the statement. See \cite{bch}, \cite{bne}.
\end{proof}

An even more conceptual result, from \cite{bb3}, \cite{bch}, \cite{bne}, is as follows:

\index{Gram matrix}
\index{inner faithful model}

\begin{theorem}
Given an inner faithful quasi-flat model 
$$\pi:C(G)\to M_K(C(X))\quad,\quad 
u_{ij}\to Proj(\xi_{ij}^x)$$
with $||\xi_{ij}^x||\in\{0,1\}$, the law of the normalized character $\chi/K$ with respect to the truncated integral $\int_G^r$ coincides with that of the Gram matrix of the vectors
$$\xi_{i_1\ldots i_r}^x=\frac{1}{\sqrt{K}}\cdot\xi^{x_1}_{i_1i_2}\otimes\xi^{x_2}_{i_2i_3}\otimes\ldots\otimes \xi^{x_r}_{i_ri_1}$$
with respect to the normalized matrix trace, and to the integration functional on $X^r$.
\end{theorem}

\begin{proof}
The moments $C_p$ of the measure that we are interested in are given by:
\begin{eqnarray*}
C_p
&=&\frac{1}{K^p}\int_G^r\left(\sum_iu_{ii}\right)^p\\
&=&\frac{1}{K^p}\sum_{i_1\ldots i_p}(T_p^r)_{i_1\ldots i_p,i_1\ldots i_p}\\
&=&\frac{1}{K^p}\cdot Tr(T_p^r)
\end{eqnarray*}

The trace on the right is given by the following formula:
$$Tr(T_p^r)=\sum_{i_1^1\ldots i_p^r}(T_p)_{i_1^1\ldots i_p^1,i_1^2\ldots i_p^2}\ldots\ldots(T_p)_{i_1^r\ldots i_p^r,i_1^1\ldots i_p^1}$$

In view of the formula in Proposition 14.27, this quantity will expand in terms of the matrices $T_p(\xi)$ constructed there. To be more precise, we have:
$$Tr(T_p^r)=\int_{X^r}\sum_{i_1^1\ldots i_p^r}T_p(\xi^{x_1})_{i_1^1\ldots i_p^1,i_1^2\ldots i_p^2}\ldots\ldots T_p(\xi^{x_r})_{i_1^r\ldots i_p^r,i_1^1\ldots i_p^1}\,dx$$

By using now the explicit formula of each $T_p(\xi)$, we have:
\begin{eqnarray*}
Tr(T_p^r)
&=&\frac{1}{K^r}\int_{X^r}\sum_{i_1^1\ldots i_p^r}<\xi_{i_1^1i_1^2}^{x_1},\xi_{i_p^1i_p^2}^{x_1}>\ldots\ldots<\xi_{i_2^1i_2^2}^{x_1},\xi_{i_1^1i_1^2}^{x_1}>\\
&&\hskip52.8mm\ldots\\
&&\hskip22.4mm<\xi_{i_1^ri_1^1}^{x_r},\xi_{i_p^ri_p^1}^{x_r}>\ldots\ldots<\xi_{i_2^ri_2^1}^{x_r},\xi_{i_1^ri_1^1}^{x_r}>dx
\end{eqnarray*}

By changing the order of the summation, we can write this formula as:
\begin{eqnarray*}
Tr(T_p^r)
&=&\frac{1}{K^r}\int_{X^r}\sum_{i_1^1\ldots i_p^r}<\xi_{i_1^1i_1^2}^{x_1},\xi_{i_p^1i_p^2}^{x_1}>\ldots\ldots<\xi_{i_1^ri_1^1}^{x_r},\xi_{i_p^ri_p^1}^{x_r}>\\
&&\hskip52.8mm\ldots\\
&&\hskip19.4mm<\xi_{i_2^1i_2^2}^{x_1},\xi_{i_1^1i_1^2}^{x_1}>\ldots\ldots<\xi_{i_2^ri_2^1}^{x_r},\xi_{i_1^ri_1^1}^{x_r}>dx
\end{eqnarray*}

But this latter formula can be written as follows:
\begin{eqnarray*}
Tr(T_p^r)
&=&K^{p-r}\int_{X^r}\sum_{i_1^1\ldots i_p^r}\frac{1}{K}<\xi_{i_1^1i_1^2}^{x_1}\otimes\ldots\otimes\xi_{i_1^ri_1^1}^{x_r}\ ,\ \xi_{i_p^1i_p^2}^{x_1}\otimes\ldots\otimes \xi_{i_p^ri_p^1}^{x_r}>\\
&&\hskip55mm\ldots\\
&&\hskip15mm\frac{1}{K}<\xi_{i_2^1i_2^2}^{x_1}\otimes\ldots\otimes \xi_{i_2^ri_2^1}^{x_r}\ ,\ \xi_{i_1^1i_1^2}^{x_1}\otimes\ldots\otimes\xi_{i_1^ri_1^1}^{x_r}>dx
\end{eqnarray*}

In terms of the vectors in the statement, and of their Gram matrix $G_r^x$, we obtain:
\begin{eqnarray*}
Tr(T_p^r)
&=&K^{p-r}\int_{X^r}\sum_{i_1^1\ldots i_p^r}<\xi_{i_1^1\ldots i_1^r}^x,\xi_{i_p^1\ldots i_p^r}^x>\ldots\ldots<\xi_{i_2^1\ldots i_2^r}^x,\xi_{i_1^1\ldots i_1^r}^x>dx\\
&=&K^{p-r}\int_{X^r}\sum_{i_1^1\ldots i_p^r}(G_r^x)_{i_1^1\ldots i_1^r,i_p^1\ldots i_p^r}\ldots\ldots(G_r^x)_{i_2^1\ldots i_2^r,i_1^1\ldots i_1^r}\,dx\\
&=&K^{p-r}\int_{X^r}Tr((G_r^x)^p)dx
\end{eqnarray*}

Summarizing, the moments of the measure in the statement are given by:
$$C_p
=\frac{1}{K^r}\int_{X^r}Tr((G_r^x)^p)dx
=\left(tr\otimes\int_{X^r}\right)\left(G_r^p\right)$$

This gives the formula in the statement of the theorem.
\end{proof}

We will see later concrete applications of the above formulae, which are far less abstract than they might seem, in the context of concrete examples of matrix models. Also, as already mentioned, there is an obvious similarity here with the computations with transfer matrices from statistical mechanics. For more on this, we refer to \cite{ba5}.

\section*{14d. Universal models}

Following \cite{bne}, let us study now the universal flat model for $C(S_N^+)$. Given a flat magic unitary $u=(u_{ij})$, we can write it, in a non-unique way, as follows:
$$u_{ij}=Proj(\xi_{ij})$$

The array $\xi=(\xi_{ij})$ is then a ``magic basis'', in the sense that each of its rows and columns is an orthonormal basis of $\mathbb C^N$, and with this being an interesting combinatorial notion, mixing linear algebra and design theory. More on such magic bases later.

\bigskip

In relation now with universal models, we are led to two spaces, as follows:

\index{magic basis}
\index{flat magic unitary}

\begin{definition}
Associated to any $N\in\mathbb N$ are the following spaces:
\begin{enumerate}
\item $X_N$, the space of all $N\times N$ flat magic unitaries $u=(u_{ij})$.

\item $K_N$, the space of all $N\times N$ magic bases $\xi=(\xi_{ij})$.
\end{enumerate}
\end{definition}

Let us recall now that the rank 1 projections $p\in M_N(\mathbb C)$ can be identified with the corresponding 1-dimensional subspaces $E\subset\mathbb C^N$, which are by definition the elements of the complex projective space $P^{N-1}_\mathbb C$. In addition, if we consider the complex sphere, $S^{N-1}_\mathbb C=\{z\in\mathbb C^N|\sum_i|z_i|^2=1\}$, we have a quotient map as follows:
$$\pi:S^{N-1}_\mathbb C\to P^{N-1}_\mathbb C\quad,\quad z\to Proj(z)$$

Observe that we have $\pi(z)=\pi(z')$ precisely when $z'=wz$, for some $w\in\mathbb T$.  Consider as well the embedding $U_N\subset(S^{N-1}_\mathbb C)^N$ given by $x\to(x_1,\ldots,x_N)$, where $x_1,\ldots,x_N$ are the rows of $x$. Finally, as before, let us call an abstract matrix stochastic/bistochastic when the entries on each row/each row and column sum up to 1. 

\bigskip

With these notations and conventions, the abstract model spaces $X_N,K_N$ from Definition 14.29 that we are interested in, and some related spaces, are as follows:

\begin{proposition}
We have inclusions and surjections as follows,
$$\begin{matrix}
K_N&\subset&U_N^N&\subset&M_N(S^{N-1}_\mathbb C)\\
\\
\downarrow&&\downarrow&&\downarrow\\
\\
X_N&\subset&Y_N&\subset&M_N(P^{N-1}_\mathbb C)
\end{matrix}$$
where $X_N,Y_N$ consist of bistochastic/stochastic matrices, and $K_N$ is the lift of $X_N$.
\end{proposition}

\begin{proof}
This follows from the above discussion. Indeed, the quotient map $S^{N-1}_\mathbb C\to P^{N-1}_\mathbb C$ induces the quotient map $M_N(S^{N-1}_\mathbb C)\to M_N(P^{N-1}_\mathbb C)$ at right, and the lift of the space of stochastic matrices $Y_N\subset M_N(P^{N-1}_\mathbb C)$ is then the rescaled group $U_N^N$, as claimed.
\end{proof}

In order to get some insight into the structure of the above spaces $X_N,K_N$, we can use some inspiration from the well-known Sinkhorn algorithm from linear algebra. Indeed, this algorithm starts with a $N\times N$ matrix having positive entries and produces, via successive averagings over rows/columns, a bistochastic matrix.  

\bigskip

In our situation, we would like to have an ``averaging'' map $Y_N\to Y_N$, whose infinite iteration lands in the model space $X_N$. Equivalently, we would like to have an ``averaging'' map $U_N^N\to U_N^N$, whose infinite iteration lands in the space $K_N$. 

\bigskip

In order to construct such averaging maps, we use the orthogonalization procedure coming from the polar decomposition, the result that we need being as follows:

\index{Sinkhorn algorithm}

\begin{proposition}
We have orthogonalization maps as follows,
$$\xymatrix@R=10mm@C=15mm{
(S^{N-1}_\mathbb C)^N\ar[r]^\alpha\ar[d]&(S^{N-1}_\mathbb C)^N\ar[d]\\
(P^{N-1}_\mathbb C)^N\ar[r]^\beta&(P^{N-1}_\mathbb C)^N}$$
where $\alpha(x)_i=Pol([(x_i)_j]_{ij})$, and $\beta(p)=(P^{-1/2}p_iP^{-1/2})_i$, with $P=\sum_ip_i$.
\end{proposition}

\begin{proof}
This is something which is routine, the idea being as follows:

\medskip

(1) Our first claim is that we have a factorization as in the statement. Indeed, pick $p_1,\ldots,p_N\in P^{N-1}_\mathbb C$, and write $p_i=Proj(x_i)$, with $||x_i||=1$. We can then apply $\alpha$, as to obtain a vector $\alpha(x)=(x_i')_i$, and then set $\beta(p)=(p_i')$, where $p_i'=Proj(x_i')$.

\medskip

(2) Our first task is to prove that $\beta$ is well-defined. Consider indeed vectors $\tilde{x}_i$, satisfying $Proj(\tilde{x}_i)=Proj(x_i)$. We have then $\widetilde{x}_i=\lambda_ix_i$, for certain scalars $\lambda_i\in\mathbb T$, and so the matrix formed by these vectors is $\widetilde{M}=\Lambda M$, with $\Lambda=diag(\lambda_i)$. It follows that $Pol(\widetilde{M})=\Lambda Pol(M)$, and so $\tilde{x}_i'=\lambda_ix_i$, and finally $Proj(\tilde{x}_i')=Proj(x_i')$, as desired.

\medskip

(3) It remains to prove that $\beta$ is given by the formula in the statement. For this purpose, observe first that, given $x_1,\ldots,x_N\in S^{N-1}_\mathbb C$, with $p_i=Proj(x_i)$ we have:
\begin{eqnarray*}
\sum_ip_i
&=&\sum_i[(\bar{x}_i)_k(x_i)_l]_{kl}\\
&=&\sum_i(\bar{M}_{ik}M_{il})_{kl}\\
&=&((M^*M)_{kl})_{kl}\\
&=&M^*M
\end{eqnarray*}

(4) We can now compute the projections $p_i'=Proj(x_i')$. Indeed, the coefficients  of these projections are given by $(p_i')_{kl}=\bar{U}_{ik}U_{il}$ with $U=MP^{-1/2}$, and we obtain, as desired:
\begin{eqnarray*}
(p_i')_{kl}
&=&\sum_{ab}\bar{M}_{ia}P^{-1/2}_{ak}M_{ib}P^{-1/2}_{bl}\\
&=&\sum_{ab}P^{-1/2}_{ka}\bar{M}_{ia}M_{ib}P^{-1/2}_{bl}\\
&=&\sum_{ab}P^{-1/2}_{ka}(p_i)_{ab}P^{-1/2}_{bl}\\
&=&(P^{-1/2}p_iP^{-1/2})_{kl}
\end{eqnarray*}

(5) An alternative proof uses the fact that the elements $p_i'=P^{-1/2}p_iP^{-1/2}$ are self-adjoint, and sum up to 1. The fact that these elements are indeed idempotents can be checked directly, via $p_iP^{-1}p_i=p_i$, because this equality holds on $\ker p_i$, and also on $x_i$.
\end{proof}

As an illustration, here is how the orthogonalization works at $N=2$:

\begin{proposition}
At $N=2$ the orthogonalization procedure for 
$$(Proj(x),Proj(y))$$
amounts in considering the vectors 
$$\frac{x+y}{\sqrt{2}}\quad,\quad \frac{x-y}{\sqrt{2}}$$
 and then rotating by $45^\circ$.
\end{proposition}

\begin{proof}
By performing a rotation, we can restrict attention to the case $x=(\cos t,\sin t)$ and $y=(\cos t,-\sin t)$, with $t\in(0,\pi/2)$. Here the computations are as follows:
\begin{eqnarray*}
M=\begin{pmatrix}\cos t&\sin t\\ \cos t&-\sin t\end{pmatrix}
&\implies&
P=M^*M=\begin{pmatrix}2\cos^2t&0\\0&2\sin^2t\end{pmatrix}\\
&\implies&P^{-1/2}=|M|^{-1}=\frac{1}{\sqrt{2}}\begin{pmatrix}\frac{1}{\cos t}&0\\0&\frac{1}{\sin t}\end{pmatrix}\\
&\implies&U=M|M|^{-1}=\frac{1}{\sqrt{2}}\begin{pmatrix}1&1\\1&-1\end{pmatrix}
\end{eqnarray*}

Thus the orthogonalization procedure replaces $(Proj(x),Proj(y))$ by the orthogonal projections on the vectors $(\frac{1}{\sqrt{2}}(1,1),\frac{1}{\sqrt{2}}(-1,1))$, and this gives the result.
\end{proof}

With these preliminaries in hand, let us discuss now the version that we need of the Sinkhorn algorithm. The orthogonalization procedure is as follows:

\begin{theorem}
The orthogonalization maps $\alpha,\beta$ induce maps as follows,
$$\xymatrix@R=10mm@C=15mm{
U_N^N\ar[r]^\Phi\ar[d]&U_N^N\ar[d]\\
Y_N\ar[r]^\Psi&Y_N}$$
which are the transposition maps on $K_N,X_N$, and which are projections at $N=2$.
\end{theorem}

\begin{proof}
It follows from definitions that $\Phi(x)$ is obtained by putting the components of $x=(x_i)$ in a row, then picking the $j$-th column vectors of each $x_i$, calling $M_j$ this matrix, then taking the polar part $x_j'=Pol(M_j)$, and finally setting $\Phi(x)=x'$. Thus:
$$\Phi(x)=Pol((x_{ij})_i)_j$$
$$\Psi(u)=(P_i^{-1/2}u_{ji}P_i^{-1/2})_{ij}$$

Thus, the first assertion is clear, and the second assertion is clear too.
\end{proof}

Our claim is that the algorithm converges, as follows:

\begin{conjecture}
The above maps $\Phi,\Psi$ increase the volume,
$$vol:U_N^N\to Y_N\to[0,1],\quad vol(u)=\prod_j|\det((u_{ij})_i)|$$
and respectively land, after an infinite number of steps, in $K_N/X_N$.
\end{conjecture}

Observe that the quantities of type $|\det(p_1,\ldots,p_N)|$ are indeed well-defined, for any $p_1,\ldots,p_N\in P^{N-1}_\mathbb C$, because multiplying by scalars $\lambda_i\in\mathbb T$ doesn't change the volume. Thus, the volume map $vol:U_N^N\to[0,1]$ factorizes through $Y_N$, as stated above.

\bigskip

As a main application of the above conjecture, the infinite iteration $(\Phi^2)^\infty:U_N^N\to K_N$ would provide us with an integration on $K_N$, and hence on the quotient space $K_N\to X_N$ as well, by taking the push-forward measures, coming from the Haar measure on $U_N^N$. 

\bigskip

In relation now with the matrix model problematics, we have:

\index{universal model}

\begin{conjecture}
The universal $N\times N$ flat matrix representation
$$\pi_N:C(S_N^+)\to M_N(C(X_N)),\quad \pi_N(w_{ij})=(u\to u_{ij})$$
is inner faithful at any $N\geq4$.
\end{conjecture}

Regarding the $N=4$ conjecture, the problem here is that of proving, as in \cite{bc2}, that the composition $C(S_4^+)\to M_4(C(X_4))\to\mathbb C$ equals the Haar integration on $S_4^+$. As for the $N\geq 5$ conjecture, the problem here is that of proving that the truncated moments $c_p^r$ converge with $r\to\infty$ to the Catalan numbers. None of these questions is trivial.

\bigskip

Still following \cite{bne}, our purpose now will be to advance towards a unification of the two conjectures formulated above. The point indeed is that when trying to approach Conjecture 13.35 with standard probabilistic tools, the estimates that are needed seem to be related to those required for approaching Conjecture 13.34. 

\bigskip

We first have the following definition, inspired by the above results:

\begin{definition}
Associated to $x\in M_N(S^{N-1}_\mathbb C)$ is the $N^p\times N^p$ matrix
$$(T_p^x)_{i_1\ldots i_p,j_1\ldots j_p}=\frac{1}{N}<x_{i_1j_1},x_{i_pj_p}><x_{i_pj_p},x_{i_{p-1}j_{p-1}}>\ldots\ldots<x_{i_2j_2},x_{i_1j_1}>$$
where the scalar products are the usual ones on $S^{N-1}_\mathbb C\subset\mathbb C^N$.
\end{definition}

The first few values of these matrices, at $p=1,2,3$, are as follows:
\begin{eqnarray*}
(T_1^x)_{ia}&=&\frac{1}{N}<x_{ia},x_{ia}>=\frac{1}{N}\\
(T_2^x)_{ij,ab}&=&\frac{1}{N}<x_{ia},x_{jb}><x_{jb},x_{ia}>=\frac{1}{N}|<x_{ia},x_{jb}>|^2\\
(T_3^x)_{ijk,abc}&=&\frac{1}{N}<x_{ia},x_{kc}><x_{kc},x_{jb}><x_{jb},x_{ia}>
\end{eqnarray*}

The interest in these matrices, in connection with Conjecture 14.34, comes from:

\begin{proposition}
For the universal model, the matrices $T_p$ are
$$T_p=\int_{K_N}T_p^xdx\,$$
where $dx$ is the measure on the model space $K_N$ coming from Conjecture 14.34.
\end{proposition}

\begin{proof}
This is a trivial statement, because by definition of $T_p$, we have:
\begin{eqnarray*}
(T_p)_{i_1\ldots i_p,j_1\ldots j_p}
&=&tr(u_{i_1j_1}\ldots u_{i_pj_p})=\int_{K_N}tr(u_{i_1j_1}^x\ldots u_{i_pj_p}^x)dx\\
&=&\int_{K_N}tr(Proj(x_{i_1j_1})\ldots Proj(x_{i_px_p}))dx\\
&=&\frac{1}{N}\int_{K_N}<x_{i_1j_1},x_{i_pj_p}>\ldots\ldots<x_{i_2j_2},x_{i_1j_1}>dx\\
&=&\int_{K_N}(T_p^x)_{i_1\ldots i_p,j_1\ldots j_p}dx
\end{eqnarray*}

Thus the formula in the statement holds indeed.
\end{proof}

In fact, the matrices $T_p^x$ are related to Conjecture 14.35 as well. Indeed, to any noncrossing partition $\pi\in NC(p)$ let us associate the following vector of $(\mathbb C^N)^{\otimes p}$: 
$$\xi_\pi=\sum_{\ker i\leq\pi}e_{i_1}\otimes\ldots\otimes e_{i_p}$$

With this convention, we have the following result, again from \cite{bne}:

\begin{proposition}
For any $x\in M_N(S^{N-1}_\mathbb C)$, the following hold:
\begin{enumerate}
\item If $\{x_{ij}\}_i$ are pairwise orthogonal then $(T_p^x)^*\xi_{||\ldots|}=\xi_{||\ldots|}$ and $T_p^x\xi_{\sqcap\hskip-0.5mm\sqcap\ldots\sqcap}=\xi_{\sqcap\hskip-0.5mm\sqcap\ldots\sqcap}$.

\item If $\{x_{ij}\}_j$ are pairwise orthogonal then $T_p^x\xi_{||\ldots|}=\xi_{||\ldots|}$ and $(T_p^x)^*\xi_{\sqcap\hskip-0.5mm\sqcap\ldots\sqcap}=\xi_{\sqcap\hskip-0.5mm\sqcap\ldots\sqcap}$.

\item If $\{x_{ij}\}_i$ or $\{x_{ij}\}_j$ are pairwise orthogonal then $<T_p^x\xi_{||\ldots|},\xi_{||\ldots|}>=N^p$.

\item We have $<T_p^x\xi_{\sqcap\hskip-0.5mm\sqcap\ldots\sqcap},\xi_{\sqcap\hskip-0.5mm\sqcap\ldots\sqcap}>=N$, without assumptions on $x$.
\end{enumerate}
\end{proposition}

\begin{proof}
Assuming that $\{x_{ij}\}_i$ are pairwise orthogonal, we have indeed:
\begin{eqnarray*}
(T_p^x\xi_{\sqcap\hskip-0.5mm\sqcap\ldots\sqcap})_{i_1\ldots i_p}
&=&\sum_j(T_p^x)_{i_1\ldots i_p,j\ldots j}\\
&=&\frac{1}{N}\sum_j<x_{i_1j},x_{i_pj}>\ldots\ldots<x_{i_2j},x_{i_1j}>\\
&=&\delta_{i_1,\ldots,i_p}
\end{eqnarray*}

Thus we have proved (1), and the proof of (2,3,4) is similar. See \cite{bne}.
\end{proof}

We have the following statement, supported by computer calculations:

\begin{conjecture}
Consider the following function, with $x\in M_N(S^{N-1}_\mathbb C)$,
$$F_p(x)=\frac{1}{N^p}||T_p^x\xi_{\sqcap\hskip-0.5mm\sqcap\ldots\sqcap}||^2$$
depending on a fixed integer $p\geq2$. Then, for any $x\in U_N^N$ we have 
$$F_p(x)\geq F_p(\Psi^2(x))$$
with equality precisely when $x\in K_N$, in which case $F_p(x)=1$. 
\end{conjecture}

This conjecture is quite interesting, in relation with the above, because by a compacity argument, this would prove that our Sinkhorn type algorithm converges. Thus, what we have here is a first step towards unifying Conjectures 14.34 and 14.35. See \cite{bne}.

\section*{14e. Exercises}

Things have been quite technical in this chapter, and as exercise on all this, which is potentially quite elementary, and is instructive and interesting as well, we have:

\begin{exercise}
Fully classify the models of type $\pi:C(S_K^+)\to M_2(\mathbb C)$.
\end{exercise}

This was something discussed in the above, and the question is that of finishing the work started there, and with all this being related to interesting combinatorics.

\chapter{Weyl matrices}

\section*{15a. Pauli matrices}

We discuss in this chapter, following \cite{bc2}, \cite{bne} and related papers, a number of key matrix models for certain quantum permutation groups $G\subset S_N^+$, which are of particular importance. These models will be, for the most, stationary, and the whole material in this chapter will lie at the interface between quantum groups and random matrices.

\bigskip

Following \cite{bc2}, let us first discuss an explicit and useful matrix model for the quantum group $S_4^+$, coming from the Pauli matrices. These matrices are as follows:
$$c_1=\begin{pmatrix}1&0\\ 0&1\end{pmatrix}\quad,\quad
c_2=\begin{pmatrix}i&0\\ 0&-i\end{pmatrix}\quad,\quad
c_3=\begin{pmatrix}0&1\\ -1&0\end{pmatrix}\quad,\quad
c_4=\begin{pmatrix}0&i\\ i&0\end{pmatrix}$$

Observe that these matrices are elements of $SU_2$. In fact, $SU_2$ consists of linear combinations of Pauli matrices, with points on the real sphere $S^3$ as coefficients:
$$SU(2)=\left\{\sum x_ic_i\,\Big| \sum|x_i|^2=1\right\}$$

The Pauli matrices multiply according to the formulae for quaternions:
$$c_2^2=c_3^2=c_4^2=-1$$
$$c_2c_3=-c_3c_2=c_4$$
$$c_3c_4=-c_4c_3=c_2$$
$$c_4c_2=-c_2c_4=c_3$$

The starting remark is that the Pauli matrices form an orthonormal basis of $M_2(\mathbb C)$, with respect to the scalar product $<a,b>=tr(ab^*)$. Moreover, the same is true if we multiply them to the left or to the right by an element of $SU_2$. Thus, we have:

\begin{proposition}
For any $x\in SU_2$ the elements 
$$\xi_{ij}=c_ixc_j$$
form a magic basis of $M_2({\mathbb C})$, with respect to the scalar product $<a,b>=tr(ab^*)$.
\end{proposition}

\begin{proof}
This follows indeed from the above discussion.
\end{proof}

We fix a Hilbert space isomorphism $M_2({\mathbb C})\simeq {\mathbb C}^4$, and we use the corresponding identification of operator algebras $B(M_2({\mathbb C}))\simeq M_4({\mathbb C})$. Associated to each element $x\in SU_2$ is the representation of $C(S_4^+)$ mapping $u_{ij}$ to the rank one projection on $c_ixc_j$:
$$\pi_x:C(S_4^+)\to M_4(\mathbb C)\quad,\quad u_{ij}\to Proj(c_ixc_j)$$

This representation depends on $x$. For getting a faithful representation, the idea is to regard all these representations as fibers of a single representation:

\begin{definition}
The Pauli representation of $C(S_4^+)$ is the map
$$\pi:C(S_4^+)\to C(SU_2, M_4({\mathbb C}))$$
mapping each standard coodinate $u_{ij}$ to the function $x\to Proj(c_ixc_j)$.
\end{definition}

Observe that we can replace here $SU_2$ by its quotient $PU_2=SO_3$. For reasons that will become clear later on, we prefer to use $SU_2$. We will prove now that $\pi$ is faithful, in the usual sense. Consider the standard linear form on the algebra on the right:
$$\int\varphi=\int_{SU_2}tr\left(\varphi(x)\right)\, dx$$

We use the following analytic formulation of faithfulness:

\begin{proposition}
The representation $\pi$ is faithful provided that
$$\int u_{i_1j_1}\ldots u_{i_kj_k}=
\int \pi_{i_1j_1}\ldots \pi_{i_kj_k}$$ for any choice of $k$ and
of various $i,j$ indices, where $\pi_{ij}=\pi(u_{ij})$.
\end{proposition}

\begin{proof}
The condition in the statement tells us that, for any product $a$ of generators $u_{ij}$, the following equality must happen:
$$\int a=\int\pi(a)$$

By linearity and density this formula holds on the whole algebra $C(S_4^+)$. In other
words, the following diagram commutes:
$$\begin{matrix}
C(S_4^+)&\ &\displaystyle{\mathop{\longrightarrow}^ {\int}}&\
&{\mathbb C}\\
\ \\ 
\pi\downarrow&\ &\ &\ &\uparrow tr\\
\ \\
C\left(SU_2,M_4({\mathbb C})\right)&\ &
\displaystyle{\mathop{\longrightarrow}^ {\int}}&\ &M_4({\mathbb
C})
\end{matrix}$$

On the other hand, we know that the algebra $C(S_4^+)$ is amenable in the discrete quantum group sense, which means that its Haar functional is faithful:
$$a\neq 0\implies\int aa^*>0$$

Assume now that we have $\pi(a)=0$. This implies $\pi(aa^*)=0$, and the commutativity of the above diagram gives $\int aa^*=0$. Thus $a=0$, and we are done.
\end{proof}

Our first goal will be to compute the integral on the right in the above statement. For this purpose, consider the canonical action of $SU_2$ on the algebra $M_2(\mathbb C)^{\otimes k}$, obtained as $k$-th tensor power of the adjoint action on $M_2(\mathbb C)$, as follows:
$$\alpha_x(a_1\otimes\ldots\otimes a_k)=xa_1x^*\otimes\ldots\otimes xa_kx^*$$

The following map will play an important role in what follows:

\begin{definition}
We define a linear map $R:M_2(\mathbb C)^{\otimes k}\to M_2(\mathbb C)^{\otimes k}$ by
$$R(c_{i_1}\otimes\ldots\otimes  c_{i_k})=\frac{1}{2}(c_{i_1}c_{i_2}^*
\otimes c_{i_2}c_{i_3}^*\otimes\ldots\otimes c_{i_k}c_{i_1}^*)$$
with the convention that at $k=1$ we have $R(c_i)=c_ic_i^*/2=1/2$.
\end{definition}

To any multi-index $i=(i_1,\ldots ,i_k)$ we associate the following element:
$$c_i=c_{i_1}\otimes\ldots\otimes c_{i_k}$$

With these notations and conventions, we have the following result:

\begin{proposition}
We have the following formula,
$$\int \pi_{i_1j_1}\ldots \pi_{i_kj_k}=<c_i,R^*ER(c_j)>$$ 
where $E$ is the expectation under the canonical action of $SU_2$.
\end{proposition}

\begin{proof}
We have indeed the following computation, where $Proj(\xi)$ denotes as usual the rank one projection onto a vector $\xi$:
\begin{eqnarray*}
\int \pi_{i_1j_1}\ldots \pi_{i_kj_k} 
&=& \int\pi(u_{i_1j_1})\ldots \pi(u_{i_kj_k})\\
&=&\int Proj(c_{i_1}xc_{j_1})\ldots Proj(c_{i_k}xc_{j_k})\\
&=&\int_{SU_2}tr\left( Proj(c_{i_1}xc_{j_1})\ldots Proj(c_{i_k}xc_{j_k})\right)\,dx
\end{eqnarray*}

We use now the following elementary formula, that we met many times in the above, valid for any sequence of norm one vectors $\xi_1,\ldots,\xi_k$ in a Hilbert space:
$$Tr\left(Proj(\xi_1)\ldots Proj(\xi_k)\right)
=<\xi_1,\xi_k><\xi_k,\xi_{k-1}>\ldots <\xi_2,\xi_1>$$

In our situation these vectors are in fact matrices, and their scalar products are given by $<\xi,\eta>=tr(\xi\eta^*)$. This gives the following formula:
\begin{eqnarray*}
\int \pi_{i_1j_1}\ldots \pi_{i_kj_k} 
&=&\frac{1}{4}\int_{SU_2}<c_{i_1}xc_{j_1},c_{i_k}xc_{j_k}>\ldots
<c_{i_2}xc_{j_2},c_{i_1}xc_{j_1}>\,dx\\
&=&\frac{1}{4}\int_{SU_2}tr(c_{i_1}xc_{j_1}c_{j_k}^*x^*c_{i_k}^*)\ldots tr(c_{i_2}xc_{j_2}c_{j_1}^*x^*c_{i_1}^*)\,dx \\
&=&\frac{1}{4}\int_{SU_2}tr(c_{i_k}^*c_{i_1}xc_{j_1}c_{j_k}^*x^*)\ldots tr(c_{i_1}^*c_{i_2}xc_{j_2}c_{j_1}^*x^*)\,dx
\end{eqnarray*}

We use now the formula $c_s^*=\pm c_s$, valid for all the Pauli matrices $c_s$. The minus signs can be rearranged, and the computation can be continued as follows:
\begin{eqnarray*}
\int \pi_{i_1j_1}\ldots \pi_{i_kj_k} 
&=&\frac{1}{4}\int_{SU_2}tr(c_{i_k}c_{i_1}^*xc_{j_1}c_{j_k}^*x^*)\ldots 
tr(c_{i_1}c_{i_2}^*xc_{j_2}c_{j_1}^*x^*)\,dx\\
&=&\frac{1}{4}\int_{SU_2}tr\left(c_{i_1}c_{i_2}^*xc_{j_2}c_{j_1}^*x^*
\otimes\ldots\otimes c_{i_k}c_{i_1}^*xc_{j_1}c_{j_k}^*x^*\right)\,dx\\
&=&\int_{SU_2}tr\left( R(c_i)\alpha_x(R(c_j)^*)\right)\,dx
\end{eqnarray*}

We can interchange the trace and integral signs, and we obtain:
\begin{eqnarray*}
\int \pi_{i_1j_1}\ldots \pi_{i_kj_k}
&=&tr\left(\int_{SU_2}R(c_i)\alpha_x(R(c_j)^*)\, dx\right)\\
&=&tr\left(R(c_i)\int_{SU_2} \alpha_x(R(c_j)^*)\,dx\right)
\end{eqnarray*}

Now since acting by group elements, then integrating, is the same as projecting onto fixed points, the computation can be continued as follows:
\begin{eqnarray*}
\int \pi_{i_1j_1}\ldots \pi_{i_kj_k}
&=&tr\left(R(c_i)ER(c_j)^*\right)\\
&=&<R(c_i),ER(c_j)>\\
&=&<c_i,R^*ER(c_j)>
\end{eqnarray*}

Thus, we are led to the conclusion in the statement.
\end{proof}

We know that the faithfulness of the Pauli representation is equivalent to a certain equality of integrals. Moreover, one of these integrals can be computed by using the Weingarten formula. As for the other integral, this can be computed as well, provided that we have enough information about the operator $R$ from Definition 15.4. 

\bigskip

We work out now a number of technical properties of this operator $R$. For this purpose, we first need to understand some aspects of the structure of the algebra of fixed points under the diagonal adjoint action of $SU_2$. We have the following result, to start with:

\begin{proposition}
The element
$$f=\sum_{i=1}^4c_i\otimes c_i^*$$
is invariant under the action of $SU_2$.
\end{proposition}

\begin{proof}
We have $c_i=-c_i^*$ for $i=2,3,4$, so the element in the statement is:
$$f=2(c_1\otimes c_1)-\sum_{i=1}^4 c_i\otimes c_i$$

Since $c_1\otimes c_1$ is invariant under $SU_2$, what is left to prove is that the sum on the right is invariant as well. This sum, viewed as a matrix, is:
$$\sum_{i=1}^4c_i\otimes c_i
=\begin{pmatrix}
 0&0&0&0\\
0&2&-2&0\\
0&-2&2&0\\
0&0&0&0
\end{pmatrix}$$

But this matrix is $4$ times the projection onto the determinant subspace of 
$\mathbb C^{2}\otimes\mathbb C^2$, which is invariant under the action of $SU_2$, and this concludes the proof. 
\end{proof}

Given a noncrossing partition $p\in NC(k)$ and a multi-index $j=(j_1,\ldots,j_k)$, we can plug $j$ into $p$ in the obvious way, and we define a Kronecker type symbol $\delta_{pj}\in\{0,1\}$ by setting $\delta_{pj}=1$ if any block of $p$ contains identical indices of $j$, and $\delta_{pj}=0$ otherwise. To any $p\in NC(k)$ we associate an element of $M_2(\mathbb C)^k$, as follows:
$$c_p=\sum_j\delta_{pj}\,c_j$$

Our next goal is to prove that $R(c_p)$ is invariant under the action of $SU_2$. We discuss first the case of the trivial partition, $0_k=\{\{1\},\ldots ,\{k\}\}$. We have:

\begin{proposition}
$R(c_{0_k})$ is invariant under the action of $SU_2$.
\end{proposition}

\begin{proof}
The partition $p=0_k$ has the particular property that $\delta_{pj}=1$ for any multi-index $j=(j_1,\ldots ,j_k)$. This gives the following formula:
\begin{eqnarray*}
R(c_{0_k})
&=&R\left(\sum_jc_j\right)\\
&=&R\left(\sum_{j_1\ldots j_k}c_{j_1}\otimes\ldots\otimes c_{j_k}\right)\\
&=&\frac{1}{2}\sum_{j_1\ldots j_k} c_{j_1}c^*_{j_2}\otimes
c_{j_2}c^*_{j_3}\otimes\ldots\otimes c_{j_k}c_{j_1}^*
\end{eqnarray*}

In ordert to compute the above quantity, to any multi-index $j=(j_1,\ldots ,j_k)$ we associate a multi-index $i=(i_1,\ldots ,i_{k-1})$ in the following way: 

\smallskip

-- $i_1$ is such that $c_{i_1}=\pm c_{j_1}c_{j_2}^*$,

\smallskip

-- then $i_2$ is such that $c_{i_2}=\pm c_{j_1}c_{j_3}^*$, 

$\vdots$

-- and so on up to $i_{k-1}$, which is such that $c_{i_{k-1}}=\pm c_{j_1}c_{j_k}^*$.

\smallskip

With this notation, we have the following formulae, where the possible dependences between the various  $\pm$ signs are not taken into account:
\begin{eqnarray*}
c_{i_1}&=&\pm c_{j_1}c_{j_2}^*\\
c_{i_1}^*c_{i_2}&=&(\pm c_{j_1}c_{j_2}^*)^*(\pm c_{j_1}c_{j_3}^*)=\pm c_{j_2}c_{j_3}^*\\
c_{i_2}^*c_{i_3}&=&(\pm c_{j_1}c_{j_3}^*)^*(\pm c_{j_1}c_{j_4}^*)=\pm c_{j_3}c_{j_4}^*\\
&\vdots&\\
c_{i_{k-2}}^*c_{i_{k-1}}&=&(\pm c_{j_1}c_{j_{k-1}}^*)^*(\pm c_{j_1}c_{j_k}^*)=\pm c_{j_{k-1}}c_{j_k}^*\\
c_{i_{k-1}}^*&=&(\pm c_{j_1}c_{j_k}^*)^*=\pm c_{j_k}c_{j_1}^*
\end{eqnarray*}

By taking the tensor product of all these formulae, we get:
$$c_{i_1}\otimes c_{i_1}^*c_{i_2}\otimes \ldots\otimes c_{i_{k-2}}^*c_{i_{k-1}} \otimes c_{i_{k-1}}^*=\pm c_{j_1}c_{j_2}^*\otimes
c_{j_2}c^*_{j_3}\otimes \ldots \otimes c_{j_k}c_{j_1}^*$$

By applying the linear map given by $a_1\otimes\ldots \otimes a_k\to a_1\ldots a_k$ to both sides we see that the sign on the right is actually $+$. That is, we have:
$$c_{i_1}\otimes c_{i_1}^*c_{i_2}\otimes \ldots\otimes c_{i_{k-2}}^*c_{i_{k-1}} \otimes c_{i_{k-1}}^*=c_{j_1}c_{j_2}^*\otimes
c_{j_2}c^*_{j_3}\otimes \ldots \otimes c_{j_k}c_{j_1}^*$$

We recognize at right the basic summand in the formula of $R(c_{0_k})$. Now it follows from the definition of $i$ that summing the right terms over all multi-indices $j=(j_1,\ldots ,j_k)$ is the same as summing the left terms over all multi-indices $i=(i_1,\ldots,i_{k-1})$, then multiplying by $4$. Thus, we obtain the following formula:
\begin{eqnarray*}
R(c_{0_k})
&=&\frac{4}{2}\sum_{i_1\ldots i_{k-1}}c_{i_1}\otimes c_{i_1}^*c_{i_2}\otimes \ldots \otimes c_{i_{k-2}}^*c_{i_{k-1}}\otimes c_{i_{k-1}}^*\\
&=&2\sum_{i_1\ldots i_{k-1}}(c_{i_1}\otimes c_{i_1}^*)_{12}(c_{i_2}\otimes c_{i_2}^*)_{23}\ldots (c_{i_{k-1}}\otimes c_{i_{k-1}}^*)_{k-1,k}\\
&=&2\left(\sum_{i_1}c_{i_1}\otimes c_{i_1}^*\right)_{12}\left(\sum_{i_2}c_{i_2}\otimes c_{i_2}^*\right)_{23}\ldots \left(\sum_{i_{k-1}}c_{i_{k-1}}\otimes c_{i_{k-1}}^*\right)_{k-1,k}\cr
&=&2f_{12}f_{23}\ldots f_{k-1,k}
\end{eqnarray*}

Now $f$ being invariant under $SU_2$, the same happens for each $f_{i,i+1}$. Since the invariants form an algebra, the above product is invariant, and this concludes the proof.
\end{proof}

We prove now that the previous result holds not only for $c_0$, but in fact for any $c_p$. We first introduce some standard combinatorial notations. The Kreweras complement of a partition $p\in NC(k)$ is constructed as follows. Consider the ordered set $\{1,\ldots ,k\}$. At right of each index $i$ we put an index $i'$, as to get the following sequence of indices:
$$1,1',2,2',\ldots ,k,k'$$

The Kreweras complement $p^c$ is then the largest noncrossing partition of the new index set $\{1',\ldots ,k'\}$, such that the union of $p$ and $p^c$ is noncrossing. We have:

\begin{definition}
For a noncrossing partition $p$ we use the notation
$$w(p)=2\,\prod_{i=1}^l\left(\prod_{p=1}^{k_i-1} f_{j_{i,p}j_{i,p+1}}\right)$$
where $\{j_{11}<\ldots <j_{1k_1}\},\ldots ,\{j_{l1}<\ldots <j_{lk_l}\}$ are the blocks of $p$, with the convention that for $k_i=1$ the product on the right is by definition $1$. 
\end{definition}

Here we use the element $f$ introoduced above, and the leg-numbering notation. Observe that the element $w(p)$ is well-defined, because the products on the right pairwise commute. In particular in the case of the partition  $1_k=\{\{1,\ldots ,k\}\}$, we have:
$$w(1_k)=2\prod_{i=1}^{k-1}f_{i,i+1}$$

Since each $f_{ij}$ is invariant under the adjoint action of $SU_2$ and since the invariants form an algebra, $w(p)$ is also invariant under the adjoint action of $SU_2$. We have:

\begin{proposition}
$R(c_p)$, with $p\in NC(k)$, is invariant under the action of $SU_2$.
\end{proposition}

\begin{proof}
The idea is to generalize the proof of the previous result. In particular we need to generalize the key formula there, namely:
$$R(c_{0_k})=2f_{12}f_{23}\ldots f_{k-1,k}$$

Our claim is that for any noncrossing partition $p$, we have:
$$R(c_p)=w(p^c)$$

As a first verification, for the trivial partition $p=0_k$ we have $p^c=1_k$, and the formula $R(c_p)=w(p^c)$ follows from the above identities. Also, for the rough partition $p=1_k$ we have $p^c=0_k$, and the claimed equality follows from:
\begin{eqnarray*}
R(c_{1_k})
&=&\sum_i R(c_i\otimes\ldots\otimes c_i)\\
&=&\frac{1}{2}\sum_ic_ic_i^*\otimes\ldots\otimes c_ic_i^*\\
&=&\frac{1}{2}\,4(1\otimes\ldots\otimes 1)\\
&=&w(0_k)
\end{eqnarray*}

In the general case, we can use the same method as before. We have:
\begin{eqnarray*}
R(c_{p})
&=&R\left(\sum_j\delta_{pj}c_j\right)\\
&=&R\left(\sum_{j_1\ldots j_k}
\delta_{pj}c_{j_1}\otimes\ldots\otimes c_{j_k}\right)\\
&=&\frac{1}{2}\sum_{j_1\ldots j_k} \delta_{pj}c_{j_1}c^*_{j_2}\otimes
c_{j_2}c^*_{j_3}\otimes\ldots\otimes c_{j_k}c_{j_1}^*
\end{eqnarray*}

As before, in the previous proof, to any multi-index $j=(j_1,\ldots ,j_k)$ we associate a multi-index $i=(i_1,\ldots ,i_k)$, in the following way: 

\smallskip

-- $i_1$ is such that $c_{i_1}=\pm c_{j_1}c_{j_2}^*$, 

\smallskip

-- $i_2$ is such that $c_{i_2}=\pm c_{j_1}c_{j_3}^*$, 

$\vdots$

-- and so on up to $i_k$, which is such that $c_{i_{k}}=\pm c_{j_k}c_{j_1}^*$.

\smallskip

With this notation, and with the observation that the product of all the above $\pm$ signs is actually $+$, the formula that we established before becomes:
$$R(c_p)=\frac{1}{2}\sum_j\delta_{pj}c_{i_1}\otimes\ldots\otimes c_{i_k}$$

Now let $l\in\{1,\,\ldots ,k\}$ and assume that it is a last element of a block of $p^c$. Let $l_1<\ldots<l_x=l$ be the ordered enumeration of the elements of this block. We have: 
$$c_{j_{l_1}}\ldots c_{j_{l_x}}=1$$

As an illustration here, if $l$ and $l+1$ are in the same block of $p$, then $\{l\}$ is a one-element block in $p^c$ and $i_l$ will be such that $c_{i_l}=\pm c_{j_l}c_{j_{l+1}}^*=c_{j_l}c_{j_l}^*=1$. The general case works by following the same argument as before. We obtain:
$$\sum_i\delta_{pi}c_{j_{l_1}}\otimes \ldots 
\otimes c_{j_{l_x}}=R(c_{0_x})$$

The point now is that the element $R(c_p)$ is the product of the above expressions, over the blocks of $p^c$. If we denote these blocks by $\{l_{11}<\ldots<l_{1x_1}\}$, \ldots , $\{l_{r1}<\ldots<l_{rx_r}\}$, we obtain, by using the leg-numbering notation, the following formula:
\begin{eqnarray*}
R(c_p)
&=&2^{1-b}\sum_j\prod_{b=1}^rR(c_{0_{x_b}})_{l_{b1}\ldots l_{bx_b}}\\
&=&2^{1-b}\sum_j\prod_{b=1}^r\omega(1_{x_b})_{l_{b1}\ldots l_{bx_b}}\\
&=&w(p^c)
\end{eqnarray*}

All this might seem a bit complicated, and to illustrate the above proof let us perform explicitly the computation in the case of the following partition:
$$p=\{\{1,5\},\{2\},\{3,4\},\{6\}\}$$ 

The Kreweras complement of this partition is then $p^c=\{\{1,2,4\},\{3\},\{5,6\}\}$, and the above method gives the following formula:
\begin{eqnarray*}
R(c_p)
&=&\frac{1}{2}\sum_{j_1j_2j_3j_4}c_{j_1}c_{j_2^*}\otimes c_{j_2}c_{j_3}^*\otimes 1\otimes c_{j_3}c_{j_1}^*\otimes c_{j_1}c_{j_4}^*\otimes c_{j_4}c_{j_1}^*\\
&=&\frac{1}{2}\,4\sum_{i_1i_2i_3}c_{i_1}\otimes c_{i_1}^*c_{i_2}\otimes 1\otimes c_{i_2}^*\otimes c_{i_3}\otimes c_{i_3}^*\\
&=&2\sum_{i_1i_2i_3}(c_{i_1}\otimes c_{i_1}^*c_{i_2}\otimes 1\otimes c_{i_2}^*\otimes 1\otimes 1)(1\otimes 1\otimes 1\otimes 1\otimes c_{i_3}\otimes c_{i_3}^*)\\
&=&2(f_{12}f_{24})f_{56}\\
&=&w(p^c)
\end{eqnarray*}

Now back to the general case, the formula $R(c_p)=w(p^c)$ that we established shows that $R(c_p)$ is a product of certain elements $f_{ij}$ obtained from $f$ by acting on the various legs of $M_2(\mathbb C)^{\otimes k}$, and we can conclude as in the proof of the previous result.
\end{proof}

We are now in position of proving our main technical result, as follows:

\begin{proposition}
$R^*ER(c_p)=c_p$.
\end{proposition}

\begin{proof}
This comes from a routine computation, as follows:

\medskip

(1) To start with, the previous result gives the following formula:
\begin{eqnarray*}
<R^*ER(c_p),c_j>
&=&<R^*R(c_p),c_j>\\
&=&<R(c_p),R(c_j)>\\
&=&\sum_i\delta_{pi}<R(c_i),R(c_j)>
\end{eqnarray*}

(2) On the other hand, we have from definitions the following formula:
$$<c_p,c_j>=\delta_{pj}$$

Since the elements $c_j$ span the ambient space, what is left to prove is:
$$\sum_i\delta_{pi}<R(c_i),R(c_j)>=\delta_{pj}$$

(3) But this can be checked by direct computation. Indeed, from the definition of the operator $R$, we get the following formula, for the above scalar products:
\begin{eqnarray*}
<R(c_i),R(c_j)>
&=&\frac{1}{4}\,<c_{i_1}c_{i_2}^*\otimes\ldots\otimes
c_{i_k}c_{i_1}^*,c_{j_1}c_{j_2}^*\otimes\ldots\otimes
c_{j_k}c_{j_1}^*>\\
&=&\frac{1}{4}\,<c_{i_1}c_{i_2}^*,c_{j_1}c_{j_2}^*>
\ldots <c_{i_k}c_{i_1}^*,c_{j_k}c_{j_1}^*>
\end{eqnarray*}

(4) In this latter formula all scalar products are $0$, $1$ and $-1$. Now assume that the indices $i_1,\ldots ,i_k$ and $j_1$ are all fixed. If the first scalar product is $\pm 1$ then $j_2$ is uniquely determined, then if the second scalar product is $\pm 1$ then $j_3$ is uniquely determined as well, and so on. Thus for all scalar products to be $\pm 1$, the multi-index $j$ is uniquely determined by the multi-index $i$, up to a possible choice of the first index $j_1$.

\medskip

(5) Moreover, each choice of the index $j_1$ leads to a multi-index $j=(j_1,\ldots,j_k)$, such that all the scalar products are $\pm 1$. Indeed, once $j_2,\ldots,j_k$ are chosen as to satisfy $c_{i_r}c_{i_{r+1}}^*=\pm c_{j_r}c_{j_{r+1}}^*$ for $r=1,\ldots,k-1$, by multiplying all these formulae we obtain $c_{i_1}c_{i_k}^*=\pm c_{j_1}c_{j_k}^*$, which shows that the last scalar product is $\pm 1$ as well.

\medskip

(6) Summarizing, given a multi-index $i=(i_1,\ldots ,i_k)$ and a number $s\in\{1,2,3,4\}$, there is a unique multi-index $j=(j_1,\ldots,j_k)$ with $j_1=s$, such that all the above scalar products are $\pm 1$. We use the notation $j=i\oplus s$, for this multi-index.

\medskip

(7) In the situation $j=i\oplus s$ we have $<R(c_i),R(c_j)>=\pm 1$, and by applying the linear map given by $a_1\otimes\ldots \otimes a_k\to a_1\ldots a_k$ to the formula $c_{i_1}c_{i_2}^*\otimes\ldots\otimes c_{i_k}c_{i_1}^*=\pm c_{j_1}c_{j_2}^*\otimes\ldots\otimes c_{j_k}c_{j_1}^*$ we see that the sign is actually $+$. Thus we have the following formula:
$$<R(c_i),R(c_j)>
=\begin{cases}
1/4\mbox{ if $j=i\oplus s$ for some $s\in\{1,2,3,4\}$}\\
0\mbox{ otherwise}
\end{cases}$$

(8) We can come back now to the missing formula. We have:
$$\sum_i\delta_{pi}<R(c_i),R(c_j)>
=\frac{1}{4}\sum_{s=1}^4\delta_{pi}$$

Our claim is that we have $\delta_{pi}=\delta_{pj}$, for any partition $p$. Indeed, this follows from the fact that for $r<s$ we have $i_r=i_s$ if and only if the product of $c_{i_t}c_{i_{t+1}}^*$ over $t=r,\ldots,s-1$ equals $\pm 1$, and a similar statement holds for the multi-index $j$. 

\medskip

(9) We can now conclude the proof. By using $\delta_{pi}=\delta_{pj}$, we get:
$$\sum_i\delta_{pi}<R(c_i),R(c_j)>
=\frac{1}{4}\sum_{s=1}^4\delta_{pj}=\delta_{pj}$$

But this is the formula that we wanted to prove, so we are done.
\end{proof}

Still following \cite{bc2}, we can now prove our main result, as follows:

\begin{theorem}
The Pauli representation of $C(S_4^+)$ is faithful.
\end{theorem}

\begin{proof}
We denote as usual by $c_1,\ldots,c_4$ the Pauli matrices, and we let $e_1,\ldots,e_4$ be the standard basis of $\mathbb C^4$. For a multi-index $i=(i_1,\ldots ,i_k)$, we set:
$$e_i=e_{i_1}\otimes\ldots\otimes e_{i_k}\quad,\quad 
c_i=c_{i_1}\otimes\ldots\otimes c_{i_k}$$

Each partition $p\in NC(k)$ creates two tensors, in the following way:
$$e_p=\sum_i\delta_{pi}\,e_i\quad,\quad 
c_p=\sum_i\delta_{pi}\,c_i$$

Consider the following $4^k\times 4^k$ matrices, with entries labeled by multi-indices $i,j$:
$$P_{ij}=tr\left(\int \pi_{i_1j_1}\ldots \pi_{i_kj_k}\right)\quad,\quad 
U_{ij}=\left(\int u_{i_1j_1}\ldots u_{i_kj_k}\right)$$

According to Proposition 15.5, we have the following formula:
$$P_{ij}=<c_i,R^*ER(c_j)>$$

Now let $\Phi$ be the linear map $(\mathbb{C}^4)^{\otimes k}\to M_2(\mathbb{C})^{\otimes k}$
given by $\Phi (e_i)=c_i$. In terms of this map, the above equation can be rephrased as follows:
$$P_{ij}=<e_i,\Phi^*R^*ER\Phi(e_j)>$$

According to our results above, it is enough to prove that we have $P=U$. But $U$ is the orthogonal projection of $(\mathbb{C}^4)^{\otimes k}$ onto the following space:
$$S_e=span\left\{ e_p\Big|p\in NC(k)\right\}$$

Thus what we have to prove is that $\Phi^*R^*ER\Phi$ is the projection onto $S_e$. But this is equivalent to proving that $R^*ER$ is the projection onto the following space: 
$$S_c=span\left\{ c_p\Big|p\in NC(k)\right\}$$

We know that $R^*ER(c_p)=c_p$, for all $p\in NC(k)$. But this implies that
$R^*ER$ restricted to $S_c$ is the identity. This vector space has dimension the Catalan number $C_k$, and this is  exactly the rank of the operator $E$. Therefore 
$S_c$ has to be the image of $R^*ER$. Now since $R^*ER$ is self-adjoint, and is the identity on its image, it is the orthogonal projection onto $S_c$, and this concludes the proof.
\end{proof}

This was for the result from \cite{bc2}. There are of course many other proofs for the above result, and generalizations. A purely algebraic proof, based on $S_4^+=SO_3^{-1}$, is explained in \cite{bb2}, and further generalizations are discussed in \cite{bbs}. We will discuss in a moment a generalization of this, from \cite{bne}, which is analytic, along the above lines.

\section*{15b. ADE models}

We know from the above, based on \cite{bc2}, that we have a stationary matrix model for the algebra $C(S_4^+)$. In view of \cite{bb2}, this suggests the following conjecture:

\begin{conjecture}
Given a quantum permutation group of $4$ points, 
$$G\subset S_4^+\simeq SO_3^{-1}$$
coming by twisting a usual ADE subgroup of the group $SO_3$,
$$H\subset SO_3$$
the restriction of the Pauli model for $C(S_4^+)$, with fibers coming from the elements of $H\subset SO_3$, has the algebra $C(G)$ as Hopf image.
\end{conjecture}

To be more precise, the main result from the previous section tells us that the conjecture holds for $G=S_4^+$ itself. Indeed, here we have $H=SO_3$, so the corresponding restriction of the Pauli model for $C(S_4^+)$ is the Pauli model itself, and this model being stationary, its Hopf image is the algebra $C(S_4^+)$ itself, as stated.

\bigskip

In general, the above conjecture does not look that scary, because the same methods used for $S_4^+$ can be used for any subgroup $G\subset S_4^+$. However, the problem is that, unless a global method in order to uniformly deal with the problem is found, this would need a case-by-case study depending on $G\subset S_4^+$, which looks quite time-consuming.

\bigskip

Finally, let us mention that such models would be extremely useful, in connection with a wide array of analytic questions. For a discussion here, see McCarthy \cite{mc2}.

\section*{15c. Weyl matrices}

Following \cite{bc2}, \cite{bne}, let us discuss now some further examples of stationary models, related to the Pauli matrices, and Weyl matrices. We first have:

\index{Weyl matrix}
\index{Pauli matrix}

\begin{definition}
Given a finite abelian group $H$, the associated Weyl matrices are
$$W_{ia}:e_b\to<i,b>e_{a+b}$$
where $i\in H$, $a,b\in\widehat{H}$, and where $(i,b)\to<i,b>$ is the Fourier coupling $H\times\widehat{H}\to\mathbb T$.
\end{definition}

As a basic example, consider the cyclic group $H=\mathbb Z_2=\{0,1\}$. Here the Fourier coupling is given by $<i,b>=(-1)^{ib}$, and so the Weyl matrices act via:
$$W_{00}:e_b\to e_b\quad,\quad 
W_{10}:e_b\to(-1)^be_b$$
$$W_{11}:e_b\to(-1)^be_{b+1}\quad,\quad 
W_{01}:e_b\to e_{b+1}$$

Thus, we have the following formulae:
$$W_{00}=\begin{pmatrix}1&0\\0&1\end{pmatrix}\quad,\quad 
W_{10}=\begin{pmatrix}1&0\\0&-1\end{pmatrix}$$
$$W_{11}=\begin{pmatrix}0&-1\\1&0\end{pmatrix}\quad,\quad
W_{01}=\begin{pmatrix}0&1\\1&0\end{pmatrix}$$

We recognize here, up to some multiplicative factors, the four Pauli matrices. Now back to the general case, we have the following well-known result:

\begin{proposition}
The Weyl matrices are unitaries, and satisfy:
\begin{enumerate}
\item $W_{ia}^*=<i,a>W_{-i,-a}$.

\item $W_{ia}W_{jb}=<i,b>W_{i+j,a+b}$.

\item $W_{ia}W_{jb}^*=<j-i,b>W_{i-j,a-b}$.

\item $W_{ia}^*W_{jb}=<i,a-b>W_{j-i,b-a}$.
\end{enumerate}
\end{proposition}

\begin{proof}
The unitarity follows from (3,4), and the rest of the proof goes as follows:

\medskip

(1) We have indeed the following computation:
\begin{eqnarray*}
W_{ia}^*
&=&\left(\sum_b<i,b>E_{a+b,b}\right)^*\\
&=&\sum_b<-i,b>E_{b,a+b}\\
&=&\sum_b<-i,b-a>E_{b-a,b}\\
&=&<i,a>W_{-i,-a}
\end{eqnarray*}

(2) Here the verification goes as follows:
\begin{eqnarray*}
W_{ia}W_{jb}
&=&\left(\sum_d<i,b+d>E_{a+b+d,b+d}\right)\left(\sum_d<j,d>E_{b+d,d}\right)\\
&=&\sum_d<i,b><i+j,d>E_{a+b+d,d}\\
&=&<i,b>W_{i+j,a+b}
\end{eqnarray*}

(3,4) By combining the above two formulae, we obtain:
\begin{eqnarray*}
W_{ia}W_{jb}^*
&=&<j,b>W_{ia}W_{-j,-b}\\
&=&<j,b><i,-b>W_{i-j,a-b}
\end{eqnarray*}

We obtain as well the following formula:
\begin{eqnarray*}
W_{ia}^*W_{jb}
&=&<i,a>W_{-i,-a}W_{jb}\\
&=&<i,a><-i,b>W_{j-i,b-a}
\end{eqnarray*}

But this gives the formulae in the statement, and we are done.
\end{proof}

With $n=|H|$, we can use an isomorphism $l^2(\widehat{H})\simeq\mathbb C^n$ as to view each $W_{ia}$ as a usual matrix, $W_{ia}\in M_n(\mathbb C)$, and hence as a usual unitary, $W_{ia}\in U_n$. Also, given a vector $\xi$, we denote by $Proj(\xi)$ the orthogonal projection onto $\mathbb C\xi$. Following \cite{bne}, we have:

\begin{proposition}
Given a closed subgroup $E\subset U_n$, we have a representation
$$\pi_H:C(S_N^+)\to M_N(C(E))$$
$$w_{ia,jb}\to[U\to Proj(W_{ia}UW_{jb}^*)]$$
where $n=|H|,N=n^2$, and where $W_{ia}$ are the Weyl matrices associated to $H$.
\end{proposition}

\begin{proof}
The Weyl matrices being given by $W_{ia}:e_b\to<i,b>e_{a+b}$, we have:
$$tr(W_{ia})=\begin{cases}
1&{\rm if}\ (i,a)=(0,0)\\
0&{\rm if}\ (i,a)\neq(0,0)
\end{cases}$$

Together with the formulae in Proposition 15.14, this shows that the Weyl matrices are pairwise orthogonal with respect to the following scalar product on $M_n(\mathbb C)$:
$$<x,y>=tr(xy^*)$$

Thus, these matrices form an orthogonal basis of $M_n(\mathbb C)$, consisting of unitaries:
$$W=\left\{W_{ia}\Big|i\in H,a\in\widehat{H}\right\}$$

Thus, each row and each column of the matrix $\xi_{ia,jb}=W_{ia}UW_{jb}^*$ is an orthogonal basis of $M_n(\mathbb C)$, and so the corresponding projections form a magic unitary, as claimed.
\end{proof}

We will need the following well-known result:

\begin{proposition}
With $T=Proj(x_1)\ldots Proj(x_p)$ and $||x_i||=1$ we have 
$$<T\xi,\eta>=<\xi,x_p><x_p,x_{p-1}>\ldots<x_2,x_1><x_1,\eta>$$
for any $\xi,\eta$. In particular, we have:
$$Tr(T)=<x_1,x_p><x_p,x_{p-1}>\ldots<x_2,x_1>$$
\end{proposition}

\begin{proof}
For $||x||=1$ we have $Proj(x)\xi=<\xi,x>x$. This gives:
\begin{eqnarray*}
T\xi
&=&Proj(x_1)\ldots Proj(x_p)\xi\\
&=&Proj(x_1)\ldots Proj(x_{p-1})<\xi,x_p>x_p\\
&=&Proj(x_1)\ldots Proj(x_{p-2})<\xi,x_p><x_p,x_{p-1}>x_{p-1}\\
&=&\ldots\\
&=&<\xi,x_p><x_p,x_{p-1}>\ldots<x_2,x_1>x_1
\end{eqnarray*}

Now by taking the scalar product with $\eta$, this gives the first assertion. As for the second assertion, this follows from the first assertion, by summing over $\xi=\eta=e_i$.
\end{proof}

Now back to the Weyl matrix models, let us first compute $T_p$. We have:

\begin{proposition}
We have the formula
\begin{eqnarray*}
(T_p)_{ia,jb}
&=&\frac{1}{N}<i_1,a_1-a_p>\ldots<i_p,a_p-a_{p-1}><j_1,b_1-b_2>\ldots<j_p,b_p-b_1>\\
&&\int_Etr(W_{i_1-i_2,a_1-a_2}UW_{j_2-j_1,b_2-b_1}U^*)\ldots tr(W_{i_p-i_1,a_p-a_1}UW_{j_1-j_p,b_1-b_p}U^*)dU
\end{eqnarray*}
with all the indices varying in a cyclic way.
\end{proposition}

\begin{proof}
By using the trace formula in Proposition 15.16, we obtain:
\begin{eqnarray*}
&&(T_p)_{ia,jb}\\
&=&\left(tr\otimes\int_E\right)\left(Proj(W_{i_1a_1}UW_{j_1b_1}^*)\ldots Proj(W_{i_pa_p}UW_{j_pb_p}^*)\right)\\
&=&\frac{1}{N}\int_E<W_{i_1a_1}UW_{j_1b_1}^*,W_{i_pa_p}UW_{j_pb_p}^*>\ldots<W_{i_2a_2}UW_{j_2b_2}^*,W_{i_1a_1}UW_{j_1b_1}^*>dU
\end{eqnarray*}

In order to compute now the scalar products, observe that we have:
\begin{eqnarray*}
<W_{ia}UW_{jb}^*,W_{kc}UW_{ld}^*>
&=&tr(W_{jb}U^*W_{ia}^*W_{kc}UW_{ld}^*)\\
&=&tr(W_{ia}^*W_{kc}UW_{ld}^*W_{jb}U^*)\\
&=&<i,a-c><l,d-b>tr(W_{k-i,c-a}UW_{j-l,b-d}U^*)
\end{eqnarray*}

By plugging these quantities into the formula of $T_p$, we obtain the result.
\end{proof}

Consider now the Weyl group $W=\{W_{ia}\}\subset U_n$, that we already met in the proof of Proposition 15.15. We have the following result, from \cite{bne}:

\index{Weyl matrix model}

\begin{theorem}
For any compact group $W\subset E\subset U_n$, the model
$$\pi_H:C(S_N^+)\to M_N(C(E))$$
$$w_{ia,jb}\to[U\to Proj(W_{ia}UW_{jb}^*)]$$
constructed above is stationary on its image.
\end{theorem}

\begin{proof}
We must prove that we have $T_p^2=T_p$. We have:
\begin{eqnarray*}
&&(T_p^2)_{ia,jb}\\
&=&\sum_{kc}(T_p)_{ia,kc}(T_p)_{kc,jb}\\
&=&\frac{1}{N^2}\sum_{kc}<i_1,a_1-a_p>\ldots<i_p,a_p-a_{p-1}><k_1,c_1-c_2>\ldots<k_p,c_p-c_1>\\
&&<k_1,c_1-c_p>\ldots<k_p,c_p-c_{p-1}><j_1,b_1-b_2>\ldots<j_p,b_p-b_1>\\
&&\int_Etr(W_{i_1-i_2,a_1-a_2}UW_{k_2-k_1,c_2-c_1}U^*)\ldots tr(W_{i_p-i_1,a_p-a_1}UW_{k_1-k_p,c_1-c_p}U^*)dU\\
&&\int_Etr(W_{k_1-k_2,c_1-c_2}VW_{j_2-j_1,b_2-b_1}V^*)\ldots tr(W_{k_p-k_1,c_p-c_1}VW_{j_1-j_p,b_1-b_p}V^*)dV
\end{eqnarray*}

By rearranging the terms, this formula becomes:
\begin{eqnarray*}
&&(T_p^2)_{ia,jb}\\
&=&\frac{1}{N^2}<i_1,a_1-a_p>\ldots<i_p,a_p-a_{p-1}><j_1,b_1-b_2>\ldots<j_p,b_p-b_1>\\
&&\int_E\int_E\sum_{kc}<k_1-k_p,c_1-c_p>\ldots<k_p-k_{p-1},c_p-c_{p-1}>\\
&&tr(W_{i_1-i_2,a_1-a_2}UW_{k_2-k_1,c_2-c_1}U^*)tr(W_{k_1-k_2,c_1-c_2}VW_{j_2-j_1,b_2-b_1}V^*)\\
&&\hskip50mm\ldots\ldots\\
&&tr(W_{i_p-i_1,a_p-a_1}UW_{k_1-k_p,c_1-c_p}U^*)tr(W_{k_p-k_1,c_p-c_1}VW_{j_1-j_p,b_1-b_p}V^*)dUdV
\end{eqnarray*}

Let us denote by $I$ the above double integral. By using $W_{kc}^*=<k,c>W_{-k,-c}$ for each of the couplings, and by moving as well all the $U^*$ variables to the left, we obtain:
\begin{eqnarray*}
I
&=&\int_E\int_E\sum_{kc}tr(U^*W_{i_1-i_2,a_1-a_2}UW_{k_2-k_1,c_2-c_1})tr(W_{k_2-k_1,c_2-c_1}^*VW_{j_2-j_1,b_2-b_1}V^*)\\
&&\hskip50mm\ldots\ldots\\
&&tr(U^*W_{i_p-i_1,a_p-a_1}UW_{k_1-k_p,c_1-c_p})tr(W_{k_1-k_p,c_1-c_p}^*VW_{j_1-j_p,b_1-b_p}V^*)dUdV
\end{eqnarray*}

In order to perform now the sums, we use the following formula:
\begin{eqnarray*}
tr(AW_{kc})tr(W_{kc}^*B)
&=&\frac{1}{N}\sum_{qrst}A_{qr}(W_{kc})_{rq}(W^*_{kc})_{st}B_{ts}\\
&=&\frac{1}{N}\sum_{qrst}A_{qr}<k,q>\delta_{r-q,c}<k,-s>\delta_{t-s,c}B_{ts}\\
&=&\frac{1}{N}\sum_{qs}<k,q-s>A_{q,q+c}B_{s+c,s}
\end{eqnarray*}

If we denote by $A_x,B_x$ the variables which appear in the formula of $I$, we have:
\begin{eqnarray*}
&&I\\
&=&\frac{1}{N^p}\int_E\int_E\sum_{kcqs}<k_2-k_1,q_1-s_1>\ldots<k_1-k_p,q_p-s_p>\\
&&(A_1)_{q_1,q_1+c_2-c_1}(B_1)_{s_1+c_2-c_1,s_1}\ldots (A_p)_{q_p,q_p+c_1-c_p}(B_p)_{s_p+c_1-c_p,s_p}\\
&=&\frac{1}{N^p}\int_E\int_E\sum_{kcqs}<k_1,q_p-s_p-q_1+s_1>\ldots<k_p,q_{p-1}-s_{p-1}-q_p+s_p>\\
&&(A_1)_{q_1,q_1+c_2-c_1}(B_1)_{s_1+c_2-c_1,s_1}\ldots (A_p)_{q_p,q_p+c_1-c_p}(B_p)_{s_p+c_1-c_p,s_p}
\end{eqnarray*}

Now observe that we can perform the sums over $k_1,\ldots,k_p$. We obtain in this way a multiplicative factor $n^p$, along with the condition:
$$q_1-s_1=\ldots=q_p-s_p$$

Thus we must have $q_x=s_x+a$ for a certain $a$, and the above formula becomes:
$$I=\frac{1}{n^p}\int_E\int_E\sum_{csa}(A_1)_{s_1+a,s_1+c_2-c_1+a}(B_1)_{s_1+c_2-c_1,s_1}\ldots(A_p)_{s_p+a,s_p+c_1-c_p+a}(B_p)_{s_p+c_1-c_p,s_p}$$

Consider now the variables $r_x=c_{x+1}-c_x$, which altogether range over the set $Z$ of multi-indices having sum 0. By replacing the sum over $c_x$ with the sum over $r_x$, which creates a multiplicative $n$ factor, we obtain the following formula:
$$I=\frac{1}{n^{p-1}}\int_E\int_E\sum_{r\in Z}\sum_{sa}(A_1)_{s_1+a,s_1+r_1+a}(B_1)_{s_1+r_1,s_1}\ldots(A_p)_{s_p+a,s_p+r_p+a}(B_p)_{s_p+r_p,s_p}$$

For an arbitrary multi-index $r$, we have the following formula:
$$\delta_{\sum_ir_i,0}=\frac{1}{n}\sum_i<i,r_1>\ldots<i,r_p>$$

Thus, we can replace the sum over $r\in Z$ by a full sum, as follows:
\begin{eqnarray*}
I
&=&\frac{1}{n^p}\int_E\int_E\sum_{rsia}<i,r_1>(A_1)_{s_1+a,s_1+r_1+a}(B_1)_{s_1+r_1,s_1}\\
&&\hskip40mm\ldots\ldots\\
&&\hskip20mm<i,r_p>(A_p)_{s_p+a,s_p+r_p+a}(B_p)_{s_p+r_p,s_p}
\end{eqnarray*}

In order to ``absorb'' now the indices $i,a$, we can use the following formula:
\begin{eqnarray*}
&&W_{ia}^*AW_{ia}\\
&=&\left(\sum_b<i,-b>E_{b,a+b}\right)\left(\sum_{bc}E_{a+b,a+c}A_{a+b,a+c}\right)\left(\sum_c<i,c>E_{a+c,c}\right)\\
&=&\sum_{bc}<i,c-b>E_{bc}A_{a+b,a+c}
\end{eqnarray*}

Thus we have the following formula:
$$(W_{ia}^*AW_{ia})_{bc}=<i,c-b>A_{a+b,a+c}$$

With this in hand, our formula becomes:
\begin{eqnarray*}
&&I\\
&=&\frac{1}{n^p}\int_E\int_E\sum_{rsia}(W_{ia}^*A_1W_{ia})_{s_1,s_1+r_1}(B_1)_{s_1+r_1,s_1}\ldots(W_{ia}^*A_pW_{ia})_{s_p,s_p+r_p}(B_p)_{s_p+r_p,s_p}\\
&=&\int_E\int_E\sum_{ia}tr(W_{ia}^*A_1W_{ia}B_1)\ldots\ldots tr(W_{ia}^*A_pW_{ia}B_p)
\end{eqnarray*}

Now by replacing $A_x,B_x$ with their respective values, we obtain:
\begin{eqnarray*}
I
&=&\int_E\int_E\sum_{ia}tr(W_{ia}^*U^*W_{i_1-i_2,a_1-a_2}UW_{ia}VW_{j_2-j_1,b_2-b_1}V^*)\\
&&\hskip30mm\ldots\ldots\\
&&tr(W_{ia}^*U^*W_{i_p-i_1,a_p-a_1}UW_{ia}VW_{j_1-j_p,b_1-b_p}V^*)dUdV
\end{eqnarray*}

By moving the $W_{ia}^*U^*$ variables at right, we obtain, with $S_{ia}=UW_{ia}V$:
\begin{eqnarray*}
I
&=&\sum_{ia}\int_E\int_Etr(W_{i_1-i_2,a_1-a_2}S_{ia}W_{j_2-j_1,b_2-b_1}S_{ia}^*)\\
&&\hskip30mm\ldots\ldots\\
&&tr(W_{i_p-i_1,a_p-a_1}S_{ia}W_{j_1-j_p,b_1-b_p}S_{ia}^*)dUdV
\end{eqnarray*}

Now since $S_{ia}$ is Haar distributed when $U,V$ are Haar distributed, we obtain:
$$I=N\int_E\int_Etr(W_{i_1-i_2,a_1-a_2}UW_{j_2-j_1,b_2-b_1}U^*)\ldots tr(W_{i_p-i_1,a_p-a_1}UW_{j_1-j_p,b_1-b_p}U^*)dU$$

But this is exactly $N$ times the integral in the formula of $(T_p)_{ia,jb}$, from Proposition 15.17. Since the $N$ factor cancels with one of the two $N$ factors that we found in the beginning of the proof, when first computing $(T_p^2)_{ia,jb}$, we are done.
\end{proof}

As an illustration for the above result, going back to \cite{bc2}, we have:

\index{Pauli model}

\begin{theorem}
We have a stationary matrix model
$$\pi:C(S_4^+)\subset M_4(C(SU_2))$$
given on the standard coordinates by the formula
$$\pi(u_{ij})=[x\to Proj(c_ixc_j)]$$
where $x\in SU_2$, and $c_1,c_2,c_3,c_4$ are the Pauli matrices.
\end{theorem}

\begin{proof}
As already explained in the comments following Definition 15.13, the Pauli matrices appear as particular cases of the Weyl matrices:
$$W_{00}=\begin{pmatrix}1&0\\0&1\end{pmatrix}\quad,\quad 
W_{10}=\begin{pmatrix}1&0\\0&-1\end{pmatrix}$$
$$W_{11}=\begin{pmatrix}0&-1\\1&0\end{pmatrix}\quad,\quad
W_{01}=\begin{pmatrix}0&1\\1&0\end{pmatrix}$$

Thus, Theorem 15.18 produces in this case the model in the statement.
\end{proof}

Observe that, since the projection $Proj(c_ixc_j)$ depends only on the image of $x$ in the quotient $SU_2\to SO_3$, we can replace the model space $SU_2$ by the smaller space $SO_3$.  This can be used in conjunction with the isomorphism $S_4^+\simeq SO_3^{-1}$, and as explained in \cite{bb2}, our model becomes in this way something more conceptual, as follows:
$$\pi:C(SO_3^{-1})\subset M_4(C(SO_3))$$

As a philosophical conclusion, to this and to some previous findings as well, no matter what we do, we always end up getting back to $SU_2,SO_3$. Thus, we are probably doing some physics here. This is indeed the case, the above computations being closely related to the standard computations for the Ising and Potts models. The general relation, however, between quantum permutations and lattice models, is not axiomatixed yet.

\section*{15d. Cocyclic models}

Following \cite{bne}, let us discuss now some interesting generalizations of the Weyl matrix models. We will need the following standard definition:

\begin{definition}
A $2$-cocycle on a group $G$ is a function $\sigma:G\times G\to\mathbb T$ satisfying:
$$\sigma(gh,k)\sigma(g,h)=\sigma(g,hk)\sigma(h,k)\quad,\quad 
\sigma(g,1)=\sigma(1,g)=1$$
The algebra $C^*(G)$, with multiplication given by $g\cdot h=\sigma(g,h)gh$, and with the involution making the standard generators $g\in C^*_\sigma(G)$ unitaries, is denoted $C^*_\sigma(G)$.
\end{definition}

As explained in \cite{bne}, we have the following general construction:

\begin{proposition}
Given a finite group $G=\{g_1,\ldots,g_N\}$ and a $2$-cocycle on it, $\sigma:G\times G\to\mathbb T$, we have a matrix model as follows, 
$$\pi:C(S_N^+)\to M_N(C(E))\quad,\quad 
w_{ij}\to[x\to Proj(g_ixg_j^*)]$$
for any closed subgroup $E\subset U_A$, where $A=C^*_\sigma(G)$. 
\end{proposition}

\begin{proof}
This is clear from definitions, because the standard generators $\{g_1,\ldots,g_N\}$ are pairwise orthogonal with respect to the canonical trace of $A$. See \cite{bne}.
\end{proof}

In order to investigate the stationarity of $\pi$, we will use:

\begin{proposition}
We have the formula
\begin{eqnarray*}
(T_p)_{i_1\ldots i_p,j_1\ldots j_p}
&=&\overline{\sigma(i_1,i_1^{-1}i_2)}\ldots\overline{\sigma(i_p,i_p^{-1}i_1)}\cdot\overline{\sigma(j_2,j_2^{-1}j_1)}\ldots\overline{\sigma(j_1,j_1^{-1}j_p)}\\
&&\frac{1}{N}\int_Etr(g_{i_1^{-1}i_2}xg_{j_2^{-1}j_1}x^*)\ldots\ldots tr(g_{i_p^{-1}i_1}xg_{j_1^{-1}j_p}x^*)dx
\end{eqnarray*}
with all the indices varying in a cyclic way.
\end{proposition}

\begin{proof}
According to the definition of $T_p$, we have the following formula:
\begin{eqnarray*}
(T_p)_{i_1\ldots i_p,j_1\ldots j_p}
&=&\left(tr\otimes\int_E\right)\left(Proj(g_{i_1}xg_{j_1}^*)\ldots Proj(g_{i_p}xg_{j_p}^*)\right)dx\\
&=&\frac{1}{N}\int_E<g_{i_1}xg_{j_1}^*,g_{i_2}xg_{j_2}^*>\ldots<g_{i_p}xg_{j_p}^*,g_{i_1}xg_{j_1}^*>dx
\end{eqnarray*}

In order to compute the scalar products, we can use the following formula:
$$g_ig_{i^{-1}k}=\sigma(i,i^{-1}k)g_k$$

Indeed, we obtain from this the following formula:
$$g_i^*g_k=\overline{\sigma(i,i^{-1}k)}g_{i^{-1}k}$$

We therefore obtain the following formula, for the above scalar products:
\begin{eqnarray*}
<g_ixg_j^*,g_kxg_l^*>
&=&tr(g_jx^*g_i^*g_kxg_l^*)\\
&=&tr(g_i^*g_kxg_l^*g_jx^*)\\
&=&\overline{\sigma(i,i^{-1}k)}\cdot\overline{\sigma(l,l^{-1}j)}\cdot tr(g_{i^{-1}k}xg_{l^{-1}j}x^*)
\end{eqnarray*}

By plugging these quantities into the formula of $T_p$, we obtain the result.
\end{proof}

We have the following result, which generalizes our previous computations:

\begin{theorem}
For any intermediate closed subgroup $G\subset E\subset U_A$, the model 
$$\pi:C(S_N^+)\to M_N(C(E))$$
constructed above is stationary on its image.
\end{theorem}

\begin{proof}
We use the formula in Proposition 15.22. Let us write $(T_p)_{ij}=\rho(i,j)(T_p^\circ)_{ij}$, where $\rho(i,j)$ is the product of $\sigma$ terms appearing there. We have:
\begin{eqnarray*}
(T_p^2)_{ij}
&=&\sum_k(T_p)_{ik}(T_p)_{kj}\\
&=&\sum_k\rho(i,k)\rho(k,j)(T_p^\circ)_{ik}(T_p^\circ)_{kj}
\end{eqnarray*}

Let us first compute the $\rho$ term. We have:
\begin{eqnarray*}
\rho(i,k)\rho(k,j)
&=&\overline{\sigma(i_1,i_1^{-1}i_2)}\ldots\overline{\sigma(i_p,i_p^{-1}i_1)}\cdot\overline{\sigma(k_2,k_2^{-1}k_1)}\ldots\overline{\sigma(k_1,k_1^{-1}k_p)}\\
&&\overline{\sigma(k_1,k_1^{-1}k_2)}\ldots\overline{\sigma(k_p,k_p^{-1}k_1)}\cdot\overline{\sigma(j_2,j_2^{-1}j_1)}\ldots\overline{\sigma(j_1,j_1^{-1}j_p)}\\
&=&\sigma(i,j)\cdot\overline{\sigma(k_2,k_2^{-1}k_1)}\cdot\overline{\sigma(k_1,k_1^{-1}k_2)}\ldots\ldots\overline{\sigma(k_1,k_1^{-1}k_p)}\cdot\overline{\sigma(k_p,k_p^{-1}k_1)}
\end{eqnarray*}

In order to compute now this quantity, observe that by multiplying the formulae $\sigma(i,i^{-1}k)g_i^*g_k=g_{i^{-1}k}$ and $\sigma(k,k^{-1}i)g_k^*g_i=g_{k^{-1}i}$ we obtain the following formula:
$$\sigma(i,i^{-1}k)\sigma(k,k^{-1}i)=\sigma(i^{-1}k,k^{-1}i)$$

Thus, our expression further simplifies, as follows:
$$\rho(i,k)\rho(k,j)=\sigma(i,j)\cdot\overline{\sigma(k_2^{-1}k_1,k_1^{-1}k_2)}\ldots\ldots\overline{\sigma(k_1^{-1}k_p,k_p^{-1}k_1)}$$

On the other hand, the $T^\circ$ term can be written as follows:
\begin{eqnarray*}
(T_p^\circ)_{ik}(T_p^\circ)_{kj}
&=&\frac{1}{N^2}\int_E\int_Etr(g_{i_1^{-1}i_2}xg_{k_2^{-1}k_1}x^*)tr(g_{k_1^{-1}k_2}yg_{j_2^{-1}j_1}y^*)\\
&&\hskip40mm\ldots\ldots\\
&&\hskip10mm tr(g_{i_p^{-1}i_1}xg_{k_1^{-1}k_1p}x^*)tr(g_{k_p^{-1}k_1}yg_{j_1^{-1}j_p}y^*)dxdy
\end{eqnarray*}

We therefore conclude that we have the following formula:
\begin{eqnarray*}
(T_p^2)_{ij}
&=&\frac{\sigma(i,j)}{N^2}\int_E\int_E\sum_{k_1\ldots k_p}
\overline{\sigma(k_2^{-1}k_1,k_1^{-1}k_2)}tr(g_{i_1^{-1}i_2}xg_{k_2^{-1}k_1}x^*)tr(g_{k_1^{-1}k_2}yg_{j_2^{-1}j_1}y^*)\\
&&\hskip50mm\vdots\\
&&\hskip10mm\overline{\sigma(k_1^{-1}k_p,k_p^{-1}k_1)}tr(g_{i_p^{-1}i_1}xg_{k_1^{-1}k_p}x^*)tr(g_{k_p^{-1}k_1}yg_{j_1^{-1}j_p}y^*)dxdy
\end{eqnarray*}

By using now $g_i^*=\overline{\sigma(i,i^{-1})}g_{i^{-1}}$, and moving as well the $x^*$ variables at left, we obtain:
\begin{eqnarray*}
(T_p^2)_{ij}
&=&\frac{\sigma(i,j)}{N^2}\int_E\int_E\sum_{k_1\ldots k_p}
tr(x^*g_{i_1^{-1}i_2}xg_{k_2^{-1}k_1})tr(g_{k_2^{-1}k_1}^*yg_{j_2^{-1}j_1}y^*)\\
&&\hskip50mm\vdots\\
&&\hskip10mm tr(x^*g_{i_p^{-1}i_1}xg_{k_1^{-1}k_p})tr(g_{k_1^{-1}k_p}^*yg_{j_1^{-1}j_p}y^*)dxdy
\end{eqnarray*}

We can compute the products of traces by using the following formula:
\begin{eqnarray*}
tr(Ag_k)tr(g_k^*B)
&=&\sum_{qs}<g_q,Ag_k><g_s,g_k^*B>\\
&=&\sum_{qs}tr(g_q^*Ag_k)tr(g_s^*g_k^*B)
\end{eqnarray*}

Thus are left with an integral involving the variable $z=xy$, which gives $T_p^\circ$.
\end{proof}

Let us discuss now the relationship with the Weyl matrices. We have:

\begin{proposition}
Given a finite abelian group $H$, consider the product $G=H\times\widehat{H}$, and endow it with its standard Fourier cocycle. 
\begin{enumerate}
\item With $E=U_n$, where $n=|H|$, the model $\pi:C(S_N^+)\to M_N(C(U_n))$ constructed above, where $N=n^2$, is the Weyl matrix model associated to $H$.

\item When assuming in addition that $H$ is cyclic, $H=\mathbb Z_n$, we obtain in this way the matrix model for $C(S_N^+)$ coming from the usual Weyl matrices.

\item In the particular case $H=\mathbb Z_2$, the model $\pi:C(S_4^+)\to M_N(C(U_2))$ constructed above is the matrix model for $C(S_4^+)$ coming from the Pauli matrices.
\end{enumerate}
\end{proposition}

\begin{proof}
All this is well-known. The general construction in Proposition 15.21 came in fact by successively generalizing $(3)\to(2)\to(1)$, and then by performing one more generalization, with $G=H\times\widehat{H}$ with its standard Fourier cocycle being replaced by an arbitrary finite group $G$, with a 2-cocycle on it. For full details here, see \cite{bne}.
\end{proof}

Regarding now the associated quantum permutation groups, in the general context of Proposition 15.21, we have the following result:

\index{cocyclic Weyl model}

\begin{theorem}
For a generalized Weyl matrix model, as in Proposition 15.21, the moments of the main character of the associated quantum group are
$$c_p=\frac{1}{N}\sum^\circ_{j_1\ldots j_p}
\int_Etr(g_{j_1}xg_{j_1}^*x^*)\ldots tr(g_{j_p}xg_{j_p}^*x^*)dx$$
where $\circ$ means that the indices are subject to the condition $j_1\ldots j_p=1$.
\end{theorem}

\begin{proof}
According to Proposition 15.21 and to Proposition 15.22, the moments of the main character are the following numbers:
\begin{eqnarray*}
c_p
&=&\frac{1}{N}\sum_{i_1\ldots i_p}\overline{\sigma(i_1,i_1^{-1}i_2)}\ldots\overline{\sigma(i_p,i_p^{-1}i_1)}\cdot\overline{\sigma(i_2,i_2^{-1}i_1)}\ldots\overline{\sigma(i_1,i_1^{-1}i_p)}\\
&&\int_Etr(g_{i_1^{-1}i_2}xg_{i_2^{-1}i_1}x^*)\ldots\ldots tr(g_{i_p^{-1}i_1}xg_{i_1^{-1}i_p}x^*)dx
\end{eqnarray*}

We can compact the cocycle part by using the following formulae:
\begin{eqnarray*}
\sigma(i_p,i_p^{-1}i_{p+1})\sigma(i_{p+1},i_{p+1}^{-1}i_p)
&=&\sigma(i_{p+1},i_{p+1}^{-1}i_p\cdot i_p^{-1}i_{p+1})\sigma(i_{p+1}^{-1}i_p,i_p^{-1}i_{p+1})\\
&=&\sigma(i_{p+1},1)\sigma(i_{p+1}^{-1}i_p,i_p^{-1}i_{p+1})\\
&=&\sigma(i_{p+1}^{-1}i_p,i_p^{-1}i_{p+1})
\end{eqnarray*}

Thus, in terms of the indices $j_1=i_1^{-1}i_2,\ldots,j_p=i_p^{-1}i_1$, which are subject to the condition $j_1\ldots j_p=1$, we have the following formula:
$$c_p=\frac{1}{N}\sum^\circ_{j_1\ldots j_p}\overline{\sigma(j_1^{-1},j_1)}\ldots\overline{\sigma(j_p^{-1},j_p)}
\int_Etr(g_{j_1}xg_{j_1^{-1}}x^*)\ldots tr(g_{j_p}xg_{j_p^{-1}}x^*)dx$$

Here the $\circ$ symbol above the sum is there for reminding us that the indices are subject to the condition $j_1\ldots j_p=1$. We can use now the following formula:
$$g_j^*=\overline{\sigma(j^{-1},j)}g_{j^{-1}}$$

We therefore obtain the following formula:
$$c_p=\frac{1}{N}\sum^\circ_{j_1\ldots j_p}
\int_Etr(g_{j_1}xg_{j_1}^*x^*)\ldots tr(g_{j_p}xg_{j_p}^*x^*)dx$$

Thus, we have obtained the formula in the statement.
\end{proof}

It is quite unclear whether the above formula further simplifies, in general. In the context of the Fourier cocycles, as in Proposition 15.24, it is possible to pass to a plain sum, by inserting a certain product of multiplicative factors $c(j_1)\ldots c(j_p)$, which equals $1$ when $j_1\ldots j_p=1$, and the computation can be finished as follows:
\begin{eqnarray*}
c_p&=&\frac{1}{N}\int_E\left(\sum_jc(j)tr(g_jxg_j^*x^*)\right)^pdx\\
&=&\frac{1}{N}\int_Etr(xx^*)dx
\end{eqnarray*}

Thus, the law of the main character of the corresponding quantum group coincides with the law of the main character of $PE$. All this suggests that the quantum group associated to a Weyl matrix model, as above, should appear as a suitable twist of $PE$. 

\bigskip

In addition, we believe that in the case where $E$ is easy these examples should be covered by a suitable projective extension of the Schur-Weyl twisting procedure.

\section*{15e. Exercises}

Things have been quite technical in this chapter, and as an instructive exercise, coming as a complement to the various questions raised above, we have:

\begin{exercise}
Develop of theory of matrix models for the ADE subgroups 
$$G\subset S_4^+$$
by suitably restricting the parameter space of the Pauli representation.
\end{exercise}

To be more precise, assuming that $G=H'$ comes as a twist of a subgroup $H\subset SO_3$, the problem is that of proving that when restricting to $H$ the model space of the Pauli representation, the Hopf image of the representation that we obtain is $C(G)$.

\chapter{Fourier models}

\section*{16a. Hadamard matrices}

Time to finish this book, and we kept the best for the end. Following \cite{bb3}, \cite{bi4} and related papers, we will discuss now the Hadamard matrix models, which came first, historically speaking, and which are of particular importance. Let us start with:

\index{Hadamard matrix}

\begin{definition}
A complex Hadamard matrix is a square matrix 
$$H\in M_N(\mathbb C)$$
whose entries are on the unit circle, and whose rows are pairwise orthogonal.
\end{definition}

Observe that the orthogonality condition tells us that the rescaled matrix $U=H/\sqrt{N}$ must be unitary. Thus, these matrices form a real algebraic manifold, given by:
$$X_N=M_N(\mathbb T)\cap\sqrt{N}U_N$$

\index{Fourier matrix}

The basic example is the Fourier matrix, $F_N=(w^{ij})$ with $w=e^{2\pi i/N}$.  With indices $i,j\in\{0,1,\ldots,N-1\}$, this matrix is as follows:
$$F_N=\begin{pmatrix}
1&1&1&\ldots&1\\
1&w&w^2&\ldots&w^{N-1}\\
1&w^2&w^4&\ldots&w^{2(N-1)}\\
\vdots&\vdots&\vdots&&\vdots\\
1&w^{N-1}&w^{2(N-1)}&\ldots&w^{(N-1)^2}
\end{pmatrix}$$

More generally, we have as example the Fourier coupling of any finite abelian group $G$, regarded via the isomorphism $G\simeq\widehat{G}$ as a square matrix, as follows:

\begin{theorem}
Given a finite abelian group $G$, with dual group $\widehat{G}=\{\chi:G\to\mathbb T\}$, consider the Fourier coupling $\mathcal F_G:G\times\widehat{G}\to\mathbb T$, given by $(i,\chi)\to\chi(i)$.
\begin{enumerate}
\item Via the standard isomorphism $G\simeq\widehat{G}$, this Fourier coupling can be regarded as a square matrix, $F_G\in M_G(\mathbb T)$, which is a complex Hadamard matrix.

\item For the cyclic group $G=\mathbb Z_N$ we obtain in this way, via the standard identification $\mathbb Z_N=\{1,\ldots,N\}$, the standard Fourier matrix, $F_N=(w^{ij})$ with $w=e^{2\pi i/N}$.

\item In general, when using a decomposition $G=\mathbb Z_{N_1}\times\ldots\times\mathbb Z_{N_k}$, the corresponding Fourier matrix is given by $F_G=F_{N_1}\otimes\ldots\otimes F_{N_k}$.
\end{enumerate}
\end{theorem}

\begin{proof}
This follows indeed from some basic facts from group theory:

\medskip

(1) With the identification $G\simeq\widehat{G}$ made our matrix is given by $(F_G)_{i\chi}=\chi(i)$, and the scalar products between the rows are computed as follows:
\begin{eqnarray*}
<R_i,R_j>
&=&\sum_\chi\chi(i)\overline{\chi(j)}\\
&=&\sum_\chi\chi(i-j)\\
&=&|G|\cdot\delta_{ij}
\end{eqnarray*}

Thus, we obtain indeed a complex Hadamard matrix.

\medskip

(2) This follows from the well-known and elementary fact that, via the identifications $\mathbb Z_N=\widehat{\mathbb Z_N}=\{1,\ldots,N\}$, the Fourier coupling here is as follows, with $w=e^{2\pi i/N}$:
$$(i,j)\to w^{ij}$$

(3) We use here the following well-known formula, for the duals of products: 
$$\widehat{H\times K}=\widehat{H}\times\widehat{K}$$

At the level of the corresponding Fourier couplings, we obtain from this:
$$F_{H\times K}=F_H\otimes F_K$$

Now by decomposing $G$ into cyclic groups, as in the statement, and by using (2) for the cyclic components, we obtain the formula in the statement.
\end{proof}

In relation with the quantum groups, the starting observation is as follows:

\begin{proposition}
If $H\in M_N(\mathbb C)$ is Hadamard, the rank one projections 
$$P_{ij}=Proj\left(\frac{H_i}{H_j}\right)$$
where $H_1,\ldots,H_N\in\mathbb T^N$ are the rows of $H$, form a magic unitary.
\end{proposition}

\begin{proof}
This is clear, the verification for the rows being as follows:
\begin{eqnarray*}
\left<\frac{H_i}{H_j},\frac{H_i}{H_k}\right>
&=&\sum_l\frac{H_{il}}{H_{jl}}\cdot\frac{H_{kl}}{H_{il}}\\
&=&\sum_l\frac{H_{kl}}{H_{jl}}\\
&=&N\delta_{jk}
\end{eqnarray*}

The verification for the columns is similar, as follows:
\begin{eqnarray*}
\left<\frac{H_i}{H_j},\frac{H_k}{H_j}\right>
&=&\sum_l\frac{H_{il}}{H_{jl}}\cdot\frac{H_{jl}}{H_{kl}}\\
&=&\sum_l\frac{H_{il}}{H_{kl}}\\
&=&N\delta_{ik}
\end{eqnarray*}

Thus, we have indeed a magic unitary, as claimed. 
\end{proof}

We can proceed now in the same way as we did with the Weyl matrices, namely by constructing a model of $C(S_N^+)$, and performing the Hopf image construction:

\begin{definition}
To any Hadamard matrix $H\in M_N(\mathbb C)$ we associate the quantum permutation group $G\subset S_N^+$ given by the fact that $C(G)$ is the Hopf image of
$$\pi:C(S_N^+)\to M_N(\mathbb C)\quad,\quad 
u_{ij}\to Proj\left(\frac{H_i}{H_j}\right)$$
where $H_1,\ldots,H_N\in\mathbb T^N$ are the rows of $H$.
\end{definition}

Summarizing, we have a construction $H\to G$, and our claim is that this construction is something really useful, with $G$ encoding the combinatorics of $H$. To be more precise, our claim is that ``$H$ can be thought of as being a kind of Fourier matrix for $G$''. 

\bigskip

This is of course quite interesting, philosophically speaking. There are several results supporting this, with the main evidence coming from the following result, coming from \cite{bbs}, which collects the basic known results regarding the construction:

\begin{theorem}
The construction $H\to G$ has the following properties:
\begin{enumerate}
\item For a Fourier matrix $H=F_G$ we obtain the group $G$ itself, acting on itself.

\item For $H\not\in\{F_G\}$, the quantum group $G$ is not classical, nor a group dual.

\item For a tensor product $H=H'\otimes H''$ we obtain a product, $G=G'\times G''$.
\end{enumerate}
\end{theorem}

\begin{proof}
All this material is standard, and elementary, as follows:

\medskip

(1) Let us first discuss the cyclic group case, where our Hadamard matrix is a usual Fourier matrix, $H=F_N$. Here the rows of $H$ are given by $H_i=\rho^i$, where:
$$\rho=(1,w,w^2,\ldots,w^{N-1})$$

Thus, we have the following formula, for the associated magic basis:
$$\frac{H_i}{H_j}=\rho^{i-j}$$

It follows that the corresponding rank 1 projections $P_{ij}=Proj(H_i/H_j)$ form a circulant matrix, all whose entries commute. Since the entries commute, the corresponding quantum group must satisfy $G\subset S_N$. Now by taking into account the circulant property of $P=(P_{ij})$ as well, we are led to the conclusion that we have:
$$G=\mathbb Z_N$$

In the general case now, where $H=F_G$, with $G$ being an arbitrary finite abelian group, the result can be proved either by extending the above proof, of by decomposing $G=\mathbb Z_{N_1}\times\ldots\times\mathbb Z_{N_k}$ and using (3) below, whose proof is independent from the rest.

\medskip

(2) This is something more tricky, needing some general study of the representations whose Hopf images are commutative, or cocommutative.

\medskip

(3) Assume that we have a tensor product $H=H'\otimes H''$, and let $G,G',G''$ be the associated quantum permutation groups. We have then a diagram as follows:
$$\xymatrix@R=45pt@C25pt{
C(S_{N'}^+)\otimes C(S_{N''}^+)\ar[r]&C(G')\otimes C(G'')\ar[r]&M_{N'}(\mathbb C)\otimes M_{N''}(\mathbb C)\ar[d]\\
C(S_N^+)\ar[u]\ar[r]&C(G)\ar[r]&M_N(\mathbb C)
}$$

Here all the maps are the canonical ones, with those on the left and on the right coming from $N=N'N''$. At the level of standard generators, the diagram is as follows:
$$\xymatrix@R=45pt@C65pt{
u_{ij}'\otimes u_{ab}''\ar[r]&w_{ij}'\otimes w_{ab}''\ar[r]&P_{ij}'\otimes P_{ab}''\ar[d]\\
u_{ia,jb}\ar[u]\ar[r]&w_{ia,jb}\ar[r]&P_{ia,jb}
}$$

Now observe that this diagram commutes. We conclude that the representation associated to $H$ factorizes indeed through $C(G')\otimes C(G'')$, and this gives the result.
\end{proof}

At a more abstract level, one interesting question is that of abstractly characterizing the magic matrices coming from the complex Hadamard matrices. We have here:

\begin{proposition}
Given an Hadamard matrix $H\in M_N(\mathbb C)$, the vectors 
$$\xi_{ij}=\frac{H_i}{H_j}$$
on which the magic unitary entries $P_{ij}$ project, have the following properties:
\begin{enumerate}
\item $\xi_{ii}=\xi$ is the all-one vector.

\item $\xi_{ij}\xi_{jk}=\xi_{ik}$, for any $i,j,k$.

\item $\xi_{ij}\xi_{kl}=\xi_{il}\xi_{kj}$, for any $i,j,k,l$.
\end{enumerate}
\end{proposition}

\begin{proof}
All these assertions are trivial, by using the formula $\xi_{ij}=H_i/H_j$.
\end{proof}

Let us call now magic basis of a given Hilbert space $H$ any square array of vectors $\xi\in M_N(H)$, all whose rows and columns are orthogonal bases of $H$. With this convention, the above observations lead to the following result, at the magic basis level:

\index{magic basis}

\begin{theorem}
The magic bases $\xi\in M_N(S^{N-1}_\mathbb C)$ coming from the complex Hadamard matrices are those having the following properties:
\begin{enumerate}
\item We have $\xi_{ij}\in\mathbb T^N$, after a suitable rescaling. 

\item The conditions in Proposition 16.6 are satisfied.
\end{enumerate}
\end{theorem}

\begin{proof}
By using the multiplicativity conditions (1,2,3) in Proposition 16.6, we conclude that, up to a rescaling, we must have $\xi_{ij}=\xi_i/\xi_j$, where $\xi_1,\ldots,\xi_N$ is the first row of the magic basis. Together with our assumption $\xi_{ij}\in\mathbb T^N$, this gives the result. 
\end{proof}

At the general level now, regarding the representation theory of the quantum groups associated to the Hadamard matrices, we have the following result, from \cite{bbs}:

\index{Tannakian category}

\begin{theorem}
The Tannakian category of the quantum group $G\subset S_N^+$ associated to a complex Hadamard matrix $H\in M_N(\mathbb C)$ is given by
$$T\in Hom(u^{\otimes k},u^{\otimes l})\iff T^\circ G^{k+2}=G^{l+2}T^\circ$$
where the objects on the right are constructed as follows:
\begin{enumerate}
\item $T^\circ=id\otimes T\otimes id$.

\item $G_{ia}^{jb}=\sum_kH_{ik}\bar{H}_{jk}\bar{H}_{ak}H_{bk}$.

\item $G^k_{i_1\ldots i_k,j_1\ldots j_k}=G_{i_ki_{k-1}}^{j_kj_{k-1}}\ldots G_{i_2i_1}^{j_2j_1}$.
\end{enumerate}
\end{theorem}

\begin{proof}
We use the Tannakian result for the Hopf image of a representation, discussed in chapter 13. With the notations here, we have the following formula:
$$Hom(u^{\otimes k},u^{\otimes l})=Hom(U^{\otimes k},U^{\otimes l})$$

The vector space on the right consists by definition of the complex $N^l\times N^k$ matrices $T$, satisfying the following relation:
$$TU^{\otimes k}=U^{\otimes l}T$$ 

If we denote this equality by $L=R$, the left term $L$ is given by:
\begin{eqnarray*}
L_{ij}
&=&(TU^{\otimes k})_{ij}\\
&=&\sum_aT_{ia}U^{\otimes k}_{aj}\\
&=&\sum_aT_{ia}U_{a_1j_1}\ldots U_{a_kj_k}
\end{eqnarray*}

As for the right term $R$, this is given by:
\begin{eqnarray*}
R_{ij}
&=&(U^{\otimes l}T)_{ij}\\
&=&\sum_bU^{\otimes l}_{ib}T_{bj}\\
&=&\sum_bU_{i_1b_1}\ldots U_{i_lb_l}T_{bj}
\end{eqnarray*}

Consider now the vectors $\xi_{ij}=H_i/H_j$. Since these vectors span the ambient Hilbert space, the equality $L=R$ is equivalent to the following equality:
$$<L_{ij}\xi_{pq},\xi_{rs}>=<R_{ij}\xi_{pq},\xi_{rs}>$$

We use now the following well-known formula, expressing a product of rank one projections $P_1,\ldots,P_k$ in terms of the corresponding image vectors $\xi_1,\ldots,\xi_k$:
$$<P_1\ldots P_kx,y>=<x,\xi_k><\xi_k,\xi_{k-1}>\ldots\ldots<\xi_2,\xi_1><\xi_1,y>$$

This gives the following formula for $L$:
\begin{eqnarray*}
<L_{ij}\xi_{pq},\xi_{rs}>
&=&\sum_aT_{ia}<P_{a_1j_1}\ldots P_{a_kj_k}\xi_{pq},\xi_{rs}>\\
&=&\sum_aT_{ia}<\xi_{pq},\xi_{a_kj_k}>\ldots<\xi_{a_1j_1},\xi_{rs}>\\
&=&\sum_aT_{ia}G_{pa_k}^{qj_k}G_{a_ka_{k-1}}^{j_kj_{k-1}}\ldots G_{a_2a_1}^{j_2j_1}G_{a_1r}^{j_1s}\\
&=&\sum_aT_{ia}G^{k+2}_{rap,sjq}\\
&=&(T^\circ G^{k+2})_{rip,sjq}
\end{eqnarray*}

As for the right term $R$, this is given by:
\begin{eqnarray*}
<R_{ij}\xi_{pq},\xi_{rs}>
&=&\sum_b<P_{i_1b_1}\ldots P_{i_lb_l}\xi_{pq},\xi_{rs}>T_{bj}\\
&=&\sum_b<\xi_{pq},\xi_{i_lb_l}>\ldots<\xi_{i_1b_1},\xi_{rs}>T_{bj}\\
&=&\sum_bG_{pi_l}^{qb_l}G_{i_li_{l-1}}^{b_lb_{l-1}}\ldots G_{i_2i_1}^{b_2b_1}G_{i_1r}^{b_1s}T_{bj}\\
&=&\sum_bG^{l+2}_{rip,sbq}T_{bj}\\
&=&(G^{l+2}T^\circ)_{rip,sjq}
\end{eqnarray*}

Thus, we obtain the formula in the statement. See \cite{bbs}.
\end{proof}

\section*{16b. Spin models}

Let us discuss now the relation with subfactor theory, and with planar algebras. As a starting point, we have the following basic observation of Popa \cite{pop}:

\index{orthogonal MASA}

\begin{proposition}
Up to a conjugation by a unitary, the pairs of orthogonal MASA in the simplest factor, namely the matrix algebra $M_N(\mathbb C)$, are as follows,
$$A=\Delta\quad,\quad 
B=H\Delta H^*$$
with $\Delta\subset M_N(\mathbb C)$ being the diagonal matrices, and with $H\in M_N(\mathbb C)$ being Hadamard.
\end{proposition}

\begin{proof}
Any maximal abelian subalgebra (MASA) in $M_N(\mathbb C)$ being conjugated to $\Delta$, we can assume, up to conjugation by a unitary, that we have, with $U\in U_N$:
$$A=\Delta\quad,\quad 
B=U\Delta U^*$$  

Now observe that given two diagonal matrices $D,E\in\Delta$, we have:
\begin{eqnarray*}
tr(D\cdot UEU^*)
&=&\frac{1}{N}\sum_i(DUEU^*)_{ii}\\
&=&\frac{1}{N}\sum_{ij}D_{ii}U_{ij}E_{jj}\bar{U}_{ij}\\
&=&\frac{1}{N}\sum_{ij}D_{ii}E_{jj}|U_{ij}|^2
\end{eqnarray*}

Thus, the orthogonality condition $A\perp B$ reformulates as follows:
$$\frac{1}{N}\sum_{ij}D_{ii}E_{jj}|U_{ij}|^2=\frac{1}{N^2}\sum_{ij}D_{ii}E_{jj}$$

Thus, we must have $|U_{ij}|=\frac{1}{\sqrt{N}}$, for any $i,j$. But this tells us that the matrix $H=\sqrt{N}U$ must be Hadamard, as claimed.
\end{proof}

Along the same lines, but at a more advanced level, we have:

\index{commuting square}

\begin{theorem}
Given a complex Hadamard matrix $H\in M_N(\mathbb C)$, the diagram formed by the associated pair of orthogonal MASA, namely
$$\xymatrix@R=35pt@C35pt{
\Delta\ar[r]&M_N(\mathbb C)\\
\mathbb C\ar[u]\ar[r]&H\Delta H^*\ar[u] }$$ is a commuting square in the sense of subfactor theory, in the sense that the expectations onto $\Delta,H\Delta H^*$ commute, and their product is the expectation onto $\mathbb C$.
\end{theorem}

\begin{proof}
The expectation $E_\Delta:M_N(\mathbb C)\to\Delta$ is the operation $M\to M_\Delta$ which consists in keeping the diagonal, and erasing the rest. Consider now the other expectation:
$$E_{H\Delta H^*}:M_N(\mathbb C)\to H\Delta H^*$$

It is better to identify this with the following expectation, with $U=H/\sqrt{N}$: 
$$E_{U\Delta U^*}:M_N(\mathbb C)\to U\Delta U^*$$

This latter expectation must be of the form $M\to UX_\Delta U^*$, with $X$ satisfying:
$$<M,UDU^*>=<UX_\Delta U^*,UDU^*>\quad,\quad\forall D\in\Delta$$

The scalar products being given by $<a,b>=tr(ab^*)$, this condition reads:
$$tr(MUD^*U^*)=tr(X_\Delta D^*)\quad,\quad\forall D\in\Delta$$

Thus $X=U^*MU$, and the formulae of our two expectations are as follows:
\begin{eqnarray*}
E_\Delta(M)&=&M_\Delta\\
E_{U\Delta U^*}(M)&=&U(U^*MU)_\Delta U^*
\end{eqnarray*}

With these formulae in hand, an elementary computation gives the result.
\end{proof}

The point now is that any commuting square $C$ produces a subfactor of the Murray-von Neumann hyperfinite ${\rm II}_1$ factor $R$. Consider indeed such a square:
$$\xymatrix@R=35pt@C35pt{
C_{01}\ar[r]&C_{11}\\
C_{00}\ar[u]\ar[r]&C_{10}\ar[u]}$$

Under suitable assumptions on the inclusions $C_{00}\subset C_{10},C_{01}\subset C_{11}$, we can perform the basic construction for them, in finite dimensions, and we obtain in this way a whole array of commuting squares, as follows:
$$\xymatrix@R=35pt@C35pt{
A_0&A_1&A_2&\\
C_{02}\ar[r]\ar@.[u]&C_{12}\ar[r]\ar@.[u]&C_{22}\ar@.[r]\ar@.[u]&B_2\\
C_{01}\ar[r]\ar[u]&C_{11}\ar[r]\ar[u]&C_{21}\ar@.[r]\ar[u]&B_1\\
C_{00}\ar[u]\ar[r]&C_{10}\ar[u]\ar[r]&C_{20}\ar[u]\ar@.[r]&B_0}$$

Here the various $A,B$ letters stand for the von Neumann algebras obtained in the limit, which are all isomorphic to the hyperfinite ${\rm II}_1$ factor $R$, and we have:

\index{Ocneanu compactness}

\begin{theorem}
In the context of the above diagram, the following happen:
\begin{enumerate}
\item $A_0\subset A_1$ is a subfactor, and $\{A_i\}$ is the Jones tower for it.

\item The corresponding planar algebra is given by $A_0'\cap A_k=C_{01}'\cap C_{k0}$.

\item A similar result holds for the ``horizontal'' subfactor $B_0\subset B_1$.
\end{enumerate}
\end{theorem}

\begin{proof}
This is something very standard, the idea being as follows:

\medskip

(1) This is something quite routine.

\medskip

(2) This is a subtle result, called Ocneanu compactness theorem.

\medskip

(3) This follows from (1,2), by flipping the diagram.
\end{proof}

Getting back now to the Hadamard matrices, we can extend our lineup of results, namely Proposition 16.9 and Theorem 16.10, with an advanced result, as follows:

\begin{theorem}
Given a complex Hadamard matrix $H\in M_N(\mathbb C)$, the diagram formed by the associated pair of orthogonal MASA is a commuting square in the sense of subfactor theory, and the associated planar algebra $P=(P_k)$ is given by the following formula,
$$T\in P_k\iff T^\circ G^2=G^{k+2}T^\circ$$
where the objects on the right are constructed as follows:
\begin{enumerate}
\item $T^\circ=id\otimes T\otimes id$.

\item $G_{ia}^{jb}=\sum_kH_{ik}\bar{H}_{jk}\bar{H}_{ak}H_{bk}$.

\item $G^k_{i_1\ldots i_k,j_1\ldots j_k}=G_{i_ki_{k-1}}^{j_kj_{k-1}}\ldots G_{i_2i_1}^{j_2j_1}$.
\end{enumerate}
\end{theorem}

\begin{proof}
We have two assertions here, the idea being as follows:

\medskip

(1) The fact that we have indeed a commuting square is something that we already know, from Theorem 16.10.

\medskip

(2) The computation of the associated planar algebra is possible thanks to the Ocneanu compactness theorem, corresponding to the formula in Theorem 16.11 (2). To be more precise, by doing some direct computations, which are quite similar to those in the proof of Theorem 16.8, we obtain the formula in the statement. See \cite{jo3}.
\end{proof}

The point now is that all the above is very similar to Theorem 16.8, and we have: 

\index{spin model}

\begin{theorem}
Let $H\in M_N(\mathbb C)$ be a complex Hadamard matrix.
\begin{enumerate}
\item The planar algebra associated to $H$ is given by $P_k=Fix(u^{\otimes k})$, where $G\subset S_N^+$ is the associated quantum permutation group.

\item The corresponding Poincar\'e series $f(z)=\sum_k\dim(P_k)z^k$ equals the Stieltjes transform $\int_G\frac{1}{1-z\chi}$ of the law of the main character $\chi=\sum_iu_{ii}$.
\end{enumerate}
\end{theorem}

\begin{proof}
This follows by comparing the quantum group and subfactor results:

\medskip

(1) As already mentioned above, this simply follows by comparing Theorem 16.8 with the subfactor computation in Theorem 16.12. For full details here, we refer to \cite{bb3}.

\medskip

(2) This is a consequence of (1), and of the Peter-Weyl type results from \cite{wo1}, which tell us that fixed points can be counted by integrating characters.
\end{proof}

Summarizing, we have now a clarification of the various quantum algebraic objects associated to a complex Hadamard matrix $H\in M_N(\mathbb C)$, the idea being that the central object, which best encodes the ``symmetries'' of the matrix, and which allows the computation of the other quantum algebraic objects as well, such as the associated planar algebra, is the associated quantum permutation group $G\subset S_N^+$. Regarding now the subfactor itself, the result here is as follows:

\begin{theorem}
The subfactor associated to $H\in M_N(\mathbb C)$ is of the form
$$A^G\subset(\mathbb C^N\otimes A)^G$$
with $A=R\rtimes\widehat{G}$, where $G\subset S_N^+$ is the associated quantum permutation group.
\end{theorem}

\begin{proof}
This is something more technical, the idea being that the basic construction procedure for the commuting squares, explained before Theorem 16.11, can be performed in an ``equivariant setting'', for commuting squares having components as follows:
$$D\otimes_GE=(D\otimes(E\rtimes\widehat{G}))^G$$

To be more precise, starting with a commuting square formed by such algebras, we obtain by basic construction a whole array of commuting squares as follows, with $\{D_i\},\{E_i\}$ being by definition Jones towers, and with $D_\infty,E_\infty$ being their inductive limits:
$$\xymatrix@R=35pt@C35pt{
D_0\otimes_GE_\infty&D_1\otimes_GE_\infty&D_2\otimes_GE_\infty\\
D_0\otimes_GE_2\ar@.[u]\ar[r]&D_1\otimes_GE_2\ar@.[u]\ar[r]&D_2\otimes_GE_2\ar@.[u]\ar@.[r]&D_\infty\otimes_GE_2\\
D_0\otimes_GE_1\ar[u]\ar[r]&D_1\otimes_GE_1\ar[u]\ar[r]&D_2\otimes_GE_1\ar[u]\ar@.[r]&D_\infty\otimes_GE_1\\
D_0\otimes_GE_0\ar[u]\ar[r]&D_1\otimes_GE_0\ar[u]\ar[r]&D_2\otimes_GE_0\ar[u]\ar@.[r]&D_\infty\otimes_GE_0}$$

The point now is that this quantum group picture works in fact for any commuting square having $\mathbb C$ in the lower left corner. In the Hadamard matrix case, that we are interested in here, the corresponding commuting square is as follows:
$$\xymatrix@R=35pt@C35pt{
\mathbb C\otimes_G\mathbb C^N\ar[r]&\mathbb C^N\otimes_G\mathbb C^N\\
\mathbb C\otimes_G\mathbb C\ar[u]\ar[r]&\mathbb C^N\otimes_G\mathbb C\ar[u] }$$ 

Thus, the subfactor obtained by vertical basic construction appears as follows:
$$\mathbb C\otimes_GE_\infty\subset\mathbb C^N\otimes_GE_\infty$$

But this gives the conclusion in the statement, with the ${\rm II}_1$ factor appearing there being by definition $A=E_\infty\rtimes\widehat{G}$, and with the remark that we have $E_\infty\simeq R$. See \cite{ba4}.
\end{proof}

There are many other things that can be said about the Hadamard matrices and their relation with operator algebras, and about the subfactors associated to the Hadamard matrices, and their generalizations, and we refer to \cite{ba3}, \cite{bb3}, \cite{jo3} and related papers.

\section*{16c. Fourier models}

Getting back now to Theorem 16.5, going beyond it is a quite delicate task. The next simplest models appear by deforming the Fourier matrices, or rather the tensor products of such matrices, $F_{G\times H}=F_G\otimes F_H$, via the following construction, due to Di\c t\u a:

\index{deformed Fourier matrix}

\begin{definition}
Given two finite abelian groups $G,H$, we consider the corresponding deformed Fourier matrix, given by the formula
$$(F_G\otimes_Q F_H)_{ia,jb}=Q_{ib}(F_G)_{ij}(F_H)_{ab}$$
and we factorize the associated representation $\pi_Q$ of the algebra $C(S_{G\times H}^+)$,
$$\xymatrix@R=40pt@C=40pt
{C(S_{G\times H}^+)\ar[rr]^{\pi_Q}\ar[rd]&&M_{G\times H}(\mathbb C)\\&C(G_Q)\ar[ur]_\pi&}$$
with $C(G_Q)$ being the Hopf image of this representation $\pi_Q$.
\end{definition}

Explicitely computing the above quantum permutation group $G_Q\subset S_{G\times H}^+$, as function of the parameter matrix $Q\in M_{G\times H}(\mathbb T)$, will be our main purpose, in what follows. In order to do so, we first have the following elementary result:

\begin{proposition}
We have a factorization as follows,
$$\xymatrix@R=40pt@C=40pt
{C(S_{G\times H}^+)\ar[rr]^{\pi_Q}\ar[rd]&&M_{G\times H}(\mathbb C)\\&C(H\wr_*G)\ar[ur]_\pi&}$$
given on the standard generators by the formulae
$$U_{ab}^{(i)}=\sum_jW_{ia,jb}\quad,\quad 
V_{ij}=\sum_aW_{ia,jb}$$
independently of $b$, where $W$ is the magic matrix producing $\pi_Q$.
\end{proposition}

\begin{proof}
With $K=F_G,L=F_H$ and $M=|G|,N=|H|$, the formula of the magic matrix $W\in M_{G\times H}(M_{G\times H}(\mathbb C))$ associated to $H=K\otimes_QL$ is as follows:
\begin{eqnarray*}
(W_{ia,jb})_{kc,ld}
&=&\frac{1}{MN}\cdot\frac{Q_{ic}Q_{jd}}{Q_{id}Q_{jc}}\cdot\frac{K_{ik}K_{jl}}{K_{il}K_{jk}}\cdot\frac{L_{ac}L_{bd}}{L_{ad}L_{bc}}\\
&=&\frac{1}{MN}\cdot\frac{Q_{ic}Q_{jd}}{Q_{id}Q_{jc}}\cdot K_{i-j,k-l}L_{a-b,c-d}
\end{eqnarray*}

Our claim now is that the representation $\pi_Q$ constructed in Definition 16.15 can be factorized in three steps, up to the factorization in the statement, as follows:
$$\xymatrix@R=60pt@C=50pt
{C(S_{G\times H}^+)\ar[rr]^{\pi_Q}\ar[d]&&M_{G\times H}(\mathbb C)\\
C(S_H^+\wr_*S_G^+)\ar[r]\ar@{.>}[rru]&C(S_H^+\wr_*G)\ar[r]\ar@{.>}[ur]&C(H\wr_*G)\ar@{.>}[u]}$$

Indeed, the construction of the map on the left is standard. Regarding the second factorization, this comes from the fact that since the elements $V_{ij}$ depend on $i-j$, they satisfy the defining relations for the quotient algebra $C(S_G^+)\to C(G)$. Finally, regarding the third factorization, observe that $W_{ia,jb}$ depends only on $i,j$ and on $a-b$. By summing over $j$ we obtain that the elements $U_{ab}^{(i)}$ depend only on $a-b$, and we are done.
\end{proof}

We have now all needed ingredients for refining Proposition 16.16, as follows:

\begin{proposition}
We have a factorization as follows,
$$\xymatrix@R=40pt@C=30pt
{C(S_{G\times H}^+)\ar[rr]^{\pi_Q}\ar[rd]&&M_{G\times H}(\mathbb C)\\&C^*(\Gamma_{G,H})\rtimes C(G)\ar[ur]_\rho&}$$
where the group on the bottom is given by:
$$\Gamma_{G,H}=H^{*G}\Big/\left<[c_1^{(i_1)}\ldots c_s^{(i_s)},d_1^{(j_1)}\ldots d_s^{(j_s)}]=1\Big|\sum_rc_r=\sum_rd_r=0\right>$$
\end{proposition}

\begin{proof}
Assume that we have a representation, as follows:
$$\pi:C^*(\Gamma)\rtimes C(G)\to M_L(\mathbb C)$$

Let $\Lambda$ be a $G$-stable normal subgroup of $\Gamma$, so that $G$ acts on $\Gamma/\Lambda$, and we can form the product $C^*(\Gamma/\Lambda)\rtimes C(G)$, and assume that $\pi$ is trivial on $\Lambda$. Then $\pi$ factorizes as:
$$\xymatrix@R=40pt@C=30pt
{C^*(\Gamma)\rtimes C(G)\ar[rr]^\pi\ar[rd]&&M_L(\mathbb C)\\&C^*(\Gamma/\Lambda)\rtimes C(G)\ar[ur]_\rho}$$

With $\Gamma=H^{*G}$, this gives the result.
\end{proof}

We have now all the needed ingredients for proving a main result, as follows:

\begin{theorem}
When $Q$ is generic, the minimal factorization for $\pi_Q$ is 
$$\xymatrix@R=45pt@C=40pt
{C(S_{G\times H}^+)\ar[rr]^{\pi_Q}\ar[rd]&&M_{G\times H}(\mathbb C)\\&C^*(\Gamma_{G,H})\rtimes C(G)\ar[ur]_\pi&}$$
where on the bottom
$$\Gamma_{G,H}\simeq\mathbb Z^{(|G|-1)(|H|-1)}\rtimes H$$
is the discrete group constructed above.
\end{theorem}

\begin{proof}
Consider the factorization in Proposition 16.17, which is as follows, where $L$ denotes the Hopf image of $\pi_Q$:
$$\theta : C^*(\Gamma_{G,H})\rtimes C(G)\to L$$

To be more precise, this morphism produces the following commutative diagram:
$$\xymatrix@R=40pt@C=40pt
{C(S_{G\times H}^+) \ar[rr]^{\pi_Q} \ar[dr]_{} \ar@/_/[ddr]_{}& & M_{G\times H}(\mathbb C) \\
& L \ar[ur]_{} & \\
& C^*(\Gamma_{G,H})\rtimes C(G) \ar@{-->}[u]_\theta \ar@/_/[uur]_{\pi}& 
}$$

The first observation is that the injectivity assumption on $C(G)$ holds by construction, and that for $f \in C(G)$, the matrix $\pi(f)$ is ``block scalar''. Now for $r \in \Gamma_{G,H}$ with $\theta(r\otimes 1)=\theta(1 \otimes f)$ for some $f \in C(G)$, we see, using the commutative diagram, that $\pi(r \otimes 1)$ is block scalar. Thus, modulo some standard algebra, we are done.
\end{proof}

Summarizing, we have computed the quantum permutation groups associated to the Di\c t\u a deformations of the tensor products of Fourier matrices, in the case where the deformation matrix $Q$ is generic. For some further computations, in the case where the deformation matrix $Q$ is no longer generic, we refer to the follow-ups of \cite{bb3}.

\bigskip

Let us compute now the Kesten measure $\mu=law(\chi)$, in the case where the deformation matrix is generic, as before. Our results here will be a combinatorial moment formula, a geometric interpretation of it, and an asymptotic result. We first have:

\begin{theorem}
We have the moment formula
$$\int\chi^p
=\frac{1}{|G|\cdot|H|}\#\left\{\begin{matrix}i_1,\ldots,i_p\in G\\ d_1,\ldots,d_p\in H\end{matrix}\Big|\begin{matrix}[(i_1,d_1),(i_2,d_2),\ldots,(i_p,d_p)]\ \ \ \ \\=[(i_1,d_p),(i_2,d_1),\ldots,(i_p,d_{p-1})]\end{matrix}\right\}$$
where the sets between square brackets are by definition sets with repetition.
\end{theorem}

\begin{proof}
According to the various formulae above, the factorization found in Theorem 16.18 is, at the level of standard generators, as follows:
$$\begin{matrix}
C(S_{G\times H}^+)&\to&C^*(\Gamma_{G,H})\otimes C(G)&\to&M_{G\times H}(\mathbb C)\\
u_{ia,jb}&\to&\frac{1}{|H|}\sum_cF_{b-a,c}c^{(i)}\otimes v_{ij}&\to&W_{ia,jb}
\end{matrix}$$

Thus, the main character of the quantum permutation group that we found in Theorem 16.18 is given by the following formula:
\begin{eqnarray*}
\chi
&=&\frac{1}{|H|}\sum_{iac}c^{(i)}\otimes v_{ii}\\
&=&\sum_{ic}c^{(i)}\otimes v_{ii}\\
&=&\left(\sum_{ic}c^{(i)}\right)\otimes\delta_1
\end{eqnarray*}

Now since the Haar functional of $C^*(\Gamma)\rtimes C(H)$ is the tensor product of the Haar functionals of $C^*(\Gamma),C(H)$, this gives the following formula, valid for any $p\geq1$:
$$\int\chi^p=\frac{1}{|G|}\int_{\widehat{\Gamma}_{G,H}}\left(\sum_{ic}c^{(i)}\right)^p$$

Consider the elements $S_i=\sum_cc^{(i)}$. With standard notations, we have:
$$S_i=\sum_c(b_{i0}-b_{ic},c)$$

Now observe that these elements multiply as follows:
$$S_{i_1}\ldots S_{i_p}=\sum_{c_1\ldots c_p}
\begin{pmatrix}
b_{i_10}-b_{i_1c_1}+b_{i_2c_1}-b_{i_2,c_1+c_2}&&\\
+b_{i_3,c_1+c_2}-b_{i_3,c_1+c_2+c_3}+\ldots\ldots&,&c_1+\ldots+c_p&\\
\ldots\ldots+b_{i_p,c_1+\ldots+c_{p-1}}-b_{i_p,c_1+\ldots+c_p}&&
\end{pmatrix}$$

In terms of the new indices $d_r=c_1+\ldots+c_r$, this formula becomes:
$$S_{i_1}\ldots S_{i_p}=\sum_{d_1\ldots d_p}
\begin{pmatrix}
b_{i_10}-b_{i_1d_1}+b_{i_2d_1}-b_{i_2d_2}&&\\
+b_{i_3d_2}-b_{i_3d_3}+\ldots\ldots&,&d_p&\\
\ldots\ldots+b_{i_pd_{p-1}}-b_{i_pd_p}&&
\end{pmatrix}$$

Now by integrating, we must have $d_p=0$ on one hand, and on the other hand:
$$[(i_1,0),(i_2,d_1),\ldots,(i_p,d_{p-1})]=[(i_1,d_1),(i_2,d_2),\ldots,(i_p,d_p)]$$

Equivalently, we must have $d_p=0$ on one hand, and on the other hand:
$$[(i_1,d_p),(i_2,d_1),\ldots,(i_p,d_{p-1})]=[(i_1,d_1),(i_2,d_2),\ldots,(i_p,d_p)]$$

Thus, by translation invariance with respect to $d_p$, we obtain:
$$\int_{\widehat{\Gamma}_{G,H}}S_{i_1}\ldots S_{i_p}
=\frac{1}{|H|}\#\left\{d_1,\ldots,d_p\in H\Big|\begin{matrix}[(i_1,d_1),(i_2,d_2),\ldots,(i_p,d_p)]\ \ \ \ \\=[(i_1,d_p),(i_2,d_1),\ldots,(i_p,d_{p-1})]\end{matrix}\right\}$$

It follows that we have the following moment formula:
$$\int_{\widehat{\Gamma}_{G,H}}\left(\sum_iS_i\right)^p
=\frac{1}{|H|}\#\left\{\begin{matrix}i_1,\ldots,i_p\in G\\ d_1,\ldots,d_p\in H\end{matrix}\Big|\begin{matrix}[(i_1,d_1),(i_2,d_2),\ldots,(i_p,d_p)]\ \ \ \ \\=[(i_1,d_p),(i_2,d_1),\ldots,(i_p,d_{p-1})]\end{matrix}\right\}$$

Now by dividing by $|G|$, we obtain the formula in the statement.
\end{proof} 

The formula in Theorem 16.19 can be further interpreted as follows:

\index{Gram matrix}

\begin{theorem}
With $M=|G|,N=|H|$ we have the formula
$$law(\chi)=\left(1-\frac{1}{N}\right)\delta_0+\frac{1}{N}law(A)$$
where the matrix
$$A\in C(\mathbb T^{MN},M_M(\mathbb C))$$
is given by $A(q)=$ Gram matrix of the rows of $q$.
\end{theorem}

\begin{proof}
According to Theorem 16.19, we have the following formula:
\begin{eqnarray*}
\int\chi^p
&=&\frac{1}{MN}\sum_{i_1\ldots i_p}\sum_{d_1\ldots d_p}\delta_{[i_1d_1,\ldots,i_pd_p],[i_1d_p,\ldots,i_pd_{p-1}]}\\
&=&\frac{1}{MN}\int_{\mathbb T^{MN}}\sum_{i_1\ldots i_p}\sum_{d_1\ldots d_p}\frac{q_{i_1d_1}\ldots q_{i_pd_p}}{q_{i_1d_p}\ldots q_{i_pd_{p-1}}}\,dq\\
&=&\frac{1}{MN}\int_{\mathbb T^{MN}}\sum_{i_1\ldots i_p}\left(\sum_{d_1}\frac{q_{i_1d_1}}{q_{i_2d_1}}\right)\left(\sum_{d_2}\frac{q_{i_2d_2}}{q_{i_3d_2}}\right)\ldots\left(\sum_{d_p}\frac{q_{i_pd_p}}{q_{i_1d_p}}\right)dq
\end{eqnarray*}

Consider now the Gram matrix in the statement, namely:
$$A(q)_{ij}=<R_i,R_j>$$

Here $R_1,\ldots,R_M$ are the rows of the following matrix:
$$q\in \mathbb T^{MN}\simeq M_{M\times N}(\mathbb T)$$

We have then the following computation:
\begin{eqnarray*}
\int\chi^p
&=&\frac{1}{MN}\int_{\mathbb T^{MN}}<R_{i_1},R_{i_2}><R_{i_2},R_{i_3}>\ldots<R_{i_p},R_{i_1}>\\
&=&\frac{1}{MN}\int_{\mathbb T^{MN}}A(q)_{i_1i_2}A(q)_{i_2i_3}\ldots A(q)_{i_pi_1}\\
&=&\frac{1}{MN}\int_{\mathbb T^{MN}}Tr(A(q)^p)dq\\
&=&\frac{1}{N}\int_{\mathbb T^{MN}}tr(A(q)^p)dq
\end{eqnarray*}

But this gives the formula in the statement, and we are done.
\end{proof}

In general, the moments of the Gram matrix $A$ are given by a quite complicated formula, and we cannot expect to have a refinement of Theorem 16.20, with $A$ replaced by a plain, non-matricial random variable, say over a compact abelian group. 

\bigskip

However, this kind of simplification does appear at $M=2$, and since this phenomenon is quite interesting, we will explain this now. We first have:

\begin{proposition}
For $F_2\otimes_QF_H$, with $Q\in M_{2\times N}(\mathbb T)$ generic, we have
$$N\int\left(\frac{\chi}{N}\right)^p=\int_{\mathbb T^N}\sum_{k\geq0}\binom{p}{2k}\left|\frac{a_1+\ldots+a_N}{N}\right|^{2k}da$$
where the integral on the right is with respect to the uniform measure on $\mathbb T^N$.
\end{proposition}

\begin{proof}
In order to prove the result, consider the following quantity, which appeared in the proof of Theorem 16.20:
$$\Phi(q)=\sum_{i_1\ldots i_p}\sum_{d_1\ldots d_p}\frac{q_{i_1d_1}\ldots q_{i_pd_p}}{q_{i_1d_p}\ldots q_{i_pd_{p-1}}}$$

We can ``half-dephase'' the matrix $q\in M_{2\times N}(\mathbb T)$ if we want to, as follows:
$$q=\begin{pmatrix}1&\ldots&1\\ a_1&\ldots&a_N\end{pmatrix}$$

Let us compute now the above quantity $\Phi(q)$, in terms of the numbers $a_1,\ldots,a_N$. Our claim is that we have the following formula:
$$\Phi(q)=2\sum_{k\geq0}N^{p-2k}\binom{p}{2k}\left|\sum_ia_i\right|^{2k}$$

Indeed, the idea is that:

\medskip

(1) The $2N^k$ contribution will come from $i=(1\ldots1)$ and $i=(2\ldots2)$.

\medskip

(2) Then we will have a $p(p-1)N^{k-2}|\sum_ia_i|^2$ contribution coming from indices of type $i=(2\ldots 21\ldots1)$, up to cyclic permutations.

\medskip

(3) Then we will have a $2\binom{p}{4}N^{p-4}|\sum_ia_i|^4$ contribution coming from indices of type $i=(2\ldots 21\ldots12\ldots21\ldots1)$.

\medskip

(4) And so on. 

\medskip

In practice now, in order to prove our claim, in order to find the $N^{p-2k}|\sum_ia_i|^{2k}$ contribution, we have to count the circular configurations consisting of $p$ numbers $1,2$, such that the $1$ values are arranged into $k$ non-empty intervals, and the $2$ values are arranged into $k$ non-empty intervals as well. Now by looking at the endpoints of these $2k$ intervals, we have $2\binom{p}{2k}$ choices, and this gives the above formula.

\medskip

Now by integrating, this gives the formula in the statement.
\end{proof}

Observe now that the integrals in Proposition 16.21 can be computed as follows:
\begin{eqnarray*}
\int_{\mathbb T^N}|a_1+\ldots+a_N|^{2k}da
&=&\int_{\mathbb T^N}\sum_{i_1\ldots i_k}\sum_{j_1\ldots j_k}\frac{a_{i_1}\ldots a_{i_k}}{a_{j_1}\ldots a_{j_k}}da\\
&=&\#\left\{i_1\ldots i_k,j_1\ldots j_k\Big|[i_1,\ldots,i_k]=[j_1,\ldots,j_k]\right\}\\
&=&\sum_{k=\sum r_i}\binom{k}{r_1,\ldots,r_N}^2
\end{eqnarray*}

We obtain in this way the following ``blowup'' result, for our measure:

\index{blowup}

\begin{proposition}
For $F_2\otimes_QF_H$, with $Q\in M_{2\times N}(\mathbb T)$ generic, we have
$$\mu=\left(1-\frac{1}{N}\right)\delta_0+\frac{1}{2N}\left(\Psi^+_*\varepsilon+\Psi^-_*\varepsilon\right)$$
where $\varepsilon$ is the uniform measure on $\mathbb T^N$, and where the blowup function is:
$$\Psi^\pm(a)=N\pm\left|\sum_ia_i\right|$$
\end{proposition}

\begin{proof}
We use the formula found in Proposition 16.21, along with the following standard identity, coming from the Taylor formula:
$$\sum_{k\geq0}\binom{p}{2k}x^{2k}=\frac{(1+x)^p+(1-x)^p}{2}$$

By using this identity, Proposition 16.20 reformulates as follows:
$$N\int\left(\frac{\chi}{N}\right)^p=\frac{1}{2}\int_{\mathbb T^N}\left(1+\left|\frac{\sum_ia_i}{N}\right|\right)^p+\left(1-\left|\frac{\sum_ia_i}{N}\right|\right)^p\,da$$

Now by multiplying by $N^{p-1}$, we obtain the following formula:
$$\int\chi^k=\frac{1}{2N}\int_{\mathbb T^N}\left(N+\left|\sum_ia_i\right|\right)^p+\left(N-\left|\sum_ia_i\right|\right)^p\,da$$

But this gives the formula in the statement, and we are done.
\end{proof}

We can further improve the above result, by reducing the maps $\Psi^\pm$ appearing there to a single one, and we are led to the following statement: 

\begin{theorem}
For $F_2\otimes_QF_H$, with $Q\in M_{2\times N}(\mathbb T)$ generic, we have
$$\mu=\left(1-\frac{1}{N}\right)\delta_0+\frac{1}{N}\Phi_*\varepsilon$$
where $\varepsilon$ is the uniform measure on $\mathbb Z_2\times\mathbb T^N$, and where the blowup map is:
$$\Phi(e,a)=N+e\left|\sum_ia_i\right|$$
\end{theorem}

\begin{proof}
This is clear indeed from Proposition 16.22.
\end{proof}

As already mentioned, the above results at $M=2$ are something quite special. In the general case, $M\in\mathbb N$, it is not clear how to construct a nice blowup of the measure.

\section*{16d. Poisson laws}

Asymptotically, things are however quite simple. Let us go back indeed to the general case, where $M,N\in\mathbb N$ are both arbitrary. The problem that we would like to solve now is that of finding the good regime, of the following type, where the measure in Theorem 16.19 converges, after some suitable manipulations:
$$M=f(K)\quad,\quad 
N=g(K)\quad,\quad 
K\to\infty$$

In order to do so, we have to do some combinatorics. Let $NC(p)$ be the set of noncrossing partitions of $\{1,\ldots,p\}$, and for $\pi\in P(p)$ we denote by $|\pi|\in\{1,\ldots,p\}$ the number of blocks. With these conventions, we have the following result:

\begin{proposition}
With $M=\alpha K,N=\beta K$, $K\to\infty$ we have:
$$\frac{c_p}{K^{p-1}}\simeq\sum_{r=1}^p\#\left\{\pi\in NC(p)\Big||\pi|=r\right\}\alpha^{r-1}\beta^{p-r}$$
In particular, with $\alpha=\beta$ we have:
$$c_p\simeq\frac{1}{p+1}\binom{2p}{p}(\alpha K)^{p-1}$$
\end{proposition}

\begin{proof}
We use the combinatorial formula in Theorem 16.19. Our claim is that, with $\pi=\ker(i_1,\ldots,i_p)$, the corresponding contribution to $c_p$ is:
$$C_\pi\simeq
\begin{cases}
\alpha^{|\pi|-1}\beta^{p-|\pi|}K^{p-1}&{\rm if}\ \pi\in NC(p)\\
O(K^{p-2})&{\rm if}\ \pi\notin NC(p)
\end{cases}$$

As a first observation, the number of choices for a multi-index $(i_1,\ldots,i_p)\in X^p$ satisfying $\ker i=\pi$ is:
$$M(M-1)\ldots (M-|\pi|+1)\simeq M^{|\pi|}$$

Thus, we have the following estimate:
$$C_\pi\simeq M^{|\pi|-1}N^{-1}\#\left\{d_1,\ldots,d_p\in Y\Big|[d_\alpha|\alpha\in b]=[d_{\alpha-1}|\alpha\in b],\forall b\in\pi\right\}$$

Consider now the following partition:
$$\sigma=\ker d$$

The contribution of $\sigma$ to the above quantity $C_\pi$ is then given by:
$$\Delta(\pi,\sigma)N(N-1)\ldots(N-|\sigma|+1)\simeq\Delta(\pi,\sigma)N^{|\sigma|}$$

Here the quantities on the right are as follows:
$$\Delta(\pi,\sigma)=\begin{cases}
1&{\rm if}\ |b\cap c|=|(b-1)\cap c|,\forall b\in\pi,\forall c\in\sigma\\
0&{\rm otherwise}
\end{cases}$$

We use now the standard fact that for $\pi,\sigma\in P(p)$ satisfying $\Delta(\pi,\sigma)=1$ we have:
$$|\pi|+|\sigma|\leq p+1$$

In addition, the equality case is well-known to happen when $\pi,\sigma\in NC(p)$ are inverse to each other, via Kreweras complementation. This shows that for $\pi\notin NC(p)$ we have:
$$C_\pi=O(K^{p-2})$$

Also, this shows that for $\pi\in NC(p)$ we have:
\begin{eqnarray*}
C_\pi
&\simeq&M^{|\pi|-1}N^{-1}N^{p-|\pi|-1}\\
&=&\alpha^{|\pi|-1}\beta^{p-|\pi|}K^{p-1}
\end{eqnarray*}

Thus, we have obtained the result.
\end{proof}

We denote by $D$ the dilation operation, given by:
$$D_r(law(X))=law(rX)$$

With this convention, we have the following result:

\begin{theorem}
With $M=\alpha K,N=\beta K$, $K\to\infty$ we have:
$$\mu=\left(1-\frac{1}{\alpha\beta K^2}\right)\delta_0+\frac{1}{\alpha\beta K^2}D_{\frac{1}{\beta K}}(\pi_{\alpha/\beta})$$
In particular with $\alpha=\beta$ we have:
$$\mu=\left(1-\frac{1}{\alpha^2K^2}\right)\delta_0+\frac{1}{\alpha^2K^2}D_{\frac{1}{\alpha K}}(\pi_1)$$
\end{theorem}

\begin{proof}
At $\alpha=\beta$, this follows from Proposition 16.24. In general now, we have:
\begin{eqnarray*}
\frac{c_p}{K^{p-1}}
&\simeq&\sum_{\pi\in NC(p)}\alpha^{|\pi|-1}\beta^{p-|\pi|}\\
&=&\frac{\beta^p}{\alpha}\sum_{\pi\in NC(p)}\left(\frac{\alpha}{\beta}\right)^{|\pi|}\\
&=&\frac{\beta^p}{\alpha}\int x^pd\pi_{\alpha/\beta}(x)
\end{eqnarray*}

When $\alpha\geq\beta$, where $d\pi_{\alpha/\beta}(x)=\varphi_{\alpha/\beta}(x)dx$ is continuous, we obtain:
\begin{eqnarray*}
c_p
&=&\frac{1}{\alpha K}\int(\beta Kx)^p\varphi_{\alpha/\beta}(x)dx\\
&=&\frac{1}{\alpha\beta K^2}\int x^p\varphi_{\alpha/\beta}\left(\frac{x}{\beta K}\right)dx
\end{eqnarray*}

But this gives the formula in the statement. When $\alpha\leq\beta$ the computation is similar, with a Dirac mass as 0 dissapearing and reappearing, and gives the same result.
\end{proof}

Let us state as well an explicit result, regarding densities:

\begin{theorem}
With $M=\alpha K,N=\beta K$, $K\to\infty$ we have:
$$\mu=\left(1-\frac{1}{\alpha\beta K^2}\right)\delta_0+\frac{1}{\alpha\beta K^2}\cdot\frac{\sqrt{4\alpha\beta K^2-(x-\alpha K-\beta K)^2}}{2\pi x}\,dx$$
In particular with $\alpha=\beta$ we have:
$$\mu=\left(1-\frac{1}{\alpha^2K^2}\right)\delta_0+\frac{1}{\alpha^2K^2}\cdot\frac{\sqrt{\frac{4\alpha K}{x}-1}}{2\pi}$$
\end{theorem}

\begin{proof}
According to the formula for the density of the free Poisson law, the density of the continuous part $D_{\frac{1}{\beta K}}(\pi_{\alpha/\beta})$ is indeed given by:
$$\frac{\sqrt{4\frac{\alpha}{\beta}-(\frac{x}{\beta K}-1-\frac{\alpha}{\beta})^2}}
{2\pi\cdot\frac{x}{\beta K}}=\frac{\sqrt{4\alpha\beta K^2-(x-\alpha K-\beta K)^2}}{2\pi x}$$

With $\alpha=\beta$ now, we obtain the second formula in the statement, and we are done.
\end{proof}

Observe that at $\alpha=\beta=1$, where $M=N=K\to\infty$, the above measure is:
$$\mu=\left(1-\frac{1}{K^2}\right)\delta_0+\frac{1}{K^2}D_{\frac{1}{K}}(\pi_1)$$

This measure is supported by $[0,4K]$. On the other hand, since the groups $\Gamma_{M,N}$ are all amenable, the corresponding measures are supported on $[0,MN]$, and so on $[0,K^2]$ in the $M=N=K$ situation. The fact that we do not have a convergence of supports is not surprising, because our convergence is in moments.

\bigskip

The above results are of course not the end of the story, because we have now to understand what happens in the case of non-generic parameters. There has been some technical work here, by Bichon and by myself, and as a sample result here, we have:

\index{free Poisson law}
\index{Marchenko-Pastur law}

\begin{theorem}
Given two finite abelian groups $G,H$, having cardinalities 
$$|G|=M\quad,\quad 
|H|=N$$ 
consider the main character $\chi$ of the quantum group associated to $\mathcal F_{G\times H}$. We have then
$$law\left(\frac{\chi}{N}\right)=\left(1-\frac{1}{M}\right)\delta_0+\frac{1}{M}\,\pi_t$$
in moments, with $M=tN\to\infty$, where $\pi_t$ is the free Poisson law of parameter $t>0$. In addition, this formula holds for any generic fiber of $\mathcal F_{G\times H}$.
\end{theorem}

\begin{proof}
We already know that the second assertion holds, as explained above. 

\bigskip

Regarding now the first assertion, our first claim is that for the representation coming from the parametric matrix $\mathcal F_{G\times H}$ we have the following formula, where $M=|G|,N=|H|$, and the sets between brackets are sets with repetitions:
$$c_p^r=\frac{1}{M^{r+1}N}\#\left\{\begin{matrix}i_1,\ldots,i_r,a_1,\ldots,a_p\in\{0,\ldots,M-1\},\\
b_1,\ldots,b_p\in\{0,\ldots,N-1\},\\
[(i_x+a_y,b_y),(i_{x+1}+a_y,b_{y+1})|y=1,\ldots,p]\\
=[(i_x+a_y,b_{y+1}),(i_{x+1}+a_y,b_y)|y=1,\ldots,p], \forall x
\end{matrix}\right\}$$

Indeed, by using the general moment formula with $K=F_G$, $L=F_H$, we have:
\begin{eqnarray*}
&&c_p^r\\
&=&\frac{1}{(MN)^r}\int_{T^r}\sum_{i_1^1\ldots i_p^r}\sum_{b_1^1\ldots b_p^r}\frac{Q^1_{i_1^1b_1^1}Q^1_{i_1^2b_2^1}}{Q^1_{i_1^1b_2^1}Q^1_{i_1^2b_1^1}}\ldots\frac{Q^1_{i_p^1b_p^1}Q^1_{i_p^2b_1^1}}{Q^1_{i_p^1b_1^1}Q^1_{i_p^2b_p^1}}\ldots\ldots\frac{Q^r_{i_1^rb_1^r}Q^r_{i_1^1b_2^r}}{Q^r_{i_1^rb_2^r}Q^r_{i_1^1b_1^r}}\ldots\frac{Q^r_{i_p^rb_p^r}Q^r_{i_p^1b_1^r}}{Q^r_{i_p^rb_1^r}Q^r_{i_p^1b_p^r}}\\
&&\hskip15mm\frac{1}{M^{pr}}\sum_{j_1^1\ldots j_p^r}\frac{K_{i_1^1j_1^1}K_{i_1^2j_2^1}}{K_{i_1^1j_2^1}K_{i_1^2j_1^1}}\ldots\frac{K_{i_p^1j_p^1}K_{i_p^2j_1^1}}{K_{i_p^1j_1^1}K_{i_p^2j_p^1}}\ldots\ldots\frac{K_{i_1^rj_1^r}K_{i_1^1j_2^r}}{K_{i_1^rj_2^r}K_{i_1^1j_1^r}}\ldots\frac{K_{i_p^rj_p^r}K_{i_p^1j_1^r}}{K_{i_p^rj_1^r}K_{i_p^1j_p^r}}\\
&&\hskip15mm\frac{1}{N^{pr}}\sum_{a_1^1\ldots a_p^r}\frac{L_{a_1^1b_1^1}L_{a_1^2b_2^1}}{L_{a_1^1b_2^1}L_{a_1^2b_1^1}}\ldots\frac{L_{a_p^1b_p^1}L_{a_p^2b_1^1}}{L_{a_p^1b_1^1}L_{a_p^2b_p^1}}\ldots\ldots\frac{L_{a_1^rb_1^r}L_{a_1^1b_2^r}}{L_{a_1^rb_2^r}L_{a_1^1b_1^r}}\ldots\frac{L_{a_p^rb_p^r}L_{a_p^1b_1^r}}{L_{a_p^rb_1^r}L_{a_p^1b_p^r}}\,dQ\\
\end{eqnarray*}

Since we are in the Fourier matrix case, $K=F_G,L=F_H$, we can perform the sums over $j,a$. To be more precise, the last two averages appearing above are respectively:
\begin{eqnarray*}
\Delta(i)&=&\prod_x\prod_y\delta(i^x_y+i^{x+1}_{y-1},i^{x+1}_y+i^x_{y-1})\\
\Delta(b)&=&\prod_x\prod_y\delta(b^x_y+b^{x+1}_{y-1},b^{x+1}_y+b^x_{y-1})
\end{eqnarray*}

We therefore obtain the following formula for the truncated moments of the main character, where $\Delta$ is the product of Kronecker symbols constructed above:
\begin{eqnarray*}
&&c_p^r\\
&=&\frac{1}{(MN)^r}\int_{T^r}\sum_{\Delta(i)=\Delta(b)=1}\frac{Q^1_{i_1^1b_1^1}Q^1_{i_1^2b_2^1}}{Q^1_{i_1^1b_2^1}Q^1_{i_1^2b_1^1}}\ldots\frac{Q^1_{i_p^1b_p^1}Q^1_{i_p^2b_1^1}}{Q^1_{i_p^1b_1^1}Q^1_{i_p^2b_p^1}}\ldots\ldots\frac{Q^r_{i_1^rb_1^r}Q^r_{i_1^1b_2^r}}{Q^r_{i_1^rb_2^r}Q^r_{i_1^1b_1^r}}\ldots\frac{Q^r_{i_p^rb_p^r}Q^r_{i_p^1b_1^r}}{Q^r_{i_p^rb_1^r}Q^r_{i_p^1b_p^r}}\,dQ
\end{eqnarray*}

Now by integrating with respect to $Q\in(\mathbb T^{G\times H})^r$, we are led to counting the multi-indices $i,b$ satisfying several conditions. First, we have the following condition:
$$\Delta(i)=\Delta(b)=1$$

We have as well the following conditions, where the sets between brackets are by definition sets with repetitions:
$$\begin{bmatrix}
i_1^1b_1^1&\ldots&i_p^1b_p^1&i_1^2b_2^1&\ldots&i_p^2b_1^1
\end{bmatrix}=
\begin{bmatrix}
i_1^1b_2^1&\ldots&i_p^1b_1^1&i_1^2b_1^1&\ldots&i_p^2b_p^1
\end{bmatrix}$$
$$\vdots$$
$$\begin{bmatrix}
i_1^rb_1^r&\ldots&i_p^rb_p^r&i_1^1b_2^r&\ldots&i_p^1b_1^r
\end{bmatrix}
=\begin{bmatrix}
i_1^rb_2^r&\ldots&i_p^rb_1^r&i_1^1b_1^r&\ldots&i_p^1b_p^r
\end{bmatrix}$$

In a more compact notation, the moment formula that we obtain in this way is therefore as follows:
$$c_p^r=\frac{1}{(MN)^r}\#\left\{i,b\Big|\Delta(i)=\Delta(b)=1,\ [i^x_yb^x_y,i^{x+1}_yb^x_{y+1}]=[i^x_yb^x_{y+1},i^{x+1}_yb^x_y],\forall x\right\}$$

Now observe that the above Kronecker type conditions $\Delta(i)=\Delta(b)=1$ tell us that the arrays of indices $i=(i^x_y),b=(b^x_y)$ must be of the following special form:
$$\begin{pmatrix}i^1_1&\ldots&i^1_p\\&\ldots\\ i^1_r&\ldots&i^r_p\end{pmatrix}=\begin{pmatrix}i_1+a_1&\ldots&i_1+a_p\\&\ldots\\ i_r+a_1&\ldots&i_r+a_p\end{pmatrix}$$

$$\begin{pmatrix}b^1_1&\ldots&b^1_p\\&\ldots\\ b^1_r&\ldots&b^r_p\end{pmatrix}=\begin{pmatrix}j_1+b_1&\ldots&j_1+b_p\\&\ldots\\ j_r+b_1&\ldots&j_r+b_p\end{pmatrix}$$

Here all the new indices $i_x,j_x,a_y,b_y$ are uniquely determined, up to a choice of $i_1,j_1$. Now by replacing $i^x_y,b^x_y$ with these new indices $i_x,j_x,a_y,b_y$, with a $MN$ factor added, which accounts for the choice of $i_1,j_1$, we obtain the following formula:
$$c_p^r=\frac{1}{(MN)^{r+1}}\#\left\{i,j,a,b\Big|\begin{matrix}[(i_x+a_y,j_x+b_y),(i_{x+1}+a_y,j_x+b_{y+1})]\\
=[(i_x+a_y,j_x+b_{y+1}),(i_{x+1}+a_y,j_x+b_y)],\forall x\end{matrix}\right\}$$

Now observe that we can delete if we want the $j_x$ indices, which are irrelevant. Thus, we obtain the announced formula. The continuation is via standard combinatorics.
\end{proof}

There are many interesting questions that are still open, regarding the computation of the spectral measure in the case where the parameter matrix $Q$ is not generic, and also regarding the computation for the deformations of the generalized Fourier matrices, which are not necessarily of Di\c t\u a type. We refer here to \cite{bb3} and related papers.

\section*{16e. Exercises} 

Congratulations for having read this book, and no exercises here. However, there are many interesting questions, as a continuation of the above, some of which were raised during our exercise sessions in this book, some other that you will find by reading the research articles referenced below, and yet some other that you will have to invent yourself. Enjoy all this, and in the hope that we will hear soon from you, new exciting things regarding $S_N^+$ and its subgroups, in relation with mathematics and physics.

\baselineskip=14pt

\printindex

\end{document}